%% file: OTERN.tex
\documentclass[reqno]{erinbook}

\usepackage{general,resistance,ifthen,text_shortcuts}%
\usepackage{cancel,xy}
\usepackage[switch*,pagewise]{lineno}
\usepackage[style=list]{glossary}
  \selectfont

\usepackage[bookmarks, colorlinks=true, pdfstartview=FitV, linkcolor=blue, citecolor=blue, urlcolor=blue]{hyperref}
\usepackage{cite}  

\newcommand{\clean}{true}  
\newcommand{\version}[2]{\ifthenelse{\equal{\clean}{true}}{#1}{{\footnotesize #2}}}
\newcommand{\linenopax}{\par}  

\newcommand{\headerquote}[2]{\scalebox{0.80}{\hstr[8]\begin{minipage}{5in}\begin{centering}{\textit{``#1'' #2}}\\ \vstr\end{centering}\end{minipage}}}

\newcommand{\q}{\quad}
\newcommand{\qq}{\qquad}

\numberwithin{equation}{chapter} \numberwithin{theorem}{chapter}

\begin{document}

  \include{title}

  \frontmatter
  
  \linenumbers
  \pagewiselinenumbers

\title{Operator theory and analysis of infinite networks}

\author{Palle E. T. Jorgensen and Erin P. J. Pearse}




\version{\vfill \eject \setcounter{tocdepth}{2}}{\setcounter{tocdepth}{2}}
{\small \tableofcontents}

\allowdisplaybreaks

\include{preface}
\include{introduction}

\mainmatter

\include{electrical-resistance-networks}

\include{potentials}

\include{energy-hilbert-space}

\include{resistance-metrics}

\include{construction-of-HE}

\include{boundary}

\include{Lap-on-HE}
\include{ell2-of-Lap-and-Trans}

\include{he-and-hd}

\include{probab-interp}

\include{examples}
\include{tree-examples}

\include{lattice-examples}

\include{magnetism}

\include{future-directions}


\begin{appendix}
\include{functional-analysis}

\include{operator-theory}

\include{road-map}

\include{biblioguide}

\end{appendix}

\addcontentsline{toc}{chapter}{List of symbols and notation}
\headerquote{Mathematics is a game played according to certain simple rules with meaningless marks on paper.}{---~D.~Hilbert}
We also attempt consistency in denoting vertices by $x,y,z$; functions on vertices by $u,v,w$; functions on edges by $\curr, J$, and denoting the beginning and end of a finite path by \ga and \gw, respectively.

\bibliographystyle{alpha}
\bibliography{OTERN}
\addcontentsline{toc}{chapter}{References}
\nocite{*} 


\end{document}

%% file: title.tex

\thispagestyle{empty}

\begin{center}
 \textsc{\Huge{{\fontencoding{OT1}\fontfamily{ppl}\fontseries{bx} \selectfont{Operator theory \\ and analysis of \\ infinite networks}}} 
  \\ \vstr[8] \\
  \scalebox{1.3}{\includegraphics{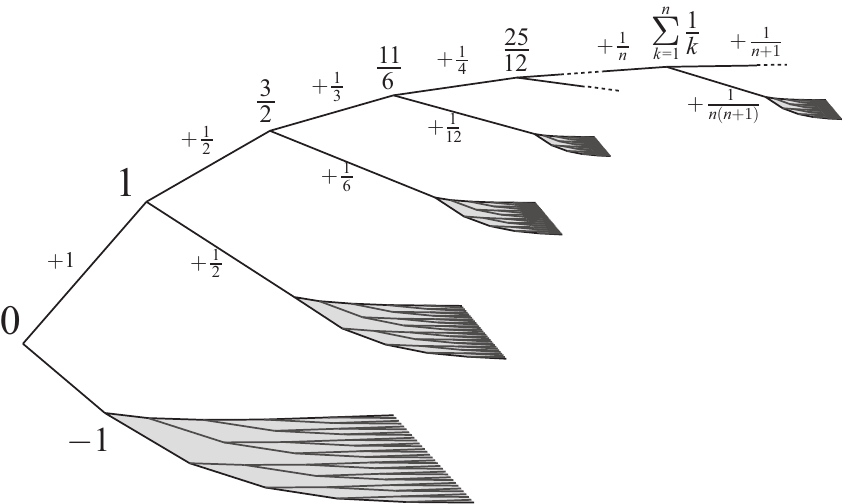}}
  \\ \vstr[8] \\
  \LARGE{\textsc{\fontencoding{OT1}\fontfamily{ppl} \selectfont Palle E. T. Jorgensen and \\ Erin P. J. Pearse}}\vstr[3]}
\end{center}




\newpage
\thispagestyle{empty}

\cleardoublepage

%% file: preface.tex

\chapter*{Preface}
\addcontentsline{toc}{section}{Preface}

\headerquote{We often hear that mathematics consists mainly of "proving theorems." Is a writer's job mainly that of ``writing sentences?''}{---~G.-C.~Rota}

\headerquote{It can be shown that a mathematical web of some kind can be woven about any universe containing several objects. The fact that our universe lends itself to mathematical treatment is not a fact of any great philosophical significance.}{---~B.~ Russell}

\section*{Note to the reader}
\addcontentsline{toc}{subsection}{Note to the reader}

  In this book, we wish to present operators in Hilbert space (with an emphasis on the theory of unbounded operators) from the vantage point of a relatively new trend, the analysis of infinite networks. This in turn involves such hands-on applications as infinite systems of resistors, and random walk on infinite graphs. Other such ``infinite'' systems include mathematical models of the internet. This new tapestry of applications offers a special appeal, and has the further advantage of bringing into play additional tools from both probability and metric geometry.
  
  While we have included some fundamentals of operator theory in the Appendices, readers will first be treated to the fundamentals of infinite networks and their operator theory. Throughout the exposition, we will make continual use of the axioms of Hilbert space, and such standard tools as the Schwarz inequality, Riesz's Lemma, projections, and the lattice of subspaces, all of which are available in any introductory functional analysis book. Readers not already familiar with this material may wish to consult the Appendices for the axioms of Hilbert space, embedding and isomorphism theorems (Appendix~\ref{sec:functional-analysis}), bounded and unbounded linear operators, the geometry of projections, infinite banded matrices, Hermitian and selfadjoint operators with dense domain, adjoint operators and their graphs, deficiency indices (Appendix~\ref{sec:operator-theory}).

Some material is motivated by deeper aspects of Hilbert space theory: the Gel'fand triple construction of Chapter~\ref{sec:the-boundary}, deficiency indices of unbounded operators in Chapter~\ref{sec:Lap-on-HE}, parallels between Kolmogorov consistency and the GNS construction in Chapter~\ref{sec:Magnetism-and-long-range-order}, and the relation of KMS states to long-range order in Chapter~\ref{sec:Magnetism-and-long-range-order}. Familiarity with these topics are not a prerequisite for this book! Conversely, we hope that the present setting allows for a smooth introduction to these areas (which may otherwise be dauntingly technical) and have correspondingly provided extensive introductory material at the relevant locations in the text.

By using the intrinsic inner product (associated to the effective resistance) we are able to obtain results which are more physically realistic than many found elsewhere in the literature. This inner product is quite different than the standard $\ell^2$ inner product for functions defined on the vertices of a graph, and holds many surprises. Many of our results apply much more generally than those already present in the literature. The next section elaborates on these rather vague remarks and highlights the advantages and differences inherent in our approach, in a variety of circumstances.

This work is uniquely interdisciplinary, and as a consequence, we have made effort to address the \emph{union} (as opposed to the intersection) of several disparate audiences: graph theory, resistance networks, spectral geometry, fractal geometry, physics, probability, unbounded operators in Hilbert space, C*-algebras, and others. It is inevitable that parts of the background material there will be unknown to some readers and so we have included the appendices to mediate this. After presenting our results at various talks, we felt that the inclusion of this material would be appreciated by most.

  The subject of operator theory enjoys periodic bursts of renewed  interest and progress, and often because of impulses and inspiration from neighboring fields. We feel that these recent trends and interconnections in discrete mathematics are ready for a self-contained presentation; a presentation we hope will help both students and researchers gain access to operator theory as well as some of its more exciting applications.

  The literature on Hilbert space and linear operators frequently breaks into a dichotomy: axiomatic vs. applied. In this book, we aim at linking the sides: after introducing a set of axioms and using them to prove some theorems, we provide examples with explicit computations. 
  For any application, there may be a host of messy choices for inner product, and often only one of them is right (despite the presence of some axiomatic isomorphisms).
  
  The most famous example of a nontrivial such isomorphism stems from the birth of quantum theory. The matrix model of Werner Heisenberg was in fierce competition with the PDE model of Erwin Schr\"{o}dinger until John von Neumann ended the dispute in 1932 by proving that the two Hilbert space models are in fact unitarily equivalent. However, despite the presence of such an axiomatic equivalence, one must still do computations in whichever one of the two models offers solutions to problems in the laboratory.

\section*{Brief overview of contents}
\addcontentsline{toc}{subsection}{Brief overview of contents}

\headerquote{Therefore, either the reality on which our space is based must form a discrete manifold or else the reason for the metric relationships must be sought for, externally, in the binding forces acting on it.}{---~G.~F.~B.~Riemann}

Among subjects in mathematics, Functional Analysis and Operator Theory are special in several respects, they are relatively young (measured in the historical scale of mathematics), and they often have a more interdisciplinary flavor. While the axiomatic side of the subjects has matured, there continues to be an inexhaustible supply of exciting applications. We focus here on a circle of interdisciplinary areas: weighted graphs and their analysis. The infinite cases are those that involve operators in Hilbert space and entail potential theory, metric geometry, probability, and harmonic analysis.  Of the infinite weighted graphs, some can be modeled successfully as systems of resistors, but the resulting mathematics has much wider implications. Below we sketch some main concepts from resistance networks.

\pgap

The following rather terse/dense sequence of paragraphs is an abstract for the reader who wishes to get an idea of the contents after just reading a page or so. A more detailed description is given in the Introduction just below.

A resistance network is a weighted graph $(G,c)$. The conductance function $c_{xy}$ weights the edges, which are then interpreted as resistors of possibly varying strengths. The effective resistance metric $R(x,y)$ is the natural notion of distance between two vertices $x,y$ in the resistance network. 

The space of functions of finite energy (modulo constants) is a Hilbert space with inner product $\mathcal E$, which we call the energy space ${\mathcal H}_{\mathcal E}$. The evaluation functionals on ${\mathcal H}_{\mathcal E}$ give rise to a reproducing kernel $\{v_x\}$ for the space. Once a reference vertex $o$ is fixed, these functions $v_x$ satisfy $\Lap v_x = \gd_x - \gd_o$, where $\Delta$ is the network Laplacian. This kernel yields a detailed description of the structure of ${\mathcal H}_{\mathcal E} = \mathcal{F}in \oplus \mathcal{H}arm$, where $\mathcal{F}in$ is the closure of the space of finitely supported functions and $\mathcal{H}arm$ is the closed subspace of harmonic functions. The energy $\mathcal{E}$ splits accordingly into a ``finite part'' expressed as a sum taken over the vertices, and an ``infinite part'' expressed as a limit of sums. Intuitively, the latter part corresponds to an integral over some sort of boundary $\textrm{bd} G$, which is developed explicitly in \S\ref{sec:Lap-on-HE}. The kernel $\{v_x\}$ also allows us to recover easily many known (and sometimes difficult) results about ${\mathcal H}_{\mathcal E}$.
As ${\mathcal H}_{\mathcal E}$ does not come naturally equipped with a natural o.n.b., we provide candidates for frames (and dual frames) when working with an infinite resistance network.

In particular, the presence of nonconstant harmonic functions of finite energy leads to 
different plausible definitions of the effective resistance metric on infinite networks. We characterize the free resistance $R^F(x,y)$ and the wired resistance $R^W(x,y)$ in terms of Neumann or Dirichlet boundary conditions on a certain operator. (In the literature, these correspond to the \emph{limit current} and \emph{minimal current}, resp.) We develop a library of equivalent formulations for each version. 
Also, we introduce the ``trace resistance'' $R^S(x,y)$, computed in terms of the trace of the Dirichlet form $\mathcal E$ to finite subnetworks. This provides a finite approximation which is more accurate from a probabilistic perspective, and gives a probabilistic explanation of the discrepancy between $R^F$ and $R^W$.

For $R=R^F$ or $R=R^W$, the effective resistance is shown to be negative semidefinite, so that it induces an inner product on a Hilbert space into which it naturally embeds. We show that for $(G,R^F)$, the resulting Hilbert space is ${\mathcal H}_{\mathcal E}$ and for $(G,R^W)$ it is $\mathcal{F}in$. Under the free embedding, each vertex $x$ is mapped to the element $v_x$ of the energy kernel; under the wired embedding it is mapped to the projection $f_x$ of $v_x$ to $\mathcal{F}in$. This establishes ${\mathcal H}_{\mathcal E}$ as the natural Hilbert space in which to study effective resistance.

We obtain an analytic boundary representation for elements of $\mathcal{H}arm$ in a sense analogous to that of Poisson or Martin boundary theory. We construct a Gel'fand triple $S \ci {\mathcal H}_{\mathcal E} \ci S'$ and obtain a probability measure $\mathbb{P}$ and an isometric embedding of ${\mathcal H}_{\mathcal E}$ into $L^2(S',\mathbb{P})$. This gives a concrete representation of the boundary in terms of the measures $(\mathbf{1} + v_{x_n})d\mathbb{P} \in S'/\mathcal{F}in$, where $\{x_n\}$ is a sequence tending to infinity.


The spectral representation for the graph Laplacian $\Delta$ on ${\mathcal H}_{\mathcal E}$ is \emph{drastically} different from the corresponding representation on $\ell^2$. Since the ambient Hilbert space ${\mathcal H}_{\mathcal E}$ is defined by the energy form, many interesting phenomena arise which are not present in $\ell^2$; we highlight many examples and explain why this occurs. In particular, we show how the deficiency indices of $\Delta$ as an operator on ${\mathcal H}_{\mathcal E}$ indicate the presence of nontrivial boundary of an resistance network, and why the $\ell^2$ operator theory of $\Delta$ does not see this. Along the way, we prove that $\Delta$ is always essentially self-adjoint on the $\ell^2$ space of functions on an resistance network, and examine conditions for the network Laplacian and its associated transfer operator to be bounded, compact, essential self-adjoint, etc. 

We consider two approaches to measures on spaces of infinite paths in an resistance network. One arises from considering the transition probabilities of a random walk as determined directly by the network, i.e., $p(x,y) = c_{xy}/\sum_{y \sim x} c_{xy}$. The other applies only to transient networks, and arises from considering the transition probabilities induced by a unit flow to infinity. The latter leads to the notion of forward-harmonic functions, for which we also provide a characterization in terms of a boundary representation.

Using our results we establish precise bounds on correlations in the Heisenberg model for quantum spin observables, and we improve earlier results of R. T. Powers. Our focus is on the quantum spin model on the rank-3 lattice, i.e., the resistance network with ${\mathbb Z}^3$ as vertices and with edges between nearest neighbors. This is known as the problem of long-range order in the physics literature, and refers to KMS states on the $C^\ast$-algebra of the model.

%% file: introduction.tex

\chapter*{Introduction}
\addcontentsline{toc}{section}{Introduction}

\headerquote{... an apt comment on how science, and indeed the whole of civilization, is a series of incremental advances, each building on what went before.}{---~Stephen~Hawking}

The subject of resistance networks has its origins in electrical engineering applications, and over decades, it has served to motivate a number of advances in discrete mathematics, such as the study of boundaries, percolation, stochastic analysis and random walk on graphs. There are already several successful schools of research, each with its own striking scientific advances, and it may be a little premature attempting to summarize the vast variety of new theorems. They are still appearing at a rapid rate in research journals!

A common theme in the study of boundaries on infinite discrete systems $X$ (weighted graphs, trees, Markov chains, or discrete groups) is the focus on a suitable subspace of functions on $X$, usually functions which are harmonic in some sense (i.e., fixed points of a given transfer operator). We are interested in the harmonic functions of finite energy, as this class of harmonic functions comes equipped with a natural inner product and corresponding Hilbert space structure. This will guide our choice of topics and emphasis, from an otherwise vast selection of possibilities.

This volume is dedicated to the construction of unified functional-analytic framework for the study of these potential-theoretic function spaces on graphs, and an investigation of the resulting structures. The primary object of study is a resistance network: a graph with weighted edges. Our foundation is the effective resistance metric as the intrinsic notion of distance, and we approach the analysis of the resistance network by studying the space of functions on the vertices which have finite (Dirichlet) energy. There is a large existing literature on this subject, but ours is unique in several respects, most of which are due to the following.
\begin{itemize}
  \item We use the effective resistance metric to find canonical Hilbert spaces of functions associated with the resistance network.
  \item We adhere to the intuition arising from the metaphor of electrical resistance networks, including Kirchhoff's Law and Ohm's Law.
  \item We apply the results of our Hilbert space construction to the isotropic Heisenberg ferromagnet and prove a theorem regarding long-range order in quantum statistical mechanics for certain lattice networks.
  \item It is known (see \cite{Lyons:ProbOnTrees} and the references therein) that the resistance metric is unique for finite graphs and not unique for certain infinite networks. 
  We are able to clarify and explain the difference in terms of certain Hilbert space structures, and also in terms of Dirichlet vs. Neumann boundary conditions for a certain operator. Additionally, we introduce trace resistance, and harmonic resistance and relate these to the aforementioned.
\end{itemize}

A large portion of this volume is dedicated to developing an operator-theoretic understanding of a certain \emph{boundary} which appears in diverse guises. The boundary appears first in Chapter~\ref{sec:energy-Hilbert-space} in a crucial but mysterious way, as the agent responsible for the misbehavior of a certain formula relating the Laplace operator to the energy form. It reappears in Chapter~\ref{sec:effective-resistance-metric} as the agent responsible for the failure of various formulations of the effective resistance $R(x,y)$ to agree for certain infinite networks. In Chapter~, we pursue the boundary directly, using tools from operator theory and stochastic integration. The pedagogical aim behind this approach is to demonstrate operator theory via a series of applications. Many examples are given throughout the book. These may serve as independent student projects, although they are not exercises in the traditional sense.

\section*{Prerequisites} 

We have endeavored making this book as accessible and self-contained as possible. Nonetheless, readers coming across various ideas for the first time may wish to consult the following books: \cite{DoSn84} (resistance networks), \cite{AlFi09,LevPerWil08} (probability), 
and \cite{DuSc88} (unbounded operators).

\section*{Detailed description of contents} 
\addcontentsline{toc}{subsection}{Detailed description of contents}

\headerquote{Mathematical science is in my opinion an indivisible whole, an organism whose vitality is conditioned upon the connection of its parts.}{---~D.~Hilbert}

\S\ref{sec:electrical-resistance-networks}\emph{ --- Electrical resistance networks.}
We introduce the resistance network as a connected simple graph $\Graph = \{\verts,\edges\}$ equipped with a positive weight function \cond on the edges. The edges $\edges \ci \verts \times \verts$ are ordered pairs of vertices, so \cond is required to be symmetric. Hence, each edge $(x,y) \in \edges$ is interpreted as a conductor with conductance $\cond_{xy}$ (or a resistor with resistance $\cond_{xy}^{-1}$. Heuristically, smaller conductances (or larger resistances) correspond to larger distances; see the discussion of \S\ref{sec:effective-resistance-metric} just below. We make frequent use of the weight that \cond defines on the vertices via $\cond(x) = \sum_{y \nbr x} \cond_{xy}$, where $y \nbr x$ indicates that $(x,y) \in \edges$. The graphs we are most interested in are infinite graphs, but we do not make any general assumptions of regularity, group structure, etc. We require that $\cond(x)$ is finite at each $x \in \verts$, but we do not generally require that the degree of a vertex be finite, nor that $\cond(x)$ be bounded.

In the ``cohomological'' tradition of von Neumann, Birkhhoff, Koopman, and others \cite{vN32c,Koo36,Koo57}, we study the resistance network by analyzing spaces of functions defined on it. These are constructed rigourously as Hilbert spaces in \S\ref{sec:vonNeumann's-embedding-thm}; in the meantime we collect some results about functions $u,v:\verts \to \bR$ defined on the vertices. The \emph{network Laplacian} (or \emph{discrete Laplace operator}) operates on such a function by taking $v(x)$ to a weighted average of its values at neighbouring points in the graph, i.e.,
\linenopax
\begin{equation}\label{eqn:intro:laplacian}
  (\Lap v)(x) := \sum_{y \nbr x} \cond_{xy}(v(x)-v(y)) = \sum_{y \nbr x} \frac{v(x)-v(y)}{ \cond_{xy}^{-1}},
\end{equation}
where $x \nbr y$ indicates that $(x,y) \in \edges$. (The rightmost expression in formula \eqref{eqn:intro:laplacian} is written so as to resemble the familiar difference quotients from calculus.) This is the usual second-difference operator of numerical analysis, when adapted to a network. There is a large literature on discrete harmonic analysis (basically, the study of the graph/network Laplacian) which include various probabilistic, combinatoric, and spectral approaches. It would be difficult to give a reasonably complete account, but the reader may find an enjoyable approach to the probabilistic perspective in \cite{Spitzer,Telcs06a}, the combinatoric in \cite{ABZ07},  the analytic in \cite{Fab06}, and the spectral in \cite{Chu01,GIL06a}. More sources are peppered about the relevant sections below. Our formulation \eqref{eqn:intro:laplacian} differs from the stochastic formulation often found in the literature, but the two may easily be reconciled; see \eqref{eqn:probablistically-reweighted-laplacian}.

Together with its associated quadratic form, the bilinear \emph{(Dirichlet) energy form}
\linenopax
\begin{align}\label{eqn:intro:energy-form}
  \energy(u,v)
  :=& \frac12 \sum_{x \in \verts} \sum_{y \nbr x}  \cond_{xy} (u(x)-u(y))(v(x)-v(y))
\end{align}
acts on functions $u,v:\verts \to \bR$ and plays a central role in the (harmonic) analysis on $(\Graph, \cond)$. (There is also the dissipation functional \diss, a twin of \energy which acts on functions defined on the edges \edges and is introduced in the following section.) The first space of functions we study on the resistance network is the \emph{domain of the energy}, that is,
\linenopax
\begin{equation}\label{eqn:intro:domE}
  \dom \energy := \{u:\verts \to \bR \suth \energy(u)<\iy\}.
\end{equation}
In \S\ref{sec:vonNeumann's-embedding-thm}, we construct a Hilbert space from the resistance metric (and show it to be a canonical invariant for $(\Graph,\cond)$ in \S\ref{sec:HE-as-an-invariant}), thereby recovering the familiar result that $\dom \energy$ is a Hilbert space with inner product \energy. (Actually, this is not quite true, as \energy is only a quasinorm; see the discussion of \S\ref{sec:vonNeumann's-embedding-thm} just below for a more accurate description.)

For finite graphs, we prove the simple and folkloric key identity which relates the energy and the Laplacian:
\linenopax
\begin{equation}\label{eqn:intro:E(u,v)=<u,Lapv>}
  \energy(u,v) = \la u, \Lap v\ra_\unwtd = \la \Lap u, v\ra_\unwtd,
  \qq u,v \in \dom \energy,
\end{equation}
where $\la u, \Lap v\ra_\unwtd = \sum_{x \in \verts} u(x) \Lap v(x)$ indicates the standard $\ell^2$ inner product. 
The formula \eqref{eqn:intro:E(u,v)=<u,Lapv>} is extended to infinite networks in Theorem~\ref{thm:E(u,v)=<u,Lapv>+sum(normals)} (see \eqref{eqn:intro:discrete-Gauss-Green} for a preliminary discussion), where a third term appears. Indeed, understanding the mysterious third term is the motivation for most of this investigation. 

\pgap

\emph{\S\ref{sec:currents-and-potentials} --- Currents and potentials on resistance networks.}
We collect several well-known and folkloric results, and reprove some variants of these results in the present context. \emph{Currents} are introduced as skew-symmetric functions on the edges; the intuition is that $\curr(x,y) = -\curr(y,x) > 0$ indicates electrical current flowing from $x$ to $y$. In marked contrast to common tradition in geometric analysis \cite{ABZ07, PaSz07}, we do not fix an orientation. For us, an orientation is a choice of one of $\{(x,y), (y,x)\}$ for each edge, and hence just a notation to be redefined as convenient. In particular, any nonvanishing current defines an orientation; one makes the choice so that \curr is a positive function. At this point we give the definition of the \emph{dissipation}, an inner product defined for functions on the edges, and its associated quadratic form:
\linenopax
\begin{equation}\label{eqn:intro:dissipation}
  \diss(\curr) = \frac12 \sum_{(x,y) \in \edges}  \cond_{xy}^{-1} \curr(x,y)^2.
\end{equation}

Most of our results in this section are groundwork for the sections to follow; several results are folkloric or obtained elsewhere in the literature. We include items which relate directly to results in later sections; the reader seeking a more well-rounded background is directed to \cite{LevPerWil08, Lyons:ProbOnTrees, Soardi94, CdV98, Bollobas} and the excellent elementary introduction \cite{DoSn84}. After establishing the Hilbert space framework of \S\ref{sec:energy-Hilbert-space}, we exploit the close relationship between the two functionals \energy and \diss, and use operators to translate a problem from the domain of one functional to the domain of the other. We also introduce Kirchhoff's Law and Ohm's Law, and in \S\ref{sec:compatibility-problem} we discuss the related \emph{compatibility problem}: every function on the vertices induces a function on the edges via Ohm's Law, but not every function on the edges comes from a function on the vertices. This is related to the fact that most currents are not ``efficient'' in a sense which can be made clear variationally (cf.~Theorem~\ref{thm:minimal-flows-are-induced-by-minimizers}) and which is important in the definition of effective resistance metric in Theorem~\ref{thm:effective-resistance-metric}. We recover the well-known fact that the dissipation of an induced current is equal to the energy of the function inducing it in Lemma~\ref{thm:energy=dissipation}; this is formalized as an isometric operator in Theorem~\ref{thm:drp-and-drp-are-partial-isometries}. We show that the equation
\linenopax
\begin{equation}\label{eqn:intro:Pot}
  \Lap v = \gd_\ga-\gd_\gw
\end{equation}
always has a solution; we call such a function a \emph{dipole}. 
In \eqref{eqn:intro:Pot} and everywhere else, we use the notation $\gd_x$ to indicate a Dirac mass at $x \in \verts$, that is,
\linenopax
\begin{equation}\label{eqn:intro:Dirac}
  \gd_x = \gd_x(y) :=
  \begin{cases}
    1, &y=x, \\
    0, &\text{else.}
  \end{cases}
\end{equation}
Proving the  existence of dipoles allows us to fill gaps in \cite{Pow76a,Pow76b} (see \S\ref{sec:motivation-and-application} just below) and extend the definition of effective resistance metric in Theorem~\ref{thm:effective-resistance-metric} to infinite dimensions.

As is discussed at length in Remark~\ref{rem:dipoles,monopoles,harmonics}, the study of dipoles, monopoles, and harmonic functions is a recurring theme of this book:
\linenopax
\begin{align*}
  \Lap v = \gd_\ga - \gd_\gw,
  \qq \Lap w = -\gd_\gw,
  \qq \Lap h = 0.
\end{align*}
As mentioned above, for any network \Graph and any vertices $x,y \in \verts$, there is a dipole in $\dom \energy$. However, $\dom \energy$ does not always contain monopoles or nonconstant harmonic functions; the existence of monopoles is equivalent to transience of the network \cite{TerryLyons}; we give a new criterion for transience in Lemma~\ref{thm:transient-iff-Fin=Ran(Lap)}. In Theorem~\ref{thm:monopoles-on-Zd}, we show that the integer lattice networks $(\bZd,\one)$ support monopoles iff $d \geq 3$, but in Theorem~\ref{thm:harmonics-on-Zd-are-linear} we show all harmonic functions on $(\bZd,\one)$ are linear and hence do not have finite energy. (Both of these results are well known; the first is a famous theorem of Polya --- we include them for the novelty of method of proof.) In contrast, the binary tree in Example~\ref{exm:monopole-on-binary-tree} support monopoles and nontrivial harmonic functions, both of finite energy (any network supporting nontrivial harmonic functions also supports monopoles, cf.~\cite[Thm.~1.33]{Soardi94}). It is apparent that monopoles and nontrivial harmonic functions are sensitive to the asymptotic geometry of $(\Graph, \cond)$.

\pgap

\emph{\S\ref{sec:energy-Hilbert-space} --- The energy Hilbert space \HE.}
We use the natural Hilbert space structure on the space of finite-energy functions (with inner product given by \energy) to reinterpret previous results as claims about certain operators, and thereby clarify and generalize results from \S\ref{sec:electrical-resistance-networks}--\S\ref{sec:currents-and-potentials}. This is the energy space \HE. 

We construct a reproducing kernel for \HE from first principles (i.e., via Riesz's Lemma) in \S\ref{sec:L_x-and-v_x}. If $o \in \verts$ is any fixed reference point, define $v_x$ to be the vector in \HE which corresponds (via Hilbert space duality) to the evaluation functional $L_{x}$:
\linenopax
\begin{align*}
  L_{x} u := u(x) - u(o).
\end{align*}
Then the functions $\{v_x\}$ form a reproducing kernel, and $v_x$ is a solution of the discrete Dirichlet problem $\Lap v_x = \gd_x - \gd_o$. Although these functions are linearly independent, they are usually neither an orthonormal basis (onb) nor a frame. However, the span of $\{v_x\}$ is dense in $\dom \energy$ and appears naturally when the energy Hilbert space is constructed from the resistance metric by von Neumann's method; cf.~\S\ref{sec:vonNeumann's-embedding-thm}. Note that the Dirac masses $\{\gd_x\}_{\verts}$, which are the usual candidates for an onb, are \textbf{not orthogonal} with respect to the energy inner product \eqref{eqn:intro:energy-form}; cf.~\eqref{eqn:energy-of-Diracs}. In fact, Theorem~\ref{thm:TFAE:Fin,Harm,Bdy} shows that $\{\gd_x\}_{\verts}$ \textbf{may not even be dense} in the energy Hilbert space! Thus, $\{v_x\}$ is the only canonical choice for a representing set for functions of finite energy.

In \S\ref{sec:The-role-of-Fin-in-HE} we use the Hilbert space structure of \HE to better understand the role of the nontrivial harmonic functions. In particular, Lemma~\ref{thm:HE=Fin+Harm} shows that we may decompose \HE into the functions of finite support (\Fin) and the harmonic functions of finite energy (\Harm):
\linenopax
\begin{equation}\label{eqn:intro:HE=Fin+Harm}
  \HE = \Fin \oplus \Harm.
\end{equation}

In \S\ref{sec:relating-energy-form-to-Laplacian}, we prove a discrete version of the Gauss-Green formula (Theorem~\ref{thm:E(u,v)=<u,Lapv>+sum(normals)}) which appears to be absent from the literature:
\linenopax
\begin{equation}\label{eqn:intro:discrete-Gauss-Green}
  \energy(u,v)
    = \sum_{x \in \verts} u(x) \Lap v(x) 
      + \sum_{x \in \bd \Graph} u(x) \dn v(x),
    \q\forall u \in \HE, v \in \MP
\end{equation}
where $\dn v(x)$ denotes the \emph{normal derivative of $v$}, and \MP is a space containing $\spn\{v_x\}$; see \S\ref{sec:relating-energy-form-to-Laplacian} for precise definitions. For the moment, both the boundary and the normal derivatives are understood as limits (and hence vanish trivially for finite graphs); we will be able to define these objects more concretely via techniques of Gel'fand in \S\ref{sec:the-boundary}.

It turns out that the boundary term (that is, the rightmost sum in \eqref{eqn:intro:discrete-Gauss-Green}) vanishes unless the network supports \emph{nontrivial harmonic functions} (that is, nonconstant harmonic functions of finite energy). More precisely, in Theorem~\ref{thm:TFAE:Fin,Harm,Bdy} we prove that there exist $u,v \in \HE$ for which $\sum_{\bd\Graph} u \dn v \neq 0$ if and only if the network is transient. That is, the random walk on the network with transition probabilities $p(x,y) = \cond_{xy}/\cond(x)$ is transient. We also give several other equivalent conditions for transience, in \S\ref{sec:More-about-monopoles}.

It is easy to prove (see Corollary~\ref{thm:nontrivial-harmonic-fn-is-not-in-L2}) that nontrivial harmonic functions cannot lie in $\ell^2(\verts)$. This is why we \textbf{do not} require $u,v \in \ell^2(\verts)$ in general, and why we stringently avoid including such a requirement in the definition of the domain of the Laplacian. Such a restriction would remove the nontrivial harmonic functions from the scope of our analysis, and we will see that they are at the core of some of the most interesting phenomena appearing on an infinite resistance network. 


\pgap

\emph{\S\ref{sec:effective-resistance-metric} --- Effective resistance metric.} 
The effective resistance metric $R$ is foundational to our study, instead of the shortest-path metric more commonly used as graph distance. The shortest-path metric on a weighted graph is usually defined to be the sum of the resistances in any shortest path between two points. The effective resistance metric is also defined via \cond, but in a more complicated way. The crucial difference is that the effective resistance metric reflects both the topology of the graph \emph{and} the weighting \cond; two points are closer together when there is more connectivity (more paths and/or paths with greater conductance) between them. The effective resistance metric is a much more accurate way to measure distance when travel from point $x$ to point $y$ can be accomplished simultaneously through many paths, for example, flow of electrical current, fluid diffusion through porous media, or data transfer over the internet.

In \S\ref{sec:Resistance-metric-on-finite-networks}, we give a multifarious definition of the effective resistance metric $R$, which may be physically characterized as the voltage drop between two vertices when electrical leads with a fixed current are applied to the two vertices. Most of these formulations appear elsewhere in the literature, but some appear to be specific to the physics literature, some to probability, and some to analysis. We collect them and prove their equivalence in Theorem~\ref{thm:effective-resistance-metric}, including a couple new formulations that will be useful in later sections.

It is somewhat surprising that when these formulas are extended to an infinite network in the most natural way, they are no longer equivalent. (Note that each of the six formulas has both a free and wired version, but some appear much less natural in one version than in the other.) Some of the formulas lead to the ``free resistance'' $R^F$ and others lead to the ``wired resistance'' $R^W$; here we follow the terminology of \cite{Lyons:ProbOnTrees}. In \S\ref{sec:Resistance-metric-on-infinite-networks}, we precisely characterize the types of extensions that lead to each, and explain this phenomenon in terms of projections in Hilbert space, Dirichlet vs. Neumann boundary conditions, and via probabilistic interpretation. Additionally, we discuss the ``trace resistance'' given in terms of the trace of the Dirichlet form \energy, and we study the ``harmonic resistance'' which is the difference between $R^F$ and $R^W$ and is not typically a metric.

\pgap

\emph{\S\ref{sec:Construction-of-HE} --- Construction of the energy space \HE.}
In \S\ref{sec:vonNeumann's-embedding-thm}, we use a theorem of von Neumann to give an isometric imbedding of the metric space $(\Graph,R)$ into \HE; cf.~Theorem~\ref{thm:R^F-embed-ERN-in-Hilbert}. 
For infinite networks, $(\Graph,R^F)$ embeds into \HE and $(\Graph,R^W)$ embeds into \Fin. 
In \S\ref{sec:HE-as-an-invariant} we discuss how this enables one to interpret \HE as an invariant of the original resistance network.

\pgap

\emph{\S\ref{sec:the-boundary} --- The boundary $\bd \Graph$ and boundary representation.}
We study the boundary $\bd \Graph$ in terms of the Laplacian by reinterpreting the boundary term of \eqref{eqn:intro:discrete-Gauss-Green} as an integral over a space which contains \HE. This gives a representation of $\bd\Graph$ as a measure space whose structure is well-studied.

In Theorem~\ref{thm:boundary-repn-for-harmonic} of \S\ref{sec:boundary-outline}, we observe that an important consequence of \eqref{eqn:intro:discrete-Gauss-Green} is the following boundary representation for the harmonic functions:
\linenopax
  \begin{equation}\label{eqn:intro:boundary-repn-for-harmonic}
    u(x) = \sum_{\bd \Graph} u \dn{{h_x}} + u(o), 
  \end{equation}
  for $u \in \Harm$, where $h_x = \Phar v_x$ is the projection of $v_x$ to \Harm; see \eqref{eqn:intro:HE=Fin+Harm}. This formula is in the spirit of Choquet theory and the Poisson integral formula and is closely related to Martin boundary theory.

Unfortunately, the sum in \eqref{eqn:intro:boundary-repn-for-harmonic} is only understood in a limiting sense and so provides limited insight into the nature of $\bd\Graph$. This motivates the development of a more concrete expression. We use a self-adjoint extension \LapS of \Lap to construct a Gel'fand triple $\Schw \ci \HE \ci \Schw'$ and a Gaussian probability measure \prob. 
Here, $\Schw := \dom(\LapS^\iy)$ is a suitable dense (Schwartz) space of ``test functions'' on the resistance network, and $\Schw'$ is the corresponding dual space of ``distributions'' (or ``generalized functions''). This enables us to identify $\bd \Graph$ as a subset of $\Schw'$, and in Corollary~\ref{thm:Boundary-integral-repn-for-harm}, we rewrite \eqref{eqn:intro:boundary-repn-for-harmonic} more concretely as
\linenopax
  \begin{equation}\label{eqn:intro-boundary-of-G-from-Minlos}
    u(x) = \int_{\Schw'} u(\gx) {h_x}(\gx) \, d\prob(\gx) + u(o),
  \end{equation}
  again for $u \in \Harm$ and with $h_x = \Phar v_x$. Thus we study the metric/measure structure of \Graph by examining an associated Hilbert space of random variables. This is motivated in part by Kolmogorov's pioneering work on stochastic processes (see \S\ref{sec:Kolmogorov-construction-of-L2(gW,prob)}) as well as on a powerful refinement of Minlos. The latter is in the context of the Gel'fand triples mentioned just above; see \cite{Nelson64} and \S\ref{sec:Gel'fand-triples-and-duality} below. 
Further applications to harmonic analysis and to physics are given in \S\ref{sec:Probabilistic-interpretation}--\S\ref{sec:magnetism}.

\pgap

\emph{\S\ref{sec:Lap-on-HE} --- The Laplacian on \HE.}
We study the operator theory of the Laplacian in some detail in \S\ref{sec:Properties-of-Lap-on-HE}, examining the various domains and self-adjoint extensions. We identify one domain for the Laplacian which allows for the choice of a particular self-adjoint extension for the constructions in \S\ref{sec:the-boundary}. Also, we give technical conditions which must be considered when the graph contains vertices of infinite degree and/or the conductance functions $\cond(x)$ is unbounded on \verts. This results in an extension of the Royden decomposition to $\HE = \Fin_1 \oplus \Fin_2 \oplus \Harm$, where $\Fin_2$ is the \energy-closure of $\spn\{\gd_x-\gd_o\}$ and $\Fin_1$ is the orthogonal complement of $\Fin_2$ within \Fin. Example~\ref{exm:Fin_2-not-dense-in-Fin} shows a case where $\Fin_2$ is not dense in \Fin.

In \S\ref{sec:defect-space}, we study the defect space of \LapV, that is, the space spanned by solutions to $\Lap u = -u$. In \S\ref{sec:boundary-form}, we relate the boundary term of \eqref{eqn:intro:discrete-Gauss-Green} to the the boundary form
\linenopax
\begin{equation}\label{eqn:intro:boundary-form}
  \bdform(u,v) :=
  \tfrac1{2\ii}\left(\la \LapV^\ad u, v\ra_\energy - \la u, \LapV^\ad v\ra_\energy\right),
  \qq u,v \in \dom(\LapV^\ad)
\end{equation}
of classical functional analysis; cf. \cite[\S{XII.4.4}]{DuSc88}. This gives a way to detect whether or not a given network has a boundary by examining the deficiency indices of \Lap. In Theorem~\ref{thm:LapV-not-ess-selfadjoint-iff-Harm=0}, we show that if \Lap fails to be essentially self-adjoint, then $\Harm \neq \{0\}$. In general, the converse does not hold: Corollary~\ref{thm:c-bdd-implies-Lap-selfadjoint-on-HE} shows that \Lap has no defect when $\deg(x)<\iy$ and $\cond(x)$ is bounded. (Thus, any homogeneous tree of degree 3 or higher with constant conductances provides a counterexample to the converse.)

In \S\ref{sec:Frames-and-dual-frames}, we study the relation between the reproducing kernel $\{v_x\}$ and the spectral properties of \Lap and its self-adjoint extensions. In particular, we examine the necessary conditions for $\{v_x\}$ to be a frame for \HE, and the relation between $v_x$ and $\gd_x$.

\pgap

\emph{\S\ref{sec:L2-theory-of-Lap-and-Trans} --- The $\ell^2$ theory of \Lap and \Trans.}
We consider some results for \Lap and \Trans as operators on $\ell^2(\unwtd)$, where the inner product is given by $\la u, \Lap v\ra_\unwtd := \sum u(x) \Lap v(x)$ and on $\ell^2(\cond)$, where the inner product is given by $\la u, \Lap v\ra_\cond := \sum \cond(x) u(x) \Lap v(x)$.

We prove that the Laplacian is essentially self-adjoint on $\ell^2(\unwtd)$ under very mild hypotheses in \S\ref{sec:self-adjointness-of-the-Laplacian}. The subsequent spectral representation allows us to give a precise characterization of the domain of the energy functional \energy in this context. In \S\ref{sec:the-transfer-operator}, we examine boundedness and compactness of \Lap and \Trans in terms of the decay properties of \cond. The space $\ell^2(\cond)$ considered in \S\ref{sec:weighted-spaces} is essentially a technical tool; it allows for a proof that the terms of the Discrete Gauss-Green formula are absolutely convergent and hence independent of any exhaustion. However, it is also interesting in its own right, and we show an interesting connection with the probabilistic Laplacian $\cond^{-1} \Lap$. Results from this section imply that \Lap is also essentially self-adjoint on \HE, subject to the same mild hypotheses as the $\ell^2(\unwtd)$ case.

The energy Hilbert space \HE contains much different information about a given infinite graph system $(\Graph, \cond)$ than does the more familiar $\ell^2$ sequence space, even when appropriate weights are assigned. In the language of Markov processes, \HE is better adapted to the study of $(\Graph,\cond)$ than $\ell^2$. One reason for this is that \HE is intimately connected with the resistance metric $R$.

\pgap

\emph{\S\ref{sec:H_energy-and-H_diss} --- \HE and \HD.}
The dissipation space \HD is the Hilbert space of functions on the edges when equipped with the dissipation inner product. We solve problems in discrete potential theory with the use of the \emph{drop operator} \drp (and its adjoint \drpa), where 
\linenopax
\begin{equation}\label{eqn:intro:drop}
  \drp v(x,y) := \cond_{xy}(v(x)-v(y)).
\end{equation}
The drop operator \drp is, of course, just an implementation of Ohm's Law, and can be interpreted as a weighted boundary operator in the sense of homology theory. The drop operator appears elsewhere in the literature, sometimes without the weighting $\cond_{xy}$; see \cite{Chu01,Telcs06a,Woess00}. However, we use the adjoint of this operator with respect to the energy inner product, instead of the $\ell^2$ inner product used by others. This approach appears to be new, and it turns out to be more compatible with physical interpretation. For example, the displayed equation preceding \cite[(2.2)]{Woess00} shows that the $\ell^2$ adjoint of the drop operator is incompatible with Kirchhoff's node law. Since the resistance metric may be defined in terms of currents obeying Kirchhoff's laws, we elect to make this break with the existing literature. Additionally, this strategy will allow us to solve the compatibility problem described in \S\ref{sec:compatibility-problem} in terms of a useful minimizing projection operator \Pdrp, discussed in detail in \S\ref{sec:Solving-potential-theoretic-problems}. Furthermore, we believe our formulation is more closely related to the (co)homology of the resistance network as a result.

We decompose \HD into the direct sum of the range of \drp and the currents which are sums of characteristic functions of cycles
\linenopax
\begin{equation}\label{eqn:intro-HD-decomp}
  \HD = \ran \drp \oplus \cl \spn{\charfn{\cycle}},
\end{equation}
where \cycle is a cycle, i.e., a path in the graph which ends where it begins. In \eqref{eqn:intro-HD-decomp} and elsewhere, we indicate the closed linear span of a set by $\cl\spn{\charfn{\cycle}} := \cl \spn\{\charfn{\cycle}\}$. From \eqref{eqn:intro:HE=Fin+Harm} (and the fact that \drp is an isometry), it is clear that the first summand of \eqref{eqn:intro-HD-decomp} can be further decomposed into weighted edge neighbourhoods $\drp\gd_x$ and the image of harmonic functions under \drp in Theorem~\ref{thm:HD=Nbd+Kir+Cyc}. After a first draught of this book was complete, we discovered that the same approach is taken in \cite{Lyons:ProbOnTrees}. One of us (PJ) recalls conversations with Raul Bott concerning an analogous Hilbert space operator theoretic approach to electrical networks; apparently attempted in the 1950s in the engineering literature. We could not find details in any journals; the closest we could come is the fascinating paper \cite{BoDu49} by Bott et al. A further early source of influence is Norbert Wiener's paper \cite{WiRo46}. 

In \S\ref{sec:Solving-potential-theoretic-problems}, we describe how \drpa solves the compatibility problem and may be used to solve a large class of problems in discrete potential theory. Also, we discuss the analogy with complex analysis.

\pgap

\emph{\S\ref{sec:Probabilistic-interpretation} --- Probabilistic interpretations.} In \cite{Lyons:ProbOnTrees,DoSn84,Telcs06a,Woess00} and elsewhere, the random walk on an resistance network is defined by the transition probabilities $p(x,y) := \cond_{xy}/\cond(x)$.  In this context, the probabilistic transition operator is $\Prob = \cond^{-1}\Trans$ and one uses the stochastically renormalized Laplacian $\Lap_\cond := \cond^{-1} \Lap$, where \cond is understood as a multiplication operator; see Definition~\ref{def:conductance-measure}. This approach also arises in the discussion of trace resistance in \S\ref{sec:trace-resistance} and allows one to construct currents on the graph as the average motion of a random walk.

As an alternative to the approach described above, we discuss a probabilistic interpretation slightly different from those typically found in the literature: we begin with a voltage potential as an initial condition, and consider the induced current \curr. The components of such a flow are called \emph{current paths} and provide a way to interpret potential-theoretic problems in a probabilistic setting. We study the random walks where the transition probability is given by $\curr(x,y)/\sum_{z \nbr x} \curr(x,z)$. We consider the harmonic functions in this context, which we call \emph{forward-harmonic functions}, and the associated \emph{forward-Laplacian} of Definition~\ref{def:forward-Laplacian}. We give a complete characterization of forward-harmonic functions as \emph{cocycles}, following \cite{Jor06}.

\pgap

\version{}{
  \marginpar{This needs to be finished once we are.}
  \emph{\S\ref{sec:cuntz} --- Cuntz and Cuntz-Krieger algebras.}

  \pgap

  \emph{\S\ref{sec:subgraphs-and-limits} --- Subgraphs and limits.}

  \pgap
}

\emph{\S\ref{sec:examples} --- Examples.} We collect an array of examples that illustrate the various phenomena encountered in the theory and work out many concrete examples. Some elementary finite examples are given in \S\ref{sec:finite-graphs} to give the reader an idea of the basics of resistance network theory. In \S\ref{sec:infinite-graphs} we move on to infinite graphs. 

\pgap

\emph{\S\ref{sec:tree-networks} --- Trees.}
When the resistance network is a tree (i.e., there is a unique path between any two vertices), the resistance distance coincides with the geodesic metric, as there is always exactly one path between any two vertices; cf.~Lemma~\ref{thm:shortest-path-bounds-resistance-distance} and the preceding discussion. When the tree has exponential growth, as in the case of homogeneous trees of degree $\geq 3$, one can always construct nontrivial harmonic functions, and monopoles of finite energy. In fact, there is a very rich family of each, and this property makes this class of examples a fertile testing ground for many of our theorems and definitions. In particular, these examples highlight the relevance and distinctions between the boundary (as we construct it), the Cauchy completion, and the graph ends of \cite{PicWoess90,Woess00}. In particular, they enable one to see how adjusting decay conditions on \cond affects these things.

\pgap

\emph{\S\ref{sec:lattice-networks} --- Integer lattices.}
The lattice resistance network $(\bZd,\cond)$ have vertices at the points of \bRd which have integer coordinates, and edges between every pair of vertices $(x,y)$ with $|x-y|=1$. The case for $\cond = \one$ is amenable to Fourier analysis, and in \S\ref{sec:simple-lattice-networks} we obtain explicit formulas for many expressions:
\begin{itemize}
  \item Lemma~\ref{thm:vx-on-Zd} gives a formula for the potential configuration functions $\{v_x\}$.
  \item Theorem~\ref{thm:R(x,y)-on-Zd} gives a formula for the resistance distance $R(x,y)$.
  \item Theorem~\ref{thm:finite-resistance-to-infinity-in-Zd} gives a formula for the resistance distance to infinity in the sense $R(x,\iy) = \lim_{y \to \iy} R(x,y)$.
  \item Theorem~\ref{thm:monopoles-on-Zd} gives a formula for the solution $w$ of $\Lap w = -\gd_o$ on \bZd; it is readily seen that this $w$ has finite energy (i.e., is a monopole) iff $d \geq 3$.
\end{itemize}
In \cite{Polya21}, P\'{o}lya proved that the random walk on this graph is transient if and only if $d \geq 3$; see \cite{DoSn84} for a nice exposition. We offer a new characterization of this dichotomy (there exist monopoles on \bZd if and only if $d \geq 3$) which we recover in this section via a new (and completely constructive) proof. In Remark~\ref{rem:HE(Zd)-lies-in-L2(Zd)} we describe how in the infinite integer lattices, functions in \HE may be approximated by functions of finite support.

\pgap

\emph{\S\ref{sec:Magnetism-and-long-range-order} --- Magnetism.}
The integer lattice networks $(\bZd,\one)$ investigated in \S\ref{sec:lattice-networks} comprise the framework of infinite models in thermodynamics and in quantum statistical mechanics. In \S\ref{sec:magnetism} we employ these formulas in the refinement of an application to the theory of the (isotropic Heisenberg) model of ferromagnetism as studied by R. T. Powers. In addition to providing an encapsulated version of the Heisenberg model, we give a commutative analogue of the model, extend certain results of Powers from \cite{Pow75,Pow76a,Pow76b,Pow78,Pow79}, and discuss the application of the resistance metric to the theory of ferromagnetism and ``long-range order''. This problem was raised initially by R. T. Powers, and may be viewed as a noncommutative version of Hilbert spaces of random variables.

Ferromagnetism in quantum statistical mechanics involves algebras of noncommutative observables and may be described with the use of states on $C^\ad$-algebras. As outlined in the cited references, the motivation for these models draw on thermodynamics; hence the notions of equilibrium states (formalized as KMS states, see \S\ref{sec:KMS-states-appendix}). These KMS states are states in the $C^\ad$-algebraic sense (that is, positive linear functionals with norm 1), and they are indexed by absolute temperature. Physicists interpret such objects as representing equilibria of infinite systems.

In the present case, we consider spin observables arranged in a lattice of a certain rank, $d=1,2,3,\dots$, and with nearest-neighbor interaction. Rigourous mathematical formulation of phase transitions appears to be a hopeless task with current mathematical technology. As an alternative avenue of enquiry, much work has been conducted on the issue of \emph{long-range order}, i.e., the correlations between observables at distant lattice points. These correlations are measured relative to states on the $C^\ad$-algebra; in this case in the KMS states for a fixed value of temperature.

While we shall refer to the literature, e.g. \cite{BrRo79,Rue69} for formal definitions of key terms from the $C^\ad$-algebraic formalism of quantum spin models, physics, and KMS states, we include a minimal amount of background and terminology from the physics literature.

\pgap

\emph{\S\ref{sec:future-directions} --- Future Directions.}
We conclude with a brief discussion of several projects which have arisen from work on the present book, as well as some promising new directions that we have not yet had time to pursue.

\pgap

\emph{Appendices.} We give some background material from functional analysis in Appendix~\ref{sec:functional-analysis}, and operator theory in Appendix~\ref{sec:operator-theory}. In Appendix~\ref{sec:navigation-aids}, we include some diagrams to help clarify the properties of the many operators and spaces we discuss, and the relations between them.

\pgap

\section*{What this book is about}
\addcontentsline{toc}{subsection}{What this book is about}

\headerquote{`Obvious' is the most dangerous word in mathematics.}{---~E.~T.~Bell}

The effective resistance metric provides the foundation for our investigations because it is the natural and intrinsic metric for a resistance network, as the work of Kigami has shown; see \cite{Kig01} and the extensive list of references by the same author therein. Moreover, the close relationship between diffusion geometry (i.e., geometry of the resistance metric) \cite{Maggioni08 ,CoifmanMaggioniSzlam, CoifmanMaggioni08, CoifmanMaggioni06} and random walks on graphs leads us to expect/hope there will be many applications of our results to several other subjects, in addition to fractals: models in quantum statistical mechanics, analysis of energy forms, interplay between self-similar measures and associated energy forms, certain discrete models arising in the study of quasicrystals (e.g., \cite{BaMo00,BaMo01}), and multiwavelets (e.g. \cite{BJMP05,DJ06, DJ07, Jor06}), among others. A general theme of these areas is that the underlying space is not sufficiently regular to support a group structure, yet is ``locally'' regular enough to allow analysis via probabilistic techniques. Consequently, the analysis of functions on such spaces is closely tied to Dirichlet energy forms and the graph Laplacian operator associated to the graph. This appears prominently in the context of the present book as follows:
\begin{enumerate}
  \item The embedding of the metric space $((\Graph, \cond), R)$ into the Hilbert space \HE of functions of finite energy, in such a way that the original metric may be recovered from the norm, i.e.,
      \linenopax
      \begin{align*}
        R(x,y) = \energy(v_x-v_y) = \|v_x-v_y\|_\energy^2,
      \end{align*}
      where $v_x \in \HE$ is the image of $x$ under the embedding.
  \item The relation of the energy form to the graph Laplacian via the equation
      \linenopax
      \begin{equation}\label{eqn:intro:aaaagain}
      \energy(u,v)
        = \sum_{x \in \verts} u(x) \Lap v (x)
          + \sum_{x \in \bd \Graph} u(x) \dn v(x),
      \end{equation}
      introduced just above in the discussion of \S\ref{sec:electrical-resistance-networks}. Each summation on the right hand side of \eqref{eqn:intro:aaaagain} is more subtle than it appears. These details for the first sum are given in Theorem~\ref{thm:E(u,v)=<u,Lapv>+sum(normals)}, and the details for the second sum are the focus of almost all of \S\ref{sec:the-boundary}. 
  \item The presence of nonconstant harmonic functions of finite energy. These are precisely the objects which support the boundary term in \eqref{eqn:intro:aaaagain} and imply $R^W(x,y) < R^F(x,y)$. They are also responsible for the boundary described in \S\ref{sec:the-boundary}.
  \item The solvability of the Dirichlet problem $\Lap w = -\gd_y$, where $\gd_y$ is a Dirac mass at the vertex $y \in \verts$. The existence of finite-energy solutions $w$ is equivalent to the transience of the random walk on the network. Such functions are called \emph{monopoles} and (via Ohm's law) they induce a \emph{unit flow to infinity} as discussed in \cite{DoSn84,Lyons:ProbOnTrees,LevPerWil08}.
\end{enumerate}

\begin{remark}[Relation to numerical analysis]\label{rem:numerical-analysis}
  In addition to uses in graph theory and electrical networks, the discrete Laplacian \Lap has other uses in numerical analysis: many problems in PDE theory lend themselves to discretizations in terms of subdivisions or grids of refinements in continuous domains. A key tool in applying numerical analysis to solving partial differential equations is discretization, and use of repeated differences; especially for using the discrete \Lap in approximating differential operators, and PDOs. See e.g., \cite{AtkinsonHan05}.

One picks a grid size \gd and then proceeds in steps:  
\begin{enumerate}
  \item Start with a partial differential operator, then study an associated discretized operator with the use of repeated differences on the \gd-lattice in \bRd.
  \item Solve the discretized problem for h fixed.
  \item As \gd tends to zero, numerical analysts evaluate the resulting approximation limits, and they bound the error terms.
\end{enumerate}

When discretization is applied to the Laplace operator in $d$ continuous variables, the result is our \Lap for the network $(\bZd,\cond)$; see \S\ref{sec:lattice-networks} for details and examples. However, when the same procedure is applied to a continuous Laplace operator on a Riemannian manifold, the discretized \Lap will be the network Laplacian on a suitable infinite network $(\Graph,\cond)$ which in general may have a much wilder geometry than \bZd.
  
This yields numerical algorithms for the solution of partial differential equations, and in the case of second order PDEs, the discretized operator is the discrete Laplacian studied in this investigation.
  
\end{remark}

\subsection*{Motivation and applications}
\label{sec:motivation-and-application}

\headerquote{A drunk man will eventually return home, but a drunk bird will lose its way in space.}{---~G.~Polya}

Applications to infinite networks of resistors serve as motivations, but our theorems have a wider scope, have other applications; and are, we believe, of independent mathematical interest. Our interest originates primarily from three sources.
\begin{enumerate}
  \item A series of papers written by Bob Powers in the 1970s which he introduced infinite systems of resistors into the resolution of an important question from quantum statistical mechanics in \cite{Pow75,Pow76a,Pow76b,Pow78, Pow79}.
  \item The pioneering work of Jun Kigami on the analysis of PCF self-similar fractals, viewing these objects as rescaled limits of networks; see \cite{Kig01}.
  \item Doyle and Snell's lovely book ``Random Walks and Electrical Networks'', which gives an excellent elementary introduction to the connections between resistance networks and random walks, including a resistance-theoretic proof of Polya's famous theorem on the transience of random walks in $\bZ^d$.
\end{enumerate}
Indeed, our larger goal is the cross-pollination of these areas, and we hope that the results of this book may be applicable to analysis on fractal spaces. A first step in this direction is given in Theorem~\ref{thm:R(x,y)-Lap-fractal}. To this end, a little more discussion of each of the above two subjects is in order.

Powers was interested in magnetism and the appearance of ``long-range order'', which is the common parlance for correlation between spins of distant particles; see \S\ref{sec:Magnetism-and-long-range-order} for a larger discussion. Consequently, he was most interested in graphs like the integer lattice \bZd (with edges between vertices of distance 1, and all resistances equal to 1), or other regular graphs that might model the atoms in a solid. Powers established a formulation of resistance metric that we adopt and extend in \S\ref{sec:effective-resistance-metric}, where we also show it to be equivalent to Kigami's formulation(s). Also, the proofs of Powers' original results on effective resistance metric contain a couple of gaps that we fill. In particular, Powers does not seem to have be aware of the possibility of nontrivial harmonic functions until \cite{Pow78}, where he mentions them for the first time. It is clear that he realized several immediate implications of the existence of such functions, but there more subtle (and just as important!) phenomena that are difficult to see without the clarity provided by Hilbert space geometry.

Powers studied an infinite graph \Graph by working with an exhaustion, that is, a nested sequence of finite graphs $\Graph_1 \ci \Graph_2 \ci \dots \ci \Graph_k \ci \Graph = \bigcup_k \Graph_k$. For example, $\Graph_k$ might be all the vertices of \bZd lying inside the ball of Euclidean radius $k$, and the edges between them. Powers used this approach to obtain certain inequalities for the resistance metric, expressing the consequences of deleting small subsets of edges from the network. Although he makes no reference to it, this approach is very analogous to Rayleigh's ``short-cut'' methods, as it is called in \cite{DoSn84}. 

Powers' use of an expanding sequence of graphs may be thought of as a ``limit in the large'' in contrast to the techniques introduces by Kigami, which may be considered ``limits in the small''. Self-similarity and scale renormalization are the hallmarks of the theory of fractal analysis as pursued by Kigami, Strichartz and others (see \cite{HKK02,Kig01, Kig03, Hut81, Str06, BHS05, Bea91, JoPed94, Jor04}, for example) but these ideas do not enter into Powers' study of resistors. One aim of the present work is the development of a Hilbert space framework suitable for the study of limits of networks defined by a recursive algorithm which introduces new vertices at each step and rescales the edges via a suitable contractive rescaling. As is known from, for example \cite{JoPed94,Jor06,Strichartz98,Str06,Teplyaev98}, there is a spectral duality between ``fractals in the large'', and ``fractals in the small''.

\pgap

\subsection*{The significance of Hilbert spaces}

\headerquote{`How large the World is!!' said the ducklings, when they found how much more room they now had compared to when they were confined inside the egg-shell. `Do you imagine this is the whole world?' asked the mother, `Wait till you have seen the garden; it stretches far beyond that to the parson's field, but I have never ventured to such a distance.'}{---~H.~C.~Andersen~(from~The~Ugly~Duckling)}

A main theme in this book is the use of Hilbert space technology in understanding metrics, potential theory, and optimization on infinite graphs, especially through finite-dimensional approximation. We emphasize those aspects that are intrinsic to \emph{infinite resistance networks,} and our focus is on \emph{analytic} aspects of graphs; as opposed to the combinatorial and algebraic sides of the subject, etc. Those of our results stated directly in the framework of graphs may be viewed as discrete analysis, yet the continuum enters via spectral theory for operators and the computation of probability of sets of infinite paths. In fact, we will display a rich variety of possible spectral types, considering the spectrum as a set (with multiplicities), as well as the associated spectral measures, and representations/resolutions.

Related issues for Hilbert space completions form a recurrent theme throughout our book. Given a resistance network, we primarily study three spaces of functions naturally associated with it: \HE, \HD, and to a lesser extent $\ell^2(\verts)$. Our harmonic analysis of functions on \Graph is studied via operators between the respective Hilbert spaces as discussed in  \S\ref{sec:H_energy-and-H_diss} and the Hilbert space completions of these three classes are used in an essential way. In particular, we obtain the boundary of the graph (a necessary ingredient of \eqref{eqn:intro:aaaagain} and the key to several mysteries) by analyzing the finite energy functions on \Graph which cannot be approximated by functions of finite support. 
\emph{However}, this metric space is naturally embedded inside the Hilbert space \HE, which is already complete by definition/construction. Consequently, the Hilbert space framework allows us to identify certain vectors as corresponding to the boundary of $(\Graph, \cond)$, and thus obtain a concrete understanding of the boundary.

However, the explicit representations of vectors in a Hilbert space completion (i.e., the completion of a pre-Hilbert space) may be less than transparent; see \cite{Yoo05}. In fact, this difficulty is quite typical when Hilbert space completions are used in mathematical physics problems. For example, in \cite{JoOl00, Jor00}, one begins with a certain space of smooth functions defined on a subset of $\bRd$, with certain support restrictions. In relativistic physics, one must deal with reflections, and there will be a separate positive definite quadratic form on each side of the ``mirror''.  As a result, one ends up with two startlingly different Hilbert space completions: a familiar $L^2$-space of functions on one side, and a space of distributions on the other. In \cite{JoOl00, Jor00}, one obtains holomorphic functions on one side of the mirror, and the space of distributions on the other side is spanned by the derivatives of the Dirac mass, each taken at the same specific point $x_0$.

It is the opinion of the authors that most interesting results of this book arise primarily from three things:
\begin{enumerate}
  \item differences between finite approximations to infinite networks, and how \& when these differences vanish in the limit, and 
  \item the phenomena that result when one works with a quadratic form whose kernel contains the constant functions, and
  \item the boundary (which is not a subset of the vertices) that naturally arises when a network supports nonconstant harmonic functions of finite energy, and how it explains other topics mentioned above.
\end{enumerate}
In classical potential theory, working modulo constant functions amounts to working with the class of functions satisfying $\|f'\|_2 < \iy$, but abandoning the $\ell^2$ requirement $\|f\|_2<\iy$. This has some interesting consequences, and the nontrivial harmonic functions play an especially important role; see Remark~\ref{rem:L2-loses-the-best-part}. What would one hope to gain by removing the $\ell^2$ condition?
\begin{enumerate}
  \item From the natural embedding of the metric space $(\Graph, R)$ into the Hilbert space \HE of functions of finite energy given by $x \mapsto v_x$, the functions $v_x$ are not generally in $\ell^2$. See Figure~\ref{fig:vx-in-Z1} of Example~\ref{exm:Z-not-bounded} for an illustration.
  \item The resistance metric does not behave nicely with respect to $\ell^2$ conditions. Several formulations of the resistance distance $R(x,y)$ involve optimizing over collections of functions which are not necessarily contained in $\ell^2$, even for many simple examples. 
  \item Corollary~\ref{thm:nontrivial-harmonic-fn-is-not-in-L2} states that nontrivial harmonic functions cannot lie in $\ell^2(\verts)$. Consequently, imposing an $\ell^2$ hypothesis removes the most interesting phenomena from the scope of study; see Remark~\ref{rem:L2-loses-the-best-part}.
\end{enumerate}
The infinite trees studied in Examples~\ref{exm:binary-tree:nontrivial-harmonic}--\ref{exm:a-forest-of-binary-trees} provide examples of these situations.

\pgap

\subsection*{Measures and measure constructions}
A reader glancing at our book will notice a number of incarnations of measures on infinite sample spaces: it may be a suitable space of paths (\S\ref{sec:path-space-of-a-general-random-walk}--\S\ref{sec:forward-harmonic-functions} and \S\ref{sec:Kolmogorov-construction-of-L2(gW,prob)}) or an analogue of the Schwartz space of tempered distributions (section \S\ref{sec:Gel'fand-triples-and-duality}). The latter case relies on a construction of ``Gel'fand triples'' from mathematical physics. The reader may wonder why they face yet another measure construction, but each construction is dictated by the problems we solve.
 Taking limits of finite subsystems is a universal weapon used with great success in a variety of applications; we use it here in the study of resistance distances on infinite graphs (\S\ref{sec:Resistance-metric-on-infinite-networks}); boundaries, boundary representations for harmonic functions (\S\ref{sec:Gel'fand-triples-and-duality}, \S\ref{sec:Dual-frames-and-the-energy-kernel}, and \S\ref{sec:path-space-of-a-general-random-walk}--\S\ref{sec:forward-harmonic-functions}); and equilibrium states and phase-transition problems in physics (\S\ref{sec:Kolmogorov-construction-of-L2(gW,prob)}--\S\ref{sec:GNS-construction}).

(1) \HE as an $L^2$ space. The central Hilbert space in this study, the energy space \HE, appears with a canonical reproducing kernel, but without any canonical basis, and there is no obvious way to see \HE as an $L^2(X,\gm)$ for some $X$ and \gm. Therefore, a major motivation for our measure constructions is just to be able to work with \HE as an $L^2$ space. In \S\ref{sec:Kolmogorov-construction-of-L2(gW,prob)}, we use a construction from probability to write $\HE = L^2(\gW,\gm)$ in a way that makes the energy kernel $\{v_x\}_{x \in \verts}$ into a system of (commuting) random variables. Here, \gW is an infinite Cartesian product of a chosen compact space \bS; one copy of \bS for each point $x \in \verts$.
  In \S\ref{sec:GNS-construction}, we use a non-commutative version of this probability technology: rather than Cartesian products, we will use infinite tensor products of $C^\ast$-algebras \sA, one for each $x \in \verts$.
  The motivation here is an application to a problem in quantum statistical mechanics.  The ``states'' on the $C^\ast$-algebra of all observables are the quantum mechanical analogues of probability measures in classical problems. Heuristically, the reader may wish to think of them as non-commutative measures; see e.g., \cite{BrRo97}.

(2) Boundary integral representation of harmonic functions. As it sometimes happens, the path to $\bd \Graph$ is somewhat circuitous: we begin with the discovery of an integral over the boundary, which leads us to understand functions on the boundary, which in turn points the way to a proper definition of the boundary itself. A closely related motivation for a measure is the formulation of an integral representation of harmonic functions $u \in \HE$: 
\linenopax
\begin{equation}\label{eqn:boundary-repn-for-harmonic-intro}
  u(x) = \int_{\Schw'} u(\gx) {h_x}(\gx) \, d\prob(\gx) + u(o).
\end{equation}
where $h_x = \Phar v_x$.
Thus the focus of \S\ref{sec:Gel'fand-triples-and-duality} is a formalization of the imprecise ``Riemann sums'' $u(x) = \sum_{\bd \Graph} u \dn{{h_x}} + u(o)$ of \S\ref{sec:relating-energy-form-to-Laplacian} as an integral of a bona fide measure. To carry this out, we construct a Gel'fand triple $\Schw \ci \HE \ci \Schw'$, where \Schw is a dense subspace of \HE and \Schw' is its dual, but with respect to a strictly finer topology. We are then able to produce a Gaussian probability measure \prob on $\Schw'$ and isometrically embed \HE into $L^2(\Schw',\prob)$. In fact, $L^2(\Schw',\prob)$ is the second quantization of \HE. However, the focus here is not on realizing \HE as an $L^2$ space (or subspace), but in obtaining the boundary integral representation of harmonic functions as in \eqref{eqn:boundary-repn-for-harmonic-intro}. Our aim is then to build formulas that allow us to compute values of harmonic functions $u \in \HE$ from an integral representation which yields $u(x)$ as an integral over $\bd G \ci \Schw'$. Note that this integration in \eqref{eqn:boundary-repn-for-harmonic-intro} is with respect to a measure depending on $x$ just as in the Poisson and Martin representations.

(3) Concrete representation of the boundary. We would like to realize $\bd \Graph$ as a measure space defined on a set of well-understood elements; this is the focus of the constructions in \S\ref{sec:the-boundary}. The goal is a measure on the space of all infinite paths in \Graph which yields the boundary $\bd\Graph$ in such a way that $\Graph \cup \bd \Graph$ is a compactification of \Graph which is compatible with the energy form \energy and the Laplace operator \Lap, and hence also the natural resistance metric on $(\Graph, \cond)$. This type of construction has been carried out with great success for the case of bounded harmonic functions (e.g., Poisson representation and the Fatou-Primalov theorem) and for nonnegative harmonic functions (e.g., Martin boundary theory), but our scope of enquiry is the harmonic functions of finite energy. Finally, we would like to use this Gaussian measure on $\Schw'$ to clarify $\bd \Graph$ as a subspace of $\Schw'$. Such a relationship is a natural expectation, as the analogous thing occurs in the work of Poisson, Choquet, and Martin.

\section*{What this book is not about}
\addcontentsline{toc}{subsection}{What this book is not about}

Many of the topics discussed in this book may appear to have been previously discussed elsewhere in the literature, but there are certain important subtleties which actually make our results quite different. This section is intended to clarify some of these.
\pgap

While there already is a large literature on electrical networks and on graphs (see e.g., \cite{CaW92, CaW07, DodziukKarp88, Dod06, DoSn84, Pow76b, CdV04, ChRi06, Chu07, FoKo07}, and the preprint \cite{Str08} which we received after the first version of this book was completed), we believe that our present operator/spectral theoretic approach suggests new questions and new theorems, and allows many problems to be solved in greater generality.

The literature on analysis on graphs breaks down into a variety of overlapping subareas, including: combinatorial aspects, systems of resistors on infinite networks, random-walk models, operator algebraic models \cite{DuJo08,Rae05}, probability on graphs (e.g., infinite particle models in physics \cite{Pow79}), Brownian motion on self-similar fractals \cite{Hut81}, Laplace operators on graphs, finite element-approximations in numerical analysis \cite{BSR08}; and more recently, use in internet-search algorithms \cite{FoKo07}. Even just the study of Laplace operators on graphs subdivides further, due to recently discovered connections between graphs and fractals generated by an iterated functions system (IFS); see e.g., \cite{Kig03, Str06}.

Other major related areas include discrete Schr\"{o}dinger operators in physics, information theory, potential theory, uses of the graphs in scaling-analysis of fractals (constructed from infinite graphs), probability and heat equations on infinite graphs, graph $C^\ad$-algebras, groupoids, Perron-Frobenius transfer operators (especially as used in models for the internet); multiscale analysis, renormalization, and operator theory of boundaries of infinite graphs (more current and joint research between the co-authors.) The motivating applications from \cite{Pow75,Pow76a,Pow76b,Pow78,Pow79} include the operator algebra of electrical networks of resistors (lattice models, $C^\ad$-algebras, and their representations), and more specifically, KMS-states from statistical mechanics. While working and presenting our results, we learned of even more such related research directions from experts working in these fields, and we are thankful to them all for taking the time to explain some aspects of them to us.

The main point here is that the related literature is \emph{vast} but our approach appears to be entirely novel and our results, while reminiscent of classical theory, are also new. We now elucidate certain specific differences.

\subsection*{Spectral theory}
The spectral theory for networks contrasts sharply with that for fractals, as is seen by considering the measures involved; they do not begin to become similar until one considers limits of networks. The spectrum of discrete Laplacians on infinite networks is typically continuous (lattices or trees provide examples, and are worked explicitly in \S\ref{sec:examples}). By contrast, in the analysis on fractals program of Kigami, Strichartz, and others, the Laplace operator has pure point spectrum; see \cite{Teplyaev98} in particular. The measures used in the analysis of networks are weighted counting measures, while the measures used in fractal analysis are based on the self-similar measures introduced by Hutchinson \cite{Hut81}. There is an associated and analogous entropy measure in the study of Julia sets; cf.~\cite{Bea91} and the recent work on Laplacians in \cite{RoTep:Basilica}.

Our approach differs from the extensive literature on spectral graph theory (see \cite{Chu01} for an excellent introduction, and an extensive list of further references) due to the fact that we eschew the $\ell^2$ basis for our investigations. We primarily study \Lap as an operator on \HE, and with respect to the energy inner product. The corresponding spectral theory is radically different from the spectral theory of \Lap in $\ell^2$. Most other work in spectral graph theory takes place in $\ell^2$, even implicitly when working with finite graphs: the adjoint of the drop operator (see Definition~\ref{def:drop-operator}) is taken with respect to the $\ell^2$ inner product and consequently violates Kirchhoff's laws. In fact, the discussion preceding \cite[(2.2)]{Woess00} shows how this version of the adjoint is incompatible with Kirhhoff's Law as mentioned in the summary of \S\ref{sec:H_energy-and-H_diss} just above. Additionally, \cite{Chu01} and others work with the spectrally renormalized Laplacian $\Lap_{\spectral} := \cond^{-1/2} \Lap \cond^{-1/2}$. However, $\Lap_{\spectral}$ is a bounded Hermitian operator (with spectrum contained in $[0,2]$) and so is unsuitable for our investigations of $\bd \Graph$ based on defect indices, etc.

As we have only encountered relatively few instances where the complete details are worked out for spectral representations in the framework of \emph{discrete analysis}, we have attempted to provide several explicit examples. These are likely folkloric, as the geometric possibilities of graphs are vast, and so is the associated range of spectral configurations. A list of recent and past papers of relevance includes \cite{Str08, Car72, Car73a, Car73b, ChRi06, Chu07, CdV99, CdV04, Jor83}, and Wigner's original paper on the semicircle law \cite{Wig55}. The present investigation also led to a spectral analysis of the binary tree from the perspective of dipoles in \cite{DuJo08}; this study discovered that the spectrum of \Lap on the binary tree is also given by Wigner's semicircle law.

There is also a literature on infinite/transfinite networks and generalized Kirchhoff laws using nonstandard analysis, etc., see \cite{Zem91,Zem97}. However, this context allows for edges with resistance 0, which we do not allow (for physical as well as theoretical reasons). One can neglect the resistance of wires in most engineering applications, but not when considering infinite networks (the epsilons add up!). The resulting theory therefore diverges rapidly from the observations of the present book; according to our definitions, all networks support currents satisfying Kirchhoff's law, and in particular, all induced currents satisfy Kirchhoff's law.

\subsection*{Operator algebras}
There are also recent papers in the literature which also examine graphs with tools from operator algebras and infinite determinants. The papers \cite{GIL06a, GIL06b, GIL06c} by Guido et al are motivated by questions for fractals and study the detection of periods in infinite graphs with the use of the Ihara zeta function, a variant of the Riemann zeta function. There are also related papers with applications to the operator algebra of groupoids \cite{Cho08, FMY05}, and the papers \cite{BaMo00, BaMo01} which apply infinite graphs to the study of quasi-periodicity in solid state physics. However, the focus in these papers is quite different from ours, as are the questions asked and the methods employed. While periods and quasi-periods in graphs play a role in our present results, they enter our picture in quite different ways, for example via spectra and metrics that we compute from energy forms and associated Laplace operators. There does not seem to be a direct comparison between our results and those of Guido et al.

\subsection*{Boundaries of graphs}
There is also no shortage of papers studying boundaries of infinite graphs: \cite{PicWoess90,Saw97,Woess00} discuss the Martin boundary, \cite{PicWoess90,Woess00} also describe the more geometrically constructed  ``graph ends'', and \cite{Car72, Car73a, Car73b} use unitary representations. There are also related results in \cite{CdV99, CdV04} and \cite{Kaimanovich98, Kaimanovich92:DirNorms, Kaimanovich92:MeasBd, KaimanovichWoess02}
\version{}{\marginpar{What are the connections? We need to determine if there is a relation between our boundary and Martin, etc, and give an example if not.}}
While there are connections to our study, the scope is different.

Martin boundary theory is really motivated by constructing a boundary for a Markov process, and the geometry/topology of the boundary is rather abstract and a bit nebulous. Additionally, one needs a Green's function, and it must satisfy certain hypotheses before the construction can proceed. Furthermore, the focus of Martin boundary theory is the nonnegative harmonic functions. Our boundary construction is more general in that it applies to any electrical network as in Definition~\ref{def:ERN} and it remains correct for all harmonic functions of finite energy, including constant functions and harmonic functions which change sign. However, it is also more restrictive in the sense that a resistance network may support functions which are bounded below but do not have finite energy. 

We should also point out that our boundary construction is related to, but different from, the ``graph ends'' introduced by Freudenthal and others. The ends of a graph are the natural discrete analogue of the ends of a minimal surface (usually assumed to be embedded in $\bR^3$), a notion which is closely related to the conformal type of the surface. Starting with the central book \cite{Woess00} by Wolfgang Woess, the following references will provide the reader with an introduction to the study of harmonic functions on infinite networks and the ends of graphs and groups: \cite{Wo86, Wo87, Wo89}, and \cite{Wo95a} on Martin boundaries, \cite{PicWoess90} on ends, \cite{Wo96} on Dirichlet problems, \cite{Wo95b} on random walk. A comparison of the examples in \S\ref{sec:lattice-networks} and \S\ref{sec:tree-networks} illustrates that varying the resistances produces dramatic changes in the topology of the boundary.

Our boundary essentially consists of infinite paths which can be distinguished by harmonic functions; 
see \S\ref{sec:bdG-as-equivalence-classes-of-paths} for details. It follows that transient networks with no nontrivial harmonic functions have exactly one boundary point (corresponding to the unique monopole). In particular, the integer lattices $(\bZd,\one)$ have precisely 1 boundary point for $d \geq 3$, and have 0 boundary points for $d=1,2$. The Martin boundary of $(\bZ^2,\one)$ consists of two points; similarly, $(\bZ^2,\one)$ has two graph ends; cf.~\cite{PicWoess90}.

\section*{General remarks}

\begin{rem}[Real- and complex-valued functions] 
  Throughout the introductory discussion of resistance networks in \S\ref{sec:electrical-resistance-networks}--\S\ref{sec:effective-resistance-metric}, we discuss collections of real-valued functions on the vertices or edges of the graph \Graph. Such objects are most natural for the heuristics of the physical model, and additionally allow for induced orientation/order and make certain probabilistic arguments possible. However, in the latter portions of this book, we need to incorporate complex-valued functions into the discussion in order to make full use of spectral theory and other methods.
\end{rem}

\begin{rem}[Symbols glossary]
  For the aid of the reader, we have included a list of symbols and abbreviations used in this document. Wherever possible, we have attempted to ensure that each symbol has only one meaning. In cases of overlap, the context should make things clear. In Appendix~\ref{sec:navigation-aids}, we also include some diagrams which we hope clarify the properties of the many operators and spaces we discuss, and the relations between them.
\end{rem}

\section*{Acknowledgements}
While working on the project, the co-authors have benefitted from interaction with colleagues and students. We thank everyone for generously suggesting improvements as our book progressed.  The authors are grateful for stimulating comments, helpful advice, and valuable references from John Benedetto, Donald Cartwright, Il-Woo Cho, Raul Curto, Dorin Dutkay, Alexander Grigor'yan, Dirk Hundertmark, Richard Kadison, Keri Kornelson, Michel Lapidus, Russell Lyons, Diego Moreira, Peter M\"{o}rters, Paul Muhly, Massimo Picardello, Bob Powers, Marc Rieffel, Karen Shuman, Sergei Silvestrov, Jon Simon, Myung-Sin Song, Bob Strichartz, Andras Telcs, Sasha Teplyaev, Elmar Teufl, Ivan Veselic, Lihe Wang, Wolfgang Woess, and Qi Zhang. The authors are particularly grateful to Russell Lyons for several key references and examples, and to Jun Kigami for several illuminating conversations and for suggesting the approach in \eqref{eqn:Schur-complement-as-sum}. Initially, the first named author (PJ) learned of discrete potential theory from Robert T. Powers at the University of Pennsylvania in the 1970s, but interest in the subject has grown exponentially since.

%% file: electrical-resistance-networks.tex

\chapter{Resistance networks}
\label{sec:electrical-resistance-networks}

\headerquote{The excitement that a gambler feels when making a bet is equal to the amount he might win times the probability of winning it}{---~B.~Pascal}

Resistance networks are the basic object of study throughout this volume; the basic idea is that a graph with weighted edges makes a good discrete model for diffusions, when the weights are interpreted as ``sizes'' or ``capacities'' for transfer, in some sense. Such a model is useful for understanding the flow of heat in perforated media, diffusion of water in porous matter, or the transfer of data through the internet. However, due to its intuitive appeal and historical precedent, we have chosen to stick predominantly with the metaphor of electricity flowing through a network of conductors. In this situation, the weights correspond to conductances (recall that conductance is the reciprocal of resistance), functions on the vertices may be interpreted as voltages, and corresponding functions on the edges of the graph may be interpreted as currents. This context also provides a natural interpretation for the energy \energy which will be central to our study: if $v$ is a function on the vertices of the graph (i.e., a voltage), then $\energy(v)$ is a number representing the potential energy of this configuration, equivalently, the power dissipated by the electrical current induced by $v$. 

\section{The \ERN model}
This section contains the basic definitions used throughout the sequel; we introduce the mathematical model of an \ERN (RN) 
  \glossary{name={ERN},description={electrical resistance network},sort=E}
as a graph \Graph whose edges are understood as conductors and whose vertices are the nodes at which these resistors are connected. The conductance data is specified by a function \cond, so that $\cond(x,y)$ is the conductance of the edge (resistor) between the vertices $x$ and $y$. With the network data $(\Graph,\cond)$ fixed, we begin the study of functions defined on the vertices. We define many basic terms and concepts used throughout the book, including the Dirichlet energy form \energy and the Laplace operator \Lap. Additionally, we prove a key identity relating \energy to \Lap for finite graphs: Lemma~\ref{thm:E(u,v)=<u,Lapv>}. In Theorem~\ref{thm:E(u,v)=<u,Lapv>+sum(normals)}, this will be extended to infinite graphs, in which case it is a discrete analogue of the familiar Gauss-Green identity from vector calculus. The appearance of a somewhat mysterious boundary term in the Theorem~\ref{thm:E(u,v)=<u,Lapv>+sum(normals)} prompts several questions which are discussed in Remark~\ref{rem:boundary-term}. Answering these questions comprises a large part of the sequel; cf.~\S\ref{sec:the-boundary}. In fact, Theorem~\ref{thm:E(u,v)=<u,Lapv>+sum(normals)} provides much of the motivation for energy-centric approach we pursue throughout our study; the reader may wish to look ahead to Remark~\ref{rem:L2-loses-the-best-part} for a preview.

\begin{defn}\label{def:graph}
  A graph $\Graph = \{\verts, \edges\}$ is given by the set of vertices \verts and the set of edges $\edges \ci \verts \times \verts$. Two vertices are \emph{neighbours} (or are \emph{adjacent}) iff there is an edge $(x,y) \in \edges$ connecting them, and this is denoted $x \nbr y$. This relation is symmetric, as $(y,x) \in \edges$ whenever $(x,y) \in \edges$. The set of neighbours of $x \in \verts$ is
\linenopax
  \begin{equation}\label{eqn:def:graph-neighbours}
    \vnbd(x) = \{y \in \verts \suth y \nbr x\}.
  \end{equation}
\end{defn}
  \glossary{name={\Graph},description={graph},sort=G,format=textbf}
  \glossary{name={\verts},description={vertices of a graph},sort=G0,format=textbf}
  \glossary{name={\edges},description={edges of a graph},sort=G1,format=textbf}
  \glossary{name={$\nbr$},description={$x \nbr y$ means $(x,y)$ is an edge},sort=~,format=textbf}

In our context, the set of edges of \Graph will be determined by the conductance function, so that all graph data is implicitly provided by \cond.

\begin{defn}\label{def:conductance}
  The \emph{conductance} $\cond_{xy}$ is a symmetric function
  \linenopax
  \begin{equation}\label{eqn:def:conductance}
    \cond:\verts \times \verts \to [0,\iy),
  \end{equation}
  in the sense that $\cond_{xy} = \cond_{yx}$.
  \glossary{name={$\cond_{xy}$},description={conductance of the edge $(x,y)$},sort=c,format=textbf}
  It is our convention that $x \not\nbr y$ if and only if $\cond_{xy}=0$; that is, there is an edge $(x,y) \in \edges$ if and only if $0 < \cond(x,y) < \iy$.
\end{defn}

Conductance is the reciprocal of resistance, and this is the origin of the name ``resistance network''. It is important to note that $\cond_{xy}^{-1}$ gives the resistance between \emph{adjacent} vertices; this feature distinguishes $\cond_{xy}^{-1}$ from the \emph{effective resistance} $R(x,y)$ discussed later, for which $x$ and $y$ need not be adjacent.

\begin{defn}\label{def:conductance-measure}
  The conductances define a measure or weighting on \verts by
  \linenopax
  \begin{equation}\label{eqn:def:cond-meas}
    \cond(x) := \sum_{y \nbr x} \cond_{xy}.
  \end{equation}
  Whenever \Graph is connected, it follows that $\cond(x) > 0$, for all $x \in \verts$. The notation \cond will also be used, on occasion, to indicate the multiplication operator $(\cond v)(x) := \cond(x) v(x)$.
  \glossary{name={$\cond(x)$},description={sum of conductances of edges incident on $x$},sort=c,format=textbf}
  \glossary{name={\cond},description={conductance (multiplication) operator},sort=c,format=textbf}
\end{defn}

\begin{defn}\label{def:paths}
  A \emph{path} \cpath from $\ga \in \verts$ to $\gw \in \verts$ is a sequence of adjacent vertices $(\ga = x_0, x_1, x_2, \dots, x_n = \gw)$, i.e., $x_i \nbr x_{i-1}$ for $i=1,\dots,n$. The path is \emph{simple} if any vertex appears at most once (so that a path is simply connected).
\end{defn}
  \glossary{name={\cpath},description={path},sort=G,format=textbf}

\begin{defn}\label{connected}
  A graph \Graph is \emph{connected} iff for any pair of vertices $\ga,\gw \in \verts$, there exists a finite path \cpath from \ga to \gw.
\end{defn}

\begin{remark}\label{rem:def:connected-ERN}
  Note that for resistors connected in series, the resistances just add, so this condition implies there is a path of finite resistance between any two points. We emphasize that \emph{all graphs and subgraphs considered in this study are connected.}
\end{remark}

At this point, the reader may wish to peruse some of the examples of \S\ref{sec:examples}.

\begin{defn}\label{def:ERN}
  An \emph{\ERN} is a connected graph $(\Graph,\cond)$ whose conductance function satisfies $\cond(x) < \iy$ for every $x \in \verts$. We interpret the edges as being defined by the conductance: $x \nbr y$ iff $c_{xy}>0$.
\end{defn}

Note that \cond need not be bounded in Definition~\ref{def:ERN}. Also, we will typically assume an RN to be simple in the sense that there are no self-loops, and there is at most one edge from $x$ to $y$. This is mostly for convenience: basic electrical theory says that two conductors $\cond^1_{xy}$ and $\cond^2_{xy}$ connected in parallel can be replaced by a single conductor with conductance $\cond_{xy} = \cond^1_{xy} + \cond^2_{xy}$. Also, electric current will never flow along a conductor connecting a node to itself. Nonetheless, such self-loops may be useful for technical considerations: one can remove the periodicity of a random walk by allowing self-loops. This can allow one to obtain a ``lazy walk'' which is ergodic, and hence amenable to application of tools like the Perron-Frobenius Theorem. See, for example, \cite{LevPerWil08, Lyons:ProbOnTrees, AlFi09}.

We will be interested in certain operators that act on functions defined on \ERNs.

\begin{defn}\label{def:graph-laplacian}
  The \emph{Laplacian} on \Graph is the linear difference operator 
  which acts on a function $v:\verts \to \bR$ by
  \linenopax
  \begin{equation}\label{eqn:def:laplacian}
    (\Lap v)(x) :
    = \sum_{y \nbr x} \cond_{xy}(v(x)-v(y)).
  \end{equation}
  \glossary{name={\Lap},description={graph Laplacian},sort=D,format=textbf}

  A \fn $v:\verts \to \bR$ is called \emph{harmonic} iff $\Lap v \equiv 0$.
\end{defn}

\begin{defn}\label{def:graph-transfer-operator}
  The \emph{transfer operator} on \Graph is the linear operator \Trans which acts on a function $v:\verts \to \bR$ by
  \linenopax
  \begin{equation}\label{eqn:def:transfer-operator}
    (\Trans v)(x): = \sum_{y \nbr x} \cond_{xy} v(y).
  \end{equation}
  Hence, the Laplacian may be written $\Lap = \cond - \Trans$, where $(\cond v)(x) := \cond(x) v(x)$.
  \glossary{name={\Trans},description={transfer operator},sort=T,format=textbf}
\end{defn}

We won't worry about the domain of \Lap or \Trans until \S\ref{sec:Lap-on-HE}. For now, consider both of these operators as defined on any function $v:\verts \to \bR$. The reader familiar with the literature will note that the definitions of the Laplacian and transfer operator given here are normalized differently than may be found elsewhere in the literature. For example, \cite{DoSn84} and other probabilistic references use
\linenopax
\begin{align}
  \Lapc := \cond^{-1} \Lap = \one - \Prob, 
  \q\text{so}\q   
  (\Lap_\cond v)(x)
  := \frac{1}{\cond(x)} \sum_{y \nbr x} \cond_{xy}(v(x)-v(y)), 
    \label{eqn:probablistically-reweighted-laplacian}
\end{align}
where $\Prob := \cond^{-1} \Trans$ is the probabilistic transition operator 
\glossary{name={\Prob},description={probabilistic transition operator},sort=P,format=textbf}
corresponding to the transition probabilities $p(x,y) = \cond_{xy}/\cond(x)$.
For another example, \cite{Chu01} and other spectral-theoretic references use
\linenopax
\begin{align}
  \Lap_\spectral := \cond^{-1/2} \Lap \cond^{-1/2} = \one - \cond^{-1/2} \Trans \cond^{-1/2}, 
  \q\text{so}\q
  (\Lap_\spectral v)(x)
  := v(x) - \sum_{y \nbr x} \frac{\cond_{xy}v(y)}{\sqrt{\cond(y)}} .
    \label{eqn:spectrally-reweighted-laplacian}\end{align}
  \glossary{name={$\Lap_\cond$},description={probabilistically renormalized Laplacian operator},sort=D,format=textbf}
  \glossary{name={$\Lap_\spectral$},description={spectrally renormalized Laplacian},sort=D,format=textbf}

However, these renormalized version are much more awkward to work with in the present context; especially when dealing with the inner product and kernels of the Hilbert spaces we shall study. Not only are \eqref{eqn:def:laplacian} and \eqref{eqn:def:transfer-operator} are better suited to the \ERN framework (as will be evinced by the operator theory developed in \S\ref{sec:energy-Hilbert-space} and succeeding sections) but both $\Lap_\cond$ and $\Lap_\spectral$ are bounded operators, and hence do not allow for the delicate spectral analysis carried out in \S\ref{sec:the-boundary}--\S{sec:Lap-on-HE}. 

\section{The energy}
\label{sec:energy}
In this section we study the relation between the energy \energy and Laplacian \Lap on finite networks, as expressed in Lemma~\ref{thm:E(u,v)=<u,Lapv>}. This formula will be used prolifically, as it also holds on infinite networks in many circumstances. In fact, a noticeable portion of \S\ref{sec:energy-Hilbert-space} is devoted to determining when this is so.

\begin{defn}\label{def:graph-energy}
  The \emph{graph energy} of an \ERN is the quadratic form defined for \fns $u:\verts \to \bR$ by
  \linenopax
  \begin{equation}\label{eqn:def:graph-energy}
    \energy(u)
    := \frac12 \sum_{x,y \in \verts} \cond_{xy}(u(x)-u(y))^2.
  \end{equation}
  \glossary{name={\energy},description={(bilinear or quadratic) Dirichlet energy form},sort=E,format=textbf}
  There is also the associated bilinear \emph{energy form}
  \linenopax
  \begin{align}\label{eqn:def:energy-form}
    \energy(u,v)
    :=& \frac12 \sum_{x,y \in \verts} \cond_{xy}(u(x)-u(y))(v(x)-v(y)).
  \end{align}
  For both \eqref{eqn:def:graph-energy} and \eqref{eqn:def:energy-form}, note that $\cond_{xy}=0$ for vertices which are not neighbours, and hence only pairs for which $x \nbr y$ contribute to the sum; the normalizing factor of $\frac12$ corresponds to the idea that each edge should only be counted once. The \emph{domain of the energy} is
  \linenopax
  \begin{equation}\label{eqn:def:energy-domain}
    \dom \energy = \{u:\verts \to \bR \suth \energy(u)<\iy\}.
  \end{equation}
\end{defn}
  \glossary{name={$\dom\energy$},description={functions of finite energy; domain of the energy form},sort=d,format=textbf}

The close relationship between the energy and the conductances is highlighted by the simple identities
\linenopax
\begin{align}\label{eqn:energy-of-Diracs}
  \energy(\gd_x) = \cond(x),
  \qq \text{and}\qq
  \energy(\gd_x,\gd_y) = -\cond_{xy},
\end{align}
where $\gd_x$ is a (unit) Dirac mass at $x \in \verts$. The easy proof is left as an exercise. A significant upshot of \eqref{eqn:energy-of-Diracs} is that the Dirac masses are not orthogonal with respect to energy.

\begin{remark}\label{rem:E(u)=0-iff-u=const}
  It is immediate from \eqref{eqn:def:graph-energy} that $\energy(u)=0$ if and only if $u$ is a constant function. The energy form is positive semidefinite, but if we work modulo constant functions, it becomes positive definite and hence an inner product. We formalize this in Definition~\ref{def:The-energy-Hilbert-space} and again in \S\ref{sec:vonNeumann's-embedding-thm}. In classical potential theory (or Sobolev theory), this would amount to working with the class of functions satisfying $\|f'\|_2 < \iy$, but abandoning the requirement that $\|f\|_2<\iy$. As a result of this, the nontrivial harmonic functions play an especially important role in this book. In particular, it is precisely the presence of nontrivial harmonic functions which prevents the functions of finite support from being dense in the space of functions of finite energy; see \S\ref{sec:The-role-of-Fin-in-HE}.

  Traditionally (e.g., \cite{Kat95,FOT94}) the study of quadratic forms would combine $\energy(u,v)$ and $\la u,v\ra_{\ell^2}$. In our context, this is counterproductive, and would eclipse some of our most interesting results. Some of our most intriguing questions for elements $v \in \HE$ involve boundary considerations, and in these cases $v$ is not in $\ell^2(\verts)$ (Corollary~\ref{thm:nontrivial-harmonic-fn-is-not-in-L2}). One example of this arises in the discrete Gauss-Green formula (Theorem~\ref{thm:E(u,v)=<u,Lapv>+sum(normals)}); another arises in study of forward-harmonic functions in \S\ref{sec:forward-harmonic-functions}.
\end{remark}

\version{}{
\begin{defn}\label{def:closed-form}
  \marginpar{Add this part after completing the proof of the lemma.}
  A symmetric form $E$ on a Hilbert space $H$ is \emph{closed} iff whenever $\{u_n\}$ is Cauchy in $E$, one has the implication
    \linenopax
    \begin{align}\label{eqn:def:closed-form}
        \lim_{n \to \iy} \|u_n\| \to 0 
        \q \implies \q
        \lim_{n \to \iy} E(u_n,u_n) = 0.
    \end{align}  
  A sequence $\{u_n\}$ is Cauchy in $E$ iff $\lim_{n,m \to \iy} E(u_n-u_m,u_n-u_m)=0$. \end{defn}

Actually, Definition~\ref{def:closed-form} is the definition for a symmetric form to be \emph{closable}, but any closable symmetric form has a closed extension, and we may always pass to this extension. Thus, Definition~\ref{def:closed-form} suffices for our purposes; see \cite[p.4]{FOT94}.

\begin{lemma}\label{thm:E-is-closed}
  \energy is a closed form.
  \begin{proof}
    Needed. Something telescopes.
    \linenopax
    \begin{align*}
      \energy
    \end{align*}
  \end{proof}
\end{lemma}
}

The following proposition may be found in \cite[\S1.3]{Str06} or \cite[Ch.~2]{Kig01}, for example.

\begin{prop}\label{prop:energy-properties}
  The following properties are readily verified:
  \begin{enumerate}
    \item $\energy(u,u) = \energy(u)$.
    \item (Polarization) $\energy(u,v) = \frac14 [\energy(u+v) - \energy(u-v)]$.
    \item (Markov property) $\energy([u]) \leq \energy(u)$, where $[u]$ is any contraction of $u$.
  \end{enumerate}
\end{prop}
For example, let $[u] := \min\{1,\max\{0,u\}\}$. The following result relates the Laplacian to the graph energy on finite networks, and can be interpreted as a relation between $\dom \energy$ and $\ell^2(\verts)$. 

\begin{lemma}\label{thm:E(u,v)=<u,Lapv>}
  Let \Graph be a finite \ERN. For $u,v \in \dom \energy$,
  \linenopax
  \begin{equation}\label{eqn:E(u,v)=<u,Lapv>}
    \energy(u,v)
    = \sum_{x \in \verts} u(x) \Lap v(x)
    = \sum_{x \in \verts} v(x) \Lap u(x).
  \end{equation}
  \begin{proof}
    Direct computation yields
    \linenopax
    \begin{align}
      \energy(u,v)
      &= \frac12 \sum_{x,y \in \verts}\cond_{xy} \left(u(x)v(x) - \vstr[2.1] u(x)v(y) - u(y)v(x) + u(y)v(y)\right) \notag \\
      &= \frac12 \sum_{x \in \verts} \cond(x)u(x)v(x) + \frac12 \sum_{y \in \verts} \cond(y)u(y)v(y) \notag \\
      &\hstr[12]- \frac12 \sum_{x \in \verts_n} u(x) \Trans v(x) - \frac12 \sum_{y \in \verts} u(y) \Trans v(y) \notag \\
      &= \sum_{x \in \verts} \cond(x) u(x)v(x) - \sum_{x,y \in \verts} u(x) \Trans v(x) \notag \\
      &= \sum_{x \in \verts} u(x) \left(\cond(x) v(x) - \Trans v(x)\right) \notag \\
      &= \sum_{x \in \verts} u(x) \Lap v(x).
    \end{align}
    Of course, the computation is identical for $\sum_{x \in \verts} v(x) \Lap u(x)$.
  \end{proof}
\end{lemma}

We include the following well-known result for completeness.
\begin{cor}\label{thm:harmonic=const-on-finite}
  On a finite \ERN, all harmonic functions of finite energy are constant.
  \begin{proof}
    If $h$ is harmonic, then $\energy(h) = \sum_{x \in \verts} h(x) \Lap h(x)= 0$. See Remark~\ref{rem:E(u)=0-iff-u=const}.
  \end{proof}
\end{cor}

Connectedness is implicit in the calculations of both Lemma~\ref{thm:E(u,v)=<u,Lapv>} and Corollary~\ref{thm:harmonic=const-on-finite}; recall that \textit{all} \ERNs considered in this work are connected.
We will extend Lemma~\ref{thm:E(u,v)=<u,Lapv>} to infinite graphs in Theorem~\ref{thm:E(u,v)=<u,Lapv>+sum(normals)}, where the formula is more complicated:
\linenopax
\begin{align*}
  \energy(u,v) = \sum_{x \in \verts} u(x) \Lap v(x) + \{\text{``boundary term''}\}.
\end{align*}
It is shown in Theorem~\ref{thm:TFAE:Fin,Harm,Bdy} that the presence of the boundary term corresponds to the transience of the random walk on the underlying network. 
In fact, one can interpret Corollary~\ref{thm:harmonic=const-on-finite} as the reason why the boundary term alluded to above vanishes on finite networks. We study the interplay between \energy and \Lap further in \S\ref{sec:the-transfer-operator}--\S\ref{sec:weighted-spaces}.

\section{Remarks and references}
\label{sec:remarks-and-references-1}

Of the cited references for this chapter, some are more specialized. However for prerequisite material (if needed), the reader may find the book \cite{DoSn84} by Doyle and Snell especially relevant. It exists in several editions, and is available for free on the arXiv (\href{http://arxiv.org/abs/math/0001057}{math/0001057}). While it is a gold mine of ideas and illuminating examples, and is accessible to undergraduates. 

We have been much inspired by Doyle and Snell's book on electrical networks \cite{DoSn84}, and by Jun Kigami's work effective resistance and discrete potential theory (especially as it pertains to on renormalization and scaling limits) \cite{Kig01, Kig03, Kig08}. We are similarly indebted to Wolfgang Woess' book \cite{Woess00}, covering probability and analysis on infinite networks, Markov chains, and especially the theory of boundaries, as developed in \cite{Wo86, Wo87, Wo89, Wo95a, Wo95b, Wo96, Thomassen90} and elsewhere. The reader will find \cite{AlFi09, Lyons:ProbOnTrees, LevPerWil08} to be excellent references for the random walks on graphs and Markov chains in general (with an emphasis on \emph{reversible} chains). The main themes in this and later chapters are also tangentially related to the fascinating work by Fan Chung on spectral theory of transfer operators on infinite graphs \cite{Chu07, ChRi06}. 


Since our first chapter serves in part as an overview of material in the book, and some results in the literature but not in our book, there are quite a number of papers and books that are appropriate to cite, and here is a partial list:
\cite{Woess00, CaWo04, Woe03, KaimanovichWoess02, Kig03, Kig01, Tera78, PrigozyWeinberg76, Sve56, Cra52, Sch91}.

In the subsequent end-of chapter sections we will discuss also pioneering work by: 
\begin{itemize}
  \item Aldous and Fill, reversible markov chains and random walks on graphs \cite{AlFi09},
  \item Benjamini, Lyons, Pemantle, Peres, and Schramm (separately or in various combinations), effective resistance, probability on trees, percolation, analysis and probability on infinite graphs \cite{Lyons:ProbOnTrees, LevPerWil08, Peres99, BLPS01, BLPS99, BLS99, LPS03, ALP99, LPS06, LyPe03, LPP96, NaPe08a, Lyo03},
  \item Cartwright, random walks, Dirichlet functions and spectrum \cite{CaSoWo93, CaW92, CaWo04, CaW07}
  \item Chung, spectral theory of transfer operators on infinite graphs \cite{Chu07, ChRi06},
  \item Doob, martingales, and probabilistic boundaries \cite{Doob53, Doob55, Doob58, Doob59},
  \item Doyle and Snell, electrical networks \cite{DoSn84}, and Doyle \cite{Doyle88},
  \item Hida, use of Hilbert space methods in stochastic integration \cite{Hida80},
  \item Kigami, effective resistance and discrete potential theory \cite{Kig01, Kig03, Kig08},
  \item Kolmogorov, foundations of probability theory \cite{Kol56},
  \item Liggett, infinite spin-models \cite{Lig93, Lig95, Lig99},
  \item von Neumann, the theory of unbounded operators, quantum mechanics, and metric geometry \cite{vN32a, vN32b, vN32c},
  \item R. T. Powers, use of resistance distance in the estimation of long-range order in quantum statistical models \cite{Pow75, Pow76a, Pow76b, Pow78, Pow79},
  \item Saloff-Coste, harmonic analysis and probability and random walks in relation to groups \cite{Saloff-Woess06, Saloff-Woess09},
  \item Manfred Schroeder, harmonic analysis and signal processing on fractals \cite{Sch91},
  \item P. M. Soardi, harmonic analysis and potential theory on infinite graphs \cite{Soardi94}, a substantial influence,
  \item Frank Spitzer, random walk \cite{Spitzer},
  \item Dan Stroock, Markov processes \cite{Stroock},
  \item A. Telcs, random walks, graphs and fractals \cite{Telcs06a, Telcs06b, Telcs03, Telcs01},
  \item G. George Yin, Qing Zhang, use of stochastic integration in renormalization theory \cite{YiZh05}.
\end{itemize}


While we present a number of theorems related in one way or the other to earlier results, the material is developed here from simple axioms and from a unifying point of view: we make use of fundamental principles in the theory of operators in Hilbert space. Using this we develop a variety of results on networks, on scaling relations, on renormalization, two spin-models, long-range order, and on discrete potential theory.  Our aim and emphasis is to develop the material from first principles: Riesz duality, reproducing kernels, metric embedding into Hilbert space, and stochastic integral models. As a bonus, we are able to point out how basic principles from operator theory lead to unification of a variety of existing results, and in some cases in their extension.

%% file: potentials.tex

\chapter[Currents and potentials]{Currents and potentials on resistance networks}
\label{sec:currents-and-potentials}

\headerquote{While electrical networks are only a different language for reversible Markov chains, the electrical point of view is useful because of the insight gained from the familiar physical laws of electrical networks.}{---~Y.~Peres}



\section{Currents on \ERNs}
\label{sec:currents}

The potential theory for an \ERN is studied via an experiment in which 1 amp of current is passed through the network, inserted into one vertex and extracted at some other vertex. The voltage drop 
measured between the two nodes is the effective resistance between them, see \S\ref{sec:effective-resistance-metric}. 

When the voltages are fixed at certain vertices, it induces a current in the network in accordance with the laws of Kirchhoff and Ohm. This \emph{induced current} is introduced formally in Definition~\ref{def:induced-current}. Induced currents are important for studying flows of minimal dissipation, and will also be useful in the study of forward-harmonic functions in \S\ref{sec:forward-harmonic-functions}. If a voltage drop of 1 volt is imposed between two vertices, the effective resistance between these two vertices is the reciprocal of the dissipation of the induced current.

In Theorem~\ref{thm:Pot-is-nonempty-by-current-flows} we show that there always exists an harmonic function satisfying the boundary conditions implied by the above described experiment, in order to fill a gap in \cite{Pow76b}. In Theorem~\ref{thm:minimal-flows-are-induced-by-minimizers} and  Theorem~\ref{thm:minimal-flows-are-induced-by-minimizers} it is shown that these harmonic functions correspond to currents which minimize energy dissipation.

\begin{defn}\label{def:currents}
  A \emph{current} is a skew-symmetric function $\curr : \verts \times \verts \to \bR$.
  \glossary{name={\curr},description={current; a skew-symmetric function on edges},sort=i,format=textbf}
\end{defn}

\begin{defn}\label{def:orientation}
  An \emph{orientation} is a subset of the edges which includes exactly one of each pair $\{(x,y), (y,x)\}$. For a given current \curr, one may pick an orientation by requiring that $\curr(x,y) > 0$ on every edge for which \curr is nonzero, and arbitrarily choosing $(x,y)$ or $(y,x)$ outside the support of \curr. We refer to this as an \emph{orientation induced by the current}; this will be used extensively in \S\ref{sec:forward-harmonic-functions} to study the forward-harmonic functions.
\end{defn}

The energy is a functional defined on functions $v:\verts \to \bR$  which give voltages between different vertices in the network. The associated notion defined on the edges of the network is the dissipation of a current.

\begin{defn}\label{def:dissipation}
  The \emph{dissipation} of a current may be thought of as the energy lost as a current flows through an \ERN. More precisely, for $\curr,\curr_1,\curr_2:\edges \to \bR$,
  \linenopax
  \begin{equation}\label{eqn:def:dissipation}
    \diss(\curr) := \frac12 \sum_{(x,y) \in \edges} \cond_{xy}^{-1} \curr(x,y)^2.
  \end{equation}
  The associated bilinear form is the \emph{dissipation form}:
  \linenopax
  \begin{equation}\label{eqn:def:dissipation-form}
    \diss(\curr_1,\curr_2) := \frac12 \sum_{(x,y) \in \edges} \cond_{xy}^{-1} \curr_1(x,y) \curr_2(x,y).
  \end{equation}
  \glossary{name={$\diss(\curr)$},description={dissipation of the current \curr},sort=d,format=textbf}
  Again, the normalizing factor of $\frac12$ corresponds to the idea that each edge only contributes once to the sum.
  The domain of the dissipation is
  \linenopax
  \begin{equation}\label{eqn:def:dissipation-domain}
    \dom \diss := \{\curr \suth \diss(\curr) < \iy \}.  
  \end{equation}
\end{defn}

\begin{remark}\label{rem:dom(diss)-is-Hilbert}
 \version{}{\marginpar{Is an orientation necessary?}}
 When an orientation $\sO$ for \Graph is chosen, it is easy to see that $\dom \diss$ is a Hilbert space under the inner product \eqref{eqn:def:dissipation-form}. Indeed, $\dom \diss = \ell^2(\sO,\cond)$.
\end{remark}

\begin{defn}\label{def:cycles}
  A \emph{cycle} \cycle is a set of $n$ edges corresponding to a sequence of vertices $(x_0, x_1, x_2, \dots, x_n=x_0) \ci \verts$, for which $(x_k,x_{k+1}) \ci \edges$ for each $k$. Denote the set of \emph{cycles} in \Graph by \Loops.
  \glossary{name={\cycle},description={cycle},sort=Q,format=textbf}
  \glossary{name={\Loops},description={set of all cycles},sort=L,format=textbf}
\end{defn}

\begin{defn}\label{def:Kirchhoff's-law}
  For physical realism, we often require that a current flow satisfy \emph{Kirchhoff's node law}, i.e., that the total current flowing into a vertex must equal the total current flowing out of a vertex:
  \linenopax
  \begin{equation}\label{eqn:Kirchhoff}
    \sum_{y \nbr x} \curr(x,y) = 0, \forall x \in \verts.
  \end{equation}
  This is indeed the version of Kirchhoff's law you would find in a physics textbook; with our convention $I(x,y)>0$ indicates that the current flows from $x$ to $y$.

  However, if we are performing the experiment described above, then there are boundary \condns at $\ga,\gw$ to take into account, and Kirchhoff's node law takes the nonhomogeneous form
  \linenopax
  \begin{equation}\label{eqn:Kirchhoff-nonhomog}
    \sum_{y \nbr x} \curr(x,y)
    = \gd_\ga - \gd_\gw
    = \begin{cases} 1, &x=\ga, \\ -1, &x=\gw, \\ 0, &\text{else}, \end{cases}
  \end{equation}
  where $\gd_x$ is the usual Dirac mass at $x \in \verts$.
\end{defn}

\begin{defn}\label{def:current-flow}
  A \emph{current flow} from \ga to \gw is a current $\curr \in \dom \diss$ that satisfies \eqref{eqn:Kirchhoff-nonhomog}. The set of all current flows is denoted $\Flo(\ga,\gw)$.
  \glossary{name={$\Flo(\ga,\gw)$},description={set of current flows from \ga to \gw},sort=f}
\end{defn}
We usually use \ga to denote the beginning of a flow and \gw to denote its end. Shortly, we will see that the currents corresponding to potentials are precisely the current flows.

\begin{remark}\label{rem:cycle-perturbation-preserves-Kirchhoff}
  Although trivial, it is important to note that the characteristic function of a current path $\charfn{\cpath}:\edges \to \{0,1\}$ trivially satisfies \eqref{eqn:Kirchhoff-nonhomog}. Also, the characteristic function of a cycle satisfies \eqref{eqn:Kirchhoff} in much the same way. As a consequence, if $\curr \in \Flo(\ga,\gw)$, then $\curr + t \charfn{\cycle} \in \Flo(\ga,\gw)$ for any $t \in \bR$ by a brief computation. In other words, perturbation on a cycle preserves the Kirchhoff condition. However, the dissipation will vary because $\diss(\charfn{\cycle}) > 0$.
  \glossary{name={$\charfn{B}$},description={characteristic function of the set $B$},sort=C}
  \glossary{name={$\charfn{B}$},description={characteristic function of the set $B$},sort=X}
\end{remark}

\section{Potential functions and their relationship to current flows.}
\label{sec:potential-functions}
From the proceeding section, it is clear that a special role is played by  functions $v : \verts \to \bR$ which satisfy the equation $\Lap v = \gd_\ga - \gd_\gw$. Such a function is the solution to a discrete Dirichlet problem, where the ``boundary'' has been chosen to be \ga and \gw (not to be confused with the boundary term discussed in Remark~\ref{rem:boundary-term}). 

\begin{defn}\label{def:dipole}
  A \emph{dipole} is a function $v \in \dom \energy$ which satisfies
  \linenopax
  \begin{equation}\label{eqn:def:dipole}
    \Lap v = \gd_\ga - \gd_\gw
  \end{equation}
  for some vertices $\ga, \gw \in \verts$. The collection of all such functions is denoted $\Pot(\ga,\gw)$. Note that when \Graph is finite, $\Pot(\ga,\gw)$ contains only a single element. This follows from Corollary~\ref{thm:harmonic=const-on-finite} because  the difference of any two solutions to \eqref{eqn:def:dipole} is harmonic.
\end{defn}
  \glossary{name={$\Pot(\ga,\gw)$},description={solutions of $\Lap v = \gd_\ga-\gd_\gw$},sort=p,format=textbf}


\begin{remark}\label{rem:def:monopole}
  The definition of a monopole that we give here is a heuristic definition; we give the precise definition in Definition~\ref{def:monopole}. A \emph{monopole} at \gw is a function $w:\verts \to \bR$ which \sats
  \linenopax
  \begin{equation}\label{eqn:monopole}
    \Lap w = k\gd_\gw,
    \qq w \in \dom \energy, k \in \bC.
  \end{equation}
  In the sequel, we are primarily concerned with monopoles $w_o$, where $o=\gw$ is some fixed vertex which acts as a point of reference or ``origin''. Also, we typically take $k=-1$, as the induced current of such a monopole is a \emph{unit flow to infinity} in the language of \cite{DoSn84}.
\end{remark}

\begin{remark}\label{rem:dipoles,monopoles,harmonics}
  The study of dipoles, monopoles, and harmonic functions is a recurring theme of this book:
  \linenopax
  \begin{align*}
    \Lap v = \gd_\ga-\gd_\gw,
    \qq \Lap w = -\gd_\gw,
    \qq \Lap h = 0.
  \end{align*}
  In Theorem~\ref{thm:Pot-is-nonempty-by-current-flows}, we will show that $\Pot(\ga,\gw)$, is nonempty for any \ga and \gw, on any network $(\Graph, \cond)$; the existence of monopoles and nontrivial harmonic functions is a much more subtle issue.

  In Corollary~\ref{thm:Pot-nonempty-via-Riesz}, we offer a more refined proof of the existence of dipoles, using Hilbert space techniques. Perhaps a more interesting question is when $\Pot(\ga,\gw)$ contains more than element; the linearity of \Lap shows immediately that any two dipoles in $\Pot(\ga,\gw)$ differ by a harmonic function. We have shown that when a connected graph is finite the only harmonic functions are constant (Corollary~\ref{thm:harmonic=const-on-finite}), and therefore $\Pot(\ga,\gw)$ consists only of a single function, up to the addition of a constant. The situation for monopoles is similar, as the difference of two monopoles at \gw is also a harmonic function.

  Not all \ERNs support monopoles; the current induced by a monopole is a finite flow to infinity and hence indicates that the random walk on the network is transient, by \cite{TerryLyons}. See also \cite{DoSn84,Lyons:ProbOnTrees} for terminology and proofs. It is well-known that for a reversible Markov chain, if the random walk started at one vertex is transient, then it is transient when started at any vertex. We give a very brief proof of this in Lemma~\ref{thm:transient-here-implies-transient-everywhere}; and a new criterion for transience in Lemma~\ref{thm:transient-iff-Fin=Ran(Lap)}.
  
  On some networks, a monopole can be understood as the limit of a sequence of dipoles $v_{x_n}$ where $\Lap v_{x_n} = \gd_{x_n} - \gd_o$ and $x_n \to \iy$. In such a situation, a monopole can be considered as a dipole where one of the Dirac masses ``sits at \iy''. However, this is not possible on all networks, as is illustrated by the binary tree in Example~\ref{exm:monopole-on-binary-tree}.
  Again, the linearity of \Lap shows immediately that any two monopoles at $\gw$ differ by a harmonic function. When these monopoles correspond to a ``distribution of dipoles at infinity'' (i.e., a limit of sums $\sum a_x v_x$ where the $v_x$'s are dipoles with $x \to \iy$ in the limit), the addition of a harmonic function transforms the distribution at infinity. It will take some work to make these ideas precise; for now the reader can consider this remark simply as a preview of coming attractions. The presence of monopoles is also extremely closely related to the existence of ``long-range order'', and the theoretical foundation of magnetism in \bRT; see \S\ref{sec:magnetism}. 

  Furthermore, it is possible for an \ERN to support monopoles but not nontrivial harmonic functions. In \S\ref{sec:lattice-networks}, we show that the integer lattice networks $(\bZd,\one)$ support monopoles (Theorem~\ref{thm:monopoles-on-Zd}). However, all harmonic functions are linear and hence do not have finite energy; cf.~Theorem~\ref{thm:harmonics-on-Zd-are-linear}. Both of these results are well-known in the literature in different contexts, and/or with different terminology.
\end{remark}

\begin{lemma}\label{thm:convexity-of-Pot,Flo}
  The dipoles $\Pot(\ga,\gw)$ and the current flows $\Flo(\ga,\gw)$ are convex sets. Furthermore, if $v \in \Pot(\ga,\gw)$, then $v+h \in \Pot(\ga,\gw)$ for any harmonic \fn $h$; similarly, if $\curr \in \Flo(\ga,\gw)$, then $\curr + J \in \Flo(\ga,\gw)$ for any \fn $J$ \satg \eqref{eqn:Kirchhoff}.
  \begin{proof}
    If $v_i \in \Pot(\ga,\gw)$, $c_i \geq 0$ and $\sum c_i = 1$, then the linearity of \Lap gives
    \linenopax
    \begin{align*}
      \Lap \left(\sum c_i v_i\right)
      = \sum c_i \Lap v_i
      = \sum c_i(\gd_\ga - \gd_\gw)
      = \gd_\ga - \gd_\gw.
    \end{align*}
    The computation for the other parts is similar.
  \end{proof}
\end{lemma}

\begin{theorem}\label{thm:energy-obtains-min}
  \energy obtains its minimum for some unique $v \in \Pot(\ga,\gw)$, and \diss obtains its minimum for some unique $\curr \in \Flo(\ga,\gw)$.
  \begin{proof}
    Each of these is a quadratic form on a convex set, by Lemma~\ref{thm:convexity-of-Pot,Flo}, so the result is an immediate application of \cite[Thm.~4.10]{Rud87} or \cite{Nel69}, e.g. To underscore the uniqueness, suppose that $\energy(v_1)=\energy(v_2)$. Then with $\ge := \inf\{\energy(v) \suth v \in \Pot(\ga,\gw)\}$, the parallelogram law gives
    \linenopax
    \begin{align*}
      \energy(v_1-v_2) &= 2\energy(v_1) + 2\energy(v_2) - 4\ge^2 = 0,
    \end{align*}
    since $\energy(v_i)=\ge$ because $v_i$ were chosen to be minimal.
  \end{proof}
\end{theorem}

\begin{defn}\label{def:ohm's-law}
  Ohm's Law ($V = RI$) appears in the present context as
  \linenopax
  \begin{equation}\label{eqn:ohm's-law}
    v(x) - v(y) = \frac1{\cond_{xy}} \curr(x,y).
  \end{equation}
\end{defn}

\begin{remark}\label{rem:Kirchhoff-vs-harmonic}
  It will shortly become evident (if it isn't already) that current flows satisfying Kirchhoff's law correspond to harmonic functions via Ohm's law 
  and that current flows satisfying the nonhomogeneous Kirchhoff's law \eqref{eqn:Kirchhoff-nonhomog} correspond to dipoles, that is, solutions of the Dirichlet problem \eqref{eqn:def:dipole} with Neumann boundary conditions. To make this precise, we need the notion of induced current given in Definition~\ref{def:induced-current} and justified by Lemma~\ref{thm:potential-induces-current}.
\end{remark}

\begin{lemma}
  \label{thm:potential-induces-current}
  \label{thm:energy=dissipation}
  Every function $v:\verts \to \bR$ induces a unique current via $I(x,y) := \cond_{xy}(v(x)-v(y))$, and the dissipation of this current is the energy of $v$:
  \linenopax
  \begin{equation}\label{eqn:energy=dissipation}
    \diss(\curr) = \energy(v).
  \end{equation}
  Moreover, if $v \in \Pot(\ga,\gw)$, then $\curr \in \Flo(\ga,\gw)$.
  \begin{proof}
    It is clear that Ohm's Law defines a current. The equality \eqref{eqn:energy=dissipation} is a very brief calculation and follows straight from the definitions; see \eqref{eqn:def:energy-form} and \eqref{eqn:def:dissipation}. A proof of \eqref{eqn:energy=dissipation} is also given in \cite{DoSn84}.

    If $v \in \Pot(\ga,\gw)$, then $\Lap v = \gd_\ga - \gd_\gw$ and
    \linenopax
    \begin{equation}\label{eqn:pot-gives-curr}
      (\gd_\ga - \gd_\gw)(x)
       = (\Lap v)(x)
       = \sum_{y \nbr x} \cond_{xy}(v(x)-v(y))
       = \sum_{y \nbr x} \curr(x,y)
    \end{equation}
    verifies the nonhomogeous form of Kirchhoff's law.
  \end{proof}
\end{lemma}

\begin{defn}\label{def:induced-current}
  Given $v \in \Pot(\ga,\gw)$, the \emph{induced current} is defined via Ohm's Law as in the statement of Lemma~\ref{thm:potential-induces-current}. That is,
  \linenopax
  \begin{equation}\label{eqn:induced-current}
    \curr(x,y) := \cond_{xy}(v(x) - v(y)).
  \end{equation}
\end{defn}

\begin{remark}\label{rem:energy=dissipation-only-for-potentials}
  Note that \eqref{eqn:energy=dissipation} holds when the current \curr is induced by $v$. It makes no sense to attempt to apply the same equality to a general current: there may be NO associated potential because of the compatibility problem described just below. Nonetheless, Theorem~\ref{thm:proj-is-minimizer} provides a way to give the identity analogous to \eqref{eqn:energy=dissipation} for general currents by using the adjoint of the operator implicit in \eqref{eqn:induced-current}.
\end{remark}

\begin{remark}\label{rem:variational-approach}
  If $\Lap v = \gd_\ga-\gd_\gw$ has a solution $v_0$, then any other solution is of the form $v=v_0+h$ where $h$ is harmonic, by linearity of \Lap. So to minimize energy, one must consider such perturbations:
  \linenopax
  \begin{align*}
    \frac{d}{dt} [\energy(v_0+th)]_{t=0} = 0
    \q\iff\q \energy(v_0,h)=0.
  \end{align*}
  Conversely, if $\energy(v,h) = 0$, then
  \linenopax
  \begin{align*}
    \energy(v+th) = \energy(v) + 2t \energy(v,h) + t^2 \energy(h) \geq \energy(v),
  \end{align*}
  shows that energy is minimized for $t = 0$. In particular, energy is minimized for $v$ which contains no harmonic component. In Lemma~\ref{thm:HE=Fin+Harm} this important principle is restated in the language of Hilbert spaces: energy is minimized for the $v$ which is orthogonal to the space of harmonic functions with respect to \energy.

  Analogous remarks hold for \curr which minimizes $\diss(\curr)$. However, note that Kirchhoff's Law is blind to conductances and so $\curr \in \Flo(\ga,\gw)$ does not imply that $\diss(\curr)$ is minimal. In the next section, we show that \emph{induced} currents are minimal with respect to \diss when they are induced by a minimal potential $v$.
\end{remark}

\section{The compatibility problem}
\label{sec:compatibility-problem}

The converse to Lemma~\ref{thm:potential-induces-current} is not always true, but a partial converse is given by Theorem~\ref{thm:minimal-flows-are-induced-by-minimizers}.
Given an electrical resistance network $(\Graph, \cond)$, one can always attempt to construct a Ohm's function by fixing the value $v(x_0)$ at some point $x_0 \in \verts$, and then applying Ohm's law to determine the value of $v$ for other vertices $x \nbr x_0$. However, this attempt can fail if the network contains a cycle (see Example~\ref{exm:one-small-cycle-network} for an example) because the existence of a cycle is equivalent to the existence of two distinct paths from one point to another. This phenomenon is worked out in detail for a simple case in Example~\ref{exm:edge-mass}.

In general, it may happen that there are two different paths from $x_0$ to $y_0$, and the net voltage drop $v(x_0) - v(y_0)$ computed along these two paths is not equal. Such a phenomenon makes it impossible to define $v$. Note that Kirchhoff's law does not forbid this, because \eqref{eqn:Kirchhoff}--\eqref{eqn:Kirchhoff-nonhomog} is expressed without reference to the conductances \cond. We refer to this as the \emph{compatibility problem}: \textbf{a general current function may not correspond to a potential}, even though every potential induces a well defined current flow (see Lemma~\ref{thm:potential-induces-current}). In this section we provide the following answer: for any current, there exists a unique associated current which \emph{does} correspond to a potential.

\begin{exm}[The Dirac mass on an edge]\label{exm:edge-mass}
  Consider a Dirac mass on an edge of the network $(\bZ^2,\one)$ as depicted in Figure~\ref{fig:edge-mass}. We use such a current here to illustrate the compatibility problem. 
  \begin{figure}
    \centering
    \scalebox{1.0}{\includegraphics{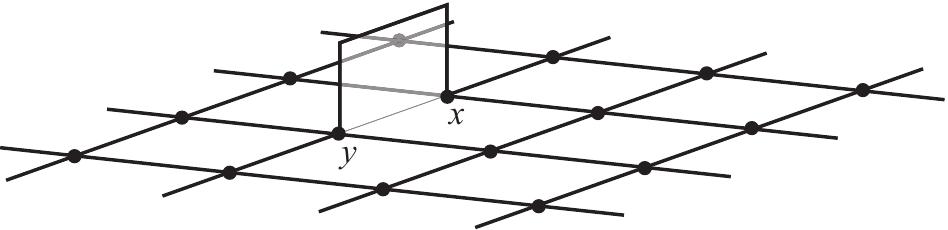}}
    \caption{\captionsize A Dirac mass on an edge of $\bZ^2$.}
    \label{fig:edge-mass}
  \end{figure}
  To find a potential corresponding to this current, consider the following dilemma: $\curr(x,y)=1$ and $\curr\equiv0$ elsewhere corresponds to a potential (up to a constant) which has $v(x)=1$ and $v(y)=0$, as in Figure~\ref{fig:edge-mass-potential}.
  \begin{figure}
    \centering
    \scalebox{1.0}{\includegraphics{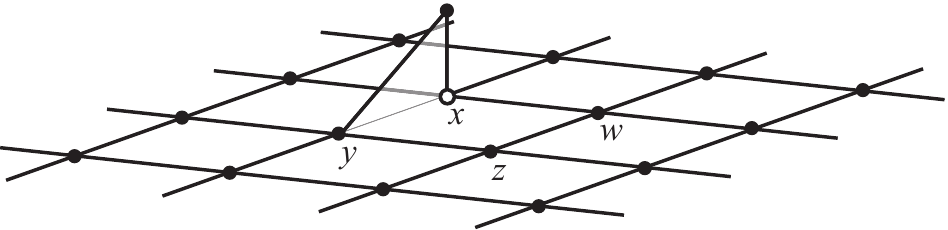}}
    \caption{\captionsize A failed attempt at constructing a potential to match Figure~\ref{fig:edge-mass}.}
    \label{fig:edge-mass-potential}
  \end{figure}
  Since $\curr(x,w)=0$, we have $v(w)=v(x)$, and since $\curr(y,z)=0$, we have $v(z)=v(y)$. But then $v(z)=1\neq0=v(z)$, contradicting the fact that $\curr(w,z)=0$! \cont
\end{exm}

\begin{defn}\label{def:cycle-condition}
  A current \curr satisfies the \emph{cycle condition} iff $\diss(\curr,\charfn{\cycle})=0$ for every cycle $\cycle \in \Loops$. (We call $\charfn{\cycle}$ a \emph{cycle}.)
\end{defn}

\begin{remark}\label{rem:noncyclic-cycles-give-path-independence}
  From the preceding discussion, it is clear that for a current satisfying the cycle condition, voltage drop between vertices $x$ and $y$ may be measured by summing the currents along any single path from $x$ to $y$, and the result will be independent of which path was chosen. In the Hilbert space interpretation of \S\ref{sec:H_energy-and-H_diss} the cycle condition is restated as ``\curr is orthogonal to cycles''. The next two results must be folklore (perhaps dating back to the 19$^\text{th}$ century?) but we include them for their relevance in \S\ref{sec:H_energy-and-H_diss}, especially the Hilbert space decomposition of Theorem~\ref{thm:HD=Nbd+Kir+Cyc} (see also Figure~\ref{fig:HE->HD->HE-decomposition}). While writing a second draft of this document, the authors discovered a similar treatment in \cite[\S9]{Lyons:ProbOnTrees}.
\end{remark}

\begin{lemma}\label{thm:induced-currents-satisfy-cycle-condition}
  \curr is an induced current if and only if \curr satisfies the cycle condition.
  \begin{proof}
    \fwd If \curr is induced by $v$, then for any $\cycle \in \Loops$, the sum
    \linenopax
    \begin{equation}\label{eqn:voltage-sum-on-cycle=0}
      \sum_{(x,y) \in \cycle} \frac1{\cond_{xy}} \curr(x,y)
      = \sum_{(x,y) \in \cycle} (v(x)-v(y))
      = 0,
    \end{equation}
    since every term $v(x_i)$ appears twice, once positive and once negative, whence $\diss(\curr,\charfn{\cycle})=0$.

    \bwd Conversely, to prove that there is such a $v$, we must show that $v(x_0)-v(y_0)$ is independence of the path from $x_0$ to $y_0$ used to compute it. In a direct analogy to basic vector calculus, this is equivalent to the fact that the net voltage drop around any closed cycle is 0.
    \linenopax
    \begin{equation*}
      \sum_{(x,y) \in \cycle} (v(x)-v(y))
      = \sum_{(x,y) \in \cycle} \frac1{\cond_{xy}} \curr(x,y)
      = \diss(\curr,\charfn{\cycle})
      = 0,
    \end{equation*}
    Now define $v$ by fixing $v(x_0)$ for some point $x_0 \in \verts$, and then coherently use $v(x) - v(y) = \frac1{\cond_{xy}}\curr(x,y)$ to compute $v$ at any other point.
  \end{proof}
\end{lemma}

The presence of cycles is not always obvious! As an exercise, we invite the reader to determine the cycles involved in Example~\ref{exm:edge-mass}.

\begin{lemma}[Resurrection of Kirchhoff's Law]
  \label{thm:Kirchhoff-iff-Pot}
  Let \curr be the current induced by $v$. Then $v \in \Pot(\ga,\gw)$ if and only if \curr satisfies the nonhomogeneous Kirchhoff's law.
  \begin{proof}
    \fwd Computing directly,
    \linenopax
    \begin{equation}\label{eqn:current-sum-at-a-point}
      \sum_{y \nbr x} \curr(x,y)
      = \sum_{y \nbr x} \cond_{xy}(v(x)-v(y))
      = \Lap v(x)
      = \gd_\ga - \gd_\gw.
    \end{equation}

    \bwd Conversely, to show $\Lap v = \gd_\ga - \gd_\gw$,
    \linenopax
    \begin{align*}
      \Lap v(x)
      &= \sum_{y \nbr x} \cond_{xy} (v(x)-v(y))
        &&\text{def of \Lap} \\
      &= \sum_{y \nbr x} \curr(x,y)
        &&\eqref{eqn:ohm's-law} \\
      &= \gd_\ga - \gd_\gw,
        &&\curr \in \Flo(\ga,\gw).
        \qedhere
    \end{align*}
  \end{proof}
\end{lemma}

\begin{cor}\label{thm:Kirchhoff-iff-Harm}
  Let \curr be the current induced by $v$. Then $v$ is harmonic if and only if \curr satisfies the homogeneous Kirchhoff's law.
  \begin{proof}
    Mutatis mutandis, this is the same as the proof of Lemma~\ref{thm:Kirchhoff-iff-Pot}.
  \end{proof}
\end{cor}

\begin{theorem}\label{thm:minimal-flows-are-induced-by-minimizers}
  \curr minimizes \diss on $\Flo(\ga,\gw)$ if and only if \curr is induced by the potential $v$ that minimizes \energy on $\Pot(\ga,\gw)$.
  \begin{proof}
    \fwd Since \curr minimizes \diss, we have
    \linenopax
    \begin{equation}\label{eqn:variation-of-diss=0-on-cycles}
      \frac{d}{dt} \left[\diss(\curr + t J)\right]_{t=0} = 0,
    \end{equation}
    for any current $J$ satisfying the homogeneous Kirchhoff's law. From Remark~\ref{rem:cycle-perturbation-preserves-Kirchhoff}, this applies in particular to $J=\charfn{\cycle}$, where \cycle is any cycle in \Loops.

    Note that $\diss(\curr,\charfn{\cycle})=0$ iff $\frac{d}{dt} \left[\diss(\curr + t \charfn{\cycle})\right]_{t=0} = 0$.
    To see this, replace \curr by $\curr + t \charfn{\cycle}$ in \eqref{eqn:def:dissipation}, differentiate $\diss(\curr + t \charfn{\cycle})$ term-by-term with respect to $t$ and evaluates at $t=0$ to obtain that \eqref{eqn:variation-of-diss=0-on-cycles} is equivalent to
    \linenopax
    \begin{equation*}
      \sum_{(x,y) \in \edges} \frac1{\cond_{xy}} \curr(x,y) \charfn{\cycle}(x,y)
      = \sum_{(x,y) \in \cycle} \frac1{\cond_{xy}} \curr(x,y)
      = 0, \qq \forall \cycle \in \Loops.
    \end{equation*}
    By Lemma~\ref{thm:induced-currents-satisfy-cycle-condition}, this shows that \curr is induced by some $v$; and by Lemma~\ref{thm:Kirchhoff-iff-Pot}, we know $v \in \Pot(\ga,\gw)$. From Lemma~\ref{thm:energy=dissipation}, it is clear that the $v$ must also be the energy-minimizing element of $\Pot(\ga,\gw)$.


    \bwd Since \curr is induced by $v \in \Pot(\ga,\gw)$, the only thing we need to check is that \curr is minimal with respect to any harmonic current (i.e. a current induced by a harmonic function); this follows from Lemma~\ref{thm:induced-currents-satisfy-cycle-condition} and the first part of the proof. If $h$ is any harmonic function on \verts, denote the induced current by $H$ as before. Then Lemma~\ref{thm:energy=dissipation} gives
    \linenopax
    \begin{align*}
      \frac{d}{dt} \left[\diss(\curr + t H)\right]_{t=0}
      = \frac{d}{dt} \left[\energy(v + t h)\right]_{t=0}
      = 0,
    \end{align*}
    by the minimality of $v$.
  \end{proof}
\end{theorem}

\begin{theorem}\label{thm:Pot-is-nonempty-by-current-flows}
  $\Pot(\ga,\gw)$ is never empty.
  \begin{proof}
    It is clear that $\Flo(\ga,\gw) \neq \es$ because one always has the characteristic \fn of a current path from \ga to \gw (since we are assuming the underlying graph is connected); see Definition~\ref{def:current-path} and Remark~\ref{rem:cycle-perturbation-preserves-Kirchhoff}. From Theorem~\ref{thm:energy-obtains-min} one sees that there is always a flow which minimizes dissipation. By Theorem~\ref{thm:minimal-flows-are-induced-by-minimizers}, this minimal flow is induced by an element of $\Pot(\ga,\gw)$.
  \end{proof}
\end{theorem}

The following result is well-known in probability (see, e.g., \cite{Stroock}), but we include it here for completeness and the novel method of proof.

\begin{cor}\label{thm:transient-here-implies-transient-everywhere}
  If the random walk on $(\Graph,\cond)$ is transient when started from $y \in \verts$, then it is transient when started from any $x \in \verts$. 
  \begin{proof}
    By \cite{TerryLyons}, the hypothesis means there is a monopole $w_y \in \dom\energy$ with $\Lap w = \gd_y$. But then by Theorem~\ref{thm:Pot-is-nonempty-by-current-flows} and the linearity of \Lap, $v + w_x$ is a monopole at $x$, for any $v \in \Pot(x,y)$.
  \end{proof}
\end{cor}

\begin{remark}\label{rem:Pot-is-not-contained-in-l2}
  There are examples for which the elements of $\Pot(\ga,\gw)$ do not lie in $\ell^2(\verts)$; see Figure~\ref{fig:vx-in-Z1} of Example~\ref{exm:infinite-lattices}.
\end{remark}

\begin{prop}\label{thm:v-min-at-a,max-at-w}
  If \Graph is finite and $v \in \Pot(\ga,\gw)$, then $v(\gw) \leq v(x) \leq v(\ga)$ for all $x \in \verts$.
  \begin{proof}
    This is immediate from the maximum principle for harmonic functions on the finite set \verts with boundary $\{\ga,\gw\}$. See \cite[\S2.1]{Lyons:ProbOnTrees}, for example, or \cite{LevPerWil08}.
  \end{proof}
\end{prop}

\begin{remark}\label{rem:v-bounded}
  In \S\ref{sec:effective-resistance-metric}, we will see that Proposition~\ref{thm:v-min-at-a,max-at-w} extends to a more general result: if $v$ is the unique element of $\Pot(\ga,\gw)$ of minimal energy, then the same conclusion follows. One way to see this is to define $u(x) = \prob_x[\gt_\ga < \gt_\gw]$ (i.e., the probability that the random walk started at $x$ reaches \ga before \gw). By Theorem~\ref{thm:wired-resistance}, $v$ is defined by $v(x) = u(x) R^W(\ga,\gw)$, where $R^W(\ga,\gw)$ is the (wired) effective resistance between \ga and \gw.
\end{remark}

\subsection{Current paths}
\label{sec:current-paths}

It is intuitively obvious that for a connected graph, current should be able to flow between any two points. Indeed, it is a basic result in graph theory that for any connected graph, one can find a path of minimal length between any two points and from Remark~\ref{rem:def:connected-ERN} we know that the resistance along such a path is finite. In Theorem~\ref{thm:Paths-is-nonempty}, we will show a stronger result: that one can always find a path along which the potential function decreases monotonically. In other words, there is always at least one ``downstream path'' between the two vertices. Somewhat surprisingly, this fact is easiest to demonstrate by an appeal to a basic fact about probability (Lemma~\ref{thm:flow-is-sum-of-its-parts}),

Our definition of an \ERN is mathematical (Definition~\ref{def:ERN}) but is motivated by engineering; modification of the conductors (\cond) will alter the associated probabilities and thus change which current flows are induced, in the sense of Definition~\ref{def:induced-current}. We are interested in quantifying this dependence. Obviously, on an infinite graph the computation of current paths involves all of \Graph, and it is not feasible to attempt to compute these paths directly.
Consequently, we feel our proof of Theorem~\ref{thm:Paths-is-nonempty} may be of independent interest.

\begin{defn}\label{def:current-path-preview}
  Let $v:\verts \to \bR$ be given and suppose we fix \ga and \gw for which $v(\ga) > v(\gw)$. Then a \emph{current path} \cpath (or simply, a \emph{path}) is an edge path from \ga to \gw with the extra stipulation that $v(x_k) < v(x_{k-1})$ for each $k=1,2,\dots,n$.
  Denote the set of all current paths by $\Paths = \Paths_{\ga,\gw}$ (dependence on the initial and terminal vertices is suppressed when context precludes confusion). Also, define $\Paths_{\ga,\gw}(x,y)$ to be the subset of current paths from \ga to \gw which pass through the edge $(x,y) \in \edges$.
  \glossary{name={\Paths},description={set of paths},sort=G,format=textbf}
\end{defn}

The following lemma is immediate from elementary probability theory, as it
represents the probability of a union of disjoint events, but it will be helpful.

\begin{lemma}\label{thm:flow-is-sum-of-its-parts}
  Suppose $(\Graph, \cond)$ is an \ERN and $v:\verts \to \bR$ satisfies
  $\Lap v = \gd_\ga - \gd_\gw$. Then if $\curr$ is the current associated to $v$ by $\curr(x,y) = \cond_{xy}(v(x)-v(y)$, then \curr satisfies 
  \begin{equation}\label{eqn:flow-is-sum-of-its-parts}
    \curr(x,y) = \sum_{\cpath \in \Paths_{\ga,\gw}(x,y)} \prob(\cpath).
  \end{equation}
\end{lemma}

The method of proof in the next proposition is a bit unusual in that it uses a probability to demonstrate existence. This result fills a hole in the proof of $\dist_\Lap(x,y) = \dist_\diss(x,y)$ in \cite{Pow76b} (recall \eqref{eqn:def:R(x,y)-Lap} and \eqref{eqn:def:R(x,y)-diss}).

\begin{theorem}\label{thm:Paths-is-nonempty}
  If $v \in \Pot(\ga,\gw)$, then $\Paths_{\ga,\gw} \neq \es$. Moreover, $v(\ga) > v(\gw)$.
  \begin{proof}
    Theorem~\ref{thm:Pot-is-nonempty-by-current-flows} ensures we can find  $v \in \Pot(\ga,\gw)$; let \curr be the current flow associated to $v$. Then $\Lap v(\ga) = 1$ implies that there is some $y \nbr \ga$ for which $\curr(\ga,y) >0$. By Lemma~\ref{thm:flow-is-sum-of-its-parts},
    \linenopax
    \begin{align*}
      \curr(\ga,y) = \sum_{\cpath \in \Paths_{\ga,\gw}(\ga,y)} \prob(\cpath) > 0,
    \end{align*}
    which implies there must exist a positive term in the sum, and hence a $\cpath \in \Paths_{\ga,\gw}$. Since we may now choose a path $\cpath \in \Paths_{\ga,\gw}$, the second claim follows. 
  \end{proof}
\end{theorem}

\section{Remarks and references}
\label{sec:Remarks-and-References-2}

Most of this chapter is just based on high-school physics, but a couple of key references for us were \cite{DoSn84}, \cite{Lyons:ProbOnTrees} and \cite{Soardi94} and the papers by Bob Powers \cite{Pow75, Pow76a, Pow76b, Pow78, Pow79}, and the early (and often overlooked) paper \cite{BoDu49} by Bott and Duffin. In addition, one may find \cite{Woess00, LPV08} 
to be useful. We especially recommend the discussion of networks and resistance distance in PowersÕ paper \cite{Pow76b}. While the main focus there is a problem for quantum spin systems on lattices, Powers develops the elementary properties of effective resistance from scratch in this paper, and adapts them to the Heisenberg model.

\begin{remark}\label{rem:preview-of-potentials-induce-currents-orthogonal-to-cycles}
  Part of the motivation for Theorem~\ref{thm:minimal-flows-are-induced-by-minimizers}
  is to fix an error in \cite{Pow76b}. The author was not apparently aware of the possibility of nontrivial harmonic functions, and hence did not see the need for taking the element of $\Pot(\ga,\gw)$ with minimal energy. This becomes especially important in Theorem~\ref{thm:effective-resistance-metric}.

  Theorem~\ref{thm:minimal-flows-are-induced-by-minimizers}
  is generalized in Theorem~\ref{thm:proj-is-minimizer} where we exploit certain operators to obtain, for any given current \curr, an associated minimal current. This minimal current is induced by a potential, even if the original is not, and provides a resolution to the compatibility problem described at the beginning of \S\ref{sec:compatibility-problem}, just above.

  In \S\ref{sec:Resolution-of-the-compatibility-problem} we revisit this scenario and show how the minimal current may be obtained by the simple application of a certain operator, once it has been properly interpreted in terms of Hilbert space theory. See Theorem~\ref{thm:proj-is-minimizer} and its corollaries in particular.
\end{remark}

\begin{remark}\label{rem:can-prove-POT-is-nonempty-via-Hilbert}
  Theorem~\ref{thm:Pot-is-nonempty-by-current-flows} fills a gap in \cite{Pow76b}. A key point is that the finite dissipation of the flow ensures the finite energy of the inducing voltage function, by Lemma~\ref{thm:energy=dissipation}. A different proof of Theorem~\ref{thm:Pot-is-nonempty-by-current-flows} is obtained in Corollary~\ref{thm:Pot-nonempty-via-Riesz} by the application of Hilbert space techniques. 
  
  Theorem~\ref{thm:Pot-is-nonempty-by-current-flows} also follows from results of \cite[\S{III.4}]{Soardi94} since the difference of two Dirac masses corresponds to a ``balanced'' flow, i.e., the same amount of current flows in as flows out.
\end{remark}

\begin{remark}\label{rem:Yamasaki}
  In Corollary~\ref{thm:transient-here-implies-transient-everywhere} we refer to \cite{TerryLyons} for the equivalence of transience with the existence of a finite flow to infinity. This is the most common reference for this result, but we should point out that a slightly different version of it (restated in terms of the Green function) appears earlier in the often overlooked paper \cite{Yamasaki79}. (These results were discovered independently.)
\end{remark}

%% file: energy-hilbert-space.tex

\chapter{The energy Hilbert space}
\label{sec:energy-Hilbert-space}

\headerquote{I would like to make a confession which may seem immoral: I do not believe in Hilbert space anymore.}{---~J.~von~Neumann}

\headerquote{I am acutely aware of the fact that the marriage between mathematics and physics, which was so enormously fruitful in past centuries, has recently ended in divorce.}{---~F.~Dyson}

A number of tools used in the theory of weighted graphs, especially for the infinite case, were first envisioned in the context of resistance networks. We will introduce them below, but it is helpful to keep in mind that they apply to a variety of other problems outside the original context of resistance networks. In particular, these concepts and tools lead to the introduction of metrics. These in turn have applications to neighboring fields; see for example Chapters~\ref{sec:effective-resistance-metric}, \ref{sec:lattice-networks}, and \ref{sec:Magnetism-and-long-range-order}.

For the analysis of resistance networks, the (Dirichlet) energy form \energy is a natural tool, and so it is helpful study the Hilbert space \HE of functions on the network where the inner product is given by \energy. While one's first instinct may be to select $\ell^2(X)$ as the preferred Hilbert space, we show in Chapter~\ref{sec:Construction-of-HE} that the energy space \HE is in some sense the natural choice. 
The relationship between the two Hilbert spaces $\ell^2(X)$ and \HE is subtle, and is explored further in Chapters~\ref{sec:Lap-on-HE}, \ref{sec:L2-theory-of-Lap-and-Trans}, and elsewhere. For example, the Laplacian operator \Lap has quite different properties depending on the choice of Hilbert space.

In this section, we study the Hilbert space \HE of (finite-energy) voltage functions, that is, equivalence classes of functions $u:\verts \to \bC$ where $u \simeq v$ iff $u-v$ is constant. On this space, the energy form is an inner product, and there is a natural reproducing kernel $\{v_x\}_{x \in \verts}$ indexed by the vertices; see Corollary~\ref{thm:vx-is-a-reproducing-kernel}. Since we work with respect to the equivalence relation defined just above, most formulas are given with respect to differences of function values; in particular, the reproducing kernel is given in terms of differences with respect to some chosen ``origin''. Therefore, for any given resistance network, we fix a reference vertex $o \in \verts$ to act as an origin. It will readily be seen that all results are independent of this choice, and this affords the convenience of working with a single-parameter reproducing kernel. When working with representatives, we typically abuse notation and use $u$ to denote the equivalence class of $u$. One natural choice is to take $u$ so that $u(o)=0$; a different but no less useful choice is to pick $k$ so that $v=0$ outside a finite set as discussed further in Definition~\ref{def:Fin}.




In Theorem~\ref{thm:E(u,v)=<u,Lapv>+sum(normals)}, we establish a discrete version of the Gauss-Green Formula which extends Lemma~\ref{thm:E(u,v)=<u,Lapv>} to the case of infinite graphs. The appearance of a somewhat mysterious boundary term 
prompts several questions which are discussed in Remark~\ref{rem:boundary-term}. Answering these questions comprises a large part of the sequel; cf.~\S\ref{sec:the-boundary}. We are able to prove in Lemma~\ref{thm:E(u,v)=<u,Lapv>-on-Fin} that this boundary term vanishes for finitely supported functions on \Graph, and in Corollary~\ref{thm:nontrivial-harmonic-fn-is-not-in-L2} that nontrivial harmonic functions cannot be in $\ell^2(\verts)$. Later, we will see that the boundary term vanishes precisely when the random walk on the network is recurrent. This is discussed further in Remark~\ref{rem:L2-loses-the-best-part} and provides the motivation for energy-centric approach we pursue throughout our study.

The energy Hilbert space \HE will facilitate our study of the resistance metric $R$ in \S\ref{sec:effective-resistance-metric}. In particular, it provides an explanation for an issue stemming from the ``nonuniqueness of currents'' in certain infinite networks; see \cite{Lyons:ProbOnTrees, Thomassen90}. This disparity leads to differences between two apparently natural extensions of the effective resistance to infinite networks, which are greatly clarified by the geometry of Hilbert space. Also, \HE presents an analytic formulation of the type problem for random walks on an \ERN: transience of the random walk is equivalent to the existence of monopoles, that is, finite-energy solutions to a certain Dirichlet problem. In fact, this approach will readily allow us to obtain explicit formulas for effective resistance on integer lattice networks in \S\ref{sec:lattice-networks}, with applications to a physics problem of \cite{Pow76b} in \S\ref{sec:Magnetism-and-long-range-order}. 

Most results in this section appeared in \cite{DGG}.

Let \one denote the constant function with value 1 and recall that $\ker \energy = \bC \one$.


\begin{defn}\label{def:H_energy}\label{def:The-energy-Hilbert-space}
  The energy form \energy is symmetric and positive definite, and its kernel is the set of constant functions on \Graph. Let \one denote the constant function with value 1. Then the \emph{energy space} $\HE := \dom \energy / \bC \one$ is a Hilbert space with inner product and corresponding norm given by
  \linenopax
  \begin{equation}\label{eqn:energy-inner-product}
    \la u, v \ra_\energy := \energy(u,v)
    \q\text{and}\q
    \|u\|_\energy := \energy(u,u)^{1/2}.
  \end{equation}
  \glossary{name={\HE},description={energy Hilbert space},sort=h,format=textbf}
  \glossary{name={$\la\cdot,\cdot\ra_\energy$},description={energy inner product},sort=<,format=textbf}
  \glossary{name={$o$},description={origin; an fixed vertex chosen arbitrarily to act as a reference point},sort=o,format=textbf}
\end{defn}
It can be checked directly that the above completion consists of (equivalence classes of) functions on \verts via an isometric embedding into a larger Hilbert space as in \cite{Lyons:ProbOnTrees,MuYaYo} or by a standard Fatou's lemma argument as in \cite{Soardi94}.
Fix a reference vertex $o \in \verts$ to act as an ``origin''. It will readily be seen that all results are independent of this choice.

\begin{remark}\label{rem:alternative-construction-of-HE}
  In \S\ref{sec:Construction-of-HE}, we provide an alternative construction of \HE via techniques of von Neumann and Schoenberg. This provides for a more explicit description of the structure of \HE and its relation to the metric geometry of $(\Graph,R)$, and shows that \HE is the natural Hilbert space in which to embed $(\Graph,R)$. However, this must be postponed until after the introduction of the effective resistance metric.
\end{remark}

\begin{remark}[Four warnings about \HE] \hfill 
  \label{rem:elements-of-HE-are-technically-equivalence-classes}
  \begin{enumerate}[(1)]
    \item \HE has no canonical o.n.b.; the usual candidates $\{\gd_x\}$ are not orthogonal and typically their span is not even dense, as we discuss further below. 
    \item \label{itm:not-Hermitian} Multiplication operators are not Hermitian; see Lemma~\ref{thm:multiplication-not-hermitian} and Remark~\ref{rem:HE-into-L2(S',prob)-gives-hermitian-multiplication}.
    \item There is no natural interpretation of \HE as an $\ell^2$-space of functions on the vertices \verts or edges \edges of $(\Graph,\cond)$. \HE does contain the embedded image of $\ell^2(\verts,\gm)$ for a certain measures \gm, but these spaces are not typically dense in \HE. Also, \HE embeds isometrically into a subspace of $\ell^2(\edges,\cond)$, but it is generally nontrivial to determine whether a given element of $\ell^2(\edges,\cond)$ lies in this subspace. \HE may also be understood as a $\ell^2$ space of random variables (see \S\ref{sec:Kolmogorov-construction-of-L2(gW,prob)}) or realized as a subspace of $L^2(S',\prob)$, where $S'$ is a certain space of distributions (see \S{sec:the-boundary}).
    \item Pointwise identities should not be confused with Hilbert space identities; see Remark~\ref{rem:Lap-defined-via-energy} and Lemma~\ref{thm:pointwise-identity-implies-adjoint-identity}.
  \end{enumerate}
  To elaborate on the last point, note that elements of \HE are technically equivalence classes of functions which differ only by a constant; this is what is meant by the notation $\dom \energy / \bR \one$. In other words, if $v_1 = v_2 + k$ for $k \in \bC$, then $v_1 = v_2$ in \HE. When working with representatives, we typically abuse notation and use $u$ to denote the equivalence class of $u$. Often, we choose $u$ so that $u(o)=0$ (occasionally without warning). A different but no less useful choice is to pick $k$ so that $v=0$ outside a finite set when $v$ is a function of finite support (see Definition~\ref{def:Fin}).
\end{remark}

As mentioned in \eqref{itm:not-Hermitian}  above, the following feature of \HE operator theory contrasts sharply with the more familiar Hilbert spaces of $L^2$ functions, where all \bR-valued functions define Hermitian multiplication operators.

\begin{lemma}\label{thm:multiplication-not-hermitian}
  If $\gf:\verts \to \bR$ and $M_\gf$ denotes the multiplication operator $M_\gf: u \mapsto \gf u$, then $M_\gf$ is Hermitian if and only if $\gf = 0$ in \HE. 
  \begin{proof}
    Choose any representatives for $u, v \in \HE$. From the formula \eqref{eqn:def:energy-form},
    \linenopax
    \begin{align*}
      \la M_\gf u, v\ra_\energy
      &= \frac12 \sum_{x,y \in \verts} \cond_{xy} (\gf(x)u(x)v(x) - \gf(x)u(x)v(y) - \gf(y)u(y)v(x) + \gf(y)u(y)v(y)).
    \end{align*}
    By comparison with the corresponding expression, this is equal to $\la u, M_\gf v\ra_\energy$ iff $(\gf(y) - \gf(x))u(y)v(x) = (\gf(y) - \gf(x))u(x)v(y)$. However, since we are free to vary $u$ and $v$, it must be the case that \gf is constant.
  \end{proof}
\end{lemma}

\begin{defn}\label{def:exhaustion-of-G}
  An \emph{exhaustion} of \Graph is an increasing sequence of finite and connected subgraphs $\{\Graph_k\}$, so that $\Graph_k \ci \Graph_{k+1}$ and $\Graph = \bigcup \Graph_k$.
  \glossary{name={$\{\Graph_k\}$},description={exhaustion of a network},sort=G,format=textbf}
\end{defn}

\begin{defn}\label{def:infinite-vertex-sum}  
  The notation
  \linenopax
  \begin{equation}\label{eqn:def:infinite-sum}
    \sum_{x \in \verts} := \lim_{k \to \iy} \sum_{x \in \Graph_k}
  \end{equation}
  is used whenever the limit is independent of the choice of exhaustion $\{\Graph_k\}$ of \Graph. We typically justify this independence by proving the sum to be absolutely convergent.
\end{defn}

\begin{remark}\label{rem:boundary-term}
  One of the main results in this section is a discrete version of the Gauss-Green theorem presented in Theorem~\ref{thm:E(u,v)=<u,Lapv>+sum(normals)}: 
  \linenopax
  \begin{equation}\label{eqn:E(u,v)=<u_0,Lapv>+sum(normals)-warmup}
    \la u, v \ra_\energy
    = \sum_{x \in \verts} \cj{u}(x) \Lap v(x)
      + \sum_{x \in \bd \Graph} \cj{u}(x) \dn v(x),
      \qq u, v \in \HE.
  \end{equation}
  This differs from the literature, where it is common to find $\energy(u,v) = \la u, \Lap v\ra_{\ell^2}$ given as a definition (of \energy or of \Lap, depending on the context), e.g. \cite{Kig01, Kig03, Str08}. 
  
  We refer to $\sum_{\bd \Graph} u \dn v$ as the ``boundary term'' by analogy with classical PDE theory. This terminology should not be confused with the notion of boundary that arises in the discussion of the discrete Dirichlet problem. In particular, the boundary discussed in \cite{Kig03} and \cite{Kig08} refers to a subset of \verts. By contrast, when discussing an \emph{infinite} network \Graph, our boundary $\bd \Graph$ is never contained in \Graph. Green's identity follows immediately from \eqref{eqn:E(u,v)=<u_0,Lapv>+sum(normals)-warmup} in the form
  \linenopax
  \begin{align}\label{eqn:Greens-identity}
    \sum_{x \in \verts} \left(\cj{u}(x) \Lap v(x) - \cj{v}(x) \Lap u(x)\right)
    = \sum_{x \in \bd \Graph} \left(\cj{v}(x) \dn u(x) - \cj{u}(x) \dn v(x)\right).
  \end{align}
  Note that our definition of the Laplace operator is the negative of that often found in the PDE literature, where one will find Green's identity written
  \linenopax
  \begin{align*}
    \int_\gW (u \Lap v - v\Lap u) = \int_{\del \gW} (u \dn{} v - v \dn{} u).
  \end{align*}

  As the boundary term may be difficult to contend with, it is extremely useful to know when it vanishes. We have several results concerning this:
  \begin{enumerate}[(i)]
    \item Lemma~\ref{thm:E(u,v)=<u,Lapv>-on-Fin} shows the boundary term vanishes when either argument of $\la u, v \ra_\energy$ has finite support,
   \item Lemma~\ref{thm:TFAE:Fin,Harm,Bdy} gives necessary and sufficient conditions on the \ERN for the boundary term to vanish for any $u,v \in \HE$,
     \item Lemma~\ref{thm:converse-to-E(u,v)=<u,Lapv>} show the boundary term vanishes when both arguments of $\la u, v \ra_\energy$ and their Laplacians lie in $\ell^2$.
  \end{enumerate}
  In fact, Lemma~\ref{thm:TFAE:Fin,Harm,Bdy} expresses the fact that it is precisely the presence of monopoles that prevents the boundary term from vanishing. 
  An example with nonvanishing boundary term is given in Example~\ref{exm:nonvanishing normal derivatives}.
\end{remark}

\section{The evaluation operator $L_x$ and the reproducing kernel $v_x$}
\label{sec:L_x-and-v_x}

\begin{defn}\label{def:L_x}
  For any vertex $x \in \verts$, define the linear evaluation operator $L_x$ on \HE by
  \linenopax
  \begin{equation}\label{eqn:def:L_x}
    L_x u := u(x) - u(o).
  \end{equation}
  \glossary{name={$L_x$},description={evaluation functional; see $v_x$},sort=L,format=textbf}
\end{defn}

\begin{lemma}\label{thm:L_x-is-bounded}
  For any $x \in \verts$, one has $|L_x u | \leq k \energy(u)^{1/2}$, where $k$ depends only on $x$.
  \begin{proof}
    Since \Graph is connected, choose a path $\{x_i\}_{i=0}^n$ with $x_0=o$, $x_n=x$ and $\cond_{x_i,x_{i-1}}>0$ for $i=1,\dots,n$. For $k = \left(\sum_{i=1}^n \cond_{x_i,x_{i-1}}^{-1} \right)^{1/2}$, the Schwarz inequality yields
    \linenopax
    \begin{align*}
      |L_x u |^2
      = |u(x)-u(o)|^2
      &= \left|\sum_{i=1}^n \sqrt{\frac{\cond_{x_i,x_{i-1}}}{\cond_{x_i,x_{i-1}}}} (u(x_i)-u(x_{i-1}))\right|^2 
      \leq k^2 \energy(u).
      \qedhere
    \end{align*}
  \end{proof}
\end{lemma}

\begin{defn}\label{def:vx}\label{def:energy-kernel}
  Let $v_x$ be defined to be the unique element of \HE for which
  \linenopax
  \begin{equation}\label{eqn:def:vx}
    \la v_x, u\ra_\energy = u(x)-u(o),
    \qq \text{for every } u \in \HE.
  \end{equation}
  This is justified by Lemma~\ref{thm:L_x-is-bounded} and the Riesz Representation Theorem.
  The family of functions $\{v_x\}_{x \in \verts}$ is called the \emph{energy kernel} because of Corollary~\ref{thm:vx-is-a-reproducing-kernel}.
  Note that $v_o$ corresponds to a constant function, since $\la v_o, u\ra_\energy = 0$ for every $u \in \HE$. Therefore, this term may be ignored or omitted.
  \glossary{name={$v_x$},description={energy kernel},sort=v,format=textbf}
\end{defn}

\begin{defn}\label{def:reproducing-kernel}
  Let $\sH$ be a Hilbert space of functions on $X$. An operator $S$ on \sH is said to have a \emph{reproducing kernel} $\{k_x\}_{x \in X} \ci \sH$ iff
  \linenopax
  \begin{equation}\label{eqn:def:reproducing-kernel}
    (S v)(x) = \la k_x, v\ra_\sH,
    \qq \forall x \in X, \forall v \in \sH.
  \end{equation}
  If $S$ is projection to a subspace $L \ci \sH$, then one says $\{k_x\}$ is a \emph{reproducing kernel for $L$}.
  If $S=\id$, then \sH is a \emph{reproducing kernel Hilbert space} with kernel $k$.
\end{defn}

\begin{theorem}[Aronszajn's Theorem {\cite{Aronszajn50}}]
  \label{thm:Aronszajn's-thm}
  Let $\{f_x\}$ be a reproducing kernel for \sH. Define a sesquilinear form on the set of all finite linear combinations of these elements by
  \linenopax
  \begin{equation}\label{eqn:Aronszajn's-form}
    \left\la \sum_x \gx_x f_x, \sum_y \gh_y f_y \right\ra
    := \sum_x \cj{\gx_x} \gh_x f_x(y).
  \end{equation}
  Then the completion of this set under the form \eqref{eqn:Aronszajn's-form} is again \sH.
\end{theorem}

\begin{cor}\label{thm:vx-is-a-reproducing-kernel} \label{thm:vx-dense-in-HE}
  $\{v_x\}_{x \in \verts}$ is a reproducing kernel for \HE. Thus, $\spn\{v_x\}$ is dense in \HE.
  \begin{proof}
    Choosing representatives with $v_x(o)=0$, it is trivial to check that $\la v_x, v_y\ra_\energy = v_x(y) = \cj{v_y(x)}$ and then apply Aronszajn's Theorem. 
  \end{proof}
\end{cor}

There is a rich literature dealing with reproducing kernels and their manifold application to both continuous analysis problems (see e.g., \cite{AD06, AL08, AAL08, BV03, Zh09}), and infinite discrete stochastic models.
One of the differences between these studies and our present work is the approach we take in Definition 4.9, i.e., the use of ``relative'' reproducing kernels.


\begin{remark}
  \label{rem:reproducing-kernel}
  Definition~\ref{def:energy-kernel} is justified by Corollary~\ref{thm:vx-is-a-reproducing-kernel}. In this book, the functions $v_x$ will play a role analogous to fundamental solutions in PDE theory; see \S\ref{sec:analogy-with-complex-functions}. 

  The functions $v_x$ are \bR-valued. This can be seen by first constructing the energy kernel for the Hilbert space of \bR-valued functions on \Graph, and then using the decomposition of a \bC-valued function $u = u_1 + \ii u_2$ into its real and imaginary parts. Alternatively, see Lemma~\ref{thm:vx-is-R-valued}.
  
  Reproducing kernels will help with many calculations and explain several of the relationships that appear in the study of \ERNs. They also extend the analogy with complex function theory discussed in \S\ref{sec:analogy-with-complex-functions}. The reader may find the references \cite{Aronszajn50,Yoo05,Jor83} to provide helpful background on reproducing kernels.
\end{remark}

\begin{remark}[Probabilistic interpretation of $v_x$]
  \label{rem:interpretation-of-v_x}
  The energy kernel $\{v_x\}$ is intimately related to effective resistance distance $R(x,y)$. In fact, $R(x,o) = v_x(x) - v_x(o) = \energy(v_x)$ and similarly, $R(x,y) = \energy(v_x-v_y)$. This is discussed in detail in \S\ref{sec:effective-resistance-metric}, but we give a brief summary here, to help the reader get a feeling for $v_x$. For a random walk (RW) started at the vertex $y$, let $\gt_x$ be the hitting time of $x$ (i.e., the time at which the random walk first reaches $x$) and define the function
  \linenopax
  \begin{align*}
    u_x(y) = \prob[\gt_x < \gt_o | \text{ RW starts at $y$}].
  \end{align*}
  Here, the RW is governed by transition probabilities  $p(x,y) = \cond_{xy}/\cond(x)$; cf.~Remark~\ref{rem:transient-iff-monopoles}. One can show that $v_x = R(x,o) u_x$ is the representative of $v_x$ with $v_x(o)=0$. Since the range of $u_x$ is $[0,1]$, one has $0 \leq v_x(y) - v_x(o) \leq v_x(x) - v_x(o) = R(x,o)$. Many other properties of $v_x$ are similarly clear from this interpretation. For example, it is easy to compute $v_x$ completely on any tree.
\end{remark}

\section{The finitely supported functions and the harmonic functions}
\label{sec:The-role-of-Fin-in-HE}

\begin{defn}\label{def:Fin}
  For $v \in \HE$, one says that $v$ has \emph{finite support} iff there is a finite set $F \ci \verts$ for which $v(x) = k \in \bC$ for all $x \notin F$. Equivalently, the set of functions of finite support in \HE is 
  \linenopax
  \begin{equation}\label{eqn:span(dx)}
    \spn\{\gd_x\} = \{u \in \dom \energy \suth u(x)=k \text{ for some $k$, for all but finitely many } x \in \verts\},
  \end{equation}
  where $\gd_x$ is the Dirac mass at $x$, i.e., the element of \HE containing the characteristic function of the singleton $\{x\}$. It is immediate from \eqref{eqn:energy-of-Diracs} that $\gd_x \in \HE$.
    \glossary{name={$\gd_x$},description={Dirac mass at the vertex $x$},sort=D}
  Define \Fin to be the closure of $\spn\{\gd_x\}$ with respect to \energy. 
\end{defn}

\begin{defn}\label{def:Harm}
  The set of harmonic functions of finite energy is denoted
  \linenopax
  \begin{equation}\label{eqn:Harm}
    \Harm := \{v \in \HE \suth \Lap v(x) = 0, \text{ for all } x \in \verts\}.
  \end{equation}
  Note that this is independent of choice of representative for $v$ in virtue of \eqref{eqn:def:laplacian}.
\end{defn}
  \glossary{name={\Harm},description={the harmonic functions of finite energy},sort=H,format=textbf}

\begin{lemma}\label{thm:<delta_x,v>=Lapv(x)}
  The Dirac masses $\{\gd_x\}_{x \in \verts}$ form a reproducing kernel for \Lap. That is, for any $x \in \verts$, one has $\la \gd_x, u \ra_\energy = \Lap u(x)$.
  \begin{proof}
    Compute $\la \gd_x, u \ra_\energy = \energy(\gd_x, u)$ directly from formula \eqref{eqn:def:energy-form}.
  \end{proof}
\end{lemma}

\begin{remark}\label{rem:Lap-defined-via-energy}
   Note that one can take the \emph{definition of the Laplacian} to be the operator $A$ defined via the equation 
   \linenopax
   \begin{align*}
     \la \gd_x, u \ra_\energy = Au(x).
  \end{align*}
  This point of view is helpful, especially when distinguishing between identities in Hilbert space and pointwise equations. For example, if $h \in \Harm$, then $\Lap h$ and the constant function \one are identified in \HE because $\la u, \Lap h \ra_\energy = \la u, \one\ra_\energy = 0$, for any $u \in \HE$. However, one should not consider a (pointwise) solution of $\Lap u(x) = 1$ to be a harmonic function.
\end{remark}

\begin{lemma}\label{thm:vx-is-dipole}
  For any $x \in \verts$, $\Lap v_x = \gd_x - \gd_o$.
  \begin{proof}
    Using Lemma~\ref{thm:<delta_x,v>=Lapv(x)}, $\Lap v_x(y) = \la \gd_y, v_x\ra_\energy = \gd_y(x) - \gd_y(o) = (\gd_x-\gd_o)(y)$.
  \end{proof}
\end{lemma}

By applying Lemma~\ref{thm:vx-is-dipole} to $v_\ga-v_\gw$, we see:

\begin{cor}\label{thm:Pot-nonempty-via-Riesz}
  The space of dipoles $\Pot(\ga,\gw)$ is nonempty.
\end{cor}

Lemma~\ref{thm:<delta_x,v>=Lapv(x)} is extremely important. Since \Fin is the closure of $\spn\{\gd_x\}$, it implies that the finitely supported functions and the harmonic functions are orthogonal. This result is called the ``Royden Decomposition'' in \cite[\S{VI}]{Soardi94} and also appears elsewhere, e.g., \cite[\S9.3]{Lyons:ProbOnTrees}.

\begin{theorem}\label{thm:HE=Fin+Harm}
  $\HE = \Fin \oplus \Harm$.
  \begin{proof}
    For all $v \in \HE$, Lemma~\ref{thm:<delta_x,v>=Lapv(x)} gives $\la \gd_x, v \ra_\energy = \Lap v(x)$. Since $\Fin = \spn\{\gd_x\}$, this equality shows $v \perp \Fin$ whenever $v$ is harmonic. Conversely, if $\la \gd_x, v \ra_\energy=0$ for every $x$, then $v$ must be harmonic. Recall that constants functions are 0 in \HE.
  \end{proof}
\end{theorem}

\begin{cor}\label{thm:Diracs-not-dense}
  $\spn\{\gd_x\}$ is dense in \HE iff $\Harm=0$.
\end{cor}

\begin{remark}\label{rem:Diracs-not-dense}
  Corollary~\ref{thm:Diracs-not-dense} is immediate from Theorem~\ref{thm:HE=Fin+Harm}, but we wish to emphasize the point, as it is not the usual case elsewhere in the literature. Part of the importance of the energy kernel $\{v_x\}$ arises from the fact that the Dirac masses are generally inadequate as a representing set for \HE. This leads to unusual consequences, e.g., one may have
  \linenopax
  \begin{align*}
    u \neq \sum_{x \in \verts} u(x) \gd_x, 
    \q \text{in \HE}.
  \end{align*}
  More precisely, $\|u - \sum_{x \in G_k} u(x) \gd_x\|_\energy$ may not tend to 0 as $k \to \iy$, for some exhaustion $\{G_k\}$.
\end{remark}

\begin{defn}\label{def:ux}
  Let $f_x = \Pfin v_x$ denote the image of $v_x$ under the projection to \Fin. Similarly, let $h_x = \Phar v_x$ denote the image of $v_x$ under the projection to \Harm. 
\end{defn}
  \glossary{name={$f_x$},description={projection of $v_x$ to \Fin},sort=f,format=textbf}
  \glossary{name={$h_x$},description={projection of $v_x$ to \Harm},sort=h,format=textbf}
  \glossary{name={\Pfin},description={projection to \Fin},sort=Pf,format=textbf}
  \glossary{name={\Phar},description={projection to \Harm},sort=Ph,format=textbf}

For future reference, we state the following immediate consequence of orthogonality.
\begin{lemma}\label{thm:repkernels-for-Fin-and-Harm}
  With $f_x = \Pfin v_x$, $\{f_x\}_{x \in \verts}$ is a reproducing kernel for \Fin, but $f_x \perp \Harm$. Similarly, with $h_x = \Phar v_x$, $\{h_x\}_{x \in \verts}$ is a reproducing kernel for \Harm, but $h_x \perp \Fin$.
\end{lemma}


\begin{remark}\label{rem:vx-vs-dx}
  The role of $v_x$ in \HE with respect $\la \cdot, \cdot \ra_\energy$ is directly analogous to role of the Dirac mass $\gd_x$ in $\ell^2$ with respect to the usual $\ell^2$ inner product. This analogy will be developed further when we show that $v_x$ is the image of $x \in \verts$ under a certain isometric embedding into \HE, in \S\ref{sec:Construction-of-HE}. It is obvious that $\gd_x \in \HE$, and the following result shows that $\gd_y$ is always in $\spn\{v_x\}$ when $\deg(y) < \iy$. However, it is not true that $v_y$ is always in $\spn\{\gd_x\}$, or even in its closure. This is discussed further in \S\ref{sec:Construction-of-HE}.
\end{remark}

\begin{lemma}\label{thm:dx-as-vx}
  For any $x \in \verts$, $\gd_x = \cond(x) v_x - \sum_{y \nbr x} \cond_{xy} v_y$.
  \begin{proof}
    Lemma~\ref{thm:<delta_x,v>=Lapv(x)} implies $\la \gd_x, u\ra_\energy = \la \cond(x) v_x - \sum_{y \nbr x} \cond_{xy} v_y, u\ra_\energy$ for every $u \in \HE$, so apply this to $u=v_z$, $z \in \verts$. Since $\gd_x, v_x \in \HE$, it must also be that $\sum_{y \nbr x} \cond_{xy} v_y \in \HE$.
  \end{proof}
\end{lemma}

\subsection{Real and complex-valued functions on \verts}
\label{sec:Real-and-complex-valued-functions-on-verts}

While we will need complex-valued functions for some later results concerning spectral theory, it will usually suffice to consider \bR-valued functions elsewhere.

\begin{lemma}\label{thm:vx-is-R-valued}\label{rem:3-repkernels}
  The reproducing kernels $v_x, f_x, h_x$ are all \bR-valued functions.
  \begin{proof}
    Computing directly,
    \linenopax
    \begin{align*}
      \la \cj{v_z}, \cj{u}\ra_\energy
      &= \frac12\sum_{x,y \in \verts} (v_z(x)-v_z(y))(\cj{u}(x)-\cj{u}(y)) 
      = \cj{\la v_z, u\ra_\energy}.
    \end{align*}
    Then applying the reproducing kernel property,
    \linenopax
    \begin{align*}
      \cj{\la v_z, u \ra_\energy}
      = \cj{u(x)-u(o)} 
      = \cj{u}(x)-\cj{u}(o) 
      = \la v_z, \cj{u}\ra_\energy.
    \end{align*}
    Thus $\la \cj{v_z}, \cj{u}\ra_\energy = \la v_z, \cj{u}\ra_\energy$ for every $u \in \Harm$, and $v_z$ must be \bR-valued. The same computation applies to $f_z$ and $h_z$.
  \end{proof}
\end{lemma}

\begin{defn}\label{def:pointwise-convergence-in-HE}
  A sequence of functions $\{u_n\} \ci \HE$ \emph{converges pointwise in \HE} iff $\exists k \in \bC$ such that $u_n(x)-u(x) \to k$, for each $x \in \verts$.
\end{defn}

\begin{lemma}\label{thm:E-convergence-implies-pointwise-convergence}
  If $\{u_n\}$ converges to $u$ in \energy, then $\{u_n\}$ converges to $u$ pointwise in \HE.
  \begin{proof}
    Define $w_n := u_n-u$ so that $\|w_n\|_\energy \to 0$. Then
    \linenopax
    \begin{align*}
      |w_n(x)-w_n(o)|
      = |\la v_x, w_n\ra_\energy|
      \leq \|v_x\|_\energy \cdot \|w_n\|_\energy
      \limas{n} 0,
    \end{align*}
    so that $\lim w_n$ exists pointwise and is a constant function.
  \end{proof}
\end{lemma}

\section{The discrete Gauss-Green formula}
\label{sec:relating-energy-form-to-Laplacian}

A key difference between our development of the relationship between the Laplace operator \Lap and the Dirichlet energy form \energy (embodied in the discrete Gauss-Green formula of Theorem~\ref{thm:E(u,v)=<u,Lapv>+sum(normals)}) is that \Lap is Hermitian but not necessarily self-adjoint in the present context. This is in sharp contrast to the literature on resistance forms \cite{Kig03}, the general theory of Dirichlet forms and probability \cite{FOT94,BouleauHirsch}, and Dirichlet spaces in potential theory \cite{Brelot,ConCor}. In fact, the ``gap'' between \Lap and its self-adjoint extensions comprises an important part of the boundary theory for $(\Graph,\cond)$, and accounts for features of the boundary terms in the discrete Gauss-Green identity of Theorem~\ref{thm:E(u,v)=<u,Lapv>+sum(normals)}.

Before completing the extension of Lemma~\ref{thm:E(u,v)=<u,Lapv>} to infinite networks, we need some definitions. 

\begin{defn}\label{def:monopole}
  A \emph{monopole} at $x \in \verts$ is an element $w_x \in \HE$ which \sats
  $\Lap w_x(y) = \gd_{xy}$, where $k \in \bC$ and $\gd_{xy}$ is Kronecker's delta. 
  When nonempty, the set of monopoles at the origin is closed and convex, so \energy attains a unique minimum here; let $w_o$ always denote the unique energy-minimizing monopole at the origin.  
  
  When \HE contains monopoles, let $\MP_x$ denote the vector space spanned by the monopoles at $x$. This implies that $\MP_x$ may contain harmonic functions; see Lemma~\ref{thm:MP-contains-spans}.
  We indicate the distinguished monopoles 
  \linenopax
  \begin{align}\label{eqn:def:monov-and-monof}
    \monov := v_x + w_o
    \q\text{and}\q
    \monof := f_x + w_o, 
  \end{align}
  where $f_x = \Pfin v_x$. (Corollary~\ref{thm:Harm-nonzero-iff-multiple-monopoles} below confirms that $\monov = \monof$ for all $x$ iff if $\Harm=0$.) 
\end{defn}

\begin{remark}\label{rem:w_O-in-Fin}
  Note that $w_o \in \Fin$, whenever it is present in \HE, and similarly that \monof is the energy-minimizing element of $\MP_x$. To see this, suppose $w_x$ is any monopole at $x$. Since $w_x \in \HE$, write $w_x = f+h$ by Theorem~\ref{thm:HE=Fin+Harm}, and get $\energy(w_x) = \energy(f) + \energy(h)$. Projecting away the harmonic component will not affect the monopole property, so $\monof = \Pfin w_x$ is the unique monopole of minimal energy. Also, $w_o$ corresponds to the projection of \one to \Gddo; see \S\ref{sec:grounded-energy-space}.
\end{remark}

\begin{defn}\label{def:LapM}
  The dense subspace of \HE spanned by monopoles (or dipoles) is
  \linenopax
  \begin{equation}\label{eqn:def:MP}
    \MP := \spn\{v_x\}_{x \in \verts} + \spn\{\monov, \monof\}_{x \in \verts}.
  \end{equation}
   
  Let \LapM be the closure of the Laplacian when taken to have the dense domain \MP. 
\end{defn}
\glossary{name={$\MP$},description={the span of the monopoles, i.e., finite linear combinations of $w_x$'s},sort=M,format=textbf}
\glossary{name={\LapM},description={the closure of the Laplacian when taken to have the dense domain $\MP$},sort=L,format=textbf}
  Note that $\MP = \spn\{v_x\}$ when there are no monopoles (i.e., when all solutions of of $\Lap w = \gd_x$ have infinite energy), and that $\MP = \spn\{\monov, \monof\}$ when there are monopoles; see Lemma~\ref{thm:MP-contains-spans}.

The space \MP is introduced as a dense domain for \Lap and for its use as a hypothesis in our main result, that is, as the largest domain of validity for the discrete Gauss-Green identity of Theorem~\ref{thm:E(u,v)=<u,Lapv>+sum(normals)}. Note that while a general monopole need not be in $\dom \LapM$ (see \cite[Ex.~13.8 or Ex.~14.39]{OTERN}), we show in Lemma~\ref{thm:pointwise-identity-implies-adjoint-identity} that it is always the case that it lies in $\dom \LapM^\ad$.

\begin{defn}\label{def:semibounded}
  A Hermitian operator $S$ on a Hilbert space \sH is called \emph{semibounded} iff
  \linenopax
  \begin{equation}\label{eqn:def:semibounded}
    \la v,Sv\ra \geq 0, \qq\text{for every } v \in \sD,
  \end{equation}
  so that its spectrum lies in some halfline $[\gk,\iy)$ and its defect indices agree.
\end{defn}

\begin{lemma}\label{thm:LapM-is-semibounded}
  \LapM is Hermitian; a fortiori, \LapM is semibounded. 
  \begin{proof}
    Suppose we have two finite sums $u = \sum a_x w_x$ and $v = \sum b_y w_y$, writing $w_x$ for \monov or \monof. We may assume that $o$ appears neither in the sum $u$ nor for $v$; see Definition~\ref{def:energy-kernel}. Then Lemma~\ref{thm:<delta_x,v>=Lapv(x)} gives
    \linenopax
    \begin{align*}
      \la u, \Lap v\ra_\energy
      = \sum \cj{a_x} b_y \la w_x, \Lap w_y\ra_\energy
      = \sum \cj{a_x} b_y\la w_x, \gd_y\ra_\energy
      = \sum \cj{a_x} b_y \Lap w_x(y)
      = \sum \cj{a_x} b_y \gd_{xy}.   
    \end{align*}
    Of course, $\la \Lap u, v\ra_\energy = \sum \cj{a_x} b_y \gd_{xy}$ exactly the same way. The argument for linear combinations from $\{v_x\}$ is similar, so \LapM is Hermitian.
    Then 
    \linenopax
    \begin{align*}
      \la u, \Lap u\ra_\energy
      = \sum_{x,y} \cj{a_x} a_y \gd_{xy}
      = \sum_{x} |a_x|^2 \geq 0
    \end{align*}
    shows \LapM is semibounded. The argument for $\{v_x\}$ is similar.
  \end{proof}
\end{lemma}

\begin{remark}[Monopoles and transience]
  \label{rem:transient-iff-monopoles}
  The presence of monopoles in \HE is equivalent to the transience of the underlying network, that is, the transience of the simple random walk on the network with transition probabilities $p(x,y) = \cond_{xy}/\cond(x)$. To see this, note that if $w_x$ is a monopole, then the current induced by $w_x$ is a unit flow to infinity with finite energy. It was proved in \cite{TerryLyons} that the network is transient if and only if there exists a unit current flow to infinity; see also \cite[Thm.~2.10]{Lyons:ProbOnTrees}. It is also clear that the existence of a monopole at one vertex is equivalent to the existence of a monopole at every vertex: consider $v_x + w_o$. The corresponding statement about transience is well-known.
\end{remark}

Since \Lap agrees with \LapM pointwise, we may suppress reference to the domain for ease of notation. When given a pointwise identity $\Lap u = v$, there is an associated identity in \HE, but the next lemma shows that one must use the adjoint.

\begin{lemma}\label{thm:pointwise-identity-implies-adjoint-identity}
  For $u,v \in \HE$, $\Lap u = v$ pointwise if and only if $v = \LapM^\ad u$ in \HE.
  \begin{proof}
    We show that $u \in \dom \LapM^\ad$ for simplicity, so let $\gf \in \spn\{v_x\}$ be given by $\gf = \sum_{i=1}^n a_i v_{x_i}$; the proof for $\gf \in \spn\{\monov,\monof\}$ is similar. Then Lemma~\ref{thm:<delta_x,v>=Lapv(x)} and Lemma~\ref{thm:vx-is-dipole} give
    \linenopax
    \begin{align*}
      \la \Lap \gf, u\ra_\energy
      = \sum_{i=1}^n a_i \la \gd_{x_i} - \gd_o, u\ra_\energy
      &= \sum_{i=1}^n a_i (\Lap u(x_i) - \Lap u(o)) .
    \end{align*}
    Since $\Lap u(x) = v(x)$ by hypothesis, this may be continued as
    \linenopax
    \begin{align*}
      \la \Lap \gf, u\ra_\energy
      &= \sum_{i=1}^n a_i (v(x_i) - v(o))       
       = \sum_{i=1}^n a_i \la v_{x_i}, v\ra_\energy 
       = \la \gf, v\ra_\energy.
    \end{align*}
    Then the Schwarz inequality gives the estimate $\left|\la \Lap \gf, u\ra_\energy\right| = \left|\la \gf, v\ra_\energy\right| \leq \|\gf\|_\energy \|v\|_\energy$, which means $u \in \dom \LapM^\ad$.
    The converse is trivial.
  \end{proof}
\end{lemma}

\begin{remark}[Monopoles give a reproducing kernel for $\ran\LapM$]
  \label{rem:monopole-definition}\label{rem:ranLap-vs-Fin}
  Lemma~\ref{thm:pointwise-identity-implies-adjoint-identity} means that 
  \linenopax
  \begin{equation}\label{eqn:def:monopole}
    \la w_x, \Lap u \ra_\energy = \la \gd_x, u \ra_\energy,
    \qq \text{ for all } u \in \dom \LapM.
  \end{equation}
  for every $w_x \in \MP_x$. Combined with Lemma~\ref{thm:<delta_x,v>=Lapv(x)}, this immediately gives  
  \linenopax
  \begin{equation}\label{eqn:monopole-as-repkernel-for-ranLap}
    \la w_x, \Lap u \ra_\energy = \Lap u(x).
  \end{equation}
 If $\{w_x\}_{x \in \verts}$ is a collection of monopoles which includes one element from each $\MP_x$, then this collection is a reproducing kernel for $\ran \LapM$. Note that the expression $\Lap u(x)$ is defined in terms of differences, so the right-hand side is well-defined even without reference to another vertex, i.e., independent of any choice of representative.
  
  As a special case, let $w_x^o$ be the representative of \monof which satisfies $w_x^o(o)=0$. Then the Green function is $g(x,y) = w_y^o(x)$, and $\{w_x^o\}_{x \in \verts \less \{o\}}$ gives a reproducing kernel for $\ran \LapM \ci \Fin$. Therefore, \MP contains an extension of the Green kernel to all of \HE.
 
In Definition~\ref{def:LapM}, we give a domain \MP for \Lap which ensures that $\ran \LapM$ contains all finitely supported functions and is thus dense in \Fin. However, even when \Lap is defined so as to be a closed operator, one may not have $\Fin \ci \ran \Lap$; in general, the containment $\ran (\opclosure{S}) \ci \opclosure{(\ran S)}$ may be strict. The operator closure $\opclosure{S}$ is done with respect to the graph norm, and the closure of the range is done with respect to \energy. 
  We note that \cite[(G.1)]{MuYaYo} claims that the Green function is a reproducing kernel for all of \Fin. In our context, at least, the Green function is a reproducing kernel only for $\ran \Lap$, where \Lap has been chosen with a suitable dense domain. In general, the containment $\ran \Lap \ci \Fin$ may be strict. In fact, it is true that $\ran \LapM^\ad \ci \Fin$, and even this containment may be strict. Note that \monof is the only element of $\MP_x$ which lies in $\opclosure{(\ran \LapM)}$, and it may not lie in ${\ran \LapM}$.

   A different choice of domain for \Lap can exacerbate the discrepancy between $\ran \Lap$ and \Fin: if one were to define \LapV to be the closure of \Lap when taken to have dense domain $\sV := \spn\{v_x\}$ (as the authors did initially), then $\ran \LapV$ is dense in $\Fin_2$, the \energy-closure of $\spn\{\gd_x-\gd_o\}$. However, it can happen that $\Fin_2$ is a proper orthogonal subspace of \Fin (the \energy-closure of $\spn\{\gd_x\}$). An example of $f \in \Fin_1 := \Fin \ominus \Fin_2$ is computed in Example~\ref{exm:Fin_2-not-dense-in-Fin}. The domain of \Lap can thus induce a refinement of the Royden decomposition:
  \linenopax
  \begin{align*}
    \HE = \Fin_1 \oplus \Fin_2 \oplus \Harm.
  \end{align*}
  See Theorem~\ref{thm:HE=Fin+Harm} and the comment preceding it.
\end{remark}
  \glossary{name={$w_x$},description={monopole at $z$},sort=w,format=textbf}
  \glossary{name={$\MP$},description={span of monopoles},sort=M,format=textbf}

Note that a monopole need not be in $\dom \LapV$; see Example~\ref{exm:Harm-notin-domLap} or Example~\ref{exm:defective-integers}. However, it is always the case that $w_x \in \dom \LapV^\ad$, which is the content of the following lemma.

\begin{remark}[Monopoles give a reproducing kernel for $\ran\LapM$]
  \label{rem:monopole-definition}\label{rem:ranLap-vs-Fin}
  Lemma~\ref{thm:pointwise-identity-implies-adjoint-identity} means that 
  \linenopax
  \begin{equation}\label{eqn:def:monopole}
    \la w_x, \Lap u \ra_\energy = \la \gd_x, u \ra_\energy,
    \qq \text{ for all } u \in \dom \LapM.
  \end{equation}
  for every $w_x \in \MP_x$. Combined with Lemma~\ref{thm:<delta_x,v>=Lapv(x)}, this immediately gives  
  \linenopax
  \begin{equation}\label{eqn:monopole-as-repkernel-for-ranLap}
    \la w_x, \Lap u \ra_\energy = \Lap u(x).
  \end{equation}
 If $\{w_x\}_{x \in \verts}$ is a collection of monopoles which includes one element from each $\MP_x$, then this collection is a reproducing kernel for $\ran \LapM$. Note that the expression $\Lap u(x)$ is defined in terms of differences, so the right-hand side is well-defined even without reference to another vertex, i.e., independent of any choice of representative.
  
  As a special case, let $w_x^o$ be the representative of \monof which satisfies $w_x^o(o)=0$. Then the Green function is $g(x,y) = w_y^o(x)$, and $\{w_x^o\}_{x \in \verts \less \{o\}}$ gives a reproducing kernel for $\ran \LapM \ci \Fin$. Therefore, \MP contains an extension of the Green kernel to all of \HE.
\end{remark}

\begin{defn}\label{def:subgraph-boundary}
  If $H$ is a subgraph of $G$, then the boundary of $H$ is
  \linenopax
  \begin{equation}\label{eqn:subgraph-boundary}
    \bd H := \{x \in H \suth \exists y \in H^\complm, y \nbr x\}.
  \end{equation}
  \glossary{name={$\bd H$},description={boundary of a subgraph},sort=b,format=textbf}
  The \emph{interior} of a subgraph $H$ consists of the vertices in $H$ whose neighbours also lie in $H$:
  \linenopax
  \begin{equation}\label{eqn:interior}
    \inn H := \{x \in H \suth y \nbr x \implies y \in H\} = H \less \bd H.
  \end{equation}
  \glossary{name={$\inn H$},description={interior of a subgraph},sort=i,format=textbf}
  For vertices in the boundary of a subgraph, the \emph{normal derivative} of $v$ is
  \linenopax
  \begin{equation}\label{eqn:sum-of-normal-derivs}
    \dn v(x) := \sum_{y \in H} \cond_{xy} (v(x) - v(y)),
    \qq \text{for } x \in \bd H.
  \end{equation}
  Thus, the normal derivative of $v$ is computed like $\Lap v(x)$, except that the sum extends only over the neighbours of $x$ which lie in $H$.
\end{defn}
  \glossary{name={$\dn v$},description={normal derivative of a function with respect to a subgraph},sort=d,format=textbf}

Definition~\ref{def:subgraph-boundary} will be used primarily for subgraphs that form an exhaustion of \Graph, in the sense of Definition~\ref{def:exhaustion-of-G}: an increasing sequence of finite and connected subgraphs $\{\Graph_k\}$, so that $\Graph_k \ci \Graph_{k+1}$ and $\Graph = \bigcup \Graph_k$. Also, recall that $\sum_{\bd \Graph} := \lim_{k \to \iy} \sum_{\bd \Graph_k}$ from Definition~\ref{def:boundary-sum}.

\begin{defn}\label{def:boundary-sum}
  A \emph{boundary sum} (or \emph{boundary term}) is computed in terms of an exhaustion $\{G_k\}$ by
  \linenopax
  \begin{equation}\label{eqn:boundary-sum}
    \sum_{\bd \Graph} := \lim_{k \to \iy} \sum_{\bd \Graph_k},
  \end{equation}
  whenever the limit is independent of the choice of exhaustion, as in Definition~\ref{def:infinite-vertex-sum}. The boundary $\bd\Graph$ is examined more closely as an object in its own right in \S\ref{sec:the-boundary}.
  \glossary{name={$\sum_{\bd \Graph}$},description={boundary sum},sort=S,format=textbf}
\end{defn}

The key point of the following result is that for $u,v$ in the specified set, the two sums are both finite. The decomposition is true for all $u,v \in \HE$ by taking limits of \eqref{eqn:<u,v>_k-decomp-2}, but is clearly meaningless if it takes the form $\iy - \iy$.
  
\begin{theorem}[Discrete Gauss-Green Formula]
  \label{thm:E(u,v)=<u,Lapv>+sum(normals)}
  If $u \in \HE$ and $v \in \MP$, then
  \linenopax
  \begin{equation}\label{eqn:E(u,v)=<u_0,Lapv>+sum(normals)}
    \la u, v \ra_\energy
    = \sum_{x \in \verts} \cj{u}(x) \Lap v(x)
      + \sum_{x \in \bd \Graph} \cj{u}(x) \dn v(x).
  \end{equation}
  \begin{proof}
    It suffices to work with \bR-valued functions and then complexify afterwards. 
    By the same computation as in Lemma~\ref{thm:E(u,v)=<u,Lapv>}, we have
    \linenopax
    \begin{align}\label{eqn:<u,v>_k-decomp-2}
      \frac12 \sum_{x,y \in G_k} \cond_{xy} (\cj{u}(x)-\cj{u}(y))(v(x)-v(y))
      &= \sum_{x \in \inn \Graph_k} \cj{u}(x) \Lap v(x)
       + \sum_{x \in \bd \Graph_k} \cj{u}(x) \dn v(x).
    \end{align}
    
    Taking limits of both sides as $k \to \iy$ gives \eqref{eqn:E(u,v)=<u_0,Lapv>+sum(normals)}. It remains to see that one of the sums on the right-hand side is finite (and hence that both are). For this part, we work just with $u$ and polarize afterwards. Note that if $v =w_z$ is a monopole, then 
    \linenopax
    \begin{align*}
      \sum_{x \in \verts} {u}(x) \Lap v(x)
      = \sum_{x \in \verts} {u}(x) \gd_z(x)
      = u(z).
    \end{align*}
    This is obviously independent of exhaustion, and immediately extends to $v \in \MP$.   \end{proof}
\end{theorem}

\begin{remark}\label{rem:DGG-holds-without-hypotheses}
  It is clear that \eqref{eqn:E(u,v)=<u_0,Lapv>+sum(normals)} remains true much more generally than under the specified conditions. Clearly, the formula holds whenever $\sum_{x \in \verts} \left|{u}(x) \Lap v(x)\right| < \iy$.
  Unfortunately, given any hypotheses more specific than this, the limitless variety of infinite networks almost always allow one to construct a counterexample; i.e. one cannot give a condition for which the formula is true for all $u \in \HE$, for all networks. To see this, suppose that $v = \sum_{i=1}^\iy a_i w_{x_i}$ with each $w_{x_i}$ a monopole at the vertex $x_i$. Then
  \linenopax
  \begin{align*}
    \sum_{x \in \verts} {u}(x) \Lap v(x)
    = \sum_{i=1}^\iy a_i {u}(x_i),
  \end{align*}
  and one would need to provide a condition on sequences $\{a_i\}$ that would ensure $\sum_{i=1}^\iy a_i {u}(x_i)$ is absolutely convergent for all $u \in \HE$. Such a hypothesis is not likely to be useful (if it is even possible to construct) and would depend heavily on the network under investigation.
  Nonetheless, the formula remains true in many specific contexts. For example, it is clearly valid whenever $v$ is a dipole, including all those in the energy kernel. We will also see that it holds for the projections of $v_x$ to \Fin and to \Harm. Consequently, for $v$ which are limits of elements in \MP, we often use this result in combination with ad hoc arguments. 

  After reading a preliminary version of this paper, a reader pointed out to us that a formula similar to \eqref{eqn:E(u,v)=<u_0,Lapv>+sum(normals)} appears in \cite[Prop~1.3]{DodziukKarp88}; however, these authors apparently do not pursue the extension of this formula to infinite networks. We were also directed towards \cite[Thm.~4.1]{Kayano88}, in which the authors give some conditions under which Lemma~\ref{thm:E(u,v)=<u,Lapv>} extends to infinite networks. The main differences here are that the scope of Kayano and Yamasaki's theorem is limited to a subset of what we call \Fin, and that Kayano and Yamasaki are interested in when the boundary term vanishes; we are more interested in when it is finite and nonvanishing; see Theorem~\ref{thm:TFAE:Fin,Harm,Bdy}, for example. Since Kayano and Yamasaki do not discuss the structure of the space of functions they consider, it is not clear how large the scope of their result is; their result requires the hypothesis $\sum_{x \in \verts} \left|{u}(x) \Lap v(x)\right| < \iy$, but it is not so clear what functions satisfy this. By contrast, we develop a dense subspace of functions on which to apply the formula. Furthermore, we show in the next chapter that these functions are relatively easy to compute.
\end{remark}

Recall that $\spn\{h_x\}$ is a dense subspace of \Harm; the following boundary representation of harmonic functions in this space is the focus of Chapter~\ref{sec:the-boundary}.

\begin{cor}[Boundary representation of harmonic functions]
  \label{thm:Boundary-representation-of-harmonic-functions}
  For all $u \in \spn\{h_x\}$, 
  \begin{align}\label{eqn:Boundary-representation-of-harmonic-functions}
    u(x) = \sum_{\bd \Graph} u \dn{h_x} + u(o). 
  \end{align}
  \begin{proof}
    Note that 
    $u(x) - u(o) 
    = \la v_x, u\ra_\energy 
    = \cj{\la u, v_x\ra_\energy} 
    = \sum_{\bd \Graph} u \dn{h_x}$ by \eqref{eqn:def:vx}.
  \end{proof}
\end{cor}

\begin{lemma}\label{thm:sumLap(u)=-sumdn(u)}
  For all $u \in \dom \LapV$, $\sum_{\verts} \Lap u = - \sum_{\bd \Graph} \dn u$. Thus, the Discrete Gauss-Green formula \eqref{eqn:E(u,v)=<u_0,Lapv>+sum(normals)} is independent of representatives.
  \begin{proof}
    On each (finite) $G_k$ in any given exhaustion,
    \linenopax
    \begin{align*}
      \sum_{x \in \inn G_k} \Lap u(x) + \sum_{x \in \bd G_k} \dn u(x)
      = \sum_{x,y \in G_k} \cond_{xy}(u(x)-u(y)) = 0,
    \end{align*}
    since each edge appears twice in the sum; once with each sign (orientation). For the second claim, we apply the formula of the first to see that the result remains true when $u$ is replaced by $u+k$: 
  \linenopax
  \begin{align*}
    \sum_{\verts} (u+k)\Lap v + \sum_{\bd G} (u+k) \dn v
    &= \sum_{\verts} u\Lap v + \sum_{\bd G} u \dn v
     +k\cancel{\left(\sum_{\verts} \Lap v + \sum_{\bd G} \dn v\right)}.
     \qedhere
  \end{align*}
  \end{proof}
\end{lemma}

\section{More about monopoles and the space \MP}
\label{sec:More-about-monopoles}

This section studies the role of the monopoles with regard to the boundary term of Theorem~\ref{thm:E(u,v)=<u,Lapv>+sum(normals)}, and provides several characterizations of transience of the network, in terms the operator-theoretic properties of \LapM. 
  
Note that if $h \in \Harm$ satisfies the hypotheses of Theorem~\ref{thm:E(u,v)=<u,Lapv>+sum(normals)}, then $\energy(h) = \sum_{\bd\Graph} h \dn h$. In Theorem~\ref{thm:TFAE:Fin,Harm,Bdy} we show that $\energy(u) = \sum_{\verts} u \Lap u$ for all $u \in \HE$ iff the network is recurrent. With respect to $\HE = \Fin \oplus \Harm$, this shows that the energy of finitely supported functions comes from the sum over \verts, and the energy of harmonic functions comes from the boundary sum. However, for a monopole $w_x$, the representative specified by $w_x(x)=0$ satisfies $\energy(w) = \sum_{\bd\Graph} w \dn w$ but the representative specified by $w_x(x) = \energy(w_x)$ satisfies $\energy(w) = \sum_{\verts} w \Lap w$. Roughly speaking, a monopole is therefore ``half of a harmonic function'' or halfway to being a harmonic function. A further justification for this comment is given by Corollary~\ref{thm:Harm-nonzero-iff-multiple-monopoles}: the proof shows that a harmonic function can be constructed from two monopoles at the same vertex. A different perspective one the same theme is given in Remark~\ref{rem:comparison-to-Royden-decomp}. The general theme of this section is the ability of monopoles to ``bridge'' the finite and the harmonic.

\begin{remark}\label{rem:transient-iff-monopoles}
  The presence of monopoles in \HE is equivalent to the transience of the underlying network, that is, the transience of the simple random walk on the network with transition probabilities $p(x,y) = \cond_{xy}/\cond(x)$. To see this, note that if $w_x$ is a monopole, then the current induced by $w_x$ is a unit flow to infinity with finite energy. It was proved in \cite{TerryLyons} that the network is transient if and only if there exists a unit current flow to infinity; see also \cite[Thm.~2.10]{Lyons:ProbOnTrees}. As mentioned in Corollary~\ref{thm:transient-here-implies-transient-everywhere}, the existence of a monopole at one vertex is equivalent to the existence of a monopole at every vertex. 
\end{remark}

\begin{lemma}\label{thm:MP-contains-spans}
  When the network is transient, \MP contains the spaces $\spn\{v_x\}$, $\spn\{f_x\}$, and $\spn\{h_x\}$, where $f_x = \Pfin v_x$ and $h_x = \Phar v_x$.
  \begin{proof}
    The first two are obvious, since $v_x = \monov - w_o$ and $f_x = \monof - w_o$ by Definition~\ref{def:monopole}. For the harmonics, note that these same identities give
    \linenopax
    \begin{align*}
      \monov - w_o = v_x = f_x + h_x = \monof - w_o + h_x,
    \end{align*}
    which implies that $h_x = \monov - \monof$. (Of course, $\monov = \monof$ when $\Harm=0$.)
  \end{proof}
\end{lemma}

\begin{cor}\label{thm:Harm-nonzero-iff-multiple-monopoles}
  $\Harm \neq 0$ iff there is more than one monopole at $x$. 
  \begin{proof}
    As usual, if this is true for any $x$, it is true for all. Suppose \HE contains a monopole $w_x \neq \monov$. Then $h := \monov - w_x$ is a nonzero harmonic function in \HE.
  \end{proof}
\end{cor}

\begin{theorem}[{\cite[Thm.~1.33]{Soardi94}}]
  \label{thm:Soardi's-harmonic-implies-transient}
  Let $u$ be a nonnegative function on a recurrent network. Then $u$ is superharmonic if and only if $u$ is constant.  
\end{theorem}

\begin{cor}\label{thm: harmonic-implies-transient}
  If $\Harm \neq 0$, then there is a monopole in \HE.
  \begin{proof}
    If $h \in \Harm$ and $h \neq 0$, then $h = h_1 - h_2$ with $h_i \in \Harm$ and $h_i \geq 0$, by \cite[Thm.~3.72]{Soardi94}. (Here, $h_i \geq 0$ means that $h_i$ is bounded below, and so we can choose a representative which is nonnegative.) Since the $h_i$ cannot both be 0, Theorem~\ref{thm:Soardi's-harmonic-implies-transient} implies the network is transient. Then by \cite[Thm.~1]{TerryLyons}, the network supports a monopole.
  \end{proof}
\end{cor} 

\begin{defn}\label{def:boundary-term-is-nonvanishing}
  The phrase ``the boundary term is nonvanishing'' indicates that \eqref{eqn:E(u,v)=<u_0,Lapv>+sum(normals)} holds with nonzero boundary sum when applied to $\la u,v\ra_\energy$, for every representative of $u$ except the one specified by $u(x) = \la u,w\ra_\energy$, for $w \in \MP_x$.
\end{defn}

Recall from Remark~\ref{rem:transient-iff-monopoles} that the network is transient iff there are monopoles in \HE.

\begin{theorem}\label{thm:TFAE:Fin,Harm,Bdy}
  The network is transient if and only if the boundary term is nonvanishing. \version{Moreover, the boundary term vanishes for the elements of $\ran \LapM$.}{}
  \begin{proof}
    \fwd If the network is transient, then as explained in Remark~\ref{rem:transient-iff-monopoles}, there is a $w \in \HE$ with $\Lap w = \gd_z$. Now let $w_z := \Pfin w$ so that for any $u \in \dom \LapV$, \eqref{eqn:E(u,v)=<u_0,Lapv>+sum(normals)}
    \linenopax
    \begin{align*}
      \la u, w_z \ra_\energy = u(z) + \sum_{\bd \Graph} u \dn{w_z}.
    \end{align*}
    It is immediate that $\sum_{\bd \Graph} u \dn{w_z} = 0$ if and only if the computation is done with the representative of $u$ specified by $u(z) = \la u, w_z \ra_\energy$.
    
    \bwd Suppose that there does not exist $w \in \HE$ with $\Lap w = \gd_z$, for any $z \in \verts$. Then $\MP = \spn\{v_x\}$ as discussed in Definition~\ref{def:monopole}.
    Therefore, it suffices to show that
    \linenopax
    \begin{align*}
      \la u, v_x\ra_\energy = \sum_{x \in \verts} u \Lap v_x,
    \end{align*}
    but this is clear because both sides are equal to $u(x)-u(o)$ by \eqref{eqn:def:vx} and Lemma~\ref{thm:vx-is-dipole}.
    
    \version{For the final claim, note that if $u \in \ran \LapM$, then \eqref{eqn:monopole-as-repkernel-for-ranLap} gives
    \linenopax
    \begin{align*}
      u(x) 
      = \la u, w_x \ra_\energy
      = \sum_{\verts} u \Lap w_x + \sum_{\bd \Graph} u \dn{w_x}
      = u(x) + \sum_{\bd \Graph} u \dn{w_x},
    \end{align*}
    so that the boundary term must vanish. }{}
  \end{proof}
\end{theorem}

\begin{remark}\label{rem:no-monopoles-in-ran(Lap)}
  It follows from Theorem~\ref{thm:TFAE:Fin,Harm,Bdy} that a monopole $w_x$ cannot lie in $\ran \LapV$. However, one can have $w_x \in \ran \LapV^\ad$, as in Example~\ref{exm:defective-integers}. 
\end{remark}

\begin{lemma}\label{thm:monopole-as-weak-limit-of-inverse}
  The network is transient if and only if there is a sequence $\{\ge_k\}$ with $\ge_k \to 0$ and $\sup_k\|(\ge_k + \Lap)^{-1} \gd_x\|_\energy \leq B < \iy$.
  \begin{proof}
    For both directions of the proof, we let $f_k := (\ge_k + \Lap)^{-1} \gd_x$.
    
    \fwd Let \LapS be any self-adjoint extension\footnote{For concreteness, one may take the Friedrichs extension, see \eqref{eqn:S_{min}<S_{max}} but this is not necessary. See also Definition~\ref{def:extn-of-Lap} and \S\ref{sec:Properties-of-Lap-on-HE} in this regard.} of \LapV, and let $E(d\gl)$ be the corresponding projection-valued measure. Then
    \linenopax
    \begin{align}\label{eqn:res(u)=spectralintegral}
      R_\ge u = (\ge + \LapS)^{-1}u = \int_0^\iy \frac1{\ge+\gl} E(d\gl)u,
    \end{align}
    where we use the notation $R_\ge := (\ge + \LapS)^{-1}$ for the resolvent. Note that $\LapS R_\ge \ci (\LapS R_\ge)^\ad = \LapS^\ad R_\ge^\ad = \LapS R_\ge$. On the other hand, $\LapS \ci \LapV^\ad$ and therefore $R_\ge \LapS \ci R_\ge \LapV^\ad$. Combining these gives $\LapS R_\ge \ci R_\ge \LapV^\ad$.
    Now we apply this and \eqref{eqn:res(u)=spectralintegral} to $u=\Lap^\ad w$ to get
    \linenopax
    \begin{align*}
      f_k = (\ge_k + \LapS)^{-1} \gd_x 
      &= (\ge_k + \LapS)^{-1} \LapV^\ad w 
       = \LapS (\ge_k + \LapS)^{-1} w 
       = \int_0^\iy \frac{\gl}{\ge_k+\gl} E(d\gl)w.
    \end{align*}
    Note that $R_\ge$ is bounded, and so $w \in \dom R_\ge$ automatically. This integral implies
    \linenopax
    \begin{align*}
      \|f_k\|_\energy^2 
      &\leq \int_0^\iy \left(\frac{\gl}{\ge_k+\gl}\right)^2 \|E(d\gl)w\|_\energy^2 
      \leq \int_0^\iy \|E(d\gl)w\|_\energy^2
      = \|w\|_\energy^2.
    \end{align*}
    Thus we have $\sup_k\|(\ge_k + \Lap)^{-1} \gd_x\|_\energy = \sup \|f_k\|_\energy \leq B =\|w\|_\energy < \iy$.
    
    \bwd We show the existence of a monopole at $x$. Since $\ge_k f_k + \Lap f_k = \gd_x$, the bound $\sup \|f_k\|_\energy \leq B$ implies that 
    \linenopax
    \begin{align*}
      \|\Lap f_k - \gd_x\|_\energy 
      = \|\ge_k f_k\| 
      \leq \ge_k B \to 0.
    \end{align*}
    Let $w$ be a weak-$\ad$ limit of $\{f_k\}$. Then for $\gf \in \dom \LapV$,
    \linenopax
    \begin{align*}
      \la \Lap \gf, w\ra_\energy
      = \lim_{k \to \iy} \la \Lap \gf, f_k\ra_\energy
      = \lim_{k \to \iy} \la \gf, \Lap f_k\ra_\energy
      = \lim_{k \to \iy} \la \gf, \gd_x - \ge_k f_k\ra_\energy
      = \la \gf, \gd_x \ra_\energy,
    \end{align*}
    so that $w$ is a monopole at $x$.
  \end{proof}
\end{lemma}

\begin{lemma}
  \label{thm:LapM-maps-into-Fin}
  \label{thm:ker(Lapadj)=Harm}
  On any network, $\opclosure{(\ran \LapM)} \ci \Fin$ and hence $\Harm \ci \ker \LapM^\ad$.   
  \begin{proof}
    If $v \in \MP$, then clearly $\LapM v \in \Fin$. To close the operator, we consider sequences $\{u_n\} \ci \MP$ which are Cauchy in \energy, and for which $\{\Lap u_n\}$ is also Cauchy in \energy, and then include $u := \lim u_n$ in $\dom \LapM$ by defining $\LapM u := \lim \LapM u_n$. Since $f_n := \LapM u_n$ has finite support for each $n$, the \energy-limit of $\{f_n\}$ must lie in \Fin. Since \Fin is closed, the first claim follows. The second claim follows upon taking orthogonal complements. 
  \end{proof}
\end{lemma}

\begin{theorem}\label{thm:transient-iff-Fin=Ran(Lap)}
 The network is transient if and only if $(\ran \LapM^\ad)^{c\ell} = \Fin$.
  \begin{proof}
    \fwd If the network is transient, we have a monopole at every vertex; see Remark~\ref{rem:transient-iff-monopoles}. Then any $u \in \spn\{\gd_x\}$ is in $\ran \LapM^\ad$ because the monopole $w_x$ is in $\dom \LapM$, 
    and so $\Fin \ci \ran \LapM^\ad$. The other inclusion is Lemma~\ref{thm:LapM-maps-into-Fin}.
    
    \bwd If $\gd_x \in \ran \LapM$ for some $x \in \verts$, then $\LapM w = \gd_x$ for $w \in \dom \LapM \ci \dom \energy$ and so $w$ is a monopole. Then the induced current $\drp w$ is a unit flow to infinity, and the network is transient, again by \cite{TerryLyons}.
  \end{proof}
\end{theorem}

\subsection{Comparison with the grounded energy space}
\label{sec:grounded-energy-space}

There are some subtleties in the relationship between \HE and \Gdd as discussed in \cite{Lyons:ProbOnTrees} and \cite{Kayano88,Kayano84,MuYaYo, Yamasaki79}, so we take a moment to give details. We have attempted to match the notation of these sources.

\begin{defn}\label{def:grounded-energy-space}
  The inner product
  \linenopax
  \begin{align*}
    \la u, v \ra_\gdd := \cj{u(o)}v(o) + \la u, v \ra_\energy.
  \end{align*}
  makes $\dom \energy$ into a Hilbert space \Gdd which we call the \emph{grounded energy space}. Let \Gddo be the closure of $\spn\{\gd_x\}$ in \Gdd and let \GHD be the space of harmonic functions in \Gdd.
\end{defn}

Throughout this section (only), we use the notation $u_o := u(o)$, for $u \in \Gdd$.

\begin{defn}\label{def:poles-and-antipoles}
  With regard to \Gdd, we define the vector subspace
  \linenopax
  \begin{align}\label{eqn:def:antipole}
    \MP_o^- := \{u \in \Gdd \suth \Lap u = -u_o \gd_o\}.
  \end{align}
  Note that $\MP_o^-$ contains the harmonic subspace
  \linenopax
  \begin{align}\label{eqn:Gdd-harmonic-subspace}
    \mathbf{HD}_o := \{u \in \Gdd \suth \Lap u = 0 \text{ and } u_o = 0\}.
  \end{align}
\end{defn}

The previous definition is motivated by the following lemma.

\begin{lemma}\label{thm:recurrent-iff-constants}\label{thm:positive-monopoles-in-Do}
  $\Gddo^\perp = \MP_o^-$ and hence $\Gdd = \Gddo \oplus \MP_o^-$. 
  \begin{proof}
    With $u_o := u(o)$, we have $u \in \Gddo^\perp$ iff $u \perp \spn\{\gd_x\}$, which means that
    \linenopax
    \begin{align}\label{eqn:<u,d_x>=0-in-Do}
      0 = \la u, \gd_x\ra_\gdd
      = u_o \gd_x(o) + \la u, \gd_x\ra_\energy
      = u_o \gd_{xo} + \Lap u(x),
      \qq \forall x \in \verts,
    \end{align}
    which means $\Lap u = -u_o\gd_o$.   \end{proof}
\end{lemma}

Let us denote the projection of \Gdd to $\Gddo$ by $P_{\Gddo}$ and the projection to $\Gddo^\perp$ by $P_{\Gddo}^\perp$.

\begin{remark}\label{rem:transient-iff-1-splits}
  The constant function $\one$ decomposes into a linear combination of two monopoles: let $v = P_{\Gddo} \one$ and $u = P_{\Gddo}^\perp \one = \one - v$, and observe that $\Lap u = -u_o \gd_o$ by Lemma~\ref{thm:recurrent-iff-constants} and that $\Lap v = \Lap(\one - u) = -\Lap u = u_o \gd_o$, so $u_o = 1-v_o$ gives $\Lap v = (1-v_o) \gd_o$.  In general, the constant function $k\one$ decomposes into $v=P_{\Gddo} k\one$ and $u=P_{\Gddo}^\perp k\one$, where
  \linenopax
  \begin{align*}
    \Lap v = (k - u_o) \gd_o
    \q\text{ and }\q
    \Lap u =  - u_o \gd_o.
  \end{align*} 
  With respect to the decomposition $\Gdd = \Gddo \oplus \MP_o^-$, given by Lemma~\ref{thm:recurrent-iff-constants}, there are two monopoles $w_o^{(1)} \in \Gddo$ and $w_o^{(2)} \in \MP_o^-$ (which may be equal) such that $\one = u_o w_o^{(1)} - u_o w_o^{(2)}$. When one passes from \Gdd to \HE by modding out constants, these components of \one add together to form (possibly constant) harmonic functions. An example of this is given in Example~\ref{exm:monopolar-decomposition}.

  Consequently, Lemma~\ref{thm:recurrent-iff-constants} yields a short proof of \cite[Exc.~9.6c]{Lyons:ProbOnTrees}: Prove that the network is recurrent iff $\one \in \Gddo$. To see this, observe that if $u$ is the projection of \one to $\Gddo^\perp$, then $u \neq 0$ iff there is a monopole. This result first appeared (in more general form) in \cite[Thm.~3.2]{Yamasaki77}.  
\end{remark}
 
\begin{remark}\label{rem:comparison-to-Royden-decomp}
  Despite the fact that Theorem~\ref{thm:HE=Fin+Harm} gives $\HE = \Fin \oplus \Harm$, note that $\Gdd \neq \Gddo \oplus \GHD$. This is a bit surprising, since $\HE = \Gdd/\bC\one$, etc., and this mistake has been made in the literature, e.g. \cite[Thm.~4.1]{Yamasaki79}. 
  The discrepancy results from the way that \one behaves with respect to $P_{\Gddo}$; this is easiest to see by considering
  \linenopax
  \begin{align*}
    \Gddo + k := \{f + k\one \suth f \in \Gddo, k \in \bC\}, \qq k \neq 0.
  \end{align*}
  If the network is transient and $f \in \Gddo + k$, $k \neq 0$, then $f = g + k \one$ for some $g \in \Gddo$, and
  \linenopax
  \begin{align*}
    f = (g + k P_{\Gddo} \one) + k P_{\Gddo}^\perp \one 
  \end{align*}
  shows $f \notin \Gddo$. Nonetheless, it is easy to check that $\Gddo+k$ is equal to the $o$-closure of $\spn{\gd_x}+k$, and hence that $(\Gddo+\bC\one)/\bC\one = \Fin$. This appears in \cite[Exc.~9.6b]{Lyons:ProbOnTrees}. Similarly, note that for a general $h \in \GHD$, one has $h = P_{\Gddo}^\perp h + k\one$, so that $h \notin \Gddo^\perp$.
\end{remark}

We conclude with a curious lemma that can greatly simplify the computation of monopoles of the form $P_{\Gddo} \one$; it is used in Example~\ref{exm:monopolar-decomposition}. In the next lemma, $u_o = u(o)$, as above.

\begin{lemma}\label{thm:parabolic-u(o)-parameter}
  Let $u \in \Gddo^\perp$. Then $u = P_{\Gddo}^\perp\one$ if and only if $u_o = \energy(u) + u_o^2 \in [0,1)$. 
  \begin{proof}
    From $\|u\|_o^2 + \|\one-u\|_o^2 = \|\one\|_o^2 = 1$, one obtains $\energy(u) - u_o + |u_o|^2 = 0$. From $\la u, \one-u\ra_o = 0$, one obtains $\energy(u) - \cj{u_o} + |u_o|^2 = 0$. Combining the equations gives $u_o = u_o = \frac12(1 \pm \sqrt{1-4\energy(u)})$, so that $u_o \in [0,1]$. However, $u_o \neq 1$ or else $\energy(u)=0$ would imply $\one \in \Gddo^\perp$ in contradiction to \eqref{eqn:<u,d_x>=0-in-Do}. The converse is clear.
  \end{proof}
\end{lemma}

\begin{remark}\label{rem:strength-of-network}
  The significance of the parameter $u_o$ is not clear. However, it appears to be related to the overall ``strength'' of the conductance of the network; we will see in Example~\ref{exm:monopolar-decomposition} that $u_o \approx 1$ corresponds to rapid growth of \cond near \iy. Also, it follows from the Remark~\ref{rem:transient-iff-1-splits} and Lemma~\ref{thm:parabolic-u(o)-parameter} that $u_o = 0$ corresponds to the recurrence. There is probably a good interpretation of $u_o$ in terms of probability and/or the speed of the random walk, but we have not yet determined it. The existence of conductances attaining maximal energy $\energy(P_{\Gddo}^\perp\one) = \frac14$ is similarly intriguing, and even more mysterious. Example~\ref{exm:monopolar-decomposition} shows that the maximum is attained on $(\bZ,c^n)$ for $c=2$.
\end{remark}

\section{Applications and extensions}
\label{sec:Applications-and-extensions}

In \S\ref{sec:More-about-Fin-and-Harm}, we use the techniques developed above to obtain new and succinct proofs of four known results, and in \S\ref{sec:Laplacian-and-its-domain} we give some useful special cases of our main result, Theorem~\ref{thm:E(u,v)=<u,Lapv>+sum(normals)}. 

\begin{defn}\label{def:limit-at-infty}
  For an infinite graph \Graph, we say $u(x)$ \emph{vanishes at \iy} iff for any exhaustion $\{\Graph_k\}$, one can always find $k$ and a constant $C$ such that $\|u(x)-C\|_\iy < \ge$ for all $x \notin \Graph_k$. One can always choose the representative of $u \in \HE$ so that $C=0$, but this may not be compatible with the choice $u(o)=0$.
\end{defn}

\begin{defn}\label{def:path-to-infinity}
  Say $\cpath = (x_0, x_1, x_2,\dots)$ is a \emph{path to \iy} iff $x_i \nbr x_{i-1}$ for each $i$, and for any exhaustion $\{G_k\}$ of \Graph, 
  \linenopax
  \begin{align}\label{eqn:def:path-to-infinity}
    \forall k, \exists N \text{ such that } n \geq N \implies x_n \notin G_k.
  \end{align}
\end{defn}

\subsection{More about \Fin and \Harm}
\label{sec:More-about-Fin-and-Harm}

The next two results are almost converse to each other, although the exact converse of Lemma~\ref{thm:Fin-vanishes-at-infinity} is false; see \cite[Fig.~10 or Ex.~14.16]{OTERN}. Lemma~\ref{thm:Fin-vanishes-at-infinity} is related to \cite[Thm.~3.86]{Soardi94}, in which the result is stated as holding almost everywhere with respect to the notion of extremal length. 

\begin{lemma}\label{thm:Fin-vanishes-at-infinity}
  If $u \in \HE$ and $u$ vanishes at \iy, then $u \in \Fin$. 
  \begin{proof}
    Let $u=f+h \in \HE$ vanish at \iy. This implies that for any exhaustion $\{G_k\}$ and any $\ge > 0$, there is a $k$ and $C$ for which $\|h(x)-C\|_\iy < \ge$ outside $G_k$. A harmonic function can only obtain its maximum on the boundary, unless it is constant, so in particular, \ge bounds $\|h(x)-C\|_\iy$ on all of $G_k$. Letting $\ge \to 0$, $h$ tends to a constant function and $u=f$.
  \end{proof}
\end{lemma}
    
\begin{lemma}\label{thm:Harm-monotonic-to-infinity}
  If $h \in \Harm$ is nonconstant, then from any $x_0 \in \verts$, there is a path to infinity $\cpath = (x_0,x_1,\dots)$, with $h(x_j) < h(x_{j+1})$ for all $j=0,1,2,\dots$.
  \begin{proof}
    Abusing notation, let $h$ be any representative of $h$. Since $h(x) = \sum_{y \nbr x} \frac{\cond_{xy}}{\cond(x)} h(y) \leq \sup_{y \nbr x} h(y)$ and $h$ is nonconstant, we can always find $y \nbr x$ for which $h(y_1) > h(x_0)$. This follows from the maximal principle for harmonic functions; cf. \cite[\S2.1]{Lyons:ProbOnTrees}, \cite[Ex.~1.12]{LevPerWil08}, or \cite[Thm.~1.35]{Soardi94}. Thus, one can inductively construct a sequence which defines the desired path \cpath. Note that \cpath is infinite, so the condition $h(x_j) < h(x_{j+1})$ eventually forces it to leave any finite subset of \verts, so Definition~\ref{def:path-to-infinity} is satisfied.
  \end{proof}
\end{lemma}

It is instructive to prove the contrapositive of Lemma~\ref{thm:Fin-vanishes-at-infinity} directly:

\begin{lemma}\label{thm:h_x-has-two-limiting-values}
  If $h \in \Harm\less\{0\}$, then $h$ has at least two different limiting values at \iy.
  \begin{proof}
    Choose $x \in \verts$ for which $h_x = \Phar v_x \in \HE$ is nonconstant. Then Lemma~\ref{thm:Harm-monotonic-to-infinity} gives a path to infinity $\cpath_1$ along which $h_x$ is strictly increasing. Since the reasoning of Lemma~\ref{thm:Harm-monotonic-to-infinity} works just as well with the inequalities reversed, we also get $\cpath_2$ to \iy along which $h_x$ is strictly decreasing. This gives two different limiting values of $h_x$, and hence $h_x$ cannot vanish at \iy.
  \end{proof}
\end{lemma}

\begin{cor}\label{thm:Harm-notin-ellP}
  If $h \in \Harm$ is nonconstant, then $h \notin \ell^p(\verts)$ for any $1 \leq p < \iy$. 
  \begin{proof}
    Lemma~\ref{thm:h_x-has-two-limiting-values} shows that no matter what representative is chosen for $h$, the sum $\|h\|_p = \sum_{x \in \verts} |h(x)|^p$ has the lower bound $\sum_{x \in F} \ge^p = \ge^p |F|$, for some infinite set $F \ci \verts$.
  \end{proof}
\end{cor}

\subsection{Special applications of the Discrete Gauss-Green formula}
\label{sec:Laplacian-and-its-domain}

In this section, we use Lemma~\ref{thm:E(u,v)=<u,Lapv>} to infinite networks to establish that \Lap is Hermitian when its domain is correctly chosen (Corollary~\ref{thm:Lap-Hermitian-on-V}), and that Lemma~\ref{thm:E(u,v)=<u,Lapv>} remains correct on infinite networks for vectors in $\spn\{v_x\}$ (Theorem~\ref{thm:E(u,v)=<u,Lapv>-on-Fin}).

\begin{lemma}\label{thm:E(u,v)=<u,Lapv>-on-spn{vx}}
  If $v \in \spn\{v_x\}$, then $\la u, v \ra_\energy = \sum_{x \in \verts} \cj{u(x)} \Lap v(x)$.    
  \begin{proof}
    It suffices to consider $v=v_x$, whence
    \linenopax
    \begin{align*}
      \sum_{\verts} u(y) \Lap v_x(y)
      = \sum_{\verts} u(y) (\gd_x - \gd_o)(y)
      = u(x) - u(o)
      = \la u, v_x\ra_\energy,
    \end{align*}
    by Lemma~\ref{thm:vx-is-dipole} and the reproducing property of Corollary~\ref{thm:vx-is-a-reproducing-kernel}.
  \end{proof}
\end{lemma}

\begin{theorem}\label{thm:uLapv-has-no-bdy-term}\label{thm:sum(Lap v)=0}
  For $u,v \in \spn\{v_x\}$, 
  \linenopax
  \begin{align}\label{eqn:uLapv-as-two-sums}
     \la u, \Lap v \ra_\energy = \sum_{x \in \verts} \cj{\Lap u(x)} \Lap v(x).
  \end{align}
  Furthermore, $\sum_{x \in \verts} \Lap u(x) = 0$ for $u \in \spn\{v_x\}$.
  \begin{proof}
    Let $u \in \spn\{v_x\}$ be given by the finite sum $u = \sum_x \gx_{x} v_{x}$. Since $v_o$ is a constant, we may assume the sum does not include $o$. Then
    \linenopax
    \begin{align}\label{eqn:Lapu(y)-is-the-yth-coord}
      \Lap u(y) 
      = \sum_x \gx_{x} \Lap v_x(y)
      = \sum_x \gx_{x} (\gd_x - \gd_o)(y)
      = \gx_{y}.
    \end{align}
    Now we have 
    \linenopax
    \begin{align*}
      \la u, \Lap u\ra_\energy
      = \sum_{x,y} \cj{\gx_x} \gx_y \la v_x, \Lap v_y \ra_\energy 
      = \sum_{x,y} \cj{\gx_x} \gx_y \la v_x, \gd_y-\gd_o \ra_\energy.
    \end{align*}
    Since it is easy to compute $\la v_x, \gd_y-\gd_o \ra_\energy = \gd_{xy}+1$ (Kronecker's delta), we have 
    \linenopax
    \begin{align}\label{eqn:uLapu-as-two-sums}
      \la u, \Lap u\ra_\energy
       = \sum_{x,y} \cj{\gx_x} \gx_y (\gd_{xy}+1)
      &= \sum_{x} |\gx_x|^2 + \left|\sum_{x} \gx_x\right|^2 \\
      &= \sum_{x} |\Lap u(x)|^2 + \left|\sum_{x} \Lap u(x)\right|^2,
    \end{align}
    by \eqref{eqn:Lapu(y)-is-the-yth-coord}.
    Since $u \in \spn\{v_x\}$, $\Lap u \in \spn\{\gd_x-\gd_o\}$ (see \eqref{eqn:Lapu(y)-is-the-yth-coord}), so that $\la u, \Lap u\ra_\energy < \iy$ and \eqref{eqn:uLapu-as-two-sums} is convergent. Therefore, $\sum_{x} \Lap u(x)$ is absolutely convergent, hence independent of exhaustion. Since 
    \linenopax
    \begin{align*}
      \sum_{x \in \verts} \Lap v_y(x) = 1 - 1 = 0
    \end{align*}
    by Lemma~\ref{thm:vx-is-dipole}, it follows that $\sum_{x} \Lap u(x)=0$, and the second sum in \eqref{eqn:uLapu-as-two-sums} vanishes. Then \eqref{eqn:uLapv-as-two-sums} follows by polarizing.
  \end{proof}
\end{theorem}

\begin{cor}\label{thm:Lap-Hermitian-on-V}\label{thm:ranLapV-is-in-ell2}
  The Laplacian \LapV is Hermitian and even semibounded on $\dom \LapV$ (see Definition~\ref{def:semibounded}) with 
  \linenopax
  \begin{equation}\label{eqn:Lap-semibounded-on-V}
    0 \leq \sum_{x \in \verts} |\Lap u(x)|^2 \leq \la u, \Lap u\ra_\energy < \iy.
  \end{equation}
  \begin{proof}
    For $u,v \in \spn\{v_x\}$, 
    two applications of Lemma~\ref{thm:uLapv-has-no-bdy-term} yield
    \linenopax
    \begin{align*}
      \la \Lap u,v\ra_\energy 
     & = \sum_{x \in \verts}  \cj{\Lap u(x)} \Lap v(x)    
      = \cj{\sum_{x \in \verts} \Lap u(x) \cj{\Lap v(x)}}
      = \cj{\la \Lap v, u\ra_\energy}.
    \end{align*} 
    This property is clearly preserved under closure of the operator.
    
    Now let $u \in \dom \LapV$ and choose $\{u_n\} \ci V$ with $\lim_{n \to \iy}\|u_n-u\|_\energy = \lim_{n \to \iy}\|\Lap u_n - \Lap u\|_\energy = 0$. Then Fatou's lemma \cite[Thm.~I.7.7]{Malliavin} yields 
  \linenopax
  \begin{align}\label{eqn:uLapVv-as-two-sums}
     \sum_{x \in \verts} |\Lap u(x)|^2
     = \sum_{x \in \verts} \lim |\Lap u_n(x)|^2
     \leq \lim_{n \to \iy}\la u_n, \Lap u_n\ra_\energy 
     = \la u, \Lap u\ra_\energy,
  \end{align}
  which gives the central inequality in \eqref{eqn:Lap-semibounded-on-V} and hence semiboundedness.
  \end{proof}
\end{cor}

\begin{remark}\label{rem:when-u-in-HE-is-in-ell2}
  The notation $u \in \ell^1$ means $\sum_{x \in \verts} |u(x)| < \iy$ and the notation $u \in \ell^2$ means $\sum_{x \in \verts} |u(x)|^2 < \iy$. When discussing an element $u$ of \HE, we say $u$ lies in $\ell^2$ if it has a representative which does, i.e., if $u+k \in \ell^2$ for some $k \in \bC$. This constant is clearly necessarily unique on an infinite network, if it exists. 
\end{remark}

The next result is a partial converse to Theorem~\ref{thm:E(u,v)=<u,Lapv>+sum(normals)}.

\begin{lemma}\label{thm:converse-to-E(u,v)=<u,Lapv>}
  If $u, v, \Lap u , \Lap v \in \ell^2$, then $\la u, v \ra_\energy = \sum_{x \in \verts} u(x) \Lap v(x)$, and $u,v \in \dom \energy$.
  \begin{proof}
    If $u, \Lap v \in \ell^2$, then $u\Lap v \in \ell^1$, and the following sum is absolutely convergent:
    \linenopax
    \begin{align*}
      \sum_{x \in \verts} \cj{u}(x) \Lap v(x)
      &= \frac12 \sum_{x \in \verts} \cj{u}(x) \Lap v(x) + \frac12\sum_{y \in \verts} \cj{u}(y) \Lap v(y) \notag \\
      &= \frac12 \sum_{x \in \verts} \sum_{y \nbr x} \cond_{xy} \cj{u}(x) (v(x)-v(y)) - \frac12\sum_{y \in \verts} \sum_{x \nbr y} \cond_{xy} \cj{u}(y) (v(x)-v(y)) \notag \\
      &= \frac12 \sum_{x \in \verts} \sum_{y \nbr x} \cond_{xy} (\cj{u}(x)-\cj{u}(y)) (v(x)-v(y)),
    \end{align*}
    which is \eqref{eqn:def:energy-form}. Absolute convergence justifies the rearrangement in the last equality; the rest is merely algebra. Substituting $u$ in for $v$ in the identity just established, $u\Lap u \in \ell^1$ shows $u \in \dom \energy$, and similarly for $v$.
  \end{proof}
\end{lemma}

\begin{remark}
  All that is required for the computation in the proof of Lemma~\ref{thm:converse-to-E(u,v)=<u,Lapv>} is that $u \Lap v \in \ell^1$, which is certainly implied by $u,\Lap v \in \ell^2$. However, this would not be sufficient to show $u$ or $v$ lies in $\dom \energy$.
\end{remark}

We will see in Theorem~\ref{thm:Fin-vanishes-at-infinity} that if $h \in \Harm$ is nonconstant, then $h+k$ is bounded away from 0 on an infinite set of vertices, for any choice of constant $k$. So the next result should not be surprising. 

\begin{cor}\label{thm:nontrivial-harmonic-fn-is-not-in-L2}
  If $h \in \HE$ is a nontrivial harmonic function, then $h$ cannot lie in $\ell^2$.
  \begin{proof}
    If $h \in \ell^2$, then $\energy(h) = \sum_{x \in \verts} \cj{h}(x) \Lap h(x) = \sum_{x \in \verts} \cj{h}(x) \cdot 0 = 0$ by Lemma~\ref{thm:converse-to-E(u,v)=<u,Lapv>}. But since $h$ is nonconstant, $\energy(h) >0$! \cont
  \end{proof}
\end{cor}

\begin{remark}[Restricting to $\ell^2$ misses the most interesting bit] \label{rem:L2-loses-the-best-part}
  When studying the graph Laplacian, some authors define $\dom \Lap = \{v \in \ell^2 \suth \Lap v \in \ell^2\}$. Our philosophy is that $\dom \energy$ is the most natural context for the study of functions on \verts, and this is motivated in detail in \S\ref{sec:vonNeumann's-embedding-thm}. Some of the most interesting phenomena in $\dom \energy$ are due to the presence of nontrivial harmonic functions, as we show in this section and the examples of \S\ref{sec:tree-networks}--\S\ref{sec:lattice-networks}. Consequently, Corollary~\ref{thm:nontrivial-harmonic-fn-is-not-in-L2} shows why one loses some of the most interesting aspects of the theory by only studying those $v$ which lie in $\ell^2$. Example~\ref{exm:binary-tree:nontrivial-harmonic} illustrates the situation of Corollary~\ref{thm:nontrivial-harmonic-fn-is-not-in-L2} on a tree network. In general, if a at least two connected components of $\Graph \less \{o\}$ are infinite, then $v_x \notin \ell^2$ for vertices $x$ in these components.
\end{remark}

\subsection{The Discrete Gauss-Green formula for networks with vertices of infinite degree}

If there are vertices of infinite degree in the network, then it does not necessary follow that $\spn\{\gd_x\} \ci \spn\{v_x\}$, or that $\spn\{\gd_x\} \ci \MP$. However, we do have the following version of Theorem~\ref{thm:E(u,v)=<u,Lapv>+sum(normals)}. When all vertices have finite degree, Theorem~\ref{thm:E(u,v)=<u,Lapv>-on-Fin} follows from  Theorem~\ref{thm:E(u,v)=<u,Lapv>+sum(normals)} by Lemma~\ref{thm:dx-as-vx}.

\begin{defn}\label{def:LapF}
  Let $\sF := \spn\{\gd_x\}_{x \in \verts}$ denote the vector space of \emph{finite} linear combinations of Dirac masses, and let \LapF be the closure of the Laplacian when taken to have the domain $\sF$.
\end{defn}
\glossary{name={$\sF$},description={the span of the energy kernel, i.e., finite linear combinations of r$v_x$'s},sort=F,format=textbf}
\glossary{name={\LapF},description={the closure of the Laplacian when taken to have the dense domain $\sF$},sort=L,format=textbf}
Note that \sF is a dense domain only when $\Harm=0$, by Corollary~\ref{thm:Diracs-not-dense}. Again, since \Lap agrees with \LapF pointwise, we may suppress reference to the domain for ease of notation. The next result extends Lemma~\ref{thm:E(u,v)=<u,Lapv>} to infinite networks.

\begin{theorem}\label{thm:E(u,v)=<u,Lapv>-on-Fin}
  If $u$ or $v$ lies in $\dom \LapF$, then $\la u, v \ra_\energy = \sum_{x \in \verts} \cj{u(x)} \Lap v(x)$.    
  \begin{proof}
    First, suppose $u \in \dom \LapF$ and choose a sequence $\{u_n\} \ci \spn\{\gd_x\}$ with $\|u_n - u\|_\energy \to 0$. From Lemma~\ref{thm:<delta_x,v>=Lapv(x)}, one has $\la \gd_x, v \ra_\energy = \Lap v(x)$, and hence 
    \[\la u_n,v\ra_\energy = \sum_{x \in \verts} u_n(x) \Lap v(x)\]
    holds for each $n$. Define $M :=  \sup\{\|u_n\|_\energy\}$, and note that $M < \iy$, since this sequence is convergent (to $\|u\|_\energy$). Moreover, $|\la u_n,v \ra_\energy| \leq M\cdot \|v\|_\energy$ by the Schwarz inequality. Since $u_n$ converges pointwise to $u$ in \HE by Lemma~\ref{thm:E-convergence-implies-pointwise-convergence}, this bound will allow us to apply Fatou's Lemma (as stated in \cite[Lemma~7.7]{Malliavin}, for example), as follows:
    \linenopax
    \begin{align*}
      \la u,v \ra_\energy
      &= \lim_{n \to \iy} \la u_n,v \ra_\energy && \text{hypothesis} \\
      &= \lim_{n \to \iy} \sum_{x \in \verts} \cj{u_n(x)} \Lap v(x) && u_n \in \spn\{\gd_x\} \\
      &= \sum_{x \in \verts} \cj{u(x)} \Lap v(x). 
    \end{align*}
    Note that the sum over \verts is absolutely convergent, as required by Definition~\ref{def:exhaustion-of-G}. 
    
    Now suppose that $v \in \dom \LapF$ and observe that this implies $v \in \Fin$ also. By Theorem~\ref{thm:HE=Fin+Harm}, one can decompose $u=f+h$ where $f = \Pfin u$ and $h=\Phar u$, and then
    \linenopax
    \begin{align*}
      \la u,v\ra_\energy = \la f,v\ra_\energy + \la h,v\ra_\energy = \la f,v\ra_\energy,
    \end{align*} 
    since $h$ is orthogonal to $v$. Now apply the previous argument to $\la f,v\ra_\energy$.
  \end{proof}
\end{theorem}

\section{Remarks and references}
\label{sec:Remarks-and-References-3}

For background material and applications of reproducing kernel Hilbert spaces, we suggest the standard references \cite{PaSc72, Aronszajn50} as well as \cite{AD06, AL08, Kai65, MuYaYo, Yoo05, Zh09, BV03, Arv97, Arv76b, ADV09}. Of the cited references for this chapter, some are more specialized. However for prerequisite material (if needed), the reader may find key facts used above on operators in Hilbert space in the books by Dunford-Schwartz \cite{DuSc88}, and Kato \cite{Kat95}. SoardiÕs book \cite{Soardi94} on potential theory is also helpful.

The space of finite-energy functions (often called Dirichlet or Dirichlet-summable functions) on a space is studied widely throughout the literature. In the context of graphs and networks, we recommend the references \cite{Soardi94} (especially Chap.~III) and \cite{Lyons:ProbOnTrees} (especially Chap.~9), and the papers \cite{Kig03, Yamasaki79, Yamasaki77, MuYaYo, CaW92, Wo96, Kayano88, Kayano84, Kayano82}, although we first learned about it from \cite{Kig01} and \cite{Str06}. Throughout most of this literature, the authors study the grounded energy space, and it is the purpose of \S\ref{sec:grounded-energy-space} to clarify the relations between
\begin{align*}
  \energy(u,v)
  \qq\text{and}\qq
  \energy(u,v) + u(o) v(o),
\end{align*}
and hence also between 
\begin{align*}
  \HE = \Fin \oplus \Harm \qq\text{and}\qq \Gdd = \Gddo \oplus \MP_o^-.
\end{align*}

\begin{remark}\label{rem:royden-decomp}
  Theorem~\ref{thm:HE=Fin+Harm}, which shows that $\HE = \Fin \oplus \Harm$, is often called the ``Royden Decomposition'' in the literature, in honour of Royden's analogous result for Riemann surfaces. In many contexts which admit a Laplace operator or suitable analogue, the ensuing decomposition into ``finite'' and ``harmonic'' function spaces is typically called the Royden decomposition, even though the actual contributions of Royden are related only in spirit. 
  
  Note that in \cite[Thm.~3.69]{Soardi94} (and see \cite[\S9.3]{Lyons:ProbOnTrees}), the author uses the grounded inner product and hence the decomposition $\Gdd = \Gddo + \GHD$ is not orthogonal. 
\end{remark}

Of course, the energy form \energy is a Dirichlet form, and the reader seeking more background on the general theory of Dirichlet forms and probability should see \cite{FOT94,BouleauHirsch}, and for Dirichlet spaces in potential theory \cite{Brelot,ConCor}. The best reference for Dirichlet forms in the present context would be Kigami's treatment of resistance forms in \cite{Kig03}. However, one should also see \cite{RS95} and the lovely volume \cite{JostKendallMoscoRocknerSturm}.

For further material on harmonic functions of finite energy, see \cite{CaW92}.

\begin{remark}\label{rem:Kigami-kernel}
  In \cite{Kig03}, Kigami constructs the Green kernel elements $g_B^x(y) = g_B(x,y)$ using potential-theoretic methods. In this context, $B$ is a nonempty finite set which is considered as the boundary of a Dirichlet problem. In the case when $B = \{o\}$, one has $g_b^x = v_x$, where $v_x$ is an element of the energy kernel as defined in Definition~\ref{def:energy-kernel}. However, the construction we give here is entirely in terms of Riesz duality and the Hilbert space structure of \HE, as opposed to discrete potential theory, and was discovered independently.
  
  While \cite{Kig01} and \cite{Kig03} are often thought to pertain specifically to self-similar fractals, there are large parts of both works which are applicable to discrete potential theory more broadly. In particular, many key properties of the resistance metric and its interrelations with the Laplacian and energy form were first worked out in \cite{Kig03}.
\end{remark}

\begin{remark}[Comparison to Kuramochi kernel]
  \label{rem:Kuramochi-kernel}
  After a first version of this book was completed, the authors were referred to \cite{MuYaYo} in which the authors construct a reproducing kernel very similar to ours, which they call the Kuramochi kernel. Indeed, the Kuramochi kernel element $k_x$ corresponds to the representative of $v_x$ which takes the value 0 at $o$. This makes the Kuramochi kernel a reproducing kernel for the space of functions
  \linenopax
  \begin{align*}
    \sD(\Graph;o) := \{u \in \dom \energy \suth u(o)=0\}.
  \end{align*} 
  As advantages of the present approach, we note that our formulation puts the Green kernel in the same space as the reproducing kernel. This will be helpful below; for example, the kernel elements $v_x$ and $f_x = \Phar v_x$ can be decomposed in terms of the Green kernel. See Definition~\ref{def:monopole} and Remark~\ref{rem:ranLap-vs-Fin}. The reader will find that we put the energy kernel to very different uses the Kuramochi kernel.
\end{remark}

%% file: resistance-metrics.tex

\chapter{The resistance metric}
\label{sec:effective-resistance-metric}

\headerquote{The further a mathematical theory is developed, the more harmoniously and uniformly does its construction proceed, and unsuspected relations are disclosed between hitherto separated branches of the science.}{---~D.~Hilbert}

We now introduce the natural notion of distance on $(\Graph, \cond)$: the resistance metric $R$. While not as intuitive as the more common shortest-path metric, it reflects the topology of the graph more accurately and may be more useful for modeling and practical applications. The effective resistance is intimately related to the random walk on $(\Graph,\cond)$, the Laplacian, and the Dirichlet energy form \cite{Kig03,Lyons:ProbOnTrees,LevPerWil08,Soardi94,Kig01,Str06,DoSn84}.

In \S\ref{sec:Resistance-metric-on-finite-networks}, we give several formulations of this metric (Theorem~\ref{thm:effective-resistance-metric}), each with its own advantages. Many of these are familiar from the literature: \eqref{eqn:def:R(x,y)-Lap} from \cite{Pow76b} and \cite[\S8]{Peres99}, \eqref{eqn:def:R(x,y)-energy} from \cite{DoSn84}, \eqref{eqn:def:R(x,y)-diss} from \cite{DoSn84,Pow76b}, \eqref{eqn:def:R(x,y)-R}--\eqref{eqn:def:R(x,y)-S} from \cite{Kig03,Kig01,Str06}. 

In \S\ref{sec:Resistance-metric-on-infinite-networks}, we extend these formulations to infinite networks. Due to the possible presence of nontrivial harmonic functions, some care must be taken when adjusting these formulations. It turns out that there are two canonical extensions of the resistance metric to infinite networks which are distinct precisely when $\Harm \neq 0$ (cf. \cite{Lyons:ProbOnTrees} and the references therein): the ``free'' resistance and the ``wired'' resistance. We are able to clarify and explain the difference in terms of the reproducing kernels for \HE and for \Fin, and also in terms of Dirichlet vs. Neumann boundary conditions; see Remark~\ref{rem:wired-vs-free-as-boundary-conditions}. We also explain the discrepancy in terms of projections in \HE and attempt to relate this to conditioning of the random walk on the network; see \S\ref{sec:probabilistic-interp-of-v_x-and-f_x} and Remark~\ref{rem:Schur-does-network-reduction}. Additionally, we introduce trace resistance and harmonic resistance and relate these to the free and wired resistances. (Note: unlike the others, harmonic resistance is not a metric.) In the limit, the trace resistance coincides with the free resistance.

\section{Resistance metric on finite networks}
\label{sec:Resistance-metric-on-finite-networks}

We make the standing assumption that the network is finite in \S\ref{sec:Resistance-metric-on-finite-networks}. However, the results actually remain true on any network for which $\Harm = 0$.

\begin{defn}\label{def:effective-resistance}
  If one amp of current is inserted into the \ERN at $x$ and withdrawn at $y$, then the \emph{(effective) resistance} $R(x,y)$ is the voltage drop between the vertices $x$ and $y$.
  \glossary{name={$R$},description={resistance metric},sort=r,format=textbf}
\end{defn}

\begin{theorem}\label{thm:effective-resistance-metric}
  The resistance $R(x,y)$ has the following equivalent formulations:
  \linenopax
  \begin{align}
    R(x,y)
    &= \dist_\Lap(x,y) := \{v(x)-v(y) \suth \Lap v = \gd_x-\gd_y\}
    \label{eqn:def:R(x,y)-Lap} \\
    &= \dist_\energy(x,y) := \{\energy(v) \suth \Lap v = \gd_x-\gd_y\} \label{eqn:def:R(x,y)-energy} \\
    &= \dist_\diss(x,y) := \min\{\diss(\curr) \suth \curr \in \Flo(x,y)\} \label{eqn:def:R(x,y)-diss} \\
    &= \dist_R(x,y) := 1/\min
    \nolimits_{v \in \dom\energy}\{\energy(v) \suth v(x)=1, v(y)=0\} \label{eqn:def:R(x,y)-R} \\
    &= \dist_\gk(x,y) := \min\nolimits_{v \in \dom\energy}\{\gk \geq 0 \suth |v(x)-v(y)|^2 \leq \gk \energy(v)\} \label{eqn:def:R(x,y)-S} \\
    &= \dist_s(x,y) := \sup_{v \in \dom \energy}\{|v(x)-v(y)|^2 \suth \energy(v) \leq 1\} \label{eqn:def:R(x,y)-sup}.
  \end{align}
  \begin{proof}
    \eqref{eqn:def:R(x,y)-Lap}$\iff$\eqref{eqn:def:R(x,y)-energy}. We may choose $v$ satisfying $\Lap v = \gd_x-\gd_y$ by Theorem~\ref{thm:Pot-is-nonempty-by-current-flows}.
    Then
    \linenopax
    \begin{align}
      \energy(v)
      = \sum_{z \in \verts} v(z) \Lap v(z)
      = \sum_{z \in \verts} v(z) (\gd_x(z) - \gd_y(z))
      = v(x) - v(y),
    \end{align}
    where first equality is justified by Theorem~\ref{thm:E(u,v)=<u,Lapv>}.

    \eqref{eqn:def:R(x,y)-energy}$\iff$\eqref{eqn:def:R(x,y)-diss}. Note that every $v \in \Pot(x,y)$ corresponds to an element $\curr \in \Flo(x,y)$ via Ohm's Law by Lemma~\ref{thm:potential-induces-current}, and $\energy(v)=\diss(\curr)$ by the same lemma. Also, this current flow is minimal by Theorem~\ref{thm:minimal-flows-are-induced-by-minimizers}.

    \eqref{eqn:def:R(x,y)-energy}$\iff$\eqref{eqn:def:R(x,y)-R}.
    Suppose that $\Lap v = \gd_x - \gd_y$. Since $\energy(v+k)=\energy(v)$ and $\Lap(v+k)=\Lap v$ for any constant $k$, we may adjust $v$ by a constant so that $v(y)=0$. Define
    \linenopax
    \begin{align*}
      u := \frac{v - v(x)}{v(x)-v(y)}
    \end{align*}
    so that $u(x)=0$ and $u(y)=1$. Observe that \eqref{eqn:def:R(x,y)-Lap} gives $\energy(v)=v(x)-v(y)$, whence
    \linenopax
    \begin{align*}
      \energy(u) 
      = \energy(v)/(v(x)-v(y))^2 
      = (v(x)-v(y))^{-1} \geq \min \energy(u).
    \end{align*}
    This shows $\energy(v) \leq [\min \energy(u)]^{-1}$ and hence $\dist_\energy \leq \dist_R$.
 
    For the other inequality, suppose $u$ minimizes $\energy(u)$, subject to $u(x)-u(y)=0$. Then by Theorem~\ref{thm:E(u,v)=<u,Lapv>} and the same variational argument as described in Remark~\ref{rem:variational-approach}, we have
    \linenopax
    \begin{align*}
      \energy(\gr,u) = \sum_{z \in \verts} \gr(z) \Lap u(z) = 0,
    \end{align*}
    for every function \gr for which $\gr(x) = \gr(y) = 0$. It follows that $\Lap u(z) = 0$ for $z \neq x,y$, and hence $\Lap u = \gx \gd_x + \gh \gd_y$. Observe that $\energy(u) = \energy(-u) = \energy(1-u)$, and so the same result follows from minimizing \energy with respect to the conditions $u(y)=1$ and $u(x)=0$. This symmetry forces $\gh = -\gx$ and we have $\Lap u = \gx \gd_x - \gx \gd_y$. Now for $v = \frac1\gx u$ one has $\Lap v = \gd_x - \gd_y$, and so
    \linenopax
    \begin{align*}
      \energy(u) = \gx^2 \energy(v) = \gx^2 (v(x)-v(y)) = \gx(u(x)-u(y)) = \gx,
    \end{align*}
    where the second equality follows by \eqref{eqn:def:R(x,y)-Lap}. Then $\gx = \frac1{\energy(v)} = \energy(u)$, whence $\dist_\energy \geq \dist_R$.

    \eqref{eqn:def:R(x,y)-R}$\iff$\eqref{eqn:def:R(x,y)-S}. Starting with \eqref{eqn:def:R(x,y)-S}, it is clear that
    \linenopax
    \begin{align*}
      \dist_\gk(x,y)
      &= \inf\{\gk \geq 0 \suth \tfrac{|v(x)-v(y)|^2}{\energy(v)} \leq \gk , v \in \dom \energy\} \\
      &= \sup\{\tfrac{|v(x)-v(y)|^2}{\energy(v)}, v \in \dom \energy, v \text{ nonconstant}\}.
    \end{align*}
    Given a nonconstant $v \in \dom \energy$, one can substitute $u:=\frac{v}{|v(x)-v(y)|}$ into the previous line to obtain
    \linenopax
    \begin{align*}
      \dist_\gk(x,y)
      &= \sup\{\tfrac{|u(x)-u(y)|^2\cancel{|v(x)-v(y)|^2}} {\energy(u)\cancel{|v(x)-v(y)|^2}}, v \in \dom \energy, v \text{ nonconstant}\} \\
      &= \sup\{\tfrac1{\energy(u)}, u \in \dom \energy, |u(x)-u(y)|=1\}\\
      &= 1/\inf\{\energy(u), u \in \dom \energy, |u(x)-u(y)|=1\}.
    \end{align*}
    Since we can always add a constant to $u$ and multiply by $\pm1$ without changing the energy, this is equivalent to letting $u$ range over the subset of $\dom\energy$ for which $u(x)=1$ and $u(y)=0$ and we have \eqref{eqn:def:R(x,y)-R}.

    \eqref{eqn:def:R(x,y)-S}$\iff$\eqref{eqn:def:R(x,y)-sup}.
    It is immediate that \eqref{eqn:def:R(x,y)-S} is equivalent to
    \linenopax
    \begin{align*}
      \sup\left\{\frac{|v(x)-v(y)|^2}{\energy(v)} \suth \energy(v) < \iy\right\}.
    \end{align*}
    For any $v \in \dom \energy$, define $w:= v/\sqrt{\energy(v)}$ so that $|w(x)-w(y)|^2 = |v(x)-v(y)|^2 \energy(v)^{-1/2}$ with $\energy(w) = 1$. Clearly then $|w(x)-w(y)|^2 \leq \dist_s(x,y)$. The other inequality is similar.
  \end{proof}
\end{theorem}

The equivalence of \eqref{eqn:def:R(x,y)-diss} and \eqref{eqn:def:R(x,y)-Lap} is shown elsewhere (e.g., see \cite[\S{II}]{Pow76b}) but the reader will find some gaps, so we have included a complete version of this proof for completeness. The terminology ``effective resistance metric'' is common in the literature (see, e.g., \cite{Kig01} and \cite{Str06}), where it is usually given in the form \eqref{eqn:def:R(x,y)-R}. The formulation \eqref{eqn:def:R(x,y)-S} will be helpful for obtaining certain inequalities in the sequel. It is also clear that $\dist_s$ of \eqref{eqn:def:R(x,y)-sup} is the norm of the operator $L_{xy}$ defined by $L_{xy} u := u(x) - u(y)$, see Lemma~\ref{thm:L_x-is-bounded} and Theorem~\ref{thm:free-resistance}. 

\begin{remark}\label{rem:minimal-energy-is-obtained}
  Taking the minimum (rather than the infimum) in \eqref{eqn:def:R(x,y)-diss}, etc, is justified by Theorem~\ref{thm:energy-obtains-min}. The same argument implies that the energy kernel on \Graph is uniquely determined.
\end{remark}

\begin{figure}
    \centering
    \includegraphics{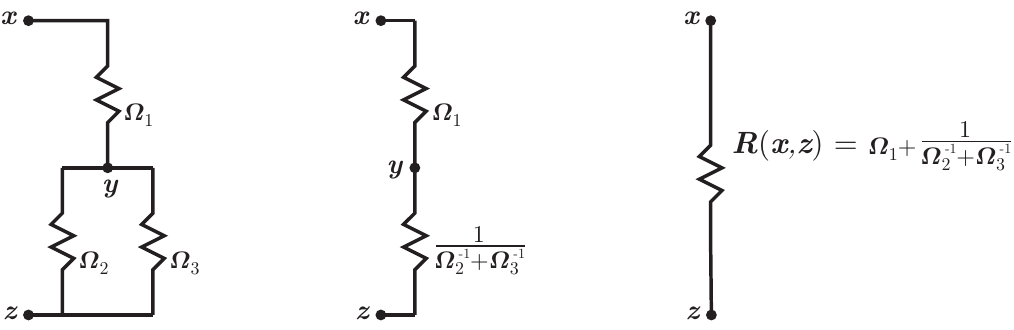}
    \caption{\captionsize Effective resistance as network reduction to a trivial network. This basic example uses parallel reduction followed by series reduction; see Remark~\ref{rem:Resistance-distance-via-network-reduction}.}
    \label{fig:network-reduction}
  \end{figure}

\begin{remark}[Resistance distance via network reduction]
  \label{rem:Resistance-distance-via-network-reduction}
  Let $H$ be a (connected) planar subnetwork of a finite network $G$ and pick any $x,y \in H$. Then $H$ may be reduced to a trivial network consisting only of these two vertices and a single edge between them via the use of three basic transformations: (i) series reduction, (ii) parallel reduction, and (iii) the $\nabla$-\textsf{Y} transform. Each of these transformations preserves the resistance properties of the subnetwork, that is, for $x,y \in G \less H$, $R(x,y)$ remains unchanged when these transformations are applied to $H$. The effective resistance between $x$ and $y$ may be interpreted as the resistance of the resulting single edge. An elementary example is shown in Figure~\ref{fig:network-reduction}. A more sophisticated technique of network reduction is given by the Schur complement construction defined in Remark~\ref{rem:Resistance-distance-via-Schur-complement}.
\end{remark}

The following result is not new (see, e.g. \cite[\S2.3]{Kig01}), but the proof given here is substantially simpler than most others found in the literature.

\begin{lemma}\label{thm:R-is-a-metric}
  $R$ is a metric.
  \begin{proof}
    Symmetry and positive definiteness are immediate from \eqref{eqn:def:R(x,y)-energy}, we use \eqref{eqn:def:R(x,y)-Lap} to check the triangle inequality. Let $v_1 \in \Pot(x,y)$ and $v_2 \in \Pot(y,z)$. By superposition, $v_3 := v_1 + v_2$ is in $\Pot(x,z)$. For $s \nbr t$, it is clear that $v_3(s)-v_3(t) = v_1(s)-v_1(t) + v_2(s)-v_2(t)$. By summing along any path from $x$ to $z$, one sees that this remains true for $s \notnbr t$, whence
    \linenopax
    \begin{align*}
      R(x,z)
      = v_3(x)-v_3(z)
      &= v_1(x)-v_1(z) + v_2(x)-v_2(z) \\
      &\leq v_1(x)-v_1(y) + v_2(y)-v_2(z)
      = R(x,y) + R(y,z),
    \end{align*}
    where the inequality follows from Proposition~\ref{thm:v-min-at-a,max-at-w}.
  \end{proof}
\end{lemma}

\section{Resistance metric on infinite networks}
\label{sec:Resistance-metric-on-infinite-networks}

There are difficulties with extending the results of the previous section to infinite networks. The existence of nonconstant harmonic functions $h \in \dom \energy$ implies the nonuniqueness of solutions to $\Lap u = f$, and hence \eqref{eqn:def:R(x,y)-Lap}--\eqref{eqn:def:R(x,y)-diss} are no longer well-defined. Two natural choices for extension lead to the free resistance $R^F$ and the wired resistance $R^W$. In this section, we attempt to explain the somewhat surprising phenomenon that one may have $R^W(x,y) < R^F(x,y)$. 
\begin{enumerate}
  \item In Theorem~\ref{thm:free-resistance}, we show how $R^F$ corresponds to choosing solutions to $\Lap u = \gd_x - \gd_y$ from the energy kernel, and how it corresponds to currents which are decomposable in terms of paths. The latter leads to a probabilistic interpretation which provides for a relation to the trace of the resistance discussed in \S\ref{sec:trace-resistance}.
  \item In Theorem~\ref{thm:wired-resistance}, we show how $R^W$ corresponds to projection to \Fin. Since this corresponds to minimization of energy, it is naturally related to capacity. 
\end{enumerate}
See also Remark~\ref{sec:probabilistic-interp-of-v_x-and-f_x}. Both of these notions are methods of specifying a \emph{unique} solutions to $\Lap u = f$ in some way. The disparity between $R^F$ and $R^W$ is thus explained in terms of boundary conditions on \Lap as an unbounded self-adjoint operator on \HE in Remark~\ref{rem:wired-vs-free-as-boundary-conditions}. For an alternative approach, see \cite{Kig03}, where the author uses a potential-theoretic formulation (axiomatic harmonic analysis) to explain the discrepancy between $R^F$ and $R^W$ in terms of domains. (This will also be apparent from our approach, see Remark~\ref{rem:comparison-to-resistance-forms}.)

To compute effective resistance in an infinite network, we will need three notions of subnetwork: free, wired, and trace. Strictly speaking, these may not actually be subnetworks of the original graph; they are networks associated to a full subnetwork. 
  Throughout this section, we use $H$ to denote a finite full subnetwork of \Graph, $\verts[H]$ to denote its vertex set, and $H^F$, $H^W$, and $H^\trc$ to denote the free, wired, and trace networks associated to $H$ (these terms are defined in other sections below). 

\begin{defn}\label{def:relative-resistance}
  If $H$ is a subnetwork of \Graph which contains $x$ and $y$, define $R_H(x,y)$ to be the resistance distance from $x$ to $y$ as computed within $H$. In other words, compute $R_H(x,y)$ by any of the equivalent formulas of Theorem~\ref{thm:effective-resistance-metric}, but extremizing over only those functions whose support is contained in $H$.
\end{defn}
  \glossary{name={$R_H(x,y)$},description={relative resistance as computed with respect to the subnetwork $H$},sort=RH,format=textbf}

We will always use the notation $\{G_k\}_{k=1}^\iy$ to denote an \emph{exhaustion} of the infinite network \Graph. Recall from Definition~\ref{def:exhaustion-of-G} that this means each $G_k$ is a finite connected subnetwork of \Graph, $G_k \ci G_{k+1}$, and $\verts = \bigcup \verts[G_k]$. Since $x$ and $y$ are contained in all but finitely many $G_k$, we may always assume that $x,y \in G_k$. Also, we assume in this section that the subnetworks are full --- this is not necessary, but simplifies the discussion and causes no loss of generality.

\begin{defn}\label{def:full-subnetwork}
  Let $\verts[H] \ci \verts$. Then the \emph{full subnetwork} on \verts[H] has all the edges of \Graph for which both endpoints lie in \verts[H], with the same conductances.
\end{defn}

\section{Free resistance}

\begin{defn}\label{def:free-resistance}
  For any subset $\verts[H] \ci \verts$, the \emph{free subnetwork} $H^F$ is just the full subnetwork $H$. That is, all edges of \Graph with endpoints in \verts[H] are edges of $H$, with the same conductances. 
  Let $R_{H^F}(x,y)$ denote the effective resistance between $x$ and $y$ as computed in $H=H^F$, as in Definition~\ref{def:relative-resistance}. 
  The \emph{free resistance} between $x$ and $y$ is then defined to be
  \linenopax
  \begin{align}\label{eqn:def:free-resistance}
    R^F(x,y) := \lim_{k \to \iy} R_{G_k^F}(x,y),
  \end{align}
  where $\{G_k\}$ is any exhaustion of \Graph.
\end{defn}
  \glossary{name={$R^F(x,y)$},description={free resistance metric},sort=RF,format=textbf}

\begin{remark}\label{rem:free-resistance}
  The name ``free'' comes from the fact that this formulation is free of any boundary conditions or considerations of the complement of $H$, in contrast to the wired and trace formulations of the next two sections. See \cite[\S9]{Lyons:ProbOnTrees} for further justification of this nomenclature.

  One can see that $R_{H^F}(x,y)$ has the drawback of ignoring the conductivity provided by all paths from $x$ to $y$ that pass through the complement of $H$. This provides some motivation for the wired and trace approaches below.
\end{remark}

\begin{defn}\label{def:L_xy}
  Fix $x, y \in \Graph$ and define the linear operator $L_{xy}$ on \HE by $L_{xy}v := v(x)-v(y)$. 
\end{defn}

\begin{remark}
  \label{rem:R(x,y) = |L{xy}|}
  Theorem~\ref{thm:free-resistance} is the free extension of Theorem~\ref{thm:effective-resistance-metric} to infinite networks; it shows that $R(x,y) = \|L_{xy}\|$ and that $R(x,o)$ is the best possible constant $k=k_x$ in Lemma~\ref{thm:L_x-is-bounded}. In the proof, we use the notation $\charfn{\cpath}$ for a current which is the characteristic function of a path, that is, a current which takes value 1 on every edge of $\cpath \in \Paths(x,y)$ and 0 on all other edges. Then $\curr = \sum \gx_\cpath \charfn{\cpath}$ indicates that \curr decomposes as a sum of currents supported on paths in \Graph.
\end{remark}
  \glossary{name={$\curr = \sum \gx_\cpath \charfn{\cpath}$},description={current decomposed as a sum over paths},sort=I,format=textbf}
  \glossary{name={\cpath},description={path},sort=G}

\begin{theorem}\label{thm:free-resistance}
  For an infinite network \Graph, the free resistance $R^F(x,y)$ has the following equivalent formulations:
  \linenopax
  \begin{align}
    R^F(x,y)
    &= (v_x(x) - v_x(y)) - (v_y(x) - v_y(y)) \label{eqn:def:R^F(x,y)-Lap} \\
    &= \energy(v_x-v_y) \label{eqn:def:R^F(x,y)-energy} \\
    &= \min \{\diss(\curr) \suth \curr \in \Flo(x,y) \text{ and } \curr = \textstyle \sum_{\cpath \in \Paths(x,y)} \gx_\cpath \charfn{\cpath}\} \label{eqn:def:R^F(x,y)-diss} \\
    &= 1 / \min\{\energy(u) \suth\vstr[2.5] u(x)=1, u(y)=0, u \in \dom \energy\} \label{eqn:def:R^F(x,y)-R} \\ 
    &= \inf\{\gk \geq 0 \suth |v(x)-v(y)|^2 \leq \gk \energy(v), v \in \dom \energy\} \label{eqn:def:R^F(x,y)-S} \\
    &= \sup\{|v(x)-v(y)|^2 \suth \energy(v) \leq 1, v \in \dom \energy\} \label{eqn:def:R^F(x,y)-sup} 
  \end{align}
  \begin{proof}
    To see that \eqref{eqn:def:R^F(x,y)-energy} is equivalent to \eqref{eqn:def:free-resistance}, fix any exhaustion of \Graph and note that 
    \linenopax
    \begin{align*}
      \energy(v_x-v_y) = \lim_{k \to \iy} \frac12 \sum_{s,t \in G_k} \cond_{st}((v_x-v_y)(s) - (v_x-v_y)(t))^2 = \lim_{k \to \iy} R_{G_k^F}(x,y),
    \end{align*}
    where the latter equality is from Theorem~\ref{thm:effective-resistance-metric}. Then for the equivalence of formulas \eqref{eqn:def:R^F(x,y)-Lap} and \eqref{eqn:def:R^F(x,y)-energy}, simply compute
    \linenopax
    \begin{align*}
      \energy(v_x-v_y)
      &= \la v_x-v_y, v_x-v_y \ra_\energy
       = \la v_x, v_x \ra_\energy - 2\la v_x, v_y \ra_\energy + \la v_y, v_y \ra_\energy
    \end{align*}
    and use the fact that $v_x$ is \bR-valued; cf.~\cite[Lemma~2.22]{DGG}.
    
    To see \eqref{eqn:def:R^F(x,y)-diss} is equivalent to \eqref{eqn:def:free-resistance}, fix any exhaustion of \Graph and define
    \linenopax
    \begin{align*}
      \Flo(x,y)\evald{H} := \{\curr \in \Flo(x,y) \suth \curr = \textstyle\sum\nolimits_{\cpath \ci H} \gx_\cpath \charfn{\cpath}\}.
    \end{align*}
    From \eqref{eqn:def:R(x,y)-diss}, it is clearly true for each $G_k$ that
    \linenopax
    \begin{align*}
      R_{G_k^F}(x,y) = \min \{\diss(\curr) \suth \curr \in \Flo(x,y) \text{ and } \curr = \textstyle \sum_{\cpath \ci G_k} \gx_\cpath \charfn{\cpath} \}.
    \end{align*}
    Since $\Flo(x,y)\evald{\Graph} = \bigcup_k \Flo(x,y)\evald{G_k}$, formula \eqref{eqn:def:R^F(x,y)-diss} follows. Note that \diss is a quadratic form on the closed convex set $\Flo(x,y)\evald{\Graph}$ and hence it attains its minimum.
    
    The equivalence of \eqref{eqn:def:R^F(x,y)-R} and \eqref{eqn:def:R^F(x,y)-sup} is \cite[Thm.~2.3.4]{Kig01}. 
    
    As for \eqref{eqn:def:R^F(x,y)-S} and \eqref{eqn:def:R^F(x,y)-sup}, they are both clearly equal to $\|L_{xy}\|$ (as described in Remark~\ref{rem:R(x,y) = |L{xy}|}) by the definition of operator norm; see \cite[\S5.3]{Rud87}, for example. To show that these are equivalent to $R^F$ as defined in \eqref{eqn:def:free-resistance}, define a subspace of \HE consisting of those voltages whose induced currents are supported in a finite subnetwork $H$ by    
    \linenopax
    \begin{equation}\label{eqn:def:HE-F-restr}
      \HE\evald[F]{H} = \{u \in \dom \energy \suth u(x)-u(y)=0 \text{ unless } x,y \in H\}. 
    \end{equation}
    This is a closed subspace, as it is the intersection of the kernels of a collection of continuous linear functionals $\|L_{st}\|$, and so we can let $Q_k$ be the projection to this subspace. Then it is clear that $Q_k \leq Q_{k+1}$ and that $\lim_{k \to \iy} \|u - Q_k u\|_\energy = 0$ for all $u \in \HE$, so
    \linenopax
    \begin{align}\label{eqn:defn-Q_k}
      R_{G_k^F}(x,y) 
      &= \| L_{xy} \|_{\HE|_{G_k} \to \bC }
       = \| L_{xy} Q_k\|,
    \end{align}
where the first equality follows from \eqref{eqn:def:R(x,y)-S} (recall that $G_k$ is finite) and therefore
    \linenopax
    \begin{align*}
      R^F(x,y) = \lim_{k \to \iy} R_{G_k^F}(x,y) 
      &= \lim_{k \to \iy} \|L_{xy} Q_k\| 
      = \left\| \lim_{k \to \iy} L_{xy} Q_k \right\| 
      = \|L_{xy}\|.
      \qedhere
    \end{align*}
  \end{proof}
\end{theorem}

  In view of the previous result, the free case corresponds to consideration of only those voltage functions whose induced current can be decomposed as a sum of currents supported on paths in \Graph. The wired case considered in the next section corresponds to considering all voltages functions whose induced current flow satisfies Kirchhoff's law in the form \eqref{eqn:Kirchhoff-nonhomog}; this is clear from comparison of \eqref{eqn:def:R^F(x,y)-diss} to \eqref{eqn:def:R^W(x,y)-diss}. See also Remark~\ref{rem:wired-vs-free-as-Kirchhoff-conditions+}.

  Formula \eqref{eqn:def:R^F(x,y)-Lap} turns out to be useful for explicit computations. Explicit formulas for the effective resistance metric on \bZd are obtained from \eqref{eqn:def:R^F(x,y)-Lap} in \cite[\S14.2]{OTERN}; compare to \cite[\S{V.2}]{Soardi94}. 

\begin{remark}\label{rem:martingales-in-Hilbert-space}
  In Theorem~\ref{thm:free-resistance}, the proofs that $R^F$ is given by \eqref{eqn:def:R^F(x,y)-diss} or \eqref{eqn:def:R^F(x,y)-S} stem from essentially the same underlying martingale argument. In a Hilbert space, a martingale is an increasing sequence of projections $\{Q_k\}$ with the martingale property $Q_k = Q_k Q_{k+1}$. Recall that conditional expectation is a projection. In this context, Doob's theorem \cite{Doob53} then states that if $\{f_k\} \ci \sH$ is such that $f_k = Q_k f_j$ for any $j \geq k$, then the following are equivalent:
\begin{enumerate}[(i)]
  \item there is a $f \in \sH$ such that $f_k = Q_k f$ for all $k$
  \item $\sup_k \|f_k\| < \iy$. 
\end{enumerate}
  For \eqref{eqn:def:R^F(x,y)-diss}, we are actually projecting to subspaces of \HD, the Hilbert space of currents introduced as the \emph{dissipation space} in \S\ref{sec:H_energy-and-H_diss}.
  In \cite[\S9.1]{Lyons:ProbOnTrees}, the free resistance $R^F(x,y)$ is defined directly  via this approach (and similarly for $R^W(x,y)$). 

  In view of the previous result, the free case corresponds to consideration of only those voltage functions whose induced current can be decomposed as a sum of currents supported on paths in \Graph. The wired case considered in the next section corresponds to considering all voltages functions whose induced current flow satisfies Kirchhoff's law \eqref{def:Kirchhoff's-law}; this is clear from comparison of \eqref{eqn:def:R^F(x,y)-diss} to \eqref{eqn:def:R^W(x,y)-diss}. See also Remark~\ref{rem:wired-vs-free-as-Kirchhoff-conditions+}.

  Formula \eqref{eqn:def:R^F(x,y)-Lap} turns out to be useful for explicit computations; we use it to obtain explicit formulas for the effective resistance metric on \bZd in Theorem~\ref{thm:R(x,y)-on-Zd}. 
\end{remark}

The following result is also a special case of \cite[Thm.~2.3.4]{Kig01}.

\begin{prop}\label{thm:free-resistance-is-a-metric}
  $R^F(x,y)$ is a metric.
  \begin{proof}
    One has $R_{G_k^F}(x,z) \leq R_{G_k^F}(x,y) + R_{G_k^F}(y,z)$ for any $k$, so take the limit.
  \end{proof}
\end{prop}

From Theorem~\ref{thm:free-resistance}, it is clear that the triangle inequality also has the formulation
\linenopax
\begin{equation*}
  \energy(v_x-v_z) \leq \energy(v_x-v_y) + \energy(v_y-v_z), \q\forall x,y,z \in \verts,
\end{equation*}
which is easily shown to be equivalent to 
\linenopax
\begin{equation*}
  v_x(z) + v_y(z) \leq v_z(z) + v_x(y), \q\forall x,y,z \in \verts,
\end{equation*}
using the convention $v_x(o)=0$.

The next result is immediate from \eqref{eqn:def:R^F(x,y)-S} and appears also in \cite{Kig03}.

\begin{cor}\label{thm:vx-is-Lipschitz}
  Every function in \HE is H\"{o}lder continuous with exponent $\frac12$.
\end{cor}

It is known from \cite{Nelson64} that the Gaussian measure of Brownian motion is supported on the space of such functions and this will be useful later; cf.~Remark~\ref{rem:boundary-of-G-from-Minlos} and the beginning of \S\ref{sec:Gel'fand-triples-and-duality}. It is somewhat subtle to determine if $R(x,\cdot)$ is in \HE.

\section{Wired resistance}

\begin{defn}\label{def:wired-resistance}
  Given a finite full subnetwork $H$ of \Graph, define the wired subnetwork $H^W$ by identifying all vertices in $\verts \less \verts[H]$ to a single, new vertex labeled \iy. Thus, the vertex set of $H^W$ is $\verts[H] \cup \{\iy\}$, and the edge set of $H^W$ includes all the edges of $H$, with the same conductances. However, if $x \in \verts[H]$ has a neighbour $y \in \verts \less \verts[H]$, then $H^W$ also includes an edge from $x$ to \iy with conductance
  \linenopax
  \begin{align}\label{eqn:cond-to-iy}
    \cond_{x\iy} := \sum_{y \nbr x, \, y \in H^\complm} \negsp[13]\cond_{xy}.
  \end{align}

  Let $R_{H^W}(x,y)$ denote the effective resistance between $x$ and $y$ as computed in $H^W$, as in Definition~\ref{def:relative-resistance}.
  The \emph{wired resistance} is then defined to be
  \linenopax
  \begin{align}\label{eqn:def:wired-resistance}
    R^W(x,y) := \lim_{k \to \iy} R_{G_k^W}(x,y),
  \end{align}
  where $\{G_k\}$ is any exhaustion of \Graph.
\end{defn}
  \glossary{name={$R^W(x,y)$},description={wired resistance metric},sort=RW,format=textbf}
  \glossary{name={$R_H(x,y)$},description={relative resistance as computed with respect to the subnetwork $H$},sort=RH}

\begin{remark}\label{rem:wired-resistance}
  The wired subnetwork is equivalently obtained by ``shorting together'' all vertices of $H^\complm$, and hence it follows from Rayleigh's monotonicity principle that $R^W(x,y) \leq R^F(x,y)$; cf. \cite[\S1.4]{DoSn84} or \cite[\S2.4]{Lyons:ProbOnTrees}. The reader will see by comparison to Theorem~\ref{thm:wired-resistance} that the wired resistance $R^W$ is also the effective resistance associated to the resistance form $(\energy,\Fin)$ of \cite{Kig03}; see Remark~\ref{rem:comparison-to-resistance-forms}. However, the wired resistance is \emph{not} related to the ``shorted resistance form'' of \cite[\S3]{Kig03} (see Prop.~3.6 in particular).

  The justification for \eqref{eqn:cond-to-iy} is that the identification of vertices in $G_k^\complm$ may result in parallel edges. Then \eqref{eqn:cond-to-iy} corresponds to replacing these parallel edges by a single edge according to the usual formula for resistors in parallel.
\end{remark}

\begin{theorem}\label{thm:wired-resistance}
  The wired resistance may be computed by any of the following equivalent formulations:
  \linenopax
  \begin{align}
    R^W(x,y)
    &= \min_v \{v(x)-v(y) \suth \Lap v = \gd_x-\gd_y, v \in \dom\energy\}
    \label{eqn:def:R^W(x,y)-Lap} \\
    &= \min_v \{\energy(v) \suth \Lap v = \gd_x-\gd_y, v \in \dom\energy\} \label{eqn:def:R^W(x,y)-energy} \\
    &= \min_\curr \{\diss(\curr) \suth \curr \in \Flo(x,y), \diss(\curr) < \iy\} \label{eqn:def:R^W(x,y)-diss} \\
    &= 1/\min\{\energy(u) \suth u(x)=1, u(y)=0, u \in \Fin\} \label{eqn:def:R^W(x,y)-R} \\
    &= \inf\{\gk \geq 0 \suth |v(x)-v(y)|^2 \leq \gk \energy(v), v \in \Fin\} \label{eqn:def:R^W(x,y)-S} \\
    &= \sup\{|v(x)-v(y)|^2 \suth \energy(v) \leq 1, v \in \Fin\} \label{eqn:def:R^W(x,y)-sup} 
  \end{align}
  \begin{proof}
    Since \eqref{eqn:def:R^W(x,y)-S} and \eqref{eqn:def:R^W(x,y)-sup} are both clearly equivalent to the norm of $L_{xy}:\Fin \to \bC$ (where again $L_{xy}u - u(x)-u(y)$ as in Remark~\ref{rem:R(x,y) = |L{xy}|}), we begin by equating them to \eqref{eqn:def:wired-resistance}. From Definition~\ref{def:Fin}, we see that 
    \linenopax
    \begin{align}\label{eqn:def:HE-W-restr}
      \HE\evald[W]{H}
      := \{u \in \HE \suth \spt u \ci H\}
    \end{align}
    is a closed subspace of \HE. Let $Q_k$ be the projection to this subspace. Then it is clear that $Q_k \leq Q_{k+1}$ and that $\lim_{k \to \iy} \| P_{\Fin}u - Q_k u\|_\energy = 0$ for all $u \in \HE$. Each function $u$ on $H^W$ corresponds to a function $\tilde u$ on \Graph whose support is contained in $H$; simply define 
    \linenopax
    \begin{align*}
    \tilde u(x) = 
    \begin{cases} 
       u(x), & x \in H, \\ 
       u(\iy_H), & x \notin H.
    \end{cases}
    \end{align*}
    It is clear that this correspondence is bijective, and that
    \linenopax
    \begin{align*}
      R_{G_k^W}(x,y) 
      &= \| L_{xy} \|_{\HE|_{G_k}^W \to \bC }
       = \| L_{xy} Q_k\|,
    \end{align*}
    where the first equality follows from \eqref{eqn:def:R(x,y)-S} (recall that $G_k$ is finite) and therefore
    \linenopax
    \begin{align*}
      R^W(x,y) = \lim_{k \to \iy} R_{G_k^W}(x,y) 
      &= \lim_{k \to \iy} \|L_{xy} Q_k\| 
      = \|L_{xy} P_{\Fin}\|,
    \end{align*}
    which is equivalent to \eqref{eqn:def:R^W(x,y)-S}.
    
    To see \eqref{eqn:def:R^W(x,y)-Lap} is equivalent to \eqref{eqn:def:R^W(x,y)-energy}, note that the minimal energy solution to $\Lap u = \gd_x - \gd_y$ lies in \Fin, since any two solutions must differ by a harmonic function. Let $u$ be a solution to $\Lap u = \gd_x - \gd_y$ and define $f = \Pfin u$. Then $f \in \Fin$ and $\Lap f = \gd_x - \gd_y$ implies 
    \linenopax
    \begin{align}
      \|f\|_\energy^2 
      = \sum_{z \in \verts} f(z) \Lap f(z) 
      = \sum_{z \in \verts} f(z) (\gd_x-\gd_y)(z)
      = f(x) - f(y).
      \label{eqn:thm:wired-resistance:eqn1}
    \end{align}

    To see \eqref{eqn:def:R^W(x,y)-Lap}$\leq$ \eqref{eqn:def:R^W(x,y)-S}, let \gk be the optimal constant from \eqref{eqn:def:R^W(x,y)-S}. If $u \in \Fin$ is the unique solution to $\Lap u = \gd_x - \gd_y$, then 
    \linenopax
    \begin{align*}
      \gk = \sup_{u \in \Fin} \left\{\frac{|u(x)-u(y)|^2}{\energy(u)}\right\}
      \geq \frac{|u(x)-u(y)|^2}{\energy(u)} 
      = u(x)-u(y),               
    \end{align*}
    where the last equality follows from $\energy(u) = u(x) - u(y)$, by the same computation as in \eqref{eqn:thm:wired-resistance:eqn1}.
    For the reverse inequality, note that with $L_{xy}$ as just above,
    \linenopax
    \begin{align*}
      \frac{|u(x)-u(y)|^2}{\energy(u)}
      = \left| L_{xy} \left(\tfrac{u}{\energy(u)^{1/2}}\right)\right|^2
      = \left| \left\la v_x-v_y, \tfrac{u}{\energy(u)^{1/2}}\right\ra_\energy\right|^2,               
    \end{align*}
    for any $u \in \Fin$. Note that Lemma~\ref{thm:wired-resistance-is-a-metric} allows one to replace $v_x$ by $f_x = \Pfin v_x$, whence 
    \linenopax
    \begin{align*}
      \frac{|u(x)-u(y)|^2}{\energy(u)}
      &\leq \energy(f_x-f_y) \energy\left(\tfrac{u}{\energy(u)^{1/2}}\right) 
      = \energy(f_x-f_y) 
    \end{align*}
    by Cauchy-Schwarz. The infimum of the left-hand side over nonconstant functions $u \in \Fin$ gives the optimal \gk in \eqref{eqn:def:R^W(x,y)-S}, and thus shows that \eqref{eqn:def:R^W(x,y)-S} $\leq$ \eqref{eqn:def:R^W(x,y)-energy}. 

    To see \eqref{eqn:def:R^W(x,y)-energy} is equivalent to \eqref{eqn:def:R^W(x,y)-diss}, recall that \curr minimizes \diss over $\Flo(x,y)$ if and only if $\curr = \drp u$ for $u$ which minimizes \energy over $\{v \in \dom \energy \suth \Lap v = \gd_x - \gd_y\}$; see \cite[Thm.~3.26]{OTERN}, for example. Apply this to $\curr = \drp f$, where $f = \Pfin u$ is the minimal energy solution to $\Lap u = \gd_x - \gd_y$. 
    
    The equivalence of \eqref{eqn:def:R^W(x,y)-R} and \eqref{eqn:def:R^W(x,y)-sup} is directly parallel to the finite case and may also be obtained from \cite[Thm.~2.3.4]{Kig01}.
  \end{proof}
\end{theorem}

\begin{remark}[$R^F$ vs. $R^W$ explained in terms of boundary conditions on \Lap]
  \label{rem:wired-vs-free-as-boundary-conditions}
  Observe that both spaces  
  \linenopax
  \begin{align*} 
    \HE\evald[F]{H} = \{u \in \HE \suth u(x)-u(y)=0 \text{ unless } x,y \in H\}
  \end{align*}
  and
  \linenopax
  \begin{align*} 
    \HE\evald[W]{H} = \{u \in \HE \suth \spt u \ci H\}
  \end{align*}
  consist of functions which have no energy outside of $H$. The difference is that if the complement of $H$ consists of several connected components, then $u \in \HE|_H^F$ may take a different constant value on each one; this is not allowed for elements of $\HE|^W_H$. Therefore, $\HE|_H^F$ corresponds to Neumann boundary conditions and $\HE|_H^W$ corresponds to Dirichlet boundary conditions. That is, from the proofs of Theorem~\ref{thm:free-resistance} and Theorem~\ref{thm:wired-resistance}, we see
  \begin{enumerate}
    \item $R_{H^F}(x,y) = u(x)-u(y)$ where $u$ is the solution to $\Lap u = \gd_x - \gd_y$ with Neumann boundary conditions on $H^\complm$, and
    \item $R_{H^W}(x,y) = u(x)-u(y)$ where $u$ is the solution to $\Lap u = \gd_x - \gd_y$ under Dirichlet boundary conditions on $H^\complm$.
  \end{enumerate}
  
\end{remark}

\begin{remark}\label{rem:wired-vs-free-as-Kirchhoff-conditions+}
  While the wired subnetwork takes into account the conductivity due to all paths from $x$ to $y$ (see Remark~\ref{rem:free-resistance}), it is overzealous in that it may also include paths from $x$ to $y$ that do not correspond to any path in \Graph (see Remark~\ref{rem:martingales-in-Hilbert-space}). On an infinite network, this leads to current flows in which some of the current travels from $x$ to \iy, and then from \iy to $y$. Consider the example of \cite{Morris2003}: let \Graph be \bZ with $\cond_{n,n+1} =1$ for each $n$. Then define $J$ by
  \linenopax
  \begin{align*}
    J(n,n-1) =
    \begin{cases}
      1, & n \neq 1 \\
      0, & n=1.
    \end{cases}
  \end{align*}
  If a unit current flow from 0 to 1 is defined to be a current satisfying $\sum_{y \nbr x} \curr(x,y) = \gd_x - \gd_y$, then $J$ is such a flow which ``passes through \iy'' (of course, $J$ certainly not of finite energy).
\end{remark}

The proof of the next result follows from the finite case, exactly as in Theorem~\ref{thm:free-resistance-is-a-metric}.
\begin{theorem}\label{thm:wired-resistance-is-a-metric}
  $R^W(x,y)$ is a metric.
\end{theorem}

\section{Harmonic resistance}
\label{sec:Harmonic-resistance}

\begin{defn}\label{def:harmonic-resistance}
  For an infinite network $(G,\cond)$ define the \emph{harmonic resistance} between $x$ and $y$ by
  \linenopax
  \begin{align}\label{eqn:def:harmonic-resistance}
    R^\hrm(x,y) := R^F(x,y) - R^W(x,y).
  \end{align}
\end{defn}

The next result is immediate upon comparing Theorem~\ref{thm:free-resistance} and Theorem~\ref{thm:wired-resistance}.

\begin{theorem}\label{thm:harmonic-resistance}
  With $h_x = \Phar v_x$ as in Remark~\ref{rem:3-repkernels}, the harmonic resistance is equal to
  \linenopax
  \begin{align}
    R^\hrm(x,y)
    &= (h_x(x) - h_x(y)) - (h_y(x) - h_y(y)) \label{eqn:def:R^ha(x,y)-Lap} \\
    &= \energy(h_x-h_y) \label{eqn:def:R^ha(x,y)-energy} \\
    &= \tfrac1{\min\{\energy(v) \suth v(x)=1, v(y)=0\}}  - \tfrac1{\min\{\energy(f) \suth f(x)=1, f(y)=0, f \in \Fin\}} \label{eqn:def:R^ha(x,y)-R} \\
    &= \inf\{\gk \geq 0 \suth |h(x)-h(y)|^2 \leq \gk \energy(h), h \in \Harm\} \label{eqn:def:R^ha(x,y)-S} \\
    &= \sup\{|h(x)-h(y)|^2 \suth \energy(h) \leq 1, h \in \Harm\} \label{eqn:def:R^ha(x,y)-sup} 
  \end{align}
\end{theorem}

  \version{}{\marginpar{Idea for precise definition: $g_o$ is constant on this portion.}}
\begin{remark}\label{rem:harmonic-resistance-not-a-metric}\label{rem:Rharm-has-no-form}
  Note that $R^\hrm$ is not the effective resistance associated to a resistance form, as in Remark~\ref{rem:comparison-to-resistance-forms}, since (RF5) may fail. If $R^\hrm$ \emph{were} the effective resistance associated to a resistance form, then \cite[Prop.~2.10]{Kig03} would imply that $R^\hrm(x,y)$ is a metric, but this can be seen to be false by considering basic examples. See Example~\ref{exm:geometric-integers}, e.g. The same remarks also apply to the  boundary resistance $R^{\bdy}(x,y)$, discussed just below. 
\end{remark}

\begin{defn}\label{def:boundary-resistance}
  For an infinite network $(G,\cond)$ define the \emph{boundary resistance} between $x$ and $y$ by
  \linenopax
  \begin{align}\label{eqn:def:boundary-resistance}
    R^{\bdy}(x,y) := \frac1{R^W(x,y)^{-1} - R^F(x,y)^{-1}}.
  \end{align}
\end{defn}

  Intuitively, some portion of the wired/minimal current from $x$ to $y$ passes through infinity; the quantity $R^{\bdy}(x,y)$ gives the voltage drop ``across infinity''; see Remark~\ref{rem:shortcut-through-infinity}. From this perspective, infinity is ``connected in parallel''. The boundary $\bd G$ in \cite{bdG} is a more rigorous definition of the set at infinity.

\begin{theorem}\label{thm:boundary-resistance}
  The boundary resistance is equal to
  \linenopax
  \begin{align}
    R^\bdy(x,y) = 
    \frac{R^W(x,y)R^F(x,y)}{R^{\hrm}(x,y)}.
    \label{eqn:def:R^hrm(x,y)-R}
  \end{align}
  In particular, the resistance across the boundary is infinite if $\Harm = 0$. 
  \begin{proof}
    From \eqref{eqn:def:harmonic-resistance} one has $R^F(x,y) = 1/(R^W(x,y)^{-1} - R^\bdy(x,y)^{-1})$, which gives
    \begin{align*}
      \frac1{\energy(v_x-v_y)}
        &=  \frac1{\energy(f_x-f_y)} - \frac1{R^\bdy(x,y)}
    \end{align*}
    by Theorem~\ref{thm:free-resistance} and Theorem~\ref{thm:wired-resistance}, and hence
    \begin{align*}
      \frac1{R^\bdy(x,y)}
        &= \frac1{\energy(f_x-f_y)} - \frac1{\energy(v_x-v_y)}.
    \end{align*}
    Now solving for $R^\bdy$ gives
    \begin{align}
      R^\bdy(x,y) 
      = \frac{\energy(f_x-f_y)\energy(v_x-v_y)}{\energy(h_x-h_y)},
    \end{align}
    and the conclusion follows from \eqref{eqn:def:R^F(x,y)-energy}, \eqref{eqn:def:R^W(x,y)-energy}, and \eqref{eqn:def:R^ha(x,y)-energy}.
  \end{proof}
\end{theorem}

\section{Trace resistance}
\label{sec:trace-resistance}

The third type of subnetwork takes into account the connectivity of the complement of the subnetwork, but does not add anything extra.  The name ``trace'' is due to the fact that this approach comes by considering the trace of the Dirichlet form \energy to a subnetwork; see \cite{FOT94}. Several of the ideas in this section were explored previously in \cite{Kig01,Kig03,Metz}; see also \cite{Kig09}.

The discussion of the trace resistance and trace subnetworks requires some definitions relating the transition operator (i.e. Markov chain) \Prob to the probability measure $\prob^{(\cond)}$ on the space of (infinite) paths in \Graph which start at a fixed vertex $a$. Such a path is a sequence of vertices $\{x_n\}_{n=0}^\iy$, where $x_0=a$ and $x_n \nbr x_{n+1}$ for all $n$. 
  
\begin{defn}\label{def:Gamma(a,b)}
  Let $\Paths(a)$ be the space of all paths \cpath beginning at the vertex $a \in \verts$, and let $\Paths(a,b) \ci \Paths(a)$ be the subset of paths that reach $b$, and that do so \emph{before} returning to $a$:
  \linenopax
  \begin{align}\label{eqn:def:Paths(a,b)}
    \Paths(a,b) := \{\cpath \in \Paths(a) \suth b=x_n \text{ for some $n$, with } x_k \neq a \text{ for } 1 \leq k \leq n\}.
  \end{align}
\end{defn}

\begin{defn}\label{def:prob^cond}
   The space $\Paths(a)$ carries a natural probability measure $\prob^{(\cond)}$ defined by 
  \linenopax
  \begin{align}\label{eqn:prob^cond}
    \prob^{(\cond)}(\cpath) := \prod_{x_i \in \cpath} p(x_{i-1},x_i),
  \end{align}
  where $p(x,y) = \cond_{xy}/\cond(x)$.
  The construction of $\prob^{(\cond)}$ comes by applying Kolmogorov consistency to the natural cylinder-set Borel topology that makes $\Paths(a)$ into a compact Hausdorff space. 
\end{defn}

\begin{defn}\label{def:hitting-time}
  Let $X_m$ be a random variable which denotes the (vertex) location of the random walker at time $m$. Then let $\gt_x$ be the \emph{hitting time} of $x$, that is, the random variable which is the expected time at which the walker first reaches $x$:
  \linenopax
  \begin{align}\label{eqn:def:hitting-time}
    \gt_x := \min\{m \geq 0 \suth X_m = x\}.
  \end{align}
  More generally, $\gt_H$ is the time at which the walker first reaches the subnetwork $H$. For a walk started in $H$, this gives $\gt_H = 0$. 
\end{defn}

\subsection{The trace subnetwork}

It is well-known that networks $\{(\Graph,\cond)\}$ are in bijective correspondence with reversible Markov processes $\{\Prob\}$; this is immediate from the \emph{detailed balance equations} which follow from the symmetry of the conductance:
\linenopax
\begin{align*}
  \cond(x) p(x,y) = \cond_{xy} = \cond_{yx} = \cond(y) p(y,x).
\end{align*}
It follows from $\Lap = \cond(\one-\Prob)$ that networks are thus in bijective correspondence with \emph{Laplacians}, if one defines a Laplacian as in \eqref{eqn:def:laplacian}. That is, a Laplacian is a symmetric linear operator which is nonnegative definite, has kernel consisting of the constant functions, and satisfies $(\Lap\gd_x)(y) \leq 0$ for $x \neq y$. In other words, every row (and column) of $\trc(\Lap,H)$ sums to 0. (This is the negative of the definition of a \emph{Laplacian} as in \cite{Kig01} and \cite{CdV98}.) 
In this section, we exploit the bijection between Laplacians and networks to define the trace subnetwork. For $\verts[H] \ci \Graph$, the idea is as follows:
\linenopax
\begin{align*}
  G \;\longleftrightarrow\; \Lap 
  \;\limmode[\text{take the trace to \verts[H]}]\; \trc(\Lap, \verts[H]) 
  \;\longleftrightarrow\; H^\trc.
\end{align*}

\begin{defn}\label{def:trace-resistance}
  The \emph{trace} of \Graph to $\verts[H]$ is the network whose edge data is defined by the trace of \Lap to \verts[H], which is computed as the Schur complement of the Laplacian of $H$ with respect to \Graph. More precisely, write the Laplacian of \Graph as a matrix in block form, with the rows and columns indexed by vertices, and order the vertices so that those of $H$ appear first:
  \linenopax
  \begin{align}\label{eqn:Lap-block-decomp}
    \Lap =
      \begin{array}{l} \scalebox{0.70}{$H$} \\ \scalebox{0.70}{$H^\complm$} \end{array}
      \negsp[12]
      \left[\begin{array}{ll} A & B^T \\ B & D \end{array}\right],
  \end{align}
  where $B^T$ is the transpose of $B$. If $\ell(G) :=\{f:\verts \to \bR\}$, the corresponding mappings are
  \linenopax
  \begin{align}
    A:&\ell(H) \to \ell(H) & B^T&:\ell(H^\complm) \to \ell(H) \notag \\
    B:&\ell(H) \to \ell(H^\complm) & D&:\ell(H^\complm) \to \ell(H^\complm).
    \label{eqn:action-of-Lap-quadrants}
  \end{align}
\glossary{name={$\trc(\Lap, \verts[H])$},description={trace of the Laplacian to a subnetwork $H$},sort=tr,format=textbf}
  It turns out that the Schur complement
  \linenopax
  \begin{equation}\label{eqn:def:Schur-complement}
    \trc(\Lap, \verts[H]) := A - B^T D^{-1} B
  \end{equation}
  is the Laplacian of a subnetwork with vertex set $\verts[H]$; cf.~\cite[\S2.1]{Kig01} and Remark~\ref{rem:trace-valid-for-subsets}.\footnote{It will be clear from \eqref{eqn:Schur-complement-as-sum} that $D^{-1}$ always exists in this context, and hence \eqref{eqn:def:Schur-complement} is always well-defined. Furthermore, the existence of the trace is given in \cite[Prop.~2.10]{Kig03}; it is known from \cite[Lem.~2.1.5]{Kig01} that $D$ is invertible and negative semidefinite.} A formula for the conductances (and hence the adjacencies) of the trace is given in Theorem~\ref{thm:edges-of-H^S}. Denote this new subnetwork by $H^\trc$.

  If $\verts[H] \ci \verts$ is finite, then for $x,y \in H$, the trace of the resistance on $H$ is denoted $R_{H^\trc}(x,y)$, and defined as in Definition~\ref{def:effective-resistance}. The \emph{trace resistance} is then defined to be
  \linenopax
  \begin{align}\label{eqn:trace-resistance}
    R^\trc(x,y) := \lim_{k \to \iy} R_{G_k^\trc}(x,y),
  \end{align}
  where $\{G_k\}$ is any exhaustion of \Graph.
\end{defn}
  \glossary{name={$R^\trc(x,y)$},description={trace resistance metric},sort=Rt,format=textbf}
  \glossary{name={$R_H(x,y)$},description={relative resistance as computed with respect to the subnetwork $H$},sort=RH}

\begin{remark}\label{rem:trace-resistance-name}
  The name ``trace'' is due to the fact that this approach comes by considering the trace of the Dirichlet form \energy on a subnetwork; see \cite{FOT94}.
\end{remark}

Recall that $\Lap = \cond - \Trans = \cond(\id - \Prob)$, where \Trans is the transfer operator and $\Prob = \cond^{-1}\Trans$ is the probabilistic transition operator defined   \glossary{name={\Prob},description={probabilistic transition operator},sort=P} \glossary{name={\Trans},description={transfer operator},sort=T}
\glossary{name={\cond},description={conductance (multiplication) operator},sort=c}
pointwise by
  \linenopax
  \begin{align}
    \Prob u(x) = \sum_{y \nbr x} p(x,y) u(y), 
    \q\text{for } 
    p(x,y) = \frac{\cond_{xy}}{\cond(x)}.
    \label{eqn:def:p(x,y)}
  \end{align}
  \glossary{name={$p(x,y)$},description={transition probability of the random walk on the network: $\cond_{xy}/\cond(x)$},sort=prob}
  The function $p(x,y)$ gives transition probabilities, i.e., the probability that a random walker currently at $x$ will move to $y$ with the next step. Since
  \linenopax
  \begin{align}\label{eqn:def:reversible}
    \cond(x) p(x,y) = \cond(y) p(y,x),
  \end{align}  
  the transition operator \Prob determines a \emph{reversible} Markov process with state space \verts; see \cite{LevPerWil08,Lyons:ProbOnTrees}. Note that the harmonic functions (i.e., $\Lap h=0$) are precisely the fixed points of \Prob (i.e., $\Prob h = h$). The proof of the next theorem requires a couple more definitions which relate \Prob to the probability measure $\prob^{(\cond)}$ on the space of paths in \Graph. Recall from Definition~\ref{def:paths} that a path is a sequence of vertices $\{x_n\}_{n=0}^\iy$, where $x_0=a$ and $x_n \nbr x_{n+1}$ for all $n$. 
  
\begin{defn}\label{def:Gamma(a,b)}
  Let $\Paths(a)$ denote the space of all paths \cpath beginning at the vertex $a \in \verts$. 
    \glossary{name={$\Paths(a)$},description={the space of all paths in $G$ beginning at the vertex $a$},sort=G,format=textbf}
Then $\Paths(a,b) \ci \Paths(a)$ consists of those paths that reach $b$, and \emph{before} returning to $a$:
  \linenopax
  \begin{align}\label{eqn:def:Paths(a,b)}
    \Paths(a,b) := \{\cpath \in \Paths(a) \suth b=x_n \text{ for some $n$, and } x_k \neq a, 1 \leq k \leq n\}.
  \end{align}
  If $a,b \in \bd H$, then we write 
  \linenopax
  \begin{align}\label{eqn:def:Paths(a,b)-outside-H}
    \Paths(a,b)\evald{H^\complm} 
    := \{\cpath \in \Paths(a,b) \suth x_i \in H^\complm, 0 < i < \gt_b\},
  \end{align}
  for the set of paths from $a$ to $b$ that do not pass through any vertex in $\verts[H]$. 
\end{defn}

\begin{remark}\label{rem:trivial-paths}
  Note that if $x,y \in \bd H$ are adjacent, then any path of the form $\cpath = (x,y,\dots)$ is trivially in $\Paths(a,b)\evald{H^\complm}$. 
\end{remark}

\begin{defn}\label{def:prob^cond}
   The space $\Paths(a)$ carries a natural probability measure $\prob^{(\cond)}$ defined by 
  \linenopax
  \begin{align}\label{eqn:prob^cond}
    \prob^{(\cond)}(\cpath) := \prod_{x_i \in \cpath} p(x_{i-1},x_i),
  \end{align}
  where $p(x,y)$ is as in \eqref{eqn:def:p(x,y)}.
  \glossary{name={\cpath},description={path},sort=G}
  \glossary{name={$\prob(\cpath)$},description={probability of a path},sort=P,format=textbf}
  \glossary{name={$\prob^{(\cond)}$},description={a measure on the space of all paths in $G$},sort=Pc,format=textbf}
  The construction of $\prob^{(\cond)}$ comes by applying Kolmogorov consistency to the natural cylinder-set Borel topology that makes $\Paths(a)$ into a compact Hausdorff space; cf. \S\ref{sec:Probabilistic-interpretation} for further discussion. 
\end{defn}

\begin{defn}\label{def:prob[a->b]}
  Let $\prob[a \to b]$ denote the probability that a random walk started at $a$ will reach $b$ before returning to $a$. That is,
  \glossary{name={$\prob[a \to b]$},description={probability of the random walk started from $a$ reaching $b$},sort=P,format=textbf}
  \linenopax
  \begin{align}\label{eqn:def:prob-a-to-b}
    \prob[a \to b] := \prob^{(\cond)}(\Paths(a,b)).
  \end{align}
  Note that this is equivalent to 
  \linenopax
  \begin{align}\label{eqn:def:prob-a-to-b-as-hitting}
    \prob[a \to b] = \prob_a[\gt_b < \gt_a] := \prob[\gt_b < \gt_a \,|\, x_0=a],
  \end{align}
  where $\gt_a$ is the \emph{hitting time} of $a$, i.e., the expected time of the first visit to $a$, \emph{after} leaving the starting point. If $a,b \in \bd H$, then we write 
  \linenopax
  \begin{align}\label{eqn:def:prob-a-to-b-outside-H}
    \prob[a \to b]\evald{H^\complm} 
    := \prob^{(\cond)}\left(\Paths(a,b)\evald{H^\complm}\right),
  \end{align}
  \glossary{name={$\complement$},description={$H^\complm$ is the complement of $H$ in $G$},sort=C}
  that is, the probability that a random walk started at $a$ will reach $b$ via a path lying outside $H$ (except for the start and end points, of course).
\end{defn}

\begin{remark}[More probabilistic notation]\label{rem:prob-notation}\label{rem:notation-restriction-to-Gk}
  The formulation in \eqref{eqn:def:prob-a-to-b-outside-H} is conditioning $\prob^{(\cond)}(\Paths(a,b))$ on avoiding $H$; the notation is intended to evoke something like ``$\prob[a \to b \,|\, \cpath \ci H^\complm]$''. However, this would not be correct because $a,b \in H$ and \cpath may pass through $H$ after $\gt_b$.

  In Theorem~\ref{thm:edges-of-H^S}, we use the following common notation as in \cite{Spitzer} or \cite{Woess00}, for example. All notations are for the random walk started at $x$.
  \begin{align*}
    \Prob^n(x,y) = p^{(n)}(x,y) =\prob_x[X_n = y] &\q\text{prob. the walk is at $y$ after $n$ steps} \\
    G(x,y) = {\textstyle\sum}_{n=0}^\iy p^{(n)}(x,y) &\q\text{exp. number of visits to $y$} \\
    f^{(n)}(x,y) = \prob_x[\gt_y = n] &\q\text{prob. the walk first reaches $y$ on the \nth step} \\
    F(x,y) = {\textstyle\sum}_{n=0}^\iy f^{(n)}(x,y) &\q\text{prob. the walk ever reaches $y$}
  \end{align*}
  Note that if the walk is killed when it reaches $y$, then $p^{(n)}(x,y) = f^{(n)}(x,y)$ because the first time it reaches $y$ is the only time it reaches $y$. Therefore, when the walk is conditioned to end upon reaching a set $S$, one has $G(x,y) = F(x,y)$ for all $y \in S$.
\end{remark}

\begin{theorem}\label{thm:edges-of-H^S}
  For $\verts[H] \ci \verts$, the conductances in the trace subnetwork $H^\trc$ are given by
  \linenopax
  \begin{align}\label{eqn:conductances-in-HS}
    \cond_{xy}^\trc = \cond_{xy} + \cond(x) \prob[x \to y]\evald{H^\complm}.
  \end{align}
  \glossary{name={$\cond_{xy}^\trc$},description={conductances in the trace subnetwork},sort=c}
  Consequently, the transition probabilities in the trace subnetwork are given by
  \linenopax
  \begin{align}\label{eqn:transition-probabilities-in-HS}
    p^\trc(x,y) = p(x,y) + \prob[x \to y]\evald{H^\complm}.
  \end{align}
\begin{proof}
  Using subscripts to indicate the block decomposition corresponding to $H$ and $H^\complm$ as in \eqref{eqn:Lap-block-decomp}, the Laplacian may be written as
  \linenopax
  \begin{align*}
    \Lap =
    \left[\begin{array}{cc} \cond_A(\one - \Prob_A) & -\cond_A\Prob_{B^T} \\ -\cond_D\Prob_{B} & \cond_D(\one - \Prob_{D}) \end{array}\right],
    \qq \text{for} \qq
    \cond =
      \begin{array}{l} \scalebox{0.70}{$H$} \\ \scalebox{0.70}{$H^\complm$} \end{array}
      \negsp[12]
      \left[\begin{array}{ll} \cond_A &  \\  & \cond_D \end{array}\right].
  \end{align*}
  Then the Schur complement is
  \linenopax
  \begin{align}
    \trc(\Lap,H)
    &= \cond_A - \cond_A \Prob_A - \cond_A \Prob_{B^T} (\id - \Prob_D)^{-1} \cond_{D}^{-1} \cond_{D} \Prob_B \notag \\
    &= \cond_A - \cond_A\left(\Prob_A + \Prob_{B^T} \left(\sum_{n=0}^\iy \Prob_D^n\right) \Prob_B \right) \notag \\
    &= \cond_A(\id - \Prob_\xX).
    \label{eqn:Schur-complement-as-sum}
  \end{align}
  Note that $\Prob_D$ is substochastic, and hence the RW has positive probability of hitting $\bd G_k$, whose vertices act as absorbing states. This means that the expected number of visits to any vertex in $H^\complm$ is finite and hence the matrix $\Prob_\xX$ has finite entries. 
  
  Meanwhile, using $\Prob_A(x,y)$ to denote the \nth[(x,y)] entry of the matrix $\Prob_A$, and $\gt_H^+$ as in Definition~\ref{def:prob[a->b]}, we have
  \linenopax
  \begin{align}
    \prob[x \to y]\evald{H^\complm}
    &= \prob^{(\cond)}\left(\Paths(x,y)\evald{H^\complm}\right) \notag \\
    &= \prob^{(\cond)}\left(\bigcup_{k=1}^\iy \{\cpath \in \Paths(x,y)\evald{H^\complm} \suth \gt_H^+ = k \}\right) \notag \\
    &= \prob^{(\cond)}\left(\{\cpath \in \Paths(x,y)\evald{H^\complm} \suth \gt_H^+ = 1 \}\right) \\
      &\hstr[10] + \sum_{k=2}^\iy \prob^{(\cond)}\left(\{\cpath \in \Paths(x,y)\evald{H^\complm} \suth \gt_H^+ = k \}\right) \notag \\
    &= \Prob_A(x,y) +\sum_{n=0}^\iy \sum_{s,t} \Prob_{B^T}(x,s) \Prob_D^n(s,t) \Prob_B(t,y) \label{eqn:prob-a-to-b-avoiding-H-computation}\\ 
    &= \Prob_\xX(x,y). \notag
  \end{align}
  To justify \eqref{eqn:prob-a-to-b-avoiding-H-computation}, note that by \eqref{eqn:action-of-Lap-quadrants}, $\Prob_D^n$ corresponds to steps taken in $H^\complm$. Therefore,
  \linenopax
  \begin{align*}
    \left(\Prob_{B^T} \left(\sum_{n=0}^\iy \Prob_D^n\right) \Prob_B\right)(x,y)
    = \Prob_{B^T} \Prob_D^0 \Prob_B(x,y) + \Prob_{B^T}\Prob_D^1 \Prob_B(x,y) + \dots
  \end{align*}
  is the probability of the random walk taking a path that steps from $x \in H$ to $H^\complm$, meanders through $H^\complm$ for any finite number of steps, and finally steps to $y \in H$. Since $y \notin H^\complm$, 
  \linenopax
  \begin{align*}
    \Prob_{B^T} \Prob_D^k \Prob_B(x,y) 
    = \prob_x[X_{k+2} = y] 
    = \prob_x[\gt_y = k+2] ,
  \end{align*}
because the walk can only reach $y$ on the last step, as in Remark~\ref{rem:prob-notation}. It follows by classical theory (see \cite{Spitzer}, for example) that the sum in \eqref{eqn:prob-a-to-b-avoiding-H-computation} is a probability (as opposed to an expectation, etc.) and justifies the probabilistic notation $\Prob_\xX$ in \eqref{eqn:Schur-complement-as-sum}.
  Note that $\Prob_A(x,y)$ corresponds to the one-step path from $x$ to $y$, which is trivially in $\Paths(x,y)\evald{H^\complm}$ by \eqref{eqn:def:Paths(a,b)-outside-H}.
  Since $\Prob_A(x,y) = p(x,y) = \cond_{xy}/\cond(x)$, the desired conclusion \eqref{eqn:conductances-in-HS} follows from combining \eqref{eqn:Schur-complement-as-sum}, \eqref{eqn:prob-a-to-b-avoiding-H-computation}, and \eqref{eqn:def:prob-a-to-b-outside-H}. Of course, \eqref{eqn:transition-probabilities-in-HS} follows immediately by dividing through by $\cond(x)$.
\end{proof}
\end{theorem}

  The authors are grateful to Jun Kigami for helpful conversations and guidance regarding the proof of Theorem~\ref{thm:edges-of-H^S}. 
 
\begin{remark}\label{rem:edges-of-HS}
  It is clear from \eqref{eqn:conductances-in-HS} that the edge sets of $\inn H$ and $\inn H^\trc$ are identical, but the conductance between two vertices $x,y \in \bd H^\trc$ is greater iff there is a path from $x$ to $y$ that does not pass through $H$. Indeed, if there is a path from $x$ to $y$ which lies entirely in $H^\complm$ except for the endpoints, then $x$ and $y$ will be adjacent in $H^\trc$, even if they were not adjacent in $H$.
\end{remark}

\begin{remark}[The trace construction is valid for general subsets of vertices]
  \label{rem:trace-valid-for-subsets}\label{rem:edges-of-HS}
  While Definition~\ref{def:trace-resistance} applies to a (connected) subnetwork of \Graph, it is essential to note that Theorem~\ref{thm:edges-of-H^S} applies to arbitrary subsets \verts[H] of \verts. 
  
  It is clear from \eqref{eqn:conductances-in-HS} that the edge sets of $\inn H$ and $\inn H^\trc$ are identical, but the conductance between two vertices $x,y \in \bd H^\trc$ is greater iff there is a path from $x$ to $y$ that does not pass through $H$. Indeed, if there is a path from $x$ to $y$ which lies entirely in $H^\complm$ except for the endpoints, then $x$ and $y$ will be adjacent in $H^\trc$, even if they were not adjacent in $H$.
\end{remark}

\begin{remark}[Resistance distance via Schur complement]
  \label{rem:Schur-does-network-reduction}
  \label{rem:Resistance-distance-via-Schur-complement} 
  A theorem of Epifanov states that every finite planar network with vertices $x,y$ can be reduced to a single equivalent conductor via the use of three simple transformations: parallel, series, and $\nabla$-\textsf{Y}; cf. \cite{Epifanov66,Truemper89} as well as \cite[\S2.3]{Lyons:ProbOnTrees} and \cite[\S7.4]{CdV98}. More precisely,
  \begin{enumerate}[(i)]
    \item Parallel. Two conductors $\cond_{xy}^{(1)}$ and $\cond_{xy}^{(2)}$ connected in parallel can be replaced by a single conductor $\cond_{xy} = \cond_{xy}^{(1)} + \cond_{xy}^{(2)}$.
    \item Series. If $z$ has only the neighbours $x$ and $y$, then $z$ may be removed from the network and the edges $\cond_{xz}$ and $\cond_{yz}$ should be replaced by a single edge $\cond_{xy} = (\cond_{xz}^{-1}+\cond_{yz}^{-1})^{-1}$.
    \item $\nabla$-\textsf{Y}. Let $t$ be a vertex whose only neighbours are $x,y,z$. Then this ``\textsf{Y}'' may be replaced by a triangle (``$\nabla$'') which does not include $t$, with conductances
      \linenopax
      \begin{align*}
        \cond_{xy} = \frac{\cond_{xt}\cond_{ty}}{\cond(t)}, \;
        \cond_{yz} = \frac{\cond_{yt}\cond_{tz}}{\cond(t)}, \;
        \cond_{xz} = \frac{\cond_{xt}\cond_{tz}}{\cond(t)}.
      \end{align*}
      This transformation may also be inverted, to replace a $\nabla$ with a \textsf{Y} and introduce a new vertex. 
  \end{enumerate}
  It is a fun exercise to obtain the series and $\nabla$-\textsf{Y} formulas by applying the Schur complement technique to remove a single vertex of degree 2 or 3 from a network. Indeed, these are both special cases of the following: let $t$ be a vertex of degree $n$, and let $H$ be the (star-shaped) subnetwork consisting only of $t$ and its neighbours. If we write the Laplacian for just this subnetwork with the \nth[t] row \& column last, then
  \linenopax
  \begin{align*}
    \Lap|_H =
    \left[\begin{array}{rrrc}
      \cond_{x_1t} & \dots & 0 & -\cond_{x_1t} \\
      \vdots & \ddots & \vdots & \vdots \\
      0 & \dots & \cond_{x_nt} & -\cond_{x_nt} \\
      -\cond_{x_1t} & \dots & -\cond_{x_nt} & \cond(t)
    \end{array}\right]
  \end{align*}
  and the Schur complement is
  \linenopax
  \begin{align*}
    \trc(\Lap|_H,H\less\{t\})
    =
    \left[\begin{array}{rrrc}
      \cond_{x_1t} & \dots & 0 \\
      \vdots & \ddots & \vdots \\
      0 & \dots & \cond_{x_nt} \\
    \end{array}\right]
    - \frac{1}{\cond(t)}
    \left[\begin{array}{rrrc}
      \cond_{x_1t} \\
      \vdots \\
      \cond_{x_nt} 
    \end{array}\right]
    \left[\begin{array}{rrrc}
      \cond_{x_1t} & \dots & \cond_{x_nt} 
    \end{array}\right],
  \end{align*}
  whence the new conductance from $x_i$ to $x_j$ is given by $\cond_{x_it}\cond_{tx_j}/\cond(t)$. It is interesting to note that the operator being subtracted corresponds to the projection to the rank-one subspace spanned by the probabilities of leaving $t$:
  \linenopax
  \begin{align*}
    \frac{1}{\cond(t)}
    \left[\begin{array}{r}
      \cond_{x_1t} \\
      \vdots \\
      \cond_{x_nt} 
    \end{array}\right]
    \left[\begin{array}{rrr}
      \cond_{x_1t} & \dots & \cond_{x_nt} 
    \end{array}\right]
    = \cond(t) |v\ra \la v|,
  \end{align*}
  using Dirac's ket-bra notation for the projection to a rank-1 subspace spanned by $v$ where
  \linenopax
  \begin{align*}
    v = \left[\begin{array}{rrr}
      p(t,x_1)  & \dots & p(t,x_n)  
    \end{array}\right].
  \end{align*}
In fact, $|v\ra \la v | = \Prob_X$, in the notation of \eqref{eqn:Schur-complement-as-sum}. In general, the trace construction (Schur complement) has the effect of probabilistically projecting away the complement of the subnetwork.

    In Remark~\ref{rem:Resistance-distance-via-network-reduction} we described how the effective resistance can be interpreted as the correct resistance for a single edge which replaces a subnetwork. The following corollary of Theorem~\ref{thm:edges-of-H^S} formalizes this interpretation by exploiting the fact that the Schur complement construction is viable for arbitrary subsets of vertices; see Remark~\ref{rem:trace-valid-for-subsets}. In this case, one takes the trace of the (typically disconnected) subset $\{x,y\} \ci \verts$; note that $\left[\begin{smallmatrix} 1 & -1 \\ -1 & 1 \end{smallmatrix}\right]$ is the Laplacian of the trivial 2-vertex network when the edge between them has unit conductance. The following result is also an extension of \cite[(2.1.4)]{Kig01} to infinite networks.
\end{remark}

\begin{cor}\label{thm:resistance-as-2-vert-network}
  Let $\verts[H] = \{x,y\}$ be any two vertices of a transient network \Graph. Then the trace resistance can be computed via
  \linenopax
  \begin{equation}\label{eqn:resistance-as-2-vert-network}
    \trc(\Lap,H)
    = \frac1{R^\trc(x,y)} \left[\begin{array}{rr} 1 & -1 \\ -1 & 1 \end{array}\right]
    = A - B^T D^{-1} B.
  \end{equation}
  \begin{proof}
    Take $H=\{x,y\}$ in Theorem~\ref{thm:edges-of-H^S}. As discussed in Remark~\ref{rem:trace-valid-for-subsets}, it is not necessary to have $x \nbr y$. Note that in this case, $\left(\Prob_{B^T} \sum_n \Prob_D^n \Prob_B\right)(x,y)$ corresponds all paths from $x$ to $y$ that consist of more than one step:
  \linenopax
  \begin{align}
    \prob[x \to y]\evald{H^\complm}
    &= \Prob_A(x,y) + \left(\Prob_{B^T} \sum_{n=0}^\iy \Prob_D^n \Prob_B\right)(x,y) 
        = p(x,y) 
    + \sum_{|\cpath| \geq 2} \prob(\cpath).
    \label{eqn:sum-over-paths-decomposition}
    \qedhere
  \end{align}
  \end{proof}
\end{cor}

\version{}{\marginpar{Alternative idea}Use Kigami's approach (\energy, Theorem~\ref{thm:resistance-as-2-vert-network}) and take limits.}

\begin{cor}\label{thm:resistance-as-path-integral}
  The trace resistance $R^\trc(x,y)$ is given by
  \linenopax
  \begin{equation}\label{eqn:resistance-as-path-integral}
    R^\trc(x,y) = \frac1{\cond(x) \prob[x \to y]}
  \end{equation}
  \begin{proof}
    Again, take $\verts[H] = \{x,y\}$.  Then    
    \linenopax
    \begin{align*}
      R^\trc(x,y)^{-1}
      = \cond_{xy}^{H^S} 
      &= \cond_{xy} + \cond(x) \prob[x \to y]\evald{H^\complm} \\
      &= \cond(x) \left(p(x,y) + \prob[x \to y]\evald{H^\complm}\right) \\
      &= \cond(x) \prob[x \to y],
    \end{align*}
    where Corollary~\ref{thm:resistance-as-2-vert-network} gives the first equality and Theorem~\ref{thm:edges-of-H^S} gives the second.
  \end{proof}
\end{cor}
   
\begin{remark}[Effective resistance as ``path integral'']
  \label{rem:resistance-as-path-integral}
  Corollary~\ref{thm:resistance-as-path-integral} may also be obtained by the more elegant (and much shorter) approach of \cite[\S2.2]{Lyons:ProbOnTrees}, where it is stated as follows: the mean number of times a random walk visits $a$ before reaching $b$ is $\prob[a \to b]^{-1} = \cond(a) R(a,b)$. We give the present proof to highlight and explain the underlying role of the Schur complement with respect to network reduction; see Remarks~\ref{rem:trace-valid-for-subsets}--\ref{rem:Resistance-distance-via-Schur-complement}.
  A key point of the present approach is to emphasize the expression of effective resistance $R(a,b)$ in terms of a \emph{sum over all possible paths from $a$ to $b$}.   By Remark~\ref{rem:wired-vs-free-as-Kirchhoff-conditions+}, it is apparent that this ``path-integral'' interpretation makes $R^\trc$ much more closely related to $R^F$ than to $R^W$, as seen by the following results.
\end{remark}

\begin{cor}[{\cite[Prop.~2.1.11]{Kig01}}]
  \label{thm:trace-invariance-of-resistance}
  Let $H_2 \ci H_1$ be finite subnetworks of a transient network \Graph. Then for $a,b \in \verts[H_2]$, one has $R_{H_1}^S(a,b) = R_{H_2}^S(a,b)$.
\end{cor}

\begin{cor}\label{thm:trace-resistance-is-free-resistance}
  On any network, $R^\trc(a,b) = R^F(a,b)$.
  \begin{proof}
    By Corollary~\ref{thm:trace-invariance-of-resistance}, it is clear that $R_{G_k^\trc}(a,b) = R_{G_{k+1}^\trc}(a,b)$ for all $k$. Meanwhile, any path from $a$ to $b$ will lie in $G_k$ for sufficiently large $k$, so it is clear by Theorem~\ref{thm:resistance-as-path-integral}, the sequence $\{R^F_{G_k}(a,b)\}_{k=0}^\iy$ is monotonically decreasing with limit $R^F(a,b) = R^\trc(a,b)$.
  \end{proof}
\end{cor}

\begin{remark}
  \version{}{\marginpar{Work out, and insert here, the corresponding formula for $R^W$}}
  Writing $[x \to y \,|\, \cpath \ci H]$ to indicate a restriction to paths from $x$ to $y$ that lie entirely in $H$, as in Remark~\ref{rem:notation-restriction-to-Gk}, one has
  \linenopax
  \begin{align*}
    R_{G_k^\trc}(x,y) 
    &= \frac1{\cond(x) \left(\prob[x \to y \,|\, \cpath \ci G_k] + \prob[x \to y  \,|\, \cpath \nsubseteq G_k]\right)}  \\
    &\leq \frac1{\cond(x) \prob[x \to y \,|\, \cpath \ci G_k]} 
    = R_{G_k^F}(x,y).
  \end{align*}
  Essentially, Corollary~\ref{thm:trace-invariance-of-resistance} is an expression of the first equality and Corollary~\ref{thm:trace-resistance-is-free-resistance} is a consequence of the inequality and how it tends to an equality as $k \to \iy$.
\end{remark}

\subsection{The shorted operator}
\label{sec:shorted-operator}

It is worth noting that the operator $D$ defined in \eqref{eqn:action-of-Lap-quadrants} is always invertible as in the discussion following \eqref{eqn:Schur-complement-as-sum}.
However, the Schur complement construction is valid more generally. As is pointed out in \cite{Butler}, the \emph{shorted operator} generalizes the Schur complement construction to positive operators on a (typically infinite-dimensional) Hilbert space \sH; see \cite{Anderson,AndersonTrapp,Krein}. In general, let $T=T^\ad$ be a positive operator so $\la \gf, T\gf\ra \geq 0$ for all $\gf \in \sH$, and let $S$ be a closed subspace of \sH. Partition $T$ analogously to \eqref{eqn:action-of-Lap-quadrants}, so that $A:S \to S$, $B:S \to S^\complm$, $B^T:S^\complm \to S$, and $D:S^\complm \to S^\complm$. 

\begin{theorem}[{\cite{AndersonTrapp}}]
  With respect to the usual ordering of self-adjoint operators, there exists a unique operator $\sS{h}(T)$ such that
  \linenopax
  \begin{align*}
    \sS{h}(T) = \sup_L\left\{L \geq 0 \suth \left[\begin{array}{cc} L & 0 \\ 0 & 0  
  \end{array}\right] \leq T\right\},
  \end{align*}
  and it is given by
  \linenopax
  \begin{align*}
    \sS{h}(T) = \lim_{\ge \to 0^+} \left(A - B^T (D + \ge)^{-1} B\right).
  \end{align*}
  In particular, the shorted operator coincides with the Schur complement, whenever the latter exists.
\end{theorem}

There is another characterization of the shorted operator due to \cite{Butler}.

\begin{theorem}[{\cite{Butler}}]
  Suppose $\{\gy_n\} \ci \sH$ is a sequence satisfying $\la\gy_n, D\gy_n\ra \leq M$ for some $M \in \bR$, and $\lim_{n \to \iy} T \left[\begin{smallmatrix} \gf \\ \gy_n \end{smallmatrix}\right] = \left[\begin{smallmatrix} \gq \\ 0 \end{smallmatrix}\right]$. Then $\sS{h}(T) \gf = \lim_{n \to \iy} \left(A \gf + B^T\gy_n\right)$.
\end{theorem}

\section{Projections in Hilbert space and the conditioning of the random walk}
\label{sec:probabilistic-interp-of-v_x-and-f_x}
  
  In Remark~\ref{rem:wired-vs-free-as-boundary-conditions}, we gave an operator-theoretic account of the difference between $R^F$ and $R^W$. The foregoing probabilistic discussions might lead one to wonder if there is a probabilistic counterpart. An alternative approach is given in \cite[App.~B]{Kig03}.
  
  On a finite network, it is well-known that 
  \begin{align}\label{eqn:vx-as-prob-on-finite}
    v_x = R(o,x) u_x,  
  \end{align}
  where $u_x(y)$ is the probability that a random walker (RW) started at $y$ reaches $x$ before $o$:
  \linenopax
  \begin{align}\label{eqn:u_x-defined-probabilistically}
    u_x(y) := \prob_y[\gt_x < \gt_o].
  \end{align}
  Here again, $\gt_x$ denotes the hitting time of $x$ as in Definition~\ref{def:hitting-time}. Note that \eqref{eqn:def:R(x,y)-energy} gives $u_x = \frac{v_x}{\energy(v_x)}$. The relationship \eqref{eqn:vx-as-prob-on-finite} is discussed in \cite{DoSn84,LevPerWil08,Lyons:ProbOnTrees}. 
  
  Theorem~\ref{thm:fx-as-prob} is a wired extension of \eqref{eqn:vx-as-prob-on-finite} to transient networks. 
  The corresponding free version appears in Conjecture~\ref{thm:v_x-as-prob-on-infinite-network}.
   
\begin{theorem}\label{thm:fx-as-prob}
  On a transient network, let $f_x$ be the representative of $\Pfin v_x$ specified by $f_x(o)=0$. Then for $x \neq o$, $f_x$ is computed probabilistically by
  \linenopax
  \begin{align}\label{eqn:prob-expr-for-fx}
    f_x(y) = R^W(o,x) \left(\prob_y[\gt_x < \gt_o] 
    \right. \\ & \hstr[5] \left. 
    + \prob_y^{G}[\gt_o = \gt_x = \iy] \prob_x^{G}[\gt_o = \iy] \lim_{k \to \iy} \tfrac{\cond(x)}{\cond(\iy_k)} \prob_{\iy_k}^{G_k^W}[\gt_{\iy_k} < \gt_{\{x,o\}}]\right).
  \end{align}
  \begin{proof}
    Fix $x,y$ and an exhaustion $\{G_k\}_{k=1}^\iy$, and suppose without loss of generality that $o,x,y \in G_1$. 
    Since $v_x=f_x$ on any finite network, the identity \eqref{eqn:vx-as-prob-on-finite} gives $f_x^{(k)} = R_{G_k^W}(o,x) \check u_x^{(k)}$, where $f_x^{(k)}$ is the unique solution to $\Lap v = \gd_x - \gd_o$ on the finite (wired) subnetwork $G_k^W$, and 
    \linenopax
    \begin{align*}
      \check u_x^{(k)}(y) := \prob_y^{G_k^W}[\gt_x < \gt_o].
    \end{align*}
    where the superscript indicates the network in which the random walk travels.
    As in the previous case, we just need to check the limit of $\check u_x^{(k)}$, for which, we have
    \begin{align}\label{eqn:fx-as-prob-derivation-1}
       \check u_x^{(k)}(y) 
       &= \prob_y^{G_k^W}[\gt_x < \gt_o \text{ \& } \gt_x < \gt_{\iy_k}] + \prob_y^{G_k^W}[\gt_x < \gt_o \text{ \& } \gt_x > \gt_{\iy_k}]
    \end{align}
    %
    The first probability in \eqref{eqn:fx-as-prob-derivation-1} is
    \begin{align*}
      \prob_y^{G_k^W}[\gt_x < \gt_o \text{ \& } \gt_x < \gt_{\iy_k}]
      &= \prob_y^{G}[\gt_x < \gt_o \text{ \& } \gt_x < \gt_{G_k^\complm}] \\
      \limas{k}&\hstr[2.2] \prob_y^{G}[\gt_x < \gt_o \text{ \& } \gt_x < \iy] 
      = \prob_y^{G}[\gt_x < \gt_o],
    \end{align*}
    where the last equality follows because $\gt_x < \gt_o$ implies $\gt_x < \iy$.
    
    The latter probability in \eqref{eqn:fx-as-prob-derivation-1} measures the set of paths which travel from $y$ to $\iy_k$ without hitting $x$ or $o$, and then on to $x$ without passing through $o$, and hence can be rewritten
    \begin{align*}
      \prob_y^{G_k^W}&[\gt_{\iy_k} < \gt_x < \gt_o]
      = \prob_y^{G_k^W}[\gt_{\iy_k} < \gt_{\{o,x\}}] \prob_{\iy_k}^{G_k^W}[\gt_x < \gt_o] \\
      &= \prob_y^{G_k^W}[\gt_{\iy_k} < \gt_{\{o,x\}}] \left(\prob_{\iy_k}^{G_k^W}[\gt_{\iy_k} < \gt_{\{x,o\}}]\prob_{\iy_k}^{G_k^W}[\gt_x < \gt_{\{\iy_k,o\}}] + \prob_{\iy_k}^{G_k^W}[\gt_x < \gt_{\{\iy_k,o\}}]\right),
    \end{align*}
    since a walk starting at $\iy_k$ may or may not return to $\iy_k$ before reaching $x$. 
    
    First, consider only those walks which do not loop back through $\iy_k$ (i.e., multiply out the above expression and take the second term) to observe
    \begin{align}
      \prob_y^{G_k^W}[\gt_{\iy_k} < \gt_{\{o,x\}}] &\prob_{\iy_k}^{G_k^W}[\gt_x < \gt_{\{\iy_k,o\}}]
      = \prob_y^{G_k^W}[\gt_{\iy_k} < \gt_{\{o,x\}}]\prob_x^{G_k^W}[\gt_{\iy_k} < \gt_o] \tfrac{\cond(x)}{\cond(\iy_k)} \label{eqn:fx-as-prob-derivation-3} \\
      &= \left(1 - \prob_y^{G_k^W}[\gt_{\{o,x\}} < \gt_{\iy_k}]\right) \left(1 - \prob_x^{G_k^W}[\gt_o < \gt_{\iy_k}]\right) \tfrac{\cond(x)}{\cond(\iy_k)} \notag \\
      &= \left(1 - \prob_y^{G}[\gt_{\{o,x\}} < \gt_{G_k^\complm}]\right) \left(1 - \prob_x^{G}[\gt_o < \gt_{G_k^\complm}]\right) \tfrac{\cond(x)}{\cond(\iy_k)} \notag \\
      \limas{k} 
      &\hstr[2.2] \left(1 - \prob_y^{G}[\gt_{\{o,x\}} < \iy]\right) \left(1 - \prob_x^{G}[\gt_o < \iy]\right) \lim_{k \to \iy} \tfrac{\cond(x)}{\cond(\iy_k)} \notag \\
      &= \prob_y^{G}[\gt_o = \gt_x = \iy] \prob_x^{G}[\gt_o = \iy] \lim_{k \to \iy} \tfrac{\cond(x)}{\cond(\iy_k)}.
      \label{eqn:fx-as-prob-derivation-2}
    \end{align}
    Note that \eqref{eqn:fx-as-prob-derivation-3} comes by reversibility of the walk, and the way probability is computed for paths from $\iy_k$ to $x$ which avoid $o$ and $\iy_k$. Since the network is transient, $\sum_{k=1}^\iy{\cond(\iy_k)^{-1}}$ is summable by Nash-William's criterion and so $\lim_{k \to \iy} \tfrac{\cond(x)}{\cond(\iy_k)} = 0$ causes \eqref{eqn:fx-as-prob-derivation-2} to vanish. 
    
    Now for walks which do loop back through $\iy_k$, the same arguments as above yield
    \begin{align*}
      \prob_y^{G_k^W}[\gt_{\iy_k} < \gt_{\{o,x\}}] &\prob_{\iy_k}^{G_k^W}[\gt_{\iy_k} < \gt_{\{x,o\}}]\prob_{\iy_k}^{G_k^W}[\gt_x < \gt_{\{\iy_k,o\}}] \\
      &\limas{k} \hstr[1] \prob_y^{G}[\gt_o = \gt_x = \iy] \prob_x^{G}[\gt_o = \iy] \lim_{k \to \iy} \tfrac{\cond(x)}{\cond(\iy_k)} \prob_{\iy_k}^{G_k^W}[\gt_{\iy_k} < \gt_{\{x,o\}}],
    \end{align*}
    and the conclusion follows.
  \end{proof}
\end{theorem}


The following conjecture expresses a free extension of \eqref{eqn:vx-as-prob-on-finite} to infinite networks. We offer an erroneous ``proof'' in the hopes that it may inspire the reader to find a correct proof. The error is discussed in Remark~\ref{rem:error-in-prob-expr-for-vx}, just below. In the statement of Conjecture~\ref{thm:v_x-as-prob-on-infinite-network}, we use the notation
\begin{align}\label{eqn:bounded-traj}
  |\cpath| < \iy
\end{align}
to denote the event that the walk is bounded, i.e., that the trajectory is contained in a finite subnetwork of \Graph.

\begin{conj}\label{thm:v_x-as-prob-on-infinite-network}
  On an infinite resistance network, let $v_x$ be the representative of an element of the energy kernel specified by $v_x(o)=0$. Then for $x \neq o$, $v_x$ is computed probabilistically by
  \linenopax
  \begin{align}\label{eqn:prob-expr-for-vx}
    v_x(y) = R^F(o,x) \prob_y[\gt_x < \gt_o \,|\, |\cpath| < \iy],
  \end{align}
  that is, the walk is conditioned to lie entirely in some finite subnetwork as in \eqref{eqn:bounded-traj}. 
  \begin{proof}[``Proof.'']
    Fix $x,y$ and suppose without loss of generality that $o,x,y \in G_1$. One can  write \eqref{eqn:vx-as-prob-on-finite} on $G_k$ as $v_x^{(k)} = R_{G_k^F}(o,x) u_x^{(k)}$. In other words, $v_x^{(k)}$ is the unique solution to $\Lap v = \gd_x - \gd_o$ on the finite subnetwork $G_k^F$. 
    Since $R^F(x,y) = \lim_{k \to \iy} R_{G_k^F}(x,y)$ by \eqref{eqn:def:free-resistance}, it only remains to check the limit of $u_x^{(k)}$. Using a superscript to indicate the network in which the random walk travels, we have 
    \linenopax
    \begin{align}\label{eqn:prob-expr-for-vx:error}
      \lim_{k \to \iy} u_x^{(k)}(y)
      &= \lim_{k \to \iy} \prob_y^{G_k^F} [\gt_x < \gt_o]
       = \lim_{k \to \iy} \prob_y^{G} [\gt_x < \gt_o \,|\, \cpath \ci G_k^F]. 
    \end{align}
    Here again, the notation $[\cpath \ci H]$ denotes the event that the random walk never leaves the subnetwork $H$, i.e., $\gt_{H^\complm} = \iy$. The events $[\cpath \ci G_k^F]$ are nested and increasing, so the limit is the union, and \eqref{eqn:prob-expr-for-vx} follows. 
    Note that $G_k^F$ is recurrent, so $\cpath \ci G_k^F$ implies $\gt_x < \iy$.
  \end{proof}
\end{conj}

\begin{remark}\label{rem:error-in-prob-expr-for-vx}
  As indicated, the argument outlined above is incomplete due to the second equality of \eqref{eqn:prob-expr-for-vx:error}. While the set of paths from $y$ to $x$ in $G_k^F$ is the same as the set of paths from $y$ to $x$ in $G$ which lie in $G_k$, the probability of a given path may differ when computed in network or the other. This happens precisely when \cpath passes through a boundary point: the transition probability away from a point in $\bd G_k$ is strictly larger in $G_k^F$ than it is in $G_k$.
\end{remark}

\section{Comparison of resistance metric to other metrics}
\label{sec:Comparison-to-other-metrics}

\subsection{Comparison to geodesic metric}
\label{sec:Comparison-to-shortest-path-metric}

On a Riemannian manifold $(\gW,g)$, the geodesic distance is
  \linenopax
  \begin{equation*}
    \distmin(x,y)
    := \inf_{\gg} \left\{\int_0^1 g(\gg'(t),\gg'(t))^{1/2} \,dt \suth \gg(0)=x, \gg(1)=y, \gg \in C^1\right\}.
  \end{equation*}

\begin{defn}\label{def:geodesic-metric}
  On $(\Graph, \cond)$, the \emph{geodesic distance} from $x$ to $y$ is
  \linenopax
  \begin{equation}\label{eqn:def:geodesic-metric}
    \distmin(x,y)
    := \inf \{r(\cpath) \suth \cpath \in \Paths(x,y)\},
  \end{equation}
  where $r(\cpath) := \left(\sum_{(x,y) \in \cpath} \cond_{xy}^{-1} \right)$ (for resistors in series, the total resistance is the sum).
\end{defn}
  \glossary{name={$\distmin(x,y)$},description={geodesic distance, shortest-path metric},sort=d,format=textbf}
  \glossary{name={$r(\cpath)$},description={total resistance of a path},sort=r,format=textbf}

\begin{remark}\label{rem:shortest-path-metric}
  Definition~\ref{def:geodesic-metric} differs from the definition of shortest path metric found in the literature on general graph theory; without weights on the edges one usually defines the shortest path metric simply as the minimal number of edges in a path from $x$ to $y$. (This corresponds to taking $\cond \equiv 1$.) Such shortest paths always exist. According to Definition~\ref{def:geodesic-metric}, shortest paths may not exist (cf.~Example~\ref{exm:one-sided-ladder-model}). Of course, even when they do exist, they are not always unique.

It should be observed that effective resistance is not a geodesic metric, in the usual sense of metric geometry;  it does not correspond to a length structure in the sense of \cite[\S2]{Burago}. 
\end{remark}

\begin{lemma}\label{thm:shortest-path-bounds-resistance-distance}
  The effective resistance is bounded above by the geodesic distance. More precisely, $R^F(x,y) \leq \distmin(x,y)$ with equality if and only if \Graph is a tree.
  \begin{proof}
    If there is a second path, then some positive amount of current will pass along it (i.e., there is a positive probability of getting to $y$ via this route). To make this precise, let $v = v_x - v_y$ and let $\cpath = (x=x_0,x_1,\dots,x_n=y)$ be any path from $x$ to $y$:
    \linenopax
    \begin{align*}
      R^F(x,y)^2
      = |v(x)-v(y)|^2
      \leq r(\cpath) \energy(v),
    \end{align*}
    by the exact same computation as in the proof of Lemma~\ref{thm:L_x-is-bounded}, but with $u=v$. The desired inequality then follows by dividing both sides by $\energy(v) = R^F(x,y)$. 
    
    The other claim follows by observing that trees are characterized by the property of having exactly one path \cpath between any $x$ and $y$ in \verts. By \eqref{eqn:def:R^F(x,y)-diss}, $R^F(x,y)$ can be found by computing the dissipation of the unit current which runs entirely along \cpath from $x$ to $y$. This means that $\curr(x_{i-1},x_i)=1$ on \cpath, and $\curr = 0$ elsewhere, so
    \linenopax
    \begin{align*}
      R^F(x,y)
      &= \diss(\curr)
      = \sum_{i=1}^n \frac1{\cond_{x_{i-1} x_i}} \curr(x_{i-1},x_i)^2 
      = \sum_{i=1}^n \frac1{\cond_{x_{i-1} x_i}}  
      = r(\cpath).
      \qedhere
    \end{align*}
  \end{proof}
\end{lemma}

This type of inequality is explicitly calculated in Example~\ref{exm:finite-cyclic-model}.

\begin{remark}\label{rem:shortcut-through-infinity}
  It is clear from the end of the proof of Lemma~\ref{thm:shortest-path-bounds-resistance-distance} that on a tree, $v_x-v_y$ is locally constant on the complement of the unique path from $x$ to $y$. However, this may not hold for $f_x-f_y$, where $f_x = \Pfin v_x$; see Example~\ref{exm:binary-tree:reproducing-kernel}. This is an example of how the wired resistance can ``cheat'' by considering currents which take a shortcut through infinity; compare \eqref{eqn:def:R^F(x,y)-diss} to \eqref{eqn:def:R^W(x,y)-diss}.
\end{remark}

\subsection{Comparison to Connes' metric}
\label{sec:Comparison-to-Connes-metric}

  The formulation of $R(x,y)$ given in \eqref{eqn:def:R(x,y)-Lap} may evoke Connes' maxim that a metric can be thought of as the inverse of a Dirac operator; cf.~\cite{Con94}. This does not appear to have a literal incarnation in the current context, but we do have the inequality of Lemma~\ref{thm:Connes'-metric} in the case when $\cond = \one$. In this formulation, $v \in \HE$ is considered as a multiplication operator defined on $u$ by
  \linenopax
  \begin{equation}\label{eqn:v-as-operator}
    (vu)(x) := v(x) u(x), \q \forall x \in \verts,
  \end{equation}
  and both $v$ and \Lap are considered as operators on $\ell^2(\verts \cap \dom \energy$. We use the commutator notation $[v,\Lap] := v\Lap-\Lap v$, and $\| [v,\Lap] \|$ is understood as the usual operator norm on $\ell^2$.

\begin{lemma}\label{thm:Connes'-metric}
  If $\cond = \one$, then for all $x,y \in \verts$ one has
  \linenopax
  \begin{align}\label{eqn:Connes'-metric}
    R(x,y) \leq \sup \{|v(x)-v(y)|^2 \suth \| [v,\Lap] \| \leq \sqrt2, v \in \dom \energy\}.
  \end{align}
  \begin{proof}
    We will compare \eqref{eqn:Connes'-metric} to \eqref{eqn:def:R(x,y)-sup}. Writing $M_v$ for multiplication by $v$, it is straightforward to compute from the definitions
    \linenopax
    \begin{align*}
      (M_v \Lap - \Lap M_v)u(x)
      &= \sum_{y \nbr x} (v(y)-v(x))u(y),
    \end{align*}
    so that the Schwarz inequality gives
    \linenopax
    \begin{align*}
      \|[M_v,\Lap]u\|_2^2
      &= \sum_{x \in \verts} \left|\sum_{y \nbr x} (v(y)-v(x))u(y)\right|^2 \\
      &\leq \sum_{x \in \verts} \left(\sum_{y \nbr x} |v(y)-v(x)|^2\right) \left(\sum_{y \nbr x} |u(y)|^2\right).
    \end{align*}
    By extending the sum of $|u(x)|^2$ to all $x \in \verts$ (an admittedly crude estimate), this gives $\|[v,\Lap]u\|_2^2 \leq 2 \|u\|_2^2 \energy(v)$, and hence $\|[v,\Lap]\|^2 \leq 2 \energy(v)$
  \end{proof}
\end{lemma}

\section{Generalized resistance metrics}
\label{sec:Generalized resistance metrics}

In this section, we describe a notion of effective resistance between probability measures, of which $R(x,y)$ (or $R^F$ and $R^W$) is a special case. This concept is closely related to the notion of total variation of measures, and hence is related to mixing times of Markov chains; cf. \cite[\S4.1]{LevPerWil08}. When the Markov chain is taken to be random walk on an ERN, the state space is just the vertices of \Graph.

\begin{defn}\label{def:resistance-distance-of-measures}
  Let \gm and \gn be two probability measures on \verts. Then the total variation distance between them is
  \linenopax
  \begin{align}\label{eqn:def:resistance-distance-of-measures}
    \dist_{\TV}(\gm,\gn)
    := 2 \sup _{A \ci \verts} |\gm(A)-\gn(A)|.
  \end{align}
\end{defn}

\begin{prop}[{\cite[Prop.~4.5]{LevPerWil08}}]
  \label{thm:resistance-distance-of-measures-as-norm}
  Let \gm and \gn be two probability measures on the state space \gW of a (discrete) Markov chain. Then the total variation distance between them is
  \linenopax
  \begin{align}\label{eqn:resistance-distance-of-measures-as-norm}
    \dist_{\TV}(\gm,\gn)
    = \sup\left\{ \left|\sum_{x \in \gW} u(x) \gm(x) - \sum_{x \in \gW} u(x) \gn(x)\right| \suth \|u\|_\iy \leq 1 \right\}.
  \end{align}
  Here, $\|u\|_\iy := \sup_{x \in \verts}|u(x)|$.
\end{prop}
  \glossary{name={$\dist_{\TV}(\gm,\gn)$},description={total variation metric},sort=d,format=textbf}

\subsection{Effective resistance between measures}
\label{sec:Effective-resistance-between-measures}

If we think of \gm as a linear functional acting on the space of bounded functions, then it is clear that \eqref{eqn:resistance-distance-of-measures-as-norm} expresses $\dist_{\TV}(\gm,\gn)$ as the operator norm $\|\gm-\gn\|$. That is, it expresses the pairing between $\gm \in \ell^1$ and $u \in \ell^\iy$. We can therefore extend $R^F$ directly (see \eqref{eqn:def:R^F(x,y)-S}--\eqref{eqn:def:R^F(x,y)-sup} and Remark~\ref{rem:R(x,y) = |L{xy}|}).

\begin{defn}
  The free resistance between two probability measures is
  \linenopax
  \begin{align}
    \dist_{R^F}(\gm,\gn) 
    := \sup\left\{ \left|\sum_{x \in \verts} u(x) \gm(x) - \sum_{x \in \verts} u(x) \gn(x)\right|^2 \suth \|u\|_\energy \leq 1 \right\}.
  \end{align} 
\end{defn}
  \glossary{name={$\dist_{R^F}(\gm,\gn)$},description={free resistance between two probability measures},sort=d,format=textbf}

It is clear from this definition (and Remark~\ref{rem:R(x,y) = |L{xy}|}) that $R^F(x,y) = \dist_{R^F}(\gd_x,\gd_y)$. This extension of $R^F$ to measures was motivated by a question of Marc Rieffel in \cite{Rieffel99}.

\subsection{Total variation spaces}
\label{sec:Total-variation-spaces}

\begin{defn}\label{def:TV-form}
  Since $\dom \energy$ is a Banach space, we may define a new pairing via the bilinear form 
  \linenopax
  \begin{align}\label{eqn:def:TV-form}
    \la u, \gm\ra_{\TV}
    := \sum_{x \in \verts} u(x) \gm(x),
  \end{align}
  where \gm is an element of 
  \linenopax
  \begin{align}\label{eqn:def:TV}
    \TV 
    := \{\gm:\verts \to \bR \suth \exists k_\gm \text{ s.t.} 
      |\la u, \gm\ra_{\TV}| \leq k_\gm \cdot \energy(u)^{1/2}, 
      \forall u \in \dom \energy\}.
  \end{align}
  \glossary{name={$\TV$},description={total variation space},sort=T,format=textbf}
  Then $\TV= \dom \la u, \cdot\ra_\TV$ is the dual of $\dom \energy$ with respect to the total variation topology induced by \eqref{eqn:def:TV-form}. Also, the norm in \TV is given by
  \linenopax
  \begin{align}\label{eqn:def:TV-norm}
    \|\gm\|_\TV
    := \inf \{k \suth |\la u, \gm\ra_{\TV}| \leq k \cdot \energy(u)^{1/2}, \forall u \in \dom \energy\}.
  \end{align}
\end{defn}

\begin{remark}
  Since \TV is a Banach space which is the dual of a normed space, the unit ball
  \linenopax
  \begin{align}\label{eqn:unit-ball-in-TV}
    \{\gm \in \TV \suth \|\gm\|_{\TV} \leq 1\}
  \end{align}
  is compact in the weak-$\star$ topology, by Alaoglu's theorem.
\end{remark}

\begin{lemma}\label{thm:Lap-maps-Fin-to-TV}
  The Laplacian \Lap maps \HE into \TV with $\|\Lap v\|_\TV \leq \|\Pfin v\|_\energy$.
  \begin{proof}
    For $u,v \in \HE$, write $v=f+h$ with $f=\Pfin v$ and $h=\Phar v$, so that 
    \linenopax
    \begin{align*}
      \left|\sum_{x \in \verts} \cj{u}(x) \Lap v(x)\right|
      &\leq \left|\sum_{x \in \verts} \cj{u}(x) \Lap f(x)\right| 
        + \cancel{\left|\sum_{x \in \verts} \cj{u}(x) \Lap h(x)\right|}
      = |\la u, f\ra_\energy|
      \leq \|u\|_\energy \|f\|_\energy,
    \end{align*}
    by Theorem~\ref{thm:E(u,v)=<u,Lapv>-on-Fin} followed by the Schwarz inequality. The mapping is contractive relative to the respective norms because $\|v\|_\energy$ is an element of the set on the right-hand side of \eqref{eqn:def:TV-norm}, and hence at least as big as the infimum, whence $\|\Lap v\|_\TV \leq \|f\|_\energy \leq \|v\|_\energy$.
  \end{proof}
\end{lemma}


\section{Remarks and references}
\label{sec:Remarks-and-References-resistance-metrics}

A key reference for this chapter is \cite{Kig03}; the relationship between the free and wired resistance can be elegantly phrased in terms of resistance forms, as we describe in the following remark. Additionally, the role of the trace resistance is apparent in Kigami's work \cite{Kig01,Kig03, Kigami95, Kigami94, Kigami93}. However, the potential-theoretic approach can be intimidating to the uninitiated, and we hope that our treatment of effective resistance from the first principles of Hilbert space theory will provide a gentle introduction, as well as new insights. As an example of this, we feel that the probabilistic proof of Theorem~\ref{thm:edges-of-H^S} (to which we are indebted to Jun Kigami) offers insight as to \emph{why} the operation of Schur complement should correspond to taking the trace.

After PowersÕ papers in the 70s (starting with \cite{Pow76b}), there has been an explosive in the interest in metric geometry, potential theory \cite{Brelot},spectral theory \cite{Chu04}, and harmonic analysis \cite{Car73a} on infinite graphs. As illustrated in \cite{Kig03} a good deal of the motivation arose from a parallel research track dealing with analysis of fractals. In addition, some of the early work was motivated by problems in statistical mechanics (see e.g., \cite{Rue69} and \cite{Rue04}, on thermodynamic formalism).

\begin{remark}[Comparison with resistance forms]
  \label{rem:comparison-to-resistance-forms}
  In \cite[Def.~2.8]{Kig03}, a \emph{resistance form} is defined as follows: let $X$ be a set and let \energy be a symmetric quadratic form on $\ell(X)$, the space of all functions on $X$, and let \sF denote the domain of \energy. Then $(\energy, \sF)$ is a resistance form iff:
  \begin{enumerate}[(RF1)]
    \item \sF is a linear subspace of $\ell(X)$ containing the constant functions and \energy is nonnegative on \sF with $\energy(u) = 0$ iff $u$ is constant.
    \item $\sF/\sim$ is a Hilbert with inner product \energy, where $\sim$ is the equivalence relation defined on \sF by $u \sim v$ iff $u-v$ is constant.
    \item For any finite subset $V \ci U$ and for any $v \in \ell(V)$, there is $u \in \sF$ such that $u \evald{V} = v$.
    \item For any $p,q \in X$, the number
      \begin{align}\label{eqn:form-resistance}
        R_{\energy,\sF}(p,q) := \sup \left\{ \tfrac{|u(p)-u(q)|^2}{\energy(u)} 
          \suth u \in \sF, \energy(u) > 0\right\}
      \end{align}
      is finite. Then $R_{\energy,\sF}$ is called the \emph{effective resistance associated to the form} $(\energy,\sF)$.
    \item If $u \in \sF$, then $\cj{u}$ defined by $\cj{u}(x) := \min\{1,\max\{0,u(x)\}\}$ (the unit normal contraction of $u$, in the language of Dirichlet forms) is also in \sF.
  \end{enumerate}
  Upon comparison of \eqref{eqn:def:R^F(x,y)-R}--\eqref{eqn:def:R^F(x,y)-S} to \eqref{eqn:def:R^W(x,y)-R}--\eqref{eqn:def:R^W(x,y)-sup}, one can see that $R^F$ is the effective resistance associated to the resistance form $(\energy,\HE)$, and that $R^W$ is the effective resistance associated to the resistance form $(\energy,\Fin)$. We are grateful to Jun Kigami for pointing this out to us. Note that the wired resistance is not related to the ``shorted resistance form'' of \cite[\S3]{Kig03} (see Prop.~3.6 in particular). See also Remark~\ref{rem:Rharm-has-no-form}.
\end{remark}

The reader will also find \cite{Soardi94} to be an good reference for effective resistance. While the sources \cite{Soardi94, DoSn84, Lyons:ProbOnTrees, LevPerWil08, Peres99} do not especially emphasize the metric aspect of effective resistance, they provide an exceptional description of the relationship between effective resistance and random walks. The books \cite{Kig01} and \cite{Str06} are also useful for understanding connections between effective resistance and the energy form and Laplace operator, on graphs and on self-similar fractals.

For different formulations effective resistance appearing in the literature, see \cite{Pow76b} and \cite[\S8]{Peres99} for \eqref{eqn:def:R(x,y)-Lap}, \cite{DoSn84} for \eqref{eqn:def:R(x,y)-energy}, \cite{DoSn84,Pow76b} for \eqref{eqn:def:R(x,y)-diss}, \cite{Kig03,Kig01,Str06} for \eqref{eqn:def:R(x,y)-R}--\eqref{eqn:def:R(x,y)-S}. 

For investigations of the ``limit current'' (corresponding to free resistance) and ``minimal current'' (corresponding to wired resistance), one should consult \cite{Soardi94} and the earlier sources \cite{Flanders71, Thomassen90, Zemanian76}.

The role of effective resistance in combinatorics (dimer configurations, percolation on finite sets, etc.) is discussed in \cite{BeKo05, Rieffel99, KeWi09}. The role of Schur complement in the trace of a resistance form appears in \cite{Kig03}, and less specifically also in \cite{Metz, KeWi09, Butler, AndersonTrapp, Anderson}, where it is sometimes called the \emph{shorted operator}. Also see \cite{KndZuLiRo} for the role of Schur complement in regard to Dirichlet-to-Neumann maps.

%% file: construction-of-HE.tex

\chapter[Schoenberg-von Neumann construction of \HE]{Schoenberg-von Neumann construction of the energy space \HE}
\label{sec:Construction-of-HE}

\headerquote{If people do not believe that mathematics is simple, it is only because they do not realize how complicated life is.}{---~John~von~Neumann}

Studying the geometry of state space $X$ through vector spaces of functions on $X$ is a fundamental idea and variations of it can be traced back in several areas of mathematics. In the setting of Hilbert space, it originates with a suggestion of B. O. Koopman \cite{Koo27,Koo36,Koo57} in the early days of ``modern'' dynamical systems, ergodic theory, and the systematic study of representations of groups. A separate impetus in 1932 were the two ergodic theorems, the $L^2$ variant due to von Neumann \cite{vN32c} and the pointwise variant due to G. D. Birkhoff. While Birkhoff's version is deeper, von Neumann's version really started a whole trend: mathematical physics, quantization \cite{vN32c}, and operator theory; especially the use of the adjoint operator and the deficiency indices which we find useful in \S\ref{sec:the-boundary}--\S\ref{sec:Lap-on-HE}; cf.~\cite{vN32a,vN32b}. Further, there is an interplay between Hilbert space on the one side, and pointwise results in function theory on the other: In fact, the $L^2$-mean ergodic theorem of von Neumann is really is a corollary to the spectral theorem in its deeper version (spectral resolution via projection-valued measures) as developed in by M. H. Stone and J. von Neumann in the period 1928-1932; cf.~\cite{vN32b} and \cite[Ch. 2]{Arveson:invitation-to-Cstar}. This legacy motivates the material in this section, as well as our overall approach.

\section{Schoenberg and von Neumann's embedding theorem}
\label{sec:vonNeumann's-embedding-thm}

In Theorem~\ref{thm:R^F-embed-ERN-in-Hilbert} we show that an resistance network equipped with resistance metric may be embedded in a Hilbert space in such a way that $R$ is induced from the inner product of the Hilbert space. As a consequence, we obtain an alternative and independent construction of the Hilbert space \HE of finite-energy functions. This provides further justification for \HE as the \emph{natural} Hilbert space for studying the metric space $(\Graph,R^F) = ((\Graph,\cond),R^F)$ and \Fin as the natural Hilbert space for studying the metric space $(\Graph,R^W)$. Although we will be interested in both $(\Graph, R^F)$ and $(\Graph, R^W)$, for brevity, we sometimes refer to both as $(\Graph, R)$ when the distinction is not important. 
  \glossary{name={$(\Graph,R)$},description={the network $(\Graph,\cond)$ as a metric space under $R$},sort=G}

It is a natural question to ask whether or not a metric space may be naturally represented as a Hilbert space, and von Neumann proved a general result which provides an answer. The reader may wish to consult \S\ref{sec:von-Neumann's-embedding-theorem} for the statement of this result (Theorem~\ref{thm:vonNeumann's-embedding-thm}) in the form it is applied below, as well as the relevant definitions. We apply Theorem~\ref{thm:vonNeumann's-embedding-thm} to the metric space $(\Graph,R)$ and to obtain a Hilbert space \sH and a natural embedding $(\Graph,R) \to \sH$. It turns out that the Hilbert space is \HE when the embedding is applied with $R=R^F$ and \Fin when applied with $R=R^W$; see Remark~\ref{rem:uniqueness-of-HE}! Therefore, the Hilbert space \HE is the natural place to study $(\Graph,R)$. The reader may find the references \cite{vN32a,Ber84,Ber96,Schoe38a} to be helpful; see also Theorem~\ref{thm:Schoenberg's-Thm}.

The following theorem is inspired by the work of von Neumann and Schoenberg \cite{Schoe38b,Ber84}, but is a completely new result. One aspect of this result that contrasts sharply with the classical theory is that the embedding is applied to the metric $R^{1/2}$ instead of $R$, for each of $R=R^F$ and $R=R^W$.

\begin{theorem}\label{thm:R^F-embed-ERN-in-Hilbert}
  $(\Graph,R^F)$ may be isometrically embedded in a Hilbert space.
  \begin{proof}
    According to Theorem~\ref{thm:vonNeumann's-embedding-thm}, we need only to check that $R^F$ is negative semidefinite (see Definition~\ref{def:negative-semidefinite}). Let $f:\verts \to \bR$ \saty $\sum_{x \in \verts} f(x) = 0$. We must show
that $\sum_{x,y \in F} \cj{f(x)} R^F(x,y) f(y) \leq 0$, for any finite subset $F \ci \verts$. From \eqref{eqn:def:R^F(x,y)-energy}, we have
    \linenopax
    \begin{align*}
      \sum_{x,y \in F} \cj{f(x)} R^F(x,y) f(y)
      &= \sum_{x,y \in F} \cj{f(x)} \energy(v_x - v_y) f(y) \\
      &= \sum_{x,y \in F} \cj{f(x)} \energy(v_x) f(y) - 2 \sum_{x,y \in F} \cj{f(x)} \la v_x, v_y \ra_\energy f(y) \\
          &\hstr[10]+ \sum_{x,y \in F} \cj{f(x)} \energy(v_x) f(y) \\
      &= -2\left\la \sum_{x \in F} f(x) v_x, \sum_{y \in F} f(y) v_y \right\ra_\energy \\
      &= -2\left\|\sum_{x \in F} f(x) v_x \right\|_\energy^2
      \leq 0.
    \end{align*}
    For the second equality, note that the first two sums vanish by the assumption on $f$.
  \end{proof}
\end{theorem}

\begin{cor}\label{thm:R^W-embed-ERN-in-Hilbert}
  $(\Graph,R^W)$ may be isometrically embedded in a Hilbert space.
  \begin{proof}
    Because the energy-minimizer in \eqref{eqn:def:R^W(x,y)-energy} is $f_x = \Pfin v_x$, we can repeat the proof of Theorem~\ref{thm:R^F-embed-ERN-in-Hilbert} with $f_x$ in place of $v_x$ to obtain the result.
  \end{proof}
\end{cor}

\begin{remark}\label{rem:uniqueness-of-HE}
  Since $R^F(x,y) =\|v_x-v_y\|_\energy^2$ by \eqref{eqn:def:R^F(x,y)-energy}, Theorem~\ref{thm:uniqueness-of-vNeu-embedding} shows that the embedded image of $(\Graph,R^F)$ is unitarily equivalent to the \energy-closure of $\spn\{v_x\}$, which is \HE. Similarly, $R^W(x,y) =\|f_x-f_y\|_\energy^2$, where $f_x := \Pfin v_x$, by \eqref{eqn:def:R^W(x,y)-energy}, whence the embedded image of $(\Graph,R^W)$ is unitarily equivalent to the \energy-closure of $\spn\{f_x\}$, which is \Fin.
\end{remark}

\begin{remark}\label{rem:norm-vs-quasinorm}
  One can choose any vertex $o \in \verts$ to act as the ``origin'' and it becomes the origin of the new Hilbert space during the construction outlined in \S\ref{sec:von-Neumann's-embedding-theorem}. As a quadratic form defined on the space of all functions $v:\verts \to \bC$, the energy is indefinite and hence allows one to define only a quasinorm. There are ways to deal with the fact that \energy does not ``see constant functions''. One possibility is to adjust the energy so as to obtain a true norm, as follows:
  \begin{equation}\label{eqn:energy-o}
    \energy^o(u,v) := \energy(u,v) + u(o)v(o).
  \end{equation}
  The corresponding quadratic form is immediately seen to be a norm; this approach is carried out in \cite{FOT94}, for example, and also occasionally in the work of Kigami. This is discussed in \S\ref{sec:grounded-energy-space} under the name ``grounded energy form'.

  We have instead elected to work ``modulo constants''; the kernel of \energy is the set of constant functions, and inspection of von Neumann's embedding theorem (cf.~\eqref{eqn:vNeu's-quotient-completion}) shows that it is precisely these functions which are ``modded out'' in von Neumann's construction. However, the constant functions resurface as multiples of the vacuum vector in the Fock space representation of \S\ref{sec:bdG-as-equivalence-classes-of-paths}.
\end{remark}

\section{\HE as an invariant of \Graph}
\label{sec:HE-as-an-invariant}
In this section, we show that \HE may be considered as an invariant of the underlying graph.

\begin{defn}\label{def:ERN-morphism}
  Let \Graph and $H$ be resistance networks with respective conductances $\cond^{\Graph}$ and $\cond^H $. A \emph{morphism of resistance networks} is a function $\gf:(\Graph, \cond^{\Graph}) \to (H, \cond^{H})$ between the vertices of the two underlying graphs for which
  \begin{equation}\label{eqn:ERN-morphism}
    \cond^H_{\gf(x) \gf(y)} = r \cond^{\Graph}_{xy}, \qq 0 < r < \iy,
  \end{equation}
  for some fixed $r$ and all $x,y \in \verts$.

  Two resistance networks are \emph{isomorphic} if there is a bijective morphism between them. 
  Note that this implies
  \begin{equation}\label{eqn:edges-in-H}
    \edges[H] = \{(\gf(x),\gf(y)) \suth (x,y) \in \edges\}.
  \end{equation}
\end{defn}

\begin{defn}\label{def:metric-space-morphism}
  A \emph{morphism of metric spaces} is a homothetic map, that is, an isometry composed with a dilation:
  \begin{equation}\label{eqn:metric-space-morphism}
    \gf:(X,d_X) \to (Y,d_Y),
    \qq d_Y(\gf(a),\gf(b)) = r d_X(a,b), \qq 0 < r < \iy,
  \end{equation}
  for some fixed $r$ and all $a,b \in X$. An \emph{isomorphism} is, of course, an invertible morphism.
\end{defn}

We allow for a scaling factor $r$ in each of the previous definitions because an isomorphism amounts to a relabeling, and rescaling is just a relabeling of lengths. More formally, an isomorphism in any category is an invertible mapping, and dilations are certainly invertible for $0 < r < \iy$.

\begin{theorem}\label{thm:R-is-a-functor}
  For each of $R=R^F,R^W$, there is a functor $\sR:(\Graph, \cond) \to ((\Graph, \cond),R)$ from the category of resistance networks to the category of metric spaces.
  \begin{proof}
    One must check that an isomorphism $\gf:(\Graph,\cond_{\Graph}) \to (H,\cond_H)$ of resistance networks induces an isomorphism of the corresponding metric spaces, so check that \gf preserves \energy. We use $x,y$ to denote vertices in \Graph and $s,t$ to denote vertices in $H$.
    \linenopax
    \begin{align}
      \la u \comp \gf, v \comp \gf\ra_\energy
      &= \sum_{x,y} \cond_{xy} (\cj{u \comp \gf}(x) - \cj{u \comp \gf}(y)) (v \comp \gf(x) - v \comp \gf(y)) \notag \\
      &= r^{-1} \sum_{x,y} \cond_{\gf(x)\gf(y)} (\cj{u(\gf(x))} - \cj{u(\gf(y))}) (v(\gf(x)) - v(\gf(y))) \notag \\
      &= r^{-1} \sum_{s,t} \cond_{st} (\cj{u}(s) - \cj{u}(t)) (v(s) - v(t)) \notag \\
      &= r^{-1} \la u \comp \gf, v \comp \gf\ra_\energy,
      \label{eqn:scaling-of-energy-under-isom}
    \end{align}
    where we can change to summing over $s,t$ because \gf is a bijection. Therefore, the reproducing kernel $\{v_x\}$ of $(\Graph,R^F)$ (or $\{\Pfin v_x\}$ of $(\Graph,R^W)$) is preserved, and hence so is the metric.
  \end{proof}
\end{theorem}

\begin{cor}\label{thm:scaling-of-Lap-under-isomorphism}
  If \gp is an isomorphism of resistance networks with scaling ratio $r$, then 
  \linenopax
  \begin{align}
    \Lap(v \comp \gf) = r^{-1} \Lap(v) \comp \gf. 
  \end{align}
  \begin{proof}
    Compute $\Lap(v \comp \gf)(x)$ exactly as in \eqref{eqn:scaling-of-energy-under-isom}. 
  \end{proof}
\end{cor}

\begin{cor}\label{thm:isomorphic-ERNs-induce-isomorphic-metric-spaces}
  An isomorphism $\gf:(\Graph, \cond^{\Graph}) \to (H, \cond^{H})$ of resistance networks induces an isomorphism of metric spaces (where the resistance networks are equipped with their respective effective resistance metrics).
\end{cor}

We use the notation $\clspn{S}$ to denote the closure of the span of a set of vectors $S$ in a Hilbert space, where the closure is taken with respect to the norm of the Hilbert space. The following theorem is just an application of Theorem~\ref{thm:uniqueness-of-vNeu-embedding} with the quadratic form $\tilde Q = \la \cdot\,,\,\cdot\ra_\sK$.

\begin{theorem}\label{thm:unitary-invariance-of-HE}
  If there is a Hilbert space $\sK = \clspn{k_x}$ for some $k:X \to \sK$ with $d(x,y)=\|k_x-k_y\|_\sK^2$, then there is a unique unitary isomorphism $U:\sH \to \sK$ and it is given by $U:\sum_x \gx_x w_x \mapsto \sum_x \gx_x k_x$.
\end{theorem}

\begin{remark}\label{rem:HE-is-an-invariant}
  Theorem~\ref{thm:unitary-invariance-of-HE} may be interpreted as the statement that there is a functor from the category of metric spaces (with negative semidefinite metrics) into the category of Hilbert spaces. However, we have avoided this formulation because the functor is not defined for the entire category of metric spaces.
  For us, it suffices to note that the composition is a functor from resistance networks to Hilbert spaces, so that $\HE = \HE(\Graph)$ is an invariant of \Graph.
\end{remark}

\begin{remark}\label{rem:first-quantization}
  To obtain a first quantization, one would need to prove that a contractive morphism between resistance networks induces a contraction between the corresponding Hilbert spaces. In other words,
  \linenopax
  \begin{align*}
    f:\Graph_1 \to \Graph_2
    \q\implies\q
    T_f:\HE(\Graph_1) \to \HE(\Graph_2)
  \end{align*}
  with $\|T_f v\|_\energy \leq \|v\|_\energy$ whenever $f$ is contractive. The authors are currently working on this endeavour in \cite{JoPea10c}. The second quantization is discussed in Remark~\ref{rem:Wiener-improves-Minlos}.
\end{remark}

\version{}{\marginpar{Add stuff about rough isometry \& rough equivalence, someday.}}

\section{Remarks and references}
\label{sec:Remarks-and-References-construction-of-HE}

Of the results in the literature of relevant to the present chapter, the references \cite{Ba04, Ale75, PaSc72, Schoe38a, Schoe38b} are especially relevant. Some of the cited references for this chapter are more specialized, but for prerequisite material (if needed), the reader may find the books \cite{Gui72} by Guichardet, \cite{Hida80} by Hida, and \cite{PaSc72} by Parthasarathy \& Schmidt especially relevant. Standard applications of a negative definite function include either the construction of an abstract Hilbert space \cite{vN55,vN32c} or else the construction of measures on a path space \cite{PaSc72,Minlos63}.

We use the terminology ``Schoenberg-von Neumann embedding'' to denote a set of general principles, both classical and modern:
\linenopax
\begin{align*}
  \xymatrix{
  \text{geometry and physics} \ar[d]^{\text{\cite{Koo57}}} \\ 
  \text{Hilbert space} \ar[d]^{\text{\cite{Nel73a}}} \\
  \text{$L^2$ space of random variables}
  }
\end{align*}

Some examples of the Schoenberg-von Neumann embedding include: Brownian motion \cite{Nelson64, Nel69}, second quantization and quantum fields \cite{Minlos63, Gro70, Hida80}, stochastic integrals \cite{Malliavin}, spin models \cite{Lig93}, quantum spin models \cite{Pow67, Pow76a, Pow76b}. 

See Chapter~\ref{sec:Magnetism-and-long-range-order} below for further details, especially Theorem~\ref{thm:Powers'-correlation-estimate}.

%% file: boundary.tex

\chapter{The boundary and boundary representation}
\label{sec:the-boundary}

\headerquote{Nature is an infinite sphere of which the center is everywhere
and the circumference nowhere.}{---~B.~Pascal}

Boundary theory is a well-established subject; the deep connections between harmonic analysis, probability, and potential theory have led to several notions of boundary; see the Remarks and References section at the end of this chapter. 

\section{Motivation and outline}
\label{sec:boundary-outline}

Recall the classical result of Poisson that gives a kernel $k:\gW \times \del \gW \to \bR$ from which a bounded harmonic function can be given via
  \begin{equation}\label{eqn:Poisson-bdy-repn}
    u(x) = \int_{\del \gW} u(y) k(x,dy),
    \qq y \in \del \gW.
  \end{equation}
The material of \S\ref{sec:the-boundary} is motivated by the following discrete analogue of the Poisson boundary representation of a harmonic function.
  \glossary{name={$k(x,dy)$},description={Poisson kernel},sort=k,format=textbf}

\begin{theorem}[Boundary sum representation of harmonic functions]
  \label{thm:boundary-repn-for-harmonic}
  For all $u \in \Harm$, and $h_x = \Phar v_x$,
  \linenopax
  \begin{equation}\label{eqn:boundary-repn-for-harmonic}
    u(x) = \sum_{\bd \Graph} u \dn{h_x} + u(o).
  \end{equation}
  \begin{proof}
    Recall from Lemma~\ref{thm:repkernels-for-Fin-and-Harm} that $\{h_x\}$ is a reproducing kernel for \Harm. Therefore, $u(x)-u(o) = \la h_x, u \ra_\energy = \cj{\la u, h_x \ra_\energy} = \sum_{\bd \Graph} u \dn{h_x}$ because $\sum_{\verts} u \Lap h_x = 0$. Note that $\cj{h_x}=h_x$ by Lemma~\ref{thm:vx-is-R-valued}.
  \end{proof}
\end{theorem}

Up to this point, the boundary sum in \eqref{eqn:boundary-repn-for-harmonic} has been understood only as a limit of sums. Comparison of \eqref{eqn:boundary-repn-for-harmonic} and \eqref {eqn:Poisson-bdy-repn} makes one optimistic that $\bd \Graph$ can be realized as some compact set which supports a ``measure'' $\dn{h_x}$, thus giving a nice representation of the boundary sum of \eqref{eqn:boundary-repn-for-harmonic} as an integral. In Corollary~\ref{thm:Boundary-integral-repn-for-harm}, we extend Theorem~\ref{thm:boundary-repn-for-harmonic} to such an integral representation.
  
 Boundary theory of harmonic functions can roughly be divided three ways: the bounded harmonic functions (Poisson theory), the nonnegative harmonic functions (Martin theory), and the finite-energy harmonic functions studied in the present book. While Poisson theory is a subset of Martin theory, the relationship between Martin theory and the study of \HE is more subtle. For example, there exist unbounded functions of finite energy; cf.~Example~\ref{exm:Unbounded-functions-of-finite-energy}. Some details are given in \cite[\S3.7]{Soardi94}.
  
  Whether the focus is on the harmonic functions which are bounded, nonnegative, or finite-energy, the goals of the associated boundary theory are essentially the same:
\begin{enumerate}[(1)]
  \item Compactify the original domain \sD by constructing/identifying a boundary $\bd \sD$. Then $\cj{\sD} = \sD \cup \bd \sD$, where the closure is with respect to some (hopefully natural) topology.
  \item Define a procedure for extending harmonic functions $u$ from \sD to $\bd \sD$. Except in the case of Poisson theory, this extension $\tilde u$ is typically a measure (or other linear functional) on $\bd \sD$; it may not be representable as a function.
  \item Obtain a kernel $\mathbbm{k}(x,\gb)$ defined on $\sD \times \bd \sD$ against which one can integrate the extension $\tilde u$ so as to recover the value of $u$ at a point in \sD:
  \begin{align*}
    u(x) = \int_{\bd \sD} \mathbbm{k}(x,\gb) \tilde u(d\gb),
    \q \forall x \in \sD,
  \end{align*}
  whenever $u$ is a harmonic functions of the given class.
\end{enumerate}
  
The difference between our boundary theory and that of Poisson and Martin is rooted in our focus on \HE rather than $\ell^2$: both of these classical theories concern harmonic functions with growth/decay restrictions. By contrast, provided they neither grow too wildly nor oscillate too wildly, elements of \HE may remain positive or even tend to infinity at infinity. See Example~\ref{exm:unbounded-harmonic} for a function $h \in \Harm$ which is unbounded.

Our boundary essentially consists of infinite paths which can be distinguished by monopoles, i.e., two paths are not equivalent iff there is a monopole $w$ with different limiting values along each path. It is an immediate consequence that recurrent networks have no boundary, and transient networks with no nontrivial harmonic functions have exactly one boundary point (corresponding to the fact that the monopole at $x$ is unique). In particular, the integer lattices $(\bZd,\one)$ each have 1 boundary point for $d \geq 3$ and 0 boundary points for $d=1,2$. In particular, the integer lattices $(\bZd,\one)$ each have 1 boundary point for $d \geq 3$ and 0 boundary points for $d=1,2$. In contrast, the Martin boundary of $(\bZ^d,\one)$ is homeomorphic to the unit sphere $S^{n-1}$ (where $S^0 = \{-1,1\}$), and each $(\bZ^d,\one)$ has only one graph ends (except for $(\bZ,\one)$, which has two); cf.~\cite[\S3.B]{PicWoess90}.

In our version of the program outlined above, we follow the steps in the order (2)-(3)-(1). A brief summary is given here; further introductory material and technical details appear at the beginning of each section: 
\begin{enumerate}
  \item[For (2),] we construct a Gel'fand triple $\Schw \ci \HE \ci \Schw'$ to extend   \glossary{name={$\Schw$},description={``Schwartz space'' of test functions (of rapid decay)},sort=S}
the energy form to a pairing $\la \cdot \,,\, \cdot \ra_\Wiener$ on $\Schw \times \Schw'$, and then use Ito integration to extend this new pairing to $\HE \times \Schw'$. This yields a suitable class of linear functionals \gx on \HE, and we can extend a function $u$ on \HE to $\tilde u$ on $\Schw'$ by duality, i.e., $\tilde u(\gx) := \la u, \gx\ra_\Wiener$. We need to expand the scope of enquiry to include $\Schw'$ because \HE will not be sufficient; no infinite-dimensional Hilbert space can support a \gs-finite probability measure, by a theorem of Nelson.
  \item[For (3),] we use the Wiener transform to isometrically embed \HE into $L^2(\Schw',\prob)$. Applying this isometry to the energy kernel $\{v_x\}$, we get a reproducing kernel $\mathbbm{k}(x,d \prob) := h_x d\prob$, where $h_x = \Phar v_x$ and $\prob$ is a version of Wiener measure. In fact, \prob is a Gaussian probability measure on $\Schw'$ whose support is disjoint from \Fin.
  \item[For (1),] we consider certain measures $\gm_x$, defined in terms of the kernel and the Wiener measure just introduced, which are supported on $\Schw'/\Fin$ and indexed by the vertices $x \in \verts$. Then elements of the boundary $\bd G$ correspond limits of sequences $\{\gm_{x_n}\}$ where $x_n \to \iy$, modulo a suitable equivalence relation. This is the content of \S\ref{sec:bdG-as-equivalence-classes-of-paths}.
\end{enumerate}
Items (2)--(3) are the content of \S\ref{sec:Gel'fand-triples-and-duality} and the main result is Theorem~\ref{thm:HE-isom-to-L2(S',P)} (and its corollaries). Due to the close relationship between the Laplacian and the random walk on a network, there are good intuitive reasons why one would expect stochastic integrals (by which we mean the Wiener transform) to be related to the boundary. ``Going to the boundary'' of $(\Graph,\cond)$ involves a suitable notion of limit, and it is a well-known principle that suitable limits of random walk yield Brownian motion realized in $L^2$-spaces of global measures (e.g., Wiener measure).

However, before this program can proceed, we need a suitable dense subspace $\Schw \ci \HE$ of ``test functions'' for the construction of a Gel'fand triple. The basic idea is to use the ``smooth functions'', that is, $u \in \HE$ for which $\Lap(\dots \Lap(u)) \in \HE$, for any number of applications of \Lap. Making this precise requires a certain amount of attention to technical details concerning the domain of \Lap, and this comprises \S\ref{sec:harmonic-functions-and-the-domain-of-Lap}. 
Caution: when studying an operator, an important subtlety is that ``the'' adjoint $\Lap^\ad$ depends on the choice of domain, i.e., the linear subspace $\dom(\Lap) \ci \sH$. We consider \Lap as an operator on a rather different Hilbert space, $\ell^2(\verts)$, in \S\ref{sec:L2-theory-of-Lap-and-Trans}.

   

Finally, in \S\ref{sec:boundary-form}, we examine the connection between the defect spaces of \Lap and $\bd G$ via the use of a boundary form akin to those of classical functional analysis.

The reader is directed to Appendix~\ref{sec:operator-theory} for a brief review of some of the pertinent ideas from operator theory; especially regarding the graph of an operator (Definition~\ref{def:graph-of-operator}) and von Neumann's theorem characterizing essential self-adjointness (Theorem~\ref{thm:essentially-self-adjointness-criterion}). Note: in several parts of this section, we use vector space ideas that are not so common when discussing Hilbert spaces; e.g. \emph{finite} linear span, and (not necessarily orthogonal) linear independence.

\begin{remark}\label{rem:comparison-to-path-space-measures}
  In \S\ref{sec:Probabilistic-interpretation} we will return to the three-way comparison of harmonic functions which are bounded, nonnegative, or finite-energy, but for a different purpose: the construction of measures on spaces of (infinite) paths in $(\Graph, \cond)$. In the case of bounded harmonic functions on $(\Graph, \cond)$, the associated probability space is derived directly as a space of infinite paths in \Graph, and the measure is constructed via the standard Kolmogorov consistency method. That is, as a projective limit constructed via cylinder sets. While the present construction is also implicitly in terms of cylinder sets (due to Minlos' framework), the reader will notice by comparison that the two probability measures and their support are quite different. As a result the respective kernels take different forms. However, both techniques yield a way to represent the values $h(x)$ for $h$ harmonic and $x \in \verts$ as an integral over ``the boundary''.
 
While Doob's martingale theory works well for harmonic functions in $L^\iy$ or $L^2$, the situation for \HE is different. The primary reason is that \HE is not immediatelly realizable as an $L^2$ space. A considerable advantage of our Gel'fand-Wiener-Ito construction is that \HE is isometrically embedded into $L^2(\Schw',\prob)$ in a particularly nice way: it corresponds to the polynomials of degree 1. See Remark~\ref{rem:Wiener-improves-Minlos}.

\end{remark}

Recall that Corollary~\ref{thm:vx-dense-in-HE} shows that $\spn\{v_x\}$ is dense in \HE and that $\{v_x\}$ is a reproducing kernel for \HE. Throughout \S\ref{sec:the-boundary}, we will implicitly be using the version of \Lap introduced in Definition~\ref{def:LapM}, which we now recall for convenience.

\begin{defn}\label{def:LapV-recalled}
  Let $\MP := \spn\{v_x\}_{x \in \verts}$ denote the vector space of \emph{finite} linear combinations of dipoles. Let \LapM be the closure of the Laplacian when taken to have the dense domain \MP. 
\end{defn}
\glossary{name={$V$},description={the span of the energy kernel, i.e., finite linear combinations of r$v_x$'s},sort=V}
\glossary{name={\LapM},description={the closure of the Laplacian when taken to have the dense domain \MP},sort=L}

Note that since \Lap agrees with \LapM pointwise, we may suppress reference to the domain for ease of notation. Recall from Corollary~\ref{thm:Lap-Hermitian-on-V} that \LapM is Hermitian and even semibounded on its domain.
We explore the properties of \LapM further, including its range, domain, and self-adjoint extensions, in \S\ref{sec:Lap-on-HE}.

\section{Gel'fand triples and duality}
\label{sec:Gel'fand-triples-and-duality}


According to the program outlined above, we would like to obtain a (probability) measure space to serve as the boundary of \Graph. It is shown in \cite{Gro67,Gro70,Minlos63} that no Hilbert space of functions \sH is sufficient to support a Gaussian measure \prob (i.e., it is not possible to have $0<\prob(\sH)<\iy$ for a \gs-finite measure). However, it \emph{is} possible to construct a \emph{Gel'fand triple} (also called a \emph{rigged Hilbert space}): a dense subspace $S$ of \sH with
\linenopax
\begin{equation}\label{eqn:Gel'fand-triple-intro}
  S \ci \sH \ci S',
\end{equation}
where $S$ is dense in \sH and $S'$ is the dual of $S$. Additionally, $S$ and $S'$ must also satisfy some technical conditions: $S$ is a Fr\'{e}chet space in its own right but realized as dense subspace in \sH, with density referring to the Hilbert norm in \sH. However, $S'$ is the dual of $S$ with respect to a Fr\'{e}chet topology defined via a specific sequence of seminorms. Finally, it is assumed that the inclusion mapping of $S$ into \sH is continuous in the respective topologies. It was Gel'fand's idea to formalize this construction abstractly using a system of nuclearity axioms \cite{GMS58, Minlos58, Minlos59}. Our presentation here is adapted from quantum mechanics and the goal is to realize $\bd \Graph$ as a subset of $S'$.

There is a concrete situation when the Gel'fand triple construction is especially natural: $\sH = L^2(\bR,dx)$ and $S$ is the \emph{Schwartz space} of functions of rapid decay. That is, each $f \in S$ is $C^\iy$ smooth functions which decays (along with all its derivatives) faster than any polynomial. In this case, $SÕ$ is the space of \emph{tempered distributions} and the seminorms defining the Fr\'{e}chet topology on $S$ are 
\linenopax
\begin{align*}
  p_m(f) := \sup \{|x^k f^{(n)}(x)| \suth x \in \bR, 0 \leq k,n \leq m\},
  \qq m=0,1,2,\dots,
\end{align*}
where $f^{(n)}$ is the \nth derivative of $f$. Then $S'$ is the dual of $S$ with respect to this Fre\'{e}chet topology. One can equivalently express $S$ as
\linenopax
\begin{align}\label{eqn:Schwartz-space-as-powers-of-Hamiltonian}
  S := \{f \in L^2(\bR) \suth (\tilde{P}^2 + \tilde{Q}^2)^n f \in L^2(\bR), \forall n\},
\end{align}
where $\tilde{P}$ and $\tilde{Q}$ are the Heisenberg operators discussed in Example~\ref{exm:banded-nonselfadjoint-matrix}. The operator $\tilde{P}^2 + \tilde{Q}^2$ is most often called the quantum mechanical Hamiltonian, but some others (e.g., Hida, Gross) would call it a Laplacian, and this perspective tightens the analogy with the present study. In this sense, \eqref{eqn:Schwartz-space-as-powers-of-Hamiltonian} could be rewritten $S := \dom \Lap^\iy$; compare to \eqref{eqn:Schwartz}.

  The duality between $S$ and $S'$ allows for the extension of the inner product on \sH to a pairing of $S$ and $S'$:
\linenopax
\begin{align*}
  \la \cdot,\cdot\ra_\sH:\sH \times \sH \to \bC
  \qq\text{to}\qq
  \la \cdot,\cdot\ra_\sH^\sim: S \times S' \to \bR.
\end{align*}
In other words, one obtains a Fourier-type duality restricted to $S$. Moreover, it is possible to construct a Gel'fand triple in such a way that $\prob(S')=1$ for a Gaussian probability measure \prob. When applied to $\sH=\HE$, the construction yields two main outcomes:
\begin{enumerate}
  \item The next best thing to a Fourier transform for an arbitrary graph.
  \item A concrete representation of \HE as an $L^2$ measure space $\HE \cong L^2(S',\prob)$.
\end{enumerate}
\version{}{This offers the exciting possibility for Fourier analysis of nonabelian groups whose Cayley graphs conform to our definition of an \ERN. The authors are currently plumbing the viability of this approach in another paper.}

As a prelude, we begin with Bochner's Theorem, which characterizes the Fourier transform of a positive finite Borel measure on the real line. The reader may find \cite{ReedSimonII} helpful for further information.
\begin{theorem}[Bochner]
  \label{thm:Bochner's-theorem}
  Let $G$ be a locally compact abelian group. Then there is a bijective correspondence $\sF:\sM(G) \to \sP\sD(\hat G)$, where $\sM(G)$ is the collection of measures on $G$, and $\sP\sD(\hat G)$ is the set of positive definite functions on the dual group of $G$. Moreover, this bijection is given by the Fourier transform
  \linenopax
  \begin{equation}\label{eqn:Bochner's-theorem}
    \sF:\gn \mapsto \gf_\gn
    \qq\text{by}\qq
    \gf_\gn(\gx) = \int_G e^{\ii \la \gx,x\ra} \,d\gn(x).
  \end{equation}
\end{theorem}
In our applications to the \ERN $(\bZd,\one)$ in \S\ref{sec:lattice-networks}, the underlying group structure allows us to apply the above version of Bochner's theorem. Specifically, in the context of group duality, Bochner's theorem characterizes the Fourier transform of a positive finite Borel measures; cf.~\cite{ReedSimonII,Ber96}.

For our representation of the energy Hilbert space \HE in the case of general \ERN, we will need Minlos' generalization of Bochner's theorem from \cite{Minlos63, Sch73}. This important result states that a cylindrical measure on the dual of a nuclear space is a Radon measure iff its Fourier transform is continuous. In this context, however, the notion of Fourier transform is infinite-dimensional, and is dealt with by the introduction of Gel'fand triples \cite{Lee96}.

\begin{theorem}[Minlos]
  \label{thm:Minlos'-theorem}
  Given a Gel'fand triple $S \ci \sH \ci S'$, Bochner's Theorem may be extended to yield a bijective correspondence between the positive definite functions on $S$ and the Radon probability measures on $S'$. Moreover, in a specific case, this correspondence is uniquely determined by the identity
  \linenopax
  \begin{align}\label{eqn:Minlos-identity}
    \int_{S'} e^{\ii \la u, \gx\ra_{\tilde \sH}} \,d\prob(\gx) 
    = e^{-\frac12 \la u,u\ra_\sH},
  \end{align}
  where $\la \cdot , \cdot \ra_\sH$ is the original inner product on \sH and $\la \cdot , \cdot \ra_{\tilde \sH}$ is its extension to the pairing on $S \times S'$. 
\end{theorem}

Formula \eqref{eqn:Minlos-identity} may be interpreted as defining the Fourier transform of \prob; the function on the right-hand side is positive definite and plays a special role in stochastic integration, and its use in quantization.

\subsection{A space of test functions \Schw on \Graph}
\label{sec:test-functions-Schw}

To apply Minlos' Theorem in the context of $(\Graph, \cond)$, we first need to construct a Gel'fand triple for \HE; we begin by identifying a certain subspace of the domain of \LapM. Recall from Definition~\ref{def:LapM} that $\MP := \spn\{v_x,\monov,\monof\}_{x \in \verts}$.

\begin{defn}\label{def:extn-of-Lap}
  Let \LapS be a self-adjoint extension of \LapM; since \LapM is Hermitian and commutes with conjugation (since \cond is \bR-valued), a theorem of von Neumann's states that such an extension exists.
\glossary{name={\LapS},description={a self-adjoint extension of \LapM},sort=L,format=textbf}
  
  Let $\LapS^p u := (\LapS\LapS\dots\LapS)u$ be the $p$-fold product of \LapS applied to $u \in \HE$. Define $\dom(\LapS^p)$ inductively by
  \linenopax
  \begin{equation}\label{eqn:domE(LapEp)}
    \dom(\LapS^p) := \{u \suth \LapS^{p-1} u \in \dom(\LapS)\}.
  \end{equation}
\end{defn}

\begin{defn}\label{def:Schw-on-G}
  The \emph{(Schwartz) space of functions of rapid decay} is 
  \linenopax
  \begin{equation}\label{eqn:Schwartz}
    \Schw := \dom(\LapS^\iy), 
  \end{equation}
\glossary{name={$\Schw$},description={``Schwartz space'' of test functions (of rapid decay)},sort=S,format=textbf}
  where $\dom(\LapS^\iy) := \bigcap_{p=1}^\iy \dom(\LapS^p)$ consists of all \bR-valued functions $u \in \HE$ for which $\LapS^p u \in \HE$ for any $p$. The space of \emph{Schwartz distributions} or \emph{tempered distributions} is the dual space $\Schw'$ of \bR-valued continuous linear functionals on \Schw.
\end{defn}

\begin{remark}\label{rem:need-for-self-adjoint-extension-of-LapV-for-Minlos}
A good choice of self-adjoint extension in Definition~\ref{def:extn-of-Lap} is the operator \LapH discussed in \S\ref{sec:harmonic-functions-and-the-domain-of-Lap}. It is critical to make the unusual step of taking a self-adjoint extension of \LapM for several reasons. Most importantly, we will need to apply the spectral theorem to extend the energy inner product $\la \cdot,\cdot\ra_\energy$ to a pairing on $\Schw \times \Schw'$. In fact, it will turn out that for $u \in \Schw, v \in \Schw'$, the extended pairing is given by $\la u,v\ra_\Wiener = \la \LapS^p u, \LapS^{-p} v \ra_\energy$, where $p$ is any integer large enough to ensure $\LapS^p u, \LapS^{-p} v \in \HE$. This relies crucially on the self-adjointness of the operator appearing on the right-hand side. Moreover, without self-adjointness, we would be unable to prove that \Schw is dense in \HE; see Lemma~\ref{thm:Schw-dense-in-HE}.
    
Additionally, the self-adjoint extensions of \LapM are in bijective correspondence with the isotropic subspaces of $\dom(\LapM^\ad)$, and we will see that these are useful for understanding the boundary of \Graph in terms of defect; see Theorem~\ref{thm:LapV-not-ess-selfadjoint-iff-Harm=0}. Recall that a subspace $\sM \ci \dom(\LapM^\ad)$ is \emph{isotropic} iff $\bdform(u,v) = 0$, $\forall u,v \in \sM$, where \bdform is as in Definition~\ref{def:boundary-form}. Since $\dom(\LapM)$ is isotropic (cf.~Theorem~\ref{thm:boundary-form-vanishes-on-dom(LapV)}), we think of \sM as a subspace of the quotient (boundary) space $B = \dom(\LapM^\ad) / \dom(\LapM)$. 
\end{remark}

\begin{remark}\label{rem:Schwartz-space-is-real-valued}
  Note that \Schw and $\Schw'$ consist of \bR-valued functions. This technical detail is important because we do not expect the integral $\int_{S'} e^{\ii \la u, \cdot \ra_{\tilde \Wiener}} \,d\prob$ from \eqref{eqn:Minlos-identity} to converge unless it is certain that $\la u, \cdot \ra$ is \bR-valued. This is the reason for the last conclusion of Lemma~\ref{thm:energy-extends}.
\end{remark}

\begin{remark}\label{rem:Schw-is-Frechet}
  Note that \Schw is dense in $\dom(\LapS)$ with respect to the graph norm, by standard spectral theory. For each $p \in \bN$, there is a seminorm on \Schw defined by
  \linenopax
  \begin{equation}\label{eqn:p-norm-on-Schw}
    \|u\|_p := \| \LapS^p u\|_\energy.
  \end{equation}
  Since $(\dom \LapS^p, \|\cdot\|_p)$ is a Hilbert space for each $p \in \bN$, the subspace \Schw is a Fr\'{e}chet space.
\end{remark}

\begin{defn}\label{def:spectral-truncation}
  Let $\charfn{[a,b]}$ denote the usual indicator function of the interval $[a,b] \ci \bR$, and let \spectrans be the spectral transform in the spectral representation of \LapS, and let $E$ be the associated projection-valued measure. Then define $E_n$ to be the \emph{spectral truncation operator} acting on \HE by
  \linenopax
  \begin{align*}
     E_n u 
     := \spectrans^\ad \charfn{[\frac1n,n]} \spectrans u 
     = \int_{1/n}^n E(dt)u.
  \end{align*}
\end{defn}


\begin{lemma}\label{thm:Schw-dense-in-HE}
  \Schw is a dense analytic subspace of \HE (with respect to \energy), and so $\Schw \ci \HE \ci \Schw'$ is a Gel'fand triple.
  \begin{proof}
    This essentially follows immediately once it is clear that $E_n$ maps \HE into \Schw. For $u \in \HE$, and for any $p=1,2,\dots$,
    \linenopax
    \begin{equation}\label{eqn:E-norm-of-Lapp(spectral-truncation)}
      \|\LapS^{p} E_n u\|_\energy^2 
      = \int_{1/n}^n \gl^{2p} \|E(d\gl)u\|_\energy^2
      \leq n^{2p} \|u\|_\energy^2,
    \end{equation}
    So $E_n u \in \Schw$. It follows that $\|u-E_n u\|_\energy \to 0$ by standard spectral theory.
  \end{proof}
\end{lemma}

\begin{theorem}\label{thm:Wiener-product-as-Lap(p)-powers}
  The energy form $\la \cdot,\cdot\ra_\energy$ extends to a pairing on $\Schw \times \Schw'$ defined by
 \linenopax
  \begin{equation}\label{eqn:energy-extends-by-ps}
    \la u,v\ra_\Wiener := \la \LapS^p u, \LapS^{-p} v \ra_\energy,
  \end{equation}
  where $p$ is any integer such that $|v(u)| \leq K \|\Lap^p u\|_\energy$ for all $u \in \Schw$.
  \glossary{name={$\la\cdot,\cdot\ra_\Wiener$},description={(Wiener) inner product on $L^2(\Schw',\prob)$},sort=<,format=textbf}
  \begin{proof}
   If $v \in \Schw'$, then there is a $C$ and $p$ such that $|\la s, v\ra_\Wiener| \leq C\|\LapS^p s\|_\energy$ for all $s \in \Schw$. Set $\gf(\LapS^p s) := \la s, v\ra_\Wiener$ to obtain a continuous linear functional on \HE (after extending to the orthogonal complement of $\spn\{\LapS^p s\}$ by 0 if necessary). Now Riesz's lemma gives a $w \in \HE$ for which $\la s, v\ra_\Wiener = \la \LapS^p s, w\ra_\energy$ for all $s \in \Schw$ and we define $\LapS^{-p} v := w \in \HE$ to make the meaning of the right-hand side of \eqref{eqn:energy-extends-by-ps} clear. 
  \end{proof}
\end{theorem}

\begin{lemma}\label{thm:energy-extends}
 The pairing on $\Schw \times \Schw'$ is equivalently given by
  \linenopax
  \begin{equation}\label{eqn:energy-extends-by-vn}
    \la u,\gx\ra_\Wiener = \lim_{n \to \iy} \gx(E_n u), 
  \end{equation}
  where the limit is taken in the topology of $\Schw'$. Moreover, $\tilde u(\gx) = \la u, \gx\ra_\Wiener$ is \bR-valued on $\Schw'$.
  \begin{proof}
    $E_n$ commutes with \LapS. This is a standard result in spectral theory, as $E_n$ and \LapS are unitarily equivalent to the two commuting operations of truncation and multiplication, respectively. Therefore,
    \linenopax
    \begin{align*}
      \gx(E_n u)
      = \la E_n u, \gx\ra_\Wiener
      = \la \LapS^p E_n s, \LapS^{-p} \gx \ra_\energy
      = \la E_n \LapS^p s, \LapS^{-p} \gx \ra_\energy
      = \la\LapS^p s, E_n \LapS^{-p} \gx \ra_\energy.
    \end{align*}
    Standard spectral theory also gives $E_n v \to v$ in \HE, so 
    \linenopax
    \begin{align*}
      \lim_{n \to \iy} \gx(E_n u)
      = \lim_{n \to \iy} \la\LapS^p s, E_n \LapS^{-p} \gx \ra_\energy
      = \la \LapS^p u,\LapS^{-p} v \ra_\energy.
    \end{align*}
    
    Note that the pairing $\la \cdot \,,\,\cdot \ra_\Wiener$ is a limit of real numbers, and hence is real.
  \end{proof}
\end{lemma}
    
    
\begin{cor}\label{thm:Spectruncation-extends-to-S'}
  $E_n$ extends to a mapping $\tilde E_n: \Schw' \to \HE$ defined via $\la u, \tilde E_n \gx \ra_\energy := \gx(E_n u)$. Thus, we have a pointwise extension of $\la \cdot \,,\,\cdot\ra_\Wiener$ to $\HE \times \Schw'$ given by
  \begin{equation}\label{eqn:energy-extends-by-truncating-on-S'}
    \la u,\gx\ra_\Wiener = \lim_{n \to \iy} \la u, \tilde E_n \gx \ra_\energy.
  \end{equation}
\end{cor}


\begin{lemma}\label{thm:vx-in-Schw}
  If $\deg(x)$ is finite for each $x \in \verts$, or if $\uBd < \iy$, then one has $v_x \in \Schw$.
  \begin{proof}
    This is immediate from the technical lemma, Lemma~\ref{thm:LapV-maps-V-to-V}, which we postpone for now.
  \end{proof}
\end{lemma}

\begin{remark}\label{rem:Lap-maps-span(vx)-into-itself}
  When the hypotheses of Lemma~\ref{thm:vx-in-Schw} are satisfied, note that $\spn\{v_x\}$ is dense in \Schw with respect to \energy, but \emph{not} with respect to the Frechet topology induced by the seminorms \eqref{eqn:p-norm-on-Schw}, nor with respect to the graph norm. One has the inclusions
  \begin{equation}\label{eqn:vx-graph-inclusions}
    \left\{\left[\begin{array}{c} v_x \\ \LapM v_x \end{array}\right]\right\}
    \ci \left\{\left[\begin{array}{c} s \\ \LapS s \end{array}\right]\right\}
    \ci \left\{\left[\begin{array}{c} u \\ \LapS u \end{array}\right]\right\}
  \end{equation}
  where $s \in \Schw$ and $u \in \HE$, with the second inclusion dense and the first inclusion not dense.
\end{remark}

We have now obtained a Gel'fand triple $\Schw \ci \HE \ci \Schw'$, and we are ready to apply the Minlos Theorem to a particularly lovely positive definite function on \Schw, in order that we may obtain a particularly nice measure on $\Schw'$. Recall that we constructed \HE from the resistance metric in \S\ref{sec:Construction-of-HE} by making use of negative definite functions. We now apply this to a famous result of Schoenberg which may be found in \cite{Ber84,ScWh49}.
\begin{theorem}[Schoenberg]
  \label{thm:Schoenberg's-Thm}
  Let $X$ be a set and let $Q: X \times X \to \bR$ be a function. Then the following are equivalent.
  \begin{enumerate}[(1)]
    \item $Q$ is negative definite.
    \item $\forall t \in \bR^+$, the function $p_t(x,y) := e^{-t Q(x,y)}$ is positive definite on $X \times X$.
    \item There exists a Hilbert space $\sH$ and a function $f:X \to \sH$ such that $Q(x,y) = \|f(x)-f(y)\|_\sH^2$.
  \end{enumerate}
\end{theorem}

In the proof of the following theorem, we apply Schoenberg's Theorem with $t=\frac12$ to the resistance metric in the form $R^F(x,y) = \|v_x-v_y\|_\energy^2$ from \eqref{eqn:def:R^F(x,y)-energy}.
The proof of Theorem~\ref{thm:HE-isom-to-L2(S',P)} also uses the notation $\Ex_\gx(f) := \int_{\Schw'} f(\gx) \,d\prob(\gx)$. 
\glossary{name={$\Ex_\gx(f)$},description={expectation of $f$ with dummy variable \gx; $\int f(\gx)\,d\gm(x)$},sort=E,format=textbf}

\begin{theorem}\label{thm:HE-isom-to-L2(S',P)}
  The Wiener transform $\sW:\HE \to L^2(\Schw',\prob)$ defined by
  \linenopax
  \begin{equation}\label{eqn:Gaussian-transform}
    \sW : v \mapsto \tilde v,
    \q \tilde v(\gx) = \la v, \gx\ra_\Wiener,
  \end{equation}
  is an isometry. The extended reproducing kernel $\{\tilde v_x\}_{x \in \verts}$ is a system of Gaussian random variables which gives the resistance distance by
  \linenopax
  \begin{equation}\label{eqn:R(x,y)-as-expectation}
    R^F(x,y) = \Ex_\gx((\tilde v_x - \tilde v_y)^2).
  \end{equation}
  Moreover, for any $u,v \in \HE$, the energy inner product extends directly as
  \linenopax
  \begin{equation}\label{eqn:Expectation-formula-for-energy-prod}
    \la u, v \ra_\energy
    = \Ex_\gx\left( \cj{\tilde{u}} \tilde{v} \right)
    = \int_{\Schw'} \cj{\tilde{u}} \tilde{v} \,d\prob.
  \end{equation}
  \begin{proof}
    Since $R^F(x,y)$ is negative semidefinite by the proof of Theorem~\ref{thm:R^F-embed-ERN-in-Hilbert}, we may apply Schoenberg's theorem and deduce that $\exp(-\tfrac12\|u-v\|_\energy^2)$ is a positive definite function on $\HE \times \HE$.     Consequently, an application of the Minlos correspondence to the Gel'fand triple established in Lemma~\ref{thm:Schw-dense-in-HE} yields a Gaussian probability measure \prob on $\Schw'$. 
    
    Moreover, \eqref{eqn:Minlos-identity} gives 
    \linenopax
    \begin{align}\label{eqn:Minlos-eqns}
      \Ex_\gx(e^{\ii\la u, \gx \ra_\Wiener}) = e^{-\frac12\|u\|_\energy^2},
    \end{align}
    whence one computes
    \linenopax
    \begin{align}\label{eqn:Minlos-expectation-integral}
      \int_{\Schw'} \left(1 + \ii\la u, \gx\ra_\Wiener - \frac12\la u, \gx\ra_\Wiener^2 + \cdots \right)\,d\prob(\gx) 
      = 1 - \frac12 \la u, u \ra_\energy + \cdots.
    \end{align}
    Now it follows that $\Ex(\tilde{u}^2) = \Ex_\gx(\la u, \gx\ra_\Wiener^2) = \|u\|_\energy^2$ for every $u \in \Schw$, by comparing the terms of \eqref{eqn:Minlos-expectation-integral} which are quadratic in $u$. Therefore, $\sW:\HE \to \Schw'$ is an isometry, and \eqref{eqn:Minlos-expectation-integral} gives
    \linenopax
    \begin{align}\label{eqn:Exp-tilde-vx=E(vx)}
      \Ex_\gx(|\tilde v_x - \tilde v_y|^2)
      = \Ex_\gx(\la v_x - v_y, \gx \ra^2)
      &= \|v_x - v_y\|_\energy^2,
    \end{align}
    whence \eqref{eqn:R(x,y)-as-expectation} follows from \eqref{eqn:def:R^F(x,y)-energy}. Note that by comparing the linear terms, \eqref{eqn:Minlos-expectation-integral} implies $\Ex_\gx(1) = 1$, so that \prob is a probability measure, and $\Ex_\gx(\la u,\gx\ra) = 0$ and $\Ex_\gx(\la u,\gx\ra^2) = \|u\|_\Wiener^2$, so that \prob is actually Gaussian.
    
    
    Finally, use polarization to compute
    \linenopax
    \begin{align*}
      \la u, v \ra_\energy
      &= \frac14 \left(\|u+v\|_\energy^2 - \|u-v\|_\energy^2\right) \\
      &= \frac14 \left(\Ex_\gx\left(\left|\tilde{u}+ \tilde{v}\right|^2\right) 
        - \Ex_\gx\left(\left|\tilde{u}-\tilde{v}\right|^2\right)\right) 
        &&\text{by \eqref{eqn:Exp-tilde-vx=E(vx)}} \\
      &= \frac14 \int_{\Schw'} \left|\tilde{u}+\tilde{v}\right|^2(\gx)
        - \left|\tilde{u}-\tilde{v}\right|^2(\gx) \,d\prob(\gx) \\
      &= \int_{\Schw'} \cj{\tilde{u}}(\gx) \tilde{v}(\gx) \,d\prob(\gx).
    \end{align*}
    This establishes \eqref{eqn:Expectation-formula-for-energy-prod} and completes the proof.
  \end{proof}
\end{theorem}

  It is important to note that since the Wiener transform $\sW : \Schw \to \Schw'$ is an isometry, the conclusion of Minlos' theorem is stronger than usual: the isometry allows the energy inner product to be extended isometrically to a pairing on $\HE \times \Schw'$ instead of just $\Schw \times \Schw'$.

\begin{remark}\label{rem:HE-into-L2(S',prob)-gives-hermitian-multiplication}
  With the embedding $\HE \to L^2(\Schw',\prob)$, we obtain a maximal abelian algebra of Hermitian multiplication operators $L^\iy(\Schw')$ acting on $L^2(\Schw',\prob)$. By contrast, see (ii) of Remark~\ref{rem:elements-of-HE-are-technically-equivalence-classes}. 
\end{remark}

\begin{remark}\label{rem:complexify-L2(S',P)}
  The reader will note that we have taken pains to keep everything \bR-valued in this chapter (especially the elements of \Schw and $\Schw'$), primarily to ensure the convergence of $\int_{\sS'} e^{\ii\la u, \gx \ra_\Wiener}\,d\prob(\gx)$ in \eqref{eqn:Minlos-eqns}. However, now that we have established the fundamental identity $\la u, v \ra_\energy = \int_{\sS'} \cj{\tilde{u}} \tilde{v} \,d\prob$ in
\eqref{eqn:Expectation-formula-for-energy-prod} and extended the pairing $\la \cdot , \cdot \ra_\Wiener$ to $\HE \times \Schw'$, we are at liberty to complexify our results via the standard decomposition into real and complex parts: $u = u_1 + \ii u_2$ with $u_i$ \bR-valued elements of \HE, etc.
\end{remark}

\begin{remark}\label{rem:Wiener-improves-Minlos}
  The polynomials are dense in $L^2(\Schw',\prob)$. More precisely, if $\gf(t_1, t_2, \dots, t_k)$ is an ordinary polynomial in $k$ variables, then 
  \linenopax
  \begin{align}\label{eqn:polynomials-in-S'}
  \gf(\gx) := \gf\left(\la u_1, \gx\ra_\Wiener, \la u_2, \gx\ra_\Wiener, \dots \vstr[2.2] \, \la u_k,\gx\ra_\Wiener\right)
  \end{align}
  is a polynomial on $\Schw'$ and 
  \linenopax
  \begin{align}\label{eqn:Poly(n)-in-S'}
    \Poly_n := \{\gf\left(\widetilde{u_1}(\gx), \widetilde{u_2}(\gx), \dots \vstr[2.2] \, \widetilde{u_k}(\gx)\right), \deg(\gf) \leq n, \suth u_j \in \HE, \gx \in \Schw'\}
  \end{align}
  is the collection of polynomials of degree at most $n$, and $\{\Poly_n\}_{n=0}^\iy$ is an increasing family whose union is all of $\Schw'$. One can see that the monomials $\la u, \gx\ra_\Wiener$ are in $L^2(\Schw',\prob)$ as follows: compare like powers of $u$ from either side of \eqref{eqn:Minlos-expectation-integral} to see that $\Ex_\gx\left(\la u, \gx\ra_\Wiener^{2n+1}\right) = 0$ and
  \linenopax
  \begin{align}\label{eqn:expectation-of-even-monomials}
    \Ex_\gx\left(\la u, \gx\ra_\Wiener^{2n}\right) 
    = \int_{\Schw'} |\la u, \gx\ra_\Wiener|^{2n} \, d\prob(\gx) 
    = \frac{(2n)!}{2^n n!} \|u\|_\energy^{2n}, 
  \end{align}
  and then apply the Schwarz inequality. 
  
  To see why the polynomials $\{\Poly_n\}_{n=0}^\iy$ should be dense in $L^2(\Schw',\prob)$ observe that the sequence $\{P_{\Poly_n}\}_{n=0}^\iy$ of orthogonal projections increases to the identity, and therefore, $\{P_{\Poly_n} \tilde u\}$ forms a martingale, for any $u \in \HE$ (i.e., for any $\tilde u \in L^2(\Schw',\prob)$).
  
  If we denote the ``multiple Wiener integral of degree $n$'' by 
  \linenopax
  \begin{align*}
    H_n := \Poly_n - \Poly_{n-1} = cl \spn\{\la u, \cdot\ra_\Wiener^n \suth u \in \HE\},
    \q n \geq 1,
  \end{align*}
  and $H_0 := \bC \one$ for a vector \one with $\|\one\|_2=1$.
  Then we have an orthogonal decomposition of the Hilbert space     
  \linenopax
  \begin{align}\label{eqn:Fock-repn-of-L2(S',P)}
    L^2(\Schw',\prob) = \bigoplus_{n=0}^\iy H_n.
  \end{align}
  See \cite[Thm.~4.1]{Hida80} for a more extensive discussion.
  
  \glossary{name={\one},description={the constant function 1, the vacuum vector},sort=1,format=textbf}
A physicist would call \eqref{eqn:Fock-repn-of-L2(S',P)} the Fock space representation of $L^2(\Schw',\prob)$ with ``vacuum vector'' \one; note that $H_n$ has a natural (symmetric) tensor product structure. Familiarity with these ideas is not necessary for the sequel, but the decomposition \eqref{eqn:Fock-repn-of-L2(S',P)} is helpful for understanding two key things:
  \begin{enumerate}[(i)]
    \item The Wiener isometry $\sW:\HE \to L^2(\Schw',\prob)$ identifies \HE with the subspace $H_1$ of $L^2(\Schw',\prob)$, in particular, $L^2(\Schw',\prob)$ is not isomorphic to \HE. In fact, it is the second quantization of \HE.
    \item The constant function \one is an element of $L^2(\Schw',\prob)$ but does not correspond to any element of \HE. In particular, constant functions in \HE are equivalent to 0, but this is not true in $L^2(\Schw',\prob)$.
  \end{enumerate}
  It is somewhat ironic that we began this story by removing the constants (via the introduction of \energy), only to reintroduce them with a certain amount of effort, much later. Item (ii) explains why it is not nonsense to write things like $\prob(\Schw') = \int_{\Schw'} \one \,d\prob = 1$, and will be helpful when discussing boundary elements in \S\ref{sec:bdG-as-equivalence-classes-of-paths}.
\end{remark}

\begin{cor}\label{thm:resistance-as-distributional-integral}
  For $e_x(\gx) := e^{\ii\la v_x, \gx\ra_\Wiener}$, one has $\Ex_\gx(e_x) = e^{-\frac12 R^F(o,x)}$ and hence
  \linenopax
  \begin{equation}\label{eqn:resistance-as-distributional-integral}
    \Ex_\gx(\cj{e_x} e_y) 
    = \int_{\Schw'} \cj{e_x(\gx)} e_y(\gx)\,d\prob
    = e^{-\frac12 R^F(x,y)}.
  \end{equation}
  \begin{proof}
    Substitute $u=v_x$ or $u=v_x-v_y$ in \eqref{eqn:Minlos-eqns} and apply Theorem~\ref{thm:free-resistance}.
  \end{proof}
\end{cor}

\begin{remark}
Remark~\ref{rem:resistance-as-path-integral} discusses the interpretation of the free resistance as the reciprocal of an integral over a path space; Corollary~\ref{thm:resistance-as-distributional-integral} provides a variation on this theme:
  \linenopax
  \begin{equation}\label{eqn:resistance-as-distributional-integral}
    R^F(x,y) = -2 \log \Ex_\gx(\cj{e_x} e_y) = -2 \log \int_{\Schw'} \cj{e_x(\gx)} e_y(\gx)\,d\prob.
  \end{equation}
  
  Observe that Theorem~\ref{thm:HE-isom-to-L2(S',P)} was carried out for the free resistance, but all the arguments go through equally well for the wired resistance; note that $R^W$ is similarly negative semidefinite by Theorem~\ref{thm:Schoenberg's-Thm} and Corollary~\ref{thm:R^W-embed-ERN-in-Hilbert}. Thus, there is a corresponding Wiener transform $\sW:\Fin \to L^2(\Schw',\prob)$ defined by
  \linenopax
  \begin{equation}\label{eqn:wired-Gaussian-transform}
    \sW : v \mapsto \tilde f,
    \qq f = \Pfin v \;\text{ and }\; \tilde f(\gx) = \la f, \gx\ra_\Wiener.
  \end{equation}
  Again, $\{\tilde f_x\}_{x \in \verts}$ is a system of Gaussian random variables which gives the wired resistance distance by $R^W(x,y) = \Ex_\gx((\tilde f_x - \tilde f_y)^2)$. 
\end{remark}

\subsection{Operator-theoretic interpretation of {bd}\Graph}
\label{sec:Operator-theoretic-interpretation-of-bdG}

Recall that we began this section with a comparison of the Poisson boundary representation
\begin{equation}\label{eqn:Poisson-bdy-repn-recall}
  u(x) = \int_{\del \gW} u(y) k(x,dy), 
  \qq u \text{ bounded and harmonic on $\gW \ci \bRd$},
\end{equation}
to the \energy boundary representation
\linenopax
\begin{equation}\label{eqn:boundary-repn-for-harmonic-recall}
  u(x) = \sum_{\bd \Graph} u \dn{h_x} + u(o), 
  \qq u \in \Harm, \text{ and } h_x = \Phar v_x.
\end{equation}

\begin{remark}
  \label{rem:abuse-of-extension-notation}
  For $u \in \Harm$ and $\gx \in \Schw'$, let us abuse notation and write $u$ for $\tilde{u}$. That is, $u(\gx) := \tilde{u}(\gx) = \la u, \gx\ra_\Wiener$. Unnecessary tildes obscure the presentation and the similarities to the Poisson kernel.
\end{remark}



\begin{cor}[Boundary integral representation for harmonic functions]
  \label{thm:Boundary-integral-repn-for-harm} \hfill \\
  For any $u \in \Harm$ and with $h_x = \Phar v_x$, 
  \linenopax
  \begin{equation}\label{eqn:integral-boundary-repn-of-h}
    u(x) = \int_{\Schw'/\Fin} u(\gx) h_x(\gx) \, d\quotprob(\gx) + u(o).
  \end{equation}
  \begin{proof}
    Starting with Lemma~\ref{thm:repkernels-for-Fin-and-Harm}, compute
    \linenopax
    \begin{align}\label{eqn:boundary-repn-for-harmonic-integral-derived}
      u(x) - u(o)
      = \la h_x, u\ra_\energy
      = \cj{\la u, h_x \ra_\energy}
      = \cj{\int_{\Schw'} \cj{u} h_x \, d\quotprob},
    \end{align}
    where the last equality comes by substituting $v=h_x$ in \eqref{eqn:Expectation-formula-for-energy-prod}; recall from Lemma~\ref{thm:vx-is-R-valued} that $\cj{h_x}=h_x$. Note that we are suppressing tildes as in Remark~\ref{rem:abuse-of-extension-notation}. 
  \end{proof}
\end{cor}

\begin{remark}[A Hilbert space interpretation of {bd}\,\Graph]
  \label{rem:boundary-of-G-from-Minlos}
  In view of Corollary~\ref{thm:Boundary-integral-repn-for-harm}, we are now able to ``catch'' the boundary between \Schw and $\Schw'$ by using \LapM and its adjoint. The boundary of \Graph may be thought of as (a possibly proper subset of) $\Squoth$. 
  Corollary~\ref{thm:Boundary-integral-repn-for-harm} suggests that $\mathbbm{k}(x,d\gx) := h_x(\gx) d\quotprob$ is the discrete analogue in \HE of the Poisson kernel $k(x,dy)$, and comparison of \eqref{eqn:boundary-repn-for-harmonic} with \eqref{eqn:integral-boundary-repn-of-h} gives a way of understanding a boundary integral as a limit of Riemann sums:
  \linenopax
  \begin{equation}\label{eqn:boundary-of-G-from-Minlos}
    \int_{\Schw'} u \, h_x \, d\quotprob
    = \lim_{k \to \iy} \sum_{\bd G_k} u(x) \dn{h_x}(x). 
  \end{equation}
  (We continue to omit the tildes as in Remark~\ref{rem:abuse-of-extension-notation}.) By a theorem of Nelson, \quotprob is fully supported on those functions which are H\"{o}lder-continuous with exponent $\ga=\frac12$, which we denote by $\Lip(\frac12) \ci \Schw'$; see \cite{Nelson64}. Recall from Corollary~\ref{thm:vx-is-Lipschitz} that $\HE \ci Lip(\frac12)$.
  See \cite{Arv76a,Arv76b,Minlos63,Nel69}.
\end{remark}

\section{The boundary as equivalence classes of paths}
\label{sec:bdG-as-equivalence-classes-of-paths}

We are finally able to give a concrete representation of elements of the boundary. We continue to use the measure \quotprob from Theorem~\ref{thm:HE-isom-to-L2(S',P)}. 
Recall the Fock space representation of $L^2(\Schw',\prob)$ discussed in Remark~\ref{rem:Wiener-improves-Minlos}:
  \linenopax
  \begin{align}\label{eqn:Fock-repn-of-L2(S',P)-as-HE}
    L^2\left(\Squot,\quotprob\right) \cong \bigoplus_{n=0}^\iy \HE^{\otimes n}.
  \end{align}
  where $\HE^{\otimes 0} := \bC \one$ for a unit ``vacuum'' vector \one corresponding to the constant function, and $\HE^{\otimes n}$ denotes the $n$-fold \emph{symmetric} tensor product of \HE with itself. Observe that \one is orthogonal to \Fin and \Harm, but is not the zero element of $L^2(\Squoth,\quotprob)$.
  
\begin{lemma}\label{thm:int(u)=0-over-bdG}
  For all $v \in \Harm$, $\int_{\Squoth} v \, d\quotprob = 0$.
  \begin{proof}
     The integral $\int_{\Squoth} v \, d\quotprob = \int_{\Squoth} \one v \, d\quotprob$ is the inner product of two elements in $L^2(\Squoth,\quotprob)$ which lie in different (orthogonal) subspaces; see \eqref{eqn:Fock-repn-of-L2(S',P)}.
  \end{proof}
\end{lemma}
Alternatively, Lemma~\ref{thm:int(u)=0-over-bdG} holds because the expectation of every odd-power monomial vanishes by \eqref{eqn:Minlos-expectation-integral}; see also \eqref{eqn:expectation-of-even-monomials} and the surrounding discussion of Remark~\ref{rem:Wiener-improves-Minlos}.

Recall that we abuse notation and write $h_x = \la h_x, \cdot \ra_\Wiener = \tilde h_x$ for elements of $\Schw'$. 

\begin{defn}\label{def:mu_x}
  Denote the measure appearing in Corollary~\ref{thm:Boundary-integral-repn-for-harm} by    
  \linenopax
  \begin{equation}\label{eqn:def:boundary-kernel}
    d\gm_x := (\one + h_x) \, d\quotprob.  
  \end{equation}
  The function \one does not show up in \eqref{eqn:integral-boundary-repn-of-h} because it is orthogonal to \Harm:
\linenopax
\begin{align*}
  \int_{\Squoth} u (\one + h_x) \, d\quotprob
  = \int_{\Squoth} u \, d\quotprob + \la u, h_x\ra_\energy
  = \la u, h_x\ra_\energy,
  \qq\text{for } u \in \Harm,
\end{align*}
where we used Lemma~\ref{thm:int(u)=0-over-bdG}. Nonetheless, its presence is necessary, 
\linenopax
\begin{align*}
  \int_{\Squoth} \one d\gm_x 
  &= \int_{\Squoth} \one (\one + h_x) \, d\quotprob 
  = \int_{\Squoth} \one \, d\quotprob  + \int_{\Squoth} \one h_x \, d\quotprob 
  = 1,
\end{align*}
again by Lemma~\ref{thm:int(u)=0-over-bdG}. 
\end{defn}

\begin{remark}\label{rem:kernel-nonnegative}
  \version{}{\marginpar{$\gm_x \geq 0$ or $h_x \geq 0$?}}
  We have shown that as a linear functional, $\gm_x[\one]=1$. It follows by standard functional analysis that $\gm_x \geq 0$ \quotprob-a.e. on $\Squoth $. Thus, $\gm_x$ is absolutely continuous with respect to \quotprob ($\gm_x \ll \quotprob$) with Radon-Nikodym derivative $\frac{d\gm_x}{d\quotprob} = \one + h_x$.
\end{remark}

\begin{defn}\label{def:harmonic-equivalence-of-paths}
  Recall that a path in \Graph is an infinite sequence of successively adjacent vertices. We say that a path $\cpath = (x_0, x_1, x_2, \dots)$ is a \emph{path to infinity}, and write $\cpath \to \iy$, iff \cpath eventually leaves any finite set $F \ci \verts$, i.e.,
  \linenopax
  \begin{align}\label{eqn:def:path-to-infinity}
    \exists N \text{ such that } n \geq N \implies x_n \notin F.
  \end{align}
  
  If $\cpath_1 = (x_0, x_1, x_2, \dots)$ and $\cpath_2 = (y_0, y_1, y_2, \dots)$ are two paths to infinity, define an equivalence relation by    
  \linenopax
  \begin{equation}\label{eqn:def:harmonic-equivalence-of-paths}
    \cpath_1 \simeq \cpath_2
    \q\iff\q
    \lim_{n \to \iy} (h(x_n)-h(y_n)) = 0, \q\text{for every } h \in \MP.
  \end{equation}
  In particular, all paths to infinity are equivalent when $\Harm = 0$.
\end{defn}
  \glossary{name={\cpath},description={path},sort=G}

If $\gb = [\cpath]$ is any such equivalence class, pick any representative $\cpath = (x_0, x_1, x_2, \dots)$ and consider the associated sequence of measures $\{\gm_{x_n}\}$. These probability measures lie in the unit ball in the weak-$\star$ topology, so Alaoglu's theorem gives a weak-$\star$ limit
  \linenopax
  \begin{equation}\label{eqn:def:prob-meas-corr-to-a-path}
    \gn_\gb := \lim_{n \to \iy} \gm_{x_n}.
  \end{equation}
For any $h \in \Harm$, this measure satisfies 
  \linenopax
  \begin{equation}\label{eqn:h(xn)-becomes-h(b)}
    h(x_{n}) = \int_{\Squoth} \tilde h \,d\gm_{x_n} 
    \limas{n} \int_{\bd \Graph} \tilde h \,d\gn_\gb.
  \end{equation}
Thus, we define $\bd \Graph$ to be the collection of all such \gb, and extend harmonic functions to $\bd \Graph$ via
  \linenopax
  \begin{equation}\label{eqn:def:h(b)}
    \tilde h(\gb) := \int_{\bd \Graph} \tilde h \,d\gn_\gb.
  \end{equation}

\begin{defn}\label{def:bounded-in-HE}
  For $u \in \HE$, denote $\|u\|_\iy := \sup_{x \in \verts} |u(x)-u(o)|$, and say $u$ is \emph{bounded} iff $\|u\|_\iy < \iy$.
\end{defn}

\begin{lemma}\label{thm:Pfin-preserves-boundedness}
  If $v \in \HE$ is bounded, then $\Pfin v$ is also bounded. 
  \begin{proof}
    Choose a representative for $v$ with $0 \leq v \leq K$. Then by Corollary~\ref{thm:Boundary-integral-repn-for-harm} and \eqref{eqn:def:boundary-kernel}, 
    \linenopax
    \begin{align*}
      \Phar v(x) = \int_{\Squoth} v(\gx) h_x(\gx) \, d\quotprob(\gx) + u(o)
           = \int_{\Squoth} v(\gx) \, d\gm_x(\gx) + u(o).
    \end{align*}
    Since $\gm_x$ is a probability measure (cf.~Remark~\ref{rem:kernel-nonnegative}), we have $\Phar v \geq 0$, and hence the finitely supported component $\Pfin v = v - \Phar v$ is also bounded. 
  \end{proof}
\end{lemma}

\begin{lemma}\label{thm:monopoles-and-dipoles-are-bounded}
  Every $v \in \MP$ is bounded. In particular, $\|v_x\|_\iy \leq R^F(o,x)$. 
  \begin{proof}
    According to Definition~\ref{def:LapM}, it suffices to check that $v_x$, \monov and \monof are bounded for each $x$. Furthermore, $\|\monof\|_\iy = \|\Pfin \monov\|_\iy \leq \|\monov\|_\iy$ by Lemma~\ref{thm:Pfin-preserves-boundedness}, and $\monov = v_x + w_o$ by definition, so it suffices to check $v_x$ and $w_o$. By \cite[Lem.~3.70]{Soardi94}, $w_o$ has a representative which is bounded, taking only values between $0$ and $w_o(o) > 0$. It remains only to check $v_x$. The following approach is taken from the ``proof'' of \cite[Conj.~3.18]{ERM}.
    
    Fix $x,y \in \verts$ and an exhaustion $\{G_k\}$, and suppose without loss of generality that $o,x,y \in G_1$. Also, let us consider the representative of $v_x$ specified by $v_x(o)=0$. On a finite network, it is well-known (see \eqref{eqn:vx-as-prob-on-finite} and the surrounding discussion) that 
  \begin{align}\label{eqn:vx-as-prob-on-finite-recalled}
    v_x = R(o,x) u_x,  
  \end{align}
  where $u_x(y)$ is the probability that a random walker (RW) started at $y$ reaches $x$ before $o$, that is, $u_x(y) := \prob_y[\gt_x < \gt_o]$, where $\gt_x$ denotes the hitting time of $x$. This idea is discussed in \cite{DoSn84,LevPerWil08,Lyons:ProbOnTrees}. 
  
    \marginpar{DOUBLE-CHECK THIS}
    Therefore, one can write \eqref{eqn:vx-as-prob-on-finite-recalled} on $G_k$ as $v_x^{(k)} = R_{G_k^F}(o,x) u_x^{(k)}$. In other words, $v_x^{(k)}$ is the unique solution to $\Lap v = \gd_x - \gd_o$ on the finite subnetwork $G_k^F$. Consequently, for every $k$ we have $v_x^{(k)}(y) \leq R_{G_k^F}(o,x)$ for all $y \in G_k$.
    Since $R^F(x,y) = \lim_{k \to \iy} R_{G_k^F}(x,y)$ by \cite[Def.~2.9]{ERM}, we have $\|v_x\|_\iy \leq R^F(o,x)$ for every $x \in \verts$. 
    %
  \end{proof}
\end{lemma}

\begin{theorem}\label{thm:bdy-point-as-linear-fnl}
  Let $\gb \in \bd \Graph$ and let $\cpath = (x_0, x_1, x_2, \dots)$ is any representative of \gb.
  Then $\gb \in \bd \Graph$ defines a continuous linear functional on \Schw via
  \linenopax
  \begin{align}\label{eqn:bdy-pt-as-meas-limit-along-a-path}
    \gb(v) := \lim_{n \to \iy} \int_{\Schw'} \tilde v \, d\gm_{x_n}, 
    \qq v \in \Schw.
  \end{align}
  In fact, the action of \gb is equivalently given by 
  \linenopax
  \begin{align}\label{eqn:bdy-pt-as-eval-limit-along-a-path}
    \gb(v) = \lim_{n \to \iy} \Phar v({x_n})- \Phar v(o), 
    \qq v \in \Schw.
  \end{align}
  \begin{proof}
    To see that \eqref{eqn:bdy-pt-as-meas-limit-along-a-path} and \eqref{eqn:bdy-pt-as-eval-limit-along-a-path} are equivalent, compute
    \linenopax
    \begin{align*}
      \int_{\Schw'} \cj{\tilde v} (\one + h_{x_n}) \,d\quotprob 
      = \cancel{\int_{\Schw'} \cj{\tilde v} \one \,d\quotprob} 
        + \int_{\Schw'} \cj{\tilde v} h_{x_n} \,d\quotprob 
      = \la v, h_{x_n} \ra_\energy 
      = \Phar v({x_n})-\Phar v(o),
    \end{align*}
    because \one is orthogonal to \HE in $L^2(\Schw',\quotprob)$; see \eqref{eqn:Fock-repn-of-L2(S',P)}. 
       
    Now, to see that \eqref{eqn:bdy-pt-as-meas-limit-along-a-path} or \eqref{eqn:bdy-pt-as-eval-limit-along-a-path} defines a {bounded} linear functional, we only need to check that $\sup_{v \in \Schw} \{\gb(v) \suth \|v\|_\energy=1\}$ is bounded, but this is the content of Lemma~\ref{thm:monopoles-and-dipoles-are-bounded}.
    Note that the equivalence relation \eqref{eqn:def:harmonic-equivalence-of-paths} ensures that the limit is independent of the choice of representative \cpath.
  \end{proof}
\end{theorem}

\begin{remark}\label{rem:bdG-as-Diracs-at-infinity}
  In light of \eqref{eqn:bdy-pt-as-eval-limit-along-a-path}, one can think of $\gn_\gb$ in \eqref{eqn:def:prob-meas-corr-to-a-path} as a Dirac mass. Thus, $\gb \in \bd \Graph$ is a boundary point, and integrating a function $f$ against $\gn_\gb$ corresponds to evaluation of $f$ at that boundary point.
\end{remark}

\section{The structure of $\Schw'$}
\label{sec:Structure-of-Schw'}

The next results are structure theorems akin to those found in the classical theory of distributions; see \cite[\S6.3]{Str03} or \cite[\S3.5]{Al-Gwaiz}. If $\HE \ci \Schw$, then Theorem~\ref{thm:structure-of-S'} would say $\Schw' = \bigcup_p \LapS^p(\HE)$ (of course, this is typically false when $\Harm \neq 0$).

\begin{theorem}\label{thm:structure-of-S'}
  The distribution space $\Schw'$ is
  \linenopax
  \begin{align}\label{eqn:structure-of-S'}
    \Schw' = \{\gx(u) = \la \LapS^p u, v\ra_\energy \suth u \in \Schw, v \in \HE, p \in \bZ^+\}.
   \end{align}
  \begin{proof}
    It is clear from the Schwarz inequality that $\gx(u) = \la \LapS^p u, v\ra_\energy$ defines a continuous linear functional on \Schw, for any $v \in \HE$ and nonnegative integer $p$. For the other direction, we use the same technique as in Lemma~\ref{thm:Wiener-product-as-Lap(p)-powers}. Observe that if $\gx \in \Schw'$, then there exists $K,p$ such that $|\gx(u)| \leq K \|\LapS^p u\|_\energy$ for every $u \in \Schw$. This implies that the map $\gx: \LapS^p u \mapsto \gx(u)$ is continuous on the subspace $Y = \spn\{\LapS^p u \suth u \in \HE, p \in \bZ^+\}$. This can be extended to all of \HE by precomposing with the orthogonal projection to $Y$. Now Riesz's lemma gives a $v \in \HE$ for which $\gx(u) = \la \LapS^p u, v\ra_\energy$. 
%
  \end{proof}
\end{theorem}

Note that $v \in \HE$ may not lie in the domain of $\LapS^p$. If it did, one would have $\la \LapS^p u, v\ra_\energy = \la u, \LapS^p v\ra_\Wiener = \la u, \LapS^p f\ra_\Wiener$, where $f = \Pfin v$. The theorem could then be written $\Schw' = \bigcup_{p=0}^\iy \LapS^p (\Fin)$. However, this turns out to have contradictory implications.

\version{}{
\begin{defn}\label{def:Lap-on-distributions}
  Extend \Lap to $\Schw'$ by defining
  \linenopax
  \begin{equation}\label{eqn:Lap-on-distributions}
    \Lap \gx(v_x) := \la \gd_x, \gx \ra_\Wiener,
  \end{equation}
  so that $\Lap\gx(v_x) = \sum_{y \nbr x} \cond_{xy} (\gx(v_x)-\gx(v_y))$ follows readily from Lemma~\ref{thm:dx-as-vx}. 
  
  Now extend \Lap to $\tilde \Lap$ defined on $\tilde v_x \in L^2(\Squot,\quotprob)$ by $\tilde \Lap (\tilde v_x)(\gx) := \widetilde{\Lap v_x} (\gx)$, so that 
  \linenopax
  \begin{equation}\label{eqn:Lap-on-L2(distns)}
    \tilde \Lap: \tilde v_x \mapsto \cond(x) \tilde v_x - \sum_{y \nbr x} \cond_{xy} \tilde v_y.
  \end{equation}
  Since $v_x \mapsto \tilde v_x$ is an isometry, it is no great surprise that
  \linenopax
  \begin{align}\label{eqn:Lap-on-distributions-inner-product}
    \la \tilde v_x, \tilde \Lap \tilde v_y\ra_{L^2}
    &= \int_{\Schw'} \tilde v_x(\gx) \tilde v_y(\Lap \gx)\,d\quotprob(\gx)
    = \la v_x, \Lap v_y\ra_\energy.
  \end{align}
\end{defn}
}

We now provide two results enabling one to recognize certain elements of $\Schw'$.

\begin{lemma}\label{thm:Schw-dist-criterion}
  A linear functional $f:\Schw \to \bC$ is an element of $\Schw'$ if and only if there exists $p \in \bN$ and $F_0,F_1,\dots F_p \in \HE$ such that
  \linenopax
  \begin{equation}\label{eqn:Schw-dist-criterion}
    f(u) = \sum_{k=0}^p \la F_k, \LapS^k u \ra_\energy, \q\forall u \in \HE.
  \end{equation}
  \begin{proof}
    By definition, $f \in \Schw'$ iff $\exists p,C<\iy$ for which $|f(u)| \leq C \|u\|_p$ for every $u \in \Schw$. Therefore, the linear functional
    \linenopax
    \begin{align*}
      \gF:\bigoplus\nolimits_{k=0}^p \dom(\LapS^k) \to \bC
      \qq\text{by}\qq
      \gF(u,\LapS u, \LapS^2 u, \dots \LapS^p u) = f(u)
    \end{align*}
    is continuous and Riesz's Lemma gives $F = (F_k)_{k=0}^p \in \bigoplus_{k=0}^p \negsp[5]\HE$ with
    \linenopax
    \begin{align*}
      f(u) &= \la F, (u,\LapS u, \dots \LapS^p u)\ra_{\bigoplus \HE}
        = \sum_{k=0}^p \la F_k, \LapS^k u \ra_{\bigoplus\negsp[5]\HE}.
        &&\qedhere
    \end{align*}
  \end{proof}
\end{lemma}

\begin{cor}\label{thm:Schw'=HE-when-Lap-bdd}
  If $\LapS:\HE \to \HE$ is bounded, then $\Schw'=\HE$.
  \begin{proof}
    We always have the inclusion $\HE \hookrightarrow \Schw'$ by taking $p=0$. If  \LapS is bounded, then the adjoint $\LapS^\ad$ is also bounded, and \eqref{eqn:Schw-dist-criterion} gives
    \linenopax
    \begin{equation}\label{eqn:Schw'=HE-when-Lap-bdd}
      f(u) = \left\la \sum_{k=0}^p (\LapS^\ad)^k  F_k, u \right\ra_{\bigoplus\negsp[5]\HE},
      \q\forall u \in \Schw.
    \end{equation}
    Since \Schw is dense in \HE by Lemma~\ref{thm:Schw-dense-in-HE}, we have $f = \sum_{k=0}^p (\LapS^\ad)^k  F_k \in \HE$.
  \end{proof}
\end{cor}

\section{Remarks and references}
\label{sec:Remarks-and-References-boundary}

Boundary theory is a well-established subject; see e.g., \cite{Brelot} and \cite{Doob59, Doob66}. The deep connections between harmonic analysis, probability, and potential theory have led to several notions of boundary and we will not attempt to give complete references. However, we recommend \cite{Saw97, Woess09} for introductory material on Martin boundary and 
\cite{Car73a, Woess00} for further information. The papers \cite{Yamasaki79, TerryLyons} and \cite{Nash-Will59} are foundational with regard to connections between energy and transience. With regard to infinite graphs and finite-energy functions, see \cite{Soardi94, SoardiWoess91, CaW92, Dod06, PicWoess90, PicWoess88, Wo86, Thomassen90}. An attractive and modern presentation especially well suited to the needs of our present chapter is \cite{CaSoWo93} by Cartwright, Soardi and Woess. An excellent book for what we need on path-space integrals is \cite{Hida80}.

The boundary representation given in Corollary~\ref{thm:Boundary-integral-repn-for-harm} above is related to a large number of analogous representations in the literature; see for example \cite{Soardi94}, \cite{Spitzer}, \cite{Stroock}, \cite[Thm.~24.7]{Woess00}, or \cite[Th.~3.1 and Thm.~4.1]{Saw97}. There are two primary differences between these more traditional approaches and the one adopted here:
\begin{enumerate}
  \item we focus on the harmonic functions of finite energy (as opposed to the nonnegative or bounded harmonic functions), and
  \item our representation is developed via Hilbert spaces. 
\end{enumerate}
In fact, the latter is made possible by the former. 
However there are no easy ways of relating say pointwise bounded functions to finite-energy functions on an infinite weighted graph. Hence Corollary~\ref{thm:Boundary-integral-repn-for-harm} does not immediately compare with analogous theorems in the literature. 

The reader may additionally wish to consult \cite{Woess00, Wo89, Wo95a, Saloff-Woess09, Saloff-Woess06, KaimanovichWoess02, Kig09b, Ign08, BrWo05, Gui72}.

%% file: Lap-on-HE.tex

\chapter{The Laplacian on \HE}
\label{sec:Lap-on-HE}

\headerquote{I have tried to avoid long numerical computations, thereby following Riemann's postulate that proofs should be given through ideas and not voluminous computations.}{---~D.~Hilbert}

We study the operator theory of the Laplacian in some detail, examining the various domains and self-adjoint extensions. One of the primary goals in \S\ref{sec:Properties-of-Lap-on-HE} is to determine when $v_x$ lies in the domain or range of \LapV; this may indicate when $v_x$ lies in the Schwartz space \Schw developed in \S\ref{sec:Gel'fand-triples-and-duality}. We also identify a particular self-adjoint extension \LapH for use in the constructions in \S\ref{sec:the-boundary}. Also, we give technical conditions which must be considered when the graph contains vertices of infinite degree and/or the conductance functions $\cond(x)$ is unbounded on \verts. 
A technical obstacle must be overcome: While $\ell^2(\verts)$ has a canonical orthonormal basis, this is not so for \HE. Instead, the analysis of \HE is carried out with the use of an independent and spanning system $\{v_x\}$ in \HE; these vectors are non-orthogonal, but this non-orthogonality is a rich source of information. 

In \S\ref{sec:boundary-form}, we relate the boundary term of \eqref{eqn:intro:discrete-Gauss-Green} to a boundary form akin to that of classical functional analysis; see Definition~\ref{def:boundary-form}. In Theorem~\ref{thm:LapV-not-ess-selfadjoint-iff-Harm=0}, we show that if \Lap fails to be essentially self-adjoint, then $\Harm \neq \{0\}$. In general, the converse does not hold: any homogeneous tree of degree 3 or higher with constant conductances provides a counterexample; cf. Corollary~\ref{thm:c-bdd-implies-Lap-selfadjoint-on-HE}.

In \S\ref{sec:Dual-frames-and-the-energy-kernel} we consider the systems $\{v_x\}$ and $\{\gd_x\}$ and a kind of spectral reciprocity between them, in terms of frame duality. In previous parts of this book, we approximated infinite networks by truncating the domain; this is the idea behind the definition of \Fin and the use of exhaustions. This approach corresponds to a restriction to $\spn\{\gd_x\}_{x \in F}$, where $F$ is some finite subset of \verts. In \S\ref{sec:Dual-frames-and-the-energy-kernel}, we consider truncations in the dual variable, i.e., restrictions to sets of the form $\spn\{v_x\}_{x \in F}$. Note that an element of this set generally will \emph{not} have finite support.

We use $\ran T$ to denote the range of the operator $T$, and $\ker T$ to denote its kernel (nullspace). We continue to use the notation from \S\ref{sec:the-boundary}: let $V := \spn\{v_x\}_{x \in \verts}$ denote the vector space of \emph{finite} linear combinations of dipoles. Then let \LapV be the closure of the Laplacian when taken to have the dense domain $V$.

\section{Properties of \Lap on \HE}
\label{sec:Properties-of-Lap-on-HE}

\begin{defn}\label{def:power-bound-preview}
  The network $(\Graph, \cond)$ \sats the \emph{Powers bound} iff $\displaystyle \|\cond\| := \sup_{x \in \verts} \cond(x) < \iy$.
\end{defn}
  \glossary{name={$\|\cond\|$},description={operator norm of \cond, also $\sup \cond(x)$},sort=M}
The Powers bound is used more in \S\ref{sec:L2-theory-of-Lap-and-Trans} (see Definition~\ref{def:power-bound} and the surrounding discussion); we include it here for use in a couple of technical lemmas.

\begin{lemma}\label{thm:Powers-bd-implies-Lap-maps-into-bdd-fns}
  If the Powers bound is satisfied, then \Lap maps \HE into $\ell^\iy(\verts)$. 
  \begin{proof}
    By Lemma~\ref{thm:<delta_x,v>=Lapv(x)} and \eqref{eqn:energy-of-Diracs},
    $|\Lap v(x)| = |\la \gd_x, v\ra_\energy|
      \leq \|\gd_x\|_\energy \cdot \|v\|_\energy
      = \cond(x)^{1/2} \|v\|_\energy.$
  \end{proof}
\end{lemma}

\begin{lemma}\label{thm:LapV-maps-V-to-V}
  If $\deg(x) < \iy$ for every $x \in \verts$, or if $\uBd < \iy$, then $\ran \LapV \ci \dom \LapV$.   
  \begin{proof}
    It suffices to show that $\LapV v_x = \gd_x-\gd_o \in \dom \LapV$ for every $x \in \verts$, and this will be clear if we show $\gd_x \in \dom \LapV$. By Lemma~\ref{thm:dx-as-vx}, $\gd_x = \cond(x) v_x - \sum_{y \nbr x} \cond_{xy} v_y$. If $\deg(x)$ is always finite, then we are done. If not, we need to see why $\sum_{y \nbr x} \cond_{xy} v_y \in \dom \LapV$ for any fixed $x \in \verts$. 
    
    Fix $x \in \verts$ and denote $\gf := \sum_{y \nbr x} \cond_{xy} v_y$ and $\gf_k := \sum_{y \in G_k} \cond_{xy} v_y$. It is clear that $\|\gf-\gf_k\|_\energy \to 0$. We next show $\left\|\LapV \gf_k - \sum_{y \nbr x} \cond_{xy} (\gd_y - \gd_o)\right\|_\energy \to 0$, from which it follows that $\{\LapV \gf_k\}$ is Cauchy, and that $\gf \in \dom \LapV$ with $\LapV \gf =  \sum_{y \nbr x} (\gd_y - \gd_o)$:
    \linenopax
    \begin{align*}
      \left\|\LapV \gf_k - \sum_{y \nbr x} \cond_{xy} (\gd_y - \gd_o)\right\|_\energy^2
      &= \left\|\sum_{y \in G_k^\complm} \cond_{xy} (\gd_y - \gd_o)\right\|_\energy^2 \\
      &\leq \left(\sum_{y \in G_k^\complm} \cond_{xy} \|\gd_y - \gd_o\|_\energy\right)^2 \\
      &\leq \cond(x) \sum_{y \in G_k^\complm} \cond_{xy} \|\gd_y - \gd_o\|_\energy^2 \\
      &= \cond(x) \left(\sum_{y \in G_k^\complm} \cond_{xy} \|\gd_y\|_\energy^2 
       - 2 \sum_{y \in G_k^\complm} \cond_{xy} \la \gd_y,\gd_o\ra_\energy
       + \sum_{y \in G_k^\complm} \cond_{xy} \|\gd_o\|_\energy^2 \right)\\
      &= \cond(x)\left(\sum_{y \in G_k^\complm} \cond_{xy} \cond(y) 
       + 2 \sum_{y \in G_k^\complm} \cond_{xy} \cond_{oy}
       + \cond(o) \sum_{y \in G_k^\complm} \cond_{xy} \right) \\
      &\leq \uBd(3\uBd + \cond(o)) \sum_{y \in G_k^\complm} \cond_{xy},
    \end{align*}
    which tends to 0 as $k$ gets large. Note that $\cond_{oy} < 1$ for $y \in G_k^\complm$ with $k$ sufficiently large. 
  \end{proof}
\end{lemma}

\subsection{Finitely supported functions and the range of \Lap}
\label{Fin-and-the-range-of-Lap}

In Remark~\ref{rem:ranLap-vs-Fin} we showed that one always has $\ran \LapV \ci \Fin$ and hence $\Harm \ci \ker \LapV^\ad$. The rest of this section is roughly an examination of the reverse containment, i.e., what conditions give $\ran \LapV = \Fin$. Determining when $\ran \LapV = \Fin$ essentially boils down to the following technical question: when is $\spn\{\gd_x - \gd_o\}$ dense in \Fin? It is curious that this never happens on a finite network (Lemma~\ref{thm:dx-are-linearly-independent}), but is often true on an infinite network. However, see Example~\ref{exm:Fin_2-not-dense-in-Fin}.

\begin{defn}\label{def:Fin2-and-Fin1}
   Let $\Fin_2$ be the \energy-closure of $\spn\{\gd_x-\gd_o\}$ and let $\Fin_1$ be the orthogonal complement of $\Fin_2$ in \Fin. This extends the decomposition $\HE = \Fin \oplus \Harm$, in some cases, to $\HE = \Fin_2 \oplus \Fin_1 \oplus \Harm$.
\end{defn}

Example~\ref{exm:Fin_2-not-dense-in-Fin} shows a situation in which $\Fin_2$ is not dense in \Fin.
 
\version{}{
\marginpar{Any ideas?}
\begin{lemma}\label{thm:deg(x)-implies-Fin1=0}
  Let $(\Graph, \cond)$ be an infinite network. If $\deg(x) < \iy$ for $x \in \Graph$, then $\Fin = \Fin_2$.
  \begin{proof}
    Needed.
  \end{proof}
\end{lemma}
}

\begin{lemma}\label{thm:powers-bd-implies-Fin1=0}
  Let $(\Graph, \cond)$ be an infinite network. If $\uBd < \iy$, then $\Fin = \Fin_2$.
  \begin{proof}
    It suffices to approximate the single Dirac mass $\gd_o$ by linear combinations of differences. For each $n$, fix $n$ vertices $\{x_k^{(n)}\}_{k=1}^n$, no two of which are adjacent. Therefore, define $\gf_n := \frac1n \sum_{k=1}^n (\gd_o - \gd_{x_k^{(n)}})$ and compute
  \linenopax
  \begin{align*}
    \|\gd_o - \gf_n\|_\energy^2
    = \left\| \frac1n \sum_{k=1}^n \gd_{x_k^{(n)}} \right\|_\energy^2
    = \frac1{n^2} \sum_{k=1}^n \left\| \gd_{x_k^{(n)}} \right\|_\energy^2
    \leq \frac1n \sup_{1 \leq k \leq n} \cond(x_k^{(n)}) 
    \leq \frac{\uBd}n \to 0,
  \end{align*}
  where the second equality comes by orthogonality; for $j \neq k$, $\gd_{x_k^{(n)}}$ and $\gd_{x_j^{(n)}}$ are not adjacent, hence $\la\gd_{x_k^{(n)}},\gd_{x_j^{(n)}}\ra_\energy=0$ by \eqref{eqn:energy-of-Diracs}. Now it is trivial to approximate $\gd_z = (\gd_z - \gd_o) + \gd_o$. 
  \end{proof}
\end{lemma}

The idea of Lemma~\ref{thm:powers-bd-implies-Fin1=0} is illustrated on the binary tree in Example~\ref{exm:Harm-notin-domLap}.

\subsection{Harmonic functions and the domain of \Lap}
\label{sec:harmonic-functions-and-the-domain-of-Lap}

Curiously, even though $\Lap h(x)=0$ pointwise for every $x \in \verts$, it may happen that $h$ is not in the domain of \Lap. Example~\ref{exm:Harm-notin-domLap} discusses a nontrivial harmonic function on the binary tree which does not appear to be in the domain of \LapV. However, harmonic functions are always in the domain of the adjoint $\LapV^\ad$ by Lemma~\ref{thm:ker(Lapadj)=Harm}.

\begin{lemma}\label{thm:Lap-extends-to-Harm-by-0}
  If \LapVc is any Hermitian extension of \LapV whose domain contains \Harm, then $\LapVc h = 0$ for any $h \in \Harm$. Moreover, $\LapVc u \in \Fin$ for any $u \in \dom \LapVc$.
  \begin{proof}
    Recall that we have the following ordering of operators: $\LapV \ci \LapVc \ci \LapV^\ad$. Since $\LapV^\ad$ is an extension of \LapVc and $\Harm \ci \LapVc$, the first claim follows immediately from Lemma~\ref{thm:ker(Lapadj)=Harm}. The second claim now follows from the first because $\la \LapVc v, h\ra_\energy = \la v, (\LapVc)^\ad h\ra_\energy = 0$ for every $h \in \Harm$, since $\LapVc \ci (\LapVc)^\ad$.
  \end{proof}
\end{lemma}

We have a partial converse of Lemma~\ref{thm:ker(Lapadj)=Harm}. Note that if $\spn\{\gd_x-\gd_o\}$ is dense in \Fin (as discussed in Lemma~\ref{thm:powers-bd-implies-Fin1=0}), then Lemma~\ref{thm:kerLapadj=spn(diff)perp} implies $\ker \LapV^\ad = \Harm.$

\begin{lemma}\label{thm:kerLapadj=spn(diff)perp}
  $\ker \LapV^\ad$ is the orthogonal complement of $\spn\{\gd_x-\gd_o\}$.
  \begin{proof}
    Suppose $u \in \ker \LapV^\ad$ so that $\LapV^\ad u = 0$. Then 
    \linenopax
    \begin{align*}
      0 = \la \LapV^\ad u, v_x \ra_\energy 
= \la u, \LapV v_x \ra_\energy
= \la u, \gd_x - \gd_o \ra_\energy.
    \end{align*}
    This shows $u$ is orthogonal to $\spn\{\gd_x-\gd_o\}$.
  \end{proof}
\end{lemma}

The Lemma~\ref{thm:powers-bd-implies-Fin1=0} gives an idea of when the hypotheses of Lemma~\ref{thm:kerLapadj=spn(diff)perp} are satisfied. In fact, a weaker hypothesis will suffice: one just needs to be able to find an infinite subset of nonadjacent vertices on which $\cond(x)$ is bounded.

\begin{defn}\label{def:LapH}
  Define \LapH to be the extension of \LapV to the domain $\dom \LapV + \Harm$ by $\LapH(v+h) := \LapV v$. By abuse of notation, let \LapH denote the closure of \LapH with respect to the graph norm; see Definition~\ref{def:graph-of-operator}.
\end{defn}

\begin{lemma}\label{thm:LapH-well-def-and-semibdd}
  \LapH is well defined, Hermitian, and semibounded. 
  \begin{proof}
    We must check that $\LapH(0)=0$, so suppose $v+h=0$ for $v \in V$ and $h \in \Harm$. Then Lemma~\ref{thm:ker(Lapadj)=Harm} gives $\LapV^\ad(v+h)=0$, whence $\LapV v = -\LapV^\ad h = 0$.
  \end{proof}
\end{lemma}

\begin{theorem}\label{thm:LapH-is-essentially-self-adjoint}
  \LapH is self-adjoint.
  \begin{proof}
    Let $w \in \HE$ satisfy $\LapH^\ad w = -w$. To see that $w=0$, note that $w \in \dom \LapH^\ad$, so $\LapV^\ad w \in \Fin$ by Lemma~\ref{thm:LapVad-maps-into-Fin}, just below. But then $w = - \LapH^\ad w = - \LapV^\ad w \in \Fin$, so
    \linenopax
    \begin{align*}
      \| w\|_\energy^2 
      = \la w, w\ra_\energy 
      = \sum_{\verts} \cj{w} \Lap w
      = -\sum_{\verts} |w|^2
      \leq 0,
    \end{align*}
    so that $w=0$ in \HE. This shows \LapH is essentially self-adjoint, but \LapH is closed by definition, so it is self-adjoint.
  \end{proof}
\end{theorem}

\begin{lemma}\label{thm:LapVad-maps-into-Fin}
  $\dom \LapH^\ad = \{w \in \dom \LapV^\ad \suth \LapV^\ad w \in \Fin\}$.
  \begin{proof}
    For purposes of this proof, it is permissible to work with \HE as a real vector space and complexify afterwards.
    
    ($\ci$) Suppose that $w \in \dom \LapH^\ad$, i.e., we have the estimate
    \linenopax
    \begin{align}\label{eqn:w-in-dom-LapHad}
      \left|\la w, \LapH(v+h)\ra_\energy\right|
      \leq C_1 \|v+h\|_\energy,
      \qq \text{ for all } v \in V \text{ and } h \in \Harm.
    \end{align}
    Then for all $t \in \bR$,
    \linenopax
    \begin{align*}
      \left|\la w, \LapH v\ra_\energy\right|^2
      \leq C_1^2 \|v+th\|_\energy^2
      \leq C_1^2 \|v\|_\energy^2 + 2t|\la v, h \ra_\energy|^2 + t^2\|h\|_\energy^2,
    \end{align*}
    for all $v \in V$ and $h \in \Harm$. This quadratic polynomial in $t$ is nonnegative, and hence its discriminant must be nonpositive, so that
    \linenopax
    \begin{align*}
      C_1^4 |\la v, h\ra_\energy|^2
      &\leq C_1^2 \|h\|_\energy^2 \left(C_1^2 \|v\|_\energy^2 - \left|\la w, \LapH v\ra_\energy\right|^2\right) \\
      \| P_{h} v\|_\energy^2 
      = \frac{|\la v, h\ra_\energy|^2}{\|h\|_\energy^2}
      &\leq \|v\|_\energy^2 - C_2 \left|\la w, \LapH v\ra_\energy\right|^2
    \end{align*}
    where $P_{h}$ is projection to the rank-1 subspace spanned by $h$ and $C_2 = \frac1{C_1}$. If we let $\{h_i\}$ be an ONB for \Harm, then 
    \linenopax
    \begin{align*}
      C_2^2 \left|\la w, \LapH v\ra_\energy\right|^2 
      &\leq \|v\|_\energy^2 - \|P_{h_1} v\|_\energy^2, \qq \text{for all } v \in V.
    \end{align*}    
    Inductively substituting $v=v-h_2$, $v = v - (h_2 + h_3)$, etc, we have
    \linenopax
    \begin{align*}
      C_2^2 \left|\la w, \LapH (v - h_2) \ra_\energy\right|^2 
      &\leq \|v - P_{h_2}v \|_\energy^2 - \|P_{h_1} v\|_\energy^2
       = \|v\|_\energy^2 - \left(\|P_{h_2}v \|_\energy^2 + \|P_{h_1} v\|_\energy^2 \right) \\    
       &\vdots \\
      C_2^2 \left|\la w, \LapH (v + {\textstyle\sum_i}h_i)\ra_\energy\right|^2 
      &\leq \| v\|_\energy^2 - {\textstyle\sum_i} \|P_{h_i} v\|_\energy^2
       = \|v\|_\energy^2 - \|\Phar v\|_\energy^2.
    \end{align*}    
    By the definition of \LapH, all the left sides are equal to $C_2^2 \left|\la w, \LapH v\ra_\energy\right|^2 = C_2^2 \left|\la w, \LapV v\ra_\energy\right|^2$. Since $\|\Pfin v\|_\energy = \|v\|_\energy^2 - \|\Phar v\|_\energy^2$, we have established
    \linenopax
    \begin{align*}
      \left|\la w, \LapV v\ra_\energy\right|
      &\leq C_3 \|\Pfin v\|_\energy, \qq \text{for all } v \in V.
    \end{align*}
    Now Riesz's lemma gives an $f \in \Fin$ such that
    \linenopax
    \begin{align*}
      \la w, \LapV v\ra_\energy 
      = \la f, \Pfin v\ra_\energy, \qq \text{for all } v \in V.
    \end{align*}
    However, orthogonality allows one to remove the projection (since the first argument is already in \Fin), whence $\la \LapV^\ad w, v\ra_\energy = \la f, v\ra_\energy$ for all $v \in V$, and so $\LapV^\ad w = f \in \Fin$.

    ($\ce$) Let $w$ be in the set on the right-hand side. To see $w \in \dom \LapH^\ad$, we need the estimate \eqref{eqn:w-in-dom-LapHad}, but
    \linenopax
    \begin{align*}
      \left|\la w, \LapH(v+h)\ra_\energy\right|
      = \left|\la w, \LapV v\ra_\energy\right| 
      = \left|\la \LapV^\ad w, v\ra_\energy\right| 
      = \left|\la \LapV^\ad w, \Pfin v\ra_\energy\right|, 
    \end{align*}
    where the last equality follows by the hypothesis $\LapV^\ad w \in \Fin$. This gives $\left|\la w, \LapH(v+h)\ra_\energy\right| \leq \|\LapV^\ad w\|_\energy \cdot \|\Pfin (v+h)\|_\energy$, but $\|\Pfin v\|_\energy = \|\Pfin (v+h)\|_\energy \leq \|v+h\|_\energy$, so \eqref{eqn:w-in-dom-LapHad} follows.
  \end{proof}
\end{lemma}

\begin{cor}\label{thm:sa-iff-dom-contains-Harm}
  A closed extension of \LapV is self-adjoint if and only if \Harm is contained in its domain.
  \begin{proof}
    It is helpful to keep in mind the operator ordering $\LapV \ci \LapH = \LapH^\ad \ci \LapV^\ad$.
    
    \fwd Let $\tilde \Lap$ be a self-adjoint extension of \LapV. If $\LapH \ci \tilde \Lap$, then the result is obvious, and if $\tilde \Lap \ci \LapH$, then again $\LapH \ci\LapH^\ad \ci (\tilde \Lap)^\ad = \tilde \Lap$, and the result is equally obvious.
    
    \bwd If $\tilde \Lap$ is a closed extension of \LapV with $\Harm \ci \dom \tilde \Lap$, then $\LapH \ci \tilde \Lap$, so
    \begin{align*}
       \LapH^{\textrm{clo}} \ci \tilde \Lap \ci (\tilde \Lap)^\ad \ci (\LapH^{\textrm{clo}})^\ad \ci \LapH^{\textrm{clo}},
    \end{align*}
    where the first inclusion holds because $\tilde \Lap$ is closed, and the last by Theorem~\ref{thm:LapH-is-essentially-self-adjoint}.
  \end{proof}
\end{cor}

\version{}{\marginpar{For when/if we have more to say ...}
In functional analysis, one commonly considers the Hilbert space obtained as the closure of the quadratic form induced by a particular operator. 
\begin{defn}\label{def:quadratic-form-space}
  Define the quadratic form $Q$ on $\dom \LapV$ by $Q(v) = \la v, \LapV v\ra_\energy$, and let $\sH_Q$ be the closure of $\dom \LapV$ with respect to this form.
\end{defn}

\begin{lemma}[{\cite[\S{XII.5.2}]{DuSc88}}]
  \label{thm:Friedrichs-identified}
  The domain of the Friedrichs extension of \LapV is $\sH_Q \cap \dom \LapV^\ad$.
\end{lemma}
}

\section{The defect space of \LapV}
\label{sec:defect-space}

Let \LapV once again denote the graph closure of the operator \Lap on the (dense) domain $V := \spn\{v_x\}$.

\begin{defn}\label{def:Def-LapV}
  Since $\LapV$ is Hermitian on its domain by Corollary~\ref{thm:Lap-Hermitian-on-V}, Definition~\ref{def:defect-index} and Theorem~\ref{thm:essentially-self-adjointness-criterion} imply that the defect space of \LapV is
  \linenopax
  \begin{align}\label{eqn:def:Def-LapV}
    \Def := \{v \in \dom \LapV^\ad \suth \LapV^\ad v = - v\}.
  \end{align}
  Observe also that $\Def^\perp = \ran(\id + \LapV)$.
\end{defn}

\begin{lemma}\label{thm:defect-iff-u+Lapu=k}
  $u$ is a defect vector of \LapV if and only if there is a constant $k$ such that $\Lap u(x) = -u(x) + k$ at each $x \in \verts$. 
  \begin{proof}
    Recall that the meaning of such a pointwise identity is that $u \in \dom \LapV^\ad$ and $\LapV^\ad u = -u + k$ in \HE; see Lemma~\ref{thm:pointwise-identity-implies-adjoint-identity}. The reverse implication is obvious; for the obverse it suffices to check the claim against the (dense) energy kernel:
    \linenopax
    \begin{align*}
      0 
      = \la v_x, \Lap^\ad u + u \ra_\energy 
      = \la \gd_x - \gd_o, u \ra_\energy + \la v_x, u \ra_\energy 
      = \Lap u (x) - \Lap u (o) + u(x) - u(o),
    \end{align*}
    by Lemma~\ref{thm:<delta_x,v>=Lapv(x)}, which proves the claim with $k = \Lap u (o) + u(o)$.
  \end{proof}
\end{lemma}

\begin{remark}[Defect vectors and the Gauss-Green formula]
  \label{rem:Defect-vectors-and-the-Gauss-Green-formula}
  We have introduced the defect space of \LapV here to alleviate any concerns regarding the convergence of $\sum_{x \in \verts} \cj{u}(x) \Lap u(x)$ in \eqref{eqn:E(u,v)=<u_0,Lapv>+sum(normals)}; the reader will note that if $u$ is a defect vector, then 
\linenopax
\begin{align*}
  \sum_{x \in \verts} \cj{u}(x) \Lap u(x)
  = - \sum_{x \in \verts} |u(x)|^2,
\end{align*}
which must equal $-\iy$, since there are no defect vectors in $\ell^2$. This is a reasonable concern, as there do exist networks with nontrivial defect; see Example~\ref{exm:integer-lattice-with-defect}. However, such defect vectors are proscribed by the hypotheses of Theorem~\ref{thm:E(u,v)=<u,Lapv>+sum(normals)}, by the following lemma.
\end{remark}

\begin{lemma}\label{thm:no-defect-in-domLapV}
  $\dom \LapV \cap \Def = 0$. 
  \begin{proof}
    Suppose $u \in \dom \LapV \cap \Def$. Note that $\LapV^\ad$ is an extension of \LapV, so such a $u$ satisfies $\LapV u = -u$. However, since \LapV is semibounded on its domain by Corollary~\ref{thm:Lap-Hermitian-on-V}, this implies
    \linenopax
    \begin{align*}
      0 \leq \la u, \LapV u\ra_\energy 
      = \la \LapV^\ad u, u\ra_\energy 
      = -\la u, u\ra_\energy = -\|u\|_\energy^2,
    \end{align*}
    whence $u=0$.
  \end{proof}
\end{lemma}

\version{}{
\begin{cor}\label{thm:boundary-sum-for-Def-perp}
  For $u \in \Def(\LapV)^\perp$, 
  \linenopax
  \begin{align}\label{eqn:boundary-sum-for-Def-perp}
    \sum_{\bd G} u \dn u = \| \Phar (\id + \Lap)^{-1} u \|_\energy^2.
  \end{align}
  \begin{proof}
    By Lemma~\ref{thm:Def-perp-is-range-of-(1+opclosure)}, we can take $u$ to be of the form $u = v + \Lap v$ with $v \in \dom\LapV$. However, we begin with the assumption that $v \in \spn\{v_x\}$ and compute
    \linenopax
    \begin{align*}
      \sum_{\verts} u \Lap u
      &= \sum_{\verts} (v + \Lap v) \Lap (v + \Lap v) \\
      &= \sum_{\verts} v \Lap v  + \sum_{\verts} v \Lap^2 v + \sum_{\verts} \Lap v \Lap v + \sum_{\verts} \Lap v \Lap^2 v \\
      &= \la v,v\ra_\energy + \la v, \Lap v\ra_\energy + \la \Lap v, v\ra_\energy + \la \Lap v, \Lap v\ra_\energy
    \end{align*}
    by Theorem~\ref{thm:uLapv-has-no-bdy-term}, since we are assuming $v \in \spn\{v_x\}$.
    \linenopax
    \begin{align*}
      \|u\|_\energy^2 =  \sum_{\verts} u \Lap u + \| \Phar v \|_\energy^2.
    \end{align*}
    $\dots$\marginpar{proof needed}
  \end{proof}
\end{cor}
}

\subsection{The boundary form}
\label{sec:boundary-form}

In this section, we relate the defect of \Lap to the boundary term of the Discrete Gauss-Green formula (Theorem~\ref{thm:E(u,v)=<u,Lapv>+sum(normals)}), thereby extending Theorem~\ref{thm:TFAE:Fin,Harm,Bdy}. The reader may find \cite[\S{XII.4.4}]{DuSc88} to be a useful reference. 

\begin{defn}\label{def:boundary-form}
  Define the \emph{boundary form}
  \linenopax
  \begin{equation}\label{eqn:boundary-form}
    \bdform(u,v) :=
    \tfrac1{2\ii}\left(\la \LapV^\ad u, v\ra_\energy - \la u, \LapV^\ad v\ra_\energy\right),
    \qq u,v \in \dom(\LapV^\ad).
  \end{equation}
  To see the significance of \bdform for the defect spaces, note that if $\LapV^\ad f = zf$ where $z \in \bC$ with $\Im z \neq 0$, then $\bdform(f,f) = (\Im z) \|f\|_\energy^2$.
\end{defn}

\begin{lemma}\label{thm:boundary-form-vanishes-on-dom(LapV)}
  The boundary form $\bdform(u,v)$ vanishes if $u$ or $v$ lies in $\dom(\LapV)$.
  \begin{proof}
    For $v \in \dom (\LapV)$, $\la \LapV^\ad u, v\ra_\energy = \la u, \LapV v\ra_\energy$ by the definition of the adjoint, and $\la u, \LapV v\ra_\energy = \la u, \LapV^\ad v\ra_\energy$ by the fact that $\LapV^\ad$ extends \LapV. Hence, both terms of \eqref{eqn:boundary-form} are equal for $u,v \in \dom (\LapV)$. The proof is identical if $u \in \dom(\LapV)$.
  \end{proof}
\end{lemma}

The following result extends Theorem~\ref{thm:TFAE:Fin,Harm,Bdy}.

\begin{theorem}\label{thm:LapV-not-ess-selfadjoint-iff-Harm=0}
  If \LapV fails to be essentially self-adjoint, then $\Harm \neq \{0\}$.
  \begin{proof}
    We prove that the boundary form $\bdform(u,v)$ vanishes identically whenever $\Harm=\{0\}$. Since the boundary sum can only be nonzero when $\Harm \neq \{0\}$, the conclusion will follow once we show that
    \linenopax
    \begin{align}\label{eqn:bd-form-as-bd-sum}
      \bdform(u,v) = \tfrac1{2\ii}\sum_{\bd \Graph} \left(\dn{\cj{u}} (\LapV^\ad v) - (\cj{\LapV^\ad u})\dn v \right).
    \end{align}
    To see this, apply Theorem~\ref{thm:E(u,v)=<u,Lapv>+sum(normals)} to obtain
    \linenopax
    \begin{align*}
      \la \LapV^\ad u, v\ra_\energy
      = \sum_{\verts} \LapV^\ad \cj{u} \LapV v + \sum_{\bd \Graph}\LapV^\ad \cj{u} \dn v
      = \sum_{\verts} \LapV \cj{u} \LapV v + \sum_{\bd \Graph}\LapV \cj{u} \dn v
    \end{align*}
    for any $u,v \in \dom (\LapV^\ad)$. The second equality follows because $\Lap^\ad = \Lap$ pointwise:
    \linenopax
    \begin{align*}
      \LapV^\ad u(x) - \LapV^\ad u(o)
      = \la v_x, \LapV^\ad u \ra_\energy
      = \la \LapV v_x, u \ra_\energy
      = \la \gd_x-\gd_o, u \ra_\energy
      = \Lap u(x) - \Lap u(o),
    \end{align*}
    where the last equality comes by Lemma~\ref{thm:<delta_x,v>=Lapv(x)}.
    Also, note that $u \in \dom (\LapV^\ad)$ implies $\LapV^\ad u \in \HE$, so that Theorem~\ref{thm:E(u,v)=<u,Lapv>+sum(normals)} applies and both terms are finite. Consequently, the two sums over \verts cancel and the theorem follows.
  \end{proof}
\end{theorem}

\begin{remark}\label{rem:defect-implies-boundary}
  There is an alternative, more elementary way to prove Theorem~\ref{thm:LapV-not-ess-selfadjoint-iff-Harm=0}. Suppose $w \neq 0$ is a nonzero defect vector with $\LapV^\ad w = \ii w$. Then we can find a representative for $w$ such that
    \linenopax
    \begin{align}\label{eqn:elementary-defect}
      \la w, w \ra_\energy
      &= \sum_{\verts} \cj{w} \Lap w + \sum_{\bd \Graph} \cj{w} \dn w 
      = \ii \sum_{\verts} |w|^2 + \Re \sum_{\bd \Graph} \cj{w} \dn w + \ii \Im \sum_{\bd \Graph} \cj{w} \dn w.
    \end{align}
    Since $\|w\|_\energy^2 = \la w, w \ra_\energy$ is real (and strictly positive, by hypothesis), this implies the boundary sum is nonzero and Theorem~\ref{thm:TFAE:Fin,Harm,Bdy} gives the existence of nontrivial harmonic functions.
    
    It also follows from \eqref{eqn:elementary-defect} that such a nonzero defect vector satisfies
    \linenopax
    \begin{align*}
      \sum_{\verts} |w|^2 = -\Im \sum_{\bd \Graph} \cj{w} \dn w > 0,
    \end{align*}
    so that $\Im \sum_{\bd \Graph} \cj{w} \dn w < 0$.
\end{remark}

\version{}{

\subsection{Self-adjointness of \Lap and the boundary of \Graph}
\label{sec:Self-adjointness-of-Lap-and-boundary-of-G}

\marginpar{No point including this until the theorem is sorted out.}
In this section, we show that \Graph has a boundary precisely when \LapV is not essentially self-adjoint. More precisely, in Theorem~\ref{thm:defect-iff-boundary} we show that \Graph has a boundary precisely when \HE contains a defect vector $w \in \Def_{\LapV}(\ii)$.
\version{}{\marginpar{Do we need to use \LapV specifically, here?}}
Recall from Definition~\ref{def:defect-index} that $\Def_{\LapV}(\ii) = \{v \in \dom \LapV^\ad \suth \LapV^\ad v = \ii v\}$, and recall from Theorem~\ref{thm:essentially-self-adjointness-criterion} that \LapV fails to be essentially self-adjoint precisely when $\Def_{\LapV}(\ii) \neq \es$. In Theorem~\ref{thm:deficiency-criterion-theorem-for-HE} we give a characterization of $\Def_{\LapV}(\ii) \neq \es$ in terms of the quadratic form $Q(\gx) := \sum_{x,y} \cj{\gx_x} \la v_x, v_y\ra_\energy \gx_y$ which may be readily computed for a given \ERN.

\begin{theorem}\label{thm:deficiency-criterion-theorem-for-HE}
  Let $w = \sum_{x \neq o} \gx_x v_x \in \HE$. Then $\LapF^\ad w = \ii w$ if and only if $\Lap \gx = \ii \gx$ and $Q(\gx) := \sum_{x,y} \cj{\gx_x} \la v_x, v_y\ra_\energy \gx_y < \iy$.
  \begin{proof}
    \fwd First, note that $\gx_x = \la \gd_x,w\ra_\energy$ by Lemma~\ref{thm:vx-are-linearly-independent}. We switch between vector and function notation $\gx(x) = \gx_x$ in the following computation:
    \linenopax
    \begin{align*}
      (\Lap \gx)(x)
      &= \cond(x) \gx(x) - \sum_{y \nbr x} \cond_{xy} \gx(y) &&\text{by \eqref{eqn:def:laplacian}} \\
      &= \cond(x) \la \gd_x, w\ra_\energy - \sum_{y \nbr x} \cond_{xy} \la \gd_y, w\ra_\energy  &&\text{Lemma~\ref{thm:vx-are-linearly-independent}} \\
      &= \left\la \cond(x) \gd_x - \sum_{y \nbr x} \cond_{xy} \gd_y, w\right\ra_{\negsp[10]\energy} &&\text{linearity, $\cond_{xy} >0$} \\
      &= \la \Lap \gd_x, w\ra_\energy &&\text{by \eqref{eqn:def:laplacian}} \\
      &= \la \gd_x, \Lap^\ad w\ra_\energy &&\text{Def. of adjoint} \\
      &= \ii \la \gd_x, w\ra_\energy &&\text{hypothesis} \\
      &= \ii \gx(x) &&\text{Lemma~\ref{thm:vx-are-linearly-independent} again}.
    \end{align*}
    To show $Q(\gx) < \iy$, let $F \ci \verts$ be finite and define $w_F = \sum_{x \in F} \gx_x v_x$. Then
    \linenopax
    \begin{align*}
      \|w_F\|_\energy^2
      = \left\la \sum_{x \in F} \gx_x v_x, \sum_{y \in F} \gx_y v_y \right\ra_{\negsp[10]\energy}
      = \sum_{x,y \in F} \cj{\gx_x} \la v_x, v_y\ra_\energy \gx_y
      = Q(\charfn{F} \gx)
      < \iy,
    \end{align*}
    where $\charfn{F}(x)$ is the usual indicator function of $F$. It only remains to see that $\|w_F\|_\energy \leq \|w\|_\energy$, but this follows from the Minimax Principle: let $M = [\la v_x, v_y\ra_\energy]_{x,y}$ and let $M_{F} = [\la v_x, v_y\ra_\energy]_{x,y \in F}$ be an $|F| \times |F|$ submatrix of $M$. Then
    \linenopax
    \begin{align*}
      \|w_F\|_\energy^2 = \la \gx, M_F \gx\ra_\unwtd
      \leq \la \gx, M \gx\ra_\unwtd = \|w\|_\energy^2,
    \end{align*}
    whence $Q(\gx) = \sup_F Q(\charfn{F} \gx) \leq \sup_F \|w_F\|_\energy^2 \leq\|w\|_\energy^2 < \iy$.

    \bwd It is important to note that there is no assumption on whether or not $\gx \in \HE$; this ensures that (ii) always has a solution. The real question is if the $w$ defined via (i) is in \HE. Let $\sF = \{F \ci \verts \suth |F|<\iy\}$ be the set of all finite subsets of vertices. For $F \in \sF$, let $w_F := \sum_{x \in F} \gx_x v_x$. Then
    \linenopax
    \begin{align*}
      \la \gd_y, v_x\ra_\energy = \gd_y(x) - \gd_y(o) = \gd_y(x)
    \end{align*}
    \marginpar{UNFINISHED!}
    shows that $w$ is of the right form. It remains to show $w \in \HE$.
 \end{proof}
\end{theorem}
}

\section{Dual frames and the energy kernel}
\label{sec:Dual-frames-and-the-energy-kernel}
\label{sec:Frames-and-dual-frames}


In previous parts of this book, we have approximated infinite networks by truncating the domain; this is the idea behind the definition of \Fin in Definition~\ref{def:Fin}, and in the use of exhaustions for various arguments (Definition~\ref{def:exhaustion-of-G}). This approach corresponds to a restriction to $\spn\{\gd_x\}_{x \in F}$, where $F$ is some finite subset of \verts. In this section, we consider truncations in the dual variable, i.e., restrictions to sets of the form $\spn\{v_x\}_{x \in F}$. This is directly analogous to the usual time/frequency duality in Fourier theory.

The energy kernel $\{v_x\}$ generally fails to be a frame for \HE, as shown by Lemma~\ref{thm:when-v_x-is-a-frame} and the ensuing remarks. However, things improve when restricting to a finite subset. We shall approach the infinite case via a compatible system of finite dual frames, one for each finite subset $F \ci \verts\less\{o\}$; see Definition~\ref{def:frames-and-dual-frames}. In Theorem~\ref{thm:dual-frames}, we show that $\{\gd_x\}_{x \in F}$ and $\{v_x\}_{x \in F}$ form a dual frame system. 

We obtain optimal frame bounds in Corollary~\ref{thm:optimal-dual-frame-bounds}. In Theorem~\ref{thm:TFAE:LapV-bdd,spectralgap,globalframebound}, we show that the boundedness of \LapV is equivalent to both the existence of a global upper frame bound (i.e., one can let $F \to \Graph$), and the existence of a spectral gap.

We begin with two lemmas whose parallels serve to underscore the theme of this section.

\begin{lemma}\label{thm:vx-are-linearly-independent}
  The vectors $\{v_x\}$ are linearly independent.
  \begin{proof}
    Suppose that we have a (finite) linear combination $\gy = \sum_{x \neq o} \gx_x v_x$, where $\gx_x \in \bC$. Then for $y \neq o$,
    \linenopax
    \begin{equation*}
      \la \gd_y, \gy \ra_\energy
      = \sum_{x \neq o} \gx_x \la \gd_y, v_x \ra_\energy
      = \sum_{x \neq o} \gx_x (\gd_y(x) - \gd_y(o))
      = \sum_{x \neq o} \gx_x \gd_y(x)
      = \gx_y.
    \end{equation*}
    If $\gy=0$, then this calculation shows $\gx_y = 0$ for each $y$.
  \end{proof}
\end{lemma}

\begin{lemma}\label{thm:dx-are-linearly-independent}
  The vectors $\{\gd_x\}$ are linearly independent.
  \begin{proof}
    Suppose that we have a (finite) linear combination $\gy = \sum_{x \neq o} \gx_x \gd_x$, where $\gx_x \in \bC$. Then for $y \neq o$,
    \linenopax
    \begin{equation*}
      \la v_y, \gy \ra_\energy
      = \sum_{x \neq o} \gx_x \la v_y, \gd_x \ra_\energy
      = \sum_{x \neq o} \gx_x (\gd_x(y) - \gd_x(o))
      = \sum_{x \neq o} \gx_x \gd_x(y)
      = \gx_y.
    \end{equation*}
    If $\gy=0$, then this calculation shows $\gx_y = 0$ for each $y$.
  \end{proof}
\end{lemma}

\begin{defn}\label{def:V(F)} \label{def:LapF}
  In this section we always let $F \ci \verts\less\{o\}$ denote a finite subset of vertices and let $V(F) = \spn\{v_x \suth x \in F\}$.
  Observe that elements of $V(F)$ do \emph{not} typically have finite support; cf. Definition~\ref{def:Fin} and Figure~\ref{fig:vx-in-Z1} of Example~\ref{exm:infinite-lattices}. Let $\Lap_{V(F)}$ denote the Laplacian when taken to have the domain $V(F)$, even though it not dense in \HE.
\end{defn}

\begin{defn}\label{def:frames-and-dual-frames}
  Denote $D(F) := \{\gd_x\}_{x \in F}$ and let \LapF be the Laplacian when taken to have this (non-dense) domain.
  
  Then $D(F)$ is a \emph{dual frame} for $V(F)$ if there are constants $0<A \leq B<\iy$ (called \emph{frame bounds}) for which
  \linenopax
  \begin{equation}\label{eqn:dual-frame}
    A \|\gy \|_\energy^2
    \leq \sum_{x \in F} |\la \gd_x,\gy\ra_\energy|^2
    \leq B\|\gy \|_\energy^2, \qq \forall \gy \in V(F).
  \end{equation}
\end{defn}

\begin{lemma}\label{thm:when-v_x-is-a-frame}
  $\{v_x\}$ is a frame for \HE if and only if $\ell^2(\verts)$ and \HE are isomorphic.
  \begin{proof}
    Since $\{v_x\}$ is a reproducing kernel, the frame inequalities take the form
    \linenopax
    \begin{align}\label{eqn:v_x-frame-inequality}
      A \|w\|_2^2 \leq \sum |w(x)|^2 \leq B \|w\|_2^2.
    \end{align}
    Each inequality indicates a (not necessarily isometric) embedding.
  \end{proof}
\end{lemma}

\begin{remark}\label{rem:when-v_x-fails-frame}
  The second inequality fails if \Lap does not have a spectral gap.  
  See also Lemma~\ref{thm:frame-for-HE}.
  \version{}{\marginpar{Is this true on the tree, too?}
  For none of the infinite graph examples in \S\ref{sec:examples}--\S\ref{sec:lattice-networks} does \Lap have a spectral gap. This is discussed further in \S\ref{sec:Frames-and-dual-frames}.}
\end{remark}

\begin{defn}\label{def:M_F}
  Define a Hermitian $|F| \times |F|$ matrix by $M_F := [\la v_x,v_y\ra_\energy]_{x,y \in F}$. Let $\gl_{min} := \min \spec(M_F)$ and $\gl_{max} := \max \spec(M_F)$.
\end{defn}


\begin{defn}\label{def:X^gy}
  For $\gy \in \HE$, define $X:\HE \to \ell(\verts)$, where $\ell(\verts)$ is the space of all functions on \verts, by
  \begin{equation}\label{eqn:def:X^gy}
    X^\gy(x) := \la \gd_x, \gy\ra_\energy.
  \end{equation}
  By Remark~\ref{rem:Lap-defined-via-energy}, $X$ is morally identical to the Laplacian when defined on all of \HE; note that $X^\psi$ may not lie in \HE.
\end{defn}

In the proof of Theorem~\ref{thm:dual-frames}, the notations $\la \cdot\,,\,\cdot\ra_\one$ and $\|\cdot\|_\one$ refer to the space $\ell^2(\one)$ discussed in \S\ref{sec:L2-theory-of-Lap-and-Trans}, that is, $\la f,g\ra_\one = \sum_{x \in \verts} \cj{f(x)} g(x)$ is the unweighted $\ell^2$ inner product, etc.

\begin{theorem}\label{thm:dual-frames}
  For any finite $F$, one has $\gl_{min} > 0$ for the minimal eigenvalue of Definition~\ref{def:M_F}, and $\{\gd_x\}_{x \in F}$ is a dual frame for $V(F)$ with frame bounds
  \linenopax
  \begin{equation}\label{eqn:dual-frame-bounds}
    \frac1{\gl_{max}} \| \gy \|_\energy^2
    \leq \sum_{x \in F} |\la \gd_x,\gy\ra_\energy|^2
    \leq \frac1{\gl_{min}} \| \gy \|_\energy^2.
  \end{equation}
  \begin{proof}
    First, to show that $\gl_{min} > 0$, we show that 0 is not in the spectrum of $M_F$. By way of contradiction, suppose $\exists\gx:F \to \bC$ such that
    \linenopax
    \begin{align*}
      M_F \gx = \sum_y \la v_x, v_y\ra_\energy \gx_y = 0.
    \end{align*}
    The vector $\gy = \sum_y \la v_x, v_y\ra_\energy \gx_y \in V(F)$ is nonzero by Lemma~\ref{thm:vx-are-linearly-independent}, and yet
    \linenopax
    \begin{align*}
      \gy(x) - \gy(o)
      = \la v_x, \gy\ra_\energy
      = \sum_y \la v_x, v_y\ra_\energy \gx_y = 0.
    \end{align*}
    Hence, \gy is constant and therefore $\gy=0$ in \HE. \cont \; So 0 is not in the spectrum of $M_F$.
    Then by \eqref{eqn:def:X^gy},
    \linenopax
    \begin{align*}
      \| \gy \|_\energy^2
      = \left\|\sum_{x \in F} \la \gd_x, \gy\ra_\energy v_x \right\|_\energy^2
      &= \sum_{x,y \in F} \cj{\la \gd_x, \gy\ra_\energy} 
\la v_x, v_y \ra_\energy \la \gd_y, \gy\ra_\energy \\
      &= \sum_{x,y \in F} \cj{X^\gy(x)} \la v_x, v_y \ra_\energy X^\gy(y) \\
      &= \la X^\gy, M_F X^\gy \ra_\one,
    \end{align*}
    whence $\gl_{min} \|X^\gy\|_\one^2 \leq \| \gy \|_\energy^2 \leq \gl_{max} \|X^\gy\|_\one^2$, and the conclusion \eqref{eqn:dual-frame-bounds} follows from $\|X^\gy\|_\one^2 = \sum_{x \in F} |\la \gd_x, \gy\ra_\energy|^2$.
  \end{proof}
\end{theorem}

\begin{cor}\label{thm:optimal-dual-frame-bounds}
  The frame bounds in \eqref{eqn:dual-frame-bounds} are optimal.
  \begin{proof}
    Let $\gx \in \spec(M_F)$ and $\gx:F \to \bC$ with $M_F \gx = \gl \gx$. The vector $\gx = \sum_{x \in F} \gx_x v_x$ is in \HE by the proposition and $\gx = \la \gd_x, \gy\ra_\energy = X^\gy(x)$ for each $x \in F$ by Lemma~\ref{thm:vx-are-linearly-independent} and \eqref{eqn:def:X^gy}. Moreover,
    \linenopax
    \begin{equation*}
      \|\gy\|_\energy^2 = \la \gx, M_F \gx\ra_2 = \gl\|\gx\|_2^2
      = \sum_{x \in F} | \la \gd_x, \gy\ra_\energy|^2.
    \end{equation*}
    We now apply this to $\gl_{min}$ and to $\gl_{min}$ and deduce the bounds are optimal.
  \end{proof}
\end{cor}

In the next lemma, we use \Lap specifically to indicate that the Laplacian is considered pointwise, and without regard to domains.

\begin{lemma}\label{thm:X-intertwines-Lap}
  $X$ represents \LapV on $\ell(\verts)$, i.e., $\Lap(X\gy) = X(\LapV \gy)$ for all $\gy \in \dom \LapV$.
  \begin{proof}
    Fix $\gy \in \dom \LapV $ and $x \in \verts$. Then
    \linenopax
    \begin{align*}
      \Lap (X\gy)(x)
      &= \cond(x) X^\gy(x) - \sum_{y \nbr x} \cond_{xy} X^\gy(y) \\
      &= \cond(x) \la \gd_x, \gy\ra_\energy
- \sum_{y \nbr x} \cond_{xy} \la \gd_y, \gy\ra_\energy \\
      &= \left\la \cond(x) \gd_x - \sum_{y \nbr x} \cond_{xy} \gd_y, \gy\right\ra_{\negsp[10]\energy} \\
      &= \la \Lap \gd_x, \gy \ra_\energy.
    \end{align*}
    Now since $\gd_x \in \dom \LapV^\ad$, we have
    $\Lap (X\gy)(x)
      = \la \LapV^\ad \gd_x, \gy \ra_\energy 
      = \la \gd_x, \LapV \gy \ra_\energy 
      = X(\LapV \gy)(x)$.
  \end{proof}
\end{lemma}

\begin{lemma}\label{thm:LapV-is-semibounded}
  For any $\gy \in V(F)$, we have $\la \gy, \Lap_{V(F)} \gy\ra_\energy = \sum_{x \in F} |X^\gy(x)|^2 + \left|\sum_{x \in F} X^\gy(x)\right|^2$.
  \begin{proof}
    Writing \Lap for $\Lap_{V(F)}$, this follows from
    \linenopax
    \begin{align}
      \la \gy, \Lap \gy\ra_\energy
      &= \sum_{x,y \in F} \cj{X^\gy(x)} X^\gy(y) \la v_x, \Lap v_y\ra_\energy \notag \\
      &= \sum_{x,y \in F} \cj{X^\gy(x)} X^\gy(y) ((\gd_y(x)-\gd_y(o))-(\gd_o(x)-\gd_o(o))) \notag \\
      &= \sum_{x,y \in F} \cj{X^\gy(x)} X^\gy(y) (\gd_y(x)+1) 
      \label{eqn:LapV-is-semibounded} \\
      &= \sum_{x \in F} \cj{X^\gy(x)} X^\gy(x) + \left(\sum_{x \in F} \cj{X^\gy(x)} \right) \left(\sum_{y \in F} X^\gy(y)\right), \notag
    \end{align}
    where \eqref{eqn:LapV-is-semibounded} follows because $o \notin F$.
  \end{proof}
\end{lemma}
Incidentally, Lemma~\ref{thm:LapV-is-semibounded} offers a proof of Theorem~\ref{thm:Lap-Hermitian-on-V}.

\begin{theorem}\label{thm:TFAE:LapV-bdd,spectralgap,globalframebound}
  \version{}{\marginpar{Do lower bound also.}}
  The following are equivalent:
  \begin{enumerate}[(i)]
    \item \LapV is a bounded operator on \HE.
    \item There is a global upper frame bound $B<\iy$ in \eqref{eqn:dual-frame-bounds}, i.e.
\begin{equation}\label{eqn:global-upper-frame-bound}
  \sum_{x \neq o} |\la \gd_x, \gy\ra_\energy|^2
  \leq B \|\gy\|_\energy^2,
  \qq\forall \gy \in \HE.
\end{equation}
    \item There is a spectral gap $\inf \spec(M_F) > 0$, where $F$ runs over the set \sF of all finite subsets of $\verts \less\{o\}$.
  \end{enumerate}
  \begin{proof}
    (i)$\implies$(ii). If \LapV is bounded, then by Lemma~\ref{thm:<delta_x,v>=Lapv(x)} followed by Corollary~\ref{thm:Lap-Hermitian-on-V}, 
    \linenopax
    \begin{align*}
      \sum_{x \neq o} |\la \gd_x, \gy\ra_\energy|^2 
      = \sum_{x \neq o} |\Lap \gy(x)|^2
      = \la \gy, \Lap \gy\ra_\energy
      \leq B \|\gy\|_\energy^2.
    \end{align*}
    
    (ii)$\implies$(i). First fix $\ge>0$. Note that $\sum_{x \in \verts} \Lap \gy(x)=0$ by Corollary~\ref{thm:sum(Lap v)=0}, so choose $F$ so that $\left|\sum_{x \in \verts} \Lap \gy(x) \right| < \ge$.
    The hypothesis of the global upper frame bound $B$ gives 
    \linenopax
    \begin{align*}
      \sum_{x \in F} |\Lap \gy(x)|^2 = \sum_{x \in F} |\la \gd_x, \Lap \gy\ra_\energy|^2 \leq B \|\gy\|_\energy^2,
    \end{align*}
    so that Lemma~\ref{thm:LapV-is-semibounded} implies
    \linenopax
    \begin{align*}
      \left|\la \gy, \Lap \gy\ra_\energy \right|
      \leq \sum_{x \in F} |\Lap \gy(x)|^2 + \left|\sum_{x \in F} \Lap \gy(x)\right|^2
      < B \|\gy\|_\energy^2 + \ge,
    \end{align*}
    and we get $\left|\la \gy, \Lap \gy\ra_\energy \right| \leq B \|\gy\|_\energy^2$ as $\ge \to 0$.
    
    (i)$\iff$(iii) Observe that \eqref{eqn:dual-frame-bounds} and Lemma~\ref{thm:optimal-dual-frame-bounds} imply that $\frac1{\gl_{min}(F)} \leq B$, and hence $\gl_{min}(F) \geq 1/B$, $\forall F \in \sF$.
    If we have an exhaustion $F_1 \ci F_2 \ci \dots \bigcup F_k = \verts\less\{o\}$, then the Minimax Theorem indicates that $\gl_{min}(F_{k+1}) \leq \gl_{min}(F_k)$ so
    \linenopax
    \begin{align*}
      \|\LapV\|^{-1}
      &= \sup \{\tfrac1B \geq 0 \suth \la \gy, \LapV \gy\ra_\energy \leq B \|\gy\|_\energy^2, \forall \gy \in V\} \\
      &= \lim_{k \to \iy} \gl_{min}(F_k).
      &&\qedhere
    \end{align*}
  \end{proof}
\end{theorem}

\begin{cor}\label{thm:dual-frames-bound-resistance-metric}
  If $\{\gd_x\}$ is a dual frame for $\{v_x\}$, then the upper and lower frame bounds $A$ and $B$ provide bounds on the free resistance metric:
  \linenopax
  \begin{equation}\label{eqn:dual-frames-bound-resistance-metric}
    \frac2B \leq R^F(x,y) \leq \frac2A.
  \end{equation}
  Note that as $F$ increases to \verts, one may have $A \to 0$ so that the upper bound tends to \iy.
  \begin{proof}
    By \eqref{eqn:def:R^F(x,y)-energy}, we are motivated to apply the frame inequalities applied to the function $v_x - v_y \in \HE$ via Theorem~\ref{thm:TFAE:LapV-bdd,spectralgap,globalframebound}:
    \linenopax
    \begin{align*}
      A\|v_x - v_y\|_\energy^2
      \leq \sum_{z \in \verts\less\{o\}} |\la \gd_z, v_x-v_y\ra_\energy|^2
      \leq B\|v_x - v_y\|_\energy^2.
    \end{align*}
   The result now follows by \eqref{eqn:def:R^F(x,y)-energy} upon observing that $\sum_{z \in \verts\less\{o\}} |\la \gd_z, v_x-v_y\ra_\energy|^2 = 2$.
    \end{proof}
  \end{cor}


\begin{lemma}\label{thm:Harm-does-not-intersect-V}
  For finite $F \ci \verts$, $\Harm \cap V(F) = \es$. A fortiori, $\LapV$ has a spectral gap.
  \begin{proof}
    Let $h = \sum_{i=1}^n c_i v_{x_i}$. If $h$ is harmonic, then
    \linenopax
    \begin{align*}
      0 = \Lap h = \sum c_i (\gd_{x_i}-\gd_o) = \sum c_i \gd_{x_i} - \gd_o\sum c_i,
    \end{align*}
    which implies $c_i =0$ for each $i$, since the Dirac masses are linearly independent vectors. The second claim follows because 0 is not in the point spectrum of \LapV on the finite-dimensional space $V$.
  \end{proof}
\end{lemma}

The symmetry of formula \eqref{eqn:<Lap-gd_x,gd_y>-formula} in $x$ and $y$ provides another proof that \LapF is Hermitian.

\begin{lemma}\label{thm:<Lap-gd_x,gd_y>-formula}
  For all $x,y \in \verts$,
  \linenopax
  \begin{equation}\label{eqn:<Lap-gd_x,gd_y>-formula}
    \la \Lap \gd_x, \gd_y \ra_\energy
    = - (\cond(x) + \cond(y)) \cond_{xy} + \sum_{z \nbr x,y} \cond_{xz} \cond_{zy}.
  \end{equation}
  \begin{proof}
    For \LapF, use $z \nbr x,y$ to denote that $z$ is a neighbour of both $x$ and $y$, and compute 
    \linenopax
    \begin{align*}
      \la \Lap \gd_x, \gd_y \ra_\energy
      &= \left\la \cond(x) \gd_x - \sum_{z \nbr x} \cond_{xz} \gd_z, \gd_y \right\ra_{\negsp[12]\energy} \\
      &= \cond(x) \la \gd_x, \gd_y \ra_\energy
      - \sum_{z \nbr x} \cond_{xz} \la \gd_z, \gd_y \ra_\energy \vstr[4]\\
      &= \cond(x) \la \gd_x, \gd_y \ra_\energy - \cond_{xy} \la \gd_y, \gd_y \ra_\energy
      - \sum_{\substack{z \nbr x \\z \neq y}} \cond_{xz} \la \gd_z, \gd_y \ra_\energy \vstr[4] \\
      &= - \cond(x) \cond_{xy} - \cond_{xy} \cond(y)
      + \sum_{\substack{z \nbr x \\z \neq y}} \cond_{xz} \cond_{zy} \\
      &= - (\cond(x) + \cond(y)) \cond_{xy}
      + \sum_{z \nbr x,y} \cond_{xz} \cond_{zy}.
      \qedhere
    \end{align*}
  \end{proof}
\end{lemma}

\begin{defn}\label{def:cond_x}
  Let $\cond_x$ be defined by $\cond_x(y) = \cond_{xy}$, so $\cond_x^2 := \cond_x \cdot \cond_x := \sum_{y \nbr x} \cond_{xy}^2$.
\end{defn}

In the next theorem, one would need to consider $\la \gf,\Lap \gf\ra_\energy$ for general $\gf \in \spn\{v_x\}$, rather than just $\gf = v_x,\gd_x$, in order to get the full spectrum $[\spec \LapV]$.

\begin{lemma}\label{thm:spec(Lap)-on-HE}
  The spectrum of \LapV satisfies 
  \linenopax
  \begin{align*}
    \inf \spec \LapV \leq \min \left\{\inf \frac{2}{R^F(x,o)}, \inf\left( \cond(x)+\frac{\cond_x^2}{\cond(x)}\right)\right\},\\
    \sup \spec \LapV \geq \max \left\{\sup \frac{2}{R^F(x,o)}, \sup\left(\cond(x)+\frac{\cond_x^2}{\cond(x)}\right)\right\}.
    \end{align*}
  \begin{proof}
    We compute the action of \Lap on certain unit vectors:
    \linenopax
    \begin{align*}
      \left\la \frac{v_x}{\|v_x\|_\energy}, \Lap \frac{v_x}{\|v_x\|_\energy}\right\ra_{\negsp[6]\energy}
      &= \frac1{\|v_x\|_\energy^2} \left((\gd_x(x)-\gd_x(o)) - \vstr[2.3](\gd_x(o) - \gd_o(o))\right)
      = \frac2{R^F(x,o)}
    \end{align*}
    and
    \linenopax
    \begin{align*}
       \left\la \frac{\gd_x}{\|\gd_x\|_\energy}, \Lap \frac{\gd_x}{\|\gd_x\|_\energy}\right\ra_{\negsp[6]\energy}
       &= \frac1{\|\gd_x\|_\energy^2} \left\la \gd_x, \Lap \gd_x \right\ra_{\energy}
       = \frac{\cond(x)^2 + \cond_x^2}{\cond(x)}
       = \cond(x) + \frac{\cond_x^2}{\cond(x)}.
    \end{align*}
    We then apply the well-known theorem that for a closed Hermitian operator $S$,
    \[[\spec S] = \{\la u, Su\ra \suth u \in \dom S, \|u\|=1\},\]
    where $[set]$ denotes the closed convex hull of $set$ in \bC.
    Note that $[\spec S] \ci \bR$. 
      \version{}{\\ \marginpar{Unfinished!} So far so good ...}
  \end{proof}
\end{lemma}

\version{}{
  \begin{cor}\label{thm:Lap-gap-iff-R-bdd}
    \marginpar{LIES?!?}
    \Lap has a spectral gap if and only if $R^F(x,o)$ is bounded.
    \begin{proof}
      Probably need $R^F(x,o)$ is bounded and $\cond(x)+\frac{\cond_x^2}{\cond(x)}$ bounded below (which is not the case if $\cond(x) \to 0$ at infinity).
    \end{proof}
  \end{cor}
}

\section{Remarks and references}
\label{sec:Remarks-and-References-Lap-on-HE}

The family of operators covered by what we here refer to as the Laplacian \Lap is large, and the literature both large and diverse; for example these operators in mathematical physics go by the name Òdiscrete Schroedinger operators.Ó Readable introductions include \cite{CdV99}, \cite{Chu01}, \cite{Dod06}, \cite{Soardi94}, and \cite{Web08}.
The reader may also find the references \cite{HKLW07, Woess00, Woe03, RS95, JostKendallMoscoRocknerSturm} to be useful.

\version{}{
\subsection{The transfer operator on \HE}
\label{sec:The-transfer-operator-on-HE} 

\begin{theorem}\label{thm:adjoint-of-transfer-operator-on-HE}
  The adjoint of $\Trans_\unwtd:\HE \to\HE$ is $\Lap_\unwtd^{-1} \Trans_\unwtd \Lap_\unwtd$.
  \begin{proof}
    Working directly with the inner product,
    \linenopax
    \begin{align*}
      \la \Lap_\unwtd^{-1} \Trans \Lap_\unwtd u, v\ra_\energy
      &= \la \Lap_\unwtd^{-1} \Trans \Lap_\unwtd u, \Lap_\unwtd v\ra_\unwtd
      = \la \Trans \Lap_\unwtd u, v\ra_\unwtd
      = \la \Lap_\unwtd u, \Trans v\ra_\unwtd
      = \la u, \Trans v\ra_\energy.
      \qedhere
    \end{align*}
  \end{proof}
\end{theorem}

Recall from Definition~\ref{def:graph-of-operator} that the graph of \Lap is
\version{}{\marginpar{Does this require $\Lap v \in \ell^2$?}}
\begin{equation}\label{eqn:graph-of-Lap}
  Graph(\Lap) := \{\left[\begin{smallmatrix}v \\ \Lap v\end{smallmatrix}\right] \suth v \in \ell^2(\unwtd)\} \ci \ell^2 \oplus \ell^2.
\end{equation}

\begin{theorem}\label{thm:embedding-Lap-graph-into-HE-is-contractive}
  \version{}{\marginpar{How is this result used? What is it for?}}
  The map $\gf:Graph(\Lap) \to \HE$ by $\gf:\left[\begin{smallmatrix}v \\ \Lap v\end{smallmatrix}\right] \mapsto v$ is contractive. If \Graph is infinite, then $\ker \gf = 0$.
  \begin{proof}
  \version{}{\marginpar{Want to apply $\gf^\ad:\HE \to Graph_\unwtd(\Lap)$? \\Claim: $\HE \cap \dom \Lap^\ad = \dom \tilde \Lap$.}}
    Observe that $v \in \dom \gf$ iff $v, \Lap v \in \ell^2(\unwtd)$, so  Lemma~\ref{thm:E(u,v)=<u,Lapv>} guarantees $v \in \HE$. Contractivity follows by combining \eqref{eqn:E(u,v)=<u,Lapv>} with the Schwarz inequality:
    \linenopax
    \begin{align*}
      \|v\|_\energy^2
      = \la v, \Lap v\ra_\unwtd
      \leq \|v\|_\unwtd \|\Lap v\|_\unwtd
      \leq \|v\|_\unwtd^2 + \|\Lap v\|_\unwtd^2.
    \end{align*}
    The map \gf can only fail to be injective if there exist elements of $\ell^2(\unwtd)$ which differ only by a constant. However, this cannot happen when \Graph is infinite, as $\sum_x (v(x)+k)$ can be finite for at most one value of $k$.
  \end{proof}
\end{theorem}

\begin{cor}\label{thm:hypoellipticity-in-HE}
  Assume $(\Graph, \cond)$ satisfies the Powers bound \eqref{eqn:def:power-bound}. Then $\|\Lapu v\|_\energy \leq \uBd \cdot \|v\|_\energy$. If \Lapu has a spectral gap $\gd := \inf \spec \Lapu >0$, then $\|v\|_\energy \leq \gd\|\Lapu v\|_\energy$.
  \begin{proof}
    \fwd Theorem~\ref{thm:E(u,v)=<u,Lapv>} gives
    \linenopax
    \begin{align*}
      \|\Lapu v\|_\energy^2
      = \la \Lapu v, \Lapu^2 v\ra_\unwtd
      = \|\Lapu^{3/2} v\|_\unwtd^2
      \leq \|\Lapu\|^2 \cdot \|\Lapu^{1/2} v\|_\unwtd^2
      = \|\Lapu\|^2 \cdot \|v\|_\energy^2.
    \end{align*}

    \bwd Let \gd denote the spectral gap as in the hypotheses. The Powers bound ensures an upper bound $K$ on the spectrum of \Lapu. Since $s \geq \frac12$ if and only if $s \leq \gd s^2$, we have
    \linenopax
    \begin{align*}
      \|v\|_\energy
      &= \|\Lapu^{1/2}v\|_\unwtd
      = \int_\gd^K s |\hat v|^2 \,d\gx
      \leq \int_\gd^K \gd s^2 |\hat v|^2 \,d\gx
      = \gd\|\Lapu v\|_\unwtd.
      \qedhere
    \end{align*}
  \end{proof}
\end{cor}

\begin{cor}\label{thm:hypoellipticity-in-L2}
  Assume $(\Graph, \cond)$ satisfies the Powers bound \eqref{eqn:def:power-bound}. Then $v \in \ell^2(\unwtd)$ if and only if $\Lap v \in \ell^2(\unwtd)$.
\end{cor}

One may also define the action of \Lap on a \fn $v \in \HE$ by
\begin{equation}\label{eqn:def:Lap-extension}
  (\Lap v)(u) := \energy(u,v) = \la u, \Lap v \ra_\unwtd.
\end{equation}

\begin{prop}\label{thm:Lap-extension}
  The formula \eqref{eqn:def:Lap-extension} extends the usual definition of \Lap as defined on $\clspn{v_x} \ci \ell^2(\unwtd)$, where $\Lap v_x = \gd_x - \gd_o$.
  \begin{proof}
    Suppose $v \in \HE$ and $u=\gd_x$. Then
    \linenopax
    \begin{align*}
      (\Lap v)(u) &= \la \gd_x, \Lap v \ra_\spectral = \Lap v(x).
      \qedhere
    \end{align*}
  \end{proof}
\end{prop}

}


\version{}{\section{Green function and operator}
\label{sec:Green-function-and-operator}

In Remark~\ref{rem:appearance-of-Greens-oper-in-vNeu-constr} we mentioned that the unitary $U$ in Theorem~\ref{thm:unitary-vN-to-H_E} is essentially the \emph{Green operator} \Wiener which acts on an element of \HE by
\linenopax
\begin{align}\label{eqn:Greens-oper-recalled}
  (\Wiener u)(x) := \sum_{x \in \verts} g(x,y) u(y),
\end{align}
where the \emph{Green function} is
\linenopax
\begin{align}\label{eqn:Greens-fn-recalled}
  g(x,y) := \la v_x, v_y\ra_\energy.
\end{align}
In combination with Corollary~\ref{thm:vx-is-Lipschitz}, this shows the Green function is Lipschitz continuous, a key result of \cite[Thm.~4.5]{Kig03} in a different setting.
}

  \version{}{\marginpar{One consequence of this is that Theorem~\ref{thm:Pot-is-nonempty-by-current-flows} can be interpreted as proving the existence of Green's functions! We don't care so much about this (yet) but others will ...}
For computations, it is often easiest to work with $g(x,y) = v_x(y)$, by selecting the representative of $g$ which vanishes at $y=o$. In this case, one can easily verify that
\linenopax
\begin{align*}
  \Lap \Wiener u(x)
  &= \sum_{x \in \verts} \Lap g(x,y) u(y)
  = \sum_{x \in \verts} \Lap v_x(y) u(y)
  = u(x) - u(o),
\end{align*}
so that $\Lap \Wiener u(x) = u(x)$ in \HE.
}

%% file: ell2-of-Lap-and-Trans.tex

\chapter{The $\ell^2$ theory of \Lap and the transfer operator}
\label{sec:L2-theory-of-Lap-and-Trans}

\headerquote{One geometry cannot be more true than another; it can only be more convenient.}{---~H.~Poincare}

This chapter is devoted to the study of the graph Laplacian \Lap, and the transfer operator \Trans, when considered as acting on the space $\ell^2(\verts)$ of square-summable functions on the vertices. This is $\ell^2(\verts) = \ell^2(\verts,\gm)$ where \gm is the counting measure, and the operators \Lap and \Trans have a profoundly different spectral theory with respect to the $\ell^2$ inner product.

\section{$\ell^2(\verts)$}

In this section, we discuss results for \Lap and \Trans when considered as operators on 
\linenopax
\begin{equation}\label{eqn:defn:l2}
  \ell^2(\unwtd) := 
  \{u:\verts \to \bC \suth {\textstyle\sum\nolimits_{x \in \verts}} |u(x)|^2 < \iy\}, 
\end{equation}
with the inner product 
\linenopax
\begin{equation}\label{eqn:defn:l2-inner-prod}
  \la u, v\ra_\unwtd := \sum_{x \in \verts} u(x) v(x).
\end{equation}
  \glossary{name={$\la\cdot,\cdot\ra_\unwtd$},description={unweighted $\ell^2$ inner product},sort=l,format=textbf}
  \glossary{name={$\ell^2(\one)$ or $\ell^2(\verts,\one)$},description={unweighted $\ell^2$ space},sort=l,format=textbf}
The constant function \one appears in the notation to specify the weight involved in the inner product, in contrast to \cond. This is necessary because we will also be interested in \Lap and \Trans as operators on 
\linenopax
\begin{equation}\label{eqn:defn:l2c}
  \ell^2(\cond) := 
  \{u:\verts \to \bC \suth {\textstyle\sum\nolimits_{x \in \verts}} \cond(x) |u(x)|^2 < \iy\},
\end{equation}
with the inner product 
\linenopax
\begin{equation}\label{eqn:defn:l2c-inner-prod}
  \la u, v\ra_\cond := \sum_{x \in \verts} \cond(x) u(x) v(x).
\end{equation}
  \glossary{name={$\la\cdot,\cdot\ra_\cond$},description={\cond-weighted $\ell^2$ inner product},sort=<,format=textbf}
  \glossary{name={$\ell^2(\cond)$ or $\ell^2(\verts,\cond)$},description={\cond-weighted $\ell^2$ space},sort=<,format=textbf}

While the pointwise definition of \Lap and \Trans remains the same on $\ell^2(\unwtd)$ and  $\ell^2(\cond)$, they are different operators with different domains and different spectra! It is important to keep in mind that in general, none of \HE, $\ell^2(\unwtd)$ or $\ell^2(\cond)$ are contained in any of the others. However, we provide some conditions under which embeddings exist in \S\ref{sec:relating-HE-and-L2}.
We give only some selected results, as this subject is well-documented elsewhere in the literature.


In \S\ref{sec:weighted-spaces}, we consider a map $J:\ell^2(\cond) \to \HE$ is the quotient map induced by the equivalence relation discussed in Remark~\ref{rem:elements-of-HE-are-technically-equivalence-classes}. It turns out that $J$ is an embedding of $\ell^2(\cond)$ into \Fin, and that its range is dense in \Fin. We will also see that \Prob is self-adjoint on $\ell^2(\cond)$, even though it is not even Hermitian on $\ell^2(\unwtd)$ or \HE except when $\cond(x)$ is constant.

\section{The Laplacian on $\ell^2(\unwtd)$}
\label{sec:self-adjointness-of-the-Laplacian}

In this section, we investigate certain properties of the Laplacian on $\ell^2(\unwtd)$, including self-adjointness and boundedness. Dealing with unbounded operators always requires a bit of care; the reader is invited to consult Appendix~\ref{sec:self-adjointness-for-unbounded-operators} to refresh on some principles of self-adjointness of unbounded operators.  Recall that for $S$ to be \emph{self-adjoint}, it must be Hermitian \emph{and} satisfy $\dom S = \dom S^\ad$, where
\linenopax
\begin{align*}
  \dom S^\ad := \{v \in \sH \suth |\la v, Su\ra| \leq K_v \|u\|, \forall u \in \dom S\}.
\end{align*}
In the unbounded case, it is not unusual for $\dom S \subsetneq \dom S^\ad$. Some good references for this section are \cite{Jor78,vN32a,Nel69, ReedSimonII,Rud91,DuSc88}.

Due to Corollary~\ref{thm:nontrivial-harmonic-fn-is-not-in-L2}, we can ignore the possibility of nontrivial harmonic functions while working in this context. Combining Theorem~\ref{thm:E(u,v)=<u,Lapv>+sum(normals)} with Theorem~\ref{thm:TFAE:Fin,Harm,Bdy}, one can relate the inner products of \HE and $\ell^2(\unwtd)$ by
\begin{equation}\label{eqn:rem:energy-as-l2-inner-prod-with-Lap}
  \la u, \Lap v\ra_\unwtd = \la u, v\ra_{\energy},
\end{equation}
for all $u,v \in \spn\{\gd_x\}$. Observe that $\spn\{\gd_x\}$ is dense in $\ell^2(\unwtd)$ with respect to \eqref{eqn:defn:l2-inner-prod}, and dense in \HE in the \energy norm when $\Harm=0$. Then \eqref{eqn:rem:energy-as-l2-inner-prod-with-Lap} immediately implies that the Laplacian is Hermitian on $\ell^2(\unwtd)$ because, again for all $u,v \in \spn\{\gd_x\}$,
\linenopax
\begin{align}\label{eqn:l2-Hermitian}
  \la u, \Lap v\ra_\unwtd = \la u,v \ra_{\energy} = \cj{\la v, u\ra}_{\energy} = \cj{\la v, \Lap u\ra}_\unwtd = \la \Lap u, v \ra_\unwtd,
\end{align}
This may seem trivial, but it turns out that \Lap is not Hermitian on $\ell^2(\cond)$; cf. Lemma~\ref{thm:Lap-not-Hermitian-on-l2(c)}.

Theorem~\ref{thm:Lap-is-bounded-selfadjoint-on-L2c} shows that if \cond is uniformly bounded \eqref{eqn:def:power-bound}, then \Lap is a bounded operator and hence self-adjoint. However, in Theorem~\ref{thm:Laplacian-is-essentially-self-adjoint} we are able to obtain a much stronger result, without assuming any bounds: the Laplacian on any \ERN is \emph{essentially self-adjoint} on $\ell^2(\unwtd)$. (Recall that \Lap is essentially self-adjoint iff it has a unique self-adjoint extension; cf.~Definition~\ref{def:essentially-self-adjoint}.) This is a sharp contrast to the case for \HE, as seen from Theorem~\ref{thm:LapV-not-ess-selfadjoint-iff-Harm=0}. In the latter parts of this section, we also derive several applications of Theorem~\ref{thm:Laplacian-is-essentially-self-adjoint}.

\subsection{The Laplacian as an unbounded operator}

We begin with the operator \Lap defined on $\spn\{\gd_x\}$, the dense domain consisting of functions with finite support. Then let \Lapu denote the closure of \Lap with respect to \eqref{eqn:defn:l2-inner-prod}, that is, its minimal self-adjoint extension to $\ell^2(\unwtd)$. Some good references for this section are \cite{vN32a,Rud91,DuSc88}.

\begin{lemma}\label{thm:Lap-is-semibounded-on-l2}
  The Laplacian \Lapu is semibounded on $\dom \Lapu$. \emph{A fortiori}, for any $u,v \in \ell^2(\unwtd)$, 
    \linenopax
    \begin{align}\label{eqn:semibounded-innerprod-l2-identity}
      \la u, \Lapu v \ra_{\unwtd}
      = \sum_{x \in \verts} \cond(x) \cj{u(x)} v(x)
      - \sum_{x,y \in \verts} \cond_{xy} \cj{u(x)} v(y)
    \end{align}
  \begin{proof}
    For any $u, v \in \Fin$, a straightforward computation shows
    \linenopax
    \begin{align}\label{eqn:semibounded-l2-identity}
      \la u, \Lapu u \ra_{\unwtd}
      = \sum_{x \in \verts} \cond(x) |u(x)|^2
      - \sum_{x,y \in \verts} \cond_{xy} \cj{u(x)} u(y),
    \end{align}
    whence the equality in \eqref{eqn:semibounded-innerprod-l2-identity} follows by taking limits and polarizing. To see that \Lapu is semibounded, apply the Schwarz inequality first with respect to $y$, then with respect to $x$, to compute
    \linenopax
    \begin{align*}
      \left|\sum_{x \in \verts} \sum_{y \nbr x} \cond_{xy} \cj{u(x)} u(y)\right|
      &\leq \sum_{x \in \verts} \sqrt{\cond(x)} |u(x)| \left[\sum_{y \nbr x} \cond_{xy} |u(y)|^2 \right]^{1/2} \\
      &\leq \left[\sum_{x \in \verts} \cond(x) |u(x)|^2\right]^{1/2} \left[ \sum_{x,y \in \verts} \cond_{xy} |u(y)|^2 \right]^{1/2} \\
      &= \sum_{x \in \verts} \cond(x) |u(x)|^2,
    \end{align*}
    so that the difference on the right-hand side of \eqref{eqn:semibounded-l2-identity} is nonnegative.
      \end{proof}
\end{lemma}

\begin{theorem}\label{thm:Laplacian-is-essentially-self-adjoint}
  If $\deg(x) < \iy$ for every $x \in \verts$, \Lapu is essentially self-adjoint on $\ell^2(\unwtd)$.
  \begin{proof}
    Lemma~\ref{thm:Lap-is-semibounded-on-l2} shows \Lapu is semibounded on $\ell^2(\unwtd)$, so by Theorem~\ref{thm:essentially-self-adjointness-criterion}, it suffices to show the implication
    \linenopax
    \begin{equation}\label{eqn:thm:essentially-self-adjointness-criterion}
      \Lapu^\ad v = -v \q\implies\q v=0, \qq v \in \dom \Lapu ^\ad.
    \end{equation}    
    Suppose that $v \in \ell^2$ is a solution to $\Lapu^\ad v = -v$. Then clearly $\Lapu^\ad v \in \ell^2$, and then by Lemma~\ref{thm:bandedness-sublemma},
    \linenopax
    \begin{align*}
      0 \leq 
       \la v, M_{\Lapu} v\ra_\unwtd
      &= \la v, \Lapu^\ad v\ra_\unwtd
      = -\la v, v\ra_\unwtd
      = -\|v\|_\unwtd^2
      \leq 0
      \q\implies\q
      v=0, 
    \end{align*}
    where $M_{\Lapu}$ is the matrix of \Lapu in the ONB $\{\gd_x\}_{x \in \verts}$. To justify the first inequality, consider that we may find a sequence $\{v_n\} \ci \Fin$ with $\|v-v_n\|_\unwtd \to 0$. Because the matrix $M_{\Lapu}$ is banded, this is sufficient to ensure that $M_{\Lapu} v_n \to M_{\Lapu} v$ and hence $(v_n,M_{\Lapu} v_n)$ converges to $(v,M_{\Lapu} v)$ in the graph norm, and so $\la v_n, M_{\Lapu} v_n\ra_\unwtd$ converges to $\la v, M_{\Lapu} v\ra_\unwtd$. Then $\la v_n, M_{\Lapu} v_n\ra_\unwtd = \energy(v_n) \geq 0$ for each $n$, and positivity is maintained in the limit (even though $\lim \energy(v_n)$ may not be finite).
  \end{proof}
\end{theorem}

See \cite{Web08} for a similar result. 
It follows from Theorem~\ref{thm:Laplacian-is-essentially-self-adjoint} that the closure of the operator \Lapu is self-adjoint on $\ell^2(\unwtd)$, and hence has a unique spectral resolution, determined by a projection valued measure on the Borel subsets of the infinite half-line $\bR_+$. This is in sharp contrast with the continuous case; in Example~\ref{exm:eigenvectors-at-infinity} we illustrate this by indicating how $\Lapu = - \frac{d^2}{dx^2}$ fails to be an essentially self-adjoint operator on the Hilbert space $L^2(\bR_+)$.

\begin{remark}\label{rem:Lap-is-banded}
  The matrix for the operator \Lapu on $\ell^2(\unwtd)$ is \emph{banded} (cf.~\S\ref{sec:banded-matrices}):
  \linenopax
  \begin{align}\label{eqn:L2-Laplacian-is-banded}
    M_{\Lapu}(x,y) = \la \gd_x,\Lapu \gd_y\ra_\unwtd
    =
    \begin{cases}
      \cond(x), &y=x,\\
      - \cond_{xy}, &y \nbr x, \\
      0, &\text{else}.
    \end{cases}
  \end{align}
  The bandedness of $M_{\Lapu}$ is a crucial element of the above proof; Example~\ref{exm:eigenvectors-at-infinity} shows how this proof technique can fail without bandedness. See also Remark~\ref{rem:Lap-is-banded} and  Example~\ref{exm:banded-nonselfadjoint-matrix} for what can go awry without bandedness.

  However, bandedness is not sufficient to guarantee essential self-adjointness. In fact, see Example~\ref{exm:banded-nonselfadjoint-matrix} for a Hermitian operator on $\ell^2$ which is not self-adjoint, despite having a \emph{uniformly banded} matrix, that is, there is some $n \in \bN$ such that each row and column has no more than $n$ nonzero entries.
  The essential self-adjointness of \Lapu in this context is likely a manifestation of the fact that the banding is geometrically/topologically local; the nonzero entries correspond to the vertex neighbourhood of a point in \verts.
\end{remark}

\subsection{The spectral representation of \Lap}
It is clear from Lemma~\ref{thm:converse-to-E(u,v)=<u,Lapv>} that $v \in \dom \energy$ whenever $v, \Lapu v \in \ell^2(\unwtd)$. However, this condition is not necessary, and the precise characterization of $\dom\energy$ is more subtle.

\begin{theorem}\label{thm:HE=Lap1/2}
  For all $u \in \ell^2(\unwtd) \cap \dom\energy$, $\|u\|_\energy = \|\hat\Lap^{1/2} \hat u\|_2$. Therefore, \HE can be characterized in terms of the spectral resolution of \Lap as
  \linenopax
  \begin{equation}\label{eqn:thm:H_energy=dom(sqrt(Lap))}
    \ell^2(\unwtd) \cap \dom\energy
    = \{v: \verts \to \bC \suth \|\widehat{\Lap}_\unwtd^{1/2} \hat v\|_2 < \iy\},
  \end{equation}
  where $\hat v$ is the image of $v$ in the spectral representation of \Lapu.
  \begin{proof}
    Theorem~\ref{thm:Laplacian-is-essentially-self-adjoint} also gives a spectral resolution
    \linenopax
    \begin{align}
      \Lap = \int \gl E(d\gl), \q E:\sB(\bR_+) \to Proj(\ell^2).
    \end{align}
    Applying the functional calculus to the Borel function $r(x) = \sqrt x$, we have
    \linenopax
    \begin{align}\label{eqn:dom-Lap1/2}
      {\Lap}^{1/2} = \int \gl^{1/2} E(d\gl),
      \qq \dom {\Lap}^{1/2} = \{v \in \ell^2 \suth \int |\gl| \cdot \|E(d\gl)v\|^2 < \iy\}.
    \end{align}
    This gives $v \in \ell^2(\unwtd) \cap \dom\energy$ if and only if $v+k \in \dom {\Lap}^{1/2}$ for some $k \in \bC$. However, $\Lap(v+k) = \Lap v$, so the same is true for 
    ${\Lap}^{1/2}$ by the functional calculus.
  \end{proof}
\end{theorem}

\begin{remark}\label{rem:domeE-not-domLap12}
  It is important to observe that $\dom\energy$ is not simply the spectral transform of $\dom \hat\Lap^{1/2} = \{\hat v \suth \hat v \in L^2$ and $\|\hat \Lap^{1/2} \hat v\|<\iy\}.$ The restriction $\hat v \in L^2$ must be removed because there are many functions of finite energy which do not correspond to $L^2$ functions. For an elementary yet important example, see Figure~\ref{fig:vx-in-Z1} of Example~\ref{exm:Z-not-bounded}. 
  Indeed, recall from Corollary~\ref{thm:nontrivial-harmonic-fn-is-not-in-L2} that no nontrivial harmonic function can be in $\ell^2$; see Example~\ref{exm:binary-tree:nontrivial-harmonic}. In this, example $v$ is equal to the constant value $1$ on one infinite subset of the graph, and equal to the constant value $0$ on another.
\end{remark}

\begin{remark}\label{rem:LapTilde-doesn't-see-constants}
  For the example of the integer lattice \bZd, Remark~\ref{rem:HE(Zd)-lies-in-L2(Zd)} shows quite explicitly why the addition of a constant to $v \in \HE$ has no effect on the spectral (Fourier) transform. In this example, one can see directly that addition of a constant $k$ before taking the transform corresponds to the addition of a Dirac mass after taking the transform. As the Dirac mass is supported where the transform of the function vanishes, it has no effect.
\end{remark}

We can also give a reproducing kernel for \Lap on $\ell^2(\unwtd)$. Recall from \eqref{eqn:def:graph-neighbours} that the \emph{vertex neighbourhood} of $x \in \verts$ is $\vnbd(x) := \{y \in \verts \suth y \nbr x\} \ci \verts$. Also recall from Definition~\ref{def:ERN} that $x \notin \vnbd(x)$ and from Definition~\ref{def:conductance-measure} that the conductance function is $\cond(x) := \sum_{y \nbr x} \cond_{xy}$.

\begin{lemma}\label{thm:Lap-has-reproducing-kernel-on-L2}
  The functions $\{\Lap \gd_x\}_{x \in \verts} = \{\cond(x) \gd_x - \cond_{(x\cdot)} \charfn{\vnbd(x)}\}_{x \in \verts}$ give a reproducing kernel for \Lap on $\ell^2(\unwtd)$.
  \begin{proof}
    Since $\la \gd_x, u \ra_\unwtd = u(x)$, the result follows by 
    \linenopax
    \begin{align*}
      \Lap v(x)
      = \cond(x) v(x) - \sum_{y \nbr x} \cond_{xy} v(y)
      = \la \cond(x) \gd_x, v\ra_\unwtd - \la \cond_{(x \cdot)}\charfn{G(x)},v\ra_\unwtd
      = \la \Lap \gd_x, v \ra_\unwtd.
    \end{align*}
    This is a recapitulation of \eqref{eqn:L2-Laplacian-is-banded}. Since $\cond(x) < \iy$, it is clear that $\Lap \gd_x \in \ell^2(\unwtd)$.
  \end{proof}
\end{lemma}

\section{The transfer operator}
\label{sec:the-transfer-operator}
\label{sec:The-transfer-operator-on-L2c}

\begin{defn}\label{def:power-bound}
  We say the graph $(\Graph, \cond)$ \sats the \emph{Powers bound} iff
  \linenopax
  \begin{equation}\label{eqn:def:power-bound}
    \|\cond\|_\iy := \sup_{x \in \verts} \cond(x)
       < \iy.
  \end{equation}
\end{defn}
  \glossary{name={$\|\cond\|$},description={operator norm of \cond, also $\sup \cond(x)$},sort=M,format=textbf}
The terminology ``Powers bound'' stems from \cite{Pow76b}, wherein the author uses this bound to study the emergence of long-range order in statistical models from quantum mechanics. Our motivation is somewhat different, and most of our results do not require such a uniform bound. However, when satisfied, it implies the boundedness of the graph Laplacian (and hence its self-adjointness) and the compactness of the associated transfer operator; see \S\ref{sec:the-transfer-operator}.

The fact that the Powers bound entails the inclusion $\ell^2(\unwtd) \ci \HE$ (see Theorem~\ref{thm:energy-bounded-by-L2}) illustrates how strong this assumption really is. 
While the Laplacian may be unbounded for infinite networks in general, Theorem~\ref{thm:Lap-is-bounded-selfadjoint-on-L2c} gives one situation in which \Lap is always bounded. To see sharpness, note that this bound is obtained in the integer lattices of Example~\ref{exm:infinite-lattices}. In particular, for $d=1$, we have $\|\Lap\| = \sup |4(\sin^2\frac t2)| = 4 = 2 \|\cond\|$.

\begin{theorem}\label{thm:Lap-is-bounded-selfadjoint-on-L2c}
  As an operator on $\ell^2(\unwtd)$, the Laplacian satisfies $\|\Lap\|_\unwtd \leq 2\uBd$, and hence is a bounded self-adjoint operator whenever the Powers bound holds. Moreover, this bound is sharp.
  \begin{proof}
    Since $\Lap = \cond - \Trans$, this is clear by the following lemma.
  \end{proof}
\end{theorem}

Recall from Definition~\ref{def:graph-transfer-operator} that the \emph{transfer operator} \Trans acts on an element of $\dom \Lap$ by
\begin{equation}\label{eqn:def:transfer-operator-on-HE}
  (\Trans v)(x) := \sum_{y \nbr x} \cond_{xy} v(y).
\end{equation}
One should not confuse \Trans with the (bounded) probabilistic transition operator $\Prob = \cond^{-1}\Trans$; recall that the Laplacian may be expressed as $\Lap = \cond - \Trans$, where \cond denotes the associated multiplication operator. Note that $\Trans = \cond - \Lap$ is Hermitian on $\ell^2(\one)$ by \eqref{eqn:l2-Hermitian}. This is a bit of a surprise, since transfer operators are not generally Hermitian. Unfortunately, $\Trans_\unwtd$ may not be self-adjoint. In fact, the transfer operator of Example~\ref{exm:banded-nonselfadjoint-matrix} is not even essentially self-adjoint; see also \cite{vN32a,Rud91,DuSc88}.

\begin{lemma}\label{thm:Transc-bounded-when-cond-is}
  $\|\Trans_\unwtd\| \leq \|\cond\|$.
  \begin{proof}
    Recall that $\Trans_\unwtd = \Trans$ pointwise. The triangle inequality and Schwarz inequality give
    \linenopax
    \begin{align*}
      |\la f, \Trans f\ra_\unwtd|
      &\leq \sum_{x \in \verts} |f(x)| \sum_{y \nbr x} \left| \sqrt{\cond_{xy}} \sqrt{\cond_{xy}} f(y) \right| \\
      &\leq \sum_{x \in \verts} |f(x)| \cond(x)^{1/2} \left(\sum_{y \nbr x} \cond_{xy}|f(y)|^2 \right)^{1/2} \\
      &\leq \left(\sum_{x \in \verts} \cond(x) |f(x)|^2 \right)^{1/2}
        \left(\sum_{x,y \in \verts} \cond_{xy}|f(y)|^2 \right)^{1/2}.
    \end{align*}
    Since the both factors above may be bounded above by $\left(\uBd \cdot \|f\|_\unwtd^2\right)^{1/2}$ (using another application of Schwarz for the one on the right), we have $|\la f, \Trans f\ra_\cond| \leq \uBd \cdot \|f\|^2_\unwtd$.
  \end{proof}
\end{lemma}

\begin{remark}\label{rem:Transc-bounded-by-cond}
  When $d=1$, Example~\ref{exm:infinite-lattices} (the simple integer lattice) shows that the bound of Corollary~\ref{thm:Transc-bounded-when-cond-is} is sharp. From the proof of Lemma~\ref{thm:Fourier-transform-of-Lap-on-Zd}, one finds that
  \linenopax
  \begin{align*}
    \|\Trans\| = \sup|2\cos t| = 2 = 1+1 = \cond(n), \; \forall n \in \bZ.
  \end{align*}
\end{remark}

\begin{defn}\label{def:cond_x2}
  Let $\cond_x$ be defined by $\cond_x(y) = \cond_{xy}$, so 
  \begin{equation}\label{eqn:cond_x^2-notation}
    \cond_x \cdot \cond_x := \sum_{y \nbr x} \cond_{xy}^2
  \end{equation}
  We denote this with the shorthand $\cond_x^2 = \cond_x \cdot \cond_x$.
\end{defn}

\begin{theorem}\label{thm:transfer-operator-is-bounded-on-L2c}
  If \cond is bounded, then $\Trans_\unwtd:\ell^2(\unwtd) \to \ell^2(\unwtd)$ is bounded and self-adjoint. If $\Trans_\unwtd$ is bounded, then $\cond_x^2$ is a bounded function of $x$.
  \begin{proof}
    \fwd The boundedness of \Trans is Lemma~\ref{thm:Transc-bounded-when-cond-is}.
    Any bounded Hermitian operator is immediately self-adjoint; see Definition~\ref{def:self-adjoint}.

    \bwd For the converse, suppose that $\cond_x^2$ is unbounded. It follows that there is a sequence $\{x_n\}_{n=1}^\iy \ci \verts$ with $\cond_{x_n}^2 \to \iy$, and a path \cpath passing through each $x_n$ exactly once.  Consider the orthonormal sequence $\left\{\gd_{x_n}\right\}$:
    \linenopax
    \begin{align*}
      \Trans_\unwtd \gd_z(x)
      &= \sum_{y \nbr x} \cond_{xy} \gd_z(y)
      = \sum_{y \nbr z} \cond_{zy} \gd_z(y), \text{ and }
      \|\Trans_\unwtd \gd_z\|_\unwtd^2 = \sum_{y \nbr z} \cond_{yz}^2.
    \end{align*}
    Then letting $z$ run through the vertices of \cpath, it is clear that  $\|\Trans_\unwtd \gd_z\|_\unwtd^2 \to \iy$.
  \end{proof}
\end{theorem}

Recall from Definition~\ref{def:limit-at-infty} that $u(x)$ \emph{vanishes at \iy} iff for any exhaustion $\{\Graph_k\}$, one can always find $k$ such that $\|u(x)\|_\iy < \ge$ for all $x \notin \Graph_k$. 

Using a nested sequence as describe in Definition~\ref{def:limit-at-infty}, it is not difficult to prove that $\Trans_\unwtd$ is always the weak limit of the finite-rank operators $\Trans_n$ defined by $\Trans_n := P_n \Trans_\unwtd P_n$, where $P_n$ is projection to $\Graph_n = \spn\{\gd_x \suth x \in G_n\}$, so that
\begin{equation}\label{eqn:Trans-is-limit-of-Transn}
  \Trans_n v(x) = \charfn{G_n(x)} (\Trans_\unwtd v)(x) = \sum_{\substack{y \nbr x \\ y \in G_n}} \cond_{xy} v(y).
\end{equation}
Norm convergence does not hold without further hypotheses (see Example~\ref{exm:nonuniformly-converging-transfer-operator}) but we do have Theorem~\ref{thm:transfer-operator-is-compact-on-L2c}, which requires a lemma.

\begin{theorem}\label{thm:transfer-operator-is-compact-on-L2c}
  If $\cond \in \ell^2$ and $\deg(x)$ is bounded on \Graph, then the transfer operator $\Trans_\unwtd:\ell^2(\unwtd) \to \ell^2(\unwtd)$ is compact. If $\Trans_\unwtd$ is compact, then $\cond_x^2$ vanishes at \iy.
  \begin{proof}
    \bwd Consider any nested sequence $\{\Graph_k\}$ of finite connected subsets of \Graph, with $\Graph = \bigcup \Graph_k$, and the restriction of the transfer operator to these subgraphs, given by $\Trans_N := P_N \Trans_\unwtd P_N$, where $P_N$ is projection to $\Graph_N$. Then for $D_N := \Trans_\unwtd - \Trans_N$, consider the operator norm
    \linenopax
    \begin{align}\label{eqn:norm-estimate-on-trans-approximation}
      \|D_N\|
      = \left\|\begin{array}{cc}
          0 & P_N \Trans_\unwtd P_N^\perp \\
          P_N^\perp \Trans_\unwtd P_N & P_N^\perp \Trans_\unwtd P_N^\perp \\
        \end{array}\right\|,
    \end{align}
    where the ONB for the matrix coordinates is given by $\{\gd_{x_k}\}_{k=1}^\iy$ for some enumeration of the vertices.
    Since $\deg(x)$ is bounded, the matrices for $\Lap_\unwtd$ and hence also $\Trans_\unwtd$ are uniformly banded; whence $D_N$ is uniformly bounded with band size $b_N$ and Lemma~\ref{thm:bandedness-gives-norm-bound} applies.
    Since the first $N$ entries of $D_{N+b_N}v$ are 0, we have
    \linenopax
    \begin{align*}
      \|D_{N+b_N}v\|^2
      = \sum_{m=N+1}^\iy \left(\sum_{k=1}^{b_m} \cond_{mn_k} \right)^2
      = \sum_{m=N+1}^\iy \cond(x_m)^2,
    \end{align*}
    which tends to 0 for $\cond \in \ell^2(\unwtd)$.

    \fwd For the converse, suppose that $\cond_x^2$ does not vanish at \iy. It follows that there is a sequence $\{x_n\}_{n=1}^\iy \ci \verts$ with $\|\cond_{x_n}^2\|_\unwtd \geq \ge > 0$, and a path \cpath passing through each of them exactly once. By passing to a subsequence if necessary, is also possible to request that the sequence satisfies
    \linenopax
    \begin{align*}
      \vnbd(x_n) \cap \vnbd(x_{n+1}) = \es, \q\forall n,
    \end{align*}
    since the sequence need not contain every point of \cpath.
    Consider the orthonormal sequence $\left\{\gd_{x_n}\right\}$. We will show that $\left\{\Trans_\unwtd \gd_{x_n}\right\}$ contains no convergence subsequence:
    \linenopax
    \begin{align*}
      \Trans_\unwtd \gd_{x_n} - \Trans_\unwtd \gd_{x_m}
      &= \sum_{y \nbr x_n} \cond_{x_n y} \gd_{x_n} - \sum_{z \nbr x_m} \cond_{x_m z} \gd_{x_m} \\
      \|\Trans_\unwtd \gd_{x_n} - \Trans_\unwtd \gd_{x_m}\|^2
      &= \|\Trans_\unwtd \gd_{x_n}\|^2 + \|\Trans_\unwtd \gd_{x_m}\|^2
       = \sum_{y \nbr x_n} \cond_{x_n y}^2 + \sum_{z \nbr x_m} \cond_{x_m z}^2 \geq 2\ge.
    \end{align*}
    There are no cross terms in the final equality by orthogonality; $x_{n+1}$ was chosen to be far enough past $x_n$ that they have no neighbours in common.
  \end{proof}
\end{theorem}

\begin{cor}\label{thm:Transc-compact-when-cond-vanishes}
  If \cond vanishes at \iy and $\deg(x)$ is bounded, then $\Trans_\unwtd$ is compact.
  \begin{proof}
    The proof of the forward direction of Theorem~\ref{thm:transfer-operator-is-compact-on-L2c} just uses the hypotheses to show that $\sup_{x,y} \cond_{xy}$ can be made arbitrarily small by restricting $x,y$ to lie outside of a sufficiently large set.
  \end{proof}
\end{cor}

\subsection{Fredholm property of the transfer operator}
A stronger form of the Theorem~\ref{thm:Pot-nonempty-by-Fredholm} was already obtained in Corollary~\ref{thm:Pot-nonempty-via-Riesz}, but we include this brief proof for its radically contrasting flavour.

\begin{defn}\label{def:Fredholm-operator}
  A \emph{Fredholm operator} $L$ is one for which the kernel and cokernel are finite dimensional. In this case, the Fredholm index is $\dim \ker L - \dim \ker L^\ad$. Alternatively, $L$ is a Fredholm operator if and only if $\hat L$ is self-adjoint in the Calkin Algebra, i.e., $L = S + K$, where $S = S^\ad$ and $K$ is compact.
\end{defn}

\begin{theorem}\label{thm:Pot-nonempty-by-Fredholm}
  If \cond vanishes at infinity, then $\Pot(\ga,\gw)$ is nonempty.
  \begin{proof}
    When the Powers bound is satisfied, the previous results show \Lap is a bounded self-adjoint operator, and \Trans is compact. Consequently, \Lap is a Fredholm operator. By the Fredholm Alternative, $\ker \Lap = 0$ if and only if $\ran \Lap = \ell^2(\unwtd)$. Modulo the harmonic functions, $\ker \Lap = 0$, so $\gd_\ga - \gd_\gw$ has a preimage in $\ell^2(\unwtd)$.
  \end{proof}
\end{theorem}

\subsection{Some estimates relating \HE and $\ell^2(\unwtd)$}
\label{sec:relating-HE-and-L2}

In this section, we make the standing assumption that the functions under consideration lie in $\HE \cap \ell^2(\unwtd)$. Strictly speaking, elements of \HE are equivalence classes, but each has a unique representative in $\ell^2$ and it is understood that we always choose this one. Our primary tool will be the identity $\energy(u,v) = \la u, \Lap v \ra_\unwtd$ from \eqref{eqn:rem:energy-as-l2-inner-prod-with-Lap}, which is valid on the intersection $\HE \cap \ell^2(\unwtd)$. For example, note that this immediately gives
\begin{equation}\label{eqn:norm(lapl)-related-to-energy}
  \la v, \Lap v \ra_\energy = \|\Lap v\|_\unwtd^2,
  \qq\text{and}\qq
  \energy(v) = \|\Lap^{1/2} v\|_\unwtd^2,
\end{equation}
where the latter follows by the spectral theorem.
Theorem~\ref{thm:Pot-is-nonempty-by-current-flows} showed that $\Pot(\ga,\gw) \neq \es$, for any choice of $\ga \neq \gw$. It is natural to ask other questions in the same vein.
\begin{enumerate}[(i)]
  \item Is $\ell^2(\unwtd) \ci \HE$? No: consider the 1-\dimnl integer lattice described in Example~\ref{exm:1-dimnl-integer-lattice}.
  \item Is $\HE \ci \ell^2(\unwtd)$? No: consider the function $f$ defined on the binary tree in Example~\ref{exm:binary-tree:nonapproximable} which takes the value $1$ on half the tree and $-1$ on the other half (and is 0 at $o$). This function has energy $\energy(f)=2$, but it is easily seen that there is no $k$ for which $f+k \in \ell^2(\unwtd)$.
  \item Does $\Lap v \in \ell^2$ imply $v \in \HE$ or $v \in \ell^2(\unwtd)$? Neither of these are true, by the example in the previous item.
  \item Is $\Pot(\ga,\gw) \ci \ell^2$? No: consider again the 1-\dimnl integer lattice, with $\ga < \gw$. Then if $v \in \Pot(\ga,\gw)$, it will be constant (and equal to $v(\ga)$) for $x_n$ to the left of \ga, and it will be constant (and equal to $v(\gw)$) for $x_n$ right of \gw.
\end{enumerate}

\version{}{
However, we have some results in these directions under further restrictions:
\begin{enumerate}[(i)]
  \item When the Powers bound \eqref{eqn:def:power-bound} is satisfied, then $\ell^2(\unwtd) \ci \HE$. See Theorem~\ref{thm:energy-bounded-by-L2}.
  \item When the Powers bound \eqref{eqn:def:power-bound} is satisfied, then $v \in \HE \iff \Lap v \in \HE$. See Theorem~\ref{thm:hypoellipticity-in-HE}.
  \item If $\Lap v \in \ell^2(\unwtd)$, then $v \in \HE$. Moreover, if \Lap has a spectral gap, then the converse holds. See Theorem~\ref{thm:hypoellipticity-in-L2}.
\end{enumerate}
}

\begin{lemma}\label{thm:energy-bounded-by-L2}
  $\|v\|_\energy \leq \|\Lap^{1/2}\| \cdot \|v\|_\unwtd$ for every $v \in \HE$. If the Powers bound \eqref{eqn:def:power-bound} is satisfied, then  $\ell^2(\unwtd) \ci \HE$.
  \begin{proof}
    Since $\|v\|_\energy^2 = \la v, \Lap v\ra_\unwtd$, this is immediate from Lemma~\ref{thm:Lap-is-bounded-selfadjoint-on-L2c}.
  \end{proof}
\end{lemma}

\begin{lemma}\label{thm:Lap-bdd-onHE-by-L2}
  If \Lap is bounded on $\ell^2(\unwtd)$, then it is bounded with respect to \energy.
  \begin{proof}
    The hypothesis implies \Lap is self-adjoint on $\ell^2$, so that one can take the spectral representation $\hat \Lap$ on $L^2(X,d\gn)$ and perform the following computation:
    \linenopax
    \begin{align*}
      \|\Lap v\|_\energy
      &= \|\Lap \Lap^{1/2} v\|_\unwtd
      \leq \|\hat \Lap\|_\iy \cdot \|\Lap^{1/2} v\|_\unwtd
      = \|\hat \Lap\|_\iy \cdot \|v\|_\energy.
      &&\qedhere
    \end{align*}
  \end{proof}
\end{lemma}

\begin{lemma}\label{thm:energy-bounded-by-positivity}
  Let $v \in \ell^2(\unwtd)$. If $v \geq 0$ (or $v \leq 0$), then $\|v\|_\energy \leq \|v\|_\unwtd$. If $v$ is bipartite and alternating, then $\|v\|_\energy \geq \|v\|_\unwtd$.
  \begin{proof}
    Both statements follow immediately from the equality
    \linenopax
    \begin{align*}
      \energy(v)
      &= \la v, \Lap v\ra_\unwtd
      = \la v, v\ra_\unwtd - \la v, \Trans v\ra_\unwtd
      = \|v\|_\unwtd - \sum_{y \nbr x} \cond_{x,y} v(x) v(y).
      \qedhere
    \end{align*}
  \end{proof}
\end{lemma}


\section{The Laplacian and transfer operator on $\ell^2(\cond)$}
\label{sec:weighted-spaces}

In \S\ref{sec:self-adjointness-of-the-Laplacian}--\ref{sec:the-transfer-operator}, we studied \Lap and \Trans as operators on the unweighted space $\ell^2(\unwtd)$. In this section, we consider the renormalized versions of these operators and attempt to carry over as many results as possible to the context of $\ell^2(\cond)$.
\linenopax
\begin{equation}\label{eqn:defn:l2c-recall}
  \ell^2(\cond) := 
  \{u:\verts \to \bC \suth {\textstyle\sum\nolimits_{x \in \verts}} \cond(x) |u(x)|^2 < \iy\},
\end{equation}
with the inner product 
\linenopax
\begin{equation}\label{eqn:defn:l2c-inner-prod-recall}
  \la u, \Lap v\ra_\cond := \sum_{x \in \verts} \cond(x) u(x) \Lap v(x).
\end{equation}

Take the operator \Lap defined on $\spn\{\gd_x\}$, the dense domain consisting of functions with finite support. Then let \Lapc denote the closure of \Lap with respect to \eqref{eqn:defn:l2c-inner-prod-recall}, that is, its minimal self-adjoint extension to $\ell^2(\cond)$. 

\begin{lemma}\label{thm:sum(uLapv)<=|u|.|v|}
  For $u \in \ell^2(\cond)$ and $v \in \HE$,
  \linenopax
  \begin{align}\label{eqn:sum(uLapv)<=|u|.|v|}
    \sum_{x \in \verts} |u(x) \Lap v(x)| \leq \sqrt2 \; \|u\|_\cond \cdot \|v\|_\energy.
  \end{align}
  \begin{proof}
    Apply the Schwarz inequality twice, first with respect to the $x$ summation, then with respect to $y$: 
    \linenopax
    \begin{align*}
      \sum_{x \in \verts} |u(x) \Lap v(x)|
      &= \sum_{x,y \in \verts} \left|\sqrt{c_{xy}} u(x)\sqrt{c_{xy}} (v(x)-v(y))\right| \\
      &\leq \sum_{y \in \verts} \left(\sum_{x \nbr y} c_{xy} |u(x)|^2\right)^{1/2} 
        \left(\sum_{x \nbr y} c_{xy}|v(x)-v(y)|^2\right)^{1/2} \\
      &\leq \left(\sum_{x,y \in \verts} c_{xy} |u(x)|^2\right)^{1/2} 
        \left(\frac22\sum_{x,y \in \verts} c_{xy}|v(x)-v(y)|^2\right)^{1/2},
    \end{align*}
    and the resulting inequality retroactively justifies the implicit initial Fubination. 
  \end{proof}
\end{lemma}

\begin{defn}\label{def:J-quotient-inclusion-map}
  For $u \in \spn\{\gd_x\} \ci \ell^2(\cond)$, define a linear functional $\gx_u$ on $\HE$ by 
  \linenopax
  \begin{align}\label{eqn:def:gx_u}
    \gx_u(v) := \sum_{x \in \verts} u(x) \Lap v(x) 
  \end{align}
  Then $\gx_u$ is continuous because $|\gx_u(v)| \leq \|u\|_\cond \cdot \|v\|_\energy$ by Lemma~\ref{thm:sum(uLapv)<=|u|.|v|}, whence Riesz's lemma gives a $w \in \Fin$ for which $\gx_u(v) = \la w, v\ra_\energy$ holds for every $v \in \Fin$. 
  Let $J:\ell^2(\cond) \to \HE$ denote the map which sends $u \mapsto w$, i.e., the map defined by
  \linenopax
  \begin{align}\label{eqn:def:J}
    \la Ju, v \ra_\energy = \gx_u(v).
  \end{align}
\end{defn}

Definition~\ref{def:J-quotient-inclusion-map} allows one to see directly that $J\gd_x = [\gd_x]$:
\linenopax
\begin{align*}
  \la J\gd_x, v\ra_\energy
  &= \sum_{y \in \verts} \gd_x(y) \Lap v(y)
   = \Lap v(x)
   = \la \gd_x, v\ra_\energy,
   \q\text{for all } v \in \HE.
\end{align*}
This idea is the reason for Definition~\ref{def:J-quotient-inclusion-map} and also Theorem~\ref{thm:l2c-in-HE}. It is also easy to see that $\|\gd_x\|_\cond^2 = \cond(x) = \|\gd_x\|_\energy^2$, although the two norms $\|\cdot\|_\cond$ and $\|\cdot\|_\energy$ are clearly different in general. In fact, if $\gf = \sum_{x \in F} \gx_x \gd_x \in \spn\{\gd_x\}$ (so $F$ is finite), then one may easily compute
\begin{align*}
  \|\gf\|_\cond^2 = \sum_{x \in F} \cond(x) |\gx|^2,
  \q\text{whereas}\q
  \|\gf\|_\energy^2 = \|\gf\|_\cond^2 - \sum_{x \in F} \sum_{y \nbr x} \cond_{xy} \gx_x \gx_y.
\end{align*}

\begin{theorem}\label{thm:l2c-in-HE}
  The map $J$ is the quotient map induced by the equivalence relation $u\simeq v$ iff $u-v=const$, and gives a continuous embedding of $\ell^2(\cond)$ into \Fin with
  \linenopax
  \begin{align}\label{eqn:l2c-in-HE}
    \|Ju\|_\energy \leq  \sqrt2 \, \|u\|_\cond,
    \qq \forall u \in \ell^2(\cond).
  \end{align}
  Furthermore, the closure of $\ran J$ with respect to \energy is \Fin.
  \begin{proof}
    The formulation of $J$ in \eqref{eqn:def:J} gives 
    \linenopax
    \begin{align*}
      \la Ju, v_x \ra_\energy
       = \gx_u(v_x)
       = \sum_{x \in \verts} u(y) \Lap v_x(y)
       = \sum_{x \in \verts} u(y) (\gd_x - \gd_o)(y)
       = u(x) - u(o).
    \end{align*}
    This shows that $J$ is the quotient map as claimed. 
    The bound \eqref{eqn:l2c-in-HE} follows immediately upon combining \eqref{eqn:sum(uLapv)<=|u|.|v|} with \eqref{eqn:def:gx_u}.
    Now let $w = Ju$ for any $u \in \ell^2(\cond)$ and apply $\gx_u$ to $v = f_x = \Pfin v_x \in \Fin$ to get
    \linenopax
    \begin{align*}
      w(x) - w(o)
      &= \la f_x, w\ra_\energy
      = \gx_u(f_x)
      = \sum_{y \in \verts} u(y) (\gd_x-\gd_o)(y) 
      = u(x) - u(o).
    \end{align*}
    Since $\{f_x\}$ is thus a reproducing kernel for any element of $\ran J$, this shows that $\ran J \ci \Fin$, and hence
    \linenopax
    \begin{align*}
      \left| \la Ju, v\ra_\energy\right| = |\gx_u(v)| \leq \sqrt2 \, \|u\|_\cond \cdot \|v\|_\energy
      \q\implies\q
      \|Ju\|_\energy \leq \sqrt2 \, \|u\|_\cond.
    \end{align*}
    The \energy-closure of $\ran J$ is equal to \Fin because $\ran J$ contains $\spn\{\gd_x\}$.
  \end{proof}
\end{theorem}

\begin{lemma}\label{thm:adjJ-and-Prob}
  The adjoint map $J^\ad:\HE \to \ell^2(\cond)$ is given by  
  \begin{align}\label{eqn:adjJ-and-Prob}
    J^\ad u = u - \Prob u,
  \end{align}
  where \Prob is the probabilistic transition operator defined in \eqref{eqn:def:p(x,y)}.
  \begin{proof}
    First let $u \in \spn\{\gd_x\}$ and $v \in \HE$ and note that $\la Ju, v\ra_\energy = \la u, \Lap v\ra_\unwtd$ by \eqref{eqn:rem:energy-as-l2-inner-prod-with-Lap}. Then
    \linenopax
    \begin{align*}
     \la u, \Lap v\ra_\unwtd
      &= \sum_{x \in \verts} u(x) \Lap v(x) 
      = \sum_{x \in \verts} u(x) \cond(x) \left(v(x) - \frac1{\cond(x)} \sum_{y \nbr x} \cond_{xy} v(y)\right),
    \end{align*}
    whence $\la Ju, v\ra_\energy = \la u, (\one-\Prob)v\ra_\cond$ on the subspace $\spn\{\gd_x\}$, which is dense in \Fin in the norm $\|\cdot\|_\energy$ and dense in $\ell^2(\cond)$ in the norm $\|\cdot\|_\cond$.
  \end{proof}
\end{lemma}

\begin{remark}\label{rem:proof-of-{thm:l2c-in-HE}-via-Harm}
  Lemma~\ref{thm:adjJ-and-Prob} provides another proof of Theorem~\ref{thm:l2c-in-HE}:
  \begin{proof}[Alternative proof of Theorem~\ref{thm:l2c-in-HE}]
    Suppose that $v \in \ran(J)^\perp \ci \HE$ so that $\la Ju, v\ra_\energy = 0$ for all $u \in \ell^2(\cond)$. Then 
    \linenopax
    \begin{align*}
     \la Ju, v\ra_\energy
      &= \la u, J^\ad v\ra_\cond
       = \la u, v-\Prob v\ra_\cond,
       \q \forall u \in \ell^2(\cond)
    \end{align*}
    by \eqref{eqn:adjJ-and-Prob}, so that $v-\Prob v = 0$ in $\ell^2(\cond)$. Then recall that $v \in \Harm$ iff $\Lap v = 0$ iff $\Prob v = v$. This shows that $\ran(J)^\perp = \Harm$, and hence $\opclosure{\ran(J)}= \ran(J)^{\perp\perp} = \Harm^\perp = \Fin$.
  \end{proof}
\end{remark}

\begin{remark}\label{rem:defs-of-Lap}
  Many authors use $J^\ad J = \one-\Prob$ as the definition of the Laplace operator on $\ell^2(\cond)$. It is intriguing to note that for $v \in \HE$, one has $v-\Prob v \in \ell^2(\cond)$, even though it is quite possible that neither $v$ nor $\Prob v$ lies in $\ell^2(\cond)$ (for an extreme example, consider $v \in \Harm$).  
  Furthermore, note that for every $v \in \HE$, one has
  \linenopax
  \begin{align}\label{eqn:v-tends-to-harmonicity-at-iy}
    \sum_{x \in \verts} \cond(x) |v(x) - \Prob v(x)|^2
    = \|J^\ad v\|_\cond^2
    \leq \sqrt2 \, \|v\|_\energy^2
    < \iy.
  \end{align}
  This obviously implies a bound $\cond(x) |v(x) - \Prob v(x)|^2 \leq B^2$, whence
  \linenopax
  \begin{align*}
    |v(x) - \Prob v(x)| \leq B \cond(x)^{-1/2}.
  \end{align*}
  Consequently, if $\{G_k\}$ is any exhaustion of \Graph, then $v \approx \Prob v$ on $G_k^\complm$ for large $k$ (in the sense of $\ell^2(\cond)$). Roughly, one can say that  any $v \in \HE$ tends to being a harmonic function at \iy, and the faster \cond grows, the better the approximation.
\end{remark}

\begin{cor}\label{thm:defective-corollary}
  If $\Lap u = -u$, then $\sum \frac1{\cond(x)} |u(x)|^2 < \iy$ and $u(x) = O(\negsp[5]\sqrt{\cond(x)})$, as $x \to \iy$. 
  \begin{proof}
    Recall that a defect vector $u$ satisfies $\Lap u = -u \in \HE$ and hence $u-\Prob u = -\frac{1}{\cond} u$. The result follows by substituting the latter into \eqref{eqn:v-tends-to-harmonicity-at-iy}. This immediately implies a bound $\frac1{\cond(x)} |u(x)|^2 \leq B$, which gives the final claim.
  \end{proof}
\end{cor}

Recall 
that \LapV denotes the closure of the Laplacian when taken to have the dense domain $V := \spn\{v_x\}_{x \in \verts \less\{o\}}$ of \emph{finite} linear combinations of dipoles. 

\begin{cor}\label{thm:c-bdd-implies-Lap-selfadjoint-on-HE}
  If $\cond(x)$ is bounded on \verts and $\deg(x) < \iy$, then \LapV is essentially self-adjoint on \HE.
  \begin{proof}
    Suppose $w \in \dom \LapV^\ad$ satisfying $\Lap^\ad w = -w$. This means there is a $K$ (possibly depending on $w$) such that $|\la w, \Lap v\ra_\energy| \leq K \|v\|_\energy$ for all $v \in V$.
        
    Then for $u \in \ell^2(\cond)$ and $v \in \spn\{v_x\}$, set $\gx_u(v) = \sum_{x \in \verts} u(x) \Lap v(x)$. As in the proof of Theorem~\ref{thm:l2c-in-HE}, $\gx_u$ extends to a continuous linear functional on \HE, so applying it to $w$ gives
    \linenopax
    \begin{align*}
      \left|\sum_{x \in \verts} u(x) w(x)\right|
      = \left|\sum_{x \in \verts} u(x) (-w(x))\right|
      = \left|\sum_{x \in \verts} u(x) \Lap w(x)\right|
      = |\gx_u(w)| 
      \leq \sqrt2 \, \|u\|_\cond \cdot \|w\|_\energy.
    \end{align*}
    However, if $\cond(x)$ is bounded by \uBd, then 
    \linenopax
    \begin{align*}
      \|u\|_\cond 
      = \left(\sum_{x \in \verts} \cond(x) |u(x)|^2 \right)^{1/2}
      \leq \uBd^{1/2} \|u\|_\unwtd.
    \end{align*}
    Combining the two displayed equations above yields the inequality 
    \linenopax
    \begin{align*}
      \left|\sum u(x) w(x)\right| 
      \leq \sqrt2 \, \uBd^{1/2} \|u\|_\unwtd \cdot \|w\|_\energy.
    \end{align*}
    This shows $u \mapsto \sum_{x \in \verts} u(x) w(x)$ is a continuous linear functional on $\ell^2(\cond)$, so that Riesz's lemma puts $w \in \ell^2(\unwtd)$. However, now that $w$ is a defect vector in $\ell^2(\unwtd)$, Theorem~\ref{thm:Laplacian-is-essentially-self-adjoint} applies, and hence $w=0$.
  \end{proof}
\end{cor}

\begin{defn}
 Let $\Lapc := \cond^{-1} \Lap = \one-\Prob = JJ^\ad$ denote the probabilistic Laplace operator on \HE, as in \eqref{eqn:probablistically-reweighted-laplacian}. Note that we abuse notation here in the suppression of the quotient map, so that $\one-\Prob$ denotes an operator on \HE and a mapping $\HE \to \ell^2(\cond)$.
\end{defn}

\begin{cor}
  For any $v \in \HE$, \Lapc is contractive on \HE and $(\one-\Prob)v \in \ell^2(\cond)$ with
  \linenopax
  \begin{align*}
   \|(\one-\Prob)v\|_\energy \leq \|v\|_\energy. 
  \end{align*}
  \begin{proof}
    Since $J$ is contractive, it follows that $J^\ad$ is contractive by basic operator theory; this is a consequence of the polar decomposition applied to $J$.
    Then $\Lapc = JJ^\ad$ is certainly continuous with $\|(\one-\Prob)v\|_\energy = \|JJ^\ad v\|_\energy \leq \sqrt2 \, \|J^\ad v\|_\cond \leq 2 \|v\|_\energy$.
  \end{proof}
\end{cor}

\begin{lemma}\label{thm:Lap-not-Hermitian-on-l2(c)}
  \Lapc is Hermitian if and only if $\cond(x)$ is a constant function on the vertices.
  \begin{proof}
   This can be seen by computing the matrix representation of \Lapc with respect to the ONB $\{\frac{\gd_x}{\sqrt{\cond(x)}}\}$, in which case the \nth[(x,y)] entry is 
    \linenopax
    \begin{align*}
      [M_{\Lapc}]_{x,y} 
      = \left\la \frac{\gd_x}{\sqrt{\cond(x)}}, \Lap \frac{\gd_y}{\sqrt{\cond(y)}} \right\ra_\cond 
      &= \sum_{z \in \verts} \cond(z)\frac{\gd_x(z)}{\sqrt{\cond(x)}} \left(\frac{\cond(y)\gd_y(z)}{\sqrt{\cond(y)}} - \sum_{t \nbr y} \cond_{ty} \frac{\gd_t(z)}{\sqrt{\cond(y)}} \right) \\
      &= \sum_{z \in \verts} \left(\frac{\cond(z)\gd_x(z)}{\sqrt{\cond(x)}} \frac{\cond(y)\gd_y(z)}{\sqrt{\cond(y)}} - \frac{\cond(z)\gd_x(z)}{\sqrt{\cond(x)}}\sum_{t \nbr y} \cond_{ty} \frac{\gd_t(z)}{\sqrt{\cond(y)}} \right) \\
      &= \cond(x)\gd_{xy} - \sqrt{\cond(x)} \sum_{t \nbr y} \cond_{ty} \frac{\gd_t(x)}{\sqrt{\cond(y)}} \\
      &= \cond(x)\gd_{xy} - \sqrt{\cond(x)} \sum_{y \nbr x} \cond_{xy} \frac{\gd_y}{\sqrt{\cond(y)}},
    \end{align*}
    which is not symmetric in $x$ and $y$. We used $\Lap \gd_y = \cond(y)\gd_y - \sum_{t \nbr y} \cond_{ty} \gd_t$, which follows easily from Lemma~\ref{thm:dx-as-vx}.
  \end{proof}
\end{lemma}

\version{}{
\begin{lemma}
\marginpar{Removed until we can fix the proof}
  Let $u \in \HE$. Then $u \in \dom \LapV^\ad$ if and only if $\Lap u \in \ell^2(\cond)$.
  \begin{proof}
    By definition, $u \in \dom \LapV^\ad$ iff there is a $K$ (possibly depending on $u$) such that
    \linenopax
    \begin{align*}
      |\la u, \Lap v\ra_\energy|
      \leq K \|v\|_\energy, \qq \forall v \in \spn\{v_x\}.
    \end{align*}
    However, 
    \linenopax
    \begin{align*}
      |\la u, \Lap v\ra_\energy|
      =
    \end{align*}
  \end{proof}
\end{lemma}
}

\begin{theorem}\label{thm:Prob-self-adjoint-in-[-1,1]}
  (i) As an operator on $\ell^2(\cond)$, $\Prob = \id - J^\ad J$ is self-adjoint with $-\id \leq \Prob \leq \id$.
  
  (ii) As an operator on $\HE$, $\Prob = \id -J J^\ad$ is self-adjoint with $-\id \leq \Prob \leq \id$.
  \begin{proof}
    (i) Since $J$ is bounded and hence closed, a theorem of von Neumann implies $J^\ad J$ is self-adjoint. Then 
    \linenopax
    \begin{align*}
      \la u, \Prob u \ra_\cond
      &= \la u, (\id - J^\ad J)u \ra_\cond
       = \la u, u \ra_\cond - \la Ju, Ju \ra_\energy
       = \|u\|_\cond^2 - \|Ju\|_\energy^2.
    \end{align*}
    Since $\|Ju\|_\energy^2 \leq 2 \|u\|_\cond^2$ by \eqref{eqn:l2c-in-HE}, this establishes $-\|u\|_\cond^2 \leq \la u, \Prob u\ra_\cond \leq \|u\|_\cond^2$.

    (ii) By the same argument as in part (i), $JJ^\ad$ is self-adjoint. Then $\|JJ^\ad\| = \|J^\ad J\|$ gives the same bound for \Prob on \HE.
  \end{proof}
\end{theorem}

\begin{remark}\label{rem:Prob-self-adjoint-by-comp}
  One can also see that \Prob is self-adjoint by independent arguments. For $\ell^2(\cond)$, we have $\la u, \Prob v\ra_\cond = \sum_{x,y} u(x) \cond_{xy} v(y) = \la \Prob u, v\ra_\cond$, and for \HE, we have
  \linenopax
  \begin{align}\label{eqn:Prob-computation-1}
    \la \Prob v_x,v_y\ra_\energy 
    &= \la v_x - \cond^{-1}(\gd_x-\gd_o),v_y\ra_\energy \notag \\
    &= \la v_x,v_y\ra_\energy-\cond^{-1}(\gd_x-\gd_o)(y) + \cond^{-1}(\gd_x-\gd_o)(o) \notag \\
    &= \la v_x,v_y\ra_\energy - \tfrac1{\cond(o)} - \tfrac{\gd_{xy}}{\cond(x)},
  \end{align}
  where $\gd_{xy}$ is the Kronecker delta.
\end{remark}

\begin{remark}\label{rem:von-Neumann's-ergodic-theorem}
  von Neumann's ergodic theorem implies that
  \linenopax
  \begin{align*}
    \lim_{N \to \iy} \frac1N \sum_{n=1}^N \Prob^n u = \Phar u.
  \end{align*}
  In general, it is difficult to know when one has the stronger result that $\Prob^n u \to \Phar u$. A Perron-Frobenius-Ruelle theorem would require an invariant measure with certain properties not satisfied in the present context. Nonetheless, one can see that $\lim_{n \to \iy} \Prob^n u$ lies in \Harm; this is shown in Theorem~\ref{thm:Prob-is-strictly-contractive-on-Fin}.
\end{remark}

\begin{lemma}\label{thm:Prob-removes-Lap(cinv)}
  For all $\gf \in \spn\{v_x\}$, one has
  \linenopax
  \begin{align}\label{eqn:Prob-removes-Lap(cinv)}
    \la \Prob \gf, \gf \ra_\energy 
    = \|\gf\|_\energy^2 - \sum_{x \in \verts} \frac{|\Lap \gf(x)|^2}{\cond(x)}.
  \end{align}
  \begin{proof}
    For some finite set $F \ci \verts$ and $\gf = \sum_{x \in F} a_x v_x$, 
    \linenopax
    \begin{align*}
      \la \Prob \gf, \gf \ra_\energy 
      &= \sum_{x,y \in F} a_x a_y \la v_x, v_y\ra_\energy 
        - \frac1{\cond(o)} \left|\sum_{x \in \verts} a_x \right|^2
        -\sum_{x \in \verts} \frac{\left| a_x \right|^2}{\cond(x)} \\
      &= \|\gf\|_\energy^2 
        - \frac1{\cond(o)} \left|\sum_{x \in \verts} \Lap \gf(x) \right|^2
        - \sum_{x \in \verts} \frac{|\Lap \gf(x)|^2}{\cond(x)},
    \end{align*}
    where we have used \eqref{eqn:Lapu(y)-is-the-yth-coord} and \eqref{eqn:Prob-computation-1} for the last step. Then Theorem~\ref{thm:sum(Lap v)=0} shows that the middle sum vanishes, and we have \eqref{eqn:Prob-removes-Lap(cinv)}.
  \end{proof}
\end{lemma}

\begin{cor}\label{thm:Lapsumming-is-contractive-on-V}
  For $\gf \in \spn\{v_x\}$, $\sum_{x \in \verts} \frac{|\Lap \gf(x)|^2}{\cond(x)} 
    \leq 2 \|\gf\|_\energy^2$.
  \begin{proof}
    Apply \eqref{eqn:Prob-removes-Lap(cinv)} to the bound $-\id \leq \Prob \leq \id$ from Theorem~\ref{thm:Prob-self-adjoint-in-[-1,1]}.
  \end{proof}
\end{cor}

\begin{theorem}\label{thm:Prob-is-strictly-contractive-on-Fin}
  \Prob is strictly contractive on \Fin. 
  \begin{proof}
    Note that no harmonic functions lie in $\spn\{v_x\}$ (by Lemma~\ref{thm:E(u,v)=<u,Lapv>-on-spn{vx}}, for example), so $\sum_{x \in \verts} \frac{|\Lap \gf(x)|^2}{\cond(x)} > 0$.
  \end{proof}
\end{theorem}

For Theorem~\ref{thm:prob-bounded-on-ell2}, we will need to consider iterates $\Prob^n$ of the probabilistic transition operator, the induced conductances $\cond_{xy}^{(n)}$, and the corresponding energy spaces $\HE(\cond^{(n)})$. From
\linenopax
\begin{align*}
  (\Prob^2 u)(x) 
  = \Prob \left(\sum_{y \nbr x} p(x,y) u(y)\right)
  = \sum_{y \nbr x} p(x,y) \sum_{z \nbr y} p(y,z) u(z)
  = \sum_{y \nbr x} \sum_{z \nbr y} \frac{\cond_{xy}}{\cond(x)} \frac{\cond_{yz}}{\cond(y)}  u(z),
\end{align*}
we take $p^2(x,y) = \sum_{z \nbr x,y} p(x,z) p(z,y)$ and define $\cond^{(2)}_{xy} := \cond(x) p^2(x,y)$. Iterating, we obtain $\cond^{(n)}$ and the spaces $ \HE^{(n)} :=\{u:\verts \to \bC \suth \energy^{(n)}(u) < \iy\}$ where
\linenopax
\begin{align*}
  \energy^{(n)}(u) := \frac12 \sum_{x,y \in \verts} \cond^{(n)}_{xy} |u(x)-u(y)|^2.
\end{align*}

\begin{lemma}\label{thm:sums-invariant-under-iteration-of-P}
  For each $n=1,2,\dots$, one has $\cond(x) = \sum_{y \in \verts} \cond^{(n)}_{xy}$.
  \begin{proof}
    From above, $\sum_{y \in \verts} \cond^{(2)}_{xy} = \cond(x) \sum_{y \in \verts} p^{(2)}(x,y) = \cond(x)$ shows that the sum at $x$ does not change. The case for general $n$ follows by iterating.
  \end{proof}
\end{lemma}

In the proof of the following theorem, we write $\Lap_n := \cond(\id - \Prob^n)$, which is \emph{not} the same as $\Lap^n$. Also, we abuse notation and write $\ell^2(\cond)$ for $J(\ell^2(\cond))$, in accordance with Theorem~\ref{thm:l2c-in-HE}.

\begin{theorem}\label{thm:prob-bounded-on-ell2}
  \Prob is densely defined on \HE with 
  \linenopax
  \begin{align}\label{eqn:prob-bounded-on-ell2}
    \|\Prob u\|_\cond \leq \|u\|_\cond + \|u\|_\energy,
    \qq \forall u \in J(\ell^2(\cond)).
  \end{align}
  In fact, $\|\Prob^{(n)} u\|_\cond \leq \|u\|_\cond + \|u\|_{\energy^{(n)}}$, for every $n \geq 1$.
  \begin{proof}
    Note that it immediately follows from Lemma~\ref{thm:sums-invariant-under-iteration-of-P} that $\ell^2(\cond) \ci \bigcap_n \HE^{(n)}$, and hence that there is a $J_n:\ell^2(\cond) \to \HE^{(n)}$ for which
    \linenopax
    \begin{align*}
      \la J_n u, v\ra_{\energy^{(n)}}
      &= \sum_{x \in \verts} u(x) \Lap_n v(x)
       = \sum_{x \in \verts} u(x) \cond(x)\left(v(x) - \Prob^n v(x)\right)
       = \la u, v \ra_\cond - \la u, \Prob^n v \ra_\cond.
    \end{align*}
    Observe that the rearrangement in the last step is justified by the convergence of $\la J_n u, v\ra_{\energy^{(n)}}$ and $\la u, v \ra_\cond$. The Schwarz inequality gives $\Prob^n v \in \ell^2(\cond)$, and 
    \linenopax
    \begin{align*}
      \left|\la u, \Prob^n v \ra_\cond\right|
      \leq \left|\la u, v \ra_\cond\right| + \left|\la J_n u, v\ra_{\energy^{(n)}}\right|
      \leq \|u\|_\cond \|v\|_\cond + \|u\|_\cond \|v\|_{\energy^{(n)}}.
    \end{align*}
    This shows that $\Prob^nv$ is in the dual of $\ell^2(\cond)$ and hence in $\ell^2(\cond)$. By Riesz's theorem, $\|\Prob u\|_\cond$ is the best constant possible in the above inequality, and so \eqref{eqn:prob-bounded-on-ell2} follows.
  \end{proof}
\end{theorem}

\section{Remarks and references}
\label{sec:Remarks-and-References-ell2-of-Lap-and-Trans}

The material above is an assortment of results on spectral theory for the operators from Chapter~\ref{sec:Lap-on-HE}. 

Of the results in the literature of relevance to the present chapter, the references \cite{Dod06, DuRo08, Sto08, DMS07, BaBa05, Chu07, ChRi06, BiLeSt07} are especially relevant. The reader may also wish to consult \cite{Nel73a}, \cite{vN32a}, and good background references include \cite{DuSc88, ReedSimonII, Arveson:spectral-theory, Chu01, Lax-Phillips}.

\version{}{
\section{The spectral Laplacian and transfer operator}
\label{sec:spectral-spaces}

One disadvantage of \Lapc is that it is not self-adjoint. In this section, we consider a renormalized version of \Lap which is defined to as to remain self-adjoint.

\subsection{The spectral transfer operator on $\ell^2(\cond)$}
\label{sec:The-spectral-transfer-operator-on-L2}
In this section, we consider the \emph{spectral Laplacian} $\Lapsp = \cond^{-1/2} \Lap \cond^{-1/2}$, given pointwise by
\linenopax
\begin{align}
  (\Lapsp v)(x) 
  := \frac1{\sqrt{\cond(x)}} \sum_{y \nbr x} 
  \cond_{xy} \left(\frac{v(x)}{\sqrt{\cond(x)}} - \frac{v(y)}{\sqrt{\cond(y)}}\right) ,
    \label{eqn:spectral-laplace-operator}
\end{align}
and the associated spectral transfer operator $\Transp = \id - \Lapsp = \cond^{-1/2} \Trans \cond^{-1/2}$. Recall that with respect to the ONB $\left\{\frac{\gd_x}{\sqrt{\cond(x)}}\right\}$, the operator \cond is the diagonal matrix with entries $\cond(x)$. 

It is not usually true that transfer operators are self-adjoint, so the following result is a bit of a surprise. 

\begin{theorem}\label{thm:spectral-transfer-operator-is-Hermitian-on-L2}
  The transfer operator is Hermitian, as an operator on on $\ell^2(\cond)$.
  \begin{proof}
    Using Theorem~\ref{thm:transfer-operator-is-Hermitian-on-L2c} for the third equality,
    \linenopax
    \begin{align*}
      \la u, \Transp v\ra_\cond
      = \la u, \cond^{-1/2} \Trans \cond^{-1/2} v \ra_\cond
      &= \la \cond^{-1/2} u, \Trans \cond^{-1/2} v \ra_\unwtd \\
      &= \la \Trans \cond^{-1/2} u, \cond^{-1/2} v \ra_\unwtd \\
      &= \la \Transp u, v \ra_\spectral.
      \qedhere
    \end{align*}
  \end{proof}
\end{theorem}

Unfortunately, $\Trans_\spectral$ may not be essentially self-adjoint, as shown by Example~\ref{exm:banded-nonselfadjoint-matrix}; see also \cite{vN32a,Rud91,DuSc88}. In view of the fact that
\linenopax
\begin{align*}
  \energy(u) = \la u, u\ra_{\energy} = \la u, \Lap u\ra_\spectral \geq 0,
\end{align*}
self-adjointness follows automatically.

\begin{theorem}\label{thm:transfer-operator-is-bounded-on-L2}
  The transfer operator $\Trans_\spectral:\ell^2(\spectral) \to \ell^2(\spectral)$ is bounded if and only if $\cond(x)$ is bounded, i.e., iff the Powers bound \eqref{eqn:def:power-bound} is satisfied.
  \begin{proof}
    Apply conjugation.
  \end{proof}
\end{theorem}

\begin{theorem}\label{thm:transfer-operator-is-compact-on-L2}
  The transfer operator $\Trans_\spectral:\ell^2(\spectral) \to \ell^2(\spectral)$ is compact if and only if $\cond(x) \limas{x} 0$ in $\ell^2$, i.e., iff $\|\cond(x)\|_2$ vanishes at \iy.
  \begin{proof}
    \bwd Consider any nested sequence $\{\Graph_k\}$ of finite connected subsets of \Graph, with $\Graph = \bigcup \Graph_k$, and the restriction of the transfer operator to these subgraphs, given by
    \linenopax
    \begin{align*}
      \Trans_n := P_n \Trans P_n,
    \end{align*}
    where $P_n$ is projection to $\Graph_n$. Then for $D_n := \Trans - \Trans_n$, consider
    \linenopax
    \begin{align}\label{eqn:norm-estimate-on-trans-approximation-L2}
      \|D_n\|
      = \left\|\begin{array}{cc}
          0 & P_n \Trans P_n^\perp \\
          P_n^\perp \Trans P_n & P_n^\perp \Trans P_n^\perp \\
        \end{array}\right\|,
    \end{align}
    where the basis for the matrix coordinates is given by $\{\gd_{x_k} \cond(x_k)^{-1/2}\}_{k=1}^\iy$ for some enumeration of the vertices.
    We consider the estimate $\|D_n\| \leq \max \cond_{xy}$, ranging over the matrix entries of the three submatrices indicated in \eqref{eqn:norm-estimate-on-trans-approximation}. Since \Lap and hence \Trans is a banded matrix, it follows that $D_n$ is, too. Denote the width of the band of the row of $D_n$ corresponding to $x_k$ by $B_{x_k}$.

    Observe that any entry $\cond_{xy}$ appearing in submatrix of \eqref{eqn:norm-estimate-on-trans-approximation} corresponding to $P_n \Trans P_n^\perp$ (i.e., upper right) must satisfy the restrictions
    \linenopax
    \begin{align*}
      \cond_{xy} \in [P_n \Trans P_n^\perp]
      \q\implies\q y > n \text{ and } n \geq x > n-B_x,
    \end{align*}
    where $B_x$ is the bandwidth mentioned just above. Since $\|\cond(x)\|_2 = \sum_y \cond_{x,y}^2 \to 0$ as $y \to \iy$, one also has $\max_y \cond_{xy} \to 0$ as $y \to \iy$. This means that taking $n$ sufficiently large ensures all entries in $[P_n \Trans P_n^\perp]$ are less than a specified $\ge>0$.

    The other two submatrices have the restrictions
    \linenopax
    \begin{align*}
      \cond_{xy} \in [P_n^\perp \Trans P_n] &\q\implies\q y \leq n \text{ and } n < x \leq n+B_x, \text { and} \\
      \cond_{xy} \in [P_n^\perp \Trans P_n^\perp] &\q\implies\q y > n \text{ and } x > n,
    \end{align*}
    but the rest of the argument is the same.

    \fwd For the converse, suppose that \cond does not vanish at \iy. It follows that there is a sequence $\{x_n\}_{n=1}^\iy \ci \verts$ with $\cond(x_n) \geq \ge > 0$, and a path \cpath passing through each of them exactly once. By passing to a subsequence if necessary, is also possible to request that the sequence satisfies
    \linenopax
    \begin{align*}
      \vnbd(x_n) \cap \vnbd(x_{n+1}) = \es, \q\forall n,
    \end{align*}
    since the sequence need not contain every point of \cpath.
    Consider the orthonormal sequence $\left\{\gd_{x_n}\right\}$. We will show that $\left\{\Trans\gd_{x_n}\right\}$ contains no convergence subsequence:
    \linenopax
    \begin{align*}
      \Trans \gd_{x_n} - \Trans \gd_{x_m} &= \sum_{k \nbr n} \gx_{nk} \gd_{x_n} - \sum_{l \nbr m} \gx_{lm} \gd_{x_m} \\
      \|\Trans \gd_{x_n} - \Trans \gd_{x_m}\|^2 &= \sum_{k \nbr n} |\gx_{nk}|^2 + \sum_{l \nbr m} |\gx_{lm}|^2 \geq 2\ge^2.
    \end{align*}
    There are no cross terms in the final equality by orthogonality; $x_{n+1}$ was chosen to be far enough past $x_n$ that they have no neighbours in common.
  \end{proof}
\end{theorem}
}

%% file: he-and-hd.tex

\chapter{The dissipation space \HD and its relation to \HE}
\label{sec:H_energy-and-H_diss}

\headerquote{Most of the fundamental ideas of science are essentially simple, and may, as a rule, be expressed in a language comprehensible to everyone.}{---~A.~Einstein}

While the vectors in \HE represent voltage differences, there is a second Hilbert space \HD which serves to complete our understanding of metric geometry for resistance networks. 
This section is about an isometric embedding \drp of the Hilbert space \HE of voltage functions into the Hilbert space \HD of current functions, and the projection \Pdrp that relates \drp to its adjoint. This dissipation space \HD will be needed for several purposes, including the resolution of the compatibility problem discussed in \S\ref{sec:compatibility-problem} and the solution of the Dirichlet problem in the energy space via its solution in the (ostensibly simpler) dissipation space; see Figure~\ref{fig:HE->HD->HE-decomposition}. The geometry of the embedding \drp is a key feature of our solution in Theorem~\ref{thm:HD=Nbd+Kir+Cyc} to a structure problem regarding current functions on graphs. Appendix~\ref{sec:partial-isometries} contains definitions of the terms isometry, coisometry, projection, initial projection, final projection, and other notions used in this section. 

In this section, we will find it helpful to use the notation $\ohm(x,y) = \cond_{xy}^{-1}$.

\begin{defn}\label{def:H_diss}
  Considering \ohm as a measure on \edges, currents comprise the Hilbert space
  \begin{equation}\label{eqn:def:H_diss}
    \HD := \{\curr:\edges \to \bC \suth \curr \text{ is antisymmetric and } \|\curr\|_\diss < \iy\},
  \end{equation}
  \glossary{name={\HD},description={dissipation Hilbert space},sort=h,format=textbf}
  where the norm and inner product are given by
  \linenopax
  \begin{align}\label{eqn:def:H_diss-form}
    \|\curr\|_\diss &:= \diss(\curr)^{1/2}
    \q\text{and}\q \\
    \la \curr_1, \curr_2\ra_\diss &:= \diss(\curr_1, \curr_2)
      = \frac12 \sum_{(x,y) \in \edges} \ohm(x,y) \cj{\curr_1(x,y)} \curr_2(x,y). \notag
  \end{align}
\end{defn}


Observe that $\HD = \ell^2(\edges,\ohm)$ but it is not true that \HE can be represented as an $\ell^2$ space in such an easy manner (but see Theorem~\ref{thm:HE-isom-to-L2(S',P)}).
As an $\ell^2$ space, $\HD$ is obviously complete. However, \HD is also blind to the topology of the underlying network and this is the reason why the space of currents is much larger than the space of potentials. This last statement is made precise in Theorem~\ref{thm:drp-and-drp-are-partial-isometries}, where it is shown that \HD is larger than \HE by precisely the space of currents supported on cycles.

The fundamental relationship between \HE and \HD is given by the following operator which implements Ohm's law. It can also be considered as a boundary operator in the sense of homology. Further motivation for the choice of symbology is explained in \S\ref{sec:analogy-with-complex-functions}.

\begin{defn}\label{def:drop-operator}
  The drop operator $\drp=\drp[\cond] : \HE \to \HD$ is defined by
  \begin{equation}\label{eqn:def:drop-operator}
    (\drp v)(x,y) := \cond_{xy}(v(x) - v(y))
  \end{equation}
  \glossary{name={\drp},description={drop operator},sort=d,format=textbf}
  and converts potential functions into currents (that is, weighted voltage drops) by implementing Ohm's law. In particular, for $v \in \Pot(\ga,\gw)$, we get $\drp v \in \Flo(\ga,\gw)$.
\end{defn}

As Lyons comments in \cite[\S9.3]{Lyons:ProbOnTrees}, thinking of the resistance $\ohm(x,y)$ as the length of the edge $(xy)$, \drp is a discrete version of directional derivative:
  \begin{equation*}
    (\drp v)(x,y) = \frac{v(x) - v(y)}{\ohm(x,y)} \approx \frac{\del v}{\del y}(x).
  \end{equation*}

\begin{lemma}\label{thm:drp-is-an-isometry}
  \drp is an isometry.
  \begin{proof}
     Lemma~\ref{thm:energy=dissipation} may be restated as follows: $\|\drp v\|_\diss^2 = \|u\|_\energy^2$.
  \end{proof}
\end{lemma}

\section{The structure of \HD}
In Theorem~\ref{thm:HD=Nbd+Kir+Cyc}, we are now able to characterize \HD by using Lemma~\ref{thm:drp-is-an-isometry} to extend Lemma~\ref{thm:HE=Fin+Harm} (the decomposition $\HE = \Fin \oplus \Harm$). First, however, we need some terminology. Whenever we consider the closed span of a set of vectors $S$ in \HE or $\HD$, we continue to use the notation $\clspn[\energy]{S}$ or $\clspn[\diss]{S}$ to denote the closure of the span in \energy or \diss, respectively.

\begin{defn}\label{def:reweighted-edge-nbd}
  Define the \emph{weighted edge neighbourhoods}
  \linenopax
  \begin{align}\label{eqn:def:reweighted-edge-nbd}
    \gh_z &= \gh_z(x,y) := \drp\gd_z =
    \begin{cases}
      \cond_{xy}, &x=z \nbr y, \\
      0, &\text{else}.
    \end{cases}
  \end{align}
  Then denote the space of all such currents by $\Nbd := \clspn[\diss]{\drp \Fin} = \clspn[\diss]{\gh_z}$. This space is called $\bigstar$ in \cite[\S2 and \S9]{Lyons:ProbOnTrees}.
\end{defn}

\begin{lemma}\label{thm:diss-of-edge}
  $\diss(\gh_z) = \deg(z)$.
  \begin{proof}
    Computing directly,
    \linenopax
    \begin{align*}
      \diss(\gh_z)
      &= \sum_{(x,y) \in \edges} \ohm(x,y) \gh_z(x,y)^2
      = \sum_{y \nbr x} 1
      = \deg(z).
      \qedhere
    \end{align*}
  \end{proof}
\end{lemma}

\begin{defn}\label{def:Kirchhoff-currents}\label{def:harmonic-currents}
  For each $h \in \Harm$, we have $\act (\drp h) = 0$ so that $\drp h$ satisfies the homogeneous Kirchhoff law by Corollary~\ref{thm:Kirchhoff-iff-Harm}. Therefore, we denote $\Kir := \drp \Harm = \ker \act$. Since the elements of \Kir are currents induced by harmonic functions, we call them \emph{harmonic currents} or \emph{Kirchhoff currents}.
\end{defn}

\begin{defn}\label{def:Cycle}
  Denote the space of cycles in \HD, that is, the closed span of the characteristic functions of cycles $\cycle \in \Loops$ by $\Cyc := \clspn[\diss]{\charfn{\cycle}}.$
\end{defn}
  
The space \Cyc is called {\large $\Diamond$} in \cite[\S2 and \S9]{Lyons:ProbOnTrees}.

We are now able to describe the structure of \HD. See \cite[(9.6)]{Lyons:ProbOnTrees} for a different proof.

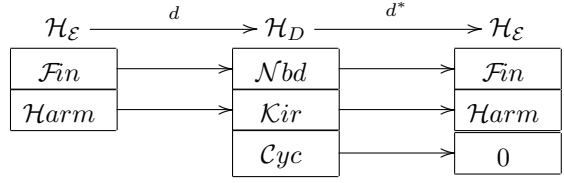
\begin{figure}
  \[\xymatrix @R=0pc @!C=5pc {
      \HE \ar^{\drp}[r] & \HD \ar^{\drpa}[r]& \HE \\
      \xyboxentry{\Fin} \ar[r]
        & \xyboxentry{\Nbd} \ar[r]
        & \xyboxentry{\Fin} \\
      \xyboxentry{\Harm} \ar[r]
        & \xyboxentry{\Kir} \ar[r]
        & \xyboxentry{\Harm} \\
      & \xyboxentry{\Cyc} \ar[r]
        & \xyboxentry{0}
    }\]
  \caption{\captionsize The action of \drp and \drpa on the orthogonal components of \HE and \HD. See Theorem~\ref{thm:HD=Nbd+Kir+Cyc} and Definition~\ref{def:projection-notation}.}
  \label{fig:HE->HD->HE-decomposition}
\end{figure}

\begin{theorem}\label{thm:HD=Nbd+Kir+Cyc}
  $\HD = \Nbd \oplus \Kir \oplus \Cyc$.
  \begin{proof}
    Lemma~\ref{thm:induced-currents-satisfy-cycle-condition} expresses the fact that $\drp \HE$ is orthogonal to \Cyc.
    Since \drp is an isometry, $\drp \HE = \drp \Fin \oplus \drp \Harm$,
    and the result follows from Theorem~\ref{thm:HE=Fin+Harm} and the definitions just above. See Figure~\ref{fig:HE->HD->HE-decomposition}.
  \end{proof}
\end{theorem}

\begin{defn}\label{def:projection-notation}
  Recall that a projection on a Hilbert space is by definition an operator satisfying $P=P^\ad=P^2$. The following notation will be used for projection operators:
  \linenopax
  \begin{align*}
    \Pfin:&\HE \to \Fin & \Phar:&\HE \to \Harm \\
    \Pdrp:&\HD \to \drp \HE & \Pdrp^\perp = \Pcyc:&\HD \to \Cyc \\
    \Pnbd:&\HD \to \Nbd & \Pkir:& \HD \to \Kir.
  \end{align*}
  Figure~\ref{fig:HE->HD->HE-decomposition} may assist the reader with seeing how these operators relate.
  \glossary{name={\Pfin},description={projection to \Fin,minimizing projection},sort=P,format=textbf}
  \glossary{name={\Phar},description={projection to \Harm},sort=P,format=textbf}
  \glossary{name={\Pdrp},description={projection to \drp\HE},sort=P,format=textbf}
  \glossary{name={\Pcyc},description={projection to \Cyc},sort=P,format=textbf}
  \glossary{name={\Pnbd},description={projection to \Nbd,minimizing projection},sort=P,format=textbf}
  \glossary{name={\Pkir},description={projection to \Kir},sort=P,format=textbf}
\end{defn}

\begin{lemma}\label{thm:drop-oper-adjoint-formula}
  The adjoint of the drop operator $\drpa : \HD \to \HE$ is given by
  \linenopax
  \begin{equation}\label{eqn:lem:drop-oper-adjoint-formula}
    (\drpa \curr)(x) - (\drpa \curr)(y) = \ohm(x,y) \Pdrp \curr(x,y).
  \end{equation}
  \begin{proof}
    Since $\Pdrp \drp = \drp$ and $\Pdrp = \Pdrp^\ad$ by definition,
    \linenopax
    \begin{align*}
      \la \drp v, \curr\ra_\diss
      &= \la \Pdrp \drp v, \curr\ra_\diss
      = \la \drp v, \Pdrp \curr\ra_\diss \\
      &= \frac12 \sum_{(x,y) \in \edges} \ohm(x,y) \cond_{xy} \cj{(v(x) - v(y))} \Pdrp \curr(x,y)
        &&\text{by } \eqref{eqn:def:drop-operator} \\
      &= \frac12 \sum_{x,y \in \verts} \cond_{xy}\cj{(v(x) - v(y))} (\drpa \curr)(x) - (\drpa \curr)(y))
         &&\text{by } \eqref{eqn:lem:drop-oper-adjoint-formula}.
         \qedhere
    \end{align*}
  \end{proof}
\end{lemma}

\begin{remark}\label{rem:adjoint-is-defined-only-up-to-constants}
  Observe that \eqref{eqn:lem:drop-oper-adjoint-formula} only defines the function $\drpa \curr$ up to the addition of a constant, but elements of \HE are equivalence classes, so this is sufficient. Also,
  \linenopax
  \begin{equation*}
    (\drpa \curr)(x) - (\drpa \curr)(y) = \ohm(x,y) \curr(x,y).
  \end{equation*}
  satisfies the same calculation as in the proof of Lemma~\ref{thm:drop-oper-adjoint-formula}. However, the compatibility problem described in \S\ref{sec:compatibility-problem} prevents this from being a well-defined operator on all of $\HD$.
  
  One can think of \drpa as a weighted boundary operator and \drp as the corresponding coboundary operator; this approach is carried out extensively in \cite{Soardi94}, although the author does not include the weight as part of his definition.
\end{remark}

\begin{theorem}\label{thm:drp-and-drp-are-partial-isometries}
  \drp and \drpa are partial isometries with initial and final projections
  \begin{equation}\label{eqn:drp-and-drp-are-partial-isometries}
    \drpa\drp = \id_{\HE}, \qq
    \drp\drpa = \Pdrp.
  \end{equation}
  Furthermore, $\drp:\HE \to \Nbd \oplus \Kir$ is unitary.
  \begin{proof}
    Lemma~\ref{thm:drp-is-an-isometry} states that \drp is an isometry; the first identity of \eqref{eqn:drp-and-drp-are-partial-isometries} follows immediately.
    The second identity of \eqref{eqn:drp-and-drp-are-partial-isometries} follows from the computation
    \linenopax
    \begin{align*}
      \drp\drpa \curr(x,y)
      &= \drp\left(\drpa\curr(x) - \drpa\curr(y)\right)
        &&\text{} \\
      &= \drp\left(\ohm(x,y) \Pdrp \curr(x,y)\right)
        &&\text{by \eqref{eqn:lem:drop-oper-adjoint-formula}} \\
      &= \cond_{xy} \left(\ohm(x,y) \Pdrp \curr(x,y) - \ohm(y,y) \Pdrp \curr(y,y)\right)
        &&\text{by \eqref{eqn:def:drop-operator}} \\
      &= \Pdrp \curr(x,y).
        &&\text{Definition~\ref{def:ERN}.}
    \end{align*}
    The last claim is also immediate from the previous computation.
  \end{proof}
\end{theorem}

We are now able to give an proof of the completeness of \HE which is independent of \S\ref{sec:vonNeumann's-embedding-thm}; see also Remark~\ref{rem:norm-vs-quasinorm}. 

\begin{lemma}\label{thm:HE=domE-is-complete}
  $\dom\energy/\{constants\}$ is complete in the energy norm.
  \begin{proof}
    Let $\{v_j\}$ be a Cauchy sequence. Then $\{\drp v_j\}$ is Cauchy in \HD by Theorem~\ref{thm:drp-and-drp-are-partial-isometries}, so it converges to some $\curr \in \HD$ (completeness of \HD is just the Riesz-Fischer Theorem). We now show that $v_j \to \drpa \curr \in \HE$:
    \linenopax
    \begin{align*}
      \energy(v_j - \drpa\curr)
      = \energy(\drpa(\drp v_j - \curr))
      \leq \diss(\drp v_j - \curr)
      \to 0,
    \end{align*}
    again by Theorem~\ref{thm:drp-and-drp-are-partial-isometries}.
  \end{proof}
\end{lemma}

\subsection{An orthonormal basis (ONB) for \HD}
\label{sec:An-ONB-for-HD}
Recall from Remark~\ref{def:orientation} that we may always choose an orientation on \edges. We use the notation $\edge = (x,y) \in \edges$ to indicate that \edge is in the orientation, and $\redge = (y,x)$ is not. For example, there is a term in the sum
\linenopax
\begin{align}\label{eqn:diss-on-oriented-edges}
  \diss(\curr) = \sum_{\edge \in \edges} \ohm(\edge) \curr(\edge)^2
\end{align}
  \glossary{name={\edge},description={oriented edge},sort=e,format=textbf}
for \edge, but there is no term for \redge (and hence no leading coefficient of $\frac12$).

\begin{defn}\label{def:ONB-for-HD}
  For $\edge = (x,y) \in \edges$, denote by $\emass$ the normalized Dirac mass on an edge:
  \begin{equation}\label{eqn:ONB-for-HD}
    \emass := \sqrt{\cond} \, \gd_{\edge}.
  \end{equation}
\end{defn}
  \glossary{name={\emass},description={normalized Dirac mass on an edge},sort=F,format=textbf}

\begin{lemma}\label{thm:ONB-for-HD}
  The weighted edge masses $\{\emass\}$ form an ONB for \HD.
  \begin{proof}
    It is immediate that every function in $\Fin(\edges)$ can be written as a (finite) linear combination of such functions. Since \HD is just a weighted $\ell^2$ space (as noted in Definition~\ref{def:H_diss}), it is clear that $\Fin(\edges)$ is dense in \HD. To check orthonormality,
    \linenopax
    \begin{align*}
      \la \emass[\edge_1], \emass[\edge_2]\ra_\diss
      = \sum_{\edge \in \edges} \ohm(\edge) \cj{\emass[\edge_1](\edge)} \emass[\edge_2](\edge)
      = \gd_{\edge_1, \edge_2},
    \end{align*}
    where $\gd_{\edge_1, \edge_2}$ is the Kronecker delta, since $\emass[\edge_1](\edge) \emass[\edge_2](\edge) = \cond_{\edge}$ iff $\edge_1=\edge_2$, and zero otherwise. (There is no term in the sum for $\redge_1$; see \S\ref{sec:An-ONB-for-HD}.) Incidentally, the same calculation verifies $\emass \in \HD$.
  \end{proof}
\end{lemma}

\begin{remark}\label{rem:edge-mass-not-in-ran(drp)}
  It would be nice if $\emass \in \drp\HE$, as this would allow us to ``pull back'' \emass to obtain a localized generating set for \HE, i.e., a collection of functions with finite support. Unfortunately, this is not the case whenever \edge is contained in a cycle, and the easiest explanation is probabilistic. If $x \nbr y$, then the Dirac mass on the edge $(x,y)$ corresponds to the experiment of passing one amp from $x$ to its neighbour $y$. However, there is always some positive probability that current will flow from $x$ to $y$ around the other part of the cycle and hence the minimal current will not be \emass; see Lemma~\ref{thm:shortest-path-bounds-resistance-distance} for a more precise statement.

  Of course, we can apply \drpa to obtain a nice result as in Lemma~\ref{thm:frame-for-HE}, however, \drpa\emass will generally not have finite support and may be difficult to compute. Nonetheless, it still has a very nice property; cf.~Lemma~\ref{thm:frame-for-HE}. In light of Theorem~\ref{thm:drp-and-drp-are-partial-isometries} and the previous paragraph, it is clear that any element $\drpa\emass$ is an element of $\Pot(x,y)$ for some $x \in \verts$ and some $y \nbr x$.
\end{remark}

\begin{lemma}\label{thm:frame-for-HE}
  \version{}{\marginpar{Is $\{\drpa\emass\}$ the same as $\{\Pot(x,y)\}_{y \nbr x}$?}}
  The collection $\{\drpa\emass\}$ is a Parseval frame for \HE.
  \begin{proof}
    The image of an ONB under a partial isometry is always a frame. That we have a Parseval frame (i.e., a tight frame with bounds $A=B=1$) follows from the fact that \drp is an isometry:
    \linenopax
    \begin{align*}
      \sum_{\edge \in \edges} |\la\drpa \emass, v \ra_\energy|^2
      = \sum_{\edge \in \edges} |\la\emass, \drp v \ra_\diss|^2
      = \|\drp v\|_\diss^2
      = \|v\|_\energy^2.
    \end{align*}
    We used Lemma~\ref{thm:ONB-for-HD} for the second equality and Theorem~\ref{thm:drp-and-drp-are-partial-isometries} for the third.
  \end{proof}
\end{lemma}

\section{The divergence operator}
\label{sec:divergence-operator}

In \S\ref{sec:Solving-potential-theoretic-problems}, we will see how \Pdrp allows one to solve certain potential-theoretic problems, but first we need an operator which enables us to study \Lap with respect to \HD rather than \HE. While the term ``divergence'' is standard in mathematic, the physics literature sometimes uses ``activity'' to connote the same idea, e.g., \cite{Pow75}--\cite{Pow79}. We like the term ``divergence'' as it corresponds to the intuition that the elements of \HD are (discrete) vector fields.

\begin{defn}\label{def:divergence-oper}
  The \emph{divergence operator} is $\act : \HD \to \HE$ given by
  \begin{equation}\label{eqn:def:divergence-oper}
    \act(\curr)(x) := \sum_{y \nbr x} \curr(x,y).
  \end{equation}
  To see that \act is densely defined, note that $\act(\gd_{(x,y)})(z) = \gd_{x}(z) - \gd_{y}(z)$, and the space of finitely supported edge functions $\Fin(\edges)$ is dense in $\ell^2(\edges,\ohm) = \HD$.
  \glossary{name={\act},description={divergence operator},sort=a,format=textbf}
\end{defn}

\begin{theorem}\label{thm:act.drp=Lap}
  $\act = \Lap\drpa$, $\act\drp = \Lap$, and $\act \Pdrp = \Lap \drpa$.
  \begin{proof}
    To compute $\Lap \drpa \curr$ for a finitely supported current $\curr \in \Fin(\edges)$, let $v := \drpa \curr$ so
    \linenopax
    \begin{align*}
      \Lap(\drpa \curr)(x)
      = \Lap v(x)
      &= \sum_{y \nbr x} \cond_{xy} (v(x) - v(y))
        &&\text{defn } \Lap\\
      &= \sum_{y \nbr x} \cond_{xy} \ohm(x,y) \Pdrp \curr(x,y)
        &&\text{defn } \drpa \\
      &= \act(\Pdrp\curr)(x)
        &&\cond_{xy} = \ohm(x,y)^{-1}.
    \end{align*}
    This establishes $\act\Pdrp =\Lap \drpa$, from which the result follows by Lemma~\ref{thm:AP=A}. The second identity follows from the first by right-multiplying by \drp and applying \eqref{eqn:drp-and-drp-are-partial-isometries}. Then the third identity follows from the second by right-multiplying by \drpa and applying \eqref{eqn:drp-and-drp-are-partial-isometries} again.
  \end{proof}
\end{theorem}

\begin{remark}\label{rem:reformulated-Theorem-{thm:act.drp=Lap}}
  Theorem~\ref{thm:act.drp=Lap} may be reformulated as follows: Let $u,v \in \HE$, and $\curr := \drp v$. Then $\Lap v = u$ if and only if $\act \curr = u$. This result will help us solve $\act\curr=w$ for general initial condition $w$ in \S\ref{sec:Solving-potential-theoretic-problems}. Also, we will see in \S\ref{sec:self-adjointness-of-the-Laplacian} that \Lap is essentially self-adjoint. In that context, the results of Theorem~\ref{thm:act.drp=Lap} have a more succinct form.
\end{remark}

\begin{cor}\label{thm:A(HD)=Fin}\label{thm:AP=A}\label{thm:ker(A)=cycles}
  The kernel of \act is $\Kir \oplus \Cyc$, whence $\act \Pnbd = \act$, $\act \Pnbd^\perp = 0$ and $\act(\HD) \ci \Fin$.
  \begin{proof}
    If $\curr \in \Kir$ so that $\curr = \drp h$ for $h \in \Harm$, then $\act \curr(x) = \Lap h = 0$ follows from Theorem~\ref{thm:act.drp=Lap}. If $\curr = \charfn{\cycle}$ for $\cycle \in \Loops$, then
    \linenopax
    \begin{align*}
      \act \curr(x)
      &= \sum_{y \nbr x} \charfn{\cycle}(x,y)
       = \sum_{(x,y) \in \cycle} \charfn{\cycle}(x,y) + \sum_{(x,y) \notin \cycle} \charfn{\cycle}(x,y) = (-1+1) + 0 = 0.
    \end{align*}
    To show $\act(\HD) \ci \Fin$, it now suffices to consider $\curr \in \Nbd$. Since $\act\gh_z = \Lap \gd_z$ by Theorem~\ref{thm:act.drp=Lap}, the result follows by closing the span.
  \end{proof}
\end{cor}

In particular, Corollary~\ref{thm:A(HD)=Fin} shows that the range of \act lies in \HE, as stated in Definition~\ref{def:divergence-oper}. The identity $\act \Pnbd = \act$ implies that the solution space $\Flo(\ga,\gw)$ is invariant under minimization; see Theorem~\ref{thm:Flo-is-closed}.

\begin{remark}\label{rem:Act-vs-Lap}
  Since \act is defined without reference to \cond, \drpa ``hides'' the measure \cond from the Laplacian. To highlight similarities with the Laplacian, recall from Definition~\ref{def:current-flow} that a current \curr satisfies the homogeneous or nonhomogeneous Kirchoff laws iff $\act\curr = 0$ or $\act\curr = \gd_\ga - \gd_\gw$, respectively. In \S\ref{sec:analogy-with-complex-functions}, we consider an interesting analogy between the previous two results and complex function theory.
\end{remark}

\begin{cor}\label{thm:Lapstar-cont-in-act}
  $\Lap^\ad \ce \drp^\ad \act^\ast$ and $\act \act^\ast = \Lap\Lap^\ast$.
  \begin{proof}
    The first follows from Theorem~\ref{thm:act.drp=Lap} by taking adjoints, and the second follows in combination with Lemma~\ref{thm:act.drp=Lap}. The inclusion is for the case when \Lap may be unbounded, in which case we must be careful about domains. When $T$ is any bounded operator, $\dom T^\ad S^\ad \ci \dom(ST)^\ad$. To see this, observe that $v \in \dom T^\ad S^\ad$ if and only if $v \in \dom S^\ad$, so assume this. Then
    \linenopax
    \begin{align*}
      |\la STu,v\ra|
      = |\la Tu,S^\ad v\ra|
      \leq \|T\| \cdot \|S^\ad v\| \cdot \|u\|
      = K_v \|u\|, \q \forall u \in \dom ST,
    \end{align*}
    for $K_v = \|T\| \cdot \|S^\ad v\|$. This shows $v \in \dom(ST)^\ad$.
  \end{proof}
\end{cor}

\begin{lemma}\label{thm:A-is-contin}
  For fixed $x \in \verts$, \act is norm continuous in \curr:
  \begin{equation}\label{eqn:thm:A-is-contin}
    |\act(\curr)(x)| \leq |\cond(x)|^{1/2} \|\curr\|_\diss.
  \end{equation}
  \begin{proof}
    Using $\cond(x) := \sum_{y \nbr x} \cond_{xy}$ as in \eqref{eqn:def:cond-meas}, direct computation yields
    \linenopax
    \begin{align*}
      |\act(\curr)(x)|^2
      = \left|\sum_{y \nbr x} \curr(x,y)\right|^2 
      &= \left|\sum_{y \nbr x} \sqrt{\cond_{xy}} \sqrt{\ohm(x,y)} \curr(x,y)\right|^2 \\
      &\leq \left|\sum_{y \nbr x} \cond_{xy}\right| \sum_{y \nbr x} \ohm(x,y) |\curr(x,y)|^2 \\
      &\leq |\cond(x)| \diss(\curr),
    \end{align*}
    where we have used the Schwarz inequality and the definitions of \cond, \diss, \act.
  \end{proof}
\end{lemma}

\begin{cor}\label{thm:pointwise-bound-on-Lap}
  For $v \in \HE$, $|\Lap v(x)| \leq \cond(x)^{1/2} \|v\|_\energy$.
  \begin{proof}
    Apply Theorem~\ref{thm:A-is-contin} to $\curr = \drp v$ and use the second claim of Theorem~\ref{thm:act.drp=Lap}.
  \end{proof}
\end{cor}

One consequence of the previous lemma is that the space of functions satisfying the nonhomogeneous Kirchhoff condition \eqref{eqn:Kirchhoff-nonhomog} is also closed, as we show in Theorem~\ref{thm:Flo-is-closed}.

In Remark~\ref{rem:reproducing-kernel}, we discussed some reproducing kernels for operators on \HE; we now introduce one for the divergence operator \act, using the weighted edge neighbourhoods $\{\gh_z\}$ of Definition~\ref{def:reweighted-edge-nbd}.

\begin{lemma}\label{thm:reproducing-kernel-for-A}
  The currents $\{\gh_z\}$ form a reproducing kernel for \act.
  \begin{proof}
    By Lemma~\ref{thm:A-is-contin}, the existence of a reproducing kernel follows from Riesz's Theorem. Since it must be of the form $\act(\curr)(z) = \la k_z, \curr \ra_\diss$, we verify
    \linenopax
    \begin{align*}
      \la \gh_z, \curr \ra_\diss
      &= \sum_{(x,y) \in \edges} \ohm(x,y) \gh_z(x,y) \curr(x,y)
      = \sum_{y \nbr z} \curr(z,y)
      = \act(\curr)(z).
      \qedhere
    \end{align*}
  \end{proof}
\end{lemma}

\section{Analogy with calculus and complex variables}
\label{sec:analogy-with-complex-functions}
The material in this book bears many analogies with vector calculus and complex function theory. Several points are obvious, like the existence and uniqueness of harmonic functions and the discrete Gauss-Green formula of Lemma~\ref{thm:E(u,v)=<u,Lapv>}. In this section, we point out a couple more subtle comparisons.

The drop operator \drp is analogous to the complex derivative
\linenopax
\begin{align*}
  \del = \frac{d}{dz} = \frac{\del}{\del z} := \frac12\left(\frac{\del}{\del x} + \frac1\ii\frac{\del}{\del y}\right),
\end{align*}
as may be seen from the discussion of the compatibility problem in \S\ref{sec:compatibility-problem}. Recall from the proof of Theorem~\ref{thm:HD=Nbd+Kir+Cyc} that Lemma~\ref{thm:induced-currents-satisfy-cycle-condition} expresses the fact
\linenopax
\begin{align*}
  \la \curr, \charfn{\cycle} \ra_\diss = 0, \;\forall \cycle \in \Loops
  \q\iff\q \exists v \in \HE \text{ such that } \drp v = \curr.
\end{align*}
This result is analogous to Cauchy's theorem: if $v$ is a complex function on an open set, then $v=f'$ (that is, $v$ has an antiderivative) if and only if every closed contour integral of $v$ is 0. Indeed, even the proofs of the two results follow similar methods.

The divergence operator \act may be compared to the Cauchy-Riemann operator
\linenopax
\begin{align}\label{eqn:Cauchy-Riemann-operator}
  \bar\del = \frac{\del}{\del \bar z} := \frac12\left(\frac{\del}{\del x} - \frac1\ii\frac{\del}{\del y}\right).
\end{align}
Indeed, in Theorem~\ref{thm:act.drp=Lap} we found that $\act \drp = \cond \Lap$, which may be compared with the classical identity $\bar\del \del = \frac14\Lap$. The Cauchy-Riemann equation $\bar \del f = 0$ characterizes the analytic functions, and $\act\curr=0$ characterizes the currents satisfying the homogeneous Kirchhoff law; see Definition~\ref{def:Kirchhoff's-law}.

In \S\ref{sec:Solving-potential-theoretic-problems}, we give a solution for the inhomogeneous equation $\act \curr = w$ when $w$ is given and satisfies certain conditions. The analogous problems in complex variables are as follows: let $W \ci \bC$ be a domain with a smooth boundary $\bd W$, and let $\bar \del$ be the Cauchy-Riemann operator \eqref{eqn:Cauchy-Riemann-operator}. Suppose that \gn is a compactly supported $(0,1)$ form in $W$. We consider the boundary value problem
\linenopax
\begin{align*}
  \bar\del f = \gn, \q\text{with}\q \bar\del\gn=0.
\end{align*}
The Bochner-Martinelli theorem states that the solution $f$ is given by the following integral representation:
\linenopax
\begin{equation}\label{eqn:Bochner-Martinelli-formula}
  f(z) = \int_{\bd W} f(\gz) \gw(d\gz,z) - \int_W \gn(\gz) \wedge \gw(d\gz,z),
\end{equation}
where \gw is the Cauchy kernel. In fact, this theorem continues to hold when $W$ is a domain in $\bC^n$, if one uses the Bochner-Martinelli kernel
\linenopax
\begin{equation}\label{eqn:Bochner-Martinelli-kernel}
  \gw(\gz,z)
  = \frac{(n-1)!}{2\gp\ii |\gz-z|^{2n}} \sum_{k=1}^n (\bar\gz_k - \bar z_k) d\bar\gz_1 \wedge d\gz_1 \wedge \dots (\hat j) \dots \wedge d\bar\gz_n \wedge d\gz_n,
\end{equation}
where $\hat j$ means that the term $d\bar\gz_j \wedge d\gz_j$ has been omitted, and where
\linenopax
\begin{align*}
  \del \gf = \sum\nolimits_k \frac{\del\gf}{\del z_k} dz_k,
  \text{ and }
  \bar\del \gf = \sum\nolimits_k \frac{\del\gf}{\del \bar z_k} d\bar z_k.
\end{align*}
Indeed, in Lemma~\ref{thm:reproducing-kernel-for-A}, we obtain a reproducing kernel for \act; this is analogous to the Bochner-Martinelli kernel $K(z,w)$; see \cite{Kyt95} for more on the Bochner-Martinelli kernel.

Theorem~\ref{thm:vx-dense-in-HE} shows that $v_x$ is analogous to the Bergman kernel, which reproduces the holomorphic functions within $L^2(\gW)$, where $\gW \ci \bC$ is a domain. Indeed, the Bergman kernel is also associated with a metric, the Bergman metric, which is defined by
\begin{equation}\label{eqn:bergman-metric}
  d_B(x,y) := \inf_\gg\int_\gg\left|\frac{\del\gg}{\del t}(t)\right|,
\end{equation}
where the infimum is taken over all piecewise $C^1$ paths \gg from $x$ to $y$; cf.~\cite{Kra01}.

\section{Solving potential-theoretic problems with operators}
\label{sec:Solving-potential-theoretic-problems}

We begin by discussing the minimizing nature of the projections \Pfin and \Pnbd. Theorem~\ref{thm:proj-is-minimizer} shows how \drpa solves the compatibility problem of \ref{sec:compatibility-problem}: Given a current flow $\curr \in \HD$, there does not necessarily exist a potential function $v \in \HE$ for which $\drp v = \curr$. Nonetheless, there is a potential function associated to \curr which satisfies $\drp v = \Pnbd \curr$, and it can be found via the minimizing projection. Consequently, Theorem~\ref{thm:proj-is-minimizer} can be seen as an analogue of Theorem~\ref{thm:minimal-flows-are-induced-by-minimizers}.

Theorem~\ref{thm:Flo-is-closed} shows that the solution space $\Flo(\ga,\gw)$ is invariant under \Pnbd. Coupled with the results of Theorem~\ref{thm:proj-is-minimizer}, this shows that if one can find any solution $\curr \in \Flo(\ga,\gw)$, one can obtain another solution to the same Dirichlet problem with minimal dissipation, namely, $\Pnbd \curr$. 

\subsection{Resolution of the compatibility problem}
\label{sec:Resolution-of-the-compatibility-problem}
In this section we relate the projections
\linenopax
\begin{align*}
  \Pfin:\HE \longrightarrow \Fin
  \q\text{and}\q
  \Pnbd:\HD \longrightarrow \Nbd = \drp \Fin
\end{align*}
of Definition~\ref{def:projection-notation} to some questions which arose in \S\ref{sec:currents-and-potentials}. The operators \Pfin and \Pnbd are \emph{minimizing projections} because they strip away excess energy/dissipation due to harmonic or cyclic functions:
\begin{itemize}
  \item If $v \in \Pot(x,y)$, then $\Pfin v$ is the unique minimizer of \energy in $\Pot(x,y)$.
  \item If $\curr \in \Flo(x,y)$, then $\Pnbd \curr$ is the unique minimizer of \diss in $\Flo(x,y)$.
\end{itemize}
In a similar sense, \Pdrp is also a minimizing projection.

Probability notions will play a key role in our solution to questions about {divergence} in electrical networks (Definition~\ref{def:divergence-oper}), as well as our solution to a potential equation. 
The divergence will be important again in \S\ref{sec:forward-harmonic-functions} where we use it to provide a foundation for a probabilistic model which is dynamic (in contrast to other related ideas in the literature) in the sense that the Markov chain is a function of a current \curr, which may vary.

\begin{theorem}\label{thm:proj-is-minimizer}
  Given $v \in \HE$, there is a unique $\curr\in \HD$ which satisfies $\drpa \curr = v$ and minimizes $\|\curr\|_\diss$. Moreover, it is given by $\Pdrp \curr$, where \curr is any solution of $\drpa \curr = v$.
  \begin{proof}
    Given $v \in \HE$, we can find some $\curr \in \HD$ for which $\drpa \curr = v$, by Theorem~\ref{thm:drp-and-drp-are-partial-isometries}.
    Then the orthogonal decomposition $\curr = \Pdrp \curr + \Pdrp^\perp \curr$ gives
    \begin{equation}\label{eqn:proj-is-minimizer}
      \|\curr\|^2_\diss
      = \|\Pdrp \curr\|^2_\diss + \|\Pdrp^\perp\curr\|^2_\diss
      \geq \|\Pdrp \curr\|^2_\diss,
    \end{equation}
    so that $\|\curr\|_\diss \geq \|\Pdrp \curr\|_\diss$ shows $\Pdrp \curr$ minimizes the dissipation norm. Finally, note that $\drpa\Pdrp \curr = \drpa \drp\drpa \curr = \drpa \curr$, by Corollary~\ref{thm:drop-oper-adjoint-formula}.
  \end{proof}
\end{theorem}

\begin{cor}\label{thm:drpadj-is-a-solution-oper}
  \drpa is a solution operator in the sense that if \curr is any element of
  $\HD$ then $\drpa \curr$ is the unique element $v \in \HE$ for which $\drp v = \Pdrp \curr$.
\end{cor}

\begin{cor}\label{thm:DPnbd=EPfin}
  $\diss\Pnbd = \energy\Pfin$. Hence for $\curr = \drp(v_x-v_y)$,
  \begin{align}\label{eqn:thm:DPnbd=EPfin}
    R^F(x,y) &= \energy(\drpa \curr)^{1/2} = \diss(\Pdrp \curr)^{1/2}, 
    \q\text{and}\q \\
    R^W(x,y) &= \min\{\diss(\curr)^{1/2} \suth \curr \in \Flo(x,y)\}. \notag
  \end{align}
  \begin{proof}
    Given $\curr \in \HD$, let $\curr_0 = \Pnbd \curr$. Then define $v$ by
    $v := \drpa\Pnbd \curr_0 = \drpa\Pnbd \curr$. Applying \drp to both sides gives $\drp v = \Pnbd \curr$ by \eqref{eqn:drp-and-drp-are-partial-isometries} (since $\Pnbd \leq \Pdrp$) so that taking dissipations and applying Lemma~\ref{thm:drp-is-an-isometry} gives $\diss(\Pnbd \curr) = \diss(\drp v) = \energy(v) = \energy(\Pfin v)$, because $\ran \drpa \Pnbd \ci \Fin$ by Theorem~\ref{thm:HD=Nbd+Kir+Cyc} and Theorem~\ref{thm:drp-and-drp-are-partial-isometries}.
  \end{proof}
\end{cor}

\begin{theorem}\label{thm:Flo-is-closed}
  For any $\ga,\gw \in \verts$, the subset $\Flo(\ga,\gw)$ is closed with respect to $\|\cdot\|_\diss$ and invariant under \Pnbd.
  \begin{proof}
    From \eqref{eqn:Kirchhoff-nonhomog} and \eqref{eqn:def:divergence-oper}, we have that $\curr \in \Flo(\ga,\gw)$ if and only $\act\curr = \gd_\ga - \gd_\gw$. Suppose that $\{\curr_n\} \ci \Flo(\ga,\gw)$ is a \seq of currents for which $\curr_n \limmode[\text{\tiny \diss}] \curr$. Then $\act \curr_n = \gd_\ga - \gd_\gw$ for every $n$, and from Lemma~\ref{thm:A-is-contin}, the inequality
    \linenopax
    \begin{align*}
      |(\act \curr_n)(x) - (\act \curr)(x)|
      \leq |\cond(x)|^{1/2} \|\curr_n - \curr\|_\diss
    \end{align*}
    gives $\act\curr(x) = \gd_\ga - \gd_\gw$. Note that $x$ is fixed, and so $\cond(x)$ is just a constant in the inequality above.

    For invariance, note that $\act\Pnbd=\act$ by Corollary~\ref{thm:AP=A}. Then $\curr \in \Flo(\ga,\gw)$ implies
    \linenopax
    \begin{align*}
      \act\Pnbd\curr &= \act\curr = \gd_\ga - \gd_\gw
      \q\implies\q \Pnbd\curr \in \Flo(\ga,\gw).
      \qedhere
    \end{align*}
  \end{proof}
\end{theorem}

Since \Pnbd is a subprojection of $\Pcyc^\perp$ and $\Pkir^\perp$, we have an easy corollary.

\begin{cor}\label{thm:Pnbd<Pdrp,Pkirperp}
  For any $\ga,\gw \in \verts$, $\Flo(\ga,\gw)$ is invariant under $\Pdrp = \Pcyc^\perp$ and $\Pkir^\perp$.
\end{cor}

\begin{remark}\label{rem:rank-Pperp-is-numcycles}
  Putting these tools together, we have obtained an extremely simple method for solving the equation $\Lap v = \gd_\ga - \gd_\gw$.
\begin{enumerate}
  \item Find any current $\curr \in \Flo(\ga,\gw)$. This is trivial; one can simply take the characteristic function of a path from  \ga to \gw.
  \item Apply \Pnbd to \curr to ``project away'' harmonic currents and cycles.
  \item Apply \drpa to $\Pnbd\curr$. Since $\Pnbd\curr \in \drp\Fin$, this only requires an application of Ohm's law in reverse as in \eqref{eqn:lem:drop-oper-adjoint-formula}.
\end{enumerate}
Then $v = \drpa\Pnbd\curr$ is the desired energy-minimizing solution (since any harmonic component is removed). As a bonus, we already obtained the current \Pdrp\curr induced by $v$.
 The only nontrivial part of the process described above is the computation of \Pnbd. For further analysis, one must understand the cycle space \Cyc of \Graph and the space \Kir of harmonic currents. We hope to make progress on this problem in a future paper, see Remark~\ref{rem:fut:rank-Pperp-is-numcycles}. 
\end{remark}

\version{}{
\subsection{Generalization of the current problem}
\label{sec:Generalization-of-the-problem}
%
%
We have already seen how \drpa and \Pnbd solve certain problems on the resistance network. Given any current \curr on the network, \drpa produces the unique associated potential function $v$ (i.e., such that $\Pdrp \curr$ is induced by $v$). Moreover, $\Pdrp\curr$ is the unique current flow associated to \curr which is physically realistic. In this section, we show how to solve resistance network problems with more general conditions by applying the tools from Hilbert space theory developed in the previous section, in this fashion. To be specific, let $X = \{x_1, x_2, \dots, x_n, \dots\} \ci \verts$ be a specified subset of vertices and consider $f = \sum_{x \in X} f(x) \gd_{x}$ as prescribing weights on this set. Then $\act\curr = f$ means that the current \curr satisfies the Kirchhoff laws
\begin{equation}\label{eqn:Kirchhoff-generalized}
  \sum_{y \nbr x} \curr(x,y) =
  \begin{cases}
    f(x), &x \in X, \\
    0, &else,
  \end{cases}
\end{equation}
in a direct generalization of \eqref{eqn:Kirchhoff-nonhomog}. A solution \curr to \eqref{eqn:Kirchhoff-generalized}, that is, to $\act\curr = f$, will represent the current flow induced by imposing the specified voltages at the appropriate vertices, see the first paragraph of \S\ref{sec:currents}.

\begin{defn}\label{def:flows-of-X}
  Let $X = \{x_1, x_2, \dots, x_n\} \ci \verts$. For $f = \sum_X f(x) \gd_{x}$, denote
  \begin{equation}\label{eqn:def:flows-of-X}
    \Flo_f := \{\curr : \edges \to \bR \suth \act\curr = f\}.
  \end{equation}
\end{defn}
  \glossary{name={$\Flo_f$},description={collection of flows specific by $f$},sort=f,format=textbf}
When $X$ is finite, the hypotheses of Lemma~\ref{thm:solving-AI=f} are automatically met, and so $\Flo_f$ is nonempty. Theorem~\ref{thm:flows-of-X-solve-dirichlet} shows that $\Flo_f$ gives solutions to a certain inhomogeneous Dirichlet problem.

If $X=\verts$, then the solution to $\act\curr=f$ is given by $\curr = \drp f$, as is shown in Lemma~\ref{thm:potential-induces-current} (in fact, this is precisely the content of \eqref{eqn:pot-gives-curr}). However, if $X \subsetneq \verts$, the following solution must be used.

\begin{lemma}\label{thm:solving-AI=f}
  For each $x \in \verts$, let $w_x$ be a solution of $\Lap w = \gd_x$. Suppose $X \ci \verts$ and let $f = \sum_{x \in X} f(x) \gd_{x}$. If $u = \sum_{x \in X} f(x) w_x \in \HE$, then $\curr=\drp u$ solves the initial value problem
  \begin{equation}\label{eqn:solving-AI=f}
    \act \curr = f.
  \end{equation}
  \begin{proof}
    Assume initially that $X$ is finite.
    For $\curr := \drp\left(\sum_{x \in X} f(x) w_x \right)$,
    \linenopax
    \begin{align*}
      \act \curr
      = \act \drp \sum_{x \in X} f(x) w_x
      = \Lap \sum_{x \in X} f(x) w_x
      = \sum_{x \in X} f(x) \gd_{x}
      = f.
    \end{align*}
    The second equality comes by Theorem~\ref{thm:act.drp=Lap}; the third comes by linearity of \Lap and the definition of $w_x$.

    If $X$ is infinite, then let $\{X_n\}$ be a nested \seq of sets increasing to $X$, i.e., $X_n \ci X_{n+1}$ and $X = \bigcup_n X_n$. Define $f_n := \sum_{x \in X_n} f(x) w_x$ so that $\|f - f_n\|_\energy \to 0$, and let $\curr_n$ be the solution obtained by the method just described, so that $\act\curr_n = f_n$. With $n>m$, we show that $\{\curr_n\}$ is Cauchy:
    \linenopax
    \begin{align*}
      \|\curr_n - \curr_m\|_\diss^2
      &= \diss\drp\left(\sum_{X_n} f(x)w_x - \sum_{X_m} f(x)w_x\right)
        &&\drp \text{ is linear} \\
      &= \energy\left(\sum_{x \in X_n \less X_m} f(x) w_x\right)
        &&\text{Cor.~\ref{thm:DPnbd=EPfin}} \\
      &\leq \energy\left(\sum_{x \in X \less X_m} f(x) w_x\right)
      = \energy(f - f_m) \limas{n} 0.
    \end{align*}
    Thus, $\curr := \lim \curr_n$ is well defined and $\act \curr = f$.
  \end{proof}
\end{lemma}

\begin{lemma}\label{thm:flow(a,w)-is-closed}
  For $f \in \HE$, $\Flo_f \ci \HD$ is invariant under \Pnbd and closed.
\end{lemma}

The proof of Lemma~\ref{thm:flow(a,w)-is-closed} is identical to that of Theorem~\ref{thm:Flo-is-closed}. Its use lies in the fact that if \curr solves $\act \curr=w$, then $\Pnbd\curr$ is also a solution, and thus one has a solution which lies in the range of \drp.

\begin{theorem}\label{thm:flows-of-X-solve-dirichlet}
  The Dirichlet problem
  \begin{equation}\label{eqn:flows-of-X-solve-dirichlet}
    \Lap v = w
  \end{equation}
  is solved by $v = \drpa \curr$, where $\curr \in \Flo_w$.
  \begin{proof}
    The set $\Flo_w$ is nonempty by Lemma~\ref{thm:solving-AI=f}, from which the result is immediate by Lemma~\ref{thm:act.drp=Lap}: $\Lap v = \Lap \drpa \curr = \act\curr = w$.
  \end{proof}
\end{theorem}

\begin{remark}\label{rem:flows-solve-Dirichlet}
  The significance of Theorem~\ref{thm:flows-of-X-solve-dirichlet} is that it allows us to solve a problem in \HE by working entirely in \HD, which is just an $\ell^2$ space and much better understood.
\end{remark}

}

\section{Remarks and references}
\label{sec:Remarks-and-References-he-and-hd}

After completing a first draught of this book, we discovered several of the results of this section in \cite{Lyons:ProbOnTrees} and \cite{Soardi94}. Both of these texts are excellent; Lyons emphasizes connections with probability and \S2 and \S9 are most pertinent to the present discussion, and Soardi emphasizes the (co)homological perspective and parallels with vector calculus. 

The subspace of currents spanned by edge neighbourhoods $\Nbd = \clspn[\diss]{\drp \Fin} = \clspn[\diss]{\gh_z}$ is called $\bigstar$ in \cite[\S2 and \S9]{Lyons:ProbOnTrees}, and the subspace of cycles $\Cyc := \clspn[\diss]{\charfn{\cycle}}$ is called {\large $\Diamond$}.

The reader may also wish to consult \cite{DuRo08, Sto08, DMS07, BaBa05, Chu07, ChRi06} with regard to the material in this chapter.

%% file: probab-interp.tex

\chapter{Probabilistic interpretations}
\label{sec:Probabilistic-interpretation}

\headerquote{From its shady beginnings devising gambling strategies and counting corpses in medieval London, probability theory and statistical inference now emerge as better foundations for scientific models, especially those of the process of thinking and as essential ingredients of theoretical mathematics, even the foundations of mathematics itself.}{---~David~Mumford}

\headerquote{God not only plays dice. He also sometimes throws the dice where they cannot be seen.}{---~S.~Hawking}

In \S\ref{sec:the-boundary}, we constructed a measure \prob on $\Schw'$, where $\Schw \ci \HE \ci \Schw'$ is a certain Gel'fand triple. In this section, we develop a different but analogous measure on the space of infinite paths in $\bd \Graph$. We carry out this construction for harmonic functions on $(\Graph, \cond)$ in \S\ref{sec:path-space-of-a-general-random-walk}, where the measure is defined in terms of transition probabilities $p(x,y) = \frac{\cond_{xy}}{\cond(x)}$ of the random walk, and the associated cylinder sets. When the random walk on $(\Graph, \cond)$ is transient, the current induced by a monopole gives a unit flow to infinity; such a current induces an orientation on the edges \edges and a new, naturally adapted, Markov chain. The state space of this new process is also \verts, but the transition probabilities are now defined by the induced current $p(x,y) = \frac{\curr(x,y)}{\act_\curr(x)}$. We call the fixed points of the corresponding transition operator the ``forward-harmonic'' functions, and carry out the analogous construction for them in \S\ref{sec:forward-harmonic-functions}. The authors are presently working to determine whether or not these measures can be readily related to each other or the measure \prob of \S\ref{sec:Gel'fand-triples-and-duality}. 


\section{The path space of a general random walk}
\label{sec:path-space-of-a-general-random-walk}

We begin by recalling some terms from \S\ref{sec:trace-resistance}, and providing some more detail. Let $\cpath = (x_0, x_1, x_2, \dots, x_n)$ be any finite path starting at $x=x_0$. The probability of a random walk started at $x$ traversing this path is
  \linenopax
  \begin{equation}\label{eqn:def:path-probability-recall}
    \prob(\cpath) := \prod_{k=1}^{n} p(x_{k-1},x_k),
  \end{equation}
  where $p(x,y) := \frac{\cond_{xy}}{\cond(x)}$ is the probability that the walk moves from $x$ to $y$ as in \eqref{eqn:def:p(x,y)}.
  \glossary{name={\cpath},description={path},sort=G}
  \glossary{name={$\prob^{(\cond)}$},description={a measure on the space of all paths in $G$},sort=Pc}
  \glossary{name={$p(x,y)$},description={transition probability of the random walk on the network: $\cond_{xy}/\cond(x)$},sort=prob}
  This intuitive notion can be extended via Kolmogorov consistency to the space of all infinite paths starting at $x$. Let $X_n(\cpath)$ denote the \nth coordinate of \cpath; one can think of \cpath as an event and $X_n$ as the random walk (a random variable), in which case 
  \linenopax
  \begin{equation}\label{eqn:def:X_n}
    X_n(\cpath) = \text{location of the random walk at time $n$.}
  \end{equation}

\begin{defn}\label{def:Paths}\label{def:prob^cond}
  Let \Paths denote the \emph{space of all infinite paths} \cpath in $(\Graph, \cond)$. Then a \emph{cylinder set} in \Paths is specified by fixing the first $n$ coordinates:
  \linenopax
  \begin{align}\label{eqn:def:cylinder-set-in-Paths(x)}
    \Paths_{(x_1,x_2,\dots,x_n)} := \{\cpath \in \Paths \suth X_k(\cpath) = x_k, k=1,\dots,n\}.
  \end{align}
  Define $\prob^{(\cond)}$ on cylinder sets by 
  \linenopax
  \begin{align}\label{eqn:def:prob^cond}
    \prob^{(\cond)}(\Paths_{(x_1,x_2,\dots,x_n)}) := \prod_{i=1}^n p(x_{i-1},x_i).
  \end{align}
\end{defn}

\begin{remark}
  It is clear from Definition~\ref{def:prob^cond} that the probability of a random walk following the finite path $\cpath = (x_0, x_1, x_2, \dots, x_n)$ is equal to the measure of the set of all infinite walks which agree with \cpath for the first $n$ steps: combining \eqref{eqn:def:path-probability-recall} and \eqref{eqn:def:prob^cond} gives $\prob^{(\cond)}(\Paths_{(x_1,x_2,\dots,x_n)}(x)) = \prob(\cpath)$. Observe that \eqref{eqn:def:prob^cond} is a conditional probability:
  \linenopax
  \begin{align}\label{eqn:cylinder-sets-as-conditional-prob}
    \prob^{(\cond)}(\Paths_{(x_1,x_2,\dots,x_n)}(x)) 
    = \prob^{(\cond)}\{\cpath \in \Paths(x) \,|\, X_k(\cpath) = x_k, k=1,\dots,n\}.
  \end{align}
\end{remark}

\begin{remark}[Kolmogorov consistency]
  \label{rem:Kolmogorov-consistency}
  We use Kolmogorov's consistency theorem to construct a measure on the space of paths beginning at vertex $x \in \verts$, see \cite[Lem.~2.5.1]{Jor06} for a precise statement of this extension principle in its function theoretic form and \cite[Exc.~2.4--2.5]{Jor06} for the method we follow here. The idea is that we consider a sequence of functionals $\{\gm^{(n)}\}$, where $\gm^{(n)}$ is defined on 
  \begin{equation}\label{eqn:def:An}
  \sA_n := \spn\{\charfn{\Paths_{(x_0,\dots,x_n)}} \suth x_i \nbr x_{i-1}, i=1,\dots,n\}.\end{equation}
Alternatively,
\begin{equation}\label{eqn:def:An-alt}
  \sA_n := \{f:\Paths \to \bR \suth f(\cpath_1)=f(\cpath_2) \text{ whenever } X_k(\cpath_1)=X_k(\cpath_2) \text{ for } k \leq n\}.
\end{equation}
That is, an element of $\sA_n$ cannot distinguish between two paths which agree for the first $n$ steps. This means that $\gm^{(n)}$ is a ``simple functional'' in the sense that it is constant on each cylinder set of level $n$:
\begin{equation}\label{eqn:def:An}
  \gm^{(n)}[f] = \sum_{x_0,\dots,x_n} a_{(x_0,\dots,x_n)} \gm^{(n)}[\charfn{\Paths_{(x_0,\dots,x_n)}}].
\end{equation}
If the functionals $\gm^{(n)}$ are mutually consistent in the sense that $\gm^{(n+1)}[f] = \gm^{(n)}[f]$, then Kolmogorov's consistency theorem gives a unique Borel probability measure on the space of all paths. More precisely, Kolmogorov's theorem gives the existence of a limit functional which is defined for functions on paths of infinite length, and this corresponds to a measure by Riesz's Theorem; see \cite{Jor06,Kol56}.
\end{remark}

In the following, we let \one denote the constant function with value equal to 1.

\begin{theorem}[Kolmogorov]\label{thm:Kolmogorov-extension-thm}
  If each $\gm^{(n)}:\sA_n \to \bC$ is a positive linear functional satisfying the consistency condition
  \begin{equation}\label{eqn:Kolmogorov-consistency}
    \gm^{(n+1)}[f] = \gm^{(n)}[f], \qq\text{for all } f \in \sA_n,
  \end{equation}
  then there exists a positive linear functional \gm defined on the space of functions on infinite paths such that
  \begin{equation}\label{eqn:Kolmogorov-functional}
    \gm[f] = \gm^{(n)}[f], \qq f \in \sA_n,
  \end{equation}
  where $f$ is considered as a function on an infinite path which is zero after the first
  $n$ edges. 
  Moreover, if we require the normalization $\gm^{(n)}[\one] = 1$, then \gm is determined uniquely.
\end{theorem} 

  We now show that $\prob^{(\cond)}$ extends to a natural probability measure on the space of infinite paths $\Paths(x)$.

\begin{theorem}\label{thm:P^cond-exists}
  For $(\Graph,\cond)$, there is a unique measure $\prob^{(\cond)}$ defined on $\Paths$ which satisfies
  \begin{equation}\label{eqn:P^cond-consistency}
    \bE[V] 
    = \int_{\Paths} V \,d\prob^{(\cond)} 
    = \int_{\Paths} V \,d\prob^{(\cond,n)} 
    = \bE^{(n)}[V], \qq \forall V \in \sA_n.
  \end{equation}
  \begin{proof}
    We must check condition \eqref{eqn:Kolmogorov-consistency} for $\gm^{(n)} = \prob^{(\cond,n)}$, defined by 
    \linenopax
    \begin{align*}
      \prob^{(\cond,n)}(\charfn{\Paths_{(x_0,\dots,x_n)}}) 
      := \prod_{i=1}^n p(x_{i-1},x_i)    
    \end{align*}
    with \eqref{eqn:def:prob^cond} in mind. 
    Think of $V \in \sA_n$ as an element of $\sA_{n+1}$ and apply $\prob^{(\cond,n+1)}$ to it:
    \linenopax
    \begin{align*}
      \prob^{(\cond,n+1)}[V]
      &= \sum_{x_0,\dots,x_{n+1}} a_{(x_0,\dots,x_{n+1})} 
        \prob^{(\cond,n+1)}(\charfn{\Paths_{(x_0,\dots,x_{n+1})}}) \\
      &= \sum_{x_0,\dots,x_n} \sum_{x_{n+1}} a_{(x_0,\dots,x_{n})} \prod_{i=1}^{n} p(x_{i-1},x_i)p(x_n,x_{n+1}) \\
      &= \sum_{x_0,\dots,x_n} a_{(x_0,\dots,x_{n})} \prod_{i=1}^{n} p(x_{i-1},x_i) \sum_{x_{n+1} \nbr x_n} p(x_n,x_{n+1}) \\
      &= \prob^{(\cond,n)} [V], \vstr[2.5]
      \qedhere
    \end{align*}
    since $\sum_{x_{n+1} \nbr x_n} p(x_n,x_{n+1}) = 1$. For the second equality, note that $f \in \sA_n$, so we can use the same constant $a$ for each $(n+2)$-tuple that begins with $(x_0,\dots,x_n)$. 
  \end{proof}
\end{theorem}

\subsection{A boundary representation for the bounded harmonic functions} 

\begin{defn}\label{def:bdd-cocycle}
  A \emph{cocycle} $V:\Paths \to \bR$ is a measurable function on the infinite path space which is independent of the first finitely many vertices in the path:
  \begin{equation}\label{eqn:def:bdd-cocycle-condition}
    V(\cpath) = V(\shift\cpath),
  \end{equation}
  where \shift is the shift operator, i.e., if $\cpath=(x_0,x_1,x_2,\dots)$, then $\shift\cpath=(x_1,x_2,x_3,\dots)$.
\end{defn}

Intuitively, a cocycle is a function on the boundary $\bd \Graph$; it depends only on the asymptotic trajectory of a path/random walk. A cocycle does not care where the random walk began, only where it goes.
\version{}{\marginpar{Make this precise}}
More precisely, a cocycle is a special kind of martingale, as we will see below. 

The goal of this section is to show that the bounded harmonic functions are in bijective correspondence with the cocycles; see Theorem~\ref{thm:harmonic-functions=cocycles}. That is, the formula 
  \begin{equation}\label{eqn:def:harmonic-vs-cocycle}
    h(x) = \Ex_x[V]Ê
  \end{equation}
spells out a bijective correspondence between functions $h \in \Harm$, and cocycles $V$ on the space of infinite paths. Our present concern is the space of all bounded harmonic functions; we will presently consider the class of finite-energy functions. A good reference for this section is \cite[Thm.~2.7.1]{Jor06}.

Note that the left hand side of \eqref{eqn:def:harmonic-vs-cocycle} involves no measure theory, in contrast to the right-hand side, where the expectation refers to the integration of cocycles $V$ against the probability measure $\prob^{(\cond)}$. The underlying Borel probability space of $\prob^{(\cond)}$ is the \gs-algebra of measurable sets generated by the cylinder sets in \Paths, i.e., by the subsets in \Paths which fix only a finite number of places (in the infinite paths). 

The condition on a measurable function $V$ on \Paths which accounts for $h$ defined by \eqref{eqn:def:harmonic-vs-cocycle} being harmonic is that $V$ is invariant under a finite left shift; cf.~\eqref{eqn:def:bdd-cocycle-condition}. It turns out that in making the integrals $\Ex_x(V)$ precise, the requirement that $V$ be measurable is a critical assumption. In fact, there is a variety of non-measurable candidates for such functions $V$ on \Paths.

\begin{defn}\label{def:path-expectation}
  For any measurable function $V:\Paths \to \bR$, we write 
  \begin{equation}\label{eqn:def:path-expectation}
    \bE_x[V] := \bE[V \,|\, X_0=x] = \int_{\Paths(x)} V(\cpath) \,d\prob^{(\cond)}
  \end{equation}
  for the expected value of $V$, conditioned on the path starting at $x$.
\end{defn}

\begin{lemma}\label{thm:harmonic-as-integral-identity}
  For $h \in \Harm$ and any $n=1,2,\dots,$
  \linenopax
  \begin{align}\label{eqn:int(h(Xn)=h}
    h(x) = \int_{\Paths_x} h \comp X_n \, d\prob^{(\cond)}.
  \end{align}
  \begin{proof}  
    By the definition of the cylinder measure $d\prob^{(\cond)}$ (Definition~\ref{def:prob^cond}),
    \linenopax
    \begin{align}\label{eqn:int(h(Xn)=h}
      \int_{\Paths_x} h \comp X_1 \, d\prob^{(\cond)}
      = \sum_{y \nbr x} p(x,y) \int_{\Paths_y} h \comp X_0 \, d\prob^{(\cond)}
      = \sum_{y \nbr x} p(x,y) h(y) \int_{\Paths_y} \, d\prob^{(\cond)}
      = \Prob h(x),
    \end{align}
    so iteration and $\Prob h = h$ gives $\int_{\Paths_x} h \comp X_n \, d\prob^{(\cond)} = h(x)$ for every $n=1,2,\dots$.
  \end{proof} 
\end{lemma}
    
\begin{theorem}\label{thm:harmonic-functions=cocycles}
  The bounded harmonic functions are in bijective correspondence with the cocycles. More precisely, if $V$ is a cocycle, then it defines a harmonic function via
  \linenopax
  \begin{equation}\label{eqn:hV=E[V]}
    h_V(x) := \bE_x[V].
  \end{equation}
  Conversely, if $h$ is harmonic, then it defines a cocycle via
  \linenopax
  \begin{equation}\label{eqn:Vh-is-martingale-limit}
    V_h(\cpath) := \lim_{n \to \iy} h(X_n(\cpath)), \qq \text{for } \prob^{(\cond)}-\text{a.e. } \cpath \in \Paths(x).
  \end{equation}
  \begin{proof}
    \fwd Recall that $\Lap = \cond - \Trans$; we will show that $\cond h_V = \Trans h_V$ whenever $V$ is a cocycle. If $\Paths_{(x,y)} := \{\cpath \in \Paths(x) \suth X_1(\cpath) = y\}$, then $\Paths(x) = \bigcup_{y \nbr x} \Paths_{(x,y)}$ is a disjoint union and
    \linenopax
    \begin{align*}
     h_V(x) 
     = \bE_x[V] 
      = \int_{\Paths(x)} V(\cpath) \,d\prob^{(\cond)} 
      = \sum_{y \nbr x} \int_{\Paths_{(x,y)}} V(\cpath) \,d\prob^{(\cond)}.
    \end{align*}
    For each $\cpath \in \Paths_{(x,y)}$, one has $\prob(\cpath) = \prob(x,y) \prob(\shift \cpath) = p(x,y) \prob(\shift \cpath)$ by \eqref{eqn:def:path-probability-recall}, whence
    \linenopax
    \begin{align*}
      \cond(x) h_V(x) 
      &= \cond(x) \sum_{y \nbr x}\int_{\Paths_{(x,y)}} \negsp[15] p(x,y) V(\shift\cpath) \,d\prob^{(\cond)} 
      = \sum_{y \nbr x} \cond_{xy} \int_{\Paths(y)} V(\cpath) \,d\prob^{(\cond)} 
     = \Trans h_V(x), 
    \end{align*}
    where the cocycle property \eqref{eqn:def:cocycle-condition} is used for the second equality.
    
    \bwd Now let $h$ be a bounded harmonic function. Since 
    \linenopax
    \begin{align*}
      \lim_{n \to \iy} h(X_n(\cpath))
      = \lim_{n \to \iy} h(X_{n+1}(\cpath))
      = \lim_{n \to \iy} h(X_n(\shift\cpath)), 
    \end{align*}
    the cocycle property \eqref{eqn:def:cocycle-condition} is obviously satisfied whenever the limit exists. Let $\gS_n$ denote the \gs-algebra generated by the cylinder sets of level $n$, and denote $X_n(\cpath) = x_n$. Then $X_{n+1}(\cpath)$ is a neighbour $y \nbr x_n$, and 
   \linenopax
    \begin{align}
      \bE[ h(X_{n+1}) \,|\, \gS_n]
      &= \bE[ h(X_{n+1}) \,|\, \gS_n] \sum_{y \nbr x_n} p(x_n,y) \notag \\
      &= \sum_{y \nbr x_n} p(x_n,y) \bE[ h(X_{n+1}) \,|\, \gS_n] \notag \\
      &= \bE\left[ \sum_{y \nbr x_n} p(x_n,y) h(X_{n+1}) \,|\, \gS_n\right] \notag \\
      &= \bE[h(X_n) \,|\, \gS_n] \vstr \notag \\
      &= h(X_n). \vstr 
    \end{align}
    Since $h$ is bounded, this shows $h(X_n)$ is a bounded martingale, whence by Doob's Theorem (cf.~\cite{Doob53}), it converges pointwise $\prob^{(\cond)}$-a.e. on \paths and \eqref{eqn:Vh-is-martingale-limit} makes $V_h$ well-defined $\prob^{(\cond)}$-a.e. on \paths.
    
    ($\leftrightarrow$) We conclude with a proof that these two constructions correspond to inverse operations. 
    If $V$ is a cocycle, we must show that $\lim_{n \to \iy} \Ex_{X_n(\cpath)}[V] = V(\cpath)$. To this end, for $A \ci \Paths$, define the conditioned measure $\prob_A := \frac{\prob^{(\cond)}(A \cap \cdot )}{\prob^{(\cond)}(A)}$, so that $d\prob_A = \frac1{\prob^{(\cond)}(A)} \charfn{A} d\prob^{(\cond)}$. Now for a fixed $\cpath \in \Paths$, let $A_n = \Paths_{(x,X_1(\cpath),\dots,X_n(\cpath))}$ be the cylinder set whose first $n+1$ coordinate agree with \cpath. Applying the measure identity $\lim \gm(A_n) = \bigcap \gm(A_n)$ for nested sets, we obtain $\lim_{n \to \iy} \prob_{A_n} = \gd_\cpath$ as a weak limit of measures. Now
    \linenopax
    \begin{align}\label{eqn:Ex{Xn}[V]=V:derivation1}
      \Ex_{X_n(\cpath)}[V]
      = \int_{\Paths_{X_n(\cpath)}} V(\gx) \, d\prob^{(\cond)}(\gx)
      = \int_{\Paths} V(\gx) \, d\prob_{A_n}(\gx)
      \limas{n} \int_{\Paths} V(\gx) \, d\gd_\cpath
      = V(\cpath).
    \end{align} 
    On the other hand, if $h$ is harmonic, we must show $\Ex_x[V_h]=h$. Then for $V_h(\cpath) := \lim_{n \to \iy} h(X_n(\cpath))$, boundedness allows us to apply the dominated convergence theorem and compute
    \linenopax
    \begin{align}\label{eqn:Ex[Vh]=h}
     \Ex_x[V_h]
     = \int_{\Paths_x} \lim_{n \to \iy} \left(h \comp X_n(\cpath)\right) \,d\prob^{(\cond)}
     = \lim_{n \to \iy} \int_{\Paths_x} h \comp X_n(\cpath) \,d\prob^{(\cond)}.
    \end{align} 
    Now the sequence on the right-hand side of \eqref{eqn:Ex[Vh]=h} is constant by Lemma~\ref{thm:harmonic-as-integral-identity}, so $\Ex_x[V_h] = h(x)$.
  \end{proof}
\end{theorem}

\begin{remark}
  The \fwd direction of the proof of Theorem~\ref{thm:harmonic-functions=cocycles} may also be computed 
  \linenopax
  \begin{align*}
    h_V(x) = \cond(x) \bE_x[V] 
    &= \cond(x) \sum_{y \nbr x} p(x,y) \bE_x[V \,|\, X_1=y]  
    = \sum_{y \nbr x} \cond_{xy} \bE_y[V] 
    = \sum_{y \nbr x} \cond_{xy} h_V(y),
  \end{align*}
  where $\bE_x[V \,|\, X_1=y] = \bE_y[V]$ because the random walk is a Markov process.   See, e.g., \cite[Prop.~9.1]{LevPerWil08}.
\end{remark}

\version{}{
\begin{defn}\label{def:space-of-bounded-paths}
  Let $\gQ_n := \{\cpath \in \Paths \suth |X_k(\cpath)| \leq n, \forall k\}$, where $|x|$ is the shortest-path distance from $o$ to $x$. Then $\gQ := \bigcup_n \gQ_n$ is the collection of all \emph{bounded paths} in \Graph.
\end{defn}

\begin{cor}\label{thm:bdd-paths-have-measure-0}
  If $\Harm \neq 0$, then $\prob^{(\cond)}(\gQ) = 0$.
  \begin{proof}
    \version{}{\marginpar{This doesn't quite work.}}
    Let $h \in \Harm$ be nonconstant. Since $X_k(\cpath) \neq X_{k+1}(\cpath)$ for any \cpath, the limit $\lim_{n \to \iy} h(X_n(\cpath)) = V(\cpath)$ cannot exist for any $\cpath \in \gQ$. The conclusion follows by \eqref{eqn:Vh-is-martingale-limit}. 
  \end{proof}
\end{cor}
}

\section{The forward-harmonic functions}
\label{sec:forward-harmonic-functions}

The current passing through a given edge may be interpreted as the expected value of the number of times that a given unit of charge passes through it. This perspective is studied extensively in the literature; see \cite{DoSn84,Lyons:ProbOnTrees} for excellent treatments. In this case, $p(x,y) = \frac{\cond_{xy}}{\cond(x)}$ helps one construct a current which is harmonic, or dissipation-minimizing. However, that is not what we do here; we are interested in studying current functions whose dissipation is finite but not necessarily minimal. In Theorem~\ref{thm:Paths-is-nonempty}, we show that the experiment always induces a ``downstream'' current flow between the selected two points; that is, a path along which the potential is strictly decreasing.

These probability notions will play a key role in our solution to questions about \emph{activity}; cf.~Definition~\ref{def:activity}. We use the forward path measure again in our representation formula (Theorem~\ref{thm:forward-harmonic-functions=cocycles}) for the class of forward-harmonic functions on \Graph. The corresponding Markov process is dynamically adapted to the network (and the charge on it). This representation is dynamic and nonisotropic, which sets it apart from other related representation formulas in the literature.

\subsection{Activity of a current and the probability of a path}
\label{sec:current-paths}
Given a (fixed) current, we are interested in computing ``how much of the current'' takes any specified path from $x$ to some other (possibly distant) vertex $y$. This will allow us to answer certain existence questions (see Theorem~\ref{thm:Paths-is-nonempty}) and provides the basis for the study of the forward-harmonic functions studied in \S\ref{sec:forward-harmonic-functions}. Note that, in contrast to \eqref{eqn:def:path-probability-recall}, the probabilistic interpretation given in Definition~\ref{def:path-probability} (and the discussion preceding it) does not make any reference to \cond. In this section we follow \cite{Pow76b} closely.

\begin{defn}\label{def:activity}
  The \emph{divergence} of a current $\curr:\edges \to \bR$ is the function on $x \in \verts$ defined by
  \begin{equation}\label{eqn:def:absolute-activity}
    \act_\curr(x) :=
    \begin{cases}
      \frac12 \sum_{y \nbr x} |\curr(x,y)|, &x \neq \ga,\gw \\
      1, &x=\ga,\gw.
    \end{cases}
  \end{equation}
  which describes the total ``current traffic'' passing through $x \in \verts$.
  Thus, \act is an operator mapping functions on \edges to functions on \verts; see \S\ref{sec:divergence-operator} for details.
\end{defn}

For convenience, we restate Definition~\ref{def:current-path-preview}:

\begin{defn}\label{def:current-path}
  Let $v:\verts \to \bR$ be given and suppose we fix \ga and \gw for which $v(\ga) > v(\gw)$. Then a \emph{current path} \cpath (or simply, a \emph{path}) is an edge path from \ga to \gw with the extra stipulation that $v(x_k) < v(x_{k-1})$ for each $k=1,2,\dots,n$.
  Denote the set of all current paths by $\Paths = \Paths_{\ga,\gw}$ (dependence on the initial and terminal vertices is suppressed when context precludes confusion). Also, define $\Paths_{\ga,\gw}(x,y)$ to be the subset of current paths from \ga to \gw which pass through the edge $(x,y) \in \edges$.
  \glossary{name={\Paths},description={set of paths},sort=G,format=textbf}
\end{defn}

Suppose we fix a source \ga and sink \gw and consider a single current path \cpath from \ga to \gw. With $\act_\curr$ defined as in \eqref{eqn:def:absolute-activity}, one
can consider $\frac{\curr(x,y)}{\act_\curr(x)}$ as the probability that a unit of charge at $x$ will pass to a ``downstream'' neighbour $y$. Note that $\curr(x,y) > 0$ and $\act_\curr \neq 0$, since we are considering an edge of our path \cpath. This allows us to define a probability measure on the path space $\Paths_{\ga,\gw}$.

\begin{defn}\label{def:path-probability}
  If $\cpath \in \Paths$ follows the vertex path $(\ga = x_0, x_1, x_2, \dots, x_n = \gw)$, the define the probability of \cpath by
  \begin{equation}\label{eqn:def:path-probability}
    \prob(\cpath) := \prod_{k=1}^{n} \frac{\curr(x_{k-1},x_k)}{\act_\curr(x_{k-1})}.
  \end{equation}
  This quantity gives the probability that a unit of charge at \ga will pass to \gw by traversing the path \cpath.
\end{defn}

\version{}{
The following result fills another hole in the proof of $\dist_\Lap(x,y) = \dist_\diss(x,y)$ (recall \eqref{eqn:def:R(x,y)-Lap} and \eqref{eqn:def:R(x,y)-diss}) in \cite{Pow76b}.

\begin{cor}\label{thm:no-infinite-current-paths-for-minimal-currents}
  \marginpar{Is this even true?}
  If \curr minimizes \diss on $\Flo(\ga,\gw)$, then all current paths flow from \ga to \gw.
  \begin{proof}
    If $\curr \in \Flo(\ga,\gw)$ is minimal, then Theorem~\ref{thm:minimal-flows-are-induced-by-minimizers} implies \curr is induced by a potential $v$ which minimizes \energy over $\Pot(\ga,\gw)$. Since orthogonality gives $\energy(v) = \energy(\Pfin v) + \energy(\Phar v)$, such a minimal $v$ must lie in \Fin. In fact, $v = f_\ga - f_\gw$, where $f_x = \Pfin v_x$. 
  \end{proof}
\end{cor}
}

\subsection{Forward-harmonic transfer operator}
\label{sec:forward-harmonic-transfer-operator}

In this section we consider the functions $h : \verts \to \bC$ which are \emph{forward-harmonic}, that is, functions which are harmonic with respect to a current \curr. We make the standing assumption that the network is transient; this guarantees the existence of a monopole at every vertex, and the induced current will be a unit flow to infinity; cf.~Corollary~\ref{thm:transient-here-implies-transient-everywhere}.

We orient the edges by a fixed unit current flow \curr to infinity, as in Definition~\ref{def:orientation}. The forward-harmonic functions functions are fixed points of a transfer operator induced by the flow which gives the value of $h$ at one vertex as a convex combination of its values at its downstream neighbours.

The main idea is to construct a measure on the space of paths beginning at vertex $x \in \verts$, and then use this measure to define forward-harmonic functions. In fact, we are able to produce all forward-harmonic functions from the class of functions which satisfy a certain cocycle condition, see Definition~\ref{def:cocycle}.

In Theorem~\ref{thm:forward-harmonic-functions=cocycles} we give an integral representation for the harmonic functions, and in Corollary~\ref{thm:current-with-sink=>trivial-harmonic-functions} we show that if \curr has a universal sink, then the only forward harmonic functions are the constants.

\begin{remark}[A current induces a direction on the resistance network.]
  If we fix a minimal current $\curr = \Pdrp \curr$, the flow gives a strict partial order on \verts and the flags in the resulting poset are the induced current paths. Thus we say $x \prec y$ iff $x$ is upstream from $y$, that is, iff there exists a current path from $x$ to $y$ in the sense of Definition~\ref{def:current-path}. Since \curr is minimal, $x \prec y$ implies $y \nprec x$ and $x \nprec x$. Transitivity is immediate upon considering the concatenation of two finite paths.
\end{remark}

\begin{defn}\label{def:set-of-current-paths}
  Given a fixed minimal current $\curr = \Pdrp \curr$, we denote the \emph{set of all
  current paths} in the resistance network $(\Graph,\cond)$ by
  \begin{equation}\label{eqn:def:set-of-current-paths}
    \paths_\curr := \{\cpath = (x_0,x_1, \dots \suth (x_i,x_{i+1}) \in \edges, x_i \prec x_{i+1}\}.
  \end{equation}
  For $n = 1, 2, \dots$ , we denote the set of all current paths of length n by
  \begin{equation}\label{eqn:def:set-of-current-paths-length-n}
    \paths^{(n)}_\curr := \{\cpath = (x_0, x_1, \dots, x_n) \suth (x_i,x_{i+1}) \in \edges, x_i \prec x_{i+1}\},
  \end{equation}
  and denote the collection of paths starting at $x$ by $\paths_\curr(x) := \{\cpath \in \paths_\curr \suth x_0=x\}$, and likewise for $\paths^{(n)}_\curr(x)$.
\end{defn}

Here, the orientation is determined by \curr, and if $\curr(x,y) = 0$ for some $(x,y) \in \edges$, then this edge will not appear in any current path, and for all practical purposes it may be considered as having been removed from \Graph for the moment.

\begin{defn}\label{def:forward-neighbours}
  When a minimal current $\curr = \Pdrp \curr$ is fixed, the set of \emph{forward neighbours} of $x \in \verts$ is
  \begin{equation}\label{eqn:def:forward-neighbours}
    \fnbrs_\curr(x) := \{y \in \verts \suth x \prec y, x \nbr y\}.
  \end{equation}
\end{defn}

\begin{defn}\label{def:forward-Laplacian}
  If $v:\verts \to \bR$, define the \emph{forward Laplacian} of $v$ by
  \begin{equation}\label{eqn:def:forward-Laplacian}
    \fLap v(x) := \sum_{y \in \fnbrs(x)} \cond_{xy}(v(x) - v(y))
  \end{equation}
  A function $h$ is forward-harmonic iff $\fLap h = 0$.
\end{defn}

For Definitions~\ref{def:set-of-current-paths}--\ref{def:forward-Laplacian}, the dependence on \curr may be suppressed when context precludes confusion. 

\begin{theorem}\label{thm:Px-exists}
  For $\curr \in H_\diss$ and $x \in \verts$, there is a unique measure $\prob_x$ defined on $\paths_\curr(x)$ which satisfies
  \begin{equation}\label{eqn:Px-consistency}
    \prob_x[f] = \prob_x^{(n)}[f], \qq f \in \sA_n.
  \end{equation}
  \begin{proof}
    We only need to check Kolmogorov's consistency condition \eqref{eqn:Kolmogorov-consistency}; see \cite{Jor06,Kol56}. For $n < m$, consider $\sA_n \ci \sA_m$ by assuming that $f$ depends only on the first $n$ edges of \cpath. (For brevity, we denote a function on $n$ edges as a function on $n+1$ vertices.) Then
    \linenopax
    \begin{align*}
      \prob_x^{(m)}[f]
      &= \int_{\paths_\curr(x)} f(\cpath) \, d\prob_x^{(m)}(\cpath) \\
      &= \int_{\paths_\curr(x)} f(x_0,x_1,x_2,\dots,x_n) \, d\prob_x^{(m)}(\cpath) \\
      &= \int_{\paths_\curr(x)} f(x_0,x_1,x_2,\dots,x_n) \, d\prob_x^{(n)}(\cpath) \\
      &= \prob_x^{(n)} [f]. \vstr[2.5]
      \qedhere
    \end{align*}
  \end{proof}
\end{theorem}

\subsection{A boundary representation for the forward-harmonic functions} We now show that the forward-harmonic functions are in bijective correspondence with the cocycles, when defined as follows.

\begin{defn}\label{def:cocycle}
  A \emph{cocycle} is a function $f:\paths_\curr \to \bR$ which is compatible with the probabilities on current paths in the sense that it satisfies
  \begin{equation}\label{eqn:def:cocycle-condition}
    f(\cpath) = f(x_0,x_1,x_2,x_3\dots) = \frac{\cond_{x_0x_1} \act_\curr(x_0)}{\cond^+(x_0) \curr(x_0,x_1)} f(x_1,x_2,x_3\dots),
  \end{equation}
  whenever $\cpath=(x_0,x_1,x_2,x_3\dots) \in \paths_\curr$ is a current path as in Definition~\ref{def:set-of-current-paths}, and $(x_0,x_1)$ is the first edge in \cpath. Also, $\cond^+(x) := \sum_{y \in \fnbrs(x)} \cond_{xy}$ is the sum of conductances of edges \emph{leaving} $x$. If the operator $m$ is given by multiplication by
  \begin{equation}\label{eqn:def:m-multiplication-operator}
    m(x,y) = \frac{\cond_{xy} \act_\curr(x)}{\cond^+(x) \curr(x,y)},
  \end{equation}
  and \gs denotes the shift operator, the cocycle condition can be rewritten $f = mf\gs$. Using $e_k = (x_{k-1}, x_k)$ to denote the edges, this gives
  \begin{equation}\label{eqn:f-as-infinite-product}
    f(e_1, e_2, \dots)
    = m(e_1) \dots m(e_n) f(e_{n+1}, e_{n+2}, \dots)
    = \prod_{k=1}^\iy m(e_k).
  \end{equation}
\end{defn}

\begin{defn}\label{def:forward-transfer-operator}
  Define the \emph{forward transfer operator} $\Trans_\curr$ induced by \curr by
  \begin{equation}\label{eqn:def:forward-transfer-operator}
    (\Trans_\curr f)(x) := \frac1{\cond^+(x)} \sum_{y \in \fnbrs(x)} \cond_{xy} f(y).
  \end{equation}
\end{defn}

\begin{lemma}\label{thm:forward-cocycles-are-harmonic}
  If $f:\paths(x) \to \bR$ is a cocycle and one defines $h_f(x) := \prob_x[f]$, it
  follows that $h_f(x)$ is a fixed point of the forward transfer operator $\Trans_\curr$.
  \begin{proof}
    With $h_f$ so defined, we conflate the linear functional $\prob_x$ with the measure associated to it via Riesz's Theorem and compute
    \linenopax
    \begin{align*}
      (\Trans_\curr h_f)(x)
      &= \frac1{\cond^+(x)} \sum_{y \in \fnbrs(x)} \cond_{xy} \prob_y[f]
        &&\text{def } \Trans, h_f \\
      &= \sum_{y \in \fnbrs(x)} \frac{\cond_{xy}}{C^+(x)}
        \int_{\paths_\curr(y)} f(\cpath) \,d\prob_y(\cpath)
        &&\prob_y \text{ as a measure} \\
      &= \sum_{y \in \fnbrs(x)} \int_{\paths_\curr(x)}
        \frac{\cond_{xy}}{\cond^+(x)} f(\gs\cpath) \,d\prob_x(\cpath)
        &&\text{change of vars} \\
      &= \sum_{y \in \fnbrs(x)} \frac{\curr(x,y)}{\act_\curr(x)} \int_{\paths_\curr(x)}
        \frac{\cond_{xy}}{\cond^+(x)} \frac{\act_\curr(x)}{\curr(x,y)} f(\gs\cpath) \,d\prob_x(\cpath)
        &&\text{just algebra} \\
      &= \sum_{y \in \fnbrs(x)} \frac{\curr(x,y)}{\act_\curr(x)}
        \int_{\paths_\curr(x)} f(\cpath) \,d\prob_x(\cpath)
        &&\text{by \eqref{eqn:def:cocycle-condition}} \\
      &= \sum_{y \in \fnbrs(x)} \frac{\curr(x,y)}{\act_\curr(x)} \prob_x(f)
        &&\prob_x \text{ as a functional} \\
      &= \prob_x(f) \vstr[2.5]
        &&\textstyle{\sum \tfrac{\curr(x,y)}{\act_\curr(x)} = 1}
    \end{align*}
    To justify the change of variables, note that if \cpath is a path starting at $x$ whose  first edge is $(x,y)$, then \gs\cpath is a path starting at $y$. Moreover, since $y$ is a downstream neighbour of $x$, every path \cpath starting at $y$ corresponds to exactly one path starting at $x$, namely, $((x,y), \cpath)$.
  \end{proof}
\end{lemma}

\begin{theorem}\label{thm:forward-harmonic-functions=cocycles}
  The forward-harmonic functions are in bijective correspondence with the cocycles. More precisely, if $f$ is a cocycle, then
  \begin{equation}\label{eqn:hf-is-P[f]}
    h_f(x) := \prob_x[f]
  \end{equation}
  is harmonic. Conversely, if $h$ is harmonic, then
  \begin{equation}\label{eqn:fh-is-martingale-limit}
    f_h(\cpath) := lim_{n \to \iy} h(X_n(\cpath)), \qq \cpath \in \paths_x
  \end{equation}
  is a cocycle, where $X_n(\cpath)$ is the \nth vertex from $x$ along the path \cpath.
  \begin{proof}
    \fwd Let $f$ be a cocycle and define $h_f$ as in \eqref{eqn:hf-is-P[f]} with $C^+(x)$ as in Definition~\ref{def:forward-transfer-operator}, compute
    \linenopax
    \begin{align*}
      \fLap h_f(x)
      &= \sum_{y \in \fnbrs(x)} \cond_{xy} (\prob_x[f] - \prob_y[f]) \\
      &= \prob_x[f] \sum_{y \in \fnbrs(x)} \cond_{xy}
        - \sum_{y \in \fnbrs(x)} \cond_{xy} \prob_y[f] \\
      &= \prob_x[f] C^+(x) - \sum_{y \in \fnbrs(x)} \cond_{xy} \prob_y[f],
    \end{align*}
    which is 0 by Lemma~\ref{thm:forward-cocycles-are-harmonic}.

    \bwd Let $h$ satisfy $\fLap h = 0$. Observe that $X_n$ is a Markov chain with transition probability $\prob_x$ at $x$. The above computations show that $h$ is then a fixed point of $\Trans_\curr$, and hence $h(X_n)$ is a bounded martingale. By Doob's Theorem (cf.~\cite{Doob53}), it converges pointwise $\prob_x$-ae on \paths and \eqref{eqn:fh-is-martingale-limit} makes $f_h$ well-defined.
    One can see that $f_h$ is a cocycle by the same arguments as in the proof of \cite[Thm.~2.7.1]{Jor06}.
  \end{proof}
\end{theorem}

\begin{cor}\label{thm:current-with-sink=>trivial-harmonic-functions}
\version{}{\marginpar{If the support of $\prob_x$ is not all of \edges, then can two distinct cocycles give the same harmonic function?}}
  If \curr has a universal sink, in other words, if all current paths \cpath end at some common point \gw, then the only forward-harmonic function is the zero function.
  \begin{proof}
    Every harmonic function comes from a cocycle, which in turn comes from a harmonic function as a martingale limit, by the previous theorem. However, formula \eqref{eqn:fh-is-martingale-limit} yields
    \linenopax
    \begin{align*}
      f_h(\cpath) := \lim_{n \to \iy} h(X_n(\cpath)) = h(\gw), \qq \forall \cpath \in \paths.
    \end{align*}
    Thus every cocycle is constant, and hence \eqref{eqn:def:cocycle-condition} implies $f_h \equiv 0$. Then \eqref{eqn:hf-is-P[f]} gives $h \equiv 0$.
  \end{proof}
\end{cor}

\version{}{\marginpar{Actually, maybe \eqref{eqn:def:cocycle-condition} doesn't imply 0 ... without further restrictions on the network or something.}}

\section{Remarks and references}
\label{sec:Remarks-and-References-probab-interp}

The probability literature is among the largest of all the subdisciplines of mathematics, and so the following list of suggested references barely begins to scratch the surface: \cite{Ign08, AHKL06, Aik05, BrWo05, ChKi04, BuZo09, IKW09, Si09, YiZh05}.

Of these references, some are more specialized. However for prerequisite material (if needed), the reader may find the following sources especially relevant: \cite{Lyons:ProbOnTrees, LevPerWil08, AlFi09, Peres99}.

%% file: examples.tex

\chapter{Examples and applications}
\label{sec:examples}

\headerquote{The art of doing mathematics consists in finding that special case which contains all the germs of generality}{---~D.~Hilbert}

\section{Finite graphs}
\label{sec:finite-graphs}

\subsection{Elementary examples}
\label{sec:elementary-examples}

\begin{exm}\label{exm:linear-network}
  Consider a ``linear'' \ERN consisting of several resistors connected in series with resistances $\ohm_i = \cond_i^{-1}$ as indicated:
  \[\xymatrix{
      \ga = x_0 \ar[r]^{\hstr[3] \ohm_1}
      & x_1 \ar[r]^{\ohm_2}
      & x_2 \ar[r]^{\ohm_3}
      & x_3 \ar[r]^{\ohm_4}
      & \dots \ar[r] ^{\ohm_n \hstr[2]}
      & x_n =\gw
    }\]
  Construct a dipole $v \in \Pot(\ga,\gw)$ on this network as follows. Let $v(x_0) = V$ be fixed. Then determine $v(x_1)$ via \eqref{eqn:def:dipole}:
  \linenopax
  \begin{align*}
    \Lap v(x_0) = \tfrac1{\ohm_1}(V-v(x_1))=1 &\implies v(x_1)=V-\ohm_1,\\
    \Lap v(x_1) = \sum_{k=1}^2 \tfrac1{\ohm_k}(v(x_{k-1})-v(x_k))=0 &\implies v(x_2)=V-\ohm_1-\ohm_2,
  \end{align*}
  and so forth. Three things to notice about this extremely elementary example are (i) $v$ is fixed by its value at one point and any other dipole on this graph can differ only by a constant, (ii) we recover the basic fact of electrical theory that the voltage drop across resistors in series is just the sum of the resistances, and (iii) all current flows are induced (this is not true of more general graphs).

  \version{}{\marginpar{What does this construction illustrate?}}
  Consider the basis $\{e_0,e_1,\dots,e_N\}$, where $e_k = \gd_{x_k}$, the unit Dirac mass at $k$. The Laplace operator for this model has the matrix
  \begin{equation}\label{exm:finite-linear-model}
    \Lap =
    \left[\begin{array}{rrrrrr}
      1 & -1 & 0 & \dots && 0 \\
      -1 & 2 & -1 & \dots && 0 \\
      0 & -1 & 2 & \dots && 0 \\
      \vdots & && \ddots \\
      0 & &\dots & -1 & 2 & -1 \\
      0 & &\dots & 0 & -1 & 1
    \end{array}\right].
  \end{equation}
  One may obtain a unitary representation on $\ell^2(\bZ_N)$ by using the diagonal matrix $U(\gz) = \diag(1,\gz,\gz^2,\dots,\gz^N)$, where $\gz := e^{2\gp\ii/(N+1)}$ is a primitive \nth[(N+1)] root of 1, so that $\gz^{-1}=\bar \gz$. It is easy to check that for any matrix $M \in \sM_{N+1}(\bC)$, one has $[U(\gz) M U(\gz)^\ad]_{j,k} = \gz^{j-k} [M]_{j,k}$. Then define
  \linenopax
  \begin{align*}\label{eqn:unitary-zeta-Lap-defined}
    \Lap(\gz)
    :=& U(\gz) \Lap U(\gz)^\ad,
  \end{align*}
  and see that $\Lap(\gz) = \Cond - U(\gz) \Trans U(\gz)^\ad$. It is clear that $\spec \Lap = \spec\Lap(\gz)$, because
  \linenopax
  \begin{align*}
    \Lap v = \gl v
    \qq\iff\qq
    \Lap(\gz) [U(\gz) v] = \gl [U(\gz) v].
  \end{align*}
  Decompose the transfer operator into the sum of two shifts, so that $\Trans = M_+ + M_-$, where $M_+$ has ones below the main diagonal and zeros elsewhere, and $M_-$ has ones above the diagonal and zeros elsewhere. Then we have $U(\gz) M_+ U(\gz)^\ad = \gz M_+$ and $U(\gz) M_- U(\gz)^\ad = \bar\gz M_-$ and $M_- = M_+^\ad$. By induction, the characteristic polynomial can be written
  \linenopax
  \begin{align*}
    p_n(x) = \det(x\id - \Trans_n) = x p_{n-1}(x) - p_{n-2}(x),
  \end{align*}
  with $p_1=x, p_2=x^2-1,p_3=x^3-2x,p_4=x^4-3x^2+1$, and corresponding Perron-Frobenius eigenvalues $\gl_1=0, \gl_2 = 1, \gl_3 = \sqrt2, \gl_4 = \phi = \frac12(1+\sqrt5)$.
  \linenopax
  \begin{align*}
    \spec \Lap_2 = \{0,1\}, \q
    \spec \Lap_3 = \{0,1,3\}, \q
    \spec \Lap_4 = \{0,2,2\pm\sqrt2\}.
  \end{align*}
\end{exm}

\begin{exm}\label{exm:one-small-cycle-network}
  The correspondence $\Pot(\ga,\gw) \to \Flo(\ga,\gw)$ described in Lemma~\ref{thm:potential-induces-current} is not bijective, i.e., the converse to the theorem is false, as can be seen from the following example. Consider the following electrical  network with resistances $\ohm_i = \cond_i^{-1}$.
  \linenopax
    \begin{align*}
    \xymatrix{
      *+[l]{\ga = x_0} \ar[r]^{\ohm_1} \ar[d]^{\ohm_2} & x_1 \ar[d]^{\ohm_3} \\
      x_2 \ar[r]^{\ohm_4} & *+[r]{x_3 = \gw}
    }
  \end{align*}
  One can verify that the following gives a current flow $\curr = \curr_t$ on the graph for any $t \in [0,1]$:
  \linenopax
    \begin{align*}
    \xymatrix{
      x_0 \ar[r]_t \ar[d]_{1-t} & x_1 \ar[d]_t \\
      x_2 \ar[r]_{1-t} & x_3
    }
  \end{align*}
  In fact, there are many flows on this network; let $\charfn{\cycle}$ be the characteristic function of the cycle $\cycle = (x_0, x_1, x_3, x_2)$
  \[\cycle =
    \begin{gathered}
    \xymatrix{
      x_0 \ar[r] & x_1 \ar[d] \\
      x_2 \ar[u] & x_3 \ar[l]
    }
    \end{gathered}
  \]
  so that $\charfn{\cycle}(e) = 1$ for each $e \in \{e_1=(x_0, x_1), e_2=(x_1, x_2), e_3=(x_2, x_3), e_4=(x_3, x_0)\}$. Then $\curr_t + \ge \charfn{\cycle}$ will be a flow for any $\ge \in \bR$. (Although this formulation seems more awkward than simply allowing $t$ to take any value in \bR, it is easier to work with characteristic functions of cycles when there are many cycles in the network.) However, there will be only one value of $t$ and \ge for which the above flow corresponds to a potential function, and that potential function is the following:
  \linenopax
  \begin{align*}
    \xymatrix{
      V \ar[r] \ar[d]
        & V - \frac{\ohm_2 + \ohm_4}{\sum_{k=1}^4 \ohm_k}\ohm_1 \ar[d] \\
      V - \frac{\ohm_1 + \ohm_3}{\sum_{k=1}^4 \ohm_k}\ohm_2 \ar[r]
        & V - \frac{(\ohm_1 + \ohm_3)(\ohm_2 + \ohm_4)}{\sum_{k=1}^4 \ohm_k}
    }
  \end{align*}
  This is the potential function which ``balances'' the flow around both sides of the square; it can be computed as in the previous example. These ideas are given formally in Theorem~\ref{thm:minimal-flows-are-induced-by-minimizers}.
\end{exm}

\begin{exm}[Finite cyclic model]\label{exm:finite-cyclic-model}
  In this case, let $\Graph_N$ have vertices given by $x_k = e^{2\gp\ii k/N}$ for $k=1,2,\dots,N$, with edges connecting each vertex to its two nearest neighbours. For example, when $N=9$,\\
  \begin{center}
    $\Graph_9$ \qq \xy \xygraph{[rrr] !P9"C"{~={42} ~*{x_{\xypolynode}}} }\endxy $=x_0$
  \end{center}
  In this case, using the same basis $\{e_0,e_1,\dots,e_N\}$, as in Example~\ref{exm:linear-network}, the Laplace operator for this model has the matrix
  \begin{equation}\label{eqn:exm:finite-cyclic-model}
    \Lap =
    \left[\begin{array}{rrrrrrrr}
      2 & -1 & 0 & \dots  & 0 & -1 \\
      -1 & 2 & -1 & \dots  & 0 & 0 \\
      0 & -1 & 2 & \dots  & 0 & 0 \\
      \vdots & && \ddots \\
      0 & 0 & 0 & \dots  & 2 & -1 \\
      -1 & 0 & 0 & \dots  & -1 & 2
    \end{array}\right].
  \end{equation}
  The Fourier transform $\sF$ is a spectral transform of \Lap that shows it to be unitarily equivalent to multiplication by $4\sin^2\left(\frac{\gp k}{N}\right)$.

  \begin{lemma}\label{thm:Lap_N-is-4sin2-etc}
    For $S := \sF^\ad \Lap \sF$ and $v \in \HE$, one has
    \linenopax
    \begin{align*}
      (Sv)(z) = (2-z-z^{-1})v(z), \qq z=\ga^k, k=0,1,2,\dots,N.
    \end{align*}
    \begin{proof}
      Denote $v_k=v(k)$ and consider the Fourier transform $\sF:\{v_k\} \mapsto \sum_{k \in \bZ_N} v_k z^k$: 
      \linenopax
      \begin{align*}
        \sF(\Lap v)(z) 
        &= \sum_{k \in \bZ_N} (2v_k - v_{k-1} - v_{k+1}) z^k \\
        &= 2\sum_{k \in \bZ_N} v_k z^k - \sum_{k \in \bZ_N} v_{k-1} z^k - \sum_{k \in \bZ_N} v_{k+1} z^k \\
        &= (2 - z - z^{-1}) \sum_{k \in \bZ_N} v_k z^k.
      \end{align*}
      This shows $\sF(\Lap v)(z) = (2-z-z^{-1})\sF(v)(z)$, so that $S$ is multiplication by $(2-z-z^{-1})$. 
    \end{proof}
  \end{lemma}

  \pgap

  In this case, $\Lap = \id - \Trans$ with $\Cond = 2\id$, so that \Cond and \Trans commute. Additionally, \Trans is the sum of two shifts and so corresponds to multiplication by $z+z^{-1} = 2\cos \gq$ on $\ell^2(\verts_N)$, where $\gq = \frac{2\gp k}N$. Consequently, $\spec \Trans_N = 2 \cos \left(\frac{2\gp k}N\right)$ and the Perron-Frobenius eigenvalue of $\Trans_N$ is $\gl_{PF}=2$, which occurs for $k=0$ and has eigenfunction $v_{PF} = [1,1,\dots,1]$. Observe that $\Lap_N$ commutes with the cyclic shift. The eigenfunctions of the shift are $v = [1,\gl,\gl^2,\dots,\gl^{N-1}]$, where $\gl \in \verts$, and hence these are the other eigenfunctions of $\Trans_N$.

  \begin{prop}
    The spectrum of $\Lap_N$ is given by
    \linenopax
    \begin{align*}
      \spec \Lap_N = \{2\left(1-\cos \left(\tfrac{2\gp k}N\right) \right)\suth k \in \bZ_N\}.
    \end{align*}
    \begin{proof}
      Let $U$ be the cyclic shift in the positive direction, i.e., $U$ has the matrix
      \linenopax
      \begin{align*}
        U =
        \left[\begin{array}{cccccc}
          0 & &&&1\\ 1 & 0 \\ & 1 & 0 \\ &&\ddots& \ddots \\ &&&1 &0
        \end{array}\right].
      \end{align*}
      The $U$ and $\Lap_N$ commute. For $n \in \bZ_n$, the Fourier transform is
      \linenopax
      \begin{align*}
        \hat v(n) = \frac1N \sum_{z \in \verts_N} z^{-n} v(z),
        \q\text{so}\q
        x_k \mapsto v(n) = \sum_{n \in \bZ_N} z^n \hat v(n).
        &\qedhere
      \end{align*}
      \version{}{\marginpar{Some comment here about a Hadamard matrix?}}
    \end{proof}
  \end{prop}
  \version{}{\marginpar{At this point, I can no longer follow the notes. What is the main point? The Fourier transform diagonalizes \Lap, but what is the application here?}}
  For $k \in \{1,2,\dots,N-1\}$, one can find $v \in \Pot(k,0)$ as follows: ``ground'' the graph with $v(0)=0$ and consider
  \linenopax
    \begin{align*}
    \xymatrix{
    \dots \ar[r]_\ga & x_1 \ar[r]_\ga & x_0=x_N & x_{N-1} \ar[l]^{1-\ga} & \dots \ar[l]^{1-\ga} & x_k \ar[l]^{1-\ga} \ar[r]_\ga & \dots
    }
  \end{align*}
  The cycle condition (the net drop of voltage around any closed cycle must be 0) yields $v(k)-v(0) = k\ga = (N-k)(1-\ga)$, and hence $\ga = \tfrac{N-k}N$. This gives $v(k) = \tfrac{k(N-k)}{N}$ and we have
  \linenopax
  \begin{align*}
    \Lap v(j) = \gd_k - \gd_0 =
    \begin{cases}
      -\tfrac{k}{N} - \tfrac{N-k}{N} = -1, &j=0 \\
      \tfrac{k}{N} - \tfrac{k}{N} = 0, &0 < j < k \\
      \tfrac{k}{N} + \tfrac{N-k}{N} = 1, &j=k \\
      \tfrac{N-k}{N} - \tfrac{N-k}{N} = 0, &k < j \leq N-1.
    \end{cases}
  \end{align*}
  This additionally shows that if the shortest path from $x$ to $y$ has length $k$, then $R(x,y) = k\frac{N-k}N < k$. Of course, there is an easier way to get $R(x,y)$. Since there are only two paths $\cpath_1,\cpath_2$ from $x$ to $y$, the laws for resistors connected in serial and parallel indicate that the entire network can be replaced by a single edge $(x,y)$ with resistance $(\ohm(\cpath_1)^{-1} + \ohm(\cpath_1)^{-1})^{-1}$. In the case of constant resistance $\ohm \equiv 1$, this becomes
  \linenopax
  \begin{align*}
    (\ohm(\cpath_1)^{-1} + \ohm(\cpath_1)^{-1})^{-1}
     = \left(\frac1k + \frac1{N-k}\right)^{-1}
     = \left(\frac{(N-k)+k}{k(N-k)}\right)^{-1}
     = k\frac{N-k}N.
  \end{align*}
\end{exm}

  \begin{figure}
    \centering
    \scalebox{0.90}{\includegraphics{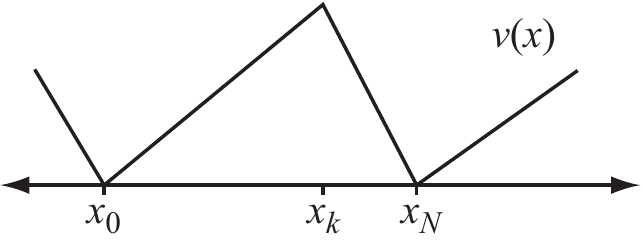}}
    \caption{\captionsize The solution $v$ as represented on \bR.}
    \label{fig:periodic-potential}
  \end{figure}

\begin{exm}\label{exm:edge-deletion-network}
  Next, it is illuminating to see how things change when edges are removed from the network. 
  \version{}{\marginpar{Add proof}
\begin{lemma}[Powers' inequality]
  \label{thm:Powers-inequality}
  Suppose \Graph is finite and $H \ci \Graph$ is obtained by removing the edges $\{e_i\}_{i=1}^n$. Then 
  \begin{align}\label{eqn:Powers-inequality}
    R_G(x,y) \leq R_H(x,y) \leq \frac{1}{\ge^2} R_G(x,y), 
  \end{align}
  for $\ge := 1- \sum_{i=1}^n |\curr(e_i)|$.
\end{lemma}
}
Consider the following example, where $\cond = \one$ and the currents are as indicated, and the second network is obtained from the first by deleting two edges, as indicated:
  \linenopax
  \begin{align*}
    \xymatrix{
      &*+[l]{\ga = x_0} \ar[dl]^{\frac{2}{11}} \ar[dd]^{\frac{4}{11}} \ar[r]^{\frac{5}{11}} & x_1 \ar[dd]^{\frac{5}{11}} \\
       x_2 \ar[dr]^{\frac{2}{11}} \\
      & x_3 \ar[r]^{\frac{6}{11}} & *+[r]{x_4 = \gw} \\
      &\hstr[5]\curr_0
    }
    \hstr[5]
    \xymatrix{
      & *+[l]{\ga = x_0} \ar[dd]^{\frac12} \ar[r]^{\frac12} & x_1 \ar[dd]^{\frac12} \\
      \\
      & x_3 \ar[r]^{\frac12} & *+[r]{x_4 = \gw} \\
      &\hstr[5]\curr_1
    }
  \end{align*}
  The dissipations are
  $R_0(\ga,\gw) = \diss(\curr_0) = \frac{10}{11}$ and $R_1(\ga,\gw) = \diss(\curr_1) = 1$.
  The set of paths from \ga to \gw that don't pass through the deleted edges contains
  only $\cpath_1 = (\ga,x_1,\gw)$ and $\cpath_2 = (\ga,x_3,\gw)$. Then
  \linenopax
  \begin{align*}
    \ge = \sum_{\cpath \in Q \less W} \prob(\cpath) = \prob(\cpath_1) + \prob(\cpath_2) = \frac{4}{11} + \frac{5}{11} = \frac9{11},
  \end{align*}
  \version{}{\marginpar{Must include!}}
  and Powers' inequality gives
  \linenopax
  \begin{align*}
    R_0(\ga,\gw) &\leq R_1(\ga,\gw) \leq \ge^{-2} R_0(\ga,\gw) \\
    \frac{10}{11} &\leq 1 \leq \frac{121}{81} \cdot \frac{10}{11} = \frac{110}{81}.
  \end{align*}
\end{exm}

\section{Infinite graphs}
\label{sec:infinite-graphs}

\begin{exm}\label{exm:edge-boundary}
  Define the projections
\linenopax
\begin{equation}\label{eqn:exhaustion-restriction-projections}
  B_n = \charfn{\Graph_n}, \qq B_n^\perp = \charfn{\Graph \less \Graph_n} = \charfn{\Graph_n^\complement}.
\end{equation}
Let us denote the \emph{edge boundary} between $\Graph_n$ and $\Graph_{n+1}$ by
\begin{equation}\label{eqn:nth-edge-boundary}
  Edge_n := \{e=(x,y)\in\edges \suth y \in \verts_n, x \in \verts \less \verts_n\}.
\end{equation}
We now consider the behavior of $\| B_n^\perp \Lap B_n\|$, where the norm is with respect to operators $\ell^2(\cond) \to \ell^2(\cond)$.
\end{exm}

\begin{lemma}\label{thm:edge-boundary-bound}
  For $C_n := \sup\{\sum_{y \nbr x} \cond_{xy}^2 \suth (x,y) \in Edge_n\}$,
  \begin{equation}\label{eqn:lem:edge-boundary-bound}
    \|B_n^\perp \Lap B_n\| \leq C_n.
  \end{equation}
  \begin{proof}
    Let $v \in \ell^2(\cond)$ and $x \in (\verts_n)^\complement$.
    Then
    \linenopax
    \begin{align*}
      (B_n^\perp \Lap B_n v)(x)
      &= \charfn{\Graph^\complement}(x) \sum_{y \nbr x} \cond_{xy} \left(\cancel{\charfn{\Graph_n} v(x)} - \charfn{\Graph_n} v(y)\right) \\
      &= - \negsp[15] \sum_{(x,y) \in Edge_n} \negsp[10] \cond_{xy} v(y), \q x \in (\verts_n)^\complement.
    \end{align*}
    Now summing over all $x \in (\verts_n)^\complement$, and hence over all edges in $Edge_n$,
    \linenopax
    \begin{align*}
      \|B_n^\perp \Lap B_n v\|_\cond^2
      &= \left|\sum_{(x,y) \in Edge_n} \negsp[20] \cond_{xy} v(y) \right|^2 \negsp[10]
      \leq \sum_{y \in \verts_n} |v(y)|^2 \negsp[10]
         \sum_{\substack{x \in (\verts_n)^\complement \\ y \nbr x}} |\cond_{xy}|^2
      \leq C_n^2 \|v\|_\cond^2,
    \end{align*}
    and hence $\|B_n^\perp \Lap B_n\|_\cond \leq C_n$.
  \end{proof}
\end{lemma}

\begin{prop}\label{thm:edge-boundary-O(n2)}
\version{}{\marginpar{Can Prop.~\ref{thm:edge-boundary-O(n2)} be applied  to the infinite Sierpinski lattice? Must compute ...}}
  If the estimate
  \begin{equation}\label{eqn:edge-boundary-O(n2)}
    \| B_n^\perp \Lap B_n\| = O(n^2), \q \text{as } n \to \iy,
  \end{equation}
  is satisfied, then \Lap has no eigenvectors at \iy.
  \begin{proof}
    Since \Lap is semibounded, this follows from \cite{Jor78}. (The reader may find information in \cite{Jor76} and \cite{Jor77} regarding the general Hermitian case.)
  \end{proof}
\end{prop}

\begin{exm}[One-sided ladder model]\label{exm:one-sided-ladder-model}
  The elementary ladder models help explain the effects of adding to the network, when the new portions added to the graph are somewhat ``peripheral'' to the original graph. Also, they provide an easy example of an infinite network in which a shortest path (see Remark~\ref{rem:shortest-path-metric}) may not exist, and then illustrate a technique for understanding what happens on a finite subgraph of interest which is embedded in an infinite graph.

Consider a graph which appears as a sideways ladder with $n$ rungs extending to the right:
  \begin{equation}\label{eqn:exm:one-sided-ladder-model}
    \xymatrix{
    *+[l]{\ga = x_0} \ar[r] \ar[d]
      & x_1 \ar[r] \ar[d] & x_2 \ar[r] \ar[d] & x_3 \ar[r] \ar[d]
      &\dots x_n \ar[d] \\
    *+[l]{\gw = y_0} & y_1 \ar[l] & y_2 \ar[l] & y_3 \ar[l] &\dots y_n \ar[l]
    }
  \end{equation}

  \version{}{\marginpar{Must add ...}
  We wish to find the minimizing current running down the leftmost edge, i.e., insert one amp at the top left corner and withdraw one amp from the bottom left corner.
 How does the flow respond as $n$ increases, and what happens as $n \to\iy$? \\

  [[To be added]]
  }
  
  The one-sided ladder model furnishes a situation where no shortest path exists, as mentioned in Remark~\ref{rem:shortest-path-metric}, if the resistances are defined as follows:
  \begin{equation}\label{eqn:exm:one-sided-ladder-model-no-shortest-path}
    \xymatrix{
    *+[l]{\ga = x_0} \ar@{-}[r]_{\frac14} \ar@{-}[d]_1
      & x_1 \ar@{-}[r]_{\frac1{16}} \ar@{-}[d]_{\frac1{4}} & x_2 \ar@{-}[r]_{\frac1{64}} \ar@{-}[d]_{\frac1{16}} & x_3 \ar@{-}[r] \ar@{-}[d]_{\frac1{64}}
      &\dots \ar@{-}[r]_{\frac1{4^n}} & x_n \ar@{-}[d]_{\frac1{4^n}} \\
    *+[l]{\gw = y_0} & y_1 \ar@{-}[l]_{\frac14} & y_2 \ar@{-}[l]_{\frac1{16}} & y_3 \ar@{-}[l]_{\frac1{64}} &\dots  \ar@{-}[l] & y_n \ar@{-}[l]_{\frac1{4^n}}
    }
  \end{equation}
  To find a path from \ga to \gw with minimal distance, one is led to consider paths stretching ever further off towards infinity. It is easy to see that the shortest path metric gives
  \linenopax
  \begin{align*}
    \distmin(\ga,\gw) = \lim \left[2\left(\frac14 + \frac1{16} + \dots + \frac1{4^n}\right) + \frac1{4^n}\right] = \frac23,
  \end{align*}
  but there is no $\cpath \in \Paths_{\ga,\gw}$ for which $\prob(\cpath) = \frac23$.
  Note that the Powers bound \eqref{eqn:def:power-bound} is violated, as $\gm(x_n) = 4^n+4^n+4^{n+1} \to \iy$.
\end{exm}

\begin{exm}[One-sided infinite ladder network]\label{exm:a,b-ladder}
  Consider the same ladder network as above, except now let the horizontal conductances grow geometrically, and let the vertical conducatances decay geometrically. More precisely, fix two positive numbers $\ga > 1 > \gb > 0$. Define the horizontal conductances between nearest neighbours by $\cond_{x_n, x_{n-1}} = \cond_{y_n, y_{n-1}} = \ga^n$, and define the vertical conductance of the ``rungs'' of the ladder by $\cond_{x_n,y_n} = \gb^n$:
  \linenopax
  \begin{equation}\label{eqn:exm:one-sided-ladder-model}
    \xymatrix{
    *+[l]{x_0} \ar@{-}[r]^{\ga} \ar@{-}[d]_{1}
      &  x_1 \ar@{-}[r]^{\ga^2} \ar@{-}[d]_{\gb} 
      &  x_2 \ar@{-}[r]^{\ga^3} \ar@{-}[d]_{\gb^2} 
      &  x_3 \ar@{-}[r]^{\ga^4} \ar@{-}[d]_{\gb^3}
      &\dots \ar@{-}[r]^{\ga^n}
      &  x_n \ar@{-}[r]^{\ga^{n+1}} \ar@{-}[d]_{\gb^n} & \dots \\
    *+[l]{y_0}\ar@{-}[r]^{\ga} 
      &  y_1 \ar@{-}[r]^{\ga^2} 
      &  y_2 \ar@{-}[r]^{\ga^3} 
      &  y_3 \ar@{-}[r]^{\ga^4} 
      &\dots \ar@{-}[r]^{\ga^n}
      &  y_n \ar@{-}[r]^{\ga^{n+1}} 
      & \dots
    }
  \end{equation}
  
  This network was suggested to us by Agelos Georgakopoulos as an example of a one-ended network with nontrivial \Harm. The function $u$ constructed below is the first example of an explicitly computed nonconstant harmonic function of finite energy on a graph with one end (existence of such a phenomenon was first proved in \cite{CaW92}). Numerical experiments indicate that this function is also bounded (and even that the sequences $\{u(x_n)\}_{n=0}^\iy$ and $\{u(y_n)\}_{n=0}^\iy$ actually converge very quickly), but we have not yet been able to prove this. Numerical evidence also suggests that \Lap is not essentially self-adjoint on this network, but we have not yet proved this, either. However, compare with the defect on the geometric integers discussed in \S\ref{sec:defect-on-geometric-integers}.
  
  This graph clearly has one end. We will show that such a network has nontrivial resistance boundary if and only if $\ga > 1$ and in this case, the boundary consists of one point for $\gb=1$, and two points for \gb such that $(1 + \frac1\ga)^2 < \ga / \gb^2$. It will be made clear that the paths $\cpath_x=(x_1,x_2,x_3,\dots)$ and $\cpath_y=(y_1,y_2,y_3,\dots)$ are equivalent in the sense of Definition~\ref{def:harmonic-equivalence-of-paths} if and only if $\gb = 1$.
  
  For presenting the construction of $u$, choose $\gb < 1$ satisfying $4\gb^2 < \ga$ (at the end of the construction, we explain how to adapt the proof for the less restrictive condition $(1 + \frac1\ga)^2 < \ga/\gb^2$). We now construct a nonconstant $u \in \Harm$ with $u(x_0) = 0$ and $u(y_0)=-1$. If we consider the flow induced by $u$, the amount of current flowing through one edge determines $u$ completely (up to a constant). Once it is clear that there are two boundary points 
  in this case, it is clear that specifying the value of $u$ at one 
  (and grounding the other) 
  determines $u$ completely.
  
  Due to the symmetry of the graph, we may abuse notation and write $n$ for $x_n$ or $y_n$, and $\check n$ for the vertex ``across the rung'' from $n$. For a function $u$ on the ladder, denote the horizontal increments and the vertical increments by  
  \linenopax
  \begin{align*}
    \gd u(n) := u(n+1) - u(n)
    \qq\text{and}\qq
    \gs u(n) := u(n) - u(\check n),
  \end{align*}
  respectively. 
  Thus, for $n \geq 1$, we can express the equation $\Lap u(n) = 0$ by
  \linenopax
  \begin{align*}
    \Lap u(n) = \ga^n \gd u(n-1) - \ga^{n+1} \gd u(n) + \gb^n \gs u(n) = 0,
  \end{align*}
  which is equivalent to 
  \linenopax
  \begin{align*}
    \gd u(n) = \frac1\ga \gd u(n-1) + \frac{\gb^n}{\ga^{n+1}} \gs u(n).  
  \end{align*}
  Since symmetry allows one to assume that $u(\check n) = 1- u(n)$, we may replace $\gs u(n)$ by $2u(n)+1$ and obtain that any $u$ satisfying
  \linenopax
  \begin{align}\label{eqn:u-harmonic-on-a,b-ladder}
    u(n+1) = u(n) + \tfrac{u(n)-u(n-1)}\ga + \tfrac2\ga\left(\tfrac\gb\ga\right)^n u(n)
     + \tfrac1\ga \left(\tfrac\gb\ga\right)^n
  \end{align}
  is harmonic. It remains to see that $u$ has finite energy.
  
  Our estimate for $\energy(u) < \iy$ requires the assumption that $\ga > 4\gb^2$, but numerical computations indicate that $u$ defined by \eqref{eqn:u-harmonic-on-a,b-ladder} will be both bounded and of finite energy, for any $\gb < 1 < \ga$. First, note that $u(1) = \frac1\ga$ and so an immediate induction using \eqref{eqn:u-harmonic-on-a,b-ladder} shows that $\gd u(n) = u(n+1)-u(n)>0$ for all $n \geq 1$, and so $u$ is strictly increasing. Since $\gb < 1 < \ga$, we may choose $N$ so that 
  \linenopax
  \begin{align*}
    n \geq N \implies \left(\frac\gb\ga\right)^n < \frac{\ga-1}2.
  \end{align*}
  Then $n \geq N$ implies 
  \linenopax
  \begin{align}\label{eqn:ladder(u)-subexponential}
    u(n+1) \leq 2u(n) + \frac1\ga,
  \end{align}
  by using \eqref{eqn:u-harmonic-on-a,b-ladder} and the fact that $u(n)$ is increasing and $\frac\gb\ga < 1$.
  Now use \eqref{eqn:u-harmonic-on-a,b-ladder} to write
  \linenopax
  \begin{align*}
    \gd u(n)
    &= \tfrac1\ga (\gd u)(n-1) + \left(\tfrac2\ga u(n) + \tfrac1\ga\right) \left(\tfrac\gb\ga\right)^n \notag \\
    &= \tfrac1{\ga^n} (\gd u)(0) + \sum_{k=0}^{n-1} \tfrac1{\ga^k} \left(\tfrac2\ga u(n-k) + \tfrac1\ga\right) \left(\tfrac\gb\ga\right)^{n-k} \notag \\
    &= \tfrac1{\ga^{n+1}} + \tfrac{\gb(1-\gb^n)}{\ga^{n+1}(1-\gb)} 
     + \tfrac2{\ga^{n+1}} \sum_{k=1}^{n} \gb^k u(k),
  \end{align*}
  where the second line comes by iterating the first, and the third by algebraic simplification. Applying the estimate \eqref{eqn:ladder(u)-subexponential} gives 
  \linenopax
  \begin{align*}
    2 \sum_{k=1}^{n} \gb^k u(k)
    &\leq 2^2 \sum_{k=1}^{n} \gb^k u(k-1) + \tfrac2\ga \sum_{k=1}^{n} \gb^k
        = 2^2 \sum_{k=2}^{n} \gb^k u(k-1) + 2\tfrac\gb\ga \cdot \tfrac{1-\gb^n}{1-\gb},
  \end{align*}
  and iterating gives
  \linenopax
  \begin{align}\label{eqn:ladder(du)-final-est}
    \gd u(n)
    \leq \frac1{\ga^{n+1}} \left(1 + \frac{\gb(1-\gb^n)}{1-\gb} + \frac{(2\gb)^n}\ga + 2\frac\gb\ga \sum_{k=0}^{n-1} 2^k \frac{\gb^k-\gb^n}{1-\gb}\right).
  \end{align}
  Now the energy $\energy(u) = \sum_{n=0}^\iy \ga^{n+1} \left(\gd u(n)\right)^2$ can be estimated by using \eqref{eqn:ladder(du)-final-est} as follows:
  \linenopax
  \begin{align*}
    \energy(u) 
    &\leq \sum_{n=0}^\iy \frac1{\ga^{n+1}} \left(1 + \frac{\gb(1-\gb^n)}{1-\gb} + \frac{(2\gb)^n}\ga + \frac{2\gb + 2\gb^{n+1} - 2^{n+2}\gb^{n+1} - 2^2\gb^{n+2} + (2\gb)^{n+2}}{\ga (1-\gb) (2\gb-1)}\right)^2
  \end{align*}
  and the condition $\ga > 4\gb^2$ ensures convergence.
  
  Note that this computations above can be slightly refined: instead of $\ga > 4\gb^2$, one need only assume that $\ga > (1 + \frac1\ga)^2\gb^2$. Then, fix $\ge>0$ for which $\ga / \gb^2 > (1 + \frac1\ga)^2 + \ge$ and choose $N$ so that $n \geq N$ implies $\left(\gb/\ga\right)^n < 1 + \frac1\ga + \ge(1 + 2\ga + \ga\ge)$. Then the calculations can be repeated, with most occurrences of $2$ replaced by $1+\frac1\ga+\ge$. 
\end{exm}

\begin{remark}\label{rem:comp-ladder-example-to-Z1}[Comparison of Example~\ref{exm:a,b-ladder} to the 1-dimensional integer lattice]
  Example~\ref{def:geometric-half-integers} shows that for $\ga > 1$, the ``geometric half-integers'' network
  \linenopax
  \begin{align*}
    \xymatrix{
      0 \ar@{-}[r]^{\ga} 
      & 1 \ar@{-}[r]^{\ga^2} 
      & 2 \ar@{-}[r]^{\ga^3} 
      & 3 \ar@{-}[r]^{\ga^4} 
      & \dots
    }
  \end{align*} 
  supports a monopole but not a harmonic function of finite energy. These conductances correspond to the biased random walk where, at each vertex, the walker has transition probabilities 
  \linenopax
  \begin{align*}
    p(n,m) = 
    \begin{cases}
      \frac1{1+\ga}, &m=n-1,  \\
      \frac{\ga}{1+\ga}, &m=n+1. 
    \end{cases}
  \end{align*}
  In particular, this is a spatially homogeneous distribution.
  In contrast, the random walk corresponding to Example~\ref{exm:a,b-ladder} has transition probabilities 
  \linenopax
  \begin{align*}
    p(n,m) = 
    \begin{cases}
      \frac1{1 + \ga + \left(\frac\gb\ga\right)^n}, &m=n-1,\\
      \frac{\ga}{1 + \ga + \left(\frac\gb\ga\right)^n}, &m=n+1, \vstr[2.5] \\
      \frac{\left(\gb/\ga\right)^n}{1 + \ga + \left(\frac\gb\ga\right)^n}, &m=\check n. \vstr[2.5]
    \end{cases}
  \end{align*}
  Thus, Example~\ref{exm:a,b-ladder} is geometrically asymptotic to the geometric half-integers. 
  
  \pgap
  
  One can even think of Example~\ref{exm:a,b-ladder} as describing the \emph{scattering theory} of the geometric half-integer model, in the sense of \cite{Lax-Phillips}. In this theory, a wave (described by a function) travels towards an obstacle. After the wave collides with the obstacle, the original function is transformed (via the ``scattering operator'') and the resulting wave travels away from the obstacle. The scattering is typically localized in some sense, corresponding to the location of the collision.
  
  To see the analogy with the present scenario, consider the current flow defined by the harmonic function $u$ constructed in Example~\ref{exm:a,b-ladder}, i.e., induced by Ohm's law: $\curr(x,y) = \cond_{xy}(u(x) - u(y))$. With $\act_{_{|\curr|}}(x) := \frac12\sum_{\{z \suth \curr(x,z) > 0\}} |\curr(x,z)|$, this current defines a Markov process with transition probabilities
  \linenopax
  \begin{align*}
    P(x,y) = \frac{\curr(x,y)}{\act_{_{|\curr|}}(x)},
    \q\text{if}\q \curr(x,y)>0,
  \end{align*}
  and $P(x,y) = 0$ otherwise; see \S\ref{sec:Probabilistic-interpretation} and also \cite{RANR}. This describes a random walk where a walker started on the bottom edge of the ladder will tend to step leftwards, but with a geometrically increasing probability of stepping to the upper edge, and then walking rightwards off towards infinity. The walker corresponds to the wave, which is scattered as it approaches the geometrically localized obstacle at the origin.
\end{remark}

\section{Remarks and references}
\label{sec:Remarks-and-References-examples}

The material of this chapter is an assortment of examples, some finite weighted graphs and others infinite. The infinite models are understood with the use of limit considerations. A good background reference is \cite{LevPerWil08}.
Additionally, the reader may find the sources \cite{Woess00, Woe03, KaimanovichWoess02, Kig03, Kig01, Lyo03, LPS03, LyPe03, HeLy03, HaJo02, Per09, FHS09, MiSh09} to be useful.

\version{}{
\begin{exm}[Two-sided ladder model]\label{exm:two-sided-ladder-model}
  Consider a graph which appears as a sideways ladder with $n$ rungs extending in either direction:
  \begin{equation}\label{eqn:exm:two-sided-ladder-model}
    \xymatrix{
    x_{-n} \ar[d] & \dots \ar[d] \ar[l] & x_{-1} \ar[d] \ar[l]
      &{\ga = x_0} \ar[l] \ar[r] \ar[d]
      & x_1 \ar[r] \ar[d] & \dots \ar[d] \ar[r] & x_n \ar[d] \\
    y_{-n} \ar[r] & \dots \ar[r] & y_{-1} \ar[r]
      &{\gw = y_0} & y_1 \ar[l] & \dots \ar[l] & y_n \ar[l]
    }
  \end{equation}
  We wish to find the minimizing current upon inserting one amp at the top center
  and withdrawing one amp from the bottom center.
  \version{}{\marginpar{Must add ...}}
  How does the flow respond as
  $n$ increases, and what happens as $n \to\iy$? \\

  [[To be added]]
\end{exm}
}

\version{}{
\pgap

\subsection{Vertices of infinite degree}
\label{sec:vertices-of-infinite degree}

As in \S\ref{sec:Properties-of-Lap-on-HE}, vertices of infinite degree and unbounded conductances make it possible that $\spn\{\gd_x - \gd_o\}$ is not dense in $\spn\{\gd_x\}$. 

\begin{exm}[The augmented integers]\label{exm:augmented-integers}
  Consider the integers with edges between nearest neighbours $(\bZ,\cond)$, where $\cond$
  \begin{equation}\label{eqn:exm:two-sided-ladder-model}
    \xymatrix{
    x_{-n} \ar[d] & \dots \ar[d] \ar[l] & x_{-1} \ar[d] \ar[l]
      &{\ga = x_0} \ar[l] \ar[r] \ar[d]
      & x_1 \ar[r] \ar[d] & \dots \ar[d] \ar[r] & x_n \ar[d] \\
    y_{-n} \ar[r] & \dots \ar[r] & y_{-1} \ar[r]
      &{\gw = y_0} & y_1 \ar[l] & \dots \ar[l] & y_n \ar[l]
    }
  \end{equation}
\end{exm}

\pgap
}

%% file: tree-examples.tex

\chapter{Infinite trees}
\label{sec:tree-networks}

\headerquote{A great discovery solves a great problem but there is a grain of discovery in the solution of any problem. Your problem may be modest; but if it challenges your curiosity and brings into play your inventive faculties, and if you solve it by your own means, you may experience the tension and enjoy the triumph of discovery.}{---~G.~P\'{o}lya}


\begin{figure}
  \centering
  \scalebox{1.0}{\includegraphics{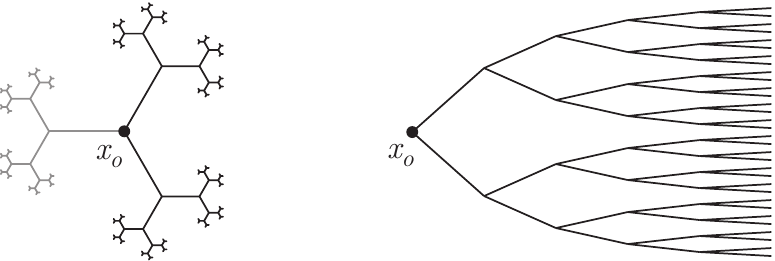}}
  \caption{\captionsize The homogeneous tree of degree 3 (left) and the binary tree from symbolic dynamics (right). The root of the tree is labeled $x_o$. If the grey branch is pruned from the homogeneous tree, the two become isomorphic.}
  \label{fig:homogeneous-vs-symbolic}
\end{figure}

The $n$-ary trees play an important role in symbolic dynamics, and they support a rich family of nontrivial harmonic functions of finite energy. These graphs are essentially homogeneous trees; the only difference is that the homogeneous tree has one more branch at the root, as can be seen from Figure~\ref{fig:homogeneous-vs-symbolic}. We use the latter examples as they are simpler yet still sufficient for our purposes, and because of our related interest in symbolic dynamics. However, almost all remarks extend to the homogeneous trees without effort; these examples are well-studied because of their close relationship with group theory (especially free groups). 
\version{}{\marginpar{Must include the case for varying \cond!}} Also, they provide an excellent testbed for studying the effects of varying \cond, and for illustrating several of our theorems.
A network whose underlying graph is a homogeneous tree always allows for the construction of a nontrivial harmonic function. In particular, \Fin is not dense in \HE by Lemma~\ref{thm:TFAE:Fin,Harm,Bdy} that these are equivalent.

\begin{remark}\label{rem:convention-negative-tree-branch}
  If the origin were removed from the binary tree, we adopt the convention that vertices in one component are ``positive'' and indices in the other are ``negative''. If the vertices are indexed with binary numbers (using the empty string $\es$ to denote the origin $o=x_\es$), then indices beginning with $1$ are positive and indices beginning with $0$ are negative.
\end{remark}

\begin{figure}
  \centering
  \scalebox{1.0}{\includegraphics{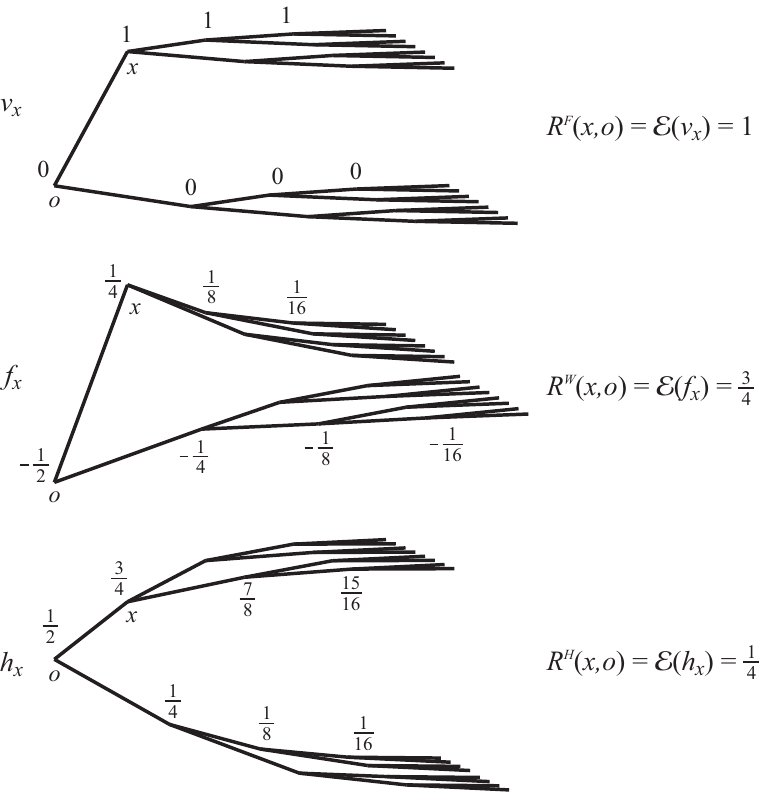}}
  \caption{\captionsize The reproducing kernel on the tree with $\cond = \one$. For a vertex $x$ which is adjacent to the origin $o$, this figure illustrates the elements $v_x$, $f_x = \Pfin v_x$, and $h_x = \Phar v_x$; see Example~\ref{exm:binary-tree:reproducing-kernel}. }
  \label{fig:tree-repkernels}
\end{figure}

\begin{figure}
  \centering
  \scalebox{1.0}{\includegraphics{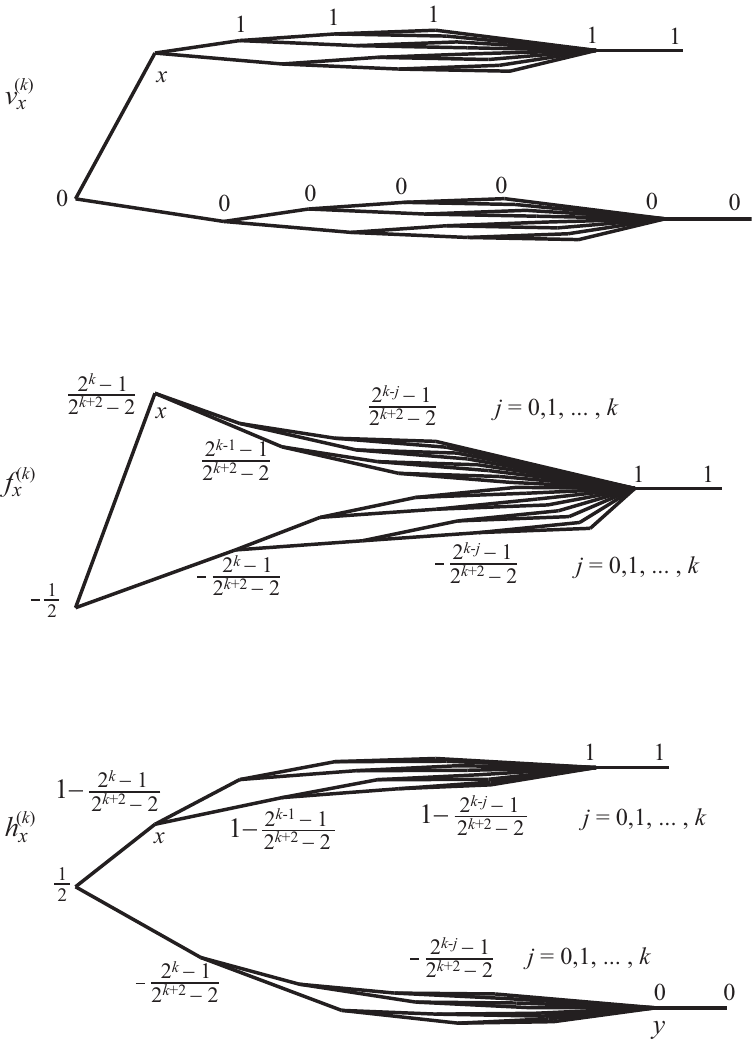}}
  \caption{\captionsize Approximants to the reproducing kernel on the tree with $\cond = \one$; see Example~\ref{exm:binary-tree:reproducing-kernel}.}
  \label{fig:tree-repkernel-approximants}
\end{figure}

\begin{exm}[The reproducing kernel on the tree]
  \label{exm:binary-tree:nontrivial-harmonic}
  \label{exm:binary-tree:reproducing-kernel}
  Let $(\sT,\one)$ be the binary tree network in Figure~\ref{fig:homogeneous-vs-symbolic} with constant conductances. Figure~\ref{fig:tree-repkernels} depicts the embedded image of a vertex $v_x$, as well as its decomposition in terms of \Fin and \Harm. We have chosen $x$ to be adjacent to the origin $o$; the binary label of this vertex would be $x_1$.

  In Figure~\ref{fig:tree-repkernels}, numbers indicate the value of the function at that vertex; artistic liberties have been taken. If vertices $s$ and $t$ are the same distance from $o$, then $|f_x(s)|=|f_x(t)|$ and similarly for $h_x$. Note that $h_x$ provides an example of a nonconstant harmonic function in \HE. Another key point is that $h_x \notin \ell^2$, see Corollary~\ref{thm:nontrivial-harmonic-fn-is-not-in-L2}. It is easy to see that $\lim_{z \to \pm\iy} h_x(z) = \frac12 \pm \frac12$, whence $h_x$ is bounded. 
  
  For $f_x = \Pfin v_x$ of Figure~\ref{fig:tree-repkernels}, the illustration of $f_x^{(k)}$ in Figure~\ref{fig:tree-repkernel-approximants} is the projection of $v_x$ (or $f_x$) to $\spn\{\gd_x \suth x \in G_k\}$, where $G_k$ consists of all vertices within $k$ steps of $o$. The lines at the right side of each figure just indicate that the function is constant on the remainder of the graph (at value 0 or 1); in particular, note that $f_x^{(k)}(y) = 0$ for every vertex $y$ which is at least $k+1$ steps from the origin. Also, observe that 
  \[\Lap f_x^{(k)} = \gd_x - \gd_y + \sum_{s \in \bd G_{k+1}^+} \frac{\gd_s}{2^{k+2}-2} - \sum_{t \in \bd G_{k+1}^-} \frac{\gd_t}{2^{k+2}-2},\]
  where $\bd G_{k+1}^+$ is the subset of $\bd G_{k+1}$ that lies on the positive branch, etc. It is interesting to note that if one were to identify all the vertices of $\bd G_{k+1} = \bd G_{k+1}^+ \cup \bd G_{k+1}^+$, then $f_x^{(k)}$ would become harmonic at this new vertex. Observe also that $h_x^{(k)} = v_x - f_x^{(k)}$ is its orthogonal complement and is harmonic everywhere except on $\bd G_{k+1}$.
\end{exm}

\begin{exm}[A function of finite energy which is not approximable by \Fin]
  \label{exm:binary-tree:nonapproximable}  
  We continue to refer to in Figure~\ref{fig:tree-repkernel-approximants}. Since $f_x^{(k)} \in \Fin$ and it is easy to see that $\|f_x-f_x^{(k)}\|_\energy \to 0$ and that $\Lap f_x = \gd_x - \gd_y$, this approximation verifies that $f_x = \Pfin v_x$. It also shows that $\min_{f \in \Fin} \|v_x - f\|_\energy = \frac14$.
\end{exm}

\begin{exm}[A monopole which is not a ``dipole at infinity'']
  \label{exm:monopole-on-binary-tree}
  Let $(\sT,\one)$ be the binary tree network in Figure~\ref{fig:homogeneous-vs-symbolic} with constant conductances.  Let $|x|$ be the number of edges in the path connecting $x$ to $o$. Define a function
  \linenopax
  \begin{equation}\label{eqn:monopole-on-binary-tree}
    w_o(x) = 1-\frac1{2^{|x|}},
  \end{equation}
  so that essentially $w_o = 2|h_x-\frac12|$ for $h_x$ of Example~\ref{exm:binary-tree:reproducing-kernel}. It is easy to check that $\Lap w_o = -\gd_o$ so that $w_o$ is a monopole at the root $o$. To see that $w_o \in \HE$,
  \begin{align*}
    \energy(w_o)
    &= 2 \sum_{n=0}^\iy 2^n \left(\left(1-\frac1{2^{n}}\right) - \left(1-\frac1{2^{n+1}}\right)\right)^2 \\
    &= 2 \sum_{n=0}^\iy 2^n \left(\frac1{2^{n+1}}\right)^2 \\
    &= \frac12 \sum_{n=0}^\iy \frac1{2^n} \\
    &= 1.
  \end{align*}

  However, $w_o$ is not a ``dipole at infinity'' in the sense that there is no sequence $\{x_n\}$ of distinct and successively adjacent vertices for which $\{v_{x_n}\}$ converges to $w_o$ (this is in contrast to the integer lattices \bZd, $d \geq 3$). Observe that $R^F(x,y)$ coincides with shortest-path distance on this network (as it does on any tree). If $\{x_n\}$ is a sequence tending to \iy (i.e., for any $N$, there is an $n$ such that $x_n$ is more than $n$ steps from $o$ for all $n \geq N$), then $\energy(v_{x_n}) = R^F(x_n,o) = n$, so that $w_o$ is not a limit of a sequence of dipoles. 
  
  Of course, since $\{v_x\}$ is dense in \HE, $w_o$ is the limit of \emph{linear combinations} of dipoles. In fact, let $\bd G_k = \{x \in \sT \suth R(x,o) = k\}$ as before. Then
  \linenopax
  \begin{align*}
    w_o(x) = \lim_{k \to \iy} \sum_{z \in \bd G_k} \frac{v_z}{2^k}.
  \end{align*}
\end{exm}

\begin{exm}[A function with nonvanishing boundary sum]
  \label{exm:nonvanishing normal derivatives}
  In Theorem~\ref{thm:E(u,v)=<u,Lapv>+sum(normals)}, we showed
  \linenopax
  \[\la u, v \ra_\energy = \lim_{k \to \iy} \sum_{x \in \inn \Graph_k} u(x) \Lap v(x)
      + \lim_{k \to \iy} \sum_{x \in \bd \Graph_k} u(x) \dn v(x).\]
  Let $w_o(x) = 1-\frac1{2^{|x|}}$ be the monopole from Example~\ref{exm:monopole-on-binary-tree}. With $G_k := \{x \suth |x| \leq k\}$, we have 
  \begin{align*}
    \dn{w_o}(x) = \left(1-\frac1{2^{k}}\right) - \left(1-\frac1{2^{k-1}}\right) = \frac1{2^{k}},
    \q \text{ for } x \in \bd G_k = \{x \suth |x| = k\}.
  \end{align*}
  Since $\Lap w = -\gd_o$, we have $\sum_{\Graph_k} w_o(x) \Lap w_o(x) = -w(o)$, for each $k$, and the energy of $w_o$ is
  \begin{align*}
    \energy(w_o) 
     = \la w_o, w_o\ra_\energy
    &= -\left(1-\frac1{2^{0}}\right)
      + \lim_{k \to \iy} \sum_{x \in \bd \Graph_k} \left(1-\frac1{2^{k}}\right)\frac1{2^{k}} = 1.
  \end{align*}
  
  For $h_x$, the harmonic function with $\energy(h_x)=\frac14$ in Example~\ref{exm:binary-tree:nontrivial-harmonic}, this becomes
  \linenopax
  \[\energy(h_x,h_x) = 
    \lim_{k \to \iy} \sum_{x \in \bd \Graph_k} h_x(x) \dn{h_x}(x)
    =\frac14.\]
  In fact, one can obtain this by computing the boundary term directly: each of the $2^{k-1}$ vertices in $\bd \Graph_k^+$ is connected by a single edge to $G_k$, and similarly for the $2^{k-1}$ vertices in $\bd \Graph_k^-$, so
  \linenopax
  \begin{align*}
    \sum_{y \in \bd \Graph_k} h_x(y) \dn{h_x}(y)
    &= 2^{k-1} \frac{2^{k+1}-1}{2^{k+1}} \cdot \frac{1}{2^{k+1}}
     + 2^{k-1} \frac{1}{2^{k+1}} \cdot \frac{-1}{2^{k+1}}
    = \frac14 \left(1-\frac1{2^k}\right).
  \end{align*}
\end{exm}

\begin{exm}[The tree supports many nontrivial harmonic functions]
  \label{exm:a-forest-of-binary-trees}
  We can use $h_x$ of Example~\ref{exm:binary-tree:nontrivial-harmonic} to describe an infinite forest of mutually orthogonal harmonic functions on the binary tree. Let $z \in \Graph$ be represented by a finite binary sequence, as discussed in Remark~\ref{rem:convention-negative-tree-branch}. Define a morphism (cf.~Definition~\ref{def:ERN-morphism}) $\gf_z:\Graph \to \Graph$ by prepending, i.e., $\gf_z(x) = zx$. This has the effect of ``rigidly'' translating the the tree so that the image lies on the subtree with root $z$. Then $h_z := h_x \comp \gf_z$ is harmonic and is supported only on the subtree with root $z$. The supports of $h_{z_1}$ and $h_{z_2}$ intersect if and only if $\Im(\gf_{z_i}) \ci \Im(\gf_{z_j})$. For concreteness, suppose it is $\Im(\gf_{z_1}) \ci \Im(\gf_{z_2})$. If they are equal, it is because $z_1=z_2$ and we don't care. Otherwise, compute the dissipation of the induced currents
  \linenopax
  \begin{align*}
    \la \drp h_{z_1}, \drp h_{z_2} \ra_\diss
    = \tfrac12 \negsp \sum_{(x,y) \in \gf_{z_1}(\edges)} \negsp
      \ohm(x,y) \drp h_{z_1}(x,y), \drp h_{z_2}(x,y).
  \end{align*}
  Note that $\drp h_{z_2}(x,y)$ always has the same sign on the subtree with root $z_1 \neq o$, but $\drp h_{z_1}(x,y)$ appears in the dissipation sum positively signed with the same multiplicity as it appears negatively signed. Consequently, all terms cancel and $0 = \la \drp h_{z_1}, \drp h_{z_2} \ra_\diss = \la h_{z_1}, h_{z_2} \ra_\energy$ shows $h_{z_1} \perp h_{z_2}$.
\end{exm}

\begin{exm}[Haar wavelets and cocycles]
  \label{exm:Haar-wavelets-and-cocycles}
  \version{}{\marginpar{Add cocycles}}
  Example~\ref{exm:a-forest-of-binary-trees} can be heuristically described in terms of Haar wavelets. Consider the boundary of the tree as a copy of the unit interval with $h_x$ as the basic Haar mother wavelet; via the ``shadow'' cast by $\lim_{n \to \pm \iy} h_x(x_n) = \pm 1$. Then $h_z$ is a Haar wavelet localized to the subinterval of the support of its shadow, etc. Of course, this heuristic is a bit misleading, since the boundary is actually isomorphic to $\{0,1\}^\bN$ with its natural cylinder-set topology.
\end{exm}

\begin{exm}[Why the harmonic functions may not be in the domain of \Lap]
  \label{exm:Harm-notin-domLap}
  We have not been able to construct an example in which we can prove that \Harm is not contained in $\dom \Lap$, but we do have the following suggestive example, motivated by Lemma~\ref{thm:powers-bd-implies-Fin1=0}.
  As before, let $V := \spn\{v_x\}_{x \in \verts}$ and let $\Lap = \LapV$ denote the closure of the Laplacian when taken to have the dense domain $V$. 
  Let $h \in \Harm$. If $h$ were an element of $\dom \LapV$, then by \eqref{eqn:def:domain-of-clo(S)}, we would have a sequence 
  \linenopax
  \begin{align}\label{eqn:def:domain-of-clo(S)}
    \dom \opclosure S := \{u \suth \lim_{n \to \iy}\|u-u_n\|_\sH = \lim_{n \to \iy}\|v-Su_n\|_\sH = 0\}
  \end{align}
  Again, think of the vertices of the tree as begin labelled by a word on $\{0,1\}$, that is, a finite binary string. If $x=\word$, then $|\word|$ is the length of the word and corresponds to the number of edges between $x$ and $o$ (i.e., shortest path distance to the root). Using $\word_1$ to denote the first coordinate of \word, define the function
  \linenopax
  \begin{align}\label{eqn:harmonic-approximants-on-tree}
    h_n := \frac1{2^n} \sum_{|\word|=n} \left(\sum_{\word_1=1} v_\word - \sum_{\word_1=0} v_\word\right).
  \end{align}
  Since $h_n$ is a (finite) sum of all the dipoles at distance $n$ from $o$, with half weighted by $2^{-n}$ and the other half weighted by $-2^{-n}$, it is clear that $h_n \in \spn\{v_x\}$. One can check that for $G_n = \{\word \suth |\word| \leq n\}$,
  \linenopax
  \begin{align*}
    h(\word) = 
    \begin{cases}
      h(\word), &|\word| \leq n, \\
      \pm 2^{-n}, &\text{else},
    \end{cases}
  \end{align*}
  whence it is immediate that $\lim_{n \to \iy} \|h_n - h\|_\energy =0$.
  One can also check that
  \linenopax
  \begin{align*}
    \Lap h_n = \frac1{2^n} \sum_{|\word|=n} \left(\sum_{\word_1=1} \gd_\word - \sum_{\word_1=0} \gd_\word \right)
  \end{align*}
  since the positive and negative weights of $\gd_0$ cancel out. If $\word \neq \word'$ but $|\word|=|\word'|$, then they cannot be neighbours, and hence $\gd_{\word}$ and $\gd_{\word'}$ are orthogonal with respect to \energy. It is then easy to compute 
  \linenopax
  \begin{align*}
    \lim_{n \to \iy} \| \Lap h_n - \Lap h\|_\energy^2
    = \lim_{n \to \iy} \| \Lap h_n\|_\energy^2
    = \lim_{n \to \iy} \frac1{2^n} \sum_{|\word|=n} \| \gd_{\word}\|_\energy^2
    = \lim_{n \to \iy} \frac1{2^n} \cdot 2^n \cdot 3
    = 3 \neq 0.
  \end{align*}
\end{exm}

\begin{exm}\label{exm:monopole-in-range-of-1+Lap}
  On the binary tree with $\cond = \one$, the monopole $w(x) = 2^{-|x|}$ can be written as $v+\Lap v$ for $v \in \dom\LapV$:
  \linenopax
  \begin{align*}
    v(x) := 2^{-|x|} - \sqrt2 \left(1 + \tfrac1{\sqrt2}\right)^{n+1}.
  \end{align*}
  We leave it to the reader to check that this $v$ satisfies the above equation and also
  \linenopax
  \begin{align*}
    \energy(v) = \tfrac{2(46\sqrt2 - 65)}{(2-\sqrt2)^2 (\sqrt2-1)}, \;
    \energy(\Lap v) = \tfrac{4(99-70\sqrt2)}{(2-\sqrt2)^4 (\sqrt2 - 1)}, 
    \; \text{and} \;
    \sum_{\verts} v \Lap v = \tfrac27 (-23 + 17 \sqrt2). 
  \end{align*}
\end{exm}

\begin{figure}
  \centering
  \scalebox{1.0}{\includegraphics{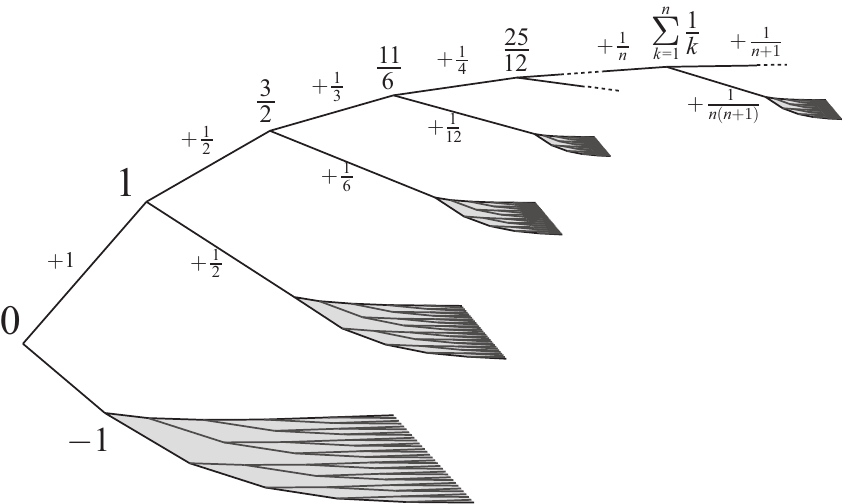}}
  \caption{\captionsize An unbounded harmonic function of finite energy. See Example~\ref{exm:unbounded-harmonic}.}
  \label{fig:unbounded-harmonic}
\end{figure}

\begin{exm}[An unbounded harmonic function of finite energy]
  \label{exm:unbounded-harmonic}
  Figure~\ref{fig:unbounded-harmonic} is a sketch of an unbounded harmonic function of finite energy on the binary tree with $\cond = \one$. To construct it, pick one ray from $o$ to \iy, and let
  \linenopax
  \begin{align*}
    h(x) = \sum_{k=1}^{|x|} \frac1k
  \end{align*}
  for $x$ along this ray. Then if $x$ is in the ray and $y \nbr x$, fix $h(y)$ so that $h$ is harmonic at $x$ (i.e., $h(y) = h(x)+\frac1{|x|(|x|+1)}$), and define $h$ along the rest of this branch by 
  \linenopax
  \begin{align*}
    h(z) = \frac{h(y)-h(x)}{2^{|y-z|}}.
  \end{align*}
  If $w$ denotes the monopole at $o$ defined by $w(x) = 2^{-|x|}$ as discussed in Example~\ref{exm:monopole-in-range-of-1+Lap} and previously, then we are essentially attaching a scaled copy of $w$ to each neighbour of the chosen ray. See Figure~\ref{fig:unbounded-harmonic}.
  
  It is clear that $h(x) \to \iy$ logarithmically along the chosen ray; the energy coming from $h(x)$ on this ray is
  \linenopax
  \begin{align*}
    \energy(h)|_{\text{ray}} = \sum_{n=1}^\iy \left(\frac1n\right)^2 = \frac{\gp^2}6.
  \end{align*}
  The energy from each branch incident upon the ray is
  \linenopax
  \begin{align*}
    \energy(h)|_{\text{branch}(n)} 
    = \left(\frac1{n(n+1)}\right)^2 + \frac1{n(n+1)}\energy(w)
    =\frac1{n^2(n+1)^2} + \frac1{n(n+1)}.
  \end{align*}
  Summing up, $\energy(h) = \energy(h)|_{\text{ray}} + \sum_{n=0}^\iy \energy(h)|_{\text{branch}(n)} = \frac{\gp^2}2 - 2$. We leave it to the reader to check that $h$ is harmonic.
\end{exm}

\begin{exm}\label{exm:Lapu-not-summable}
  On the binary tree with $\cond = \one$, the function $u_\gx(x) = \gx^{-|x|}$ has energy
  \linenopax
  \begin{align*}
    \energy(u_\gx) = 2\frac{(1-\gx)^2}{1-2\gx^2} < \iy, 
    \q\text{for } \gx \in [0,\tfrac1{\sqrt2}).
  \end{align*}
  Moreover, the Discrete Gauss-Green formula applies to this example with
  \linenopax
  \begin{align*}
    \sum_{\verts} u_\gx \Lap u_\gx
    = (1-\gx)\left(2 + (2\gx-1)\tfrac{2\gx}{1-2\gx^2}\right), 
    \q\text{for } \gx \in [0,\tfrac1{\sqrt2}).
  \end{align*}
  However, for $\gx \in [0,\tfrac12)$, one also has $\sum_{\verts} \Lap u_\gx = (1-\gx) \left(2 + (2\gx-1)\tfrac{1}{1-2\gx} \right)$. Thus, for $\frac12 < \gx < \frac1{\sqrt2}$, one has $\sum \Lap u_\gx = \iy$, even though $\energy(u_\gx) < \iy$ and $\sum u_\gx \Lap u_\gx < \iy$.
\end{exm}

\section{Remarks and references}
\label{sec:Remarks-and-References-tree-examples}

The infinite trees offer an especially attractive source of examples, and there are theories devoted to them; see \cite{Car72, Car73a, Car73b}. Here we merely scratch the surface. Excellent books and introductions are \cite{Lyons:ProbOnTrees, Peres99, LevPerWil08}.
The reader may also find the sources \cite{Dha98, DfdeG04, DLP09, HLL08, NaPe08a, NaPe08b} to be helpful, and the preprint \cite{Kig09b} deals specifically with trees in our context, and in relation to self-similar fractals.

%% file: lattice-examples.tex

\chapter{Lattice networks}
\label{sec:lattice-networks}

\headerquote{Observe also (what modern writers almost forgot, but some older writers, such as Euler and Laplace, clearly perceived) that the role of inductive evidence in mathematical investigation is similar to its role in physical research.}{---~G.~P\'{o}lya}

The integer lattices $\bZd \ci \bRd$ are some of the most widely-studied infinite graphs and have an extensive literature; see \cite{DoSn84,Telcs06a}, for example. We begin with some results for the simple lattices; in \S\ref{sec:nonsimple-integer-lattices} we consider the case when \cond is nonconstant. Because the case when $\cond = \one$ is amenable to Fourier analysis, we are able to compute many explicit formulas for many expressions, including $v_x$ and $R(x,y)$. For $d \geq 3$, we even compute $R(x,\iy) = \lim_{y \to \iy} R(x,y)$ in Theorem~\ref{thm:finite-resistance-to-infinity-in-Zd} and give a formula for the monopole $w$. There is a small amount of overlap here with the results of \cite[\S{V.2}]{Soardi94}, where the focus is more on solving Poisson's equation $\Lap u = f$. In \S\ref{sec:magnetism} we employ our formulas in the refinement of an application to the isotropic Heisenberg model of ferromagnetism.

In the present context, we may choose canonical representatives when working pointwise: given $u \in \HE$, we use the representative which tends to 0 at infinity. We take this as a standing assumption for this section, as it allows us to use the Fourier transform without ambiguity or unnecessary technical details. To see that this is justified, note that $\ell^2(\cond)$ is dense\footnote{Technically, the embedded image of $J:\ell^2(\cond) \to \Fin$ is dense in \Fin; see Definition~\ref{def:J-quotient-inclusion-map}.} in \Fin by Theorem~\ref{thm:l2c-in-HE}, and hence dense in \HE for these examples, as it is well-known that there are no nonconstant harmonic functions of finite energy on the integer lattices (we provide a proof in Theorem~\ref{thm:harmonics-on-Zd-are-linear} for completeness). Clearly, $\cond=\one$ implies that $\ell^2(\cond) = \ell^2(\one)$ and that all elements of $\ell^2(\cond)$ vanish at \iy. 

\begin{remark}\label{rem:lattices-and-numerical-analysis}
  As mentioned in Remark~\ref{rem:numerical-analysis}, one of the applications of the present investigation is to numerical analysis. Discretization of the real line amounts to considering a graph which is a scaled copy of the integers $\Graph_\epsilon = (\epsilon\bZ, \frac1\epsilon \one)$ where $\epsilon \bZ = \{\epsilon n \suth n\in \bZ\}$. After finding the solution to a given problem, as a function of the parameter $\epsilon$, one lets $\epsilon \to 0$. Let $x_n$ denote the vertex at $\epsilon n$. 
  \linenopax
    \begin{align*}
    \xymatrix{
      \dots \ar@{-}[r]_{1/\epsilon}
      & x_{-2} \ar@{-}[r]_{1/\epsilon} 
      & x_{-1} \ar@{-}[r]_{1/\epsilon} 
      & x_{0} \ar@{-}[r]_{1/\epsilon} 
      & x_{1} \ar@{-}[r]_{1/\epsilon} 
      & x_{2} \ar@{-}[r]_{1/\epsilon} 
      & x_{3} \ar@{-}[r]_{1/\epsilon} & \dots
    }
  \end{align*} 
  The \emph{difference operator} $D$ acts on a function on this network by $Df(x_n) := f(x_n)-f(x_{n+1})$. The adjoint of $D$ with respect to $\ell^2$ is $D^\ad f(x_n) = f(x_n) - f(x_{n-1})$. Then $D^\ad D = \Lap$.
\end{remark}

\section{Simple lattice networks}
\label{sec:simple-lattice-networks}

\begin{exm}[Simple integer lattices]\label{exm:infinite-lattices}
  The lattice network $(\bZd, \cond)$, with an edge between any two vertices which are one unit apart is called \emph{simple} or \emph{translation-invariant} when $\cond=\one$. The term ``simple'' originates in the literature on random walks.
\end{exm}

\begin{figure}
    \centering
    \scalebox{1.0}{\includegraphics{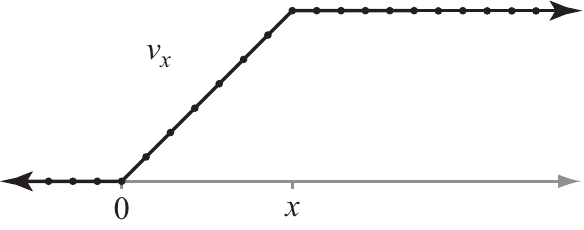}}
    \caption{\captionsize The function $v_x$, a solution to $\Lap v = \gd_x - \gd_0$ in $(\bZ,\one)$.}
    \label{fig:vx-in-Z1}
  \end{figure}

  One may compute the energy kernel directly using \eqref{eqn:def:R(x,y)-Lap}, that is, by finding a solution $v_x$ to $\Lap v = \gd_x - \gd_0$ as depicted in Figure~\ref{fig:vx-in-Z1}. Then $R(o,x) = v_x(x)-v_x(o) = x-0 = x$, which is unbounded as $x \to \iy$. This also provides an example of a function $v_x \in \HE$ for which $v_x \notin \ell^2(\cond)$, as discussed in \S\ref{sec:relating-HE-and-L2} and elsewhere.
  \linenopax
  In Lemma~\ref{thm:vx-on-Zd} we obtain a general formula for $v_x$ on $(\bZd, \one)$. Figure~\ref{fig:fejer-in-Z1} of Example~\ref{exm:Z-not-bounded} shows how this compares to Figure~\ref{fig:vx-in-Z1}.
  
  To see how the function $v=v_{x_1} = \charfn{[1,\iy)}$ may be approximated by elements of \Fin, define
  \linenopax
  \begin{equation}
    u_n(x_k) =
    \begin{cases}
      1-\frac kn, &1 \leq k \leq n, \\
      0, &\text{else}.
    \end{cases}
  \end{equation}
  The reader can verify that $u_n$ minimizes $\energy(v-u)$ over the set of $u$ for which $\spt(u) \ci [1,n-1]$ and that
  \linenopax
  \begin{align*}
    \energy(v-u_n)
    &= (1-(1-\tfrac1n))^2 + \sum_{k=1}^{n-1} \left((1-\tfrac kn)-(1-\tfrac{k-1}n)\right)^2 + (1-(1-0))^2 \\
    &= \frac1{n^2} + \frac{n-1}{n^2} \limas{n} 0.
  \end{align*}
  The fact that $v \in \clspn{\Fin}$ but $\lim_{k \to \iy} v(x_{k}) \neq \lim_{k \to \iy} v(x_{-k})$ reflects that $(\bZ,\one)$ has two \emph{graph ends}, unlike the other integer lattices; cf. \cite{PicWoess90}.
  Therefore, $(\bZ,\one)$ also provides an example of a network with more than one end which does not support nontrivial harmonic functions.

An explicit formula is given for the potential configuration functions $\{v_x\}$ on the simple $d$-dimensional lattice in Lemma~\ref{thm:vx-on-Zd}. By combining this formula with the dipole formulation of resistance distance
\linenopax
\begin{equation*}
  R(x,y) = v(x)-v(y), \qq \text{for } v=v_x-v_y,
  \text{from \eqref{eqn:def:R(x,y)-Lap},}
\end{equation*}
we are able to compute an explicit formula for resistance distance on the translation-invariant lattice network \bZd in Theorem~\ref{thm:R(x,y)-on-Zd}. Results in this section exploit the Fourier duality $\bZd \simeq \bTd$; \cite{Rud62} is a good reference. We are using Pontryagin duality of locally compact abelian groups; as an additive group of rank $d$, the discrete lattice \bZd is the dual of the $d$-torus \bTd. Conversely, \bTd is the compact group of unitary characters on \bZd (the operation in \bTd is complex multiplication). This duality is the basis for our Fourier analysis in this context. For convenience, we view \bTd as a $d$-cube, i.e., the Cartesian product of $d$ period intervals of length $2\gp$. In this form, the group operation in \bTd is written additively and the Haar measure on \bTd is normalized with the familiar factor of $(2\gp)^{-d}$.

In \cite{Polya21}, P\'{o}lya proved that the random walk on the simple integer lattice is transient if and only if $d \geq 3$; see \cite{DoSn84} for a nice introduction and a proof using \ERNs. In the present context, this can be reformulated as the statement that there exist monopoles on \bZd if and only if $d \geq 3$. We offer a new characterization of this dichotomy, which we recover in Theorem~\ref{thm:monopoles-on-Zd} via a new (and completely constructive) proof. In Remark~\ref{rem:HE(Zd)-lies-in-L2(Zd)} we show that in the infinite integer lattices, functions in \HE may be approximated by functions of finite support.

Sometimes P\'{o}lya's result is restated: the resistance to infinity is finite if and only if $d \geq 3$. There is an ambiguity in this statement which is specific to the nature of resistance metric. One interpretation is that one can construct a \emph{unit flow to infinity}; this is the terminology of \cite{DoSn84} for a current with $\act(\curr) = \gd_x$ and it is clear that this is the induced current of a monopole. Probabilistically, this definition may be rephrased: for a random walk beginning at $x \in \verts$, the expected hitting time of the sphere of (shortest-path) radius $n$ remains bounded as $n \to \iy$. This approach interprets ``infinity'' as the ``set of all points at infinity''. 

By contrast, we prove a much stronger result for the simple lattice networks \bZd in  Theorem~\ref{thm:finite-resistance-to-infinity-in-Zd}, where we show $\lim_{y \to \iy} R(x,y)$ is bounded as $y \to \iy$, for any $x \in \verts$. To see the strength of this result, note that the simple ($\cond = \one$) homogeneous trees of degree $d \geq 3$ have finite resistance to infinity, even though $\lim_{y \to \iy} R(x,y) = \iy$ for any $x \in \verts$, and any choice of $y \to \iy$. This is discussed further in Example~\ref{exm:binary-tree:nontrivial-harmonic} of the previous section. The heuristic explanation is that the resistance distance between two places is much smaller when there is high connectivity between them; there is much more connectivity between $x$ and the ``set of all points at infinity'' than between $x$ and a single ``point at infinity''.

In the next result, we obtain the Fourier transform of the Laplacian; we recently noticed that this corresponds almost identically to the inverse Fourier transform $H$ of the ``potential kernel'' of \cite[\S{V.2}]{Soardi94}.

\begin{lemma}\label{thm:Fourier-transform-of-Lap-on-Zd}
  On the \ERN $(\bZd, \one)$, the spectral (Fourier) transform of \Lap is multiplication by $S(t) = S(t_1,\dots,t_d) = 4 \sum_{k=1}^d \sin^2\left(\frac{t_k}2\right)$.
  \begin{proof}
    Each point $x$ in the lattice \bZd has $2d$ neighbours, so we need to find the $L^2(\bT^d)$ Fourier representation of
    \linenopax
    \begin{equation}\label{eqn:Fourier-transform-of-Lap-on-Zd}
      \Lap v(x) = (2d\id - T) v(x) = 2d v(x) - \sum_{k=1}^d v(x_1,\dots,x_k \pm 1,\dots x_d).
    \end{equation}
    Here, $t = (t_1, \dots, t_d) \in \bT^d$ and $x = (x_1, \dots, x_d) \in \bZ^d$. The \kth entry of $t$ can be written $t_k = t \cdot \ge_k$ where $\ge_k = [0,\dots,0,1,0,\dots,0]$ has the 1 in the \kth slot. Then moving one step in the lattice by $x \mapsto x+\ge_k$ corresponds to $e^{\ii x \cdot t} \mapsto e^{\ii t_k} e^{\ii x \cdot t}$ under the Fourier transform, and
    \linenopax
    \begin{align*}
      \widehat{\Lap v}(t) &= \left(2d - \sum_{k=1}^d (e^{\ii t_k} + e^{-\ii t_k} )\right) \hat v(t) \\
      &= 2\left(\sum_{k=1}^d (1-\cos(t_k))\right) \hat v(t) \\
      &= 4 \sum_{k=1}^d \sin^2\left(\frac{t_k}2\right) \hat v(t).
      \qedhere
    \end{align*}
  \end{proof}
\end{lemma}

\begin{lemma}\label{thm:vx-on-Zd}
  Let $\{v_x\}_{x \in \bZ^d}$ be the potential configuration on the integer lattice $\bZd$ with $\cond=\one$. Then for $y \in \bZd$,
  \linenopax
  \begin{equation}\label{eqn:vx-on-Zd}
    v_x(y) = \frac1{(2\gp)^d} \int_{\bT^d} \frac{\cos((x-y)\cdot t) - \cos(y\cdot t)}{S(t)} \,dt.
  \end{equation}
  \begin{proof}
    Under the Fourier transform, Lemma~\ref{thm:Fourier-transform-of-Lap-on-Zd} indicates that the equation $\Lap v_x = \gd_x - \gd_o$ becomes $S(t) \hat v_x = e^{\ii x\cdot t} - 1$, whence
    \linenopax
    \begin{align}\label{eqn:vx-on-Zd-as-exponential}
      v_x(y) = \frac1{(2\gp)^d} \int_{\bT^d} e^{-\ii y \cdot t} \frac{e^{\ii x\cdot t} - 1}{S(t)} \,dt.
    \end{align}
    Since we may assume $v_x$ is \bR-valued, the result follows.
  \end{proof}
\end{lemma}

The following result is well-known in the literature (cf.~\cite{DoSn84,Nash-Will59}, e.g.), but usually stated in terms of the current flow induced by the monopole.

\begin{theorem}\label{thm:monopoles-on-Zd}
  The network $(\bZd,\one)$ has a monopole
  \linenopax
  \begin{equation}\label{eqn:monopoles-on-Zd}
    w(x) = -\frac1{(2\gp)^d} \int_{\bTd} \frac{\cos(x \cdot t)}{S(t)}\,dt
  \end{equation}
  if and only if $d \geq 3$, in which case the monopole is unique.
  \begin{proof}
    As in the proof of Lemma~\ref{thm:vx-on-Zd}, we use the Fourier transform to solve $\Lap w = -\gd_o$ by converting it into $S(t) \hat w(t) = -1$.
    This gives \eqref{eqn:monopoles-on-Zd}, and since $\frac{\cos t}{S(t)} \approx \frac1{S(t)} \in L^1(\bTd)$ for $t \approx 0$, the integral is finite iff $d \geq 3$ by the same argument as in the proof of Theorem~\ref{thm:finite-resistance-to-infinity-in-Zd}; see \eqref{eqn:exm:infinite-lattices:spherical-integral}.
    It remains to check that $w \in \HE$. Note that it follows from Theorem~\ref{thm:harmonics-on-Zd-are-linear} that the boundary term of \eqref{eqn:E(u,v)=<u_0,Lapv>+sum(normals)} vanishes, and hence we may compute the energy for $d \geq 3$ via
    \linenopax
    \begin{align}\label{eqn:energy-of-monopole-on-Zd}
      \|w\|_\energy &= \| \Lap^{1/2} w\|_2
      = \int_{\bTd} S(t) \hat w(t)^2\,dt
      = \int_{\bTd} \frac1{S(t)} \,dt
      < \iy.
    \end{align}

    Uniqueness is an immediate corollary of the previous theorem; if $w'$ were another, then $\Lap(w-w') = \gd_o - \gd_o = 0$ and $w-w'$ is constant by  Theorem~\ref{thm:harmonics-on-Zd-are-linear}.
  \end{proof}
\end{theorem}

\begin{remark}
  \version{}{\marginpar{FINISH THIS}}
  Upon comparing \eqref{eqn:monopoles-on-Zd} to \eqref{eqn:vx-on-Zd}, it is easy to see why all networks support finite-energy dipoles: the numerator in the integral for the monopole is of the order of 1 for $t \approx 0$, while the corresponding numerator for the dipole is $o(t)$ for $t \approx 0$.  
\end{remark}

\begin{theorem}\label{thm:R(x,y)-on-Zd}
  Resistance distance on the integer lattice $(\bZd,\one)$ is given by
  \linenopax
  \begin{equation}\label{eqn:R(x,y)-on-Zd}
    R(x,y) = \frac1{(2\gp)^d} \int_{\bT^d} \frac{\sin^2((x-y)\cdot \frac t2)} {\sum_{k=1}^d \sin^2\left(\frac{t_k}2\right)} \,dt.
  \end{equation}
  \begin{proof}
    Let $\{v_x\}_{x \in \bZd}$ be the potential configuration on \bZd. Then $v_x-v_y \in \Pot(x,y)$, so by \eqref{eqn:def:R(x,y)-energy} we may use \eqref{eqn:def:R^F(x,y)-Lap}to compute the resistance distance via $R(x,y) = v_x(x) + v_y(y) - v_x(y) - v_y(x)$, since $R^F=R^W$ on \bZd. Using $e_x = e^{\ii x \cdot t}$, substitute in the terms from \eqref{eqn:vx-on-Zd-as-exponential} of Lemma~\ref{thm:vx-on-Zd}:
    \linenopax
    \begin{align}
      R(x,y)
      &= \frac1{(2\gp)^d} \int_{\bTd} \frac{\cj{e_x}(e_x-1) + \cj{e_y}(e_y-1) - \cj{e_x}(e_y-1) - \cj{e_y}(e_x-1)}{S(t)} \,dt \notag \\
      &= \frac1{(2\gp)^d} \int_{\bTd} \frac{1-\cancel{\cj{e_x}} + 1 - \cancel{\cj{e_y}} - e_{y-x} + \cancel{\cj{e_x}} - \cj{e_{y-x}} + \cancel{\cj{e_y}}}{S(t)} \,dt \notag \\
      &= \frac1{(2\gp)^d} \int_{\bTd} \frac{2-2\cos((x-y) \cdot t)}{S(t)} \,dt, \label{eqn:vx-on-Zd-cos}
    \end{align}
    and the formula follows by the half-angle identity.
  \end{proof}
\end{theorem}

\begin{cor}\label{thm:neighbours-in-Zd-are-1/d-apart}
  If $y \nbr x$, then $R(x,y) = \frac1d$ on $(\bZd,\one)$.
  \begin{proof}
    The symmetry of $(\bZd,\one)$ indicates that the distance from $x$ to its neighbour will not depend on which of the $2d$ neighbours is chosen. For $k=1,2,\dots,d$, let $y_k$ be a neighbour of $x$ in the \kth direction. Then  \eqref{eqn:vx-on-Zd} gives
    \begin{equation}\label{eqn:nbr-distance-in-Zd}
      R(x,y_k)
      = \frac1{(2\gp)^d} \int_{\bT^d} \frac{\sin^2(\frac{t_k}2)} {\sum_{k=1}^d \sin^2\left(\frac{t_k}2\right)} \,dt.
    \end{equation}
    Thus, $\sum_{k=1}^d R(x,y_k) = 1$ and $R(x,y_k) = R(x,y_j)$ gives the result.
  \end{proof}
\end{cor}

\begin{theorem}\label{thm:finite-resistance-to-infinity-in-Zd}
  The metric space $((\bZd,\one),R)$ is bounded if and only if $d \geq 3$, in which case
  \linenopax
  \begin{align}\label{eqn:limR(x,y)-on-Zd}
    \lim_{y \to \iy} R(x,y)
    = \frac2{(2\gp)^d} \int_{\bT^d} \frac{1} {S(t)}\,dt
    \qq\text{for } d \geq 3.
  \end{align}
  \begin{proof}
    This result hinges upon the convergence properties of the integrand for $R(x,y)$ as computed in Lemma~\ref{thm:R(x,y)-on-Zd}. In particular, to see that $1/S(t) \in L^1(\bT^d)$ one only needs to check for $t \approx 0$, where
    \linenopax
    \begin{align*}
      \frac{1}{S(t)} = O\left(\frac{1}{\sum t_k^2}\right), \qq \text{as } t \to 0.
    \end{align*}
    Switching to spherical coordinates, $1/S(\gr) = O\left(\gr^{-2}\right)$, as $\gr \to 0$, and one requires
    \linenopax
    \begin{equation}\label{eqn:exm:infinite-lattices:spherical-integral}
      \int_0^1 |\gr^{-2}| \gr^{d-1} \, dS_{d-1} < \iy,
    \end{equation}
    where $dS_{d-1}$ is the usual $(d-1)$-\dimnl spherical measure. Of course, \eqref{eqn:exm:infinite-lattices:spherical-integral} holds precisely when $-2+d-1>-1$, i.e., when $d > 2$. Similarly, the function $\cos((x-y) \cdot t)/S(t) \in L^1(\bT^d)$ iff $d \geq 3$. Therefore, \eqref{eqn:vx-on-Zd-cos} gives
    \linenopax
    \begin{align*}
      R(x,y)
      = \frac2{(2\gp)^d} \int_{\bT^d} \frac{1} {S(t)}\,dt
       -\frac2{(2\gp)^d} \int_{\bT^d} \frac{\cos((x-y) \cdot t)}{S(t)} \,dt,
      \qq\text{for } d \geq 3.
    \end{align*}
    Now replace $y$ with a sequence of vertices tending to infinity as in Definition~\ref{def:path-to-infinity}. By the Riemann-Lebesgue lemma, the second integral vanishes and for any such $y \to \iy$, we have \eqref{eqn:limR(x,y)-on-Zd}. Note that this is independent of $x \in \verts$, as one would expect from the translational invariance of the network, since $\cond = \one$.
  \end{proof}
\end{theorem}

\begin{defn}\label{def:distance-to-infty}
  Denote $R_\iy := \lim_{y \to \iy} R(o,y)$, as it is clear from the previous result that the limit does not depend on the choice of $y$.
\end{defn}

\begin{cor}\label{thm:proximity-of-infinity}
  For $d \geq 3$, there exists $x \in \bZd$ for which $R(o,x) > R_\iy$.
  \begin{proof}
     From \eqref{eqn:R(x,y)-on-Zd}, it is clear that $R(o,x) \leq R(o,\iy)$ if and only if
    \linenopax
    \begin{align*}
      \frac1{(2\gp)^d} \int_{\bTd} \frac{1-\cos(x \cdot t)}{S(t)}\,dt
      \leq \frac1{(2\gp)^d} \int_{\bTd} \frac{1}{S(t)}\,dt,
    \end{align*}
    which is equivalent to
    \linenopax
    \begin{align}\label{eqn:black-hole}
      \frac1{(2\gp)^d} \int_{\bTd} \frac{\cos(x \cdot t)}{S(t)} \geq 0.
    \end{align}
    However, \eqref{eqn:black-hole} cannot hold for all $x \in \bZd$, as such an inequality would mean that all Fourier coefficients of $w$ are nonnegative, in violation of Heisenberg's uncertainty principle.
  \end{proof}
\end{cor}

\begin{remark}\label{rem:black-hole}
  Corollary~\ref{thm:proximity-of-infinity} leads to the paradoxical conclusion that given $x \in \verts$, there may be a $y$ which is ``further from $x$ than infinity''. This is the case for $d=3$; numerical computation of \eqref{eqn:limR(x,y)-on-Zd} gives
  \begin{equation}\label{eqn:R_iy-in-Z3}
    \lim_{y \to \iy} R(x,y) \approx 0.5054620038965394,
    \qq \text{in } \bZ^3,
  \end{equation}
  and for $y=(1,1,1)$,
  \linenopax
  \begin{equation}\label{eqn:R(x,y)>R(x,infty)}
    R\left(o,y\right) \approx 0.5334159062457338.
  \end{equation}
  In fact, numerical computations indicate the following extremely bizarre situation:
  \linenopax
  \begin{align*}
    R(x,y_{2k}) < R(x,\iy) < R(x,y_{2k+1}), \qq \text{for } y_n := (n,n,0).
  \end{align*}
\end{remark}

\begin{figure}
  \centering
  \scalebox{0.70}{\includegraphics{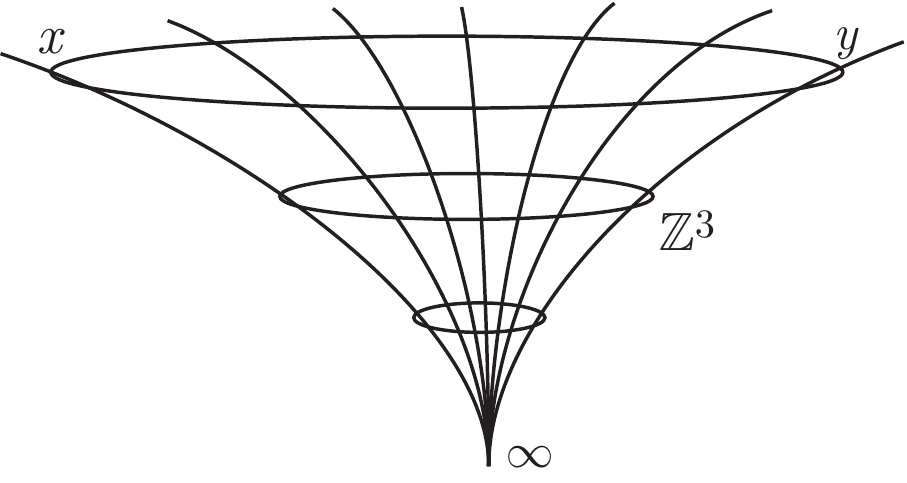}}
  \caption{\captionsize In $\bZ^3$, it may happen that $R(x,y) > R(x,\iy)$, where $R(x,\iy) = \lim_{z \to \iy} R(x,z)$. This phenomenon is represented here schematically as a ``black hole''.}
  \label{fig:Z3-black-hole}
\end{figure}

\begin{remark}\label{rem:Bochner-on-Zd}
  An application of Bochner's Theorem (see Theorem~\ref{thm:Bochner's-theorem}) yields a unique Radon probability measure \prob on \bTd such that
  \linenopax
  \begin{align*}
    \int_{\bTd} e^{\ii t \cdot x} \,d\prob(t) = e^{-\frac12 R(o,x)},
    \q \forall x \in \bZd.
  \end{align*}
\end{remark}



\begin{cor}\label{thm:vx-in-L2(Zd)-for-big-d}
  For $(\bZd,\one)$, $v_x \in \ell^2(\bZd)$ if and only if $d \geq 3$.
  \begin{proof}
    By computations similar to those in the proof of Theorem~\ref{thm:finite-resistance-to-infinity-in-Zd}, one can see that in absolute values, the integrand $\left|(e^{\ii x\cdot t} - 1)/S(t)\right|$ of \eqref{eqn:vx-on-Zd-as-exponential} is in $L^2(\bTd)$ if and only if $d \geq 3$, in which case Parseval's theorem applies.
  \end{proof}
\end{cor}

\begin{cor}\label{thm:w-in-L2(Zd)-for-big-d}
  For $(\bZd,\one)$, the monopole $w \in \ell^2(\bZd)$ if and only if $d \geq 5$.
  \begin{proof}
    The proof is almost identical to that of Corollary~\ref{thm:vx-in-L2(Zd)-for-big-d}, except that the integrand is $1/S(t)$, which is in $L^2(\bTd)$ if and only if $d \geq 5$.
  \end{proof}
\end{cor}

\begin{exm}\label{exm:Z-not-bounded}
  To see why $R$ is not bounded on $(\bZ,\one)$, one can evaluate \eqref{eqn:vx-on-Zd} explicitly via the Fej\'{e}r kernel:
  \linenopax
  \begin{align*}
    R(x,y)
    &= \frac1{(2\gp)} \int_{\bT} \frac{\sin^2((x-y)\frac t2)} {\sin^2\left(\frac{t}2\right)} \,dt \\
    &= \frac1{(2\gp)} \int_{\bT} |x-y|\frac{\sin^2((x-y)\frac t2)} {|x-y| \sin^2\left(\frac{t}2\right)} \,dt \\
    &= |x-y| \frac1{(2\gp)} \int_{\bT} F_N(t) \,dt \\
    &= |x-y|.
  \end{align*}
  where $F_N(t)$ is the Fej\'{e}r kernel with $N=|x-y|$; see Figure~\ref{fig:fejer-in-Z1}. Of course, this was to be expected because $R$ coincides with shortest path metric on trees.
  \begin{figure}
    \centering
    \scalebox{0.90}{\includegraphics{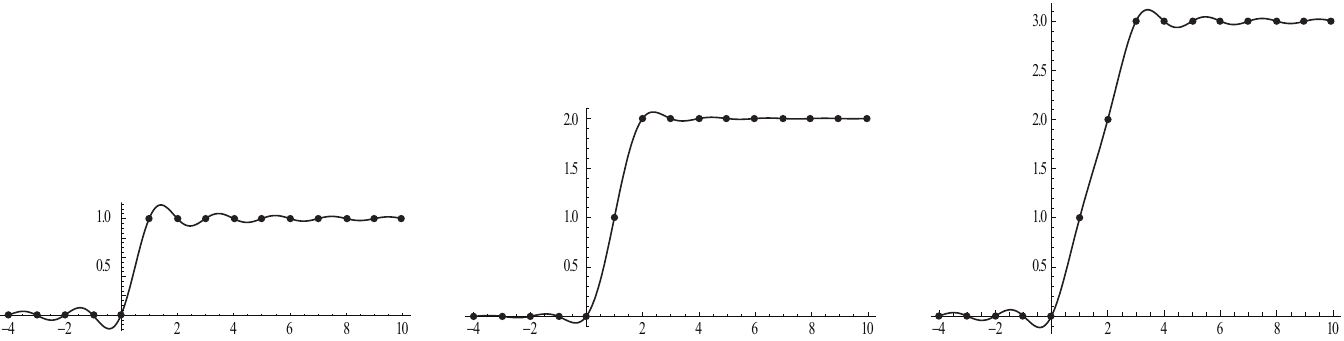}}
    \caption{\captionsize The function $v_x$, for $x=1,2,3$ in \bZ, as obtained from the Fej\'{e}r kernel. See Example~\ref{exm:Z-not-bounded}.}
    \label{fig:fejer-in-Z1}
  \end{figure}
\end{exm}

The following result is well known; we include it for completeness and the novelty of the proof.
\begin{theorem}\label{thm:harmonics-on-Zd-are-linear}
  $h$ is a harmonic function on $(\bZd,\one)$ if and only if $h$ is linear (or affine). Consequently, $\HE = \clspn{\Fin}$ for \bZd.
  \begin{proof}
    \version{}{\marginpar{Use proof from Rem.~\ref{rem:HE(Zd)-lies-in-L2(Zd)}.}}
    From $\Lap h = 0$, the Fourier transform gives $S(t) h(t)=0$. By the formula of Lemma~\ref{thm:Fourier-transform-of-Lap-on-Zd} for $S(t)$, this means $\hat h$ can only be supported at $t=0$ and hence that $\hat h$ is a distribution which is a linear combination of derivatives of the Dirac mass at $t=0$; see \cite{Rud91} for this structure theorem from the theory of distributions.

    Denoting this by $\hat h(t) = P(\gd_o)$, where $\gd_o$ is the Dirac mass at $t=0$, and $P$ is some polynomial. The inverse Fourier transform gives $h(x) = P(x)$. If the degree of $P$ is 2 or higher, then $\Lap h$ will have a constant term of $-2d$ (cf.~\eqref{eqn:Fourier-transform-of-Lap-on-Zd}) and hence cannot vanish identically.

    It is clear that a linear function on \bZd has infinite energy; consequently \Harm is empty on this network and the second conclusion follows.
  \end{proof}
\end{theorem}

\begin{exm}[Nontrivial harmonic functions on $\bZd \cup \bZd$]
  \label{exm:Nontrivial-harmonic-functions-on-2Zd}
  Consider the disjoint union of two copies of \bZd, with $\cond = \one$ and $d \geq 3$. Now connect the origins $o_1,o_2$ of the two lattices with a single edge of conductance $\cond_{o_1o_2} = 1$. Let $w_1 \in \HE$ be a monopole on the first copy of \bZd, as ensured by Theorem~\ref{thm:monopoles-on-Zd}. We may assume $w_1$ is normalized so that $w(o_1)=1$, and then extend $w_1$ to the rest of the network by letting $\tilde w_1(x) = 0$ for all $x$ in the second copy of \bZd. Similarly, let $w_2$ be a function which is a monopole on the second \bZd, satisfies $w(o_2)=1$, and extend it to $\tilde w_2$ by defining $\tilde w_2(x) = 0$ for $x$ in the first copy of \bZd. Now one can check that $\Lap \tilde w_1 = \Lap \tilde w_2 = -\gd_{o_2}$. Note that the unit drop in $\tilde w_1$ across the edge $\cond_{o_1o_2}$ moves the Dirac mass of $\Lap w_1$ to the second copy of \bZd. Now define 
  \linenopax
  \begin{align}\label{eqn:harmonic-functions-on-2Zd}
    h := w_1 - w_2.
  \end{align}
  It is easy to check that $h \in \HE$ and that $h \in \Harm$.
\end{exm}

\begin{cor}\label{thm:monopole-on-Zd-vanishes-at-infty}
\label{thm:monopole-on-Zd-as-resistance-distance}
  If $w$ is the monopole on $(\bZd,\one)$, $d \geq 3$, then
  \begin{equation}\label{eqn:monopole-on-Zd-as-resistance-distance}
    w(x) = \tfrac12 (R(o,x) - R(o,\iy)),
  \end{equation}
  and consequently $\lim_{x \to \iy} w(x) =0$.
  \begin{proof}
    Subtract \eqref{eqn:limR(x,y)-on-Zd} from \eqref{eqn:vx-on-Zd-cos} and compare to \eqref{eqn:monopoles-on-Zd}. For the latter statement, one can take the limit of \eqref{eqn:monopole-on-Zd-as-resistance-distance} as $x \to \iy$ directly or apply the Riemann-Lebesgue lemma to \eqref{eqn:monopoles-on-Zd}.
  \end{proof}
\end{cor}

\begin{cor}\label{thm:energy-of-monopole-on-Zd}
  If $w$ is the monopole on $(\bZd,\one)$, then 
  \begin{equation}\label{eqn:energy-of-monopole-on-Zd}
    \energy(w) = \frac12 \lim_{y \to \iy} R(x,y).
  \end{equation}
  \begin{proof}
    Compare \eqref{eqn:energy-of-monopole-on-Zd} to \eqref{eqn:limR(x,y)-on-Zd} and note that $\Harm=\{0\}$, so $w$ is unique.
  \end{proof}
\end{cor}


\begin{remark}\label{rem:HE(Zd)-lies-in-L2(Zd)}
  For $(\bZd,\one)$, it is instructive to work out directly why \Fin is dense in \HE. That is, let us suppose that the boundary term vanishes for every $v \in \HE$, and use this to prove that every function which is orthogonal to \Fin must be constant (and hence 0 in \HE). This shows that \Fin is dense in \HE in the energy norm.

  \begin{proof}[``Proof'']
    If $v \in \HE$, then $\|v\|_\energy = \la v, \Lap v\ra_\cond < \iy$, the Fourier transform sends $v \mapsto \hat v(t) = \sum_\bZ v_n e^{\ii n \cdot t}$ and
    \linenopax
    \begin{align}\label{eqn:thm:HE(Zd)-lies-in-L2(Zd):orthog-to-Dirac}
      \la v, \Lap v\ra_\cond
      \mapsto (2\gp)^{-d} \int_{\bTd} \hat v(t) S(t) \hat v(t) \,dt
      <\iy,
    \end{align}
    where $S(t) = 4\sum_{k=1}^d \sin^2\left(\tfrac{t_k}2\right)$, as in Lemma~\ref{thm:Fourier-transform-of-Lap-on-Zd}.
    Then note that the Schwarz inequality gives
    \linenopax
    \begin{align*}
      \left(\int_{\bTd} S(t) \hat v(t) \,dt\right)^2
      \leq \int_{\bTd} S(t) \,dt \int_{\bTd} S(t) \hat v(t)^2 \,dt,
    \end{align*}
    so that $S(t) \hat v(t) \in L^1(\bTd)$.
    From the other hypothesis, $v \perp \Fin$ means that $\la \gd_x, v \ra_\energy = 0$ for each $x \in \verts$, whence Parseval's equation gives
    \linenopax
    \begin{align*}
      0 &= \la \gd_{x_m}, v \ra_\energy
      = \la \gd_{x_m}, \Lap v \ra_\cond
      \mapsto (2\gp)^{-d} \int_{\bTd} e^{\ii m \cdot t} S(t) \hat v(t) \,dt
      = 0, \forall m.
    \end{align*}
    This implies that $S(t) \hat v(t) =0$ in $L^1(\bTd)$, and hence $\hat v$ can only be supported at $t=0$. From Schwartz's theory of distributions, this means
    \linenopax
    \begin{align*}
      \hat v(t) = f_0(t) + c_0 \gd_0(t) + c_1 D \gd_0(t) + c_2 D^{(2)} \gd_0(t) + \dots,
    \end{align*}
    where $f_0$ is an $L^1$ function and all the other terms are derivatives of the Dirac mass at $t=0$ ($D^{(2)}$ is a differential operator of rank 2, etc.).

    If $\hat v$ is just a function, then it is 0 a.e. and we are done. If the distribution $\gd_0(t)$ is a component of $\hat v$, then $\sF^{-1}(\gd_0) = \one$, which is zero in \HE. In one dimension, the distribution $\gd_0'(t)$ cannot be a component of $\hat v$ because  $\sF^{-1}(\gd_0')(x_m) = m$, and this function does not have finite energy (the computation of the energy picks up a term of 1 on every edge of the lattice \bZd). The computation is similar for higher derivatives of $\gd_0$, but they diverge even faster. For higher dimensions, note that $D_1 \gd_{(0,0)} = D \gd_0 \otimes \gd_0$ and $\energy(D \gd_0 \otimes \gd_0) = \energy(D \gd_0) \energy(\gd_0)$ (this is a basic fact about quadratic forms on a Hilbert space), and so this devolves into same argument as for the 1-\dimnl case.
  \end{proof}
\end{remark}

Remark~\ref{rem:HE(Zd)-lies-in-L2(Zd)} does not hold for general graphs; see Example~\ref{exm:binary-tree:nontrivial-harmonic}. Also, the end of the proof shows why $\tilde{\Lap}^{1/2} (\widehat{v+k}) = \tilde{\Lap}^{1/2} \hat v$, as mentioned in Remark~\ref{rem:LapTilde-doesn't-see-constants}; the addition of a constant corresponds on the Fourier side to the addition of a Dirac mass outside the support of \gc.

The case of $(\bZd, \one)$, for $d=1$ is a tree and hence very simple with $R(x,y) = |x-y|$, and for $d \geq 3$ may be fairly well understood by the formulas given above. However, the case $d=2$ seems to remain a bit mysterious. It appears that $R(x,y) \approx \log(1+|x-y|)$; we now give two results in this direction.

\begin{remark}\label{rem:asymptotic-neighobur-difference-R(x,y)}
  From Theorem~\ref{thm:finite-resistance-to-infinity-in-Zd} it is clear that for $d \geq 3$, if $y_n \nbr z_n$ and both tend to \iy, one has $\lim_{n \to \iy} (R(x,y_n)-R(x,z_n)) =0$. In fact, this remains true in $\bZ^2$ but not \bZ. For \bZ,
  \linenopax
  \begin{align*}
    y_n \nbr z_n \q\implies\q |R(x,y_n)-R(x,z_n)|=1 \limas{n} 1 \neq 0.
  \end{align*}
  A little more work is required for $\bZ^2$, where we work with $x=o$ for simplicity: 
  \linenopax
  \begin{align*}
    R(o,z) - R(o,y)
    &= \frac1{(2\gp)^2} \int_{\bT^2} \frac{\cos(y \cdot t) - \cos(y \cdot t + t_k)} {S(t)} \,dt \\
    &= \frac1{(2\gp)^2} \int_{\bT^2} \frac{\cos(y \cdot t)(1-\cos t_k) + \sin(y \cdot t) \sin t_k} {S(t)} \,dt.
  \end{align*}
  One can check that $\frac{1-\cos t_k}{S(t)},\frac{\sin t_k}{S(t)} \in L^1(\bT^2)$ by converting to spherical coordinates and making the estimate
  \linenopax
  \begin{align*}
    \int_{\bT^2} \frac{\gr}{\gr^2} \gr \, d\gr \, d\gq < \iy.
  \end{align*}
  Now the Riemann-Lebesgue Lemma shows that $R(o,y) - R(o,z)$ tends to 0 as $y$ (and hence also $z$) tends to \iy.
\end{remark}

\begin{theorem}\label{thm:gradR(x,y)-vanishes-at-infty}
  On $(\bZ^2, \one)$, the gradient of $R$ vanishes at infinity, i.e.,
  \begin{equation}\label{eqn:gradR(x,y)-vanishes-at-infty}
    \lim_{y \to \iy} \nabla R(x,y) = 0.
  \end{equation}
\end{theorem}

\begin{theorem}\label{thm:R(x,y)-in-Z2-as-contour-integral}
  On $(\bZ^2, \one)$, resistance distance is given by
  \begin{equation}\label{eqn:R(x,y)-in-Z2-as-contour-integral}
    R(x,y) = \int_0^x A(o,t) \,dt + \int_0^\iy B(x,t) \,dt,
  \end{equation}
  \begin{equation}\label{eqn:A1-and-A2}
    \text{where}\q
    A(s,t) := \frac{\del}{\del t_1}
    B(s,t) := \frac{\del}{\del t_2}
  \end{equation}
\end{theorem}

\section{Noncompactness of the transfer operator}

\begin{exm}[\Trans may not be the uniform limit of finite-rank operators]
  \label{exm:nonuniformly-converging-transfer-operator}
  Let \Graph be the integers \bZ with edges only between vertices of distance one apart (as in Example~\ref{exm:infinite-lattices} with $d=1$), with $\cond \equiv 1$. Then the transfer operator $T := \gs^+ + \gs^-$ consists of the sum of two unilateral shifts, for which the finite truncations (as described just above) are the banded matrices
  \begin{equation}\label{eqn:exm:nonuniformly-converging-transfer-operator}
    \Trans_N =
    \left[\begin{array}{rrrrrrrrr}
      0 &1 &0 &0 &\dots &0 &0 &0 \\
      1 &0 &1 &0 &\dots &0 &0 &0 \\
      0 &1 &0 &1 &\dots &0 &0 &0 \\
      0 &0 &1 &0 &\dots &0 &0 &0 \\
      \vdots &\vdots &\vdots &\vdots & \ddots &\vdots &\vdots &\vdots \\
      0 &0 &0 &0 &\dots &0 &1 &0 \\
      0 &0 &0 &0 &\dots &1 &0 &1 \\
      0 &0 &0 &0 &\dots &0 &1 &0
    \end{array}\right].
  \end{equation}
  Then consider the vectors
  \begin{equation}\label{eqn:exm:nonuniformly-converging-transfer-operator:vectors}
    \gx_n := (\underbrace{0,\dots,0}_{n\text{ zeros}}, 1, \tfrac12, \tfrac13, \tfrac14, \dots), \qq
    \|\gx_n\|^2 = 2\sum_{k=1}^\iy \frac1{k^2} = \frac{\gp^2}3.
  \end{equation}
  Then $\Trans_N$ does not converge to \Trans uniformly, because for $n=N$,
  \linenopax
  \begin{align*}
    \la \gx_n, (\Trans - \Trans_n)\gx_n\ra_\cond
    = \sum_{|k|>n} \gx_k \gx_{k+1}
    &= \sum_{|k|>n} \frac1{(k-n)(k-n+1)}
    = \sum_{|k| \geq 1} \frac1{k(k+1)} \\
    &\geq \sum_{|k| \geq 1} \left(\frac1{k+1}\right)^2
    \approx \frac{\gp^2}6,
  \end{align*}
  which is bounded away from 0 as $n \to \iy$.
\end{exm}

\subsection{The Paley-Wiener space $H_s$}
The transfer operator is not compact in \HE, as Example~\ref{exm:nonuniformly-converging-transfer-operator} shows. However, by introducing the correct weights we can obtain a compact operator, i.e., the transfer operator is compact when considered as acting on the correct Hilbert space. To this end, we make the identification between $\gx \in \ell^2(\bZ)$ and $f(z) = \sum_{n \in \bZ} \gx_n z^n \in L^2(\bT)$ via Fourier series, so that we may use analytic continuation and introduce the following spaces.

\begin{defn}\label{def:s-norm}
  For an function $f \in L^2(\bT)$ given by $f(z) = \sum_{n \in \bZ} \gx_n z^n$, we define
  \begin{equation}\label{eqn:def:s-norm}
    \|f\|_s := \left(\sum_{n \in \bZ} e^{s|n|} |\gx_n|^2 \right)^2,
  \end{equation}
  and consider the space.
  \begin{equation}\label{eqn:Hs-space}
    H_s := \{f : \bT \to \bT \suth \|f\|_s < \iy\}.
  \end{equation}
  For $s=0$ we recover good old $\ell^2(\cond)$, but for $s>0$, we have the subspace of $\ell^2(\cond)$ which consists of those functions with an analytic continuation to the annulus $\{z \suth 1 - s < |z| < 1 + s\}$ about \bT. In general, we have $H_s \ci L^2(\bT) \ci H_{-s} = H_s^\ast$.
\end{defn}

\begin{theorem}\label{thm:Trans-is-compact-on-Hs}
  The transfer operator is a compact operator on $H_s$.
  \begin{proof}
    Using $\Lap = \cond-\Trans$, we show that there exist solutions to $\Lap v = \gd_\ga - \gd_\gw$ by construction, using spectral theory. The Laplacian may be represented as the infinite symmetric banded matrix
    \linenopax
    \begin{align*}
      \left[\begin{array}{cccccc}
        \dots &\cond(x_1) &-\cond(x_1, y_1) &\dots \\
        \dots &-\cond(x_1, y_1) &\cond(x_2) &-\cond(x_2, y_2) &\dots \\
        \dots & &-\cond(x_2, y_2) &\cond(x_3) &-\cond(x_3, y_3) &\dots \\
        \dots & & &\ddots &\ddots &\ddots
      \end{array}\right]
    \end{align*}
    The symmetry is immediate from the symmetry of $\cond_{xy}$, of course.

    Using the same notation as in Example~\ref{exm:nonuniformly-converging-transfer-operator}, we must check that $\Trans_N \to \Trans$ uniformly in $H_s$, so we examine $D_N := \Trans - \Trans_N$:
\version{}{\marginpar{How to justify \eqref{eqn:thm:Trans-is-compact-on-Hs:deriv1}? This is supposed to be for any old \gx, right?}}
    \linenopax
    \begin{align}
      \la \gx, D_N \gx \ra
      &= \sum_{|n| \geq N} e^{s|n|} \cj{\gx_n} \gx_{n+1} \notag \\
      &\leq \left(\sum_{|n| \geq N} e^{s|n|} |\gx_{n}|^2\right)^{1/2}
        \left(\sum_{|n| \geq N} e^{s|n|} |\gx_{n+1}|^2\right)^{1/2} \notag \\
      &= \left(\sum_{|n| \geq 0} e^{s|n+N|} |\gx_{n+N}|^2\right)^{1/2}
        \left(\sum_{|n| \geq 0} e^{s|n+N|} |\gx_{n+N+1}|^2\right)^{1/2} \notag \\
      &= e^{-s/2} \left(e^{sN} |\gx_N|^2 + \sum_{|n| \geq N+1} e^{sn} |\gx_{n}|^2\right)^{1/2}
        \left(\sum_{|n| \geq N+1} e^{sn} |\gx_{n}|^2\right)^{1/2}
        \label{eqn:thm:Trans-is-compact-on-Hs:deriv1} \\
      \frac{\la \gx, D_N \gx\ra}{\|\gx\|_s}
      &= e^{-s/2} \frac{\left(\sum_{|n| \geq N+1} e^{sn} |\gx_{n}|^2\right)^{1/2}}{\left(e^{sN} |\gx_N|^2 + \sum_{|n| \geq N+1} e^{sn} |\gx_{n}|^2\right)^{-1/2}}
      \limas{N}&\; 0. \notag
    \end{align}
\version{}{\marginpar{The hat here indicates passing to the Calkin Alg, i.e, modding by compacts, right?}}
    Then \Trans is compact in $H_s$, and in fact, \Trans is trace class! Then $\hat{cpt} = 0$, so $\Lap = \id - \Trans$ implies $\hat{\Lap} = \hat{\Cond}$.
  \end{proof}
\end{theorem}

\section{Non-simple integer lattice networks}
\label{sec:nonsimple-integer-lattices}

In this section, we illustrate some of the phenomena that may occur on integer lattices when the conductances are allowed to vary. 
\version{}{\marginpar{Add a decay example for $\lim \Graph$.}}
Many of these examples serve to demonstrate certain definitions or general properties discussed in previous sections.

\begin{exm}[Symmetry of the graph vs. symmetry of the network]
  \label{exm:2-dimnl-integer-lattice}
  Consider the 2-\dimnl integer lattice; the case $d=2$ in Example~\ref{exm:infinite-lattices}, and think of these points as living in the complex plane, so each vertex is $m+n\ii$, where $m$ and $n$ are integers and $\ii=\sqrt{-1}$. It is possible to define the conductances in such a way that a function $v(z)$ has finite energy, but $v(\ii z)$ does not (this is just precomposing with a symmetry of the graph: rotation by by $\frac\gp2$). However, $v(z)$ is in $\ell^2(\one)$ if and only if $v(\ii z)$ is in $\ell^2(\one)$. Thus, $\ell^2(\one)$ does not see the graph.

  Define the conductances by
  \linenopax
  \begin{align*}
    \cond_{xy} =
    \begin{cases}
      1, & y=x+1, \\
      2^{|\Im(y)|}, & y=x+\ii,
    \end{cases}
  \end{align*}
  so that the conductances of horizontal edges are all 1 and the conductances of vertical edges grow like $2^k$. Now consider the function
  \linenopax
  \begin{align*}
    v(z) = v(x+\ii y) =
    \begin{cases}
      2^{-|Re(x)|}, & y=0, \\
      0, & y \neq 0,
    \end{cases}
  \end{align*}
  When computing the energy $\energy(v)$, the only contributing terms are the edges along the real axis, and the edges immediately adjacent to the real axis:
  \linenopax
  \begin{align*}
    \energy(v(z))
    &= \mathrm{horizontal} + \mathrm{vertical} \\
    &= 2(1/2 + 1/4 + 1/8 + ...) + 4(1/2 + 1/4 + 1/8 + ...)
    = 6,
  \end{align*}
  which is finite. However,
  \begin{align*}
    \energy(v(\ii z))
    = 2(1 + 1 + 1 + ...) + 4(1/2 + 1/4 + 1/8 + ...) = \iy.
  \end{align*}
\end{exm}

\begin{exm}[An example where $\ell^2 \nsubseteq \HE$]
  \label{exm:1-dimnl-integer-lattice}
  Let $\bZ$ have $\cond_{n-1,n} = n$. Consider $\ell^2(\verts,\gn)$ where \gn is the counting measure. The Dirac functions $\gd_{x_k}$ satisfy $\|\gd_{x_k}\|_\unwt = 1$, so $\{\gd_{x_k}\}$ is a bounded \seq in $\ell^2(\unwtdsp)$. However, the Laplacian is
  \linenopax
  \begin{align*}
    \Lap_\unwt = \left[
    \begin{array}{rrr}
      \ddots \\ -n & 2n+1 & -(n+1) \\ &&\ddots
    \end{array}\right]
  \end{align*}
  and $\energy(\gd_{x_k}) = \la \gd_{x_k}, \Lap \gd_{x_k}\ra = 2k+1 \limas{k} \iy$. So we cannot have the bound $\|v\|_\energy \leq K \|v\|_\unwt$, for any constant $K$.

  \version{}{\marginpar{Why doesn't unitary equivalence provide an $\ell^2_\cond$ counterexample?}}
  This is ``corrected'' by using the measure \cond instead. In this case, $\|\gd_k\|_\cond = 2k+1$ so that $\{\gd_k\}$ is not bounded and we must use $\{\gd_k/\sqrt{\cond(k)}\}$. But then the Laplacian is
  \linenopax
  \begin{align*}
    \Lap_\cond = \left[
    \begin{array}{rrr}
      \ddots \\ -\tfrac{n}{2n+1} & 1 & -\tfrac{n+1}{2n+1} \\ &&\ddots
    \end{array}\right]
  \end{align*}
  and $\energy\left(\frac{\gd_{x_k}}{\sqrt{\cond(x_k)}}\right) = \frac1{\cond(x_k)} \energy(\gd_x) = 1$.
\end{exm}

\begin{exm}\label{exm:Unbounded-functions-of-finite-energy}
  It is quite possible to have unbounded functions of finite energy. Consider the network $(\Graph, \cond) = (\bZ,\one)$ with vertices at each integer and unit conductances to nearest nearest neighbours. Then it is simple to show that $u(n) = \sum_{i=1}^n \frac1n$ and $v(n) = \log |1+n|$ are unbounded and have finite energy --- use the identity $\sum_{n=1}^\iy \frac1{n^2}=\frac{\gp^2}6$. For $v$, note that $\log|1+n| - \log|1+(n-1)| = \log\left|1+\frac1n\right| \leq \frac1n$.
\end{exm}

\begin{exm}[An unbounded function with finite energy]
  \label{exm:unbounded-finite-energy-fn}
  Let $\bZ$ have $\cond_{n-1,n} = \frac1{n^2}$. Then the function $f(n)=n$ is clearly unbounded, but
  \linenopax
  \[\energy(f)
    = \sum \frac1{n^2}(f(n)-f(n-1))^2
    = \sum \frac1{n^2}
    = \frac{\gp^2}3 < \iy.\]
  Conclusion: it is possible to have unbounded functions of finite energy if \cond decays sufficiently fast. 

\version{}{\marginpar{tree example}  It is natural to wonder if there exist unbounded harmonic functions of finite energy. We do not know the answer but have so far failed to produce any. The equation
  \linenopax
  \begin{align*}
    \energy(h) = \sum_{\bd \Graph} h \dn h < \iy
  \end{align*}
  leads us to expect it is not possible.}
\end{exm}

  We've seen that there are no harmonic functions of finite energy on $(\bZd,\cond)$, when $\cond = \one$. However, the situation is very different when \cond is not bounded.
  
\begin{theorem}\label{thm:harmonic-functions-on-summable-resistance-integers}
  $\Harm \neq 0$ for $(\bZ,\cond)$ iff $\sum \cond_{xy}^{-1} < \iy$. In this case, \Harm is spanned by a single bounded function. 
  \begin{proof}
    \fwd Fix $u(0)=0$, define $u(1) = \frac1{\cond_{01}}$ and let $u(n)$ be such that
    \linenopax
    \begin{align}\label{eqn:harmonic-recipe-on-Z1}
      u(n) - u(n-1) = \frac1{\cond_{n-1,n}}, 
      \q\forall n.
    \end{align}
    Now $u$ is harmonic:
    \linenopax
    \begin{align*}
      \Lap u(n) 
      &= \cond_{n-1,n}(u(n)-u(n-1)) - \cond_{n,n+1}(u(n+1)-u(n)) \\
      &= \cond_{n-1,n}\frac1{\cond_{n-1,n}} - \cond_{n,n+1}\frac1{\cond_{n,n+1}}
       = 0, 
    \end{align*}
    and $u$ is of finite energy
    \linenopax
    \begin{align*}
      \energy(u)
      &= \sum_{n \in \bZ} \cond_{n-1,n}(u(n)-u(n-1))^2
       = \sum_{n \in \bZ} \frac1{\cond_{n-1,n}}
       < \iy.
    \end{align*}
    Note that once the values of $u(0)$ and $u(1)$ are fixed, all the other values of $u(n)$ are determined by \eqref{eqn:harmonic-recipe-on-Z1}. Therefore, \Harm is 1-dimensional.
    
    \bwd If $\Lap u(n) = \cond_{n-1,n}(u(n)-u(n-1)) - \cond_{n,n+1}(u(n+1)-u(n)) = 0$ for every $n$, then
    \linenopax
    \begin{align*}
      \cond_{n-1,n}(u(n)-u(n-1)) = \cond_{n,n+1}(u(n+1)-u(n)) = a,
    \end{align*}
    for some fixed $a$ (the amperage of a sourceless current). Then
    \linenopax
    \begin{align}\label{eqn:finite-energy-harmonic-implies-bounded-on-Z}
      \energy(u) 
      = \sum_{n \in \bZ} \cond_{n-1,n}(u(n)-u(n-1))^2
      = a^2 \sum_{n \in \bZ} \frac{1}{\cond_{n-1,n}}
      < \iy,
    \end{align}
    since $u \in \Harm \ci \HE$. Note that \eqref{eqn:finite-energy-harmonic-implies-bounded-on-Z} implies $u$ is bounded: $\energy(u) = a \sum_{n \in \bZ}(u(n)-u(n-1))$ and $\sum_{n \geq 1} (u(n)-u(n-1)) = \lim_{n \to \iy} u(n) - u(0)$. The function $u$ is monotonic because it is harmonic, so the sum is absolutely convergent.
  \end{proof}
\end{theorem}

We will now explore a specific example of this kind of network, where explicit computations are tractable.

  \begin{figure}
    \centering
    \scalebox{1.0}{\includegraphics{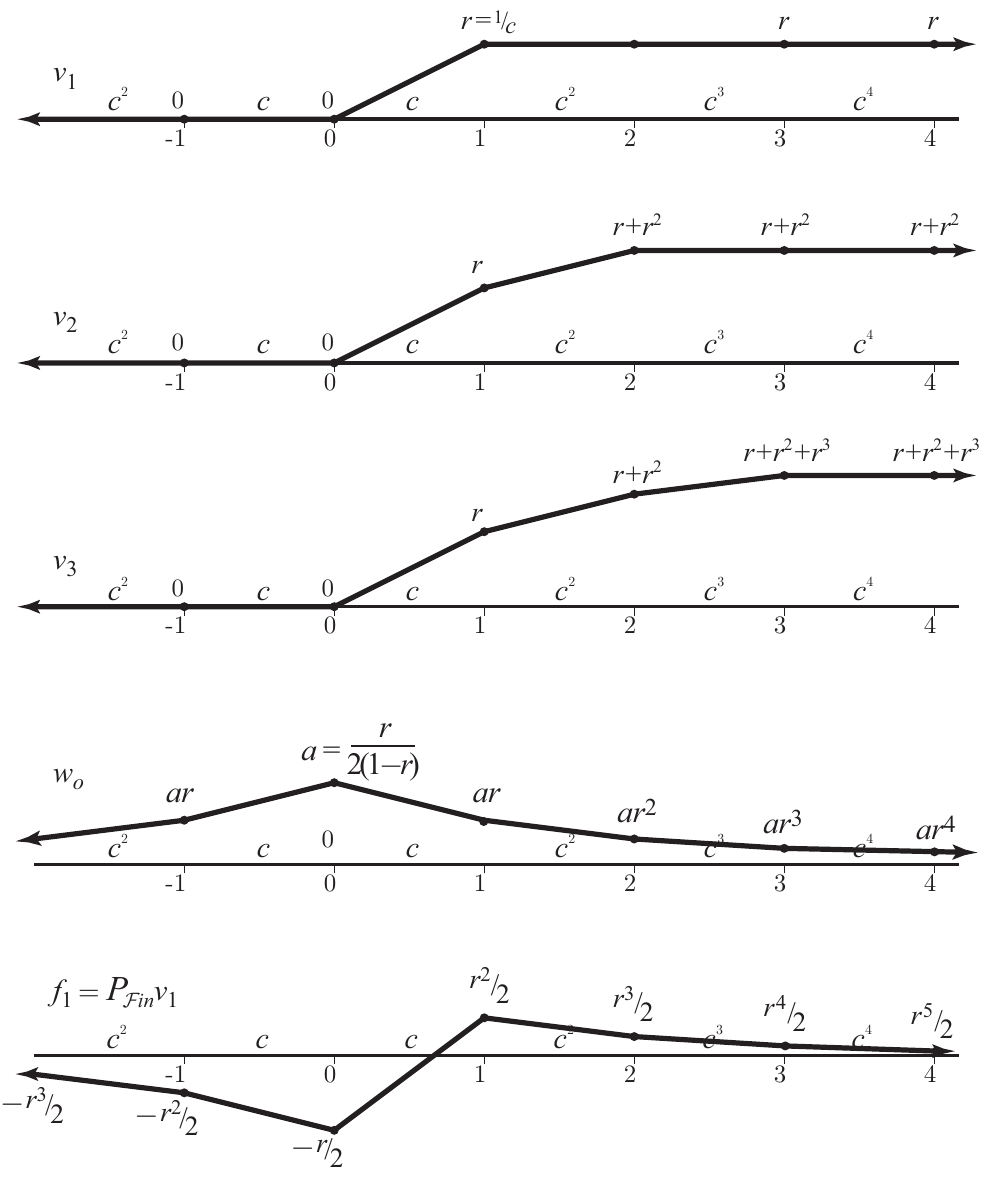}}
    \caption{\captionsize The functions $v_1$, $v_2$, and $v_3$ on $(\bZ,c^n)$. Also, the monopole $w_o$ and the projection $f_1 = \Pfin v_1$. See Lemma~\ref{thm:repkernels-on-geometric-integers}.}
    \label{fig:vx-on-Zc}
  \end{figure}
  
\begin{exm}\label{exm:geometric-integers}\label{def:geometric-integers}
  For a fixed constant $c>1$, let $(\bZ,c^n)$ denote the network with integers for vertices, and with geometrically increasing conductances defined by 
  \linenopax
    \begin{align*}
    \cond_{n-1,n} = c^{\max\{|n|,|n-1|\}}, 
  \end{align*}
  so that the network under consideration is
  \linenopax
    \begin{align*}
    \xymatrix{
      \dots \ar@{-}[r]^{c^3}
      & -2 \ar@{-}[r]^{c^2} 
      & -1 \ar@{-}[r]^{c} 
      & 0 \ar@{-}[r]^{c} 
      & 1 \ar@{-}[r]^{c^2} 
      & 2 \ar@{-}[r]^{c^3} 
      & 3 \ar@{-}[r]^{c^4} 
      & \dots
    }
  \end{align*} 
  We fix $o=0$.
\end{exm}

\begin{lemma}\label{thm:monopole-on-geometric-integers}\label{thm:repkernels-on-geometric-integers}
  On $(\bZ,c^n)$, the energy kernel is given by
  \linenopax
  \begin{align*}
    v_n(k) = 
    \begin{cases}
      0, &k \leq 0, \\
      \frac{1-r^{k+1}}{1-r}, &1 \leq k \leq n, \\
      \frac{1-r^{n+1}}{1-r}, &k \geq n,
    \end{cases}
    n > 0,
  \end{align*}
  and similarly for $n < 0$.
  Furthermore, the function $w_o(n) = ar^{|n|}$, $a:= \frac{r}{2(1-r)}$, defines a monopole, and $h(n) = \operatorname{sgn}(n) (1-w_o(n))$ defines an element of \Harm.
  \begin{proof}
    It is easy to check that $\Lap w_o(0) = 2c(a-ar) = 1$, and that $\Lap w_o(n) = c^n(ar^n-ar^{n-1}) + c^{n+1}(ar^n-ar^{n+1}) = 0$ for $n \neq 0$. The reader may check that $\energy(w_o) = \frac r{2(1-r)}$ so that $w_o \in \HE$. The computations for $v_x$ and $h$ are essentially the same. See Figure~\ref{fig:vx-on-Zc}.
  \end{proof}
\end{lemma}
  
  In Figure~\ref{fig:vx-on-Zc}, one can also see that $f_1 = \Pfin v_1$ induces a current flow of 1 amp from 1 to 0, with  $\frac{1+r}2$ amps flowing down the 1-edge path from 1 to 0, and the remaining current of $\frac{1-r}2$ amps flowing down the ``pseudo-path'' from 1 to $+\iy$ and then from $-\iy$ to 0.

  \begin{figure}
    \centering
    \scalebox{1.0}{\includegraphics{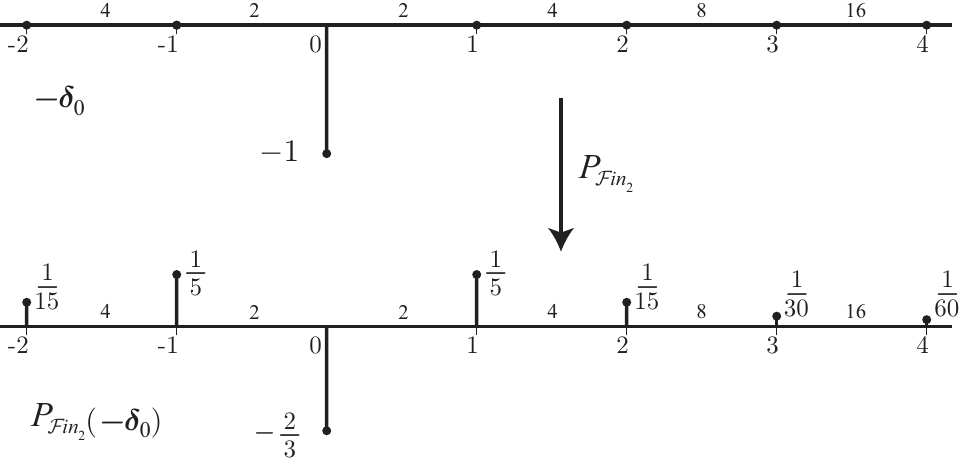}}
    \caption{\captionsize The projection of the Dirac mass $-\gd_o$ onto $\Fin_2$; see Example~\ref{exm:Fin_2-not-dense-in-Fin} and also Lemma~\ref{thm:monopole-on-geometric-integers} and Lemma~\ref{thm:repkernels-on-geometric-integers}.}
    \label{fig:Dirac-projected-to-Fin2}
  \end{figure}

\begin{exm}[Geometric half-integer model]\label{def:geometric-half-integers}
  It is also interesting to consider $(\bZ_+,c^n)$, as this network supports a monopole, but has $\Harm = 0$. 
  \linenopax
  \begin{align*}
    \xymatrix{
      0 \ar@{-}[r]^{c} 
      & 1 \ar@{-}[r]^{c^2} 
      & 2 \ar@{-}[r]^{c^3} 
      & 3 \ar@{-}[r]^{c^4} 
      & \dots
    }
  \end{align*} 
  As in Lemma~\ref{thm:monopole-on-geometric-integers}, it is straightforward to check that $w_o(n) = ar^{|n|}$, $a:= \frac{r}{(1-r)}$, defines a monopole on the geometric half-integer model $(\bZ_+,c^n)$. However, it is also easy to check by induction that $\Harm = 0$ for this model.
  
  For $k=2,3,\dots$, the network $(\bZ_+,k^n)$ can be thought of as the ``projection'' of the homogeneous tree of degree $k$ $(\sT_k, \tfrac1k \one)$ under a map which sends $x$ to $n \in \bZ$ iff there are $n$ edges between $x$ and $o$.
\end{exm}

\begin{exm}[Decomposition in \Gdd]\label{exm:monopolar-decomposition}
  In Remark~\ref{rem:comparison-to-Royden-decomp}, we discussed the Hilbert space \Gdd and its inner product $\la u, v \ra_\gdd := u(o)v(o) + \la u, v \ra_\energy$. Since $(\bZ_+,c^n)$ and $(\bZ,c^n)$ are both transient for $c>1$ (but only the latter contains harmonic functions), it is interesting to consider $P_{\Gddo}\one$ for these models (see Remark~\ref{rem:transient-iff-1-splits}). The projections $v = P_{\Gddo}\one$ and $u = \one-v = P_{\Gddo}^\perp\one$ on $(\bZ,c^n)$ are given by
  \linenopax
  \begin{align}\label{eqn:v=P_Do1-on-Z}
    v(x) = 2 - 2a + ar^{|x|} 
    \qq\text{and}\qq
    u(x) = 2a - 1 - a r^{|x|}, 
  \end{align}
  where with $a = \frac1{3-2c}$ and $r=c^{-1}$, and one can check that $v \in \MP_o^+$ and $u \in \MP_o^-$; see Definition~\ref{def:poles-and-antipoles} and Lemma~\ref{thm:positive-monopoles-in-Do}. In particular, $\Lap v = (1-v_o)\gd_0$ and $\Lap u = -u_o\gd_0$ (as usual, $o=0$). Now consider the representative of $w \in \MP_o$ given by
  \linenopax
  \begin{align}\label{eqn:monopolar-charfn}
    w(x) = (2-a) \charfn{[-\iy,0]} + \left(1 + 2a(r^{|x|} - c)\right)\charfn{[1,\iy)} \,.
  \end{align}
  A straightforward computation shows that $w = v + h$ with $h \in \GHD_o$. 
  
  The function $v = P_{\Gddo}\one$ was computed for $(\bZ,c^n)$ in \eqref{eqn:v=P_Do1-on-Z} by using the formula $\energy(u) = u_o - u_o^2$, from Lemma~\ref{thm:parabolic-u(o)-parameter}, where $u := P_{\Gddo}^\perp \one = \one - v$ and $u_o = u(o)$. For a general network $(\Graph,\cond)$, this formula implies that $(u_o, \energy(u))$ lies on a parabola with $u_o \in [0,1)$ and maximum at $(\frac12,\frac14)$. From \eqref{eqn:v=P_Do1-on-Z}, it is clear that the network $(\bZ,c^n)$ provides an example of how $u_o = 1 - \frac1{2c-1}$ can take any value in $[0,1)$. Note that $c=1$ corresponds to $\energy(u)=0$, which is the recurrent case.
\end{exm}
  
\begin{exm}[Star networks]\label{exm:star}
  Let $(\sS_m,c^n)$ be a network constructed by conjoining $m$ copies of $(\bZ_+,c^n)$ by identifying the origins of each; let $o$ be the common origin. 
\end{exm}

\begin{exm}[$\Fin_2$ not dense in \Fin]
  \label{exm:Fin_2-not-dense-in-Fin}
  On $(\bZ,2^n)$, $\Fin_2 = \cl {\spn\{\gd_x-\gd_o\}}$ is not dense in $\Fin$ (see Definition~\ref{def:Fin2-and-Fin1}). To illustrate this, we compute the projection of $-\gd_o$ to $\Fin_2$. This may be accomplished by computing 
  \linenopax
  \begin{align*}
    u_n := \left[\text{projection of } -\gd_o \text{ to } \spn\{\gd_x-\gd_o \suth |x| \leq n\}\right],
  \end{align*}
  and then taking the limit as $n \to \iy$. The result is depicted in Figure~\ref{fig:Dirac-projected-to-Fin2}.
  We leave the computation of the case of general geometric conductance $(\bZ,c^n)$ as an exercise.
\end{exm}

\section{Defect spaces}
\label{sec:defect-on-geometric-integers}

We will construct a defect vector $u \in \HE$ satisfying $\Lap u = -u$ on $(\bZ,c^n)$, $c>1$, the 1-dimensional integer lattice with geometrically growing conductances. We do this in two stages: (i) construct a defect vector on $(\bZ_+,c^n)$, and (ii) combine two copies of this defect vector to obtain an example on $(\bZ,c^n)$. 

\begin{exm}[Defect on the positive integers]
  \label{exm:defective-positive-integers}\label{exm:integer-lattice-with-defect}
  We consider $(\bZ_+,\cond)$ where 
  \linenopax
  \begin{align*}
    \cond_{n-1,n} = c^n, \q &n \geq 1,
  \end{align*}
  for some fixed $c>1$. Thus, the network under consideration is
  \linenopax
    \begin{align*}
    \xymatrix{
      0 \ar@{-}[r]^{c} 
      & 1 \ar@{-}[r]^{c^2} 
      & 2 \ar@{-}[r]^{c^3} 
      & 3 \ar@{-}[r]^{c^4} & \dots
    }
  \end{align*} 
  Now recursively define a system of polynomials in $r=1/c$ by
  \linenopax
  \begin{align*}
    \left[\begin{array}{cc}
      p_n \\ q_n
    \end{array}\right] 
    =
    \left[\begin{array}{cc}
      1 & 1 \\ r^n & 1+r^n
    \end{array}\right] 
    \cdots
    \left[\begin{array}{cc}
      1 & 1 \\ r^2 & 1+r^2
    \end{array}\right] 
    \left[\begin{array}{cc}
      1 & 1 \\ r & 1+r
    \end{array}\right] 
    \left[\begin{array}{cc}
      0 \\ 1
    \end{array}\right] 
  \end{align*}
  We will show that $u(n) := q_n$ satisfies $\Lap u = -u$ and has finite energy. It will be helpful to note that 
  \linenopax
  \begin{align}\label{eqn:p_n-is-diff}
    p_n = c^n (u(n) - u(n-1)),
  \end{align}
  and hence
  \linenopax
  \begin{align*}
    p_{n+1} = p_n + q_n,
    \q\text{and}\q
    q_{n+1} = q_n + r^{n+1} p_{n+1}.
  \end{align*}
  Now, $\Lap u = -u$ because
  \linenopax
  \begin{align*}
     \Lap u(n) 
     = p_n - p_{n+1}
     = -q_n = -u(n).
  \end{align*}
  We will need the following lemma to show that $u \in \HE$.

\begin{lemma}\label{thm:geometric-bounds-on-p,q}
  There is an $m$ such that $p_n \leq n^m$ and $q_n \leq (n+1)^m - n^m$ for $n \in \bZ_+$.   
  \begin{proof}
    We prove both bounds simultaneously by induction, so assume both bounds hold for $n$ and prove
    \linenopax
    \begin{align*}
      p_{n+1} &\leq (n+1)^m, \q \text{and} \\
      q_{n+1} &\leq (n+2)^m - (n+1)^m.
    \end{align*}
    The estimate for $p_{n+1} = p_n + q_n$ is immediate from the inductive hypotheses.
    For the $q_{n+1}$ estimate, choose an integer $m$ so that 
    \linenopax
    \begin{align*}
      m(m-1) 
      \geq \max\{t^2 r^t \suth t \geq 0\} 
      = \left(\frac{2}{e \log c}\right)^2. 
    \end{align*}
    Then $(n+1)^2 r^{n+1} \leq m(m-1)$ for all $n$, so
    \linenopax
    \begin{align*}
      2 + r^{n+1}
      \leq 2 + \frac{m(m-1)}{(n+1)^2}
      \leq \left(\frac{n}{n+1}\right)^m + \left(\frac{n+2}{n+1}\right)^m,
    \end{align*}
    by using the binomial theorem to expand $\left(\frac{n}{n+1}\right)^m = \left(1-\frac{1}{n+1}\right)^m$ and $\left(\frac{n+2}{n+1}\right)^m = \left(1+\frac{1}{n+1}\right)^m$. Multiplying by $(n+1)^m$ gives
    \linenopax
    \begin{align*}
      \left((n+1)^m - n^m\right) + r^{n+1}(n+1)^m &\leq (n+2)^m - (n+1)^m,
    \end{align*}
    which is sufficient because the left side is an upper bound for $q_{n+1} = q_n + r^{n+1} p_{n+1}$.
  \end{proof}
\end{lemma}

\begin{lemma}\label{thm:defect-vector-has-finite-energy}
  The defect vector $u(n) := q_n$ has finite energy and is bounded. 
  \begin{proof}
    Applying Lemma~\ref{thm:geometric-bounds-on-p,q} to the formula for \energy yields 
    \linenopax
    \begin{align*}
      \energy(u)
      &= \sum_{n=1}^\iy c^n(u(n)-u(n-1))^2
      = \sum_{n=1}^\iy r^n p_n^2
      \leq \sum_{n=1}^\iy r^n n^{2m}
      = \operatorname{Li}_{-2m}(r) < \iy,
    \end{align*}
    since a polylogarithm indexed by a negative integer is continuous on $\bR$, except for a single pole at 1 (but recall that $r \in (0,1)$).
  \end{proof}
\end{lemma}

  Lemma~\ref{thm:defect-vector-has-finite-energy} ensures that the defect vector is bounded; in the example in Figure~\ref{fig:defect-vector}, the defect vector has a limiting value of $\approx 4.04468281$, although the function value does not exceed 4 until $x=10$. The first few values of the function are
  \linenopax
  \begin{align*}
    u &= \left[ \tfrac{3}{2}, \tfrac{17}{8}, \tfrac{173}{64}, \tfrac{3237}{1024}, \tfrac{114325}{32768}, \tfrac{7774837}{2097152}, \tfrac{1032268341}{268435456}, \tfrac{270040381877}{68719476736}, \tfrac{140010315667637}{35184372088832}, \dots\right] \\
      &\approx \left[{1.5}, {2.125}, {2.7031}, {3.1611}, {3.4889}, {3.7073}, {3.8455}, \
{3.9296}, {3.9793}, {4.0080}, \dots \right]
  \end{align*}
  
  \begin{figure}
    \centering
    \scalebox{1.0}{\includegraphics{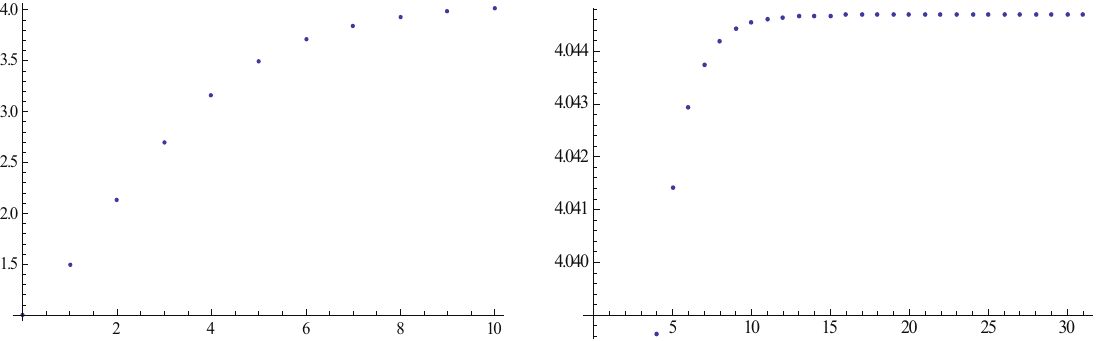}}
    \caption{\captionsize A Mathematica plot of the defect vector $u$ on $(\bZ_+,2^n)$; see Example~\ref{exm:defective-positive-integers} and Lemma~\ref{thm:defect-vector-has-finite-energy}. The left plot shows $u(x)$ for $x=0,1,\dots,10$, and the plot on the right shows data points for $u(x)$, $x=10,11,12,\dots$.}
    \label{fig:defect-vector}
  \end{figure}

  While we are unable to provide a nice closed-form formula for the defect vector, we can provide generating functions for it, using the $p_n = p_n(r)$ and $q_n = q_n(r)$ obtained just above. Define
  \linenopax
  \begin{align*}
    P(x) = \sum_{n=0}^\iy p_n(r) x^n
    \q\text{and}\q
    Q(x) = \sum_{n=0}^\iy q_n(r) x^n.
  \end{align*}
  Multiplying both sides of $p_{n+1} = p_n + q_n$ by $x^{n+1}$ and summing from $n=0$ to \iy, 
  \linenopax
  \begin{align}\label{eqn:genfun-for-P-1}
    P(x) = xP(x) + xQ(x),
  \end{align}
  where we have used the fact that $p_0=0$.
  Meanwhile, multiplying both sides of $q_{n+1} = q_n + r^{n+1}p_{n+1}$ by $x^{n+1}$ and summing from $n=0$ to \iy, 
  \linenopax
  \begin{align}\label{eqn:genfun-for-Q-1}
    Q(x) -1 = xQ(x) + P(rx).
  \end{align}
  Write \eqref{eqn:genfun-for-P-1} in the form $(1-x)P(x) = xQ(x)$ and substituting in $(1-x)Q(x) = 1 + P(rx)$ from \eqref{eqn:genfun-for-Q-1}, to get $1 + P(rx) = (1-x)Q(x) = \tfrac{(1-x)^2}{x} P(x)$ or
  \linenopax
  \begin{align*}
    P(x) 
    = \tfrac{x}{(1-x)^2} + \tfrac{x}{(1-x)^2}P(rx) 
    &= \tfrac{x}{(1-x)^2} + \tfrac{x(rx)}{(1-x)^2(1-rx)^2} + \tfrac{x(rx)}{(1-x)^2(1-rx)^2}P(r^2x) \\
    \cdots 
    &= \sum_{n=0}^\iy \prod_{k=0}^n \frac{r^k x}{(1-r^kx)^2} 
    = \sum_{n=0}^\iy\frac{r^{n(n+1)/2} x^n}{\prod_{k=0}^n (1-r^kx)^2}.
  \end{align*}
  Note that $r \in (0,1)$, so $P(r^kx) \limas{k} P(0) = 0$, again since $p_0=0$. Now \eqref{eqn:genfun-for-P-1} gives $Q(x) = \frac{1-x}x P(x)$, whence
  \linenopax
  \begin{align*}
    Q(x) = \sum_{n=0}^\iy\frac{r^{n(n+1)/2} x^{n-1}}{\prod_{k=1}^n (1-r^kx)^2}. 
  \end{align*}
\end{exm}

\begin{exm}[Defect on the integers]
  \label{exm:defective-integers}
  We consider $(\bZ,\cond)$ as in Definition~\ref{def:geometric-integers}:
  \linenopax
  \begin{align*}
    \xymatrix{
      \dots \ar@{-}[r]^{c^3}
      & {-2} \ar@{-}[r]^{c^2} 
      & {-1} \ar@{-}[r]^{c} 
      & {0} \ar@{-}[r]^{c} 
      & {1} \ar@{-}[r]^{c^2} 
      & {2} \ar@{-}[r]^{c^3} 
      & {3} \ar@{-}[r]^{c^4} & \dots
    }
  \end{align*} 
  
  Proceeding as in Example~\ref{exm:defective-positive-integers}, one uses $\Lap u(0) = -u(0)$ to compute
  \linenopax
  \begin{align*}
    2c(u(0) - u(1)) = -u(0)
    \q\implies\q
    u(1) = \left(1 + \tfrac 1{2c}\right) u(0),
  \end{align*}
  and obtain the initial values $p_1 = \frac12$ and $q_1 = 1+\frac r2$. 
  Therefore, for \bZ we instead use the polynomials defined by  
  \linenopax
  \begin{align*}
    \left[\begin{array}{cc}
      p_n \\ q_n
    \end{array}\right] 
    =
    \left[\begin{array}{cc}
      1 & 1 \\ r^n & 1+r^n
    \end{array}\right] 
    \cdots
    \left[\begin{array}{cc}
      1 & 1 \\ r^2 & 1+r^2
    \end{array}\right] 
    \left[\begin{array}{cc}
      1 & 1 \\ r & 2+r
    \end{array}\right] 
    \left[\begin{array}{cc}
      0 \\ \tfrac12
    \end{array}\right] 
  \end{align*}
  The other computations are essentially identical to those for $(\bZ_+,c^n)$.
\end{exm}

\version{}{\marginpar{Add $w_x \in \ran \LapV^\ad$}}

\section{Remarks and references}
\label{sec:Remarks-and-References-lattice-examples}

The infinite lattices offer a second attractive family of examples; and they are especially relevant for lattice-spin models in physics, as discussed in Chapter~\ref{sec:Magnetism-and-long-range-order}. The book \cite{Soardi94} by Soardi is a nice introduction to the subject, and a classical introductory reference is \cite{Spitzer}. Of the results in the literature of relevance to the present chapter, the references \cite{GovN51, GoSt06, ScZi09, Martinelli99, ChLi07, Lig99, Lig95, Lig93} are especially relevant. 

The geometric integers of Example~\ref{exm:defective-positive-integers} came about from our desire to apply von NeumannÕs theory of unbounded operators and their deficiency indices \cite{vN32a, vN32b, vN32c, DuSc88} to the metric geometry of infinite weighted graphs $(\Graph, \cond)$. Starting with $(\Graph, \cond)$ there are two natural Hilbert spaces $\ell^2(\verts)$ (where \verts is the vertex-set) and the energy Hilbert space \HE. An intriguing aspect of \S\ref{sec:defect-on-geometric-integers} is that the boundary features of $(\Graph, \cond)$ deriving from deficiency indices cannot be accounted for with the use of the more naive of the two Hilbert spaces $\ell^2(\verts)$; \HE is forced upon us. 

The geometric integers discussed in Example~\ref{exm:defective-positive-integers} is called a \emph{weighted linear graph} in \cite{AlFi09} and is studied in conjunction with birth and death processes; see the references cited therein.

%% file: magnetism.tex

\chapter[Magnetism and long-range order]{Application to magnetism and long-range order}
\label{sec:Magnetism-and-long-range-order}

\headerquote{Physics is becoming too difficult for the physicists.}{---~D.~Hilbert}

\headerquote{For a physicist mathematics is not just a tool by means of which phenomena can be calculated, it is the main source of concepts and principles by means of which new theories can be created.}{---~F.~Dyson}

The integer lattice examples studied in \S\ref{sec:lattice-networks} may be applied to the theory of ferromagnetism. In \S\ref{sec:Kolmogorov-construction-of-L2(gW,prob)}, we construct a Hilbert space $L^2(\gW,\prob)$ with a probability measure, following techniques of Kolmogorov. Since $L^2(\gW,\prob) \cong \HE$, this provides a concrete realization of \HE as a probability space and a commutative analogue/precursor of the Heisenberg spin model developed in \S\ref{sec:magnetism}. In \S\ref{sec:GNS-construction} we carry out the GNS construction \cite{Arveson:invitation-to-Cstar} to obtain a Hilbert space $\sH_\gf$. Again, this will be useful for \S\ref{sec:magnetism}, where we recall Powers' approach (using a \gb-KMS state \gf) and show how our results may be used to obtain certain refinements of Powers' results.

In \cite{Pow76a}, Powers made the first connection between the two seemingly unrelated ideas: resistance distances in electrical networks, and a problem from statistical mechanics. Even in the physics literature, one often distinguishes between quantum statistical model as emphasizing (a) physics, or (b) rigorous mathematics. The literature for (a) is much larger than it is for (b); in fact, the most basic questions (phase transition and long-range order) are notoriously difficult for (b). Powers was concerned with long-range order in ferromagnetic models from quantum statistical mechanics, especially Heisenberg models. The notion of long-range order in these models depends on a chosen Hamiltonian $H$, and a $C^\ad$-algebra \xA of local observables for such models. These objects and ideas are discussed in more detail in \S\ref{sec:long-range-order-in-ERNs} and \S\ref{sec:magnetism}. The reader may also find some information on \gb-KMS states in \S\ref{sec:KMS-states-appendix}.

While we shall refer to the literature, e.g. \cite{BrRo79,Rue69} for formal definitions of the key terms from the $C^\ad$-algebra formalism of quantum spin models, physics, KMS states and the like, we present a minimal amount of background and terminology from the mathematical physics literature so our presentation is agreeable to a mathematical audience. A brief discussion of KMS states is given in \S\ref{sec:KMS-states-appendix} and the reader may wish to peruse the general GNS construction is outlined in \S\ref{sec:GNS-construction} before reading \S\ref{sec:GNS-construction}.

Here we turn to a non-commutative version of the infinite Cartesian products that went into the probabilistic constructs used in sections 7 and 11 above. This is dictated directly from quantum physics: Think of an algebra of observables placed on each vertex point in an infinite graph, each algebra non-abelian because of quantum mechanics. The infinite graphs here may represent sites from a solid state model, or a spin-model for magnetization. To achieve our purpose, we will use infinite tensor products of $C^\ad$-algebras, one for each point $x \in \verts$. This is dictated by our application to quantum statistical mechanics. In quantum physics, the entity that corresponds to probability measures in classical problems are however ``states'' on the algebra of all the quantum mechanical observables, a $C^\ad$-algebra, but the $C^\ad$-algebra for the entire system will be an infinite tensor product $C^\ad$-algebra. To gain intuition, the reader may wish to think of states as non-commutative measures, and hence non-commutative probability theory (see e.g., \cite{BrRo97}.)

\section{Kolmogorov construction of $L^2(\gW,\prob)$}
\label{sec:Kolmogorov-construction-of-L2(gW,prob)}
As a prelude to the quantum-mechanical model, we first give a probabilistic model, that is a model for classical particles, which serves to illustrate the main themes. In particular, long-range order appears in this setting as an estimate on correlations (in the sense of probability).

We consider a Brownian motion on $(\Graph,\cond)$ as a system of Gaussian random variables, again indexed by \verts. For these (commutating) random variables, we will show the correlations are given by the resistance distance $R(x,y)$. This result is extended to the noncommutative setting via the GNS construction in \S\ref{sec:GNS-construction}.

In \HE, we don't really have an algebra of functions, so first we make one. Since $\Ex(v_x,v_y) := \la v_x, v_y\ra_\energy$ is a positive definite form \version{}{\marginpar{Mention why this is good enough to produce consistency \`{a} la Theorem~\ref{thm:Kolmogorov-extension-thm}.}} $\verts \times \verts \to \bC$, we can follow Kolmogorov's construction of a measure on the space of functions on \verts. Denoting the Riemann sphere by $\bS^2 = \bC \cup \{\iy\}$, this produces a probability measure \prob on
\linenopax
\begin{equation}\label{eqn:def:gW=prod(C)}
  \gW = \prod_{x \in \verts} \bS^2,
\end{equation}
the space of all functions on \verts. Also, we define
\linenopax
\begin{equation}\label{eqn:def:tilde-vx}
  \tilde v_x: \gW \to \bC \q\text{by}\q
  \tilde v_x (f) := f(x)-f(o).
\end{equation}
Since \HE is a Hilbert space, \HE is its own dual, and we can think of $v_x$ as an element of \HE or the function on \HE defined by $\la v_x, \cdot\ra_\energy$. In the latter sense, $\tilde v_x$ is an extension of $v_x$ to \gW; observe that for $u \in \HE$,
\linenopax
\begin{equation}\label{eqn:tildevx-repkernel-extn}
  \tilde v_x(u) = \la v_x, u \ra_\energy = u(x)-u(o).
\end{equation}
Thus we have a Hilbert space $L^2(\gW,\prob)$ which contains as a dense subalgebra the algebra generated by $\{\tilde v_x\}$.

Another consequence of the Kolmogorov construction is that
\linenopax
\begin{equation}\label{eqn:expectation-as-inner-prod}
  \Ex(\cj{\tilde v_x}\, \tilde v_y) = \int \cj{\tilde v_x}\, \tilde v_y \,d\prob = \la v_x, v_y \ra_\energy.
\end{equation}

\begin{lemma}\label{thm:L2(prob)=HE}
  $L^2(\gW,\prob)$ is unitarily equivalent to \HE.
  \begin{proof}
    The mapping $\tilde v_x \mapsto v_x$ extends by linearity to an isometric isomorphism:
    \linenopax
    \begin{align*}
      \Ex\left(\left|\sum_{x \in \verts} c_x \tilde v_x\right|^2\right)
        = \left\|\sum_{x \in \verts} c_x v_x \right\|_\energy^2
      \iff \sum_{x,y \in \verts} \negsp[8] \bar c_x \Ex(\cj{\tilde v_x}\, v_y) c_y
      = \sum_{x,y \in \verts} \negsp[8] \bar c_x \la v_x, v_y\ra_\energy c_y,
    \end{align*}
    which is true by \eqref{eqn:expectation-as-inner-prod}.
  \end{proof}
\end{lemma}

Observe that one recognizes $\gy \in L^2(\gW,\prob)$ as corresponding to a \emph{finite} linear combination $\sum c_x v_x$ if and only if there is a finite subset $F \ci \verts$ and a function $u:F \to \bC$ with $\spt u = F$ such that $\gy(\gw) = u(\gw)$.
By Riesz's Lemma, integration with respect to \prob is given by a positive linear functional $\gf_\prob$, i.e., the expectation is
\linenopax
\[\Ex(f) = \int f \,d\prob = \gf_\prob(f).\]
Since \prob is a probability measure, we even have $\gf_\prob(\one) = 1$. Consequently, $\gf_\prob$ corresponds to a \emph{state} in the noncommutative version (cf.~Definition~\ref{def:state}); see Remark~\ref{rem:commutative-noncommutative-dictionary} and the table of Figure~\ref{fig:classical-vs-quantum-dictionary}.

\begin{exm}[Application of Lemma~\ref{thm:L2(prob)=HE} to the integer lattice network $(\bZd,\one)$]\label{exm:L2(prob)=HE(Zd,1)}
  Observe that Bochner's Theorem (Theorem~\ref{thm:Bochner's-theorem}) gives
  \linenopax
  \begin{align*}
    \Ex_\gx(\cj{e^{\ii x \cdot \gx}} e^{\ii y \cdot \gx})
    &= \int_{\bRd} e^{\ii(y-x) \gx} \;d\prob(\gx) = \la v_x, v_y\ra_\energy.
  \end{align*}
  Thus, we are obliged to set $\tilde v_x(\gx) := e^{\ii x \cdot \gx}$ for $\gx \in \bRd$, whence the mapping $e^{\ii x \cdot \gx} \mapsto v_x(\gx)$ extends by Lemma~\ref{thm:L2(prob)=HE} and $L^2(\bRd,\prob) \cong \HE$. A particularly striking feature of this example is that one sees that translation-invariance of the underlying network causes the \emph{a priori} infinite-dimensional lion \gW to devolve into the finite-dimensional lamb \bRd. The duals of abelian groups are much tamer!
\end{exm}

\section{The GNS construction}
\label{sec:GNS-construction}

The GNS construction takes a $C^\ad$-algebra \xA and a state $\gf:\xA \to \bC$ (see Definition~\ref{def:state} just below), and builds Hilbert space $\sH_\gf$ and a representation $\gp:\xA \to B(\sH_\gf)$. The main point is that even though a $C^\ad$-algebra can be defined axiomatically and without reference to any Hilbert space, one can always think of a $C^\ad$-algebra as an algebra of operators on some Hilbert space. The GNS construction stands for Gel'fand, Naimark, and Segal, and the literature on this construction is extensive; we include a sketch of the proof, but point the reader to \cite[\S1.6]{Arveson:invitation-to-Cstar} (for newcomers) and \cite[\S2.3.2]{BrRo79} (for details).

Following the overview of the general GNS construction, we explain how the GNS construction provides a noncommutative analogue of the Kolmogorov model discussed in the previous section. The Heisenberg model is built within the representation of a certain $C^\ad$-algebra, and we will need this framework to describe Powers' results concerning magnetism. We also provide some of the background material relevant to the applications to the theory of magnetism and long-range order discussed in \S\ref{sec:Magnetism-and-long-range-order}; see also the excellent references \cite{Arv76a,Arveson:invitation-to-Cstar,Arv76b,BrRo79}.

\begin{defn}\label{def:state}
  A \emph{state} on a $C^\ad$-algebra \xA is a linear functional $\gf:\xA \to \bC$ which satisfies $\gf(A^\ad A) \geq 0$ and $\gf(\one)=1$.
\end{defn}

\begin{theorem}[GNS construction]
  \label{thm:GNS-construction}
  Given a $C^\ad$-algebra \xA, a unit vector $\one \in \xA$ and a state \gf, there exists
  \begin{enumerate}
    \item a Hilbert space: $\sH_\gf$, $\la \cdot,\cdot \ra_\gf$,
    \item a representation $\gp:\xA \to B(\sH_\gf)$ given by $A \mapsto \gp(A)(\cdot)$, and
    \item a cyclic vector (the \emph{ground state}\footnote{\gz is called the ground state because when \gf is a KMS state built from the Hamiltonian $H$, one has $H\gz=0$, i.e., the energy of \gz is 0}): $\gz = \gz_\gf \in \sH_\gf$, $\|\gz\|_\gf=1$,
  \end{enumerate}
  for which $\gf(A) = \la \gz, \gp(A)\gz \ra_\gf$, $\forall A \in \xA$.
  \begin{proof}[Sketch of proof]
    For (1), define $\la A,B \ra_\gf := \gf(A^\ad B)$.\footnote{This is why physicists make the inner product linear in the second variable.}
    Define the kernel of \gf in the nonstandard fashion
    \[\ker(\gf) := \{A \in \xA \suth \gf(A^\ad A) = 0\}.\]
    Intuitively, think $\gf(A^\ad A) \leftrightarrow \int |f|^2$. Then one has a Hilbert space by taking the completion
    \[\sH_\gf = \left(\xA / \ker(\gf)\right)^\sim.\]

    For (2), show that the multiplication operator $\gp(A): B \mapsto AB$ is a bounded linear operator on \xA. This follows from the computation
    \linenopax
    \begin{align*}
      \|\gp(A) B \|_\gf \leq \|A\|_{C^\ad} \|B\|_\gf
      &\iff \|\gp(A) B \|_\gf^2 \leq \|A\|_{C^\ad}^2 \|B\|_\gf^2 \\
      &\iff \gf((AB)^\ad AB) \leq \|A\|_{C^\ad}^2 \gf(B^\ad B),
    \end{align*}
    which is true because $\gf_B(A) := \frac{\gf(B^\ad A B)}{\gf(B^\ad B)}$ is a state, $\|A^\ad A\| = \|A\|^2$, and $|\gf_B(A)| \leq \|A\|$ for every $A \in \xA$.

    For (3), start with $\one \in \xA$. Then $\gz = \gz_\gf$ is the image of \one under the embedding
    \linenopax
    \[\xymatrix{
        \xA \ar[rr]^-{\; projection \;}
        && \frac{\xA}{\ker(\gf)} \ar[rr]^-{\; completion \;}
        && \sH_\gf = \left(\frac{\xA}{\ker(\gf)}\right)^\sim
      }\]
    During this composition, \one is transformed as follows: $\one \mapsto \one + \ker(\gf) \mapsto \gz_\gf$. Finally, to verify the condition relating (1),(2),(3), use $[\cdot]_\gf$ to denote an equivalence class in the quotient space and then
    $\la \gz_\gf, \gp(A) \gz_\gf \ra_\gf
      = \la[\one]_\gf, A \cdot [\one]\ra_\gf
      = \gf(\one A \one)
      = \gf(A)$.
  \end{proof}
\end{theorem}

\begin{remark}\label{rem:commutative-simplification}
  When \xA is a \emph{commutative algebra of functions}, it turns out that $\gp(f)(\cdot)$ is multiplication by $f$, in which case the notation is a bit heavy handed:
  \linenopax
  \[\gp(f) \one_\prob = f \cdot \one = f, \qq \one_\prob \in L^2(\gW,\prob).\]
  For the noncommutative case, things are different and the full notation is really necessary. (Note that $\one_\prob$ really does depend on \prob, in the same way that the unit in $L^2(X,\gd_o)$ is different from the unit of $L^2(X,dx)$).
\end{remark}

\begin{remark}\label{rem:contractivity-of-gp(A)}
  In general, the resulting representation $\gp: \xA \to B(\sH_\gf)$ is a contractive injective homomorphism, so that $\|\gp(A)\| \leq \|A\|$. However, when \xA is simple (as is the case in our setting), then \gp is actually an isometry.
\end{remark}

\begin{figure}
  \begin{center}
    \begin{tabular}{|l|l|}
      \hline
      probabilistic/classical & quantum \vstr[2.3] \\ \hline
      space $\gW = \prod_{x \in \verts} \bS$ & $C^\ad$-algebra $\xA = \bigotimes_{x \in \verts} \xA_x$ \vstr[2.3] \\
      Gaussian measure \prob & state \gf (or KMS state \gw) \\
      probability space $L^2(\gW,\prob)$ & Hilbert space $\sH_\gf = GNS(\xA,\gf)$ \\
      function $\tilde v \in L^2(\gW,\prob)$ & observable $\gs:\verts \to \xA$, $\gs_x \in \xA_x$ \\
      constant function \one & ground state \gz \\
      embedding $W:\HE \to L^2(\gW,\prob)$ & representation $\gp_\gf:\xA \to \sB(\sH_\gf)$ \\
      $v \mapsto \tilde v(\one)$ & $A \mapsto \gp_\gf(A) \gz$ \\
      expectation $\Ex(v) = \int v \,d\prob$
        & measurement $\gf(\gs) = \la \gz, \gp_\gf(\gs) \gz \ra_\gf$ \\
      covariance $\Ex(\bar{\tilde v}_x \tilde v_y) = \int \bar{\tilde v}_x \tilde v_y d\prob$ & correlation $\gf(\gs_x^\ad \cdot \gs_y) = \la \gs_x, \gs_y \ra_\gf$ \\ \hline
    \end{tabular}
  \end{center}
  \caption{\captionsize A ``dictionary'' between the classical and quantum aspects of this problem. In this table, $\bS$ is the Riemann sphere (the one-point compactification of \bC) and $H$
   is the Hamiltonian discussed by Powers. The notation $\gs_x \cdot \gs_y$ is explained in \eqref{eqn:gs.gs-notation}. This table is elaborated upon in Remark~\ref{rem:commutative-noncommutative-dictionary}.}
  \label{fig:classical-vs-quantum-dictionary}
\end{figure}

\begin{remark}(Kolmogorov construction vs. GNS construction)
  \label{rem:commutative-noncommutative-dictionary}
  In Figure~\ref{fig:classical-vs-quantum-dictionary} we present a table which gives an idea of how analogous ideas match up in the commutative and noncommutative models on the same \ERN $(\Graph, \cond)$. The titles of the two columns in the table refer to nature of the corresponding random variables. In both columns, the variables are indexed by vertices.

  In the left column, \gW is simply a (commutative) family of measurable functions; the collection of all measurable functions on \verts. One can think of this as the tensor product of 1-dimensional algebras \bC. On the right, the variables are quantum observables, so noncommuting self-adjoint elements in a $C^\ad$-algebra.  For infinite \ERNs, the $C^\ad$-algebra $\xA = \xA(\Graph)$ in question is built as an infinite tensor product of finite-dimensional $C^\ad$-algebras $\xA_x$, $x \in \verts$. More specifically, \xA is the inductive limit $C^\ad$-algebra of algebras $\xA(F)$ where $F$ ranges over all finite subsets in \verts, and where each $\xA(F) = \bigotimes_{x \in F} \xA_x$ is a finite tensor product.

  A general element $A \in \xA$ does not have a direct analogue in the left-hand column, as the structure is much richer in the noncommutative case. The observable \gs is a particularly simple type of element of \xA; it is one which can be represented as a single element of $\prod_{x \in \verts} \xA_x$. A general $A \in \xA$ can only be represented as a sum of such things; cf.~\cite{BrRo79}.


  One can think of \gW as $L^\iy(\gW,\prob)$, as the latter is generated by the coordinate functions. Let $X_x$ be the random variable
  \[X_x: \gW \to \bC \qq\text{by}\qq X_x(\gw) = \gw(x).\]
  (Recall that \gw is any measurable function on \verts.) Then $X_x$ corresponds to $\tilde v_x = \la v_x, \cdot \ra_\energy$, as can be seen by considering the reproducing kernel property \eqref{eqn:tildevx-repkernel-extn}, when representatives of \HE are chosen so that $v(o)=0$. In this sense, $L^\iy(\gW)$ is the commutative version of \xA. Recall from Stone's Theorem that every abelian von Neumann algebras is $L^\iy(X)$ for some measurable space $X$. In a similar vein, The Gel'fan-Naimark Theorem states that every abelian $C^\ad$-algebra is $C(X)$ for some compact Hausdorff space $X$. Thus, $C(\gW)$ is a dense subalgebra of $L^2(\gW,\prob)$ in the same way that \xA is a dense subalgebra of $\sH_\gf$.

  One key point is that a state is the noncommutative version of a probability measure, in the $C^\ad$-algebraic framework of quantum statistical mechanics. In the abelian case, the Gaussian measure \prob is unique, while in the quantum statistical case, the states \gf typically are not unique. In fact, when requiring \gf to be a \emph{KMS state} as in the next section, then a ``phase transition'' is precisely the situation of multiple \gb-KMS states corresponding to the same value of \gb. KMS states are equilibrium states, so a phase transition is when more than one equilibrium state (e.g., liquid and vapour) are simultaneously present; see \S\ref{sec:KMS-states-appendix}. The physicist will recognize the table entry $\gf(A) = \la \gz, \gp_\gf(A) \gz \ra_\gf$ as the transition probability from the ground state to the excited state $A$.

  In the Heisenberg model of ferromagnetism, the graph is $(\bZd,\one)$, and the support of the associated Gaussian measure \prob is \bRd. Thus, the associated probability space is $(\bRd,\prob)$ and the Hilbert space is $L^2(\bRd,\prob)$. As a result, the random variables are $L^2$-functions on \bRd, obtained by extension from \bZd. See Example~\ref{exm:L2(prob)=HE(Zd,1)}.
\end{remark}

The use of \gb-KMS states is actually a crucial hypothesis in Powers' Theorem (Theorem~\ref{thm:Powers'-correlation-estimate}), although this detail is obscured in the present exposition. The technical definition of a KMS state is not critical for the main exposition of Powers' problem, but his results (and ours) would be unobtainable without this assumption. (The reasons for this are somewhat involved, but hinge upon the stability of KMS states as equilibria.) Consequently, we include a discussion in \S\ref{sec:KMS-states-appendix} outlining some key features of these states.

Powers did not consider the details of the spectral representation in GNS representation for the KMS states. More precisely, Powers did not consider the explicit function representation (with multiplicity) of the resistance metric and the graph Laplacian \Lap as it acts on the energy Hilbert space \HE. The prior literature regarding \Lap has focused on $\ell^2$, as opposed to the drastically different story for \HE.

\section{Magnetism and long-range order in \ERNs}
\label{sec:magnetism}
\label{sec:long-range-order-in-ERNs}

Following \cite{Pow75}, we apply Theorem~\ref{thm:GNS-construction} to a \gb-KMS state \gw and the $C^\ad$-algebra $\xA = \bigotimes_{x \in \verts} \xA_x$ described in Remark~\ref{rem:commutative-noncommutative-dictionary}. The inverse temperature \gb will be fixed throughout the discussion.

It is known that the translation-invariant ferromagnetic models do not have long-range order in \bZd when $d=1,2$. Powers suggested that it happens for $d=3$. Below we supply detailed estimates which bear out Powers' expectations. We can now make more precise the allusion which begins this section: Powers was the first to make a connection between
\begin{enumerate}[(i)]
  \item the resistance metric $R(x,y)$, and
  \item estimates of \gw-correlations between observables localized at distant vertices $x,y \in \verts$.
\end{enumerate}
Precise estimates for (ii) are called ``\textbf{long-range order}''; the Gestalt effect of this phenomenon is \emph{magnetism}.

The Hamiltonian $H$ appearing in Definition~\ref{def:KMS-state} as part of the definition of a \gb-KMS state is a formal infinite sum over the edges \edges of the network, where the terms in the sum are weighted with the conductance function \gm. (An explicit formula appears just below in \eqref{eqn:Hamiltonian}.) The Hamiltonian $H$ then induces a one-parameter unitary group of automorphisms $\{\ga_t\}_{t \in \bR}$ (as in Definition~\ref{def:Hamiltonian-unitary-group}) describing the dynamics in the infinite system; and a KMS state \gw refers to this automorphism group. As mentioned above, the KMS states \gw are indexed by the inverse temperature \gb. The intrepid reader is referred to the books \cite{BrRo79,BrRo97,Rue69,Rue04} for details.

In this section, we discuss an application to the spin model of the isotropic Heisenberg ferromagnet. Let $\Graph = \bZ^3$ and $\ohm \equiv 1$. We consider each vertex $x$ (or ``lattice site'') to be a particle whose spin is given by an observable $\gs_x$ which lies in the finite-dimensional $C^\ad$-algebra $\xA_x$. The $C^\ad$-algebra $\xA = \bigotimes_x \xA_x$ describes the entire system. For the case when the particles are of spin $\frac12$, an element $\gs_x \in \xA$ is expressed in terms of the three Pauli matrices:
\linenopax
\begin{equation}\label{eqn:Pauli-matrices}
  \gs_{x_1} = \left[\begin{array}{rr} 0 & 1 \\ 1 & 0\end{array}\right], \q
  \gs_{x_2} = \left[\begin{array}{rr} 0 & -\ii \\ \ii & 0\end{array}\right], \q
  \gs_{x_3} = \left[\begin{array}{rr} 1 & 0 \\ 0 & -1\end{array}\right],
\end{equation}
for any $x \in \verts$. Interaction in this isotropic Heisenberg model is given in terms of the Hamiltonian
\linenopax
\begin{equation}\label{eqn:Hamiltonian}
  H = \tfrac12 \sum_{x,y \in \verts} \cond_{xy}(\one - \gs_x\cdot\gs_y),
\end{equation}
where
\begin{equation}\label{eqn:gs.gs-notation}
  \gs_x\cdot\gs_y = \gs_{x_1} \otimes \gs_{y_1} + \gs_{x_2} \otimes \gs_{y_2} + \gs_{x_3} \otimes \gs_{y_3}.
\end{equation}
More precisely, $\gw(\id - \gs_x\cdot\gs_y)$ gives the amount of energy that would be required to interchange the spins of the particle at $x$ and the particle at $y$ when measured in state \gw, and
\linenopax
\begin{align*}
  \gw(H) = \tfrac12 \sum_{x,y \in \verts} \cond_{xy} \gw(\one - \gs_x \cdot \gs_y)
\end{align*}
is the weighted sum of all such interactions.
The Hamiltonian $H$ may be translated by time $t$ into the future by $A \mapsto \ga_t(A):= e^{\ii t H} A e^{-\ii t H} \in Aut(\xA)$.

Motivated by \eqref{eqn:Powers'-general-estimate}, Powers conjectures the following estimate for \gb-KMS states in \cite{Pow76a}: there exists a constant $K$ (independent of \Graph) for which 
\linenopax
\begin{align*}
  \gw(\one - \gs_x \cdot \gs_y) \leq K \gb^{-1} R(x,y).
\end{align*}
The following result appears in \cite{Pow76b}.

\begin{theorem}[Powers]
  \label{thm:Powers'-correlation-estimate}
  Let \gw be a \gb-KMS state and let $H$ be the Hamiltonian of \eqref{eqn:Hamiltonian}. Then
  \linenopax
  \begin{equation}\label{eqn:Power-KMS-bound}
    \gw(\one - \gs_x \cdot \gs_y) \leq \gw(H) R(x,y),
  \end{equation}
\end{theorem}
After obtaining a bound for $\gw(H)$, the author notes that in $\bZ^3$, the ``resistance between $o$ and infinity is finite'' and uses this to show that $\gw(\one - \gs_x \cdot \gs_y) = 1 - \gw(\gs_x \cdot \gs_y)$ is bounded. The interpretation is that correlation between the spin states of $x$ and $y$ remains positive, even when $x$ is arbitrarily far from $y$, and this is ``long-range order'' manifesting as magnetism. We offer the following improvement.

\begin{lemma}\label{thm:improvement-on-Powers}
  If \gw is a \gb-KMS state, then
  \linenopax
  \begin{equation}\label{eqn:improvement-on-Powers}
    1 - \liminf_{y \to \iy} \gw(\gs_x \gs_y) \leq \left(\frac1{(2\gp)^3} \int_{\bTd} \frac{dt}{2 \sum_{k=1}^3 \sin^2(\frac{t_k}2)}\right) \gw(H).
  \end{equation}
  \begin{proof}
    The identity $\lim_{y \to \iy} R(x,y) = (2\gp)^{-3} \int_{\bTd} (2 \sum_{k=1}^3 \sin^2(\frac{t_k}2))^{-1}\,dt$ is shown in Theorem~\ref{thm:finite-resistance-to-infinity-in-Zd}; a computer gives the numerical approximation
    \linenopax
    \begin{align*}
      \lim_{y \to \iy} R(x,y) \approx 0.505462
    \end{align*}
    for this integral. While the limit may not exist on the left-hand side of \eqref{eqn:Power-KMS-bound}, we can certainly take the $\limsup$, whence the result follows.
  \end{proof}
\end{lemma}

\begin{remark}[Long-range order]
  \label{rem:Long-range-order}
  In the model of ferromagnetism described above, consider the collection of spin observables $\{\gs_x\}$ located at vertices $x \in \verts$ as a system of non-commutative random variables. One interpretation of the previous results is that in KMS-states, the correlations between pairs of vertices $x,y \in \verts$ are asymptotically equal to the resistance distance $R(x,y)$. As mentioned just above, the idea is that correlation between the spin states of $x$ and $y$ remains positive, even when $x$ is arbitrarily far from $y$.

  One interpretation of this result is that magnetism can only exist in dimensions 3 and above, or else $R(x,y)$ is unbounded and the estimate \eqref{eqn:improvement-on-Powers} fails. A different interpretation of the existence of magnetism is the existence of multiple \gb-KMS states for a fixed temperature $T=1/k\gb$. This more classical view is quite different.
\end{remark}

\begin{remark}\label{rem:HE/KMS-parallel}
  We leave it to the reader to ponder the enticing parallel:
  \linenopax
  \begin{align*}
    |u(x)-u(y)|^2 &\leq \tfrac12 R(x,y) \sum\nolimits_{x,y} \cond_{xy} (u(x)-u(y))^2, \q \forall u \in \HE \q (\text{Cor.~\ref{thm:vx-is-Lipschitz}}) \\
    |\gf(\one - \gs_x \cdot \gs_y)| &\leq \tfrac12 R(x,y) \sum\nolimits_{x,y} \cond_{xy} (\gf(\one - \gs_x \cdot \gs_y)), \q \forall \gf \in\{\text{KMS states}\}
  \end{align*}
\end{remark}

\section{KMS states}
\label{sec:KMS-states-appendix}

While the rigorous definitions provided in this mini-appendix are not absolutely essential for understanding the Heisenberg model of ferromagnetism, they may help the reader understand what a \gb-KMS state is, and hence have a better feel for the discussion in the previous section. We suggest the references \cite{BrRo79,BrRo97,Rue69,Rue04} for more details.

\begin{defn}\label{def:Hamiltonian-unitary-group}
  Define $\ga: \bR \to Aut(\xA)$ by $\ga_t(A) = e^{-\ii t H} A e^{\ii t H}$, for all $t \in \bR$ and $A \in \xA$, where $H$ is a Hamiltonian (as in \eqref{eqn:Hamiltonian} below, for example). This unitary group accounts for time evolution of the system, i.e.,
  \linenopax
  \begin{align*}
    \la \gy(t), A \gy(t) \ra = \la \gy(0), \ga_t(A) \gy(0) \ra
  \end{align*}
  shows that measuring the time-evolved observable $\ga_t(A)$ in the (ground) state $\gy_0 = \gy(0)$ is the same as measuring the observable $A$ in the time-evolved state $\gy(t)$.
\end{defn}

\begin{defn}\label{def:KMS-state}
  Let \gf be a state as in Definition~\ref{def:state}. We say \gf is a \emph{KMS state} iff for all $A,B \in \xA$, there is a function $f$ with:
  \begin{enumerate}
    \item $f$ is bounded and analytic on $\{z \in \bC \suth 0 < \Im z < \gb\}$ and continuous up to the boundary of this region;
    \item $f(t) = \gf(A \ga_t(B))$, for all $t \in \bR$; and
    \item $f(t+\ii\gb) = \gf(\ga_t(B) A)$, for all $t \in \bR$.
  \end{enumerate}
  Note that $f$ depends on $A$ and $B$. This definition is roughly saying that there is an analytic continuation from the graph of $\gf(A \ga_t(B))$ to the graph of $\gf(\ga_t(B)A)$, where both are considered as functions of $t \in \bR$.
\end{defn}

\begin{defn}\label{def:beta-KMS-state}
  If \xA is finite-dimensional, then
  \linenopax
  \begin{equation}\label{eqn:beta-KMS}
    \gf(A) = \gf_\gb(A) := \frac{trace(e^{-\gb H}A)}{trace(e^{-\gb H})}
  \end{equation}
  defines $\gf_\gb$ uniquely. In this case, $\gf=\gf_\gb$ is called a \emph{\gb-KMS state}.
\end{defn}

\begin{remark}\label{rem:Powers'-general-estimate}
  Let \gd be the infinitesimal generator of the flow $\ga:\bR \to Aut(\xA)$ so that $\ga_t = e^{t\gd}$. For all \gb-KMS states \gw, Powers established the following a priori estimate in \cite{Pow76a}:
  \linenopax
  \begin{align}\label{eqn:Powers'-general-estimate}
    \left|\gw([A,B])\right|^2 
    \leq \frac\gb2 \gw(A^\ad A + AA^\ad) \gw(-\ii[B^\ad,\gd(B)]) 
  \end{align}
  for all $A,B \in \xA$ and $B \in \dom \gd$.
\end{remark}

\begin{remark}[Long-range order vs. phase transitions]
  \label{rem:Long-range-order-vs.-phase-transitions}
  It is excruciatingly important to notice that when \xA is infinite-dimensional, formula \eqref{eqn:beta-KMS} becomes meaningless, as was discovered by Bob Powers in his Ph.D. Dissertation \cite{Pow67}; see also \cite{BrRo97}. The reason for this is somewhat subtle: KMS states should really be formulated in terms the representation of \xA obtained via GNS construction (see \S\ref{sec:GNS-construction}). Thus, each occurrence of $A$ in Definition~\ref{def:beta-KMS-state} should be replaced by $\gp_\gf(A)$ if we are being completely honest. However, in the finite-dimensional case, one can use the identity representation and recover \eqref{eqn:beta-KMS} as it reads above.  Unfortunately, the von Neumann algebras generated by KMS states are almost always type III, i.e., the double commutant $\gp_\gf(\xA)''$ typically \emph{does not have a trace} (even though the $C^\ad$-algebra \xA always does). The von Neumann algebra is the weak-$\ad$ closure of the representation (obtained via GNS construction) of the $C^\ad$-algebra; this connection is expressed in the notation of \S\ref{sec:GNS-construction} by the identity
  \linenopax
  \begin{align*}
    \gf( A \ga_t(B))
    &= \la \gp_\gf(A^\ad) \gz_\gf, \gp_\gf(\ga_t(B)) \gz_\gf \ra_{\sH_\gf} \\
    &= \la \gp_\gf(A^\ad) \gz_\gf, e^{-\ii t \sH_\gf} \gp_\gf(B) e^{\ii t \sH_\gf} \gz_\gf \ra_{\sH_\gf},
  \end{align*}
  where now $\sH_\gf$ in the exponent is an unbounded self-adjoint operator in the Hilbert space of the GNS representation derived from the state \gf as in \S\ref{sec:GNS-construction}.

  As a consequence of the lack of trace described just above, there is no uniqueness for $\gf_\gb$ in general, and this has an important physical interpretation in terms of phase transitions. The parameter \gb is inverse temperature: $\gb = 1/kT$ where $k$ is Boltzmann's constant and $T$ is temperature in degrees Kelvin. Whenever \gb is a number for which the set of \gb-KMS states contains more than one element, one says that \gb corresponds to a phase transition; i.e., $T=1/k\gb$ is a temperature at which more than one equilibrium state can exist. Conversely, ``when the system is heated, all is vapor,'' and we expect that the equilibrium state \gf is then unique for $\gb = 1/kT \approx 0$. The lowest $T$ for which multiplicity exceeds 1 is called the critical temperature; it is found experimentally but rigorous results are hard to come by. Indeed, the phase-transition problem in rigorous models is notoriously extremely difficult. Instead the related long-range order problem (as described just below) is thought to be more amenable to computations.
\end{remark}

To get a feel for why KMS states must exist, consider the following construction. Suppose we begin with a finite set $F \ci \verts$ and the corresponding truncated Hamiltonian
\linenopax
\begin{align*}
  H_F := \tfrac12 \sum_{x,y \in F} \cond_{xy} \gf(\one - \gs_x \cdot \gs_y) \in \xA(F).
\end{align*}
Observe that $\xA(F)$ is finite-dimensional; for spin observables with spin $s$, for example, $\dim(\xA_x) = 2s+1$ for each $x$. Here, $\xA_x$ is a subalgebra of the matrices $\sM_{2s+1}(\bC)$. Consequently,
\linenopax
\begin{align*}
  \gf_\gb^F(A) := \frac{trace(e^{-\gb H_F}A)}{trace(e^{-\gb H_F})}
\end{align*}
is a well-defined and unique \gb-KMS state. If we now let $F \to \verts$, then $\gf_\gb$ is a \gb-KMS state also. However, $\gf_\gb$ exists as a weak-$\ad$ limit and hence is not unique!

\section{Remarks and references}
\label{sec:Remarks-and-References-magnetism}

The material in this chapter is based primarily on papers by Powers. Our presentation draws on the resistance estimates derived in the previous chapter for lattice models. The best introduction to this chapter is the paper \cite{Pow76b}, and the books \cite{Rue69} and \cite{BrRo97}. Of the results in the literature of relevance to the present chapter, the references \cite{CoMa070, Con07, Han96, Lig99, Lig95, Lig93} are especially relevant. See also Powers \cite{Pow76a, Pow75}. 

Of the work on infinite spin systems, we are influenced by profound work of Thomas Liggett (infinite spin-models \cite{Lig93, Lig95, Lig99}), and by Robert T. Powers: his use of resistance distance in the estimation of long-range order in quantum statistical models \cite{Pow75, Pow76a, Pow76b, Pow78, Pow79}. The work of Liggett is the classical case and it is more directly connected with estimation of metric distances for statistical models. In contrast, Powers deals with quantum statistical lattice spin-models, and in this case the role of the weighted graphs and their resistance metrics is more subtle, see e.g., Theorem~\ref{thm:Powers'-correlation-estimate} above.

%% file: future-directions.tex

\chapter{Future directions}
\label{sec:future-directions}

\headerquote{The bottom line for mathematicians is that the architecture has to be right. In all the mathematics that I did, the essential point was to find the right architecture. It's like building a bridge. Once the main lines of the structure are right, then the details miraculously fit. The problem is the overall design.}{---~-F.~Dyson}

\headerquote{An expert is a man who has made all the mistakes, which can be made, in a very narrow field.}{---~N.~Bohr}

\begin{remark}\label{rem:fut:boundary}
  We have done some groundwork in \S\ref{sec:the-boundary} for the formal construction of the boundary of an infinite resistance network, however there is much more to be done. The development of this boundary theory is currently underway in \cite{bdG}, where we make explicit the connections between our boundary, Martin boundary, and the theory of graph ends. As in this book, the notions of dipoles, monopoles, and harmonic functions play key roles.
\end{remark}

\begin{remark}\label{rem:fut:fractals}
  In \cite{JoPea10b}, we attempt to apply some results of the present investigation to the theory of fractal analysis.\footnote{Finally! If you remember from the introduction, this was our initial aim!} For now, we just show that the resistance distance as defined by \eqref{eqn:def:R(x,y)-Lap} extends to the context of analysis on \emph{PCF self-similar fractals}. The reader is referred to the definitive text \cite{Kig01} and the excellent tutorial \cite{Str06} for motivation and definitions.

  Suppose that \sF is a post-critically finite (PCF) self-similar set with an approximating sequence of graphs $\Graph_1, \Graph_m, \dots \Graph_m, \dots$ with $\Graph = \bigcup_m \verts_m$ and \sF is the closure of \Graph in resistance metric (which is equivalent to closure in Euclidean metric; see \cite[(1.6.10)]{Str06}). The definition of PCF can be found in \cite[Def.~1.3.4 and Def.~1.3.13]{Kig01}. In the following proof, the subscript $m$ indicates that the relevant quantity is computed on the corresponding \ERN $(\Graph_m,R_m)$. For example, $\Pot_m(x,y)$ is the set of dipoles on $\Graph_m$ (cf.~Definition~\ref{eqn:def:dipole}) and $\energy_m(u)$ is the appropriately renormalized energy of a function $u:\verts_m \to \bR$ (cf.~\cite[(1.3.20)]{Str06}).

  \begin{theorem}\label{thm:R(x,y)-Lap-fractal}
    For $x,y \in \sF$, the resistance distance is given by
    \linenopax
    \begin{equation}\label{eqn:R(x,y)-Lap-fractal}
      \min \{v(x)-v(y) \suth v \in \dom \energy, \Lap v = \gd_x-\gd_y\}.
    \end{equation}
    \begin{proof}
      By the definition of \sF, it suffices to consider the case when $x,y$ are junction points, that is, $x,y \in \Graph_m$ for some $m$. Then the proof follows for general $x,y \in \sF$ by taking limits.

      For $x,y \in \Graph_m$, let $v=v_m$ denote the element of $\Pot_m(x,y)$ of minimal energy; the existence and uniqueness of $v_m$ is justified by the results of \S\ref{sec:potential-functions}. From Theorem~\ref{thm:effective-resistance-metric} we have $R_m(x,y) = \energy(v_m)$.
      Next, apply the harmonic extension algorithm to $v_m$ to obtain $v_{m+1}$ on $\Graph_{m+1}$. By \cite[Lem.~1.3.1]{Str06},
      \linenopax
      \begin{align*}
        v_m(x)-v_m(y) = R_m(x,y) = \energy(v_m) = \energy(v_{m+1}) = \dots = \energy(\tilde v),
      \end{align*}
      where $\tilde v$ is the harmonic extension of $v$ to all of \Graph. It is clear by construction and the cited results that $\tilde v$ minimizes \eqref{eqn:R(x,y)-Lap-fractal}. Note that we do not need to worry about the possible appearance of nontrivial harmonic functions, as $\tilde v$ is constructed as a limit of functions with finite support.
    \end{proof}
  \end{theorem}

  This theorem offers a practical improvement over the formulation of resistance metric as found in the literature on fractals in a couple of respects:
  \begin{enumerate}
    \item \eqref{eqn:R(x,y)-Lap-fractal} provides a formula (or at least, an equation to solve) for the explicit function which gives the minimum in \cite[(1.6.1) or (1.6.2)]{Str06}.
    \item One can compute $v=v_m$ on $\Graph_m$ by basic methods, i.e., Kirchhoff's law and the cycle condition. To find $R(x,y)$, one need only evaluate $v(x)-v(y)$, and this may be done without even fully computing $v$ on all of $\Graph_m$.
  \end{enumerate}
\end{remark}

\version{
\begin{remark}\label{rem:spectral-reciprocity}
  The authors have uncovered a form of spectral reciprocity relating the Laplacian to the matrix $[\la v_x,v_y\ra_\energy]$. This topic is currently under investigation in \cite{SRAMO}.
\end{remark}
}{}

\begin{remark}\label{rem:quantum-graphs}
  The metric graphs and their analysis presented in this volume are ubiquitous, and we can not do justice to the vast literature. However, the application to quantum communication appears especially intriguing, and we refer to the following papers for detail: \cite{NeBr08,DuBr07,Fab06,GHB05,HCDB07}.
  \begin{quote}
    ``$\dots$ the computational power of an important class of quantum states called graph states, representing resources for measurement-based quantum computation, is reflected in the expressive power of (classical) formal logic languages defined on the underlying mathematical graphs.'' \q from \cite{NeBr08}.
  \end{quote}
  \emph{Quantum graphs} (also called \emph{cable systems} in \cite{Kig03} and \emph{graph refinements} in \cite{Telcs06a}) are essentially a refinement of \ERNs where the edges are replaced by intervals and functions are allowed to vary continuously for different values of $x$ in a single edge.
\end{remark}

\begin{remark}\label{rem:fut:rank-Pperp-is-numcycles}
  As noted in Remark~\ref{rem:rank-Pperp-is-numcycles}, the rank of $\Pdrp^\perp$ is an invariant related to the space of cycles in \Graph. However, this object is rather a blunt tool, and it would interesting to see if one can obtain a more refined analysis by applying extensions of the techniques Terras and Stark, as in  \cite{GIL06b}, for example.
\end{remark}



%% file: functional-analysis.tex

\chapter{Some functional analysis}
\label{sec:functional-analysis}

Since this presentation addresses disparate audiences, we found it helpful to organize tools from functional analysis and the theory of unbounded operators in appendix sections. The reader may find the references \cite{Rud91, KadisonRingroseI, DuSc88} to be helpful.


The magic of Hilbert space resides in the following innocent-looking axioms:
\begin{enumerate}
  \item A complex vector space \sH.
  \item A complex-valued function on $\sH \times \sH$, denoted by $\la \cdot,\cdot\ra$ and satisfying
  \begin{enumerate}[(a)]
    \item For every $v \in \sH$, $\la v, \cdot\ra : \sH \to \bC$ is linear.
    \item For every $v_1,v_2 \in \sH$, $\la v_1,v_2\ra = \cj{\la v_2,v_1\ra}$.
    \item For every $v \in \sH$, $\la v,v\ra \geq 0$, with equality if and only if $v=0$.
  \end{enumerate}
  \item Under the norm defined by $\|v\|_\sH := \la v,v\ra^{1/2}$, \sH is complete.
\end{enumerate}

\begin{exm}[Square-summable sequences]
  For a \bC-valued function $v$ on the integers, let $\|v\|_2 := \left(\sum_{n \in \bZ} |v(n)|^2\right)^{1/2}$. Then
  \[\ell^2(\bZ) := \{v:\bZ \to \bC \suth \|v\|_2 < \iy\}\]
  is a Hilbert space.
\end{exm}

\begin{exm}[Classical $L^2$-spaces]
  For a measurable \bC-valued function $v$ on a measure space $(X,\gm)$, let $\|v\|_2 := \left(\int_X |v(x)|^2 \,d\gm(x) \right)^{1/2}$. Then
  \[L^2(\gm) := \{v:X \to \bC \suth \|v\|_2 < \iy\}\]
  is a Hilbert space.
\end{exm}

\section{von Neumann's embedding theorem}
\label{sec:von-Neumann's-embedding-theorem}

\begin{theorem}[von Neumann]\label{thm:vonNeumann's-embedding-thm}
  Suppose $(X,d)$ is a metric space. There exists a Hilbert space $\sH$ and an embedding $w:(X,d) \to \sH$ sending $x \mapsto w_x$ and satisfying
  \linenopax
  \begin{equation}\label{eqn:vonNeumann's-embedding-thm}
    d(x,y) = \|w_x - w_y\|_\sH
  \end{equation}
  if and only if $d^2$ is negative semidefinite.
\end{theorem}

\begin{defn}\label{def:negative-semidefinite}
  A function $d:X \times X \to \bR$ is \emph{negative semidefinite} iff for any $f:X \to \bR$ satisfying $\sum_{x \in X} f(x) = 0$, one has
  \linenopax
  \begin{equation}\label{eqn:negative-semidefinite}
    \sum_{x,y \in F} f(x) d^2(x,y) f(y) \leq 0,
  \end{equation}
  where $F$ is any finite subset of $X$.
\end{defn}

von Neumann's theorem is constructive, and provides a method for obtaining the embedding, which we briefly describe, continuing in the notation of Theorem~\ref{thm:vonNeumann's-embedding-thm}.

\begin{FlatList}
  \item Schwarz inequality. If $d$ is a negative semidefinite function on $X \times X$, then define a positive semidefinite bilinear form on functions $f,g:X \to \bC$ by
      \linenopax
      \begin{align} \label{eqn:vNeu-inner-prod}
        Q(f,g) = \la f,g\ra_Q := - \sum_{x,y} f(x) d^2(x,y) g(y).
      \end{align}
      One obtains a quadratic form $Q(f) := Q(f,f)$, and checks that the generalized Schwarz inequality holds $Q(f,g)^2 \leq Q(f) Q(g)$ by elementary methods.
  \item The kernel of $Q$. Denote the collection of finitely supported functions on $X$ by $\Fin(X)$ and define
      \linenopax
      \begin{align}\label{eqn:def:Fin0(X)}
        \Fin_0(X) := \{f \in \Fin(X) \suth {\textstyle \sum}_x f(x) = 0\}.
      \end{align}
      The idea is to complete $\Fin_0(X)$ with respect to $Q$, but first one needs to identify functions that $Q$ cannot distinguish. Define
      \linenopax
      \begin{align}\label{eqn:def:ker(Q)}
        \ker Q = \{f \in \Fin_0(X) \suth Q(f)=0\}.
      \end{align}
      It is easy to see that $\ker Q$ will be a subspace of $\Fin_0(X)$.
  \item Pass to quotient. Define $\tilde Q$ to be the induced quadratic form on the quotient space $\Fin_0(X)/\ker Q$. One may then verify that $\tilde Q$ is \emph{strictly positive definite} on the quotient space. As a consequence, $\|\gf\|_{\sH_{vN}} := -\tilde Q(\gf)$ will be a bona fide norm.
  \item Complete. When the quotient space is completed with respect to $\tilde Q$, one obtains a Hilbert space
      \linenopax
      \begin{align}\label{eqn:vNeu's-quotient-completion}
        \sH_{vN} &:= \left(\frac{\Fin_0(X)}{\ker Q}\right)^\sim,
        \q\text{with}\q
        \la \gf, \gy \ra_{\sH_{vN}} = -\tilde Q(\gf,\gy).
      \end{align}
  \item Embed $(X,d)$ into $\sH_{vN}$. Fix some point $o \in X$ to act as the origin; it will be mapped to the origin of $\sH_{vN}$ under the embedding. Then define
      \linenopax
      \begin{align*}
        w:(X,d) \to \sH_{vN} \q\text{by}\q
        x \mapsto w_x := \tfrac1{\sqrt2}(\gd_x - \gd_o).
      \end{align*}
      Now $w$ gives an embedding of $(X,d)$ into the Hilbert space $H_{vN}$, and
      \linenopax
      \begin{align}\label{eqn:embedding-check-app}
        \|w_x - w_y\|_{vN}^2
        &= \la w_x - w_y, w_x - w_y \ra_{vN}  \\
        &= \la w_x, w_x\ra_{vN} - \la w_x, w_y \ra_{vN} + \tfrac12 \la w_y, w_y \ra_{vN} \notag \\
        &= d^2(x,o) + \left(d^2(x,y) - d^2(x,o) - d^2(y,o)\right) + d^2(y,o) \notag \\
        &= d^2(x,y),\notag
      \end{align}
      which verifes \eqref{eqn:vonNeumann's-embedding-thm}. The third equality follows by three computations of the form
      \linenopax
      \begin{align}\label{eqn:crossterm-sums-of-madness}
        \la w_x, w_y\ra_{vN}
        &= -\sum_{a,b} w_x(a) d^2(a,b) w_y(b) \notag \\
        &= -\sum_{a,b} d^2(a,b)  \tfrac1{\sqrt2}(\gd_x(a)-\gd_o(a))(\gd_y(b)-\gd_o(b)) \notag \\
        &= \dots = d^2(x,y) - d^2(x,o) - d^2(y,o),
      \end{align}
      noting that $d(a,a)=0$, etc.
\end{FlatList}

von Neumann's theorem also has a form of uniqueness which may be thought of as a universal property.

\begin{theorem}\label{thm:uniqueness-of-vNeu-embedding}
  If there is another Hilbert space \sK and an embedding $k:\sH \to \sK$, with $\|k_x-k_y\|_\sK = d(x,y)$ and $\{k_x\}_{x \in X}$ dense in \sK, then there exists a unique unitary isomorphism $U:\sH \to \sK$.
  \begin{proof}
    We show that $U:w \mapsto k$ by $U(\sum \gx_x w_x) = \sum \gx_x k_x$ is the required isometric isomorphism. Let $\sum \gx_x = 0$. It is conceivable that $U$ fails to be well-defined because of linear dependency; we show this is not the case:
    \linenopax
    \begin{align}
      \left\|\sum_{x \in X} \gx_x w_x\right\|^2
      &= \sum_{x,y \in X} \gx_x \tilde Q(w_x, w_y) \gx_y \notag \\
      &= \sum_{x,y \in X} \gx_x \left(d^2(x,y) - d^2(x,o) - d^2(y,o) \right) \gx_y \qq\text{by \eqref{eqn:crossterm-sums-of-madness}} \notag \\
      &= \sum_{x,y \in X} \gx_x \gx_y d^2(x,y)
        - \sum_{x \in X} \gx_x d^2(x,o) \cancel{\sum_{y \in X} \gx_y}
        - \sum_{y \in X} \gx_y d^2(y,o) \cancel{\sum_{x \in X} \gx_x} \notag \\
      &= -\sum_{x,y \in X} \gx_x \gx_y d(x,y), \label{eqn:uniqueness-of-vNeu-embedding}
    \end{align}
    since $\sum_x \gx_x = 0$ by choice of \gx. However, the same computation may be applied to $k$ with the same result; note that \eqref{eqn:uniqueness-of-vNeu-embedding} does not depend on $w$. Thus, $\|w\|_\sH = \|k\|_\sK$ and $U$ is an isometry. Since it is an isometry from a dense set in \sH to a dense set in \sK, we have an isomorphism and are finished.
  \end{proof}
\end{theorem}

\version{}{\marginpar{Is this paragraph correct/ useful?}}
The importance of using $\Fin_0(X)$ in the above construction is that the finitely supported functions $\Fin(X)$ are in duality with the bounded functions $B(X)$ via
\linenopax
\begin{align}\label{eqn:vNeu-duality-bracket}
  \la f,\gb\ra := \sum_{x \in X} f(x) \gb(x) < \iy
  \qq f \in \Fin(X), \gb \in B(X).
\end{align}
The constant function $\gb_1 := \one$ is a canonical bounded function. With respect to the pairing in \eqref{eqn:vNeu-duality-bracket}, its orthogonal complement is
\linenopax
\begin{equation*}
  \gb_1^\perp = \{\gf \suth \la \gf,\gb_1\ra = 0\} = \Fin_0(X).
\end{equation*}

\section{Remarks and references}
\label{sec:Remarks-and-References-App-A}

The material here is collected to help make our presentation self-contained. The reader may find the references \cite{Rud91, KadisonRingroseI, DuSc88, ReedSimonII, Arveson:spectral-theory, vN55} to be helpful. The further references \cite{BrRo79}, and \cite{Kat95} may also be useful.

%% file: operator-theory.tex

\chapter{Some operator theory}
\label{sec:operator-theory}

\headerquote{Anyone who attempts to generate random numbers by deterministic means is, of course, living in a state of sin.}{---~J.~von~Neumann}

\begin{defn}\label{def:adjoint}
  If $S:\sH \to \sH$ is an operator on a Hilbert space \sH, its \emph{adjoint} is the operator satisfying $\la S^\ad u, v\ra = \la u, Sv\ra$ for every $v \in \dom S$. The restriction to $v \in \dom S$ becomes significant only when $S$ is unbounded, in which case one sees that the domain of the adjoint is defined by
  \begin{equation}\label{eqn:def:adjoint-domain}
    \dom S^\ad := \{u \in \sH \suth |\la u, S v \ra| \leq k\|v\|, \forall v \in \dom S\},
  \end{equation}
  where the constant $k=k_u \in \bC$ may depend on $u$.

  An operator $S$ is said to be \emph{self-adjoint} iff $S=S^\ad$ \textbf{and} $\dom S = \dom S^\ad$. It is often the equality of domains that is harder to check.
\end{defn}

\version{}{
\begin{lemma}\label{thm:adjoint-of-contractive-is-contractive}
    \marginpar{Something about applying polar decomp to $S$. This may no longer be necessary.}
  If $S$ is a contractive operator, then its adjoint $S^\ad$ is also contractive.
  \begin{proof}
  \end{proof}
\end{lemma}
}

\section{Projections and closed subspaces}
\label{sec:projections-and-closed-subspaces}
Let \sH be a complex (or real) Hilbert space, and define
\linenopax
\begin{align*}
  \sB(\sH) &:= \{A:\sH \to \sH \suth A \text{ is bounded and linear}\} \\
  \sB(\sH)_{sa} &:= \{A \in \sB(\sH) \suth A=A^\ad\} \\
  Cl(\sH) &:= \{V \ci \sH \suth V \text{ is a closed linear space}\}.
\end{align*}

\begin{defn}\label{def:projection}
  An operator $P$ on a Hilbert space \sH is a \emph{projection} iff it \sats $P = P^2 = P^\ad$. Denote the space of projections by $Proj(\sH)$.
\end{defn}

\begin{theorem}[Projection Theorem]\label{thm:Projection-Theorem}
  There is a bijective correspondence between the set $Cl(\sH)$ of all closed subspaces $V \ci \sH$ and the set of all projections $P$ acting on \sH.
\end{theorem}
Specifically, if a subspace $V \in Cl(\sH)$ is given, there is a unique projection $P \in Proj(\sH)$ with range $L$. Conversely, if a projection $P$ is given, set $V := P\sH$. Moreover, the mapping
\linenopax
\begin{align*}
  Proj(\sH) \to Cl(\sH), \qq\text{by}\qq P \mapsto V = P\sH
\end{align*}
is a lattice isomorphism. The ordering in $Cl(\sH)$ is defined by containment; and for a pair of projections $P$ and $Q$ we say that $P \leq Q$ iff $P = PQ$. This ordering coincides with the usual order on the self-adjoint elements of the algebra:
\linenopax
\begin{align*}
  A \leq B \text{ in } \sB(\sH) \qq\text{iff}\qq \la v,Av\ra \leq v, Bv \ra, \;\forall v \in \sH,
\end{align*}
and induces an ordering on $Proj(\sH)$. Since $Proj(\sH)$ is a lattice, it follows that $Cl(\sH)$ is a lattice as well.

\begin{theorem}[{$\negsp[4]$\cite[Prop.~2.5.2]{KadisonRingroseI}}]\label{thm:3defs-for-projn-order}
  For projections $P,Q$, the following are equivalent:
  \begin{enumerate}[(i)]
    \item $P\sH \ci Q\sH$.
    \item $P = PQ$.
    \item $\la v,Pv\ra \leq \la v, Qv \ra$, $\forall v \in \sH$.
    \item $\|Pv\| \leq \|Qv \|$, $\forall v \in \sH$.
  \end{enumerate}
\end{theorem}

\section{Partial isometries}
\label{sec:partial-isometries}
Let \sH and \sK be two complex (or real) Hilbert spaces, and let $L: \sH \to \sK$ be a bounded linear operator.

\begin{defn}\label{def:partial-isometry}
  $L$ is a \emph{partial isometry} if one (all) of the following equivalent conditions is satisfied:
  \begin{enumerate}[(i)]
    \item $L^\ad L$ is a projection in \sH (the \emph{initial projection}).
    \item $LL^\ad $ is a projection in \sK (the \emph{final projection}).
    \item $LL^\ad L = L$.
    \item $L^\ad L L^\ad = L^\ad$.
  \end{enumerate}
  In this case, we say that $P_i := L^\ad L$ is the \emph{initial projection} and $P_f := LL^\ad$ is the \emph{final projection}. Moreover, $P_i$ is the \projn onto $\ker(L)^\perp$ and $P_f$ is the \projn onto the closed subspace $\ran L$.
\end{defn}

\begin{theorem}\label{thm:def:partial-isometry}
  The four conditions of Definition~\ref{def:partial-isometry} are equivalent. Consequently, the initial and final projections satisfy $L P_i = L = P_f L$.
  \begin{proof}
    \version{Define $P := L^\ad L$. Compute that $(LP-L)^\ad(LP-L) = 0$ to deduce $LP=L$; the reader can fill in the rest.}
    {(i)$\iff$(iii). Define $P := L^\ad L$. First, suppose $P$ is a \projn, so that $P = P^2 = P^\ad$. Then
    \linenopax
    \begin{align*}
      (LP-L)^\ad(LP-L)
      &= P(L^\ad L)P - P(L^\ad L) - (L^\ad L)P + L^\ad L \\
      &= P^3 - P^2 - P^2 + P
      = 0,
    \end{align*}
    which means $LP=L$. Now, assume $LP=L$. Then
    \linenopax
    \begin{align*}
      P^2 = L^\ad LP = L^\ad L = P,
      \q\text{and}\q
      P^\ad = (L^\ad L)^\ad = L^\ad L = P,
    \end{align*}
    so that $P$ is a \projn.

    (iii)$\iff$(iv). This is immediate by taking adjoints.

    (ii)$\iff$(iv). \emph{Mutatis mutandis}, this is the same proof as (i)$\iff$(iii).}
  \end{proof}
\end{theorem}

\begin{defn}\label{def:isometry}
  The operator $L$ is an \emph{isometry} iff $P_i = \id_\sH$, the identity operator on \sH. The operator $L$ is a \emph{coisometry} iff $P_f = \id_\sK$, the identity operator on \sK. It is clear that $L$ is an isometry if and only if $L^\ad$ is a coisometry.
\end{defn}

\section{Self-adjointness for unbounded operators}
\label{sec:self-adjointness-for-unbounded-operators}
Throughout this section, we use \sD to denote a dense subspace of \sH. Some good references for this material are \cite{vN32a,Rud91,DuSc88}.

\begin{defn}\label{def:Hermitian}
  An operator $S$ on \sH is called \emph{Hermitian} (or \emph{symmetric} or \emph{formally self-adjoint}) iff
  \begin{equation}\label{eqn:def:Hermitian}
    \la u,Sv\ra = \la Su,v\ra, \qq\text{for every } u,v \in \sD.
  \end{equation}
  In this case, the spectrum of $S$ lies in \bR.
\end{defn}

\begin{defn}\label{def:self-adjoint}
  An operator $S$ on \sH is called \emph{self-adjoint} iff it is Hermitian and $\dom S = \dom S^\ad$.
\end{defn}

\begin{defn}\label{def:semibounded}
  An operator $S$ on \sH is called \emph{semibounded} iff
  \begin{equation}\label{eqn:def:semibounded}
    \la v,Sv\ra \geq 0, \qq\text{for every } v \in \sD.
  \end{equation}
  The spectrum of a semibounded operator lies in some halfline $[\gk,\iy)$ and the defect indices of a semibounded operator always agree (see Definition~\ref{def:defect-index}). The graph Laplacian \Lap considered in much of this book falls into this class.
\end{defn}

\begin{defn}\label{def:bounded}
  An operator $S$ on \sH is called \emph{bounded} iff there exists $k \in \bR$ \st
  \begin{equation}\label{eqn:def:bounded}
    |\la v, Sv\ra| \leq k\|v\|^2, \qq\text{for every } v \in \sD.
  \end{equation}
  The spectrum of a bounded operator lies in a compact subinterval of \bR.  Bounded Hermitian operators are automatically self-adjoint. When $(\Graph, \ohm)$ satisfies the Powers bound \eqref{eqn:def:power-bound}, the transfer operator falls into this class.
\end{defn}

\begin{defn}\label{def:graph-of-operator}
  For an operator $S$ on the Hilbert space \sH, the \emph{graph} of $S$ is
  \linenopax
  \begin{equation}\label{eqn:def:graph-of-operator}
    \opGraph(S) := \{\left[\begin{smallmatrix}v \\Av\end{smallmatrix}\right] \suth v \in \sH\} \ci \sH \oplus \sH,
  \end{equation}
  with the norm
  \linenopax
    \begin{align}\label{eqn:def:graph-norm-of-operator}
    \left\|\left[\begin{smallmatrix}v \\Sv\end{smallmatrix}\right] \right\|_{Graph}^2
    := \|v\|_\sH^2 + \|S v\|_\sH^2
  \end{align}
  and the corresponding inner product. The operator $S$ is \emph{closed} iff $\opGraph(S)$ is closed in $\sH \oplus \sH$ or \emph{closable} if the closure of $\opGraph(S)$ is the graph of an operator. In this case, the corresponding operator is $\opclosure S$, the \emph{closure} of $S$. The domain of $\opclosure S$ is therefore defined
  \linenopax
    \begin{align}\label{eqn:def:domain-of-clo(S)}
    \dom \opclosure S := \{u \suth \lim_{n \to \iy}\|u-u_n\|_\sH = \lim_{n \to \iy}\|v-Su_n\|_\sH = 0\}
  \end{align}
  for some $v \in \sH$ and Cauchy sequence $\{u_n\} \ci \dom S$. Then one defines $\opclosure{S}u := v$.
  If $S$ is Hermitian, then $\opclosure{S}$ will also be Hermitian, but it will not be self-adjoint in general.
\end{defn}

\begin{remark}\label{rem:closable-iff-adjoint-densely-defined}
  It is important to observe that an operator $S$ is closable if and only if $S^\ad$ has dense domain. However, this is clearly satisfied when $S$ is Hermitian with dense domain, since then $\dom S \ci \dom S^\ad$.
\end{remark}

\begin{defn}\label{def:operator-graph-rotation}
  \version{}{\marginpar{There is probably a better name?}}
  Suppose that $S$ is a linear operator on \sH with a dense domain $\dom S$. Define the \emph{graph rotation operator} $\xG:\sH \oplus \sH \to \sH \oplus \sH$ by $\xG(u,v) := (-v,u)$.
  It is easy to show that the graph of $S^\ad$ is
  \begin{equation}\label{eqn:operator-graph-rotation}
    \opGraph(S^\ad) = (\xG(\opGraph(S)))^\perp.
  \end{equation}
\end{defn}

For any semibounded operator $S$ on a Hilbert space, there are unique self-adjoint extensions $S_{\mathrm{min}}$ (the Friedrichs extension) and $S_{\mathrm{max}}$ (the Krein extension) such that
\begin{equation}\label{eqn:S_{min}<S_{max}}
  S \ci \opclosure S \ci S_{\mathrm{min}} \ci \tilde{S} \ci S_{\mathrm{max}},
\end{equation}
 where $\tilde{S}$ is any non-negative self-adjoint extension of $S$. \emph{For general unbounded operators, these inclusions may all be strict.} In \eqref{eqn:S_{min}<S_{max}}, $A \ci B$ means graph containment, i.e., it means $\opGraph(A) \ci \opGraph(B)$, where $\opGraph$ is as in Definition~\ref{def:graph-of-operator}. The case when $S_{\mathrm{min}} = S_{\mathrm{max}}$ is particularly important.

\begin{defn}\label{def:essentially-self-adjoint}
  An operator is defined to be \emph{essentially self-adjoint} iff it has a unique self-adjoint extension. An operator is essentially self-adjoint if and only if it has defect indices 0,0 (see Definition~\ref{def:defect-index}). A self-adjoint operator is trivially essentially self-adjoint.
\end{defn}

\begin{theorem}\label{thm:Hermitian-properties}\emph{\cite{vN32a,Rud91,DuSc88}}
  Let $S$ be a Hermitian operator.
  \begin{enumerate}
    \item $S$ is closable, its closure $\opclosure S$ is Hermitian, and $S^\ad = (\opclosure S)^\ad$.
    \item Every closed Hermitian extension $T$ of $S$ satisfies
      \begin{align*}
        S \ci T \ci T^\ad \ci S^\ad.
      \end{align*}
    \item $S$ is essentially self-adjoint if and only if $\dom (\opclosure{S}) = \dom S^\ad$.
    \item $S$ is essentially self-adjoint precisely when both its defect indices are 0.
    \item $S$ has self-adjoint extensions iff $S$ has equal defect indices.
  \end{enumerate}
\end{theorem}

The Hermitian operator $S := QPQ$ of Example~\ref{exm:banded-nonselfadjoint-matrix} has defect indices 1,1, and yet is not even semibounded.

\begin{defn}\label{def:defect-index}
  Let $S$ be an operator with adjoint $S^\ad$. For $\gl \in \bC$, define
  \linenopax
  \begin{equation}\label{eqn:defect-space}
    \Def_\gl(S)
    := \ker(\gl - S^\ad)
    = \{v \in \dom S^\ad \suth S^\ad v = \gl v\}.
  \end{equation}
  Then $\Def_\gl(S)$ is the \emph{defect space} of $S$ corresponding to \gl. Elements of $\Def_\gl(S)$ are called \emph{defect vectors}. The number $\dim \Def_\gl(S)$ is constant in the connected components of the resolvent set $\bC \less \gs(S)$ and is called the \emph{defect index} of the component containing \gl.
\end{defn}

Note that if $S$ is Hermitian, then its resolvent set can have at most two connected components. Further, if $S$ is semibounded, then its resolvent set can have only one connected component, and we have only one defect index to compute: the dimension of $\Def(S) = \Def_{-1}(S)$. These facts explain the two consequences of the following theorem, which can be found in \cite{vN32a,Rud91,DuSc88}.

\begin{theorem}[von Neumann]
  \label{thm:essentially-self-adjointness-criterion}
  For a Hermitian operator $S$ on \sH, one has
  \linenopax
  \begin{align*}
    \dom S^\ad = 
    \dom \opclosure{S} 
      \oplus \{v \in \sH \suth S^\ad v = \ii v\}
      \oplus \{v \in \sH \suth S^\ad v = -\ii v\},
  \end{align*}
  where the orthogonality of the direct sum on the right-hand side is with respect to the graph inner product \eqref{eqn:def:graph-norm-of-operator} (not the inner product of \sH).
  Consequently, $S$ is essentially self-adjoint if and only if
  \linenopax
  \begin{equation}\label{eqn:thm:essentially-self-adjointness-criterion-Herm}
    S^\ad v = \pm \ii v \q\implies\q v=0, \qq v \in \sH.
  \end{equation}
  For a semibounded Hermitian operator $S$ on \sH, one has
  \linenopax
  \begin{align*}
    \dom S^\ad = \dom \opclosure{S} \oplus \{v \in \sH \suth S^\ad v = -v\},
  \end{align*}
  where again the orthogonality of the direct sum on the right-hand side is with respect to the graph inner product \eqref{eqn:def:graph-norm-of-operator}.
  Consequently, $S$ is essentially self-adjoint if and only if
  \linenopax
  \begin{equation}\label{eqn:thm:essentially-self-adjointness-criterion-semi}
    S^\ad v = -v \q\implies\q v=0, \qq v \in \sH.
  \end{equation}
\end{theorem}
A solution $v$ of \eqref{eqn:thm:essentially-self-adjointness-criterion-Herm} or \eqref{eqn:thm:essentially-self-adjointness-criterion-semi} is called a \emph{defect vector} (as in Definition~\ref{def:defect-index}) or an \emph{vector at \iy}.
The idea of the proof in von Neumann's theorem is to obtain the essential self-adjointness of a Hermitian operator $S$ by using the following stratagem: an unbounded function applied to a bounded self-adjoint operator is an unbounded self-adjoint operator. In this case, the function is $f(x) = \gl - x^{-1}$. If we can see that $(\gl - S_{\mathrm{min}})^{-1}$ is bounded and self-adjoint, then
\linenopax
\begin{align*}
  f((\gl - S_{\mathrm{min}})^{-1}) = S_{\mathrm{min}}
\end{align*}
is an unbounded self-adjoint operator.
First, note that
\linenopax
\begin{align*}
  \clspn {\ran (\bar{\gl} - S)} = \ran(\bar \gl - S_{\mathrm{min}}) = \ker(\gl - S^\ad)^\perp = \Def_\gl(S)^\perp.
\end{align*}
Note that if $\gl \in \operatorname{res}(S)$ and $(\gl - S_{\mathrm{min}})^{-1} = \int_\bR x \,E(dx)$ with projection-valued measures $E:\sB(\bR) \to Proj(\sH)$, then
\linenopax
\begin{align*}
  S_{\mathrm{min}} = \int_\bR (\gl - x^{-1})\, E(dx).
\end{align*}
This will show $S_{\mathrm{min}}$ is self-adjoint; if $\Def_\gl(S)=0$ for ``enough'' \gl, then $S_{\mathrm{min}}$ is self-adjoint and hence $S$ is essentially self-adjoint.

\begin{lemma}\label{thm:bounded-and-Hermitian-implies-self-adjoint}
  If $S$ is bounded and Hermitian, then it is essentially self-adjoint.
  \begin{proof}
    $S$ is bounded iff it is everywhere defined, by the Hellinger-Toeplitz Theorem. Since $S^\ad$ is also everywhere defined in this case, it is clear the two operators have the same domain.
  \end{proof}
\end{lemma}

\begin{lemma}\label{thm:Def-perp-is-range-of-(1+opclosure)}
  For an operator $S$ which is semibounded but not necessarily closed, 
  \linenopax
  \begin{align}\label{eqn:Def-perp-is-range-of-opclosure}
    \Def(S)^\perp = \ran(1+\opclosure{S}). 
  \end{align}
  \begin{proof}
    Recall that $\Def(S) = \Def_{-1}(S)$. General theory gives
    \linenopax
    \begin{align*}
      \Def(S)^\perp 
      = (\ker 1+S^\ad)^\perp  
      = (\ran (1+S))^{\perp \perp}
      = (\ran (1+S))^{\textrm{clo}}.
    \end{align*}
    It remains to check that $\ran (1+\opclosure{S}) = (\ran (1+S))^{\textrm{clo}}$. If $S$ is not semibounded, one may have only the containment $(\ci)$: note that $u \in \ran(1 +\opclosure{S})$ iff $u = v + Sv$ for $v \in \dom \opclosure{S}$. Then $v = \lim v_n$ for some $v_n \in \dom S$, and the containment is clear. For $(\ce)$, let $u \in \opclosure{\ran(1+S)}$ so that $u = \lim u_n$, where $u_n = v_n Sv_n$ for $v_n \in \dom S$, and note that
    \linenopax
    \begin{align}\label{eqn:1+S-is-contractive}
      \|v_n\|^2 \leq \|v_n\|^2 + \la Sv_n,v_n\ra + \la v_n, Sv_n\ra + \|Sv_n\|^2 
      = \|v_n + Sv_n\|^2
      =\|u_n\|^2,
    \end{align}
    where the inequality uses the fact that $S$ is semibounded. By passing to a subsequence if necessary, \eqref{eqn:1+S-is-contractive} implies $\|v_n - v_m\| \leq \|u_n - u_m\|$, whence $\{u_n\}$ Cauchy implies $\{v_n\}$ is also Cauchy and hence has a limit $v$. Therefore
    \linenopax
    \begin{align*}
      S v_n = u_n - v_n \limas{n} u-v,
    \end{align*}
    which allows one to define $\opclosure S v = u-v$ and see $u \in \ran(1+\opclosure{S})$.
  \end{proof}
\end{lemma}

\begin{exm}[The defect of the Laplacian on $(0,\iy)$]
  \label{exm:eigenvectors-at-infinity}
  Probably the most basic example of defect vectors (and how an Hermitian operator can fail to be essentially self-adjoint) is provided by the Laplace operator $\Lap = -\frac{d^2}{dx^2}$ on the Hilbert space $\sH = L^2(0,\iy)$. Example~\ref{exm:banded-nonselfadjoint-matrix} gives an even more striking (though less simple) example. We take \Lap as having the dense domain
  \linenopax
  \begin{align*}
    \sD = \{f \in C_0^\iy(0,\iy) \suth f^{(k)}(0) = \lim_{x \to \iy} f^{(k)}(x) = 0, k=0,1,2,\dots\}.
  \end{align*}
  We always have $\la u,\Lap u\ra \geq 0$ for $u \in \sD$. However, $\Lap^\ad u = - u$ is satisfied by $e^{-x} \in \sH \less \sD$. To see this, take any test function $\gf \in \sD$ and compute the weak derivative via integration by parts (applied twice):
  \linenopax
  \begin{align*}
    \la e^{-x}, \Lap \gf \ra
    = \left\la e^{-x}, -\tfrac{d^2}{dx^2} \gf \right\ra
    = - (-1)^2 \left\la \tfrac{d^2}{dx^2} e^{-x}, \gf \right\ra
    = - \la e^{-x}, \gf \ra.
  \end{align*}
  Thus the domain of $\Lap^\ad$ is strictly larger than $\tilde{\Lap}$, and so $\tilde{\Lap}$ fails to be self-adjoint.
  
  \pgap

  One might try the approximation argument used to prove essential self-adjointness of the Laplacian on $\ell^2(\cond)$ in  Theorem~\ref{thm:Laplacian-is-essentially-self-adjoint}: let $\{v_n\} \ci \sD$ be a sequence with $\|v_n - e^{-x}\|_\cond \to 0$. Since $\Lap^\ad$ agrees with \Lap when restricted to \sD,
  \linenopax
  \begin{align*}
    \la v_n, \Lap v_n \ra_\cond
    = \la v_n, \Lap^\ad v_n \ra_\cond
    \to \la e^{-x}, \Lap^\ad e^{-x} \ra_\cond.
  \end{align*}
  However, there are two mistakes here. First, one does not have convergence unless the original sequence is chosen so as to approximate $v$ in the nonsingular quadratic form
  \linenopax
  \begin{align*}
    \la u,v \ra_{\Lap^\ad} := \la u,v \ra + \la u, \Lap^\ad v \ra.
  \end{align*}
  Second, one \emph{cannot} approximate $e^{-x}$ with respect to this nonsingular quadratic form by elements of \sD. In fact, $e^{-x}$ is orthogonal to \sD in this sense:
  \linenopax
  \begin{align*}
    \la \gf, e^{-x}\ra_{\Lap^\ad}
    &= \la \gf, e^{-x}\ra_\cond + \la \gf, \Lap e^{-x}\ra_\cond
    = \la \gf, e^{-x}\ra_\cond - \la \gf, e^{-x}\ra_\cond =0.
  \end{align*}

  Alternatively, observe that von Neumann's theorem (Theorem~\ref{thm:essentially-self-adjointness-criterion}), a general element in the domain of $\Lap^\ad$ is $v+\gf_\iy$, where $v \in \dom \Lap$ and $\gf_\iy$ is in the defect space.

  Suppose we have an exhaustion $\{\sH_k\}$ with $\sH_k \ci \sH_{k+1} \ci \sD$ and $\sH = \bigvee \sH_k$ (this notation indicates closed linear span of the union).
  \linenopax
  \begin{align*}
    \sH_{n,m} = \clspn{x^k e^{-(x+1/x)} \suth -n \leq k \leq m, \text{for } n,m \in \bN}.
  \end{align*}
  \version{}{\marginpar{Notes say ``the axiom that is required is that \Lap shifts the space, i.e., \Lap maps one $\sH_k$ into a bigger one''. What axiom for what? And shouldn't it be $\Lap \sH_{n,m} \ci \sH_{n+4,m}$?}}
  We have that $\Lap \sH_{n,m} \ci \sH_{n+3,m}$, since
  \linenopax
  \begin{align*}
    \tfrac{d^2}{dx^2} x^k e^{-(x+1/x)}
    &= x^k e^{-(x+1/x)} \left(\tfrac1{x^4} + \tfrac{2(k-1)}{x^3} + \tfrac{-2 -k + k^2}{x^2} - \tfrac{2k}x + 1\right).
  \end{align*}
  However, $\Lap^\ad u = - u$ is satisfied by $e^{-x} \in \sH \less \sD$.
\end{exm}

\section{Banded matrices}
\label{sec:banded-matrices}

\begin{defn}\label{def:banded-matrix}
  Consider the matrix $M_S$ corresponding to an operator $S$ in some ONB $\{b_x\}$, so the entries of $M_S$ are given by
  \linenopax
  \begin{equation}\label{eqn:def:banded-matrix}
    M_S(x,y) := \la b_x, Sb_y\ra.
  \end{equation}
  We say $M_S$ is a \emph{banded matrix} iff every row and column contains only finitely many nonzero entries. A fortiori, $M_S$ is \emph{uniformly banded} if no row or column has more than $N$ nonzero entries, for some $N \in \bN$.
\end{defn}
With $M_S$ defined as in \eqref{eqn:def:banded-matrix}, it is immediate that $M_S$ is Hermitian whenever $S$ is:
\linenopax
\begin{align}\label{eqn:banded-matrix-is-Hermitian}
  M_S(x,y) = \la b_x, S b_y\ra = \la S b_x, b_y\ra = \cj{\la b_y, Sb_x\ra} = \cj{M_S(y,x)}.
\end{align}

Banded matrices are of interest in the present context because the graph Laplacian is always a banded matrix in virtue of the fact that each vertex has only finitely many neighbours; recall the form of $M_\Lap$ given in \eqref{eqn:L2-Laplacian-is-banded}. Since $\Lap = \cond - \Trans$, the transfer operator \Trans is also banded. In general, the bandedness of an operator does not imply the operator is self-adjoint. In fact, see Example~\ref{exm:banded-nonselfadjoint-matrix} for a Hermitian operator on $\ell^2$ which is not self-adjoint, despite having a \emph{uniformly banded} matrix. However, this property does make it much easier to compute the adjoint.

\begin{lemma}\label{thm:bandedness-sublemma}
  Let $S$ be an unbounded Hermitian operator on \sH with dense domain of definition $\sD = \dom S \ci \sH$. Suppose that the matrix $M_S$ defines as in \eqref{eqn:def:banded-matrix} is banded with respect to the ONB $\{b_x\}_{x \in X}$, and define $\hat v(x) := \la b_x, v\ra$. Then $v \in \dom S^\ad$ and $S^\ad v = w$ if and only if $\hat v, \hat w \in \ell^2(X)$ and $\hat w(x)$ is
  \linenopax
  \begin{align}
    \hat w(x) = \sum_{y \in X} M_S(x,y) \hat v(y).
      \label{eqn:lem:banded-representation-of-Lap-adjoint:claim}
  \end{align}
  Thus, $S^\ad$ is represented by the banded matrix $M_S$.
  \begin{proof}
    \fwd  To see the form of $\hat w$ in \eqref{eqn:lem:banded-representation-of-Lap-adjoint:claim},
    \linenopax
    \begin{align*}
      \sum_{y \in X} M_S(x,y) \hat v(y)
      = \sum_{y \in X} \la b_x, S b_y\ra \la b_y, v \ra
      &= \sum_{y \in X} \la S b_x, b_y\ra \la b_y, v \ra &&\text{Hermitian} \\
      &= \la S b_x, v\ra  &&\text{Parseval} \\
      &= \la b_x, w\ra  && S^\ad v = w.
    \end{align*}
    where the last equality is possible since $v \in \dom S^\ad$. It is the hypothesis of bandedness that guarantees all these sums are finite, and hence meaningful.

    Conversely, first note that it is the hypothesis of bandedness which makes the sum in \eqref{eqn:lem:banded-representation-of-Lap-adjoint:claim} finite, ensuring $\hat w$ is well-defined. suppose \eqref{eqn:lem:banded-representation-of-Lap-adjoint:claim} holds, and that $\hat v, \hat w \in \ell^2(X)$. To show $v \in \dom S^\ad$, we must find a constant $K < \iy$ for which $|\la v, S u\ra| \leq K \|u\|$ for every $u \in \sD$.
    \linenopax
    \begin{align*}
      \la v, S u \ra
      &= \sum_{x \in X} \la v, b_x\ra\la b_x, S u\ra
        &&\text{Parseval} \\
      &= \sum_{x \in X} \la v, b_x\ra\la S b_x, u\ra \\
      &= \sum_{x \in X} \cj{\hat v(x)} \left\la \sum\nolimits_y M_S(y,x) b_y, u\right\ra \\
      &= \sum_{y \in X} \sum_{x \in X} \cj{\hat v(x)} \cj{M_S(y,x)} \hat u(y)
        &&\text{by \eqref{eqn:banded-matrix-is-Hermitian}}\\
      |\la v, S u\ra|^2
      &\leq \sum_{y \in X} \left|\sum_{x \in X} \cj{\hat v(x)} M_S(x,y) \right|^2 \sum_{y \in X} |\hat u(y)|^2
      &&\text{by Schwarz},
    \end{align*}
    and see that we can take $K = \|\hat w\|^2$.
  \end{proof}
\end{lemma}

\begin{lemma}\label{thm:bandedness-gives-norm-bound}
  Let $A$ be an operator on $\ell^2(\bZ)$ whose matrix $M_A$ is uniformly banded, with all bands having no more than \gb nonzero entries. Then
  \linenopax
    \begin{equation}\label{eqn:thm:bandedness-gives-norm-bound}
    \|A\| \leq \gb \sup_{x,y} |a_{xy}|.
  \end{equation}
  \begin{proof}
    The Schwarz inequality gives
    \linenopax
    \begin{align*}
      \|A\| \leq \max \left\{\sup\nolimits_x \left(\sum\nolimits_y |a_{xy}|^2\right)^{1/2}, \sup\nolimits_y \left(\sum\nolimits_x |a_{xy}|^2\right)^{1/2}\right\}
    \end{align*}
    However, uniform bandedness gives
    \linenopax
    \begin{align*}
      \sup_x \left(\sum_y |a_{xy}|^2\right)^{1/2} \leq \gb \sup_x \max_y |a_{xy}|,
    \end{align*}
    and similarly for the other term.
  \end{proof}
\end{lemma}

\version{}{
\begin{prop}\label{thm:Lap-has-finite-propagation}
 \marginpar{Generalize this to banded matrices}
  If each vertex of \Graph has finite degree, then \Lap has finite propagation. In other words, $v \in \spn\{\gd_x\}$ implies $\Lap v \in \spn\{\gd_x\}$.
  \begin{proof}
    Denote the set of vertices which are adjacent to the support of $v$ by
    \linenopax
    \begin{align*}
      nbd(v) := \{x \in \verts \suth x \nbr y \in \spt v\}.
    \end{align*}
    Then for $x \in \verts \less nbd(v)$, the sum defining $\Lap v(x)$ contains only 0s, which implies $\spt \Lap v \ci nbd(v)$.
  \end{proof}
\end{prop}
}

\begin{exm}[Two operators which are each self-adjoint, but whose product is not essentially self-adjoint]
  \label{exm:banded-nonselfadjoint-matrix}
  In \S\ref{sec:Gel'fand-triples-and-duality}, we discussed the Schwartz space $S$ of functions of rapid decay, and its dual $S'$, the space of tempered distributions; cf. \eqref{eqn:Schwartz-space-as-powers-of-Hamiltonian}. If we use the ONB for $L^2(\bR)$ consisting of the Hermite polynomials, then the operators $\tilde P:f(x) \mapsto \frac1\ii \frac{d}{dx}f(x)$ and $\tilde Q:f(x) \mapsto xf(x)$ have the following matrix form:
  \linenopax
  \begin{align}
    P &= \frac12
    \left[\begin{array}{ccccccc}
      0 &1  \\
      1 &0 &\sqrt2 \\
       &\sqrt2 &0 &\sqrt3 \\
       & &\sqrt3 &\ddots &\ddots \\
       & & &\ddots & 0 &\sqrt n \\
       & & & &\sqrt n &0 &\ddots \\
       & & & & &\ddots &\ddots\\
    \end{array}\right],
    \label{eqn:exm:Heisenberg's-P} \\
    Q &= \frac1{2\ii}
    \left[\begin{array}{rrrrrrr}
      0 &1  \\
      -1 &0 &\sqrt2 \\
       &-\sqrt2 &0 &\sqrt3 \\
       & &-\sqrt3 &\ddots &\ddots \\
       & & &\ddots & 0 &\sqrt n \\
       & & & &-\sqrt n &0 &\ddots \\
       & & & & &\ddots &\ddots\\
    \end{array}\right].
    \label{eqn:exm:Heisenberg's-Q}
  \end{align}
  $P$ and $Q$ are Heisenberg's matrices, and they satisfy the canonical commutation relation $PQ-QP = \frac{\ii}{2}\id$. $P$ and $Q$ provide examples of Hermitian operators on $\ell^2(\bZ)$ which are each essentially self-adjoint, but for which $T=QPQ$ is not essentially self-adjoint. In fact, $T$ has defect indices 1,1 (cf.~Definitions~\ref{def:defect-index} and \ref{def:essentially-self-adjoint}). These indices are found directly by solving the the defect equation
  \linenopax
  \begin{align*}
    T^\ad f = QPQ f = \pm \ii f
    \q\implies\q
    x(x f)' = \pm f
  \end{align*}
  to obtain the $C^\iy$ solutions
  \begin{align*}
    f_+(x) =
    \begin{cases}
      \frac{e^{-1/x}}x, & x > 0,\\
      0, &x \leq 0,
    \end{cases}
    \qq
    f_-(x) =
    \begin{cases}
      \frac{e^{1/x}}x, & x < 0,\\
      0, &x \geq 0.
    \end{cases}
  \end{align*}
  Thus there is a 1-\dimnl space of solutions to each defect equation, and the defect indices are 1,1. To see that $P$ and $Q$ are actually self-adjoint, one can observe that $P$ generates the unitary group $f(x) \mapsto f(x+t)$ and $Q$ generates the unitary group $e^{\ii tx}$. Therefore, $P$ and $Q$ are self-adjoint by Stone's theorem; see \cite{DuSc88}.
\end{exm}

\section{Remarks and references}
\label{sec:Remarks-and-References-App-B}

Representations for unbounded operators by infinite matrices were suggested early in quantum mechanics by Heisenberg. Since then, they have served as sources of other applications, as well as the theory of operators in Hilbert space; see especially \cite{vN55}. 
The reader may find the references \cite{Rud91, KadisonRingroseI, DuSc88, ReedSimonII, Arveson:spectral-theory, vN55} to be helpful.

%% file: road-map.tex

\chapter{Navigation aids for operators and spaces}
\label{sec:navigation-aids}

\section{A road map of spaces}
\label{sec:road-map}

Each arrow represents an embedding.
  \begin{equation*}
    \xymatrix{
    \spn\{\gd_x\} \ar[r] & \ell^2(\cond) \ar[r] & \Fin \ar[d] \\
    \spn\{v_x\} \ar[r]   & \Schw \ar[r] & \HE \ar[r] & \Schw'\\
        &   & \Harm \ar[u] \ar[r] & \frac{\Schw'}{\Fin} \ar[u]
    }
  \end{equation*}
  
\pgap
  
\section{A summary of the operators on various Hilbert spaces}
\label{sec:summary-of-Lap-and-Trans}
\hfill \\

\begin{tabular}{r|c|c|c}
   & \cond unbdd & \cond bdd & $\cond=\one$ \\ \hline
  \Lap on $\HE$ & unbdd, Herm, poss. defect & unbdd, Herm, ess. s.-a. & unbdd, Herm, ess. s.-a. \\
  \Lap on $\ell^2(\one)$ & unbdd, Herm, ess. s.-a. & & bdd, s.-a. \\
  \Lap on $\ell^2(\cond)$ & unbdd, non-Herm & non-Herm & bdd, s.-a.\\
  \Trans on $\HE$ & unbdd, Herm, poss. defect & unbdd, Herm, ess. s.-a. & unbdd, Herm, ess. s.-a. \\
  \Trans on $\ell^2(\one)$ & unbdd, Herm, ess. s.-a. & & bdd, s.-a. \\
  \Trans on $\ell^2(\cond)$ & unbdd, non-Herm & non-Herm & bdd, s.-a.\\
\end{tabular}

%% file: biblioguide.tex

\chapter{A guide to the bibliography}
\label{sec:biblio-guide}

We have endeavored to offer a self-contained presentation, but readers looking to brush up on preliminaries may find a number of terms and ideas used here in \cite{AlFi09}, \cite{DoSn84}, \cite{LevPerWil08} and \cite{YiZh05}. These sources have an primarily probabilistic point of view, but resistance networks (weighted graphs) are essentially equivalent to reversible Markov chains; see the quote by Peres at the start of Chapter~\ref{sec:currents-and-potentials}. Consequently, these books cover the fundamentals on electrical resistance networks, as well as such related topics as Reversible Markov Chains; Hitting and Convergence Time, and Flow Rate, Parameters; Special Graphs and Trees; Symmetric Graphs and Chains; $L^2$ Techniques for Bounding Mixing Times; Randomized Algorithms; Continuous State, Infinite State and Random Environment; Interacting Particles; and Markov Chain Monte Carlo.

\begin{enumerate}
  \item \label{first}
  Analysis on weighted graphs, including electrical resistance networks (ERN).
    \begin{enumerate}[\ref{first}.1]
      \item Metrics: R. Alexander, Y. Bartal, D. Burago, I. Benjamini, Y. Colin de Verdier, G. Sved, C. M. Cramlet, J. Dodziuk, P. G. Doyle, J. L. Snell, J. Kigami, R. Lyons, M. A. Rieffel.
      \item Graph theory: G. Polya, B. Bollobas, Fan Chung, R. Diestel, P. G. Doyle, J. L. Snell, C. Thomassen, G. V. Epifanov, M. L. Lapidus, L. Hartman, I. Raeburn.
      \item Potential theory: R. Aharoni, C. Constantinescu, J. Dodziuk, X. W. C. Faber, J. Kigami, V. S. Guliyev, R. Lyons, W. Woess, P. M. Soardi.
      \item Dirichlet forms:  N. Bourleau, F. Hirsch, D. I. Cartwright, W. Woess, J. Dodziuk, M. Fukushima, Y. Oshima, V. A. Kaimanovich, Z. Kuramochi, M. Yamasaki, M. Rockner, B. Schmuland.
      \item Boundaries (Martin and others) and discrete analogues of Greens-Gauss-Stokes theorems: H. Aikawa, S. Brofferio, W. Woess, N. Budarina, C. Z. Chen, R. Diestel, Y. Peres, J. L. Doob, I. Ignatiouk-Robert, Z. Kuramochi, H. L. Royden, S. A. Sawyer, M. Silhavy, P. M. Soardi.
      \item Trees: P. Cartier, R. Permantle, Y. Colin de Verdier, M. D. Esposti, S. Isola, R. Diestel, R. Diestel, R. Lyons, Y. Peres, R. Froese, W. Woess, P. M. Soardi.
      \item Fractals, limits,  and renormalization: A. F. Beardon, M. Barnsley, J. Hutchinson, O. Bratteli, D. Dutkay, L. Baggett, J. Packer, A. Connes, R. L. Devaney, R. Duits, J. de Graaf, J. Kigami, R. S. Strichartz, D. Guido, T. Isola, M. L. Lapidus, M. Golubitsky, B. M. Hambly, T. Kumagai, P. E. T. Jorgensen, E. Pearse, S. Pedersen.
    \end{enumerate}

  \item Tools from the theory of operators in Hilbert space.
    \begin{enumerate}[2.1]
      \item Reproducing kernel Hilbert space: D. Alpay, E. Nelson, N. Aronszajn, D. Larson, T. Kailath.
      \item Unbounded Hermitian operators and their extensions; J. von Neumann,  M. Krein, K. O. Friedrichs, E. Nelson, P. D. Lax, R. S. Phillips, N. Dunford, J. T. Schwartz, M. Reed, B. Simon, P. E. T. Jorgensen.
      \item Shorted operators: W. N.  Anderson, G. E.  Trapp,  C. A.  Butler, T. D.  Morley, V. Metz.
      \item Spectral theory: W. B. Arveson, Y. Colin de Verdier, J. Dodziuk, J. Kigami, R. S. Strichartz, D. Guido, T. Isola, M. L. Lapidus, D. Sarason, P. M. Soardi.
      \item Transfer operators:  V. Baladi, Fan Chung, D. I. Cartwright, W. Woess, L. Saloff-Coste,  D. Dutkay.
      \item Quadratic forms, and perturbation theory: C. Z. Chen, M. Fukushima, Y. Oshima, J. Kigami, T. Kato.
      \item Laplacians: F. Chung, R. M. Richardson, J. Dodziuk, J. L. Doob, X. W. C. Faber, W. Woess, L. G. Rogers, A. Teplyaev, P. M. Soardi, A. Weber.
      \item Fock space methods: T. Hida, W. B.  Arveson, L. Gross, A. Guichardet.
    \end{enumerate}

  \item Stochastic integrals.
    \begin{enumerate}[3.1]
      \item Gelfand triples: I. M. Gelfand, N. Wiener, R. A. Minlos, L. Gross, A. Guichardet, T. Hida.
      \item Schoenberg-von Neumann embeddings: I. J. Schoenberg, J. von Neumann,  B. O. Koopman, J. L. Doob, L. Gross, K. R. Parthasarathy, K. Schmidt.
    \end{enumerate}

  \item Probabilistic methods.
    \begin{enumerate}[4.1]
      \item Random walk models, Markov processes: D. Aldous, J. A. Fill, G. Polya, F. Spitzer, D. W. Stroock, C. A. Nash-Williams, A. I.  Aptekarev, Fan Chung, Y. Peres, J. L. Doob, R. Lyons, T. Lyons, Y. Peres, A. Telcs, R. Pemantle, I. Ignatiouk-Robert, B. Morris, M. A. Picardello, L. Saloff-Coste, G. G. Yin, Q. Zhang.
      \item Potential theory: M. Brelot, R. Lyons, W. Woess, P. Malliavin, A. N. Kolmogorov.
      \item Martingale techniques: J. L. Doob, E. Nelson, W. Woess.
      \item Diffusion models: R. R. Coifman, I. G. Kevrekidis, M. Maggioni, R. Lyons, Y. Peres, W. Woess.
    \end{enumerate}

  \item Tools from other areas of mathematics.
    \begin{enumerate}[5.1]
      \item $C^\ast$-algebras: E. Andruchow, W. B.  Arveson, O. Bratteli, D. W. Robinson, R. V. Kadison, J. A. Ball, I. Cho, A. Connes, P. S. Muhly.
      \item Harmonic analysis: C. Berg, P. Cartier, J. Dodziuk, J. L. Doob, R. S. Strichartz, F. P. Greenleaf, G. Olaffson, W. Woess, A. M. Kytmanov, P. M. Soardi.
      \item Functional analysis: L. Schwartz, M. A.  Al-Gwaiz, H. Aikawa, D. L. Cohn, S. G. Krantz, R. R. Phelps, W. Rudin.
      \item Numerical analysis: K. Atkinson, S. C. Brenner, H. H. Goldstine.
    \end{enumerate}

  \item Applications.
    \begin{enumerate}[6.1]
      \item Classical statistical models: G. Polya, T. M. Liggett, Y. Peres, K. Handa, L. Hartman.
      \item Quantum statistical mechanics. Equilibrium states and long-range order: R. T. Powers, H. Schwetlick, E. Nelson.
      \item Engineering: R. Bott, R. J. Duffin, P. G. Doyle, J. L. Snell, T. Kailath.
      \item Physics (other): M. Baake, R. V.  Moody, W. Dur, H. J. Briegel, A. Dhar, E. Formenti, O. Guhne, J. Kellendonk, D. Ruelle.
    \end{enumerate}
\end{enumerate}